\RequirePackage{fix-cm}
\documentclass[a4paper,10pt]{amsart}
\usepackage{amsmath,amssymb,amsthm,mathtools,amsfonts,amscd,comment,enumitem,etoolbox,textcomp,gensymb,mathdots,stmaryrd,accents}
\usepackage{array,booktabs}
\usepackage{longtable}
\usepackage{tikz}
\usepackage{xcolor}

\usepackage[nopatch=footnote]{microtype}
\usepackage{multicol}
\usepackage{verbatim}
\usepackage{hyperref}
\usepackage{cleveref}
\usepackage{graphicx}
\ifdefined\XeTeXversion
  \usepackage{dsfont}
  \let\mathbbm\mathds %
\else
  \usepackage{bbm} %
\fi

\theoremstyle{definition}

\def\del{\varsigma}
\DeclareMathOperator{\im}{im}
\newcommand{\cM}{\mathcal M}
\newcommand{\cO}{\mathcal O}
\newcommand{\Zloc}{Z_{\lambda,V}}
\newcommand{\eps}{\varepsilon}
\newcommand{\Id}{\operatorname{Id}}
\newcommand{\ev}{\operatorname{ev}}
\newcommand{\coev}{\operatorname{coev}}
\newcommand{\fg}{\mathrm{fg}}

\newcommand\one{\mathbf{1}}
\usetikzlibrary{calc, intersections, arrows.meta, decorations.markings}
\usepackage{amstext,mathrsfs}

\usepackage{pdflscape}
\usepackage{setspace}%
\usepackage[all]{xy}%
\usepackage{indentfirst}%
\usepackage{ytableau}
\usepackage{color}%
\usepackage{colortbl}
\usepackage{dsfont}

\DeclareFontFamily{U}{stmry}{}
\DeclareFontShape{U}{stmry}{m}{n}{%
  <-6>stmary5 <6-7>stmary6 <7-8>stmary7
  <8-9>stmary8 <9-10>stmary9 <10->stmary10}{}
\SetSymbolFont{stmry}{bold}{U}{stmry}{m}{n}
\DeclareFontFamily{U}{rsfs}{\skewchar\font127 }
\DeclareFontShape{U}{rsfs}{m}{n}{%
  <-7>rsfs5 <7-8>rsfs7 <8->rsfs10}{}
\usetikzlibrary{shapes.geometric, arrows}
\usetikzlibrary{calc}
\def\arraystretch{1.2} %
\numberwithin{equation}{section}
 \definecolor{darkred}{HTML}{993333}
\definecolor{darkblue}{RGB}{20,40,120}
\newcommand{\clr}{darkblue}

\definecolor{GRSWRed}{RGB}{174,48,48}
\definecolor{diagramred}{RGB}{174,48,48}
\tikzset{
  grsw strand/.style={line width=.78pt,line cap=round,line join=round},
  grsw red strand/.style={draw=GRSWRed,line width=1.05pt,
    line cap=round,line join=round},
  grsw morphism/.style={-,densely dashed,line width=.62pt,>=stealth},
  grsw morphism reversed/.style={-,densely dashed,line width=.62pt,>=stealth},
  grsw endpoint/.style={circle,draw,line width=.66pt,fill=white,
    inner sep=0pt,minimum size=3.6pt},
  grsw minus/.style={circle,draw,line width=.66pt,fill=white,
    inner sep=.2pt,minimum size=5pt},
  every picture/.append style={>=stealth,line cap=round,line join=round},
  thick/.style={line width=.78pt},
  very thick/.style={line width=1.05pt}
}

\newcommand*\getscale[1]{%
  \begingroup
    \pgfgettransformentries{\scaleA}{\scaleB}{\scaleC}{\scaleD}{\whatevs}{\whatevs}%
    \pgfmathsetmacro{#1}{sqrt(abs(\scaleA*\scaleD-\scaleB*\scaleC))}%
    \expandafter
  \endgroup
  \expandafter\def\expandafter#1\expandafter{#1}%
}

\def\dim{\text{dim}}
\def\Ob{\text{Ob}}
\def\Vect{\text{Vect}}

\def\op{\text{op}}
\def\Id{\text{Id}}
\def\Z{\mathbb{Z}}
\def\N{\mathbb{N}}
\def\id{\text{id}}
\def\Obj{\text{Obj}}
\def\R{\mathbb{R}}
\def\k{\Bbbk}
\def\P{\mathcal{P}}
\def\pcl{P^+_{cl,k}}
\def\ssum{\textstyle\sum\limits}
\def\dim{\operatorname{dim}}
\def\Hom{\operatorname{Hom}}
\def\End{\operatorname{End}}
\def\max{\operatorname{max}}
\def\ev{\operatorname{ev}}
\def\deg{\operatorname{deg}}
\def\res{\operatorname{res}}
\def\ch{\operatorname{char}}

\newtheorem{Theorem}{Theorem}[section] %
\newtheorem{theorem*}{Theorem}
\newtheorem{Lemma}[Theorem]{Lemma}
\newtheorem{Cor}[Theorem]{Corollary}
\newtheorem{Prop}[Theorem]{Proposition}

\newtheorem{Defn}[Theorem]{Definition}

\newtheorem{rem}[Theorem]{Remark}
\newtheorem{Assumption}[Theorem]{Assumption}
\theoremstyle{definition}
\newtheorem{THEOREM}{Theorem}
\def\theTHEOREM{\Alph{THEOREM}}

\newcommand{\rclr}{rgb:red,4}
\newcommand{\wdot}{ node[circle, draw, fill=white, thick, inner sep=0pt, minimum width=4pt]{}}
\newcommand{\bdot}{ node[circle, draw, fill=\clr, thick, inner sep=0pt, minimum width=4pt]{}}
\newcommand{\rdot}{ node[circle, draw=darkred, fill=darkred, thick, inner sep=0pt, minimum width=4pt]{}}
\newcommand{\bd}[1]{\mathbf{#1}}
\newcommand{\ob}[1]{\mathsf{#1}}
\newcommand{\up}{\uparrow}
\newcommand{\down}{\downarrow}
\newcommand{\sdim}{\text{sdim}}
\newcommand{\C}{\mathcal{C}}
\newcommand{\OBC}{\mathcal{OBC}}
\newcommand{\AOB}{\mathcal{AOB}}
\newcommand{\COB}{\mathcal{COB}}
\newcommand{\AOBC}{\mathcal{AOBC}}
\newcommand{\unit}{\mathds{1}}
\newcommand{\OB}{\mathcal{OB}}
\newcommand{\B}{\mathcal{B}}
\newcommand{\CB}{\mathcal{CB}}
\newcommand{\AB}{\mathcal{AB}}
\newcommand{\fgl}{\mathfrak{gl}}
\newcommand{\fq}{\mathfrak{q}}
\newcommand{\fh}{\mathfrak{h}}
\newcommand{\fn}{\mathfrak{n}}
\newcommand{\fb}{\mathfrak{b}}
\newcommand{\wrd}{\langle\up,\down\rangle}
\definecolor{darkblue}{HTML}{111199}
\newcommand{\fulldot}{
    \begin{tikzpicture}[color=\clr]
    \draw (0,0) \bdot;
    \end{tikzpicture}
}
\newcommand{\emptydot}{
    \begin{tikzpicture}[color=\clr]
    \draw (0,0) \wdot;
    \end{tikzpicture}
}
\newcommand{\undot}[1]{\operatorname{undot}({#1})}
\newcommand{\p}[1]{|{#1}|}

\def\clock{\begin{tikzpicture}[baseline=-.9mm]
\filldraw[white] (0,0) circle (1.72mm);
\draw[-] (0,-0.18) to[out=180,in=-102] (-.178,0.02);
\draw[-] (-0.18,0) to[out=90,in=180] (0,0.18);
\draw[-] (0.18,0) to[out=-90,in=0] (0,-0.18);
\draw[-] (0,0.18) to[out=0,in=90] (0.18,0);
\end{tikzpicture}\,}
\newcommand{\cE}{\mathsf{E}}
\newcommand{\cF}{\mathsf{F}}
\newcommand{\excise}[1]{}
\newcommand{\thal}{{}^{\theta}\!\alpha}
\newcommand{\So}{\mathsf{S}}
\def\iLam{\Lambda^{\!\jmath}}
\def\iH{\H^{\jmath}}
\tikzset{anchorbase/.style={>=To,baseline={([yshift=-0.5ex]current bounding box.center)}}}
\tikzset{ %
    centerzero/.style={>=To,baseline={([yshift=-0.5ex]#1)}},
    centerzero/.default={0,0}
}
\tikzset{wipe/.style={white,line width=4pt}}
\tikzset{rotateup/.style={anchor=south, rotate=90, inner sep=.8mm}}
\tikzset{rotatedown/.style={anchor=north, rotate=-90, inner sep=.8mm}}

\tikzset{gcolor/.style={green!60!black,semithick}}
\tikzset{pinhead/.style={thin,gray,fill=white!40!white}}
\tikzset{iQ/.style={semithick,black}}
\tikzset{bulb/.style={thin,black,fill=black}}
\tikzset{Fakebubble/.style={thin,black,fill=white!20!white}}
\tikzset{internalbubble/.style={thin,black,fill=blue!30!white}}
\tikzset{Internalbubble/.style={thin,black,fill=green}}
\tikzset{Q/.style={thin,black}}
\tikzset{Bulb/.style={thin,black,fill=white}}
\tikzset{fakebubble/.style={semithick,black,fill=white!20!white}}
\newcommand{\Zge}{\mathbb Z_{\ge0}}

\newcommand{\Pin}[3]{
    \path (#1) node[inner sep=1.6pt] (x) {} to (#2) node[rectangle,rounded corners,draw,pinhead,inner sep=2.5pt](y) {$\color{black}\scriptstyle#3$};
    \draw[Triangle Cap-,thick,gray] (x)--(y);
    \opendot{#1}
}

\newcommand{\Pinpin}[4]{
\path (#1) node[inner sep=1.6pt] (x) {} to (#3) node[rectangle,rounded corners,draw,pinhead,inner sep=2.5pt](y) {$\color{black}\scriptstyle#4$};
    \draw[Triangle Cap-,thick,gray] (x)--(#2)--(y);
\opendot{#2}; \opendot{#1}
}

\newcommand{\Pinpinpin}[5]{
\path (#1) node[inner sep=1.6pt] (x) {} to (#4) node[rectangle,rounded corners,draw,pinhead,inner sep=2.5pt](y) {$\color{black}\scriptstyle#5$};
    \draw[Triangle Cap-,thick,gray] (x)--(#2)--(#3)--(y);
\opendot{#3}; \opendot{#2}; \opendot{#1}
}

\newcommand{\Pinpinpinmid}[5]{
\path (#1) node[inner sep=1pt] (x) {} to (#3) node[inner sep=1pt] (z) {} to (#4) node[rectangle,rounded corners,draw,pinhead,inner sep=2.5pt](y) {$\color{black}\scriptstyle#5$};
\draw[-,thick,gray] (#2)--(y);
\draw[Triangle Cap-,thick,gray] (x)--(#2);
\draw[Triangle Cap-,thick,gray] (z)--(#2);
\opendot{#3}; \opendot{#2}; \opendot{#1}
}

\newcommand{\pin}[3]{
    \path (#1) node[inner sep=1.6pt] (x) {} to (#2) node[rectangle,rounded corners,draw,pinhead,inner sep=2.5pt](y) {$\color{black}\scriptstyle#3$};
    \draw[Triangle Cap-,thick,gray] (x)--(y);
    \closeddot{#1}
}

\newcommand{\pinpin}[4]{
\path (#1) node[inner sep=1.6pt] (x) {} to (#3) node[rectangle,rounded corners,draw,pinhead,inner sep=2.5pt](y) {$\color{black}\scriptstyle#4$};
    \draw[Triangle Cap-,thick,gray] (x)--(#2)--(y);
\closeddot{#2}; \closeddot{#1}
}

\newcommand{\pinpinpin}[5]{
\path (#1) node[inner sep=1.6pt] (x) {} to (#4) node[rectangle,rounded corners,draw,pinhead,inner sep=2.5pt](y) {$\color{black}\scriptstyle#5$};
    \draw[Triangle Cap-,thick,gray] (x)--(#2)--(#3)--(y);
\closeddot{#3}; \closeddot{#2}; \closeddot{#1}
}

\newcommand{\pinpinpinmid}[5]{
\path (#1) node[inner sep=1pt] (x) {} to (#3) node[inner sep=1pt] (z) {} to (#4) node[rectangle,rounded corners,draw,pinhead,inner sep=2.5pt](y) {$\color{black}\scriptstyle#5$};
\draw[-,thick,gray] (#2)--(y);
\draw[Triangle Cap-,thick,gray] (x)--(#2);
\draw[Triangle Cap-,thick,gray] (z)--(#2);
\closeddot{#3}; \closeddot{#2}; \closeddot{#1}
}

\newcommand{\circledbar}[2]{
\getscale{\scalefactor};
\node[circle,draw,bulb,inner sep=.85pt] at (#1) {$\hspace{0.06pt}\color{white}\scriptscriptstyle \pmb#2$};\draw[line width=.8pt,white] ($(#1)+({-.06/ \scalefactor},{.07 / \scalefactor})$) to ($(#1)+({.06/ \scalefactor},{.07 / \scalefactor})$)}

\newcommand{\circled}[2]{\node[circle,draw,bulb,inner sep=.85pt] at (#1) {$\hspace{0.06pt}\color{white}\scriptscriptstyle\pmb #2$}}

\newcommand{\cross}[2][45]{
\getscale{\scalefactor};
\draw[line width=1.3pt] ($(#2)+({-.12*cos(#1)/\scalefactor},{-.12*sin(#1)/\scalefactor})$) -- ($(#2)+({.12*cos(#1)/\scalefactor},{.12*sin(#1)/\scalefactor})$);
\draw[line width=1.3pt] ($(#2)+({-.12*sin(#1)/\scalefactor},{.12*cos(#1)/\scalefactor})$) -- ($(#2)+({.12*sin(#1)/\scalefactor},{-.12*cos(#1)/\scalefactor})$);
}

\newcommand\multreddot[3]{
    \closeddot[lightgray]{#1};
    \draw (#1) node[anchor=#2] {$\color{teal!30!white}\scriptstyle{#3}$}
}

\newcommand\braidup{to[out=up,in=down]}
\newcommand\braiddown{to[out=down,in=up]}

\colorlet{darkg}{green!50!black}

\newcommand\objlabel[1]{$\color{black}\scriptstyle{#1}$}
\newcommand\dotlabel[1]{$\color{teal}\scriptstyle{#1}$}
\newcommand\regionlabel[1]{$\color{gray}\scriptstyle{#1}$}
\newcommand\strandlabel[1]{$\color{olive}\scriptstyle{#1}$}
\newcommand\stringlabel[2]{\node at (#1) {\strandlabel{#2}}}
\newcommand\thickness[2]{\node at (#1) {$\color{blue}\scriptstyle{#2}$}}
\newcommand\catlabel[2]{\node at (#1) {$\color{darkg}\scriptstyle{#2}$}}
\newcommand\botlabel[1]{node[inner sep=0.5pt,anchor=north] {\strandlabel{#1}}}
\newcommand\toplabel[1]{node[anchor=south,inner sep=0.5pt] {\strandlabel{#1}}}
\newcommand\dottybubblelabel[2]{\node at (#1) {$\color{teal}\scriptstyle{#2}$}}
\newcommand\bubblelabel[2]{\node at (#1) {$\scriptstyle{#2}$}}
\newcommand\region[2]{\node at (#1) {\regionlabel{#2}}}
\newcommand\strand[2]{\node at (#1) {\strandlabel{#2}}}

\newcommand\closeddot[2][black]{
    \node at (#2) {\color{#1}$\scriptstyle{\bullet}$}
}

\newcommand\opendot[2][black]{
    \node at (#2) {\color{white} $\scriptstyle{\bullet}$};\node at (#2) {\color{#1}$\scriptstyle{\circ}$}
}

\newcommand\multcloseddot[3]{
    \closeddot{#1};
    \draw (#1) node[anchor=#2] {\dotlabel{#3}}
}

\newcommand\multopendot[3]{
    \opendot{#1};
    \draw (#1) node[anchor=#2] {\dotlabel{#3}}
}

\newcommand{\teleporter}[3][black]{
\draw[Triangle Cap-Triangle Cap,ultra thick,blue!40!white] (#2)--(#3);
\opendot{#2}; \opendot{#3}}

\newcommand{\Teleporter}[3][black]{
\draw[Triangle Cap-Triangle Cap,ultra thick,green] (#2)--(#3);
\opendot{#2}; \opendot{#3}}

\newcommand{\bentteleporter}[4][black]{
\draw[Triangle Cap-Triangle Cap,ultra thick,blue!40!white] (#2)..controls#4..(#3);
\opendot{#2}; \opendot{#3}
}

\newcommand\bubble[2][black]{%
\draw[-,thin,#1] (#2)++(0,.2) arc(90:-270:0.2);
}
\newcommand\clockwisebubble[2][black]{%
\draw[-to,thin,#1] (#2)++(0,.2) arc(90:-270:0.2);
}
\newcommand\anticlockwisebubble[2][black]{%
  \draw[-to,thin,#1] (#2)++(0,.2) arc(90:450:0.2);
}
\newcommand\filledclockwisebubble[2][black]{%
\draw[-to,#1,Fakebubble] (#2)++(0,.2) arc(90:-270:0.2);
}
\newcommand\filledanticlockwisebubble[2][black]{%
  \draw[-to,#1,Fakebubble] (#2)++(0,.2) arc(90:450:0.2);
}
\newcommand{\M}{\mathcal M}
\newcommand{\wt}{\operatorname{wt}}
\newcommand{\ZZ}{\mathbb Z}
\newcommand{\NN}{\mathbb N}
\newcommand{\Iv}{I^{\nu}}
\newcommand{\Xv}{X_{\nu}}
\newcommand{\Yv}{Y^{\nu}}
\newcommand{\Coeff}[2]{\left[#1\right]_{#2}}

\newcommand{\Bzero}[3]{B^{(#3)}_{#1,#2}}
\newcommand{\Bplus}[2]{B^{+}_{#1,#2}}
\newcommand{\Bminus}[2]{B^{-}_{#1,#2}}
\newcommand{\Om}[3]{\Omega^{(#3)}_{#1,#2}}

\newcommand\anticlockwiseinternalbubble[2][black]{%
 \getscale{\scalefactor};
 \draw[-to,#1,internalbubble] ($(#2)+({-.12 / \scalefactor},0)$)  arc(-180:180:{.12/\scalefactor});
}
\newcommand\clockwiseinternalbubble[2][black]{%
 \getscale{\scalefactor};
 \draw[-to,#1,internalbubble] ($(#2)+({.12 / \scalefactor},0)$)  arc(0:-360:{.12/\scalefactor});
}

\newcommand\babyanticlockwiseinternalbubble[2][black]{%
 \getscale{\scalefactor};
 \draw[-to,#1,internalbubble] ($(#2)+({-.1 / \scalefactor},0)$)  arc(-180:180:{.1/\scalefactor});
}
\newcommand\babyclockwiseinternalbubble[2][black]{%
 \getscale{\scalefactor};
 \draw[-to,#1,internalbubble] ($(#2)+({.1 / \scalefactor},0)$)  arc(0:-360:{.1/\scalefactor});
}

\newcommand\Anticlockwiseinternalbubble[2][black]{%
 \getscale{\scalefactor};
 \draw[-to,#1,Internalbubble] ($(#2)+({-.12 / \scalefactor},0)$)  arc(-180:180:{.12/\scalefactor});
}
\newcommand\Clockwiseinternalbubble[2][black]{%
 \getscale{\scalefactor};
 \draw[-to,#1,Internalbubble] ($(#2)+({.12 / \scalefactor},0)$)  arc(0:-360:{.12/\scalefactor});
}

\newcommand\Babyanticlockwiseinternalbubble[2][black]{%
 \getscale{\scalefactor};
 \draw[-to,#1,Internalbubble] ($(#2)+({-.1 / \scalefactor},0)$)  arc(-180:180:{.1/\scalefactor});
}
\newcommand\Babyclockwiseinternalbubble[2][black]{%
 \getscale{\scalefactor};
 \draw[-to,#1,Internalbubble] ($(#2)+({.1 / \scalefactor},0)$)  arc(0:-360:{.1/\scalefactor});
}

\renewcommand*{\thealphatheorem}{\Alph{alphatheorem}}

\theoremstyle{definition}

\numberwithin{equation}{section}
\allowdisplaybreaks

\newcommand{\lcap}{
\begin{tikzpicture}[baseline = 3pt, scale=0.5]
        \draw[-,thick] (1,0) to[out=up, in=right] (0.53,0.5) to[out=left, in=right] (0.47,0.5);
        \draw[-,thick] (0.49,0.5) to[out=left,in=up] (0,0);
\end{tikzpicture}
}
\newcommand{\lcup}{
\begin{tikzpicture}[baseline = 6pt, scale=0.5]
        \draw[-,thick] (1,1) to[out=down, in=right] (0.53,0.5) to[out=left, in=right] (0.47,0.5);
        \draw[-,thick] (0.49,0.5) to[out=left,in=down] (0,1);
\end{tikzpicture}
}
\newcommand{\lcapl}{
\begin{tikzpicture}[baseline = 3pt, scale=0.5]
        \draw[<-,thick] (1,0) to[out=up, in=right] (0.53,0.5);
         \draw[-,thick]  (0.53,0.5) to[out=left, in=right] (0.47,0.5);
        \draw[-,thick] (0.49,0.5) to[out=left,in=up] (0,0);
\end{tikzpicture}
}
\newcommand{\lcupl}{
\begin{tikzpicture}[baseline = 6pt, scale=0.5]
        \draw[<-,thick] (1,1) to[out=down, in=right] (0.53,0.5);
        \draw[-,thick](0.53,0.5) to[out=left, in=right] (0.47,0.5);
        \draw[-,thick] (0.49,0.5) to[out=left,in=down] (0,1);
\end{tikzpicture}
}
\newcommand{\lcuplL}{
\begin{tikzpicture}[baseline = 6pt, scale=0.5]
        \draw[<-,thick] (1,1) to[out=down, in=right] (0.53,0.5);
        \draw[-,thick](0.53,0.5) to[out=left, in=right] (0.47,0.5);
        \node at (-0.15, 0.85) {$\color{darkblue}\scriptstyle\bullet$};
        \draw[-,thick] (0.49,0.5) to[out=left,in=down] (0,1);
\end{tikzpicture}
}
\newcommand{\lcuplr}{
\begin{tikzpicture}[baseline = 6pt, scale=0.5]
        \draw[<-,thick] (1,1) to[out=down, in=right] (0.53,0.5);
        \draw[-,thick](0.53,0.5) to[out=left, in=right] (0.47,0.5);
        \node at (1.15, 0.85) {$\color{darkblue}\scriptstyle\bullet$};
        \draw[-,thick] (0.49,0.5) to[out=left,in=down] (0,1);
\end{tikzpicture}
}
\newcommand{\lcapo}{
\begin{tikzpicture}[baseline = 3pt, scale=0.5]
        \draw[->,thick] (1,0) to[out=up, in=right] (0.53,0.5) to[out=left, in=right] (0.47,0.5);
        \draw[-,thick] (0.49,0.5) to[out=left,in=up] (0,0);
\end{tikzpicture}
}
\newcommand{\lcupo}{
\begin{tikzpicture}[baseline = 6pt, scale=0.5]
        \draw[-,thick] (1,1) to[out=down, in=right] (0.53,0.5) to[out=left, in=right] (0.47,0.5);
        \draw[->,thick] (0.49,0.5) to[out=left,in=down] (0,1);
\end{tikzpicture}
}
\newcommand{\rcap}{
\begin{tikzpicture}[baseline = 3pt, scale=0.5]
        \draw[<-,thick] (1,0) to[out=up, in=right] (0.53,0.5) to[out=left, in=right] (0.47,0.5);
        \draw[-,thick] (0.49,0.5) to[out=left,in=up] (0,0);
\end{tikzpicture}
}
\newcommand{\rcup}{
\begin{tikzpicture}[baseline = 6pt, scale=0.5]
        \draw[<-,thick] (1,1) to[out=down, in=right] (0.53,0.5) to[out=left, in=right] (0.47,0.5);
        \draw[-,thick] (0.49,0.5) to[out=left,in=down] (0,1);
\end{tikzpicture}
}
\newcommand{\swap}{
\begin{tikzpicture}[baseline = 3pt, scale=0.4]
        \draw[-,thick] (0,0) to[out=up, in=down] (1,1);
        \draw[-,thick] (1,0) to[out=up, in=down] (0,1);
\end{tikzpicture}
}
\newcommand{\swapo}{
\begin{tikzpicture}[baseline = 3pt, scale=0.5]
        \draw[->,thick] (0,0) to[out=up, in=down] (1,1);
        \draw[->,thick] (1,0) to[out=up, in=down] (0,1);
\end{tikzpicture}
}

\newcommand{\swapr}{
\begin{tikzpicture}[baseline = 3pt, scale=0.5]
        \draw[<-,thick] (0,0) to[out=up, in=down] (1,1);
        \draw[<-,thick] (1,0) to[out=up, in=down] (0,1);
\end{tikzpicture}
}
\newcommand{\rswap}{
\begin{tikzpicture}[baseline = 3pt, scale=0.5]
        \draw[->,thick] (0,0) to[out=up, in=down] (1,1);
        \draw[<-,thick] (1,0) to[out=up, in=down] (0,1);
\end{tikzpicture}
}

\newcommand{\lswap}{
\begin{tikzpicture}[baseline = 3pt, scale=0.5]
        \draw[<-,thick] (0,0) to[out=up, in=down] (1,1);
        \draw[->,thick] (1,0) to[out=up, in=down] (0,1);
\end{tikzpicture}
}
\newcommand{\dswap}{
\begin{tikzpicture}[baseline = 3pt, scale=0.5]
        \draw[<-,thick] (0,0) to[out=up, in=down] (1,1);
        \draw[<-,thick] (1,0) to[out=up, in=down] (0,1);
\end{tikzpicture}
}
\newcommand{\caup}{
\begin{tikzpicture}[baseline = 5pt, scale=0.75]
                 \draw[<-,thick] (0.75,1) to[out=down,in=left] (1,0.6) to[out=right,in=down] (1.25,1);
          \draw[<-,thick] (1.25,0) to[out=up,in=right] (1,0.4) to[out=left,in=up] (0.75,0);\end{tikzpicture}}

\newcommand{\xd}{
\begin{tikzpicture}[baseline = 3pt, scale=0.5, color=\clr]
\draw[->,thick] (0,0) to[out=up, in=down] (0,1);
\draw(0,0.5) \bdot;
 \draw (0.35,0.4) node{\tiny $2$};
\end{tikzpicture}
}
\newcommand{\xdx}{
\begin{tikzpicture}[baseline = 3pt, scale=0.5, color=\clr]
\draw[<-,thick] (0,0) to[out=up, in=down] (0,1);
\draw(0,0.5) \bdot;
 \draw (0.35,0.4) node{ \tiny$ 2$};\end{tikzpicture}
}
\newcommand{\xdz}{
\begin{tikzpicture}[baseline = 3pt, scale=0.5, color=\clr]
\draw[->,thick] (0,0) to[out=up, in=down] (0,1);
\draw(0,0.5) \bdot;
\end{tikzpicture}
}
\newcommand{\xli}{
\begin{tikzpicture}[baseline = 3pt, scale=0.5, color=\clr]
        \draw[<-,thick] (0,0) to[out=up, in=down] (0,1);
\end{tikzpicture}
}
\newcommand{\sli}{
\begin{tikzpicture}[baseline = 3pt, scale=0.5, color=\clr]
        \draw[->,thick] (0,0) to[out=up, in=down] (0,1);
\end{tikzpicture}
}
\newcommand{\xdxz}{
\begin{tikzpicture}[baseline = 3pt, scale=0.5, color=\clr]
\draw[<-,thick] (0,0) to[out=up, in=down] (0,1);
\draw(0,0.5) \bdot;
 \end{tikzpicture}
}

\newcommand{\obg}{\begin{tikzpicture}[baseline = 10pt, scale=0.5, color=\clr] \draw[-,thick] (0,0.5)to[out=up,in=down](0,1.2);
    \end{tikzpicture}
    }
\def\bigdt{{\color{white}\bullet}\hspace{.1mm}\!\!\!\pmb\circ}
\def\dt{{\color{white}\bullet}\!\!\!\circ}
\def\bull{{\scriptstyle\bullet}}

\newcommand{\xdotk}{
\begin{tikzpicture}[baseline = 3pt, scale=0.5, color=\clr]
        \draw[-,thick] (0,0) to[out=up, in=down] (0,1);
            \draw (0,0.4)\bdot;
        \draw (0.5,0.4) node{$k$};
\end{tikzpicture}
}
\newcommand{\xdot}{\begin{tikzpicture}[baseline = -1mm, color=\clr]
	\draw[-] (0.08,-.3) to (0.08,.4);
      \node at (0.08,0.15) {$\dt$};
\end{tikzpicture}}
\newcommand{\cev}[1]{\reflectbox{\ensuremath{\vec{\reflectbox{\ensuremath{#1}}}}}}
\newcommand{\xdota}{
\begin{tikzpicture}[baseline = 3pt, scale=0.5, color=\clr]
        \draw[-,thick] (0,0) to[out=up, in=down] (0,1);
            \draw (0,0.4)\bdot;
        \draw (0.5,0.4) node{$a$};
\end{tikzpicture}
}
\tikzset{strand/.style={-,thick}}
\newcommand{\hf}{\diamond}
\newcommand{\darkg}{\color{green!70!black}}
\tikzset{darkg/.style={green!70!black}}

\tikzset{
	cplus/.style={
		circle,
		draw=black,
		fill=white,
		line width=0.2pt,
		inner sep=0.025pt,
		minimum size=.2pt,
		font=\fontsize{5pt}{6pt}\selectfont
	}
}
\newcommand{\midscript}[1]{{\fontsize{6pt}{7pt}\selectfont$#1$}}

\tikzset{
  blackstrand/.style={black,thick,line cap=round,line join=round},
  redstrand/.style={darkred,thick,line cap=round,line join=round},
  guide/.style={black,densely dotted,thin},
  smalldot/.style={circle,fill=black,inner sep=1.25pt}
}

\newcommand{\BlackCup}{\mathord{\begin{tikzpicture}[baseline=-.55ex,scale=.48]
  \draw[blackstrand] (-.42,.55) to[out=-90,in=180] (0,0)
    to[out=0,in=-90] (.42,.55);
\end{tikzpicture}}}
\newcommand{\RedCup}{\mathord{\begin{tikzpicture}[baseline=-.55ex,scale=.48]
  \draw[redstrand] (-.42,.55) to[out=-90,in=180] (0,0)
    to[out=0,in=-90] (.42,.55);
\end{tikzpicture}}}
\newcommand{\BlackCap}{\mathord{\begin{tikzpicture}[baseline=-.15ex,scale=.48]
  \draw[blackstrand] (-.42,0) to[out=90,in=180] (0,.55)
    to[out=0,in=90] (.42,0);
\end{tikzpicture}}}
\newcommand{\RedCap}{\mathord{\begin{tikzpicture}[baseline=-.15ex,scale=.48]
  \draw[redstrand] (-.42,0) to[out=90,in=180] (0,.55)
    to[out=0,in=90] (.42,0);
\end{tikzpicture}}}
\newcommand{\BlackCross}{\mathord{\begin{tikzpicture}[baseline=-.45ex,scale=.48]
  \draw[blackstrand] (-.4,0)--(.4,.72);
  \draw[blackstrand] (.4,0)--(-.4,.72);
\end{tikzpicture}}}
\newcommand{\RedCross}{\mathord{\begin{tikzpicture}[baseline=-.45ex,scale=.48]
  \draw[redstrand] (-.4,0)--(.4,.72);
  \draw[redstrand] (.4,0)--(-.4,.72);
\end{tikzpicture}}}
\newcommand{\TwoVertical}{\mathord{\begin{tikzpicture}[baseline=-.45ex,scale=.48]
  \draw[blackstrand] (-.25,0)--(-.25,.72);
  \draw[blackstrand] (.25,0)--(.25,.72);
\end{tikzpicture}}}
\newcommand{\RedTwoVertical}{\mathord{\begin{tikzpicture}[baseline=-.45ex,scale=.48]
  \draw[redstrand] (-.25,0)--(-.25,.72);
  \draw[redstrand] (.25,0)--(.25,.72);
\end{tikzpicture}}}
\newcommand{\RedDot}{\mathord{\begin{tikzpicture}[baseline=-.45ex,scale=.48]
  \draw[redstrand] (0,0)--(0,.72);
  \fill[darkred] (0,.36) circle (1.7pt);
\end{tikzpicture}}}
\newcommand{\BlackDot}{\mathord{\begin{tikzpicture}[baseline=-.45ex,scale=.48]
  \draw[blackstrand] (0,0)--(0,.72);
  \fill[black] (0,.36) circle (1.7pt);
\end{tikzpicture}}}
\newcommand{\CupCap}{\mathord{\begin{tikzpicture}[baseline=-.5ex,scale=.42]
  \draw[blackstrand] (-.35,.62) to[out=-90,in=180] (0,.32)
    to[out=0,in=-90] (.35,.62);
  \draw[blackstrand] (-.35,0) to[out=90,in=180] (0,.30)
    to[out=0,in=90] (.35,0);
\end{tikzpicture}}}

\newcommand{\LabeledCross}[2]{\mathord{\begin{tikzpicture}[baseline=-.5ex,scale=.62]
  \draw[blackstrand] (-.48,0)--(.48,.95);
  \draw[blackstrand] (.48,0)--(-.48,.95);
  \draw[guide] (-.25,.72)--(.18,.72);
  \node[anchor=west,scale=.68] at (.18,.72) {$#1$};
  \draw[guide] (-.18,.23)--(.25,.23);
  \node[anchor=west,scale=.68] at (.25,.23) {$#2$};
\end{tikzpicture}}}

\newcommand{\VerticalLabel}[1]{\mathord{\begin{tikzpicture}[baseline=-.5ex,scale=.62]
  \draw[blackstrand] (-.32,0)--(-.32,.95);
  \draw[blackstrand] (.32,0)--(.32,.95);
  \draw[guide] (-.20,.48)--(.20,.48);
  \node[anchor=west,scale=.68] at (.20,.48) {$#1$};
\end{tikzpicture}}}

\newcommand{\CrossDotLeft}{\mathord{\begin{tikzpicture}[baseline=-.45ex,scale=.48]
  \draw[blackstrand] (-.4,0)--(.4,.72);
  \draw[blackstrand] (.4,0)--(-.4,.72);
  \fill[black] (-.18,.20) circle (1.7pt);
\end{tikzpicture}}}
\newcommand{\CrossDotRight}{\mathord{\begin{tikzpicture}[baseline=-.45ex,scale=.48]
  \draw[blackstrand] (-.4,0)--(.4,.72);
  \draw[blackstrand] (.4,0)--(-.4,.72);
  \fill[black] (.18,.20) circle (1.7pt);
\end{tikzpicture}}}

\newcommand{\RedCapCross}{\mathord{\begin{tikzpicture}[baseline=-.5ex,scale=.55]
  \draw[redstrand] (-.42,.76) to[out=90,in=180] (0,1.06)
    to[out=0,in=90] (.42,.76);
  \draw[redstrand] (-.42,.76)--(.42,0);
  \draw[redstrand] (.42,.76)--(-.42,0);
\end{tikzpicture}}}

\newcommand{\BlackCapCross}{\mathord{\begin{tikzpicture}[baseline=-.5ex,scale=.55]
  \draw[blackstrand] (-.42,.76) to[out=90,in=180] (0,1.06)
    to[out=0,in=90] (.42,.76);
  \draw[blackstrand] (-.42,.76)--(.42,0);
  \draw[blackstrand] (.42,.76)--(-.42,0);
\end{tikzpicture}}}

\newcommand{\BlackCurl}{\mathord{\begin{tikzpicture}[baseline=-.55ex,scale=.53]
  \draw[blackstrand] (-.30,0) to[out=90,in=-90] (.25,.80)
    to[out=90,in=0] (-.05,1.10) to[out=180,in=90] (-.34,.78)
    to[out=-90,in=90] (.30,0);
\end{tikzpicture}}}

\newcommand{\BraidLeft}[1]{\mathord{\begin{tikzpicture}[baseline=-.5ex,scale=#1]
  \draw[blackstrand] (.45,.6)--(-.45,-.6);
  \draw[blackstrand] (.45,-.6)--(-.45,.6);
  \draw[blackstrand] (0,-.6) to[out=90,in=-90] (-.45,0)
    to[out=90,in=-90] (0,.6);
\end{tikzpicture}}}
\newcommand{\BraidRight}[1]{\mathord{\begin{tikzpicture}[baseline=-.5ex,scale=#1]
  \draw[blackstrand] (.45,.6)--(-.45,-.6);
  \draw[blackstrand] (.45,-.6)--(-.45,.6);
  \draw[blackstrand] (0,-.6) to[out=90,in=-90] (.45,0)
    to[out=90,in=-90] (0,.6);
\end{tikzpicture}}}
\definecolor{diagramred}{RGB}{205,35,35}
\tikzset{
  strand/.style={line width=1.05pt,line cap=round,line join=round},
  blackstrand/.style={strand,draw=black},
  redstrand/.style={strand,draw=diagramred},
  overblack/.style={blackstrand,preaction={draw=white,line width=3.4pt}},
  overred/.style={redstrand,preaction={draw=white,line width=3.4pt}},
  slide/.style={densely dashed,-{Stealth[length=1.6mm,width=1.2mm]},line width=.65pt},
  lab/.style={font=\scriptsize,inner sep=.7pt,fill=white,fill opacity=.88,text opacity=1}
}

\newcommand{\pt}[3]{\draw[draw=#1,fill=white,line width=.8pt] (#2,#3) circle[radius=.075];}
\newcommand{\dmath}[1]{\mathord{\vcenter{\hbox{#1}}}}

\newcommand{\dCross}[1]{\dmath{\begin{tikzpicture}[x=.72em,y=.72em]
  \draw[#1strand] (-.72,-.95) .. controls (-.32,-.35) and (.32,.35) .. (.72,.95);
  \draw[over#1] (.72,-.95) .. controls (.32,-.35) and (-.32,.35) .. (-.72,.95);
  \pt{#1}{-.72}{-.95}\pt{#1}{.72}{-.95}\pt{#1}{-.72}{.95}\pt{#1}{.72}{.95}
\end{tikzpicture}}}
\newcommand{\BCross}{\dCross{black}}
\newcommand{\RCross}{\dCross{red}}

\newcommand{\dVert}[1]{\dmath{\begin{tikzpicture}[x=.72em,y=.72em]
  \draw[#1strand] (-.48,-.95)--(-.48,.95);\draw[#1strand] (.48,-.95)--(.48,.95);
  \pt{#1}{-.48}{-.95}\pt{#1}{.48}{-.95}\pt{#1}{-.48}{.95}\pt{#1}{.48}{.95}
\end{tikzpicture}}}
\newcommand{\BVert}{\dVert{black}}
\newcommand{\RVert}{\dVert{red}}

\newcommand{\dCup}[1]{\dmath{\begin{tikzpicture}[x=.72em,y=.72em]
  \draw[#1strand] (-.62,.42) .. controls (-.35,-.58) and (.35,-.58) .. (.62,.42);
  \pt{#1}{-.62}{.42}\pt{#1}{.62}{.42}
\end{tikzpicture}}}
\newcommand{\dCap}[1]{\dmath{\begin{tikzpicture}[x=.72em,y=.72em]
  \draw[#1strand] (-.62,-.42) .. controls (-.35,.58) and (.35,.58) .. (.62,-.42);
  \pt{#1}{-.62}{-.42}\pt{#1}{.62}{-.42}
\end{tikzpicture}}}
\newcommand{\BCup}{\dCup{black}}\newcommand{\Bcap}{\dCap{black}}
\newcommand{\RCup}{\dCup{red}}\newcommand{\Rcap}{\dCap{red}}

\newcommand{\dBubble}[1]{\dmath{\begin{tikzpicture}[x=.78em,y=.78em]
  \draw[#1strand] (0,-.75) .. controls (-.9,-.75) and (-.9,.75) .. (0,.75)
                    .. controls (.82,.75) and (.82,-.25) .. (.36,-.25);
  \draw[draw=#1,fill=white,line width=.8pt] (.36,-.25) circle[radius=.075];
\end{tikzpicture}}}
\newcommand{\RBubble}{\dBubble{red}}

\newcommand{\dCapCross}[1]{\dmath{\begin{tikzpicture}[x=.76em,y=.76em]
  \draw[#1strand] (-.62,.45) .. controls (-.38,1.2) and (.38,1.2) .. (.62,.45);
  \draw[#1strand] (-.62,.45) .. controls (-.28,.05) and (.28,-.28) .. (.68,-.92);
  \draw[over#1] (.62,.45) .. controls (.28,.05) and (-.28,-.28) .. (-.68,-.92);
  \pt{#1}{-.62}{.45}\pt{#1}{.62}{.45}\pt{#1}{-.68}{-.92}\pt{#1}{.68}{-.92}
\end{tikzpicture}}}
\newcommand{\dCrossCup}[1]{\dmath{\begin{tikzpicture}[x=.76em,y=.76em]
  \draw[#1strand] (-.68,.92) .. controls (-.28,.28) and (.28,-.05) .. (.62,-.45);
  \draw[over#1] (.68,.92) .. controls (.28,.28) and (-.28,-.05) .. (-.62,-.45);
  \draw[#1strand] (-.62,-.45) .. controls (-.38,-1.2) and (.38,-1.2) .. (.62,-.45);
  \pt{#1}{-.68}{.92}\pt{#1}{.68}{.92}\pt{#1}{-.62}{-.45}\pt{#1}{.62}{-.45}
\end{tikzpicture}}}
\newcommand{\RCapCross}{\dCapCross{red}}
\newcommand{\RCrossCup}{\dCrossCup{red}}

\newcommand{\dCapCrossCup}[1]{\dmath{\begin{tikzpicture}[x=.78em,y=.78em]
  \draw[#1strand] (-.58,.55) .. controls (-.35,1.18) and (.35,1.18) .. (.58,.55);
  \draw[#1strand] (-.58,.55) .. controls (-.24,.16) and (.24,-.16) .. (.58,-.55);
  \draw[over#1] (.58,.55) .. controls (.24,.16) and (-.24,-.16) .. (-.58,-.55);
  \draw[#1strand] (-.58,-.55) .. controls (-.35,-1.18) and (.35,-1.18) .. (.58,-.55);
  \pt{#1}{-.58}{.55}\pt{#1}{.58}{.55}\pt{#1}{-.58}{-.55}\pt{#1}{.58}{-.55}
\end{tikzpicture}}}
\newcommand{\RCapCrossCup}{\dCapCrossCup{red}}

\newcommand{\dDottedCross}[2]{\dmath{\begin{tikzpicture}[x=.72em,y=.72em]
  \draw[#1strand] (-.72,-.95) .. controls (-.32,-.35) and (.32,.35) .. (.72,.95);
  \draw[over#1] (.72,-.95) .. controls (.32,-.35) and (-.32,.35) .. (-.72,.95);
  \fill[#1] #2 circle[radius=.085];
  \pt{#1}{-.72}{-.95}\pt{#1}{.72}{-.95}\pt{#1}{-.72}{.95}\pt{#1}{.72}{.95}
\end{tikzpicture}}}
\newcommand{\RDotCrossA}{\dDottedCross{red}{(-.22,.18)}}
\newcommand{\RDotCrossB}{\dDottedCross{red}{(.22,-.18)}}

\newcommand{\dLabCross}[5]{\dmath{\begin{tikzpicture}[x=.78em,y=.78em]
  \draw[blackstrand] (-.74,-1.05) .. controls (-.32,-.36) and (.32,.36) .. (.74,1.05);
  \draw[overblack] (.74,-1.05) .. controls (.32,-.36) and (-.32,.36) .. (-.74,1.05);
  \pt{black}{-.74}{-1.05}\pt{black}{.74}{-1.05}\pt{black}{-.74}{1.05}\pt{black}{.74}{1.05}
  \draw[slide] (-.22,.70)--(-.22,.10);\draw[slide] (.22,.70)--(.22,.10);
  \node[lab,anchor=east] at (-.15,.70) {$#2$};
  \node[lab,anchor=west] at (.15,.78) {$#3$};
  \node[lab,anchor=east] at (-.15,-.72) {$#4$};
  \node[lab,anchor=west] at (.15,-.72) {$#5$};
\end{tikzpicture}}}

\newcommand{\dCoeffCross}[2]{\dmath{\begin{tikzpicture}[x=.78em,y=.78em]
  \draw[blackstrand] (-.74,-1.05) .. controls (-.32,-.36) and (.32,.36) .. (.74,1.05);
  \draw[overblack] (.74,-1.05) .. controls (.32,-.36) and (-.32,.36) .. (-.74,1.05);
  \pt{black}{-.74}{-1.05}\pt{black}{.74}{-1.05}\pt{black}{-.74}{1.05}\pt{black}{.74}{1.05}
  \draw[slide] (-.36,.46)--(.22,.46);\draw[slide] (-.36,-.46)--(.22,-.46);
  \node[lab,anchor=west] at (.18,.55) {$#1$};
  \node[lab,anchor=west] at (.18,-.55) {$#2$};
\end{tikzpicture}}}

\newcommand{\dCoeffVert}[1]{\dmath{\begin{tikzpicture}[x=.78em,y=.78em]
  \draw[blackstrand] (-.46,-1.05)--(-.46,1.05);\draw[blackstrand] (.46,-1.05)--(.46,1.05);
  \pt{black}{-.46}{-1.05}\pt{black}{.46}{-1.05}\pt{black}{-.46}{1.05}\pt{black}{.46}{1.05}
  \draw[slide] (-.28,0)--(.28,0);\node[lab,anchor=west] at (.23,.10) {$#1$};
\end{tikzpicture}}}

\newcommand{\dCupCap}[2]{\dmath{\begin{tikzpicture}[x=.78em,y=.78em]
  \draw[blackstrand] (-.62,.52) .. controls (-.35,.08) and (.35,.08) .. (.62,.52);
  \draw[blackstrand] (-.62,-.52) .. controls (-.35,-.08) and (.35,-.08) .. (.62,-.52);
  \node[lab] at (0,.55) {$#1$};\node[lab] at (0,-.55) {$#2$};
\end{tikzpicture}}}

\newcommand{\dSnake}[1]{\dmath{\begin{tikzpicture}[x=.72em,y=.72em]
  \draw[#1strand] (-.78,-1.05) .. controls (-.92,.15) and (-.78,.92) .. (-.20,.92)
                    .. controls (.32,.92) and (.15,.18) .. (.58,-.02)
                    .. controls (.92,-.23) and (.72,-.92) .. (.15,-.92);
  \draw[over#1] (-.20,-1.05) .. controls (.10,-.52) and (.18,.22) .. (.70,1.05);
  \pt{#1}{-.78}{-1.05}\pt{#1}{-.20}{-1.05}\pt{#1}{.70}{1.05}\pt{#1}{-.78}{1.05}
\end{tikzpicture}}}
\newcommand{\BSnake}{\dSnake{black}}

\newcommand{\coeff}[2]{\bigl[#1\bigr]_{#2}}
\newcommand{\DeltaX}{\Delta}
\newcommand{\OmegaD}{\Omega}

\tikzset{
  zerovtx/.style={font=\scriptsize,inner sep=.15pt,fill=white},
  labelline/.style={densely dotted,thin,-{Stealth[length=2.4pt,width=2.1pt]}},
  diagramlabel/.style={font=\scriptsize,inner sep=.4pt},
  strand/.style={thick,line cap=round,line join=round}
}

\newcommand{\zpt}[4]{%
  \node[zerovtx,text=#1] at (#2,#3) {$\scriptstyle #4$};}

\newcommand{\TwoCross}[1]{\dmath{\begin{tikzpicture}
  [baseline=-.5ex,scale=1.08,transform shape]
  \draw[strand,#1] (.45,.6) -- (-.45,-.6);
  \draw[strand,#1] (.45,-.6) -- (-.45,.6);
  \zpt{#1}{-.45}{-.6}{0}\zpt{#1}{.45}{-.6}{0}
  \zpt{#1}{-.45}{.6}{0}\zpt{#1}{.45}{.6}{0}
\end{tikzpicture}}}

\newcommand{\LabeledVert}[1]{\dmath{\begin{tikzpicture}
  [baseline=-.5ex,scale=1.00,transform shape]
  \draw[strand] (-.45,-.6) -- (-.45,.6);
  \draw[strand] (.45,-.6) -- (.45,.6);
  \draw[labelline] (-.38,0) -- (.38,0);
  \node[diagramlabel,anchor=west] at (.48,0) {$#1$};
  \zpt{black}{-.45}{-.6}{0}\zpt{black}{.45}{-.6}{0}
  \zpt{black}{-.45}{.6}{0}\zpt{black}{.45}{.6}{0}
\end{tikzpicture}}}

\newcommand{\LabeledCupCap}[2]{\dmath{\begin{tikzpicture}
  [baseline=-.5ex,scale=1.05,transform shape]
  \draw[strand] (-.45,.6) to[out=-90,in=180] (0,.2);
  \draw[strand] (0,.2) to[out=0,in=-90] (.45,.6);
  \draw[strand] (-.45,-.6) to[out=90,in=180] (0,-.2);
  \draw[strand] (0,-.2) to[out=0,in=90] (.45,-.6);
  \draw[labelline] (-.40,.48) -- (.02,.48);
  \node[diagramlabel,anchor=east] at (-.48,.48) {$#1$};
  \draw[labelline] (-.02,-.48) -- (.40,-.48);
  \node[diagramlabel,anchor=west] at (.48,-.48) {$#2$};
  \zpt{black}{-.45}{-.6}{0}\zpt{black}{.45}{-.6}{0}
  \zpt{black}{-.45}{.6}{0}\zpt{black}{.45}{.6}{0}
\end{tikzpicture}}}

\newcommand{\EigenPair}[3]{\dmath{\begin{tikzpicture}
  [baseline=-.5ex,scale=1.00,transform shape]
  \draw[strand] (-.45,-.6) -- (-.45,.6);
  \draw[strand] (.45,-.6) -- (.45,.6);
  \draw[labelline] (-.38,0) -- (.38,0);
  \node[diagramlabel,anchor=west] at (.48,0) {$#1$};
  \node[diagramlabel] at (-.45,-.78) {$#2$};
  \node[diagramlabel] at (.45,-.78) {$#3$};
\end{tikzpicture}}}

\newcommand{\localpiece}[4]{%
  \pgfmathsetmacro{\ya}{#4}%
  \pgfmathsetmacro{\yb}{#4+1.20}%
  \ifnum#2=1
    \ifcase#3
      \draw[strand,#1] (-.75,\ya) -- (0,\yb);
      \draw[strand,#1] (0,\ya) -- (-.75,\yb);
    \or
      \draw[strand,#1] (-.75,\ya) -- (-.75,\yb);
      \draw[strand,#1] (0,\ya) -- (0,\yb);
    \or
      \draw[strand,#1] (-.75,\ya)
        to[out=90,in=180] (-.375,\ya+.28)
        to[out=0,in=90] (0,\ya);
      \draw[strand,#1] (-.75,\yb)
        to[out=-90,in=180] (-.375,\yb-.28)
        to[out=0,in=-90] (0,\yb);
    \fi
    \draw[strand,#1] (.75,\ya) -- (.75,\yb);
  \else
    \draw[strand,#1] (-.75,\ya) -- (-.75,\yb);
    \ifcase#3
      \draw[strand,#1] (0,\ya) -- (.75,\yb);
      \draw[strand,#1] (.75,\ya) -- (0,\yb);
    \or
      \draw[strand,#1] (0,\ya) -- (0,\yb);
      \draw[strand,#1] (.75,\ya) -- (.75,\yb);
    \or
      \draw[strand,#1] (0,\ya)
        to[out=90,in=180] (.375,\ya+.28)
        to[out=0,in=90] (.75,\ya);
      \draw[strand,#1] (0,\yb)
        to[out=-90,in=180] (.375,\yb-.28)
        to[out=0,in=-90] (.75,\yb);
    \fi
  \fi}

\newcommand{\ThreeStackL}[4]{\dmath{\begin{tikzpicture}
  [baseline=-.5ex,scale=1.00,transform shape]
  \localpiece{#1}{1}{#2}{-1.80}
  \localpiece{#1}{2}{#3}{-.60}
  \localpiece{#1}{1}{#4}{.60}
  \foreach \x in {-.75,0,.75}{
    \zpt{#1}{\x}{-1.80}{0}\zpt{#1}{\x}{1.80}{0}}
\end{tikzpicture}}}

\newcommand{\ThreeStackR}[4]{\dmath{\begin{tikzpicture}
  [baseline=-.5ex,scale=1.00,transform shape]
  \localpiece{#1}{2}{#2}{-1.80}
  \localpiece{#1}{1}{#3}{-.60}
  \localpiece{#1}{2}{#4}{.60}
  \foreach \x in {-.75,0,.75}{
    \zpt{#1}{\x}{-1.80}{0}\zpt{#1}{\x}{1.80}{0}}
\end{tikzpicture}}}

\newcommand{\RedBraidL}{\ThreeStackL{darkred}{0}{0}{0}}
\newcommand{\RedBraidR}{\ThreeStackR{darkred}{0}{0}{0}}

\newcommand{\ctype}[1]{\ifcase#1\mathsf X\or\mathsf I\or\mathsf C\fi}
\newcommand{\termcell}[3]{%
  \begin{minipage}[t][5.15cm][t]{.30\textwidth}
    \vspace{0pt}\centering
    $\ThreeStackL{black}{#1}{#2}{#3}$\par
    \vspace{.45em}
    {\small
      $\mathcal D_{\ctype{#1}\ctype{#2}\ctype{#3}}$}\par
    \vspace{.15em}
    {\scriptsize
      $(d^{12}_{\ctype{#1}}d^{23}_{\ctype{#2}}
        u^{12}_{\ctype{#3}})$}
  \end{minipage}}

\newcommand{\ThreeCrossVert}{\dmath{\begin{tikzpicture}
  [baseline=-.5ex,scale=1.02,transform shape]
  \draw[strand] (-.75,-.7) -- (0,.7);
  \draw[strand] (0,-.7) -- (-.75,.7);
  \draw[strand] (.75,-.7) -- (.75,.7);
  \foreach \x in {-.75,0,.75}{
    \zpt{black}{\x}{-.7}{0}\zpt{black}{\x}{.7}{0}}
\end{tikzpicture}}}

\newcommand{\ThreeIdentity}{\dmath{\begin{tikzpicture}
  [baseline=-.5ex,scale=1.02,transform shape]
  \foreach \x in {-.75,0,.75}{
    \draw[strand] (\x,-.7) -- (\x,.7);
    \zpt{black}{\x}{-.7}{0}\zpt{black}{\x}{.7}{0}}
\end{tikzpicture}}}

\newcommand{\ThreeCupCap}{\dmath{\begin{tikzpicture}
  [baseline=-.5ex,scale=1.02,transform shape]
  \draw[strand] (-.75,-.7)
    to[out=90,in=180] (-.375,-.32)
    to[out=0,in=90] (0,-.7);
  \draw[strand] (-.75,.7)
    to[out=-90,in=180] (-.375,.32)
    to[out=0,in=-90] (0,.7);
  \draw[strand] (.75,-.7) -- (.75,.7);
  \foreach \x in {-.75,0,.75}{
    \zpt{black}{\x}{-.7}{0}\zpt{black}{\x}{.7}{0}}
\end{tikzpicture}}}

\newcommand{\minusnode}[1]{%
  \node[circle,draw,fill=white,line width=0.45pt,
        inner sep=0.25pt,minimum size=5.2pt]
        at #1 {\ensuremath{\scriptscriptstyle -}};%
}

\newcommand{\leftcomposite}{%
\begin{tikzpicture}[baseline=-0.5ex,scale=0.8,>=stealth]
  \draw[thick]
    (-0.38,1.35)
    .. controls (-0.28,1.22) and (-0.12,1.05) .. (0,0.94)
    .. controls (0.12,0.83) and (0.34,0.56) .. (0.30,0.34)
    .. controls (0.26,0.14) and (0.08,0.10) .. (0,0.10);
  \draw[thick]
    (0.38,1.35)
    .. controls (0.28,1.22) and (0.12,1.05) .. (0,0.94)
    .. controls (-0.12,0.83) and (-0.34,0.56) .. (-0.30,0.34)
    .. controls (-0.26,0.14) and (-0.08,0.10) .. (0,0.10);

  \draw[->,densely dashed] (0,0.02) -- (0,-0.25);

  \draw[thick]
    (-0.38,-1.15)
    .. controls (-0.38,-0.55) and (0.38,-0.55) ..
    (0.38,-1.15);

  \node at (-0.40, 1.54) {\ensuremath{\scriptstyle 0}};
  \node at ( 0.40, 1.54) {\ensuremath{\scriptstyle 0}};
  \node at (-0.40,-1.34) {\ensuremath{\scriptstyle 0}};
  \node at ( 0.40,-1.34) {\ensuremath{\scriptstyle 0}};

  \node[anchor=west] at (0.48,-0.10)
    {\ensuremath{\scriptstyle
      \frac{r_1(2x_1)-r_1(-2x_2)}{x_1+x_2}}};
\end{tikzpicture}%
}

\newcommand{\crosscupcapdots}{%
\begin{tikzpicture}[baseline=-0.5ex,scale=0.8]
  \draw[thick]
    (-0.38,1.35)
    .. controls (-0.28,1.22) and (-0.12,1.05) .. (0,0.94)
    .. controls (0.12,0.83) and (0.34,0.56) .. (0.30,0.34)
    .. controls (0.26,0.14) and (0.08,0.10) .. (0,0.10);
  \draw[thick]
    (0.38,1.35)
    .. controls (0.28,1.22) and (0.12,1.05) .. (0,0.94)
    .. controls (-0.12,0.83) and (-0.34,0.56) .. (-0.30,0.34)
    .. controls (-0.26,0.14) and (-0.08,0.10) .. (0,0.10);

  \draw[thick]
    (-0.38,-1.15)
    .. controls (-0.38,-0.55) and (0.38,-0.55) ..
    (0.38,-1.15);

  \minusnode{(0.30,0.45)}
  \minusnode{(-0.33,-0.72)}

  \node at (-0.40, 1.54) {\ensuremath{\scriptstyle 0}};
  \node at ( 0.40, 1.54) {\ensuremath{\scriptstyle 0}};
  \node at (-0.40,-1.34) {\ensuremath{\scriptstyle 0}};
  \node at ( 0.40,-1.34) {\ensuremath{\scriptstyle 0}};
\end{tikzpicture}%
}

\newcommand{\cupcapdots}[1]{%
\begin{tikzpicture}[baseline=-0.5ex,scale=0.8]
  \draw[thick]
    (-0.38,1.20)
    .. controls (-0.38,0.46) and (0.38,0.46) ..
    (0.38,1.20);

  \draw[thick]
    (-0.38,-1.15)
    .. controls (-0.38,-0.40) and (0.38,-0.40) ..
    (0.38,-1.15);

  \ifx#1L
    \draw[white,line width=2.2pt] (-0.38,1.20) -- (-0.38,0.78);
    \draw[densely dashed,thick] (-0.38,1.20) -- (-0.38,0.78);
    \minusnode{(-0.38,0.78)}
  \else
    \draw[white,line width=2.2pt] (0.38,1.20) -- (0.38,0.78);
    \draw[densely dashed,thick] (0.38,1.20) -- (0.38,0.78);
    \minusnode{(0.38,0.78)}
  \fi

  \draw[white,line width=2.2pt] (-0.38,-0.76) -- (-0.38,-1.15);
  \draw[densely dashed,thick] (-0.38,-0.76) -- (-0.38,-1.15);
  \minusnode{(-0.38,-0.76)}

  \node at (-0.40, 1.39) {\ensuremath{\scriptstyle 0}};
  \node at ( 0.40, 1.39) {\ensuremath{\scriptstyle 0}};
  \node at (-0.40,-1.34) {\ensuremath{\scriptstyle 0}};
  \node at ( 0.40,-1.34) {\ensuremath{\scriptstyle 0}};
\end{tikzpicture}%
}
\newcommand{\cupcapminus}{%
\begin{tikzpicture}[baseline=-0.5ex,scale=0.8]
  \draw[thick]
    (-0.38,1.20)
    .. controls (-0.38,0.46) and (0.38,0.46) ..
    (0.38,1.20);

  \draw[thick]
    (-0.38,-1.15)
    .. controls (-0.38,-0.40) and (0.38,-0.40) ..
    (0.38,-1.15);

  \draw[white,line width=2.2pt] (-0.38,1.20) -- (-0.38,0.78);
  \draw[densely dashed,thick] (-0.38,1.20) -- (-0.38,0.78);
  \minusnode{(-0.38,0.78)}

  \draw[white,line width=2.2pt] (-0.38,-0.76) -- (-0.38,-1.15);
  \draw[densely dashed,thick] (-0.38,-0.76) -- (-0.38,-1.15);
  \minusnode{(-0.38,-0.76)}

  \node at (-0.40, 1.39) {\ensuremath{\scriptstyle 0}};
  \node at ( 0.40, 1.39) {\ensuremath{\scriptstyle 0}};
  \node at (-0.40,-1.34) {\ensuremath{\scriptstyle 0}};
  \node at ( 0.40,-1.34) {\ensuremath{\scriptstyle 0}};
\end{tikzpicture}%
}

\newcommand{\cupcaparrow}{%
\begin{tikzpicture}[baseline=-0.5ex,scale=0.8,>=stealth]
  \draw[thick]
    (-0.38,1.20)
    .. controls (-0.38,0.48) and (0.38,0.48) ..
    (0.38,1.20);

  \draw[->,densely dashed] (0,0.48) -- (0,-0.25);

  \draw[thick]
    (-0.38,-1.15)
    .. controls (-0.38,-0.40) and (0.38,-0.40) ..
    (0.38,-1.15);

  \node at (-0.40, 1.39) {\ensuremath{\scriptstyle 0}};
  \node at ( 0.40, 1.39) {\ensuremath{\scriptstyle 0}};
  \node at (-0.40,-1.34) {\ensuremath{\scriptstyle 0}};
  \node at ( 0.40,-1.34) {\ensuremath{\scriptstyle 0}};
\end{tikzpicture}%
}
\newcommand{\annotatedcupcap}[2]{%
\begin{tikzpicture}[baseline=-0.5ex,scale=0.8,>=stealth]
  \draw[thick]
    (-0.38,1.25)
    .. controls (-0.38,0.48) and (0.38,0.48) ..
    (0.38,1.25);

  \draw[->,densely dashed] (0,0.48) -- (0,-0.20);

  \draw[thick]
    (-0.38,-1.20)
    .. controls (-0.38,-0.38) and (0.38,-0.38) ..
    (0.38,-1.20);

  \draw[->,densely dashed] (-0.20,-0.82) -- (0.20,-0.82);
  \node[anchor=west] at (0.52,-0.82)
    {\ensuremath{\scriptstyle \lowerfactor}};

  \node[anchor=west] at (0.64,0.88)
    {\ensuremath{\scriptstyle #1}};
  \node[anchor=west] at (0.88,-0.53)
    {\ensuremath{\scriptstyle #2}};

  \node at (-0.40, 1.44) {\ensuremath{\scriptstyle 0}};
  \node at ( 0.40, 1.44) {\ensuremath{\scriptstyle 0}};
  \node at (-0.40,-1.39) {\ensuremath{\scriptstyle 0}};
  \node at ( 0.40,-1.39) {\ensuremath{\scriptstyle 0}};
\end{tikzpicture}%
}

\newcommand{\cupcapwithfactors}{%
\begin{tikzpicture}[baseline=-0.5ex,scale=0.8,>=stealth]
  \draw[thick]
    (-0.38,1.25)
    .. controls (-0.38,0.48) and (0.38,0.48) ..
    (0.38,1.25);

  \draw[->,densely dashed] (-0.20,0.92) -- (0.20,0.92);
  \node[anchor=west] at (0.52,0.92)
    {\ensuremath{\scriptstyle \upperfactor}};

  \draw[->,densely dashed] (0,0.48) -- (0,-0.20);

  \draw[thick]
    (-0.38,-1.20)
    .. controls (-0.38,-0.38) and (0.38,-0.38) ..
    (0.38,-1.20);

  \draw[->,densely dashed] (-0.20,-0.82) -- (0.20,-0.82);
  \node[anchor=west] at (0.52,-0.82)
    {\ensuremath{\scriptstyle \lowerfactor}};

  \node at (-0.40, 1.44) {\ensuremath{\scriptstyle 0}};
  \node at ( 0.40, 1.44) {\ensuremath{\scriptstyle 0}};
  \node at (-0.40,-1.39) {\ensuremath{\scriptstyle 0}};
  \node at ( 0.40,-1.39) {\ensuremath{\scriptstyle 0}};
\end{tikzpicture}%
}
\newcommand{\topfactor}{\frac{1}{1-x_1+x_2}}
\newcommand{\negtopfactor}{-\frac{1}{1-x_1+x_2}}
\newcommand{\lowerfactor}{A_V(y)}
\newcommand{\upperfactor}{B_V(y)}
\newcommand{\bottomfactor}{\frac{1}{1+x_1-x_2}}
\newcommand{\middlefactor}{r(x_1,x_2)}
\newcommand{\middlecorrection}{\frac{1-r(x_1,x_2)}{x_1-x_2}}

\newcommand{\topzeros}{%
  \node at (-0.92,1.69) {$\scriptstyle 0$};
  \node at ( 0.00,1.69) {$\scriptstyle 0$};
  \node at ( 1.47,1.69) {$\scriptstyle 0$};}
\newcommand{\bottomzeros}{%
  \node at (-0.92,-1.69) {$\scriptstyle 0$};
  \node at ( 0.00,-1.69) {$\scriptstyle 0$};
  \node at ( 1.47,-1.69) {$\scriptstyle 0$};}

\newcommand{\crossone}[3]{%
  \draw[thick]
    (-0.90,#1) .. controls (-0.69,{#1-0.19}) and (-0.20,{#2+0.19}) .. (0,#2);
  \draw[thick]
    (0,#1) .. controls (-0.20,{#1-0.19}) and (-0.69,{#2+0.19}) .. (-0.90,#2);
  \draw[thick] (1.45,#1)--(1.45,#2);
  \draw[->,densely dashed] (-0.59,{#1-0.27})--(-0.18,{#1-0.27});
  \node[anchor=west] at (0.10,{#1-0.27}) {$\scriptstyle #3$};}

\newcommand{\crosstwo}[3]{%
  \draw[thick] (-0.90,#1)--(-0.90,#2);
  \draw[thick]
    (0,#1) .. controls (0.31,{#1-0.19}) and (1.14,{#2+0.19}) .. (1.45,#2);
  \draw[thick]
    (1.45,#1) .. controls (1.14,{#1-0.19}) and (0.31,{#2+0.19}) .. (0,#2);
  \draw[->,densely dashed] (0.55,{#2+0.28})--(1.02,{#2+0.28});
  \node[anchor=west] at (1.20,{#2+0.28}) {$\scriptstyle #3$};}

\newcommand{\parallelone}[3]{%
  \draw[thick] (-0.90,#1)--(-0.90,#2);
  \draw[thick] (0,#1)--(0,#2);
  \draw[thick] (1.45,#1)--(1.45,#2);
  \draw[->,densely dashed] (-0.68,{(#1+#2)/2})--(-0.21,{(#1+#2)/2});
  \node[anchor=west] at (0.10,{(#1+#2)/2}) {$\scriptstyle #3$};}

\newcommand{\paralleltwo}[3]{%
  \draw[thick] (-0.90,#1)--(-0.90,#2);
  \draw[thick] (0,#1)--(0,#2);
  \draw[thick] (1.45,#1)--(1.45,#2);
  \draw[->,densely dashed] (0.44,{(#1+#2)/2})--(1.05,{(#1+#2)/2});
  \node[anchor=west] at (1.22,{(#1+#2)/2}) {$\scriptstyle #3$};}

\newcommand{\turnone}[4]{%
  \draw[thick]
    (-0.90,#1) .. controls (-0.90,{#1-0.43}) and (0,{#1-0.43}) .. (0,#1);
  \draw[thick]
    (-0.90,#2) .. controls (-0.90,{#2+0.43}) and (0,{#2+0.43}) .. (0,#2);
  \draw[thick] (1.45,#1)--(1.45,#2);
  \draw[->,densely dashed] (0,{#1-0.12})--(0,{#2+0.12});
  \node[anchor=west] at (0.12,{#1-0.13}) {$\scriptstyle #3$};
  \node[anchor=west] at (0.12,{(#1+#2)/2-0.10}) {$\scriptstyle #4$};}

\newcommand{\turntwo}[4]{%
  \draw[thick] (-0.90,#1)--(-0.90,#2);
  \draw[thick]
    (0,#1) .. controls (0,{#1-0.43}) and (1.45,{#1-0.43}) .. (1.45,#1);
  \draw[thick]
    (0,#2) .. controls (0,{#2+0.43}) and (1.45,{#2+0.43}) .. (1.45,#2);
  \draw[->,densely dashed] (1.45,{#1-0.12})--(1.45,{#2+0.12});
  \node[anchor=west] at (1.57,{#1-0.13}) {$\scriptstyle #3$};
  \node[anchor=west] at (1.57,{(#1+#2)/2-0.10}) {$\scriptstyle #4$};}

\newcommand{\upperstage}[1]{%
  \ifcase#1\relax
  \or\crossone{1.50}{0.50}{\topfactor}%
  \or\parallelone{1.50}{0.50}{\negtopfactor}%
  \or\turnone{1.50}{0.50}{B_N(x)}{P_N(x_1,x_2)}%
  \fi}

\newcommand{\middlestage}[1]{%
  \ifcase#1\relax
  \or\crosstwo{0.50}{-0.50}{\middlefactor}%
  \or\paralleltwo{0.50}{-0.50}{\middlecorrection}%
  \or\turntwo{0.50}{-0.50}{q_V(x_1,x_2)}{A_V(x)}%
  \fi}

\newcommand{\lowerstage}[1]{%
  \ifcase#1\relax
  \or\crossone{-0.50}{-1.50}{\bottomfactor}%
  \or\parallelone{-0.50}{-1.50}{\bottomfactor}%
  \or\turnone{-0.50}{-1.50}{q_N(x_1,x_2)}{A_N(x)}%
  \fi}

\newcommand{\braidterm}[3]{%
\begin{tikzpicture}[
  baseline=(current bounding box.center),scale=0.92,>=stealth,
  line cap=round,line join=round]
  \upperstage{#1}
  \middlestage{#2}
  \lowerstage{#3}
  \topzeros
  \bottomzeros
\end{tikzpicture}%
}

\newcommand{\postcomposition}{%
\begin{tikzpicture}[baseline=-0.5ex,scale=0.76,>=stealth,line cap=round]
  \draw[thick] (-0.75,0.75)--(-0.75,-0.75);
  \draw[thick] (0,0.75)--(0,-0.75);
  \draw[thick] (0.75,0.75)--(0.75,-0.75);
  \draw[->,densely dashed] (-0.60,0.32)--(-0.13,0.32);
  \node[anchor=west] at (0.02,0.32)
    {$\scriptscriptstyle 1-x_1+x_2$};
  \draw[->,densely dashed] (0.12,-0.30)--(0.60,-0.30);
  \node[anchor=west] at (0.77,-0.30)
    {$\scriptscriptstyle 1+x_1-x_2$};
\end{tikzpicture}%
}
\newcommand{\afactor}{\frac{1}{1-x_1+x_2}}
\newcommand{\nafactor}{-\frac{1}{1-x_1+x_2}}
\newcommand{\bfactor}{\frac{1}{1+x_1-x_2}}
\newcommand{\rfactor}{r(x_1,x_2)}
\newcommand{\cfactor}{\frac{1-r(x_1,x_2)}{x_1-x_2}}

\newcommand{\threezerostop}{%
  \node at (-0.92,1.57) {$\scriptstyle 0$};
  \node at (0,1.57) {$\scriptstyle 0$};
  \node at (1.47,1.57) {$\scriptstyle 0$};}
\newcommand{\threezerosbottom}{%
  \node at (-0.92,-1.57) {$\scriptstyle 0$};
  \node at (0,-1.57) {$\scriptstyle 0$};
  \node at (1.47,-1.57) {$\scriptstyle 0$};}

\newcommand{\residualterm}[1]{%
\biggl[
\vcenter{\hbox{%
\begin{tikzpicture}[scale=0.82,>=stealth,line cap=round,line join=round]
  \draw[thick] (-0.80,1.25)--(-0.80,0.30)
    .. controls (-0.80,-0.15) and (0,-0.15) .. (0,0.30)--(0,1.25);
  \draw[thick] (-0.80,-1.20)
    .. controls (-0.80,-0.55) and (0,-0.55) .. (0,-1.20);
  \draw[thick] (1.10,1.25)--(1.10,-1.20);
  \node[circle,draw,fill=white,inner sep=0.2pt,minimum size=5pt]
    at (0,0.02) {$\scriptscriptstyle -$};
  \node[circle,draw,fill=white,inner sep=0.2pt,minimum size=5pt]
    at (-0.80,-0.70) {$\scriptscriptstyle -$};
  \node[anchor=west] at (0.12,-0.05) {$\scriptstyle #1$};
  \draw[->,densely dashed] (-0.52,-0.88)--(0.74,-0.88);
  \node[anchor=west] at (0.84,-0.88) {$\scriptstyle 1+x_1-x_2$};
  \node at (-0.82,1.43) {$\scriptstyle0$};
  \node at (0,1.43) {$\scriptstyle0$};
  \node at (1.12,1.43) {$\scriptstyle0$};
  \node at (-0.82,-1.38) {$\scriptstyle0$};
  \node at (0,-1.38) {$\scriptstyle0$};
  \node at (1.12,-1.38) {$\scriptstyle0$};
\end{tikzpicture}}}
\biggr]u^{-1}%
}

\newcommand{\cupdiagram}[2]{%
\begin{tikzpicture}[baseline=-0.5ex,scale=0.82,>=stealth,line cap=round]
  \draw[thick] (-0.50,0.65)
    .. controls (-0.50,-0.05) and (0.50,-0.05) .. (0.50,0.65);
  \draw[->,densely dashed] (-0.28,0.35)--(0.28,0.35);
  \node[anchor=west] at (0.62,0.35) {$\scriptstyle #1$};
  \node[anchor=west] at (0.62,-0.02) {$\scriptstyle #2$};
  \node at (-0.52,0.84) {$\scriptstyle0$};
  \node at (0.52,0.84) {$\scriptstyle0$};
\end{tikzpicture}%
}

\hypersetup{hidelinks}

\title[Integral Brauer categorification]{Affine and cyclotomic Brauer categorification via $\imath$-Kac--Moody $2$-categories: the integral type
 $\operatorname{AIII}$ case}

 \author{Mengmeng  Gao}
\email{g19920119@163.com}
\author{Hebing  Rui}
\email{hbrui@tongji.edu.cn}
\author{Linliang Song}
\email{llsong@tongji.edu.cn}
\author{Jun Wu}
\email{wujunstu@tongji.edu.cn}
\address{School of Mathematical Sciences, Tongji University, Shanghai, 200092, China.}
\thanks{M. Gao is  partially supported by NSFC (grant No. 12301038), H. Rui is partially supported   by NSFC (grant No.  11971351), and  L. Song is  partially supported   by NSFC (grant No.  12071346). }

\date{}

\begin{document}
\begin{abstract}Continuing our work on the half-integral case, we construct
generalized nilpotent $2$-representations of the integral
quasi-split type $\operatorname{AIII}$ $\imath$-Kac--Moody
$2$-category from affine Brauer actions with dot spectrum
exactly $\mathbb Z$ on locally Schurian $\mathbb C$-linear
categories. Conversely, nilpotent $2$-representations on such
categories yield affine Brauer actions with spectrum contained
in $\mathbb Z$.
At the fixed vertex $0$, we directly define corrected generators
and verify the nil-Brauer and mixed relations.
For $\mathbf u$-admissible integral parameters at arbitrary
positive level, we prove an isomorphism between the locally
unital algebras of the cyclotomic Brauer category and a
cyclotomic quotient of the principal $2$-representation.
This gives compatible $\mathbb Z$-gradings on the category and
its endomorphism algebras, providing a Brauer analogue of the
graded realization of type~$A$ cyclotomic Hecke algebras.
Over algebraically closed fields of odd characteristic $p$,
 we also obtain both constructions with spectrum
contained in $\mathbb F_p$, for the mixed affine datum
restricted to its spectral weight subgroup.
\end{abstract}
\maketitle
\enlargethispage{8pt}

\par\medskip
\setlength{\baselineskip}{14.2pt}

\tableofcontents

\section{Introduction}

\subsection{Main results}

Ariki's theorem identifies decomposition numbers of cyclotomic Hecke
algebras in characteristic zero with canonical-basis coefficients
\cite{Ariki-dec}. Quiver Hecke algebras give graded algebraic models
for quantum groups \cite{KL-rings,R-2cat}, and the
Brundan--Kleshchev isomorphism identifies blocks of cyclotomic Hecke
algebras, including their degenerate versions, with cyclotomic quiver
Hecke algebras \cite{BK-klr}. Brundan--Kleshchev also described their
graded decomposition numbers in characteristic zero in terms of
canonical bases \cite{BK-gr}. These results motivate the search for
analogous constructions for cyclotomic Brauer algebras.

Cyclotomic Brauer algebras, also called cyclotomic Nazarov--Wenzl
algebras \cite{AMR-cyc}, arise in higher Schur--Weyl duality for
classical Lie algebras. Rui--Song introduced the affine Brauer
category, proved basis theorems for it and its admissible cyclotomic
quotients, and related them to parabolic categories $\mathcal O$
in types $B$, $C$, and $D$ \cite{RS-cyc}. In addition to dots and
crossings, the category has cups and caps, which must be included
in the categorical action.

The corresponding quantum-group structures are quantum symmetric
pairs. Bao--Wang developed their canonical-basis theory
\cite{BW-ihes}, and Bao--Shan--Wang--Webster constructed a
categorification in quasi-split type $\operatorname{AIII}$
\cite{BSWW-icate}. Brundan--Wang--Webster introduced the nil-Brauer
category \cite{BWW24} and the $\imath$-Kac--Moody $2$-categories
for quasi-split symmetric types \cite{BWW25}. We use these
categories to give an integral counterpart of the half-integral
construction in \cite{GRS}.

Let $\mathcal{AB}$ be the affine Brauer category and let
$\mathfrak U^\imath$ be the integral quasi-split type
$\operatorname{AIII}$ $\imath$-Kac--Moody $2$-category, with
\[
 I=\mathbb Z,\qquad \tau(i)=-i.
\]
For an affine Brauer action on a locally Schurian category with
dot spectrum $\mathbb Z$, the generating functor decomposes as
\[
 E=\bigoplus_{i\in\mathbb Z}E_i.
\]
We use these generalized-eigenspace functors to construct an action
of $\mathfrak U^\imath$, and conversely recover an affine Brauer
action under the nilpotence and finiteness conditions below.
The affine constructions are ungraded. The cyclotomic isomorphism
in Theorem~\ref{thm:intro-B} transports the grading of the
$2$-category to the cyclotomic Brauer algebras.

The notions of nilpotent and generalized nilpotent
$2$-representations are recalled in Definition~2.5.
Our first main result is as follows.

\begin{THEOREM}\label{thm:intro-A}
Let $\mathcal M$ be a locally Schurian $\mathbb C$-linear category.
\begin{itemize}
\item[(1)]
Suppose that $\mathcal M$ admits an affine Brauer categorification
and that the spectrum of the dot endomorphism on the generating
endofunctor is exactly $\mathbb Z$. Then $\mathcal M$ carries a
generalized nilpotent $2$-representation of the integral quasi-split
type $\operatorname{AIII}$ $2$-category $\mathfrak U^\imath$.
Moreover, the full subcategory of finitely generated objects is a
nilpotent $2$-subrepresentation. In particular, if $\mathcal M$ is
locally finite Abelian, then the resulting $2$-representation is
nilpotent.
\item[(2)]
Suppose that $\mathcal M$ carries a $2$-representation of
$\mathfrak U^\imath$. For each weight $\lambda$, let
$\mathcal N_\lambda\subseteq\mathcal M_\lambda$ be a full additive
subcategory of finitely generated ambient objects, closed under
isomorphisms, and suppose that all $E_i$ preserve
$\mathcal N=\bigoplus_\lambda\mathcal N_\lambda$.
If every object of $\mathcal N$ has finite color support and
nilpotent dot endomorphisms, then $\mathcal N$ admits an affine
Brauer action whose dot spectrum is contained in $\mathbb Z$.
\end{itemize}
\end{THEOREM}

Part~\textup{(1)} is Theorem~5.16, and part~\textup{(2)} is
Corollary~8.32. Theorem~8.31 also gives the reverse construction
on an entire locally Schurian nilpotent $2$-representation.
The forward statement assumes the full spectrum $\mathbb Z$,
whereas the reconstructed dot in part~\textup{(2)} need not realize
every integer. We do not assert that the two constructions are
mutually inverse for arbitrary $2$-representations.

For the cyclotomic application, fix a positive integer $a$ and
$\mathbf u=(u_1,\ldots,u_a)\in\mathbb Z^a$ and the corresponding
$\mathbf u$-admissible cyclotomic Brauer category
$\mathcal{CB}_{\mathbf u}$, with locally unital path algebra
$B_{\mathbf u}$. The root multiplicities of
\[
 m(u)=\prod_{r=1}^a(u-u_r)
\]
determine an initial $\imath$-weight $\kappa$ and a cyclotomic
quotient $\mathcal Q(m)$ of the principal $2$-representation.
Let $Q_m$ be the locally unital algebra of the full subcategory
of $\mathcal Q(m)$ on the objects $E_{\mathbf i}1_\kappa$,
where $\mathbf i$ runs over finite sequences in $I$.

\begin{THEOREM}\label{thm:intro-B}
For $\mathbf u$-admissible integral cyclotomic parameters, there is
an isomorphism of locally unital $\mathbb C$-algebras
\[
 Q_m\cong B_{\mathbf u}.
\]
Consequently, the additive idempotent completions of this full
subcategory and $\mathcal{CB}_{\mathbf u}$ are equivalent, and
every $\mathbf u$-admissible integral cyclotomic Brauer algebra
inherits a natural $\mathbb Z$-grading.
\end{THEOREM}

This follows from Theorem~9.8 and Corollary~9.9.
The isomorphism identifies the generalized-eigenspace idempotents
and the dot, crossing, cup, and cap generators.
Transporting the grading of $Q_m$ along this ungraded isomorphism
gives a $\mathbb Z$-grading on every morphism space of
$\mathcal{CB}_{\mathbf u}$, compatible with composition, and
hence on its fixed-rank endomorphism algebras.
This provides a Brauer analogue of the graded realization of
type~$A$ cyclotomic Hecke algebras.

The proof follows the method of our half-integral paper
\cite[Section~8, Theorem~8.6]{GRS}: the two constructions descend
with matching cyclotomic parameters, and recovery of the generators
and spectral idempotents gives mutually inverse homomorphisms.
For the analogous cyclotomic isomorphism in the Heisenberg setting,
see \cite[Theorem~5.19]{BSW-heisenberg}.

Ehrig--Stroppel realized certain cyclotomic Nazarov--Wenzl algebras
in parabolic category $\mathcal O$ of type $D$ \cite{ES-NW}
and constructed gradings on ordinary Brauer algebras with integral
parameter using level-two cyclotomic algebras \cite{ES-grading}.
Theorem~\ref{thm:intro-B} treats admissible integral parameters at
arbitrary level, with the grading induced by the
$\imath$-Kac--Moody $2$-category.

Theorem~9.10 also gives actions on the subcategories of parabolic
categories $\mathcal O$ whose simple composition factors have
highest weights with integral $\varepsilon$-coordinates.
The integral $2$-category acts in type $B$, and the even component
of the half-integral $2$-category acts in types $C$ and $D$.

Over an algebraically closed field $k$ of odd characteristic $p$,
the images of both $\mathbb Z$ and $\tfrac12+\mathbb Z$ are
$\mathbb F_p$. The resulting cyclic datum contains both the
fixed vertex $0$ and the exceptional adjacent pair
$\{-\tfrac12,\tfrac12\}$.
Write $\mathfrak U^\imath_{p,\mathrm{sp}}$ for the full
$2$-subcategory on the spectral weight subgroup defined in
Section~10. Corollaries~10.3 and~10.5 construct a nilpotent
$2$-representation of this $2$-category from an affine Brauer
action with spectrum contained in $\mathbb F_p$ on a locally
Schurian $k$-linear category, and conversely.
For $p>3$, the proof combines the zero-color calculations below
with the half-integral relations of \cite{GRS}.
In characteristic $3$, both vertices of the exceptional pair
are also adjacent to $0$. Subsection~10.5 adjusts the zero-color
normalization and verifies the corresponding mixed braid
relations, completing the construction in every odd
characteristic.

\subsection{Comparison with the half-integral case}

In \cite{GRS}, we constructed affine Brauer actions and
half-integral $\imath$-Kac--Moody $2$-representations in both
directions, using the even component $\mathfrak U^\imath_+$,
and identified the corresponding admissible cyclotomic algebras
and their gradings. Theorems~\ref{thm:intro-A} and~\ref{thm:intro-B}
give the integral counterparts. In the half-integral case, the
involution $i\mapsto-i$ on $\tfrac12+\mathbb Z$ has no fixed point.
The exceptional relation occurs at the adjacent pair
$\{-\tfrac12,\tfrac12\}$. In the integral case, the vertex $0$
is fixed, and the corresponding generators satisfy nil-Brauer
relations \cite{BWW24}.

Both constructions use generalized eigenspaces and rational
functions of the dot. Relations between opposite colors also
involve cups and caps. We refer to \cite{GRS} for the nonexceptional
calculations that apply unchanged. The formulas at $0$, and their
compatibility with the other colors, require separate proofs.

The weakly triangular representation theory and its
decategorification were developed in \cite{RS-brauer,GRS-tri}.
The cyclotomic isomorphism and grading provide a setting for a
graded refinement of this representation theory. Constructing
a categorical bar involution, identifying indecomposable
projectives with $\imath$-canonical basis elements, and
interpreting decomposition numbers as a Brauer analogue of
Ariki's theorem are goals of subsequent work.

\subsection{The zero vertex and the two constructions}

The functor $E_0$ is self-adjoint. Its dot, crossing, cup, and cap
must satisfy the nil-Brauer relations, which cannot be obtained
by substituting $i=0$ into the nonzero-color formulas.

For the forward construction, the local factors of the dot's
minimal polynomial and the bubble identities in
\cite{GRS-block,SW,GRS} determine the normalizing functions.
Their expansions are evaluated at nilpotent shifted dots.
The factor $\Gamma$ in Definition~4.4 makes the zero-color
formulas compatible with the other colors.
We verify the nil-Brauer and mixed relations for these generators.
For the zero-color braid relation, we reduce both sides to
Brauer diagrams and compare all coefficients, including the
bubble contributions.

For the reverse construction, we take a $2$-subrepresentation
satisfying Theorem~\ref{thm:intro-A}(2). We sum its color functors
and define the affine Brauer generators using the
$\imath$-Kac--Moody generators. The finiteness and nilpotence
conditions make the sums and power-series evaluations finite
on each finitely generated object. The resulting dot has
integral spectrum.

\subsection{Related recent work}
Brundan--Savage--Webster \cite{BSW-Brauer} study Brauer and
$\imath$-quantum categorifications in the locally finite
Abelian and Karoubian settings. They establish mutually
inverse constructions
\cite[Corollary~8.3 and Theorem~9.9]{BSW-Brauer}
and equivalences between the corresponding deformed
cyclotomic quotients \cite[Theorem~10.8]{BSW-Brauer}.
 In the cyclotomic setting,
\cite[Lemma~10.6]{BSW-Brauer} extends the local
finite-dimensionality result of
\cite[Proposition~8.5]{GRS}.

The present paper continues \cite{GRS} by developing
the integral case through objectwise constructions
on locally Schurian categories. It treats the fixed
vertex and the mixed affine case, and establishes
the corresponding cyclotomic algebra isomorphisms
and gradings.

\subsection{History of the project}

This paper forms part of the project initiated in 2022,
as described in \cite[p.~2]{GRS}.
The integral generators at the fixed vertex~$0$ were
constructed and their nil-Brauer relations verified
in June 2026. The $\Gamma$-correction terms for
compatibility with the other colors were introduced
in August 2026.

Our half-integral paper was first submitted to arXiv
on 7 September 2026. Its September 2026 version
describes the integral case as a separate work
in preparation \cite[pp.~4 and~76]{GRS}.

\subsection{Organization of the paper}

Sections~2--3 give the categorical setup and the
generalized-eigenspace and weight decompositions.
Sections~4--7 construct and verify the forward action.
Section~8 gives the reverse construction, Section~9 proves
the cyclotomic isomorphisms and grading results, and Section~10
treats the mixed affine case in every odd characteristic,
with characteristic $3$ handled separately in Subsection~10.5.
The appendices contain the diagrammatic reductions and
coefficient calculations.

\medskip
\noindent\textbf{Declaration on the use of artificial intelligence.}

AI assistance in the original mathematical construction was
limited to ChatGPT's suggestion of the $\Gamma$-correction terms
in Definition~4.4. All defining relations, including those
involving the corrected generators, were verified by the authors.
AI tools also assisted with language editing, formatting, and
supplementary symbolic consistency checks.
The authors take responsibility for all mathematical statements
and proofs.

\section{The integral type \texorpdfstring{$\operatorname{AIII}$}{AIII}
\texorpdfstring{$\imath$}{i}--Kac--Moody 2-category}

Fix an algebraically closed field $\Bbbk$ of characteristic different
from $2$. We specify the characteristic-zero integral datum in
Subsection~\ref{subsec:integral-datum}; the mixed affine datum over a
field of characteristic $p\ge3$ is treated separately in
Section~\ref{sec:mixed-affine}.

\subsection{Root data}
We first recall the root datum used by Brundan--Wang--Webster
\cite[Sections~2.1 and~3.1]{BWW25}, restricting to symmetric Cartan
matrices and $d_i=1$ for all $i\in I$.
Let $A=(a_{ij})_{i,j\in I}$ be a symmetric generalized Cartan
matrix, and let $\tau:I\to I$ be an involution satisfying
$a_{\tau(i),\tau(j)}=a_{ij}$ for all $i,j\in I$.
Thus $(I,\tau)$ is a quasi-split Satake datum.
We choose a decomposition
\begin{equation}\label{dec-I}
I=I_1\sqcup I_0\sqcup I_{-1},
\end{equation}
where $I_0$ is the set of $\tau$-fixed vertices, $I_1$ is a set of
representatives for the $\tau$-orbits of size two, and
$I_{-1}=\tau(I_1)$. Choose integers $\boldsymbol\varsigma=(\varsigma_i)_{i\in I}$ such
that
\begin{equation}\label{bold-des}
\begin{aligned}
\varsigma_i+\varsigma_{\tau(i)}&=-a_{i,\tau(i)} &&(i\in I),\\
\varsigma_i&\geq0 &&(i\ne\tau(i)).
\end{aligned}
\end{equation}
In particular, $\varsigma_i=-1$ at a fixed vertex, and
$\varsigma_i=\varsigma_{\tau(i)}=0$ when $a_{i,\tau(i)}=0$.

Let
$(X,Y,\Pi,\Pi^\vee)
$
be a realization of $A$, where $X$ and $ Y$ are the weight  and coweight lattices, respectively, and
$\Pi=\{\alpha_i\mid i\in I\}\subseteq X$ and
$\Pi^\vee=\{\alpha_i^\vee\mid i\in I\}\subseteq Y$ are  the sets of simple roots and  simple  coroots.
Thus, $X$ and $Y$ are free abelian groups equipped with a non-degenerate  pairing
\begin{equation}\label{pair-ikm1}
\langle\ ,\ \rangle:Y\times X\longrightarrow\mathbb Z
\end{equation}
such that
$\langle \alpha_i^\vee,\alpha_j\rangle=a_{ij}$ for all $i,j\in I$.

 We assume that there are involutions, again
denoted by $\tau$ and $\tau^*$, on $X$ and $Y$, respectively, such
that
\[
\tau(\alpha_i)=\alpha_{\tau(i)},\qquad
\tau^*(\alpha_i^\vee)=\alpha^\vee_{\tau(i)},\qquad
\langle\tau^*(h),\mu\rangle=\langle h,\tau(\mu)\rangle
\quad(h\in Y,\ \mu\in X).
\]
The $\imath$-weight and $\imath$-coweight lattices are
\begin{equation}\label{i-coweight}
X^\imath:=X/\operatorname{im}(\operatorname{id}+\tau)
\ \ \text{and}\ \
Y^\imath:=\ker(\operatorname{id}+\tau^*).
\end{equation}
Write $\bar\alpha_j$ for the image of $\alpha_j$ in $X^\imath$, and use $\lambda+\alpha_j$ as shorthand for $\lambda+\bar \alpha_j$.
For $\lambda\in X^\imath$, choose a lift $\widehat\lambda\in X$ and
set
$$
\lambda_i:=(\alpha_i^\vee-\alpha^\vee_{\tau(i)})(\widehat\lambda).$$
This is independent of the choice of $\widehat\lambda$, and
$$
(\lambda+\alpha_j)_i
=
\lambda_i+a_{ij}-a_{i,\tau(j)}.
$$
If $i=\tau(i)$, then $\lambda_i=0$, but the parity
$\epsilon_i(\lambda)\in\{0,1\}$, defined by
\[
\epsilon_i(\lambda)\equiv
\langle\alpha_i^\vee,\widehat\lambda\rangle\pmod2,
\]
is well defined. Indeed, changing the lift by $(1+\tau)\mu$ changes
this pairing by $2\langle\alpha_i^\vee,\mu\rangle$.

Besides the datum above, the general definition involves KLR
parameters $Q_{ij}(u,v)$, normalization homomorphisms
$c_i\colon X\to\Bbbk^\times$, and parameters
$\zeta=(\zeta_i)_{i\in I}$ satisfying the compatibility conditions
in \cite[(3.12)--(3.18)]{BWW25}.
We give the integral choices below and use them in
Definition~\ref{def-2ikac}. The mixed affine choices are given in
Section~\ref{subsec:mixed-affine-datum}.

\subsection{Quasi-split type $\operatorname{AIII}$ data and standing assumptions}
\label{subsec:integral-datum}
We now specialize the general datum to the quasi-split type $\operatorname{AIII}$ cases relevant to our applications.
The underlying Dynkin diagram is either the bi-infinite diagram
of type $A_\infty$ or a cyclic affine diagram of type
$A_{e-1}^{(1)}$, where $e>2$. The Cartan matrix is given by
$$
a_{ij}=
\begin{cases}
2,&i=j,\\
-1,&i\text{ and }j\text{ are joined by an edge},\\
0,&\text{otherwise}.
\end{cases}
$$
The involution $\tau$ is a reflection of the corresponding Dynkin
diagram.
In the two infinite cases, two vertices $i$ and $j$ are adjacent if
and only if
$
|i-j|=1$.
In the affine cases, the index set is taken modulo $e$, and two
vertices $i$ and $j$ are adjacent if and only if
$i-j\equiv\pm 1\pmod e$.

For a $\tau$-orbit $\mathcal O$ in $I$, we call the corresponding
vertex of the folded diagram a \emph{$0$-type vertex} if
$\mathcal O=\{i\}$, a \emph{$\frac12$-type vertex} if
$\mathcal O=\{i,\tau(i)\}$ has two adjacent elements, and an
\emph{$A$-type vertex} if it has two nonadjacent elements.

\begin{Defn}\label{qsta3} Up to isomorphism, the infinite and affine quasi-split type $\operatorname{AIII}$  data
which occur below fall into the following three cases.
\begin{itemize}
\item[(1)] The half-integral infinite case:
$I=\frac12+\mathbb Z$ and
$\tau(i)=-i$.
The orbit
$\{-\frac12,\frac12\}
$ is the unique $\frac12$-type vertex, and every other vertex of the
folded diagram is of $A$-type.

\item[(2)]
The integral infinite case:
$I=\mathbb Z$ and
$\tau(i)=-i$,
with two vertices adjacent when their difference is $\pm1$. The vertex
$0$ is the unique $0$-type vertex, and every other vertex of the
folded diagram is of $A$-type.

\item[(3)] The mixed affine case. Let $e\geq 3$ be odd and take
$
I=\mathbb Z/e\mathbb Z$ and
$
\tau(i)=-i$.
The vertex $0$ is the unique $0$-type vertex. The orbit
$\{\frac{e-1}{2},\frac{e+1}{2}\}
$ is the unique $\frac12$-type vertex. All remaining vertices of the
folded diagram are of $A$-type.

\end{itemize}
\end{Defn}

For the mixed affine construction in Section~\ref{sec:mixed-affine},
we take $e=p=\operatorname{char}\Bbbk$, where $p$ is odd.
The case $p>3$ is treated first, while $p=3$ is considered
separately in Subsection~\ref{subsec:mixed-characteristic-three}.

\begin{Defn} \label{i-para} For each of these three cases, we use the following standard choice of
$\boldsymbol\varsigma$:
\begin{itemize}\item [(1)] At a $0$-type vertex $i=\tau(i)$, let $\varsigma_i=-1$.
\item [(2)] At an $A$-type vertex, set
$
\varsigma_i=\varsigma_{\tau(i)}=0
$.\item [(3)] At a $\frac12$-type vertex $\{i,\tau(i)\}$, choose one of its two
vertices as the representative in $I_1$ and set
$\varsigma_i=1$ and $\varsigma_{\tau(i)}=0$.
 \end{itemize}
\end{Defn}
We fix these choices whenever a datum is specified. Changing the
parameters gives the isomorphisms of \cite[Lemma~3.6]{BWW25} when
its weight-shift and scalar conditions are satisfied; such a change
is not automatic for an arbitrary realization or a restricted weight
subgroup.
The half-integral infinite case $I=\frac12+\mathbb Z$
was treated in \cite{GRS}. We now specialize to the integral infinite
case. Except where other assumptions are stated explicitly, we use
this datum through Section~9:
\begin{equation} \label{integral-index} I=\mathbb Z,
\tau(i)=-i \text{ for all $i\in \mathbb Z$, and $I_1=\mathbb Z_+$,}\end{equation}
where $\mathbb Z_+$ is the set of positive integers.

Unless otherwise stated, we assume that $A=(a_{i,j})_{i,j\in I}$ is the Cartan matrix of type $A_\infty$ with vertex set $I$ in \eqref{integral-index}, and that $\Bbbk$ is an algebraically closed field of characteristic zero. We identify each $i\in I=\mathbb Z$ with $i1_\Bbbk$.

From now on, we use the standard realization of integral type
$A_\infty$: $X=\bigoplus_{i\in\mathbb Z}\mathbb Z\Lambda_i$,
$Y=\bigoplus_{i\in\mathbb Z}\mathbb Z\alpha_i^\vee$, with
$\langle\alpha_i^\vee,\Lambda_j\rangle=\delta_{ij}$,
$\alpha_i=2\Lambda_i-\Lambda_{i-1}-\Lambda_{i+1}$, and
$\tau(\Lambda_i)=\Lambda_{-i}$. In this realization,
\[
X^\imath\cong(\mathbb Z/2\mathbb Z)\bar\Lambda_0
\oplus\bigoplus_{i>0}\mathbb Z\bar\Lambda_i,
\qquad \bar\Lambda_{-i}=-\bar\Lambda_i.
\]
In particular,
$\alpha_0=(1+\tau)(\Lambda_0-\Lambda_1)$, so
$\bar\alpha_0=0$.
Thus the folded diagram has one $0$-type vertex and otherwise only
$A$-type vertices. The new arguments are concentrated at the
$0$-type vertex, or equivalently at the eigenvalue $0$.

Let $Q$ be the quiver with vertex set $I=\mathbb Z$ and arrows
$i\to i+1$ for all $i\in\mathbb Z$. Write $\#(i\to j)$ for the
number of arrows from $i$ to $j$. Then
$\#(i\to j)+\#(j\to i)=-a_{ij}$ and
$\#(i\to j)=\#(\tau(j)\to\tau(i))$ for all $i\ne j$.
The parameters of Definition~\ref{i-para} can be
written in terms of this quiver as
\[
\del_i=\begin{cases}
-1,&i=\tau(i),\\
\#(i\to\tau(i)),&i\ne\tau(i).
\end{cases}
\]
For the integral datum, $i$ and $-i$ are nonadjacent when $i\ne0$.
Thus $\del_0=-1$ and $\del_i=0$ for $i\ne0$, in agreement with
the choice of $\boldsymbol\varsigma$ fixed above.

For each $i$, let
$c_i:X \rightarrow \Bbbk^\times$ be a group homomorphism
such that
\begin{equation}\label{ci1}
\begin{cases} c_i(\alpha_j)  = (-1)^{\#(j \rightarrow i)}, & \text{for all $i,j \in I$;}\\
c_{\tau i}(\tau(\lambda))  = (-1)^{\alpha_i^\vee(\lambda)} c_i(\lambda), &
\text{for all $\lambda \in X$ and $i \in I$ with $i \neq \tau i$.}\end{cases}\end{equation}

The characters in \eqref{ci1} exist. Indeed, write
$c_i(\Lambda_j)=(-1)^{f_i(j)}$. The first condition in \eqref{ci1}
becomes
$f_i(j-1)+f_i(j+1)\equiv\#(j\to i)\pmod 2$,
which can be solved in both directions on $\mathbb Z$.
Choose such solutions for $i\geq0$, and, for $i>0$, define
$c_{-i}$ by the second condition in \eqref{ci1}. The symmetry of
the quiver then gives the first condition for $c_{-i}$.

For $i\in I$ and $\lambda\in X^\imath$, choose a lift
$\hat\lambda\in X$ and set
\begin{equation}\label{gammai}
\gamma_i(\lambda)=\begin{cases}
c_i\big(\hat\lambda-\tau(\hat\lambda)\big),&i\ne\tau(i),\\
(-1)^{\epsilon_i(\lambda)},&i=\tau(i),
\end{cases}
\qquad
\zeta_i=\begin{cases}
\pm2^{\del_i},&i\ne\tau(i),\\
-\tfrac12,&i=\tau(i).
\end{cases}
\end{equation}
The value of $\gamma_i(\lambda)$ is independent of the lift:
$(1-\tau)(1+\tau)=0$ for $i\ne\tau(i)$, and the fixed-vertex case
uses the parity defined above. Choose opposite signs for
$\zeta_i$ and $\zeta_{\tau(i)}$ when $i\ne\tau(i)$.
The KLR polynomials and the associated rational functions are
\begin{equation} \label{Qij} Q_{i,j}(x,y):=\begin{cases}
 (x-y)^{\#(i \rightarrow j)}(y-x)^{\#(j \rightarrow i)}
&\text{if $i \neq j$}\\
0&\text{if $i=j$,}
\end{cases}\ \ \
R_{i,j}(x,y) :=
\begin{cases}
Q_{i,j}(x,y)&\text{if $i \neq j$}\\
1/(x-y)^2&\text{if $i=j$,}
\end{cases}
\end{equation}
Using the corrected sign in \cite[(3.24), footnote~2]{BWW25v2}, set
\begin{equation}\label{iQ}
Q^\imath_{i,j}(x,y):=
(-1)^{\delta_{i,\tau j}} Q_{i,j}(x,y),
\end{equation}

We use formal Laurent series in $u^{-1}$, expanded at $u=\infty$.
For $F(u)=\sum_{r\le N}F_r u^r$, write
\[
[F(u)]_{u:r}=F_r,\qquad
[F(u)]_{u:\ge s}=\sum_{r=s}^N F_r u^r,\qquad
[F(u)]_{u:<s}=\sum_{r<s}F_r u^r.
\]
If $x$ denotes the dot on a strand, the two resolvent labels mean
$(u-x)^{-1}=\sum_{r\ge0}x^r u^{-r-1}$ and
$(u+x)^{-1}=\sum_{r\ge0}(-x)^r u^{-r-1}$:
\begin{equation}\label{idgf}
\begin{tikzpicture}[baseline=-0.5ex, scale=0.7,draw=darkred]
		\draw[thick] (0,-0.6) -- (0,-0.2);
		\draw[thick] (0,0.2) -- (0,0.5);
		\draw (0,0) circle (0.13);
		\node at (0,0) {$-$};
	\end{tikzpicture}
:=
\begin{tikzpicture}[baseline=-0.5ex, scale=0.7,draw=darkred]
		\draw[thick] (0,-0.6) -- (0,0.6);
		\pin{0,0}{1.3,0}{(u-x)^{-1}};
	\end{tikzpicture}\
, \quad
\begin{tikzpicture}[baseline=-0.5ex, scale=0.7,draw=darkred]
		\draw[thick] (0,-0.6) -- (0,-0.2);
		\draw[thick] (0,0.2) -- (0,0.5);
		\draw (0,0) circle (0.13);
		\node at (0,0) {$+$};
	\end{tikzpicture}
:=
\begin{tikzpicture}[baseline=-0.5ex, scale=0.7,draw=darkred]
		\draw[thick] (0,-0.6) -- (0,0.6);
		\pin{0,0}{1.3,0}{(u+x)^{-1}};
	\end{tikzpicture}\ ,
\end{equation}
In \eqref{idgf}, the dot marks the variable at which the
displayed rational function is evaluated; it does not indicate an
additional multiplication by the dot.

These choices satisfy the parameter conditions of \cite{BWW25}. For $i\ne j$, put
$t_{i,j}:=Q_{i,j}(1,0)$. Then
$t_{i,j}=c_i(\alpha_j)$; also
$t_{i,i}=c_i(\alpha_i)=1$. Hence the BWW parameters
$r_{i,j}:=c_i(\alpha_j)/t_{i,j}$ equal $1$ for all $i,j\in I$.
The second condition in \eqref{ci1} is \cite[(3.12)]{BWW25}.
For $i\ne0$, the vertices $i$ and $-i$ are nonadjacent, so
$Q_{-i,i}(x,-x)=1$. The quiver symmetry gives
$Q_{-i,-j}(x,y)=Q_{i,j}(-x,-y)$. Finally,
$\zeta_0=-\frac12$ and $\zeta_i\zeta_{-i}=-1$ for $i\ne0$.
These identities verify the remaining conditions in
\cite[(3.12)--(3.18)]{BWW25}.

In the following formulas and in Definition~\ref{def-2ikac}, a
circle with an index $n$ inside denotes an auxiliary generator,
not a bubble with $n$ dots. Ordinary dotted bubbles are composites
of cups, caps, and dots.
\begin{align}\label{ifakebubblepolynomial}
\begin{tikzpicture}[baseline=-1mm,iQ,draw=darkred]
\draw[fakebubble,draw=darkred] (-0.25,0) arc(180:-180:0.25);
\node at (-0.42,0) {\strandlabel{\tau i}};
\region{0.4,0}{\lambda};
\bubblelabel{0,0}{u};
\end{tikzpicture}
&:=
\sum_{n = 0}^{\del_i-\lambda_i}
\begin{tikzpicture}[baseline=-1mm,iQ,draw=darkred]
\draw[-] (-0.25,0) arc(180:-180:0.25);
\node at (-0.42,0) {\strandlabel{\tau i}};
\region{0.4,0}{\lambda};
\dottybubblelabel{0,0}{n};
\end{tikzpicture}\
u^{\del_i-\lambda_i-n} (\text{which is 0 if $i=\tau i$}),
\\\label{ibubblegeneratingfunction}
\begin{tikzpicture}[baseline=-1mm,iQ,draw=darkred]
\draw[-] (-0.25,0) arc(180:-180:0.25);
\node at (-0.42,0) {\strandlabel{\tau i}};
\region{0.95,0}{\lambda};
\node at (.55,0) {$(u)$};
\end{tikzpicture}
&:=
\begin{dcases}
\begin{tikzpicture}[baseline=-1mm,iQ,draw=darkred]
\draw[-] (-0.25,0) arc(180:-180:0.25);
\node at (-0.42,0) {\strandlabel{\tau i}};
\region{0.6,0}{\lambda};
 \node[cplus] at (.25,0) {$-$};
\end{tikzpicture}+
\begin{tikzpicture}[baseline=-1mm,iQ,draw=darkred]
\draw[fakebubble,draw=darkred] (-0.25,0) arc(180:-180:0.25);
\node at (-0.42,0) {\strandlabel{\tau i}};
\region{0.4,0}{\lambda};
\bubblelabel{0,0}{u};
\node at (0,.3) {$\phantom{.}$};
\node at (0,-.3) {$\phantom{.}$};
\end{tikzpicture}&\text{if $i \neq \tau i$}\\
\begin{tikzpicture}[baseline=-1mm,iQ,draw=darkred]
\draw[-] (-0.25,0) arc(180:-180:0.25);
\node at (-0.42,0) {\strandlabel{\tau i}};
\region{0.6,0}{\lambda};
 \node[cplus] at (.25,0) {$-$};
\end{tikzpicture}-\frac{1}{2u} \id_{\one_\lambda}&\text{if $i  =\tau i$.}
\end{dcases}
\end{align}
\begin{rem}
We use the presentation and unnormalized bubble convention of
\cite[Definition~3.4 and (3.23)--(3.24)]{BWW25}, with the corrected
sign from \cite[(3.24), footnote~2]{BWW25v2} in \eqref{iQ}.
Let $B_i(\lambda;u)$ denote the series in
\eqref{ibubblegeneratingfunction}, whose circle is labeled $\tau i$.
The series in \cite[(3.27)]{BWW25v2} is
$(\zeta_i\gamma_i(\lambda))^{-1}B_i(\lambda;u)$.
For $i=\tau(i)$, that version includes
$-\frac{1}{2u}\mathrm{id}_{\one_\lambda}$ in its fake-bubble series;
we include this term in \eqref{ibubblegeneratingfunction} instead.
Thus \cite[(3.30)]{BWW25v2} is equivalent to \eqref{ibubblerel}.
The rational functions also differ: with $t_{i,i}=1$, the
second-version convention is $R^{(2)}_{i,j}=t_{i,j}^{-1}R_{i,j}$.
By \eqref{ci1} and \eqref{gammai},
\[
 \frac{\gamma_i(\lambda-\bar\alpha_j)}{\gamma_i(\lambda)}
 =(-1)^{a_{\tau i,j}}\frac{t_{\tau i,j}}{t_{i,j}}.
\]
Substituting these identities into \cite[(3.32)]{BWW25v2}
cancels the scalar factors and gives \eqref{ibubslide}.
\end{rem}

We write $E_i$ in place of $B_i$.

Diagrams are read from bottom to top. Horizontal juxtaposition
represents composition of $1$-morphisms, with the rightmost factor
acting first. A strand changes color from $i$ to $\tau(i)$ at a cup
or cap. We label one region by its weight; all other region labels
are determined by the strand colors. Relations with omitted labels
hold for every compatible labeling. Unless the marked positions are explicitly numbered or specified
in the text, the arguments of a polynomial attached to a dashed
connector are ordered by increasing vertical coordinate and,
at equal height, from left to right. Thus a horizontal connector
is read from left to right. This convention applies to straight
and curved connectors alike. All contacts with solid strands,
including interior contacts, are counted. The ordering is local
to each connector; explicitly assigned variables are retained.

\begin{Defn}\label{def-2ikac}
The $\imath$-Kac--Moody $2$-category $\mathfrak U^\imath$ is the
underlying ungraded strict $\Bbbk$-linear $2$-category of
\cite[Definition~3.4]{BWW25}, with the sign convention in
\eqref{iQ}. It has object set $X^\imath$,
generating 1-morphisms
    $$
E_i \one_\lambda = \one_{\lambda -\bar\alpha_i} E_i
: \lambda \rightarrow \lambda-\bar\alpha_i
$$
for $\lambda \in X^\imath$ and $i \in I$,
with identity 2-endomorphisms
denoted by unoriented strings
$\begin{tikzpicture}[iQ,anchorbase,draw=darkred]
\draw[-] (0,-0.2)\botlabel{i} -- (0,0.2);
\region{0.2,0}{\lambda};
\region{-0.4,0}{\lambda-\bar\alpha_i};
\end{tikzpicture}$,
and generating 2-morphisms for all $i,j \in I$ and $\lambda \in X^\imath$ as follows
\[ \begin{tikzpicture}[iQ,centerzero,draw=darkred]
\draw[-] (0,-0.3) \botlabel{i} -- (0,0.3) \toplabel{i};
\closeddot{0,0};
\region{0.2,0}{\lambda};
\end{tikzpicture}\ \,
\colon E_i \one_\lambda \Rightarrow E_i \one_\lambda,
\quad \begin{tikzpicture}[iQ,centerzero,draw=darkred]
\draw[-] (-0.25,-0.15) \botlabel{\tau i} to [out=90,in=90,looseness=3](0.25,-0.15) \botlabel{i};
\region{0.45,0.1}{\lambda};
\node at (0,.3) {$\phantom.$};
\node at (0,-.4) {$\phantom.$};
\end{tikzpicture}\hspace{-.03pt}
\colon E_{\tau i} E_{i} \one_\lambda \Rightarrow \one_\lambda, \quad \begin{tikzpicture}[iQ,centerzero,draw=darkred]
\draw[-] (-0.25,0.15) \toplabel{\tau i} to[out=-90,in=-90,looseness=3] (0.25,0.15) \toplabel{i};
\region{0.45,-.1}{\lambda};
\node at (0,.2) {$\phantom.$};\node at (0,-.3) {$\phantom.$};
\end{tikzpicture}
\colon \one_\lambda \Rightarrow E_{\tau i} E_i \one_\lambda \]
\[\begin{tikzpicture}[baseline=-1mm,iQ,draw=darkred]
\draw[-] (-0.25,0) arc(180:-180:0.25);
\node at (-0.42,0) {\strandlabel{\tau i}};
\region{0.45,0}{\lambda};
\dottybubblelabel{0,0}{n};
\node at (0,-.3) {$\phantom{.}$};
\end{tikzpicture}
\colon \one_\lambda \Rightarrow\one_\lambda
\ \text{ for }0 \leq n \leq \del_i-\lambda_i, \quad    \begin{tikzpicture}[iQ,centerzero,scale=.9,draw=darkred]
\draw[-] (-0.3,-0.3) \botlabel{i} -- (0.3,0.3) \toplabel{i};
\draw[-] (0.3,-0.3) \botlabel{j} -- (-0.3,0.3) \toplabel{j};
\region{0.35,0}{\lambda};
\end{tikzpicture}
\;\,\colon E_{i} E_{j} \one_\lambda \Rightarrow E_j E_i \one_\lambda\]
An empty range $0\le n\le\del_i-\lambda_i$ contributes no bubble
generators. In \eqref{ifakebubblepolynomial}, a sum with negative
upper limit is also understood to be zero.
The generating $2$-morphisms are subject to the following relations:

\begin{align}
\label{ibubblerel}
\left[\!\!
\begin{tikzpicture}[baseline=-1mm,iQ,draw=darkred]
\draw[-] (-0.25,0) arc(180:-180:0.25);
\node at (-0.42,0) {\strandlabel{\tau i}};
\region{0.95,0}{\lambda};
\node at (.55,0) {$(u)$};
\end{tikzpicture}\!\!\right]_{\!u:\geq \del_i-\lambda_i}&= \zeta_i \gamma_i(\lambda) u^{\del_i-\lambda_i} \id_{\one_\lambda},\\
\label{iinfgrass}\begin{tikzpicture}[baseline=-1mm,iQ,draw=darkred]
\draw[-] (-0.25,0) arc(180:-180:0.25);
\node at (-0.42,0) {\strandlabel{\tau i}};
\node at (.55,0) {$(u)$};
\region{1,0}{\lambda};
\end{tikzpicture}\!\!
\begin{tikzpicture}[baseline=-1mm,iQ,draw=darkred]
\draw[-] (-0.25,0) arc(180:-180:0.25);
\node at (-0.35,0) {\strandlabel{i}};
\node at (.7,0) {$(-u)$};
\end{tikzpicture}\!
 &= - R_{i,\tau i}(u,-u) \id_{\one_\lambda},
 \\
\label{ibubslide}
\begin{tikzpicture}[anchorbase,iQ,draw=darkred]
\draw[-] (-0.25,0) arc(180:-180:0.25);
\node at (-.42,0) {\strandlabel{\tau i}};
\node at (.55,0) {$(u)$};
\draw[-] (1,-0.5)\botlabel{j} to (1,0.5);
\pin{1,0}{1.95,0}{R_{i, j}(u,x)};
\end{tikzpicture}
&=\begin{tikzpicture}[anchorbase,iQ,draw=darkred]
\draw[-] (-0.7,-0.5)\botlabel{j} to (-0.7,0.5);
\pin{-.7,0}{-1.85,0}{R_{\tau i, j}(-u,x)};
\draw[-] (-0.25,0) arc(180:-180:0.25);
\node at (-0.42,0) {\strandlabel{\tau i}};
\node at (.55,0) {$(u)$};
\end{tikzpicture}\\
 \label{izigzag}\begin{tikzpicture}[iQ,centerzero,scale=1.2,draw=darkred]
\draw[-] (-0.3,0.4) -- (-0.3,0) arc(180:360:0.15) arc(180:0:0.15) -- (0.3,-0.4) \botlabel{i};
\end{tikzpicture}
&=
\begin{tikzpicture}[iQ,centerzero,scale=1.2,draw=darkred]
\draw[-] (0,-0.4)\botlabel{i} -- (0,0.4);
\end{tikzpicture}\ =
\begin{tikzpicture}[iQ,centerzero,scale=1.2,draw=darkred]
\draw[-] (-0.3,-0.4) \botlabel{i}-- (-0.3,0) arc(180:0:0.15) arc(180:360:0.15) -- (0.3,0.4);
\end{tikzpicture}\ ,
\\\label{ipivotal}
\begin{tikzpicture}[iQ,baseline=-3pt,scale=1,draw=darkred]
\draw[-] (-.3,.25)\toplabel{i} to [out=-90,in=-90,looseness=3] (.3,.25);
\closeddot{-.27,.03};
\end{tikzpicture}\ &=
-\
\begin{tikzpicture}[iQ,baseline=-3pt,scale=1,draw=darkred]
\draw[-] (-.3,.25)\toplabel{i} to [out=-90,in=-90,looseness=3] (.3,.25);
\closeddot{.27,.03};
\end{tikzpicture}
\ ,\qquad
\begin{tikzpicture}[iQ,centerzero,scale=1,draw=darkred]
\draw[-] (-.3,-.25) \botlabel{i} to [out=90,in=90,looseness=3] (.3,-.25);
\closeddot{-.27,-.03};
\end{tikzpicture}\ =
-\begin{tikzpicture}[iQ,centerzero,scale=1,draw=darkred]
\draw[-] (-.3,-.25)\botlabel{i} to [out=90,in=90,looseness=3] (.3,-.25);
\closeddot{.27,-.03};
\end{tikzpicture}\ ,\\
\begin{tikzpicture}[iQ,baseline=-5pt,scale=1.2,draw=darkred]
\draw[-] (-0.25,0.15)\toplabel{i}  to[out=-90,in=-90,looseness=3] (0.25,0.15);
\draw[-] (0.3,-0.4)to[out=up,in=down] (0,0.15)\toplabel{j};
\end{tikzpicture}&=
\ \begin{tikzpicture}[iQ,baseline=-5pt,scale=1.2,draw=darkred]
\draw[-] (-0.25,0.15) \toplabel{i} to[out=-90,in=-90,looseness=3] (0.25,0.15);
\draw[-] (-0.3,-0.4) to[out=up,in=down] (0,0.15)\toplabel{j};
\end{tikzpicture}\ ,
\qquad
\begin{tikzpicture}[iQ,baseline=.5pt,scale=1.2,draw=darkred]
\draw[-] (-0.25,-0.15)  to[out=90,in=90,looseness=3] (0.25,-0.15)\botlabel{i};
\draw[-] (-0.3,0.4) \braiddown (0,-0.15)\botlabel{j};
\end{tikzpicture}
=\ \begin{tikzpicture}[iQ,baseline=.5pt,scale=1.2,draw=darkred]
\draw[-] (-0.25,-0.15)  to[out=90,in=90,looseness=3] (0.25,-0.15)\botlabel{i};
\draw[-] (0.3,0.4) \braiddown (0,-0.15)\botlabel{j};
\end{tikzpicture}\ ,\label{ipitchfork}\\
\label{idotslide}
\begin{tikzpicture}[iQ,centerzero,scale=1.2,draw=darkred]
\draw[-] (-0.3,-0.3) \botlabel{i} -- (0.3,0.3);
\draw[-] (0.3,-0.3) \botlabel{j} -- (-0.3,0.3);
\closeddot{-0.15,-0.15};
\end{tikzpicture}
-
\begin{tikzpicture}[iQ,centerzero,scale=1.2,draw=darkred]
\draw[-] (-0.3,-0.3) \botlabel{i} -- (0.3,0.3);
\draw[-] (0.3,-0.3) \botlabel{j} -- (-0.3,0.3);
\closeddot{0.15,0.15};
\end{tikzpicture}
&= \delta_{i,j}\
\begin{tikzpicture}[iQ,centerzero,scale=1.2,draw=darkred]\draw[-] (-0.2,-0.3) \botlabel{i} -- (-0.2,0.3);
\draw[-] (0.2,-0.3) \botlabel{j} -- (0.2,0.3);
\end{tikzpicture}-
\delta_{i,\tau j}\
\begin{tikzpicture}[iQ,centerzero,scale=1.2,draw=darkred]
\draw[-] (-0.2,-0.3) \botlabel{i} to [looseness=2,out=90,in=90] (0.2,-0.3);
\draw[-] (0.2,0.3) to [looseness=2,out=-90,in=-90] (-0.2,0.3)\toplabel{j};
\end{tikzpicture}=
\begin{tikzpicture}[iQ,centerzero,scale=1.2,draw=darkred]
\draw[-] (-0.3,-0.3) \botlabel{i} -- (0.3,0.3);
\draw[-] (0.3,-0.3) \botlabel{j} -- (-0.3,0.3);
\closeddot{-0.15,0.15};
\end{tikzpicture}
-
\begin{tikzpicture}[iQ,centerzero,scale=1.2,draw=darkred]
\draw[-] (-0.3,-0.3) \botlabel{i} -- (0.3,0.3);
\draw[-] (0.3,-0.3) \botlabel{j} -- (-0.3,0.3);
\closeddot{0.15,-0.15};
\end{tikzpicture}\ ,\\\label{icurl}
\begin{tikzpicture}[anchorbase,scale=1.1,iQ,draw=darkred]
\draw[-] (0,-0.5)\botlabel{i} to[out=up,in=180] (0.3,0.2) to[out=0,in=up] (0.45,0) to[out=down,in=0] (0.3,-0.2) to[out=180,in=down] (0,0.5);
\end{tikzpicture}
&=
\left[\
\begin{tikzpicture}[anchorbase,scale=1.1,iQ,draw=darkred]
\draw[-] (-0.8,-0.5)\botlabel{i} -- (-0.8,0.5);
\node[cplus] at (-.8,0) {$+$};
\draw[fakebubble,draw=darkred] (-0.4,0) arc(180:-180:0.2);
\node at (-0.48,0) {\strandlabel{i}};
\bubblelabel{-.2,0}{u};
\end{tikzpicture}\
\right]_{u:-1},\\
\label{iquadratic}
\begin{tikzpicture}[iQ,centerzero,scale=1.2,draw=darkred]
\draw[-] (-0.2,-0.4) \botlabel{i} to[out=45,in=down] (0.15,0) to[out=up,in=-45] (-0.2,0.4);
\draw[-] (0.2,-0.4) \botlabel{j} to[out=135,in=down] (-0.15,0) to[out=up,in=225] (0.2,0.4);
\end{tikzpicture}
&=
\begin{tikzpicture}[iQ,centerzero,scale=1.4,draw=darkred]
\draw[-] (-0.2,-0.3) \botlabel{i} -- (-0.2,0.3);
\draw[-] (0.2,-0.3) \botlabel{j} -- (0.2,0.3);
\pinpin{0.2,0}{-0.2,0}{-.95,0}{Q^\imath_{i,j}(x,y)};
\end{tikzpicture}
+
\delta_{i,\tau j}
\left[\begin{tikzpicture}[iQ,centerzero,scale=1.4,draw=darkred]
\draw[-] (-0.2,-0.3) \botlabel{i} to [looseness=2.2,out=90,in=90] (0.2,-0.3);
\draw[-] (0.2,0.3) to [looseness=2.2,out=-90,in=-90] (-0.2,0.3)\toplabel{i};
 \node[cplus] at (-.17,.14) {$-$};
 \node[cplus] at (-.17,-.14) {$-$};
\draw[fakebubble,draw=darkred] (0.4,0) arc(180:-180:0.16);
\node at (0.3,0) {\strandlabel{j}};
\bubblelabel{0.56,0}{u};
\end{tikzpicture}\right]_{u:-1}\!\!,
\\\label{ibraid}
\begin{tikzpicture}[iQ,centerzero,scale=1.1,draw=darkred]
\draw[-] (-0.4,-0.4) \botlabel{i} -- (0.4,0.4);
\draw[-] (0,-0.4) \botlabel{j} to[out=135,in=down] (-0.32,0) to[out=up,in=225] (0,0.4);
\draw[-] (0.4,-0.4) \botlabel{k} -- (-0.4,0.4);
\end{tikzpicture}
\,-\,
\begin{tikzpicture}[iQ,centerzero,scale=1.1,draw=darkred]
\draw[-] (-0.4,-0.4) \botlabel{i} -- (0.4,0.4);
\draw[-] (0,-0.4) \botlabel{j} to[out=45,in=down] (0.32,0) to[out=up,in=-45] (0,0.4);
\draw[-] (0.4,-0.4) \botlabel{k} -- (-0.4,0.4);
\end{tikzpicture}
&= \delta_{i,k}
\begin{tikzpicture}[iQ,centerzero,scale=1.1,draw=darkred]
\draw[-] (-0.28,-0.4) \botlabel{i} -- (-0.28,0.4);
\draw[-] (0,-0.4) \botlabel{j} -- (0,0.4);
\draw[-] (0.28,-0.4) \botlabel{k} -- (0.28,0.4);
\pinpinpin{-.28,0}{0,0}{.28,0}{1.7,0}{
\frac{Q^\imath_{i,j}(x,y)-Q^\imath_{i,j}(z,y)}{x-z}};
\end{tikzpicture}
+
\delta_{i,\tau j}\delta_{j,\tau k}\!
\left[\begin{tikzpicture}[iQ,centerzero,scale=1.4,draw=darkred]
\draw[-] (-0.2,-0.3) \botlabel{i} to [looseness=2.2,out=90,in=90] (0.2,-0.3);
\draw[-] (0.2,0.3) to [looseness=2.2,out=-90,in=-90] (-0.2,0.3)\toplabel{i};
\draw[-] (.84,-.3) \botlabel{k} to (.84,.3);
 \node[cplus] at (.84,0) {$-$};
 \node[cplus] at (-.17,-.14) {$-$};
 \node[cplus] at (-.17,.14) {$-$};
\draw[fakebubble,draw=darkred] (0.33,0) arc(180:-180:0.16);
\node at (0.25,0) {\strandlabel{j}};
\bubblelabel{0.49,0}{u};
\end{tikzpicture}
-
\begin{tikzpicture}[iQ,centerzero,scale=1.4,draw=darkred]
\draw[-] (-0.22,-0.3) \botlabel{j} to [looseness=2.2,out=90,in=90] (0.18,-0.3);
\draw[-] (0.18,0.3) to [looseness=2.2,out=-90,in=-90] (-0.22,0.3)\toplabel{j};
\draw[-] (-.9,-.3) \botlabel{i} to (-.9,.3);
\node[cplus] at (-.9,0) {$+$};
\node[cplus] at (.15,.14) {$+$};
\node[cplus] at (.15,-.14) {$+$};
\draw[fakebubble,draw=darkred] (-0.63,0) arc(180:-180:0.16);
\node at (-0.68,0) {\strandlabel{i}};
\bubblelabel{-0.47,0}{u};
\end{tikzpicture}
\right]_{\!u:-1}\notag\\
&\hspace{-10mm}
-\delta_{i,\tau j}
\begin{tikzpicture}[iQ,centerzero,scale=1.4,draw=darkred]
\draw[-] (-0.3,0.3)  -- (0.3,-0.3)\botlabel{k};
\draw[-] (0,-0.3) [out=90,in=90,looseness=2.5] to  (-0.3,-0.3)\botlabel{i};
\draw[-] (0.3,0.3) [out=-90,in=-90,looseness=2.5] to (0,0.3) \toplabel{j};
\pinpinpinmid{-.08,-.1}{0,0}{.08,.1}{1.4,0}{
\frac{Q^\imath_{j,k}(x,y)-Q^\imath_{j, k}(z,y)}{x-z}};
\end{tikzpicture}
-\delta_{j, \tau k}
\begin{tikzpicture}[iQ,centerzero,scale=1.4,draw=darkred]
\draw[-] (0.3,0.3)  -- (-0.3,-0.3)\botlabel{i};
\draw[-] (0,-0.3) \botlabel{j} [out=90,in=90,looseness=2.5] to  (0.3,-0.3);
\draw[-] (-0.3,0.3) \toplabel{k}[out=-90,in=-90,looseness=2.5] to (0,0.3) ;
\pinpinpinmid{.08,-.1}{0,0}{-.08,.1}{1.4,0}{
\frac{Q^\imath_{j, i}(x,y)-Q^\imath_{j, i}(z,y)}{x-z}};
\end{tikzpicture}.
\end{align}
\end{Defn}

Our affine Brauer constructions use this ungraded $2$-category.
We do not assume a grading on the affine Brauer category or on its
module categories. In Section~9, the natural grading of the
$2$-category in \cite[Table~3]{BWW25} induces a grading on its
cyclotomic quotients; the algebra isomorphism then transports that
grading to the corresponding cyclotomic Brauer algebras.

\subsection{Nilpotent and generalized nilpotent $2$-representations}

A $2$-representation of $\mathfrak U^\imath$ is a strict
$\Bbbk$-linear $2$-functor to the $2$-category of $\Bbbk$-linear
categories, $\Bbbk$-linear functors, and natural transformations.
The actions considered here are ungraded.
Write $\mathcal R_\lambda$ for the category assigned to the weight
$\lambda$, and use $E_i$ and the diagrammatic generators for their
images under this $2$-functor.

We use the categorical conventions of \cite[Section~2]{BSW-heisenberg}
and \cite{GRS}. A locally finite Abelian category is an essentially
small $\Bbbk$-linear Abelian category whose objects have finite length
and whose morphism spaces are finite-dimensional. A locally Schurian
category is a category equivalent to the category of locally
finite-dimensional left modules over a locally finite-dimensional
locally unital algebra $H$. Explicitly, there is a family $\{e_a\}$
of pairwise orthogonal idempotents such that the algebra and its
modules satisfy
\[
\begin{gathered}
H=\bigoplus_{a,b}e_aHe_b,\qquad \dim_\Bbbk e_aHe_b<\infty,\\
V=\bigoplus_a e_aV,\qquad \dim_\Bbbk e_aV<\infty.
\end{gathered}
\]
In this module realization, $V$ is finitely generated if it is a
quotient of a finite direct sum of modules $He_a$.
We call the $2$-representation locally finite Abelian or locally
Schurian when all its weight categories belong to the corresponding
class. Its total category is
\[
\mathcal R=
\begin{cases}
\displaystyle\bigoplus_{\lambda\in X^\imath}\mathcal R_\lambda,
 &\text{in the locally finite Abelian case},\\[2pt]
\displaystyle\prod_{\lambda\in X^\imath}\mathcal R_\lambda,
 &\text{in the locally Schurian case}.
\end{cases}
\]
The functors $E_i:\mathcal R_\lambda\to
\mathcal R_{\lambda-\bar\alpha_i}$ extend to functors on
$\mathcal R$ by shifting the weight components.

A $2$-subrepresentation consists of full additive subcategories
$\mathcal N_\lambda\subseteq\mathcal R_\lambda$, closed under
isomorphisms and stable under all generating $1$-morphisms. Fullness
ensures that the generating $2$-morphisms restrict to this family.
The subcategories need not themselves be locally Schurian or locally
finite Abelian. When referring to finitely generated objects of
$\mathcal N_\lambda$, we always mean finitely generated in the
ambient category $\mathcal R_\lambda$.

The following definition is the integral version of
\cite[Definition~3.6]{GRS}. The nilpotent terminology is adapted
from \cite[Section~5.1]{BSW-heisenberg}; the distinction between
nilpotent and generalized nilpotent actions is the one used in
\cite{GRS}.

\begin{Defn}\label{def:integral-nilpotence}
Let $\mathcal R=(\mathcal R_\lambda)_{\lambda\in X^\imath}$ be a
locally Schurian or locally finite Abelian $2$-representation.
The $2$-representation is \emph{nilpotent} if:
\begin{itemize}
\item[(1)] for every $\lambda\in X^\imath$ and every
$V\in\mathcal R_\lambda$, we have $E_iV=0$ for all but finitely
many $i\in I$ (finite color support);
\item[(2)] for every $i\in I$, $\lambda\in X^\imath$, and finitely
generated $V\in\mathcal R_\lambda$, the dot endomorphism of $E_iV$
is nilpotent.
\end{itemize}
It is \emph{generalized nilpotent} if condition~\textup{(2)} holds
and condition~\textup{(1)} is required only for finitely generated
$V$. The same conditions define these notions for a
$2$-subrepresentation, with finite generation understood in the
ambient weight category. For the mixed affine datum in
Section~\ref{sec:mixed-affine}, the same definition is used with its
vertex set and weight subgroup in place of $I$ and $X^\imath$.
\end{Defn}

The nilpotence exponent in condition~\textup{(2)} may depend on both
$i$ and $V$; no uniform exponent is required.
Every nilpotent $2$-representation is generalized nilpotent.
The two notions agree when all objects in the weight categories
are finitely generated, in particular in the locally finite Abelian
case. They also agree when $I$ is finite, since finite color support
is then automatic.

At the $0$-type vertex, the relations involving only the color $0$
are governed by the nil-Brauer category. We therefore recall its
definition for later use.

\begin{Defn}\label{nil}\cite[Definition~2.1]{BWW24}
For $t\in\Bbbk$, the nil-Brauer category $\mathcal{NB}_t$ is the
strict $\Bbbk$-linear monoidal category, with its grading forgotten,
generated by one object $B$ and four morphisms: a dot $\begin{tikzpicture}[baseline=-0.5ex,scale=.55]
\draw[-,thick,darkred,] (0,-0.5) -- (0,0.5);
		\fill (0,0) circle (0.1);
	\end{tikzpicture}
	: B\rightarrow B$,
    a crossing $\begin{tikzpicture}[baseline=-0.5ex,scale=.5]
\draw[-,thick,darkred,] (-0.4,-0.5) -- (0.4,0.5);
\draw[-,thick,darkred,] (0.4,-0.5) -- (-0.4,0.5);
	\end{tikzpicture}
	: B\otimes  B \rightarrow  B \otimes  B$,
    a cap $\begin{tikzpicture}[baseline=-0.5ex,scale=.5]
\draw[-,thick,darkred,] (-0.4,-0.3) .. controls (-0.4,0.3) and (0.4,0.3) .. (0.4,-0.3);
	\end{tikzpicture}
	: B \otimes B\rightarrow  \mathbbm 1$,
and a cup $\begin{tikzpicture}[baseline=-0.5ex,scale=.5]
\draw[-,thick,darkred,](-0.4,0.4) .. controls (-0.4,-0.3) and (0.4,-0.3) .. (0.4,0.41);
	\end{tikzpicture}
	: \mathbbm 1 \rightarrow  B \otimes  B$.  Thus, its objects are $B^{ \otimes  n}$, $n\in \mathbb N$ with $ B^{ \otimes  0}=\mathbbm 1$.
These morphisms   are subject to the following defining relations:
\begin{multicols}{2}
\begin{enumerate}[label=(\arabic*),ref=\arabic*,leftmargin=*,itemsep=.5ex,topsep=0pt]
    \item
		\label{nil:zigzags-invertibility}
		$\begin{tikzpicture}[baseline = 10pt, scale=0.4, transform shape]
			\draw[-,thick,darkred,] (2,0) to[out=up, in=down] (2,1) to[out=up, in=right] (1.5,1.5) to[out=left,in=up] (1,1);
			\draw[-,thick,darkred] (1,1) to[out=down,in=right] (0.5,0.5) to[out=left,in=down] (0,1) to[out=up,in=down] (0,2);
		\end{tikzpicture}
		~=~
		\begin{tikzpicture}[baseline = 10pt, scale=0.4,  transform shape]
			\draw[-,darkred,thick] (0,0) to (0,1);
			\draw[-,darkred,thick] (0,1) to (0,2);
		\end{tikzpicture}
		~=~
		\begin{tikzpicture}[baseline = 10pt, scale=0.4, transform shape]
			\draw[-,darkred,thick] (2,2) to[out=down, in=up] (2,1) to[out=down, in=right] (1.5,0.5) to[out=left,in=down] (1,1);
			\draw[-,darkred,thick] (1,1) to[out=up,in=right] (0.5,1.5) to[out=left,in=up] (0,1) to[out=down,in=up] (0,0);
		\end{tikzpicture}$,
        \item $\begin{tikzpicture}[baseline = 4pt]
  \draw[-,thick,darkred] (0.2,0.2) to[out=90,in=0] (0,.4);
  \draw[-,thick,darkred] (0,0.4) to[out=180,in=90] (-.2,0.2);
\draw[-,thick,darkred] (-.2,0.2) to[out=-90,in=180] (0,0);
  \draw[-,thick,darkred] (0,0) to[out=0,in=-90] (0.2,0.2);
\end{tikzpicture}=t\,\mathrm{id}_{\mathbbm1}$
\item
		\label{nil:right-cup-down-cross}
        $\begin{tikzpicture}[baseline = 5pt, scale=0.5]
			\draw[-,darkred,thick] (0,0) to[out=up,in=down] (1,1) to[out=up,in=right] (0.5,1.5) to[out=left,in=up] (0,1) to[out=down,in=up] (1,0);
		\end{tikzpicture}=0$.
	\item  %
		\label{nil:brauer-relation-4}
$\begin{tikzpicture}[baseline = 10pt, scale=0.5, transform shape]
			\draw[-,darkred,thick] (1,0) to[out=up,in=right] (0.5,1.5) to[out=left,in=up] (0,1) to[out=down,in=up] (0,0);
			\draw[-,darkred,thick] (0.5,0) to[out=up,in=down] (1.3,1.5);
		\end{tikzpicture}
		~=~
		\begin{tikzpicture}[baseline = 10pt, scale=0.5, transform shape]
			\draw[-,darkred,thick] (1,0) to[out=up,in=right] (0.5,1.5) to[out=left,in=up] (0.2,1) to[out=down,in=up] (0.2,0);
			\draw[-,darkred,thick] (0.7,0) to[out=up,in=down] (0,1.5);
		\end{tikzpicture}$.
\item%
		\label{nil:AB-relations}
 $\begin{tikzpicture}[baseline = 7.5pt, scale=0.4]
			\draw[-,darkred,thick] (0,0) to[out=up, in=down] (1,2);
			\draw[-,darkred,thick] (0,2) to[out=up, in=down] (0,2.2);
			\draw[-,darkred,thick] (1,0) to[out=up, in=down] (0,2);
			\draw[-,darkred,thick] (1,2) to[out=up, in=down] (1,2.2);
			\node at (0,1.9) {$\scriptstyle\bullet$};
		\end{tikzpicture}
		~-~
		\begin{tikzpicture}[baseline = 7.5pt, scale=0.4]
			\draw[-,darkred,thick] (0,0) to[out=up, in=down] (1,2);
			\draw[-,darkred,thick] (0,0) to[out=up, in=down] (0,-0.2);
			\draw[-,darkred,thick] (1,0) to[out=up, in=down] (0,2);
			\draw[-,darkred,thick] (1,0) to[out=up, in=down] (1,-0.2);
			\node at (1,0.1){ $\scriptstyle\bullet$};
		\end{tikzpicture}
		~=~
		-\begin{tikzpicture}[baseline = 10pt, scale=0.4, transform shape]
			\draw[-,darkred,thick] (2,2) to[out=down,in=right] (1.5,1.5) to[out=left,in=down] (1,2);
			\draw[-,darkred,thick] (2,0) to[out=up, in=right] (1.5,0.5) to[out=left,in=up] (1,0);
		\end{tikzpicture}
		~+~
		\begin{tikzpicture}[baseline = 7.5pt, scale=0.4, transform shape]
			\draw[-,darkred,thick] (0,0) to[out=up, in=down] (0,2);
			\draw[-,darkred,thick] (1,0) to[out=up, in=down] (1,2);
		\end{tikzpicture}$,
  \item  %
		\label{nil:dot-cup-cap}
\begin{tikzpicture}[baseline=-0.5mm]
		\draw[-,darkred,thick]
		(0,0) to[out=up,in=left] (0.28,0.28) to[out=right,in=up] (0.56,0);
		\draw[-,darkred,thick] (0,-0.2) -- (0,0);
		\draw[-,darkred,thick] (0.56,0) -- (0.56,-0.2);

		\node at (0,0) {$\scriptstyle{\bullet}$};
	\end{tikzpicture}
	~=~-\:
	\begin{tikzpicture}[baseline=-0.5mm]
		\draw[-,darkred,thick]
		(0,0) to[out=up,in=left] (0.28,0.28) to[out=right,in=up] (0.56,0);
		\draw[-,darkred,thick] (0.56,0) -- (0.56,-0.2);
		\draw[-,darkred,thick] (0,-0.2) -- (0,0);

		\node at (0.56,0) {$\scriptstyle{\bullet}$};
	\end{tikzpicture},
  \tikzset{strand/.style={-,thick}}

    \item\label{nil:square-zero}
     $
	\tikzset{strand/.style={-,thick}}
	\begin{tikzpicture}[baseline=8pt,scale=0.5]
		\draw[strand,darkred,] (0,0) to[out=90,in=-90] (1,1);
		\draw[darkred,strand] (1,1) to[out=90,in=-90] (0,2);
		\draw[darkred,strand] (1,0) to[out=90,in=-90] (0,1);
		\draw[darkred,strand] (0,1) to[out=90,in=-90] (1,2);
	\end{tikzpicture}
	\;=\; 0$
    \item$\tikzset{strand/.style={-,thick}}
	 \begin{tikzpicture}[baseline=8pt,scale=0.5]
		\draw[darkred,strand] (0,0) to[out=90,in=-90] (2,2);
		\draw[darkred,strand] (2,0) to[out=90,in=-90] (0,2);
		\draw[darkred,strand] (1,0) to[out=90,in=-90] (0,1)
		to[out=90,in=-90] (1,2);
	\end{tikzpicture}
	\;=\;
	\begin{tikzpicture}[baseline=8pt,scale=0.5]
		\draw[darkred,strand] (0,0) to[out=90,in=-90] (2,2);
		\draw[darkred,strand] (2,0) to[out=90,in=-90] (0,2);
		\draw[darkred,strand] (1,0) to[out=90,in=-90] (2,1)
		to[out=90,in=-90] (1,2);
	\end{tikzpicture}$.
    \end{enumerate}
\end{multicols}
\end{Defn}
Only $t\in\{0,1\}$ occurs in our applications. Since
$\bar\alpha_0=0$, the $1$-morphism $E_0\one_\lambda$ is an
endomorphism of $\lambda$. The coefficient of $u^{-1}$ in
\eqref{ibubblerel}, using \eqref{ibubblegeneratingfunction}, gives
the value of the undotted zero-color circle as
\[
\frac{1-\gamma_0(\lambda)}2\,\mathrm{id}_{\one_\lambda}
=\epsilon_0(\lambda)\,\mathrm{id}_{\one_\lambda}.
\]
The zero-color relations therefore give a strict $\Bbbk$-linear
monoidal functor
\[
\mathcal{NB}_{\epsilon_0(\lambda)}
\longrightarrow \operatorname{End}_{\mathfrak U^\imath}(\lambda),
\qquad B\longmapsto E_0\one_\lambda,
\]
sending the dot, crossing, cup, and cap to the corresponding
zero-color generators; see \cite[Section~3.3]{BWW25}.
Thus the nil-Brauer relations describe the part involving only
color $0$. Relations involving both $0$ and nonzero colors must
also be checked when constructing an action of the full
$2$-category.

\section{Affine Brauer actions and integral spectral data}
\label{sec:affine-spectral}
We recall affine Brauer actions, their generalized eigenspaces,
and the integral weight decomposition, using the normalization
of \cite{GRS}.

Throughout this section, $\Bbbk$ is algebraically closed of
characteristic different from $2$. From Assumption~\ref{assumi}
onward, we assume $\operatorname{char}\Bbbk=0$ and $I=\mathbb Z$.

 \subsection{The affine Brauer category}
 All diagrams below
are read from bottom to top, and vertical juxtaposition denotes composition.

 \begin{Defn}\label{C defn}
Let  $\mathcal {AB}$ be the $\Bbbk$-linear strict monoidal category generated by a single object, denoted by
$\begin{tikzpicture}[baseline = 10pt, scale=0.5]
	\draw[-,thick] (0,0.5)to[out=up,in=down](0,1.2);
\end{tikzpicture}$.
The objects are indexed by $\mathbb N=\{0,1,2,\ldots\}$, with $n$ corresponding to $
\begin{tikzpicture}[baseline = 10pt, scale=0.5]
	\draw[-,very thick] (0,0.5)to[out=up,in=down](0,1.2);
\end{tikzpicture}^{\otimes n}$.
In particular, $0$ corresponds to the unit object.
The generating morphisms are $\begin{tikzpicture}[baseline=-0.5ex,scale=.5]
		\draw[ thick] (-0.4,0.4) .. controls (-0.4,-0.3) and (0.4,-0.3) .. (0.4,0.41);
	\end{tikzpicture}
	: 0 \to 2$,
	$	\begin{tikzpicture}[baseline=-0.5ex,scale=.5]
		\draw[ thick] (-0.4,-0.3) .. controls (-0.4,0.3) and (0.4,0.3) .. (0.4,-0.3);
	\end{tikzpicture}
	: 2 \to 0$,
$	\begin{tikzpicture}[baseline=-0.5ex,scale=.5]
		\draw[ thick] (-0.4,-0.5) -- (0.4,0.5);
		\draw[ thick] (0.4,-0.5) -- (-0.4,0.5);
	\end{tikzpicture}
	: 2 \to 2$, and $
	\begin{tikzpicture}[baseline=-0.5ex,scale=.55]
		\draw[ thick] (0,-0.5) -- (0,0.5);
		\fill (0,0) circle (0.1);
	\end{tikzpicture}
	: 1 \to 1 $.
These morphisms   are subject to the following defining relations:
\begin{multicols}{2}
\begin{enumerate}[label=(\arabic*),ref=\arabic*,leftmargin=*,itemsep=.5ex,topsep=0pt]
    \item
		\label{affine:zigzags-invertibility}
		$\begin{tikzpicture}[baseline = 10pt, scale=0.4, transform shape]
			\draw[-,thick] (2,0) to[out=up, in=down] (2,1) to[out=up, in=right] (1.5,1.5) to[out=left,in=up] (1,1);
			\draw[-,thick] (1,1) to[out=down,in=right] (0.5,0.5) to[out=left,in=down] (0,1) to[out=up,in=down] (0,2);
		\end{tikzpicture}
		~=~
		\begin{tikzpicture}[baseline = 10pt, scale=0.4,  transform shape]
			\draw[-,thick] (0,0) to (0,1);
			\draw[-,thick] (0,1) to (0,2);
		\end{tikzpicture}
		~=~
		\begin{tikzpicture}[baseline = 10pt, scale=0.4, transform shape]
			\draw[-,thick] (2,2) to[out=down, in=up] (2,1) to[out=down, in=right] (1.5,0.5) to[out=left,in=down] (1,1);
			\draw[-,thick] (1,1) to[out=up,in=right] (0.5,1.5) to[out=left,in=up] (0,1) to[out=down,in=up] (0,0);
		\end{tikzpicture}$,
	 \item
		\label{affine:right-cup-down-cross}
		$\begin{tikzpicture}[baseline = 5pt, scale=0.5]
			\draw[-,thick] (0,1) to[out=down,in=left] (0.5,0.35) to[out=right,in=down] (1,1);
		\end{tikzpicture}
		~=~
		\begin{tikzpicture}[baseline = 5pt, scale=0.5]
			\draw[-,thick] (0,1) to[out=down,in=up] (1,0) to[out=down,in=right] (0.5,-0.5) to[out=left,in=down] (0,0) to[out=up,in=down] (1,1);
		\end{tikzpicture}$, \hskip.9cm
		$\begin{tikzpicture}[baseline = 5pt, scale=0.5]
			\draw[-,thick] (0,0) to[out=up,in=left] (0.5,0.65) to[out=right,in=up] (1,0);
		\end{tikzpicture}
		~=~
		\begin{tikzpicture}[baseline = 5pt, scale=0.5]
			\draw[-,thick] (0,0) to[out=up,in=down] (1,1) to[out=up,in=right] (0.5,1.5) to[out=left,in=up] (0,1) to[out=down,in=up] (1,0);
		\end{tikzpicture}$,
	\item  %
		\label{affine:brauer-relation-4}
		$\begin{tikzpicture}[baseline = 10pt, scale=0.5, transform shape]
			\draw[-,thick] (1,0) to[out=up,in=right] (0.5,1.5) to[out=left,in=up] (0,1) to[out=down,in=up] (0,0);
			\draw[-,thick] (0.5,0) to[out=up,in=down] (1.3,1.5);
		\end{tikzpicture}
		~=~
		\begin{tikzpicture}[baseline = 10pt, scale=0.5, transform shape]
			\draw[-,thick] (1,0) to[out=up,in=right] (0.5,1.5) to[out=left,in=up] (0.2,1) to[out=down,in=up] (0.2,0);
			\draw[-,thick] (0.7,0) to[out=up,in=down] (0,1.5);
		\end{tikzpicture}$, \hskip.5cm $
		\begin{tikzpicture}[baseline = 10pt, scale=0.5, transform shape]
			\draw[-,thick]  (2,1.5) to[out=down, in=right] (1.5,0) to[out=left,in=down] (1,1.5);
			\draw[-,thick] (0.7,0) to[out=up,in=down] (1.5,1.5);
		\end{tikzpicture}
		~=~
		\begin{tikzpicture}[baseline = 10pt, scale=0.5, transform shape]
			\draw[-,thick]  (2,1.5) to[out=down, in=right] (1.5,0) to[out=left,in=down] (1,1.5);
			\draw[-,thick] (2.3,0)to[out=up,in=down](1.5,1.5);
		\end{tikzpicture}$,

    \item%
		\label{affine:AB-relations}
		$\begin{tikzpicture}[baseline = 7.5pt, scale=0.4]
			\draw[-,thick] (0,0) to[out=up, in=down] (1,2);
			\draw[-,thick] (0,2) to[out=up, in=down] (0,2.2);
			\draw[-,thick] (1,0) to[out=up, in=down] (0,2);
			\draw[-,thick] (1,2) to[out=up, in=down] (1,2.2);
			\node at (0,1.9) {$\scriptstyle\bullet$};
		\end{tikzpicture}
		~-~
		\begin{tikzpicture}[baseline = 7.5pt, scale=0.4]
			\draw[-,thick] (0,0) to[out=up, in=down] (1,2);
			\draw[-,thick] (0,0) to[out=up, in=down] (0,-0.2);
			\draw[-,thick] (1,0) to[out=up, in=down] (0,2);
			\draw[-,thick] (1,0) to[out=up, in=down] (1,-0.2);
			\node at (1,0.1){ $\scriptstyle\bullet$};
		\end{tikzpicture}
		~=~
		-\begin{tikzpicture}[baseline = 10pt, scale=0.4, transform shape]
			\draw[-,thick] (2,2) to[out=down,in=right] (1.5,1.5) to[out=left,in=down] (1,2);
			\draw[-,thick] (2,0) to[out=up, in=right] (1.5,0.5) to[out=left,in=up] (1,0);
		\end{tikzpicture}
		~+~
		\begin{tikzpicture}[baseline = 7.5pt, scale=0.4, transform shape]
			\draw[-,thick] (0,0) to[out=up, in=down] (0,2);
			\draw[-,thick] (1,0) to[out=up, in=down] (1,2);
		\end{tikzpicture}$,

  \item  %
		\label{affine:dot-cup-cap}
		$\begin{tikzpicture}[baseline = 10pt, scale=0.4]
			\draw[-,thick] (2,0) to[out=up, in=down] (2,1) to[out=up, in=right] (1.5,1.5) to[out=left,in=up] (1,1);
			\draw[-,thick] (1,1) to[out=down,in=right] (0.5,0.5) to[out=left,in=down] (0,1) to[out=up,in=down] (0,2);
			\node at ( 1,1) {$\scriptstyle\bullet$};
		\end{tikzpicture}
		~=~
		-
		\begin{tikzpicture}[baseline = 10pt, scale=0.4]
			\draw[-,thick] (0,0) to (0,2);
			\node at (0,1) {$\scriptstyle\bullet$};
		\end{tikzpicture}
		~=~
		\begin{tikzpicture}[baseline = 10pt, scale=0.4]
			\draw[-,thick] (2,2) to[out=down, in=up] (2,1) to[out=down, in=right] (1.5,0.5) to[out=left,in=down] (1,1);
			\draw[-,thick] (1,1) to[out=up,in=right] (0.5,1.5) to[out=left,in=up] (0,1) to[out=down,in=up] (0,0);
			\node at ( 1,1) {$\scriptstyle\bullet$};
		\end{tikzpicture}$.

 \tikzset{strand/.style={-,thick}}

    \item\label{affine:symmetric-group}
    $
	\tikzset{strand/.style={-,thick}}
	\begin{tikzpicture}[baseline=8pt,scale=0.5]
		\draw[strand] (0,0) to[out=90,in=-90] (1,1);
		\draw[strand] (1,1) to[out=90,in=-90] (0,2);
		\draw[strand] (1,0) to[out=90,in=-90] (0,1);
		\draw[strand] (0,1) to[out=90,in=-90] (1,2);
	\end{tikzpicture}
	\;=\;
	\begin{tikzpicture}[baseline=8pt,scale=0.5]
		\draw[strand] (0,0) -- (0,2);
		\draw[strand] (1,0) -- (1,2);
	\end{tikzpicture}$,
    \item$\tikzset{strand/.style={-,thick}}
	\begin{tikzpicture}[baseline=8pt,scale=0.5]
		\draw[strand] (0,0) to[out=90,in=-90] (2,2);
		\draw[strand] (2,0) to[out=90,in=-90] (0,2);
		\draw[strand] (1,0) to[out=90,in=-90] (0,1)
		to[out=90,in=-90] (1,2);
	\end{tikzpicture}
	\;=\;
	\begin{tikzpicture}[baseline=8pt,scale=0.5]
		\draw[strand] (0,0) to[out=90,in=-90] (2,2);
		\draw[strand] (2,0) to[out=90,in=-90] (0,2);
		\draw[strand] (1,0) to[out=90,in=-90] (2,1)
		to[out=90,in=-90] (1,2);
	\end{tikzpicture}
	$.
\end{enumerate}
\end{multicols}
\end{Defn}
The affine Brauer category $\AB$ in Definition~\ref{C defn}	 differs from that of \cite[Definition 1.2]{RS-cyc} by the sign of the dot generator. It was proved in \cite[Lemma~2.2(3)]{GRS} that
\begin{equation}\label{AB relations1} \begin{tikzpicture}[baseline = 7.5pt, scale=0.4]
			\draw[-,thick] (0,0) to[out=up, in=down] (1,2);
			\draw[-,thick] (0,2) to[out=up, in=down] (0,2.2);
			\draw[-,thick] (1,0) to[out=up, in=down] (0,2);
			\draw[-,thick] (1,2) to[out=up, in=down] (1,2.2);
			\node at (0,0.1) {$\scriptstyle\bullet$};
		\end{tikzpicture}
		~-~
		\begin{tikzpicture}[baseline = 7.5pt, scale=0.4]
			\draw[-,thick] (0,0) to[out=up, in=down] (1,2);
			\draw[-,thick] (0,0) to[out=up, in=down] (0,-0.2);
			\draw[-,thick] (1,0) to[out=up, in=down] (0,2);
			\draw[-,thick] (1,0) to[out=up, in=down] (1,-0.2);
			\node at (1,1.9){ $\scriptstyle\bullet$};
		\end{tikzpicture}
		~=~
		-\begin{tikzpicture}[baseline = 10pt, scale=0.4, transform shape]
			\draw[-,thick] (2,2) to[out=down,in=right] (1.5,1.5) to[out=left,in=down] (1,2);
			\draw[-,thick] (2,0) to[out=up, in=right] (1.5,0.5) to[out=left,in=up] (1,0);
		\end{tikzpicture}
		~+~
		\begin{tikzpicture}[baseline = 7.5pt, scale=0.4, transform shape]
			\draw[-,thick] (0,0) to[out=up, in=down] (0,2);
			\draw[-,thick] (1,0) to[out=up, in=down] (1,2);
		\end{tikzpicture}.\end{equation} This implies the following result.
  \begin{Lemma}\label{anti} \cite[Lemma~2.6]{GRS}
 There is a $\Bbbk$-linear monoidal isomorphism $\sigma: \AB\rightarrow \AB^{\text{op}}$ that fixes  objects and the  two  generating morphisms $   \begin{tikzpicture}[baseline = 2pt, scale=0.4]
 	\draw[-,thick] (0,0) to (0,1);
 	\node at (0,.5) {$\scriptstyle\bullet$};
 \end{tikzpicture} ,
 \swap$, and swaps $\lcap$ and $\lcup$  of $\AB$.
 \end{Lemma}

\subsection{Module categories and generalized eigenspaces}
Let $\mathcal A$ be a strict $\Bbbk$-linear monoidal category.
A left $\mathcal A$-module category is a $\Bbbk$-linear category
$\mathcal M$ equipped with a strict $\Bbbk$-linear monoidal functor
$\mathrm M:\mathcal A\to\operatorname{End}_{\Bbbk}(\mathcal M)$,
where the monoidal product on the target is composition:
$\mathrm M(a\otimes b)=\mathrm M(a)\circ\mathrm M(b)$.
Throughout this section, $\mathcal M$ is locally Schurian and
carries an $\AB$-module-category structure. All categories and
functors in this section are ungraded.

\begin{Lemma} \cite[Lemma~4.2]{GRS}\label{sweet1}
Suppose that $\mathcal M$ is a module category over the  affine Brauer category $\AB$. Let $E=\mathrm M( \begin{tikzpicture}[baseline = 10pt, scale=0.5]
	\draw[-,thick] (0,0.5)to[out=up,in=down](0,1.2);
\end{tikzpicture}) $.
The cup and cap make $E$ biadjoint to itself. In particular,
$E$ is sweet, meaning that it admits a biadjoint. It is exact,
preserves finitely generated objects, and satisfies
$\dim_{\Bbbk}\operatorname{End}_{\mathcal M}(EV)<\infty$
whenever $V$ is finitely generated.
\end{Lemma}
\begin{proof}
The zigzag relations give the two adjunctions. Hence $E$ is exact
and preserves all colimits that exist. To check finite generation,
write $EV$ as the directed union of its finitely generated subobjects
$U_a$. If $V$ is finitely generated, the adjunction unit
$V\to E(EV)=\varinjlim_a EU_a$ factors through some $EU_a$.
Under adjunction this gives a map $EV\to U_a$ whose composite
with $U_a\hookrightarrow EV$ is the identity. Thus $U_a=EV$.
Finally, in a module realization of a locally Schurian category,
$\operatorname{Hom}_{\mathcal M}(W,V')$ is finite-dimensional
whenever $W$ is finitely generated: take a surjection onto $W$
from a finite direct sum of the modules $He_a$ and use
$\operatorname{Hom}_H(He_a,V')=e_aV'$. Apply this with
$W=V'=EV$.
\end{proof}

When drawing natural transformations, we omit $\mathrm{M}$ from the notation and simply write
\begin{tikzpicture}[baseline = 10pt, scale=0.6]
  \draw[-,thick] (0,0.5) to[out=up, in=down] (0,1.2);
\end{tikzpicture}
for $E$. We use the same graphical notation for morphisms in $\AB$ and their images under $\mathrm{M}$.
 Furthermore, for any morphism $g$ in $ \mathcal{AB}$ and any object $V$ in $ \mathcal M$, the morphism
$g_V $ is drawn as $\begin{tikzpicture}[baseline = 9pt,scale=0.4,inner sep=0pt, minimum width=11pt]
        \draw[-,thick] (0,0) to (0,2);
        \draw (0,1) node[circle,draw,thick,fill=white]{$g$};
        \draw[-,darkg,thick] (1,0) to (1,2);
      \node at (1,-.4) {$\darkg\scriptstyle{V}$};
    \end{tikzpicture}$.
   Usually, the rightmost factor corresponding to $V$ will be suppressed in the diagrams.

Set
$\begin{tikzpicture}[baseline=0, scale=0.7]
	\draw[thick] (0,-0.3) -- (0,0.3);
    \node at (0,0.1){$\scriptstyle\bullet$};
	\node[right] at (0.1,0.12) {$\scriptstyle{n}$};
\end{tikzpicture}
\;=\; \begin{tikzpicture}[baseline=0, scale=0.7]
	\draw[thick] (0,-0.3) -- (0,0.3);
	\node at (0,0.1){$\scriptstyle\bullet$};
\end{tikzpicture} \circ \begin{tikzpicture}[baseline=0, scale=0.7]
	\draw[thick] (0,-0.3) -- (0,0.3);
	\node at (0,0.1){$\scriptstyle\bullet$};
    \node[right] at (0.1,0.12) {$\scriptstyle{n-1}$};
\end{tikzpicture}$  for all positive  integers $n$, where
\begin{tikzpicture}[baseline=0, scale=0.7]
	\draw[thick] (0,-0.3) -- (0,0.3);
	\node at (0,0.1){$\scriptstyle\bullet$};
    \node[right] at (0.1,.12) {$\scriptstyle{0}$};\end{tikzpicture}
    is defined to be    the identity morphism.
Later, we also write   $\begin{tikzpicture}[baseline=0, scale =0.7]
\draw [thick](0,-0.3) -- (0,0.3);
\node at (0,0.1){$\scriptstyle\bullet$};
\node[right] at (0.1,0.12) {$x^n$};\end{tikzpicture}
$ for
 \begin{tikzpicture}[baseline=0, scale=0.7]
	\draw[thick] (0,-0.3) -- (0,0.3);
   \node at (0,0.1){$\scriptstyle\bullet$};
	\node[right] at (0.1,0.12) {$\scriptstyle{n}$};
\end{tikzpicture}.
This convention  allows  us to
represent a $\Bbbk$-linear combination of monomials by  labeling dots with   polynomials in
$x$.

For a commutative $\Bbbk$-algebra $R$, write $R((u^{-1}))$
for the ring of formal Laurent series at $u=\infty$. Define
\begin{equation} \label{ou1} O(u) = u-\frac{1}{2}+\sum_{r\ge 0}
\begin{tikzpicture}[baseline = 1.25mm]
		\draw[-] (0,0.4) to[out=180,in=90] (-.2,0.2);
		\draw[-] (0.2,0.2) to[out=90,in=0] (0,.4);
		\draw[-] (-.2,0.2) to[out=-90,in=180] (0,0);
		\draw[-] (0,0) to[out=0,in=-90] (0.2,0.2);
		\node at (0.2,0.2) {$\scriptstyle\bullet$};
		\node at (0.45,0.2) {$\scriptstyle{\frac{x^r}{u^r}}$};
\end{tikzpicture} \in  \End_{\AB}(0)((u^{-1})).\end{equation}

\begin{Lemma} \label{ou}
The following identities hold in $\End_{\AB}(0)((u^{-1}))$
and $\End_{\AB}(1)((u^{-1}))$, respectively; see
\cite[Lemmas~4.1--4.2]{GRS-block} and the formulation over
$\Bbbk$ in \cite[Lemma~2.9]{GRS}.
\begin{multicols}{2}
\begin{enumerate}[label=(\arabic*),ref=\arabic*,leftmargin=*,itemsep=.5ex,topsep=0pt]
\item  $
    O (u) O (-u)=(\frac{1}{2}-u)(\frac{1}{2}+u)$,
    \item $\mathord{
\begin{tikzpicture}[baseline = -1mm]
	\draw[-] (0.08,-.3) to (0.08,.3);
      \node at (-0.4,0) {$O \scriptstyle (u)$};
\end{tikzpicture}
  }=
\mathord{
\begin{tikzpicture}[baseline = -1mm]
	\draw[-] (0.08,-.3) to (0.08,.3);
      \node at (0.6,0) {$O \scriptstyle (u)$};
      \node at (.08,0) {$\scriptstyle\bullet$};
      \node at (-1.2,0) {$\scriptstyle{\frac{(u+x)^2-1}{(u-x)^2-1}\frac{(u-x)^2}{(u+x)^2}}$};
\end{tikzpicture}}$.\end{enumerate}
\end{multicols}
\end{Lemma}

Let $V$ be a finitely generated object of $\mathcal M$.
By Lemma~\ref{sweet1}, the dot endomorphism $x_V=
\begin{tikzpicture}[baseline=-0.5ex,scale=1]
            \draw[thick] (0,-0.3) -- (0,0.3);
            \fill (0,0) circle (0.06);
        \end{tikzpicture}:EV\to EV$
has a monic minimal polynomial $m_V(u)\in\Bbbk[u]$.
We set $m_V(u)=1$ when $EV=0$. For every $i\in\Bbbk$, let
$\eps_i(V)$ denote the multiplicity of $u-i$ in $m_V(u)$.
Only finitely many of these multiplicities are nonzero.

\begin{Defn}\label{subf}\cite[Definition~4.4]{GRS}
For $i\in\Bbbk$, define $E_i$ on objects as follows.
\begin{itemize}
\item[(1)] If $V$ is finitely generated, set
\[
E_iV=\ker\bigl((x_V-i\operatorname{id}_{EV})^N\bigr)
\qquad(N\geq\eps_i(V),\ N\geq1).
\]
This kernel is independent of $N$ in the indicated range.
\item[(2)] For arbitrary $V$, let $V_a$ run through its finitely
generated subobjects and set
$E_iV=\varinjlim_a E_iV_a$.
\end{itemize}
\end{Defn}
Every object is the directed union of its finitely generated
subobjects. Since $E$ is exact and preserves these unions, all
$EV_a$ can be viewed as subobjects of $EV$. Naturality of the dot
shows that the inclusions preserve generalized eigenspaces. In a
module realization, every element of $EV$ belongs to some $EV_a$,
and consequently
\[
E_iV=\bigcup_{N\geq1}\ker\bigl((x_V-i\operatorname{id}_{EV})^N\bigr)
\quad\text{inside }EV.
\]
This description proves independence of the chosen directed union
and shows that every morphism $V\to W$ induces a morphism
$E_iV\to E_iW$. Primary decomposition on finitely generated objects,
followed by passage to these unions, gives a natural decomposition
$E=\bigoplus_{i\in\Bbbk}E_i$. For each finitely generated $V$, only
finitely many summands $E_iV$ are nonzero, and each is finitely
generated because it is a direct summand of $EV$.

The dot--cup and dot--cap relations force the cup and cap components
with colors $i,j$ to vanish unless $i+j=0$. Indeed, on $E_iE_jV$
for finitely generated $V$, the sum of the two dot endomorphisms
is $(i+j)$ times the identity plus a nilpotent endomorphism, so it
is invertible when $i+j\ne0$. The assertion extends to arbitrary
objects by the directed unions above. Restricting the zigzag
identities therefore gives a biadjunction between $E_i$ and $E_{-i}$;
see \cite[Lemmas~4.7--4.8]{GRS}. Thus $E_i$ is exact, preserves
finitely generated objects, and preserves all filtered colimits
that exist.

We call
\[
I=\{i\in\Bbbk\mid E_iL\ne0\text{ for some simple object }L\in\mathcal M\}
\]
the spectrum of the dot endomorphism.

\begin{Lemma}\label{endoe-dec}
We have $E=\bigoplus_{i\in I}E_i$.
\end{Lemma}
\begin{proof}
It remains to show that a nonzero $E_i$ is nonzero on a simple
object. By preservation of directed unions, $E_iV\ne0$ for some
finitely generated $V$. The nonzero finitely generated object
$E_iV$ has a simple quotient $L$. Adjunction gives a nonzero map
$V\to E_{-i}L$. The object $E_{-i}L$ is therefore nonzero and
finitely generated, so it has a simple quotient $L'$.
The other adjunction gives a nonzero map $L\to E_iL'$.
Thus $i\in I$. This is the argument of \cite[Lemma~4.10]{GRS}.
\end{proof}

The construction above does not require a particular dot spectrum.
The half-integral case $I=\frac12+\mathbb Z$ was treated in
\cite{GRS}. Here we impose the following assumption for the integral case.

\begin{Assumption}\label{assumi} $\Bbbk$ is an algebraically closed field with characteristic zero and
$I=\mathbb Z$.
\end{Assumption}

\begin{Defn} \label{mini-local} For a non-zero  finitely generated object $V$ in $\cM$, define  $n_V(u)=
\cO_V(u) m_V(u)$, where $\mathcal O_V(u):=
\mathord{
\begin{tikzpicture}[baseline = -1mm]
     \node at (0.08,0) {$\scriptstyle O(u)$};
 	\draw[-,darkg,thick] (0.48,.2) to (0.48,-.22);
     \node at (0.48,-.37) {$\darkg\scriptstyle{V}$};
\end{tikzpicture}
}$, and $m_V(u)$ is the minimal polynomial of the dot endomorphism $
        \begin{tikzpicture}[baseline=-0.5ex,scale=1]
            \draw[thick] (0,-0.3) -- (0,0.3);
            \fill (0,0) circle (0.06);
        \end{tikzpicture} : EV \to EV$.  \end{Defn}
By \cite[Lemma~4.15]{GRS}, $n_V(u)$ is a monic polynomial in $u$
whose coefficients belong to the subalgebra of
$\operatorname{End}_{\mathcal M}(V)$ generated by the endomorphisms
obtained by evaluating dotted bubbles on $V$. The relation between
these coefficients and the local ring $Z_{\lambda,V}$ will be
explained in Subsection~\ref{subsec:polynomial-representatives}.

\begin{Lemma} \label{barn-dec} \cite[Lemma~4.16]{GRS} Suppose that  $L$ is a simple object.    Then   $n_L(u)=((-1)^{a_L} u-\frac {1}{2}) m_L(-u)$, where
 $\deg m_{L}=a_{L}$.
   \end{Lemma}

\subsection{Weights for the integral  folded datum}
We use the standard realization in Subsection~\ref{subsec:integral-datum}
and the folded lattices in \eqref{i-coweight}. Thus
\begin{equation}\label{ynu}Y^\imath=\bigoplus_{i>0}\mathbb Z\alpha_i^{\theta,\vee},
\qquad
\alpha_i^{\theta,\vee}:=\alpha_i^\vee-\alpha_{-i}^\vee.
\end{equation}
We write $\overline{\lambda}$ for the image of $\lambda\in X$
in $X^\imath$. The pairing in \eqref{pair-ikm1} induces a pairing
\[
\langle\ ,\ \rangle:Y^\imath\times X^\imath\longrightarrow\mathbb Z.
\]
In particular, for $i,j>0$, we have
\begin{equation}\label{pair-3}
\left\langle
\alpha_i^{\theta,\vee},\overline{\Lambda}_j
\right\rangle
=\delta_{ij},
\qquad
\left\langle
\alpha_i^{\theta,\vee},\overline{\Lambda}_0
\right\rangle
=0.
\end{equation}
Consequently, the pairing is non-degenerate on the left, and its
right radical is $\langle\overline{\Lambda}_0\rangle
\cong\mathbb Z/2\mathbb Z$. It therefore induces a non-degenerate
pairing between $Y^\imath$ and
$X^\imath/\langle\overline{\Lambda}_0\rangle$.

The defining relations in $X^\imath$ are
$\bar\Lambda_{-i}=-\bar\Lambda_i$ for $i>0$
and $2\bar\Lambda_0=0$.
Hence every $\lambda\in X^\imath$ has a unique expression
\begin{equation}\label{lambdacoe}
\lambda
=
\sum_{i>0}\lambda_i\overline{\Lambda}_i
+
\epsilon_0(\lambda)\overline{\Lambda}_0,\end{equation}
where $
\lambda_i
=
\left\langle
\alpha_i^{\theta,\vee},\lambda
\right\rangle\in\mathbb Z$ and
$
\epsilon_0(\lambda)\in\{0,1\}$,
and only finitely many $\lambda_i$ are nonzero. The coefficient
$\epsilon_0(\lambda)$ records the parity at the $0$-type vertex.

\begin{Defn}\label{def-wei}
Let $L$ be a simple object of $\cM$. Define
\begin{itemize}\item[(1)] $\operatorname{wt}_L
:=
\sum_{i\in I}
\bigl(\eps_{i}(L)-\eps_{-i}(L)\bigr)\Lambda_i
\in X$.
\item [(2)]
$\operatorname{wt}_L^\imath
:=
-\sum_{i\in I\setminus\{0\}}
\eps_{-i}(L)\overline{\Lambda}_i
-
\eps_0(L)\overline{\Lambda}_0
\in X^\imath
$
\end{itemize}
 Here $\eps_i(L)$ is the multiplicity of $u-i$ in the minimal
polynomial $m_L(u)$ of the dot endomorphism
$
x_L:EL\longrightarrow EL$.
\end{Defn}
Since $\bar \Lambda_i=-\bar \Lambda_{-i}$, (2) is equivalent to
$\operatorname{wt}_L^\imath
=
\sum_{i>0}
\bigl(\eps_{i}(L)-\eps_{-i}(L)\bigr)
\overline{\Lambda}_i
-
\eps_0(L)\overline{\Lambda}_0$.
Here $\eps_0(L)$ is the full multiplicity of the zero root, whereas
$\epsilon_0(\operatorname{wt}_L^\imath)\in\{0,1\}$ is its parity.

Dotted bubbles act by central natural transformations. Let
$B$ be the image of $\operatorname{End}_{\AB}(0)$ in
$\operatorname{End}(\operatorname{Id}_{\mathcal M})$.
We use \emph{central block} for a nonzero generalized central-character
subcategory for this action: for a character $\chi:B\to\Bbbk$,
it consists of the objects whose simple subquotients all have
character $\chi$. The proof below constructs this decomposition
and identifies its central blocks with the folded weights that occur.

\begin{Lemma}\label{weight-pro}
Suppose that $L$ and $K$ are two simple objects of $\cM$. Then the
following statements hold.
\begin{itemize}\item [(1)]
For every $i>0$,
$\langle
\alpha_i^{\theta,\vee},\operatorname{wt}_L^\imath
\rangle
=
\left\langle
\alpha_i^\vee,\operatorname{wt}_L
\right\rangle$.

\item [(2)]
If $L$ and $K$ lie in the same central block of $\cM$, then
$\operatorname{wt}_L^\imath=\operatorname{wt}_K^\imath$.

\item [(3)]
If $K$ is a simple subquotient of $E_iL$ for some $i\in I$, then
$
\operatorname{wt}_K^\imath
=
\operatorname{wt}_L^\imath-\overline{\alpha}_i
=
\operatorname{wt}_L^\imath+\overline{\alpha}_{-i}
$.

\item [(4)]
There is a central block decomposition
$
\cM=\prod_{\lambda\in X^\imath}\cM_\lambda,
$
where $\cM_\lambda$ is the Serre subcategory of $\cM$ consisting of
the objects $V$ such that every simple subquotient $L$ of $V$
satisfies
$
\operatorname{wt}_L^\imath=\lambda
$. The nonzero $\mathcal M_\lambda$ are precisely the central blocks
for the dotted-bubble action.
If $\mathcal M$ is locally finite Abelian, the product is replaced
by the direct sum.

\item [(5)]
For every $i\in I$, the functor $E_i$ maps $\cM_\lambda$ into
$
\cM_{\lambda-\overline{\alpha}_i}
=
\cM_{\lambda+\overline{\alpha}_{-i}}
$.

\item [(6)]
If $\operatorname{wt}_L^\imath=\operatorname{wt}_K^\imath$,
then
$
\cO_L(u)=\cO_K(u)
$.
\end{itemize}
\end{Lemma}

\begin{proof}
Suppose that $i>0$. By \eqref{pair-ikm1}, \eqref{ynu},
\eqref{pair-3}, and Definition~\ref{def-wei}, we have
$$\langle
\alpha_i^{\theta,\vee},\operatorname{wt}_L^\imath
\rangle
=
\eps_{i}(L)-\eps_{-i}(L)=
\langle
\alpha_i^\vee,\operatorname{wt}_L
\rangle.
$$
This proves (1).

Suppose that $L$ and $K$ lie in the same central block. Dotted
bubbles then act on them by the same scalars, so
$\cO_L(u)=\cO_K(u)$. By Lemma~\ref{barn-dec},
$$
\cO_T(u)
=
\left(
(-1)^{a_T}u-\frac12
\right)
\frac{m_T(-u)}{m_T(u)},
\qquad
T\in\{L,K\}.
$$
After cancelling the common power of $u$ in the numerator and the
denominator, we have
\begin{equation}\label{limi}
\left.\frac{m_T(-u)}{m_T(u)}\right|_{u=0}
=
(-1)^{\eps_0(T)}.
\end{equation}
Evaluating the resulting regular functions at $u=0$ in the equality
$\cO_L(u)=\cO_K(u)$, we obtain
$
(-1)^{\eps_0(L)}
=
(-1)^{\eps_0(K)}$.
Thus $\eps_0(L)$ and $\eps_0(K)$ have the same parity, and hence
\begin{equation}\label{parlk}
\eps_0(L)\overline{\Lambda}_0
=
\eps_0(K)\overline{\Lambda}_0.
\end{equation}
For every $i>0$, the factor
$
(-1)^{a_T}u-\frac12
$
has neither a zero nor a pole at $u=i$. Therefore,
\begin{equation}\label{dif-equal}
\operatorname{ord}_{u=i}\cO_T(u)
=
\eps_{-i}(T)-\eps_i(T).
\end{equation}

Since $\mathcal O_L(u)=\mathcal O_K(u)$, equation~\eqref{dif-equal}
gives $\eps_i(L)-\eps_{-i}(L)
=\eps_i(K)-\eps_{-i}(K)$ for every $i>0$.
Together with \eqref{parlk} and Definition~\ref{def-wei}, this gives
$\operatorname{wt}^{\imath}_L
=\operatorname{wt}^{\imath}_K$.

For (6), assume that $\operatorname{wt}_L^\imath=\operatorname{wt}_K^\imath$. Definition~\ref{def-wei} implies that $\eps_0(K)$ and $\eps_0(L)$ have the same parity. On the other hand, for  every $i>0$, we have
\begin{equation}\label{dif-equ}
\eps_{-i}(L)-\eps_i(L)
=-\langle\alpha_i^{\theta,\vee},\operatorname{wt}_L^\imath\rangle
=\eps_{-i}(K)-\eps_i(K).\end{equation}
Summing \eqref{dif-equ} and reducing modulo $2$, we obtain
$a_L-\eps_0(L)\equiv
a_K-\eps_0(K)\pmod2$.
Since $\eps_0(L)\equiv\eps_0(K)\pmod2$,
it follows that $a_L\equiv a_K\pmod2$.
Therefore,
\[
\cO_K(u)=\left((-1)^{a_K} u-\frac12\right)
\frac{m_K(-u)}{m_K(u)}
=\left((-1)^{a_L} u-\frac12\right)
\frac{m_L(-u)}{m_L(u)}
=\cO_L(u),
\]
where the second equality follows from \eqref{dif-equ}, together
with $\eps_0(L)\equiv\eps_0(K)\pmod2$ and
$a_K\equiv a_L\pmod2$. This proves (6).

Suppose that $K$ is a simple subquotient of $E_iL$.
For $i\in I$, put
\[
R_i(u):=
\frac{((u+i)^2-1)(u-i)^2}
     {((u-i)^2-1)(u+i)^2}.
\]
By Lemma~\ref{ou}(2), each coefficient of
$\mathcal O_{E_iL}(u)-\mathcal O_L(u)R_i(u)$ is a
polynomial in the nilpotent endomorphism $x-i$ with zero
constant term. These coefficients are central
endomorphisms, so they induce nilpotent endomorphisms
of $K$. Since $K$ is simple, those endomorphisms vanish.
It follows that $\mathcal O_K(u)=\mathcal O_L(u)R_i(u)$.

If $i=0$, then $R_0(u)=1$. The argument proving (2)
therefore gives $\operatorname{wt}^{\imath}_K
=\operatorname{wt}^{\imath}_L$. Since $\bar\alpha_0=0$,
this proves (3) for $i=0$.

Suppose that $i\ne0$. The identity just obtained reads
\begin{equation}\label{kl-equal}
\cO_K(u)=\cO_L(u)
\frac{(u+i)^2-1}{(u-i)^2-1}
\frac{(u-i)^2}{(u+i)^2}.
\end{equation}
Since $i\in\mathbb Z$, the function $R_i(u)$ is regular
and nonzero at $u=\pm\frac12$. Neither $m_T(u)$ nor
$m_T(-u)$ vanishes at these points, for $T\in\{L,K\}$.
By Lemma~\ref{barn-dec}, the zero of $\mathcal O_T(u)$ at one of
these two points is therefore determined by
$(-1)^{a_T}u-\frac12$. Comparing the orders in \eqref{kl-equal}
gives $a_K\equiv a_L\pmod2$. Canceling this common
linear factor yields

 \begin{equation}\label{mueql}
\frac{m_K(-u)}{m_K(u)}
=\frac{(u+i)^2-1}{(u-i)^2-1}
 \frac{(u-i)^2}{(u+i)^2}
 \frac{m_L(-u)}{m_L(u)}.
\end{equation}
Now, \eqref{mueql} implies that $
\operatorname{wt}_K
=
\operatorname{wt}_L-\alpha_i+\alpha_{-i}$. For every $j>0$, part~(1) gives
$$
\left\langle
\alpha_j^{\theta,\vee},
\operatorname{wt}_K^\imath
-
\operatorname{wt}_L^\imath
+
\overline{\alpha}_i
\right\rangle=
\left\langle
\alpha_j^\vee,
-\alpha_i+\alpha_{-i}
\right\rangle
+
\left\langle
\alpha_j^\vee-\alpha_{-j}^\vee,
\alpha_i
\right\rangle
 =
-a_{ji}+a_{j,-i}+a_{ji}-a_{-j,i}
=0.$$
Here we have used $
a_{j,-i}=a_{-j,i}$  to obtain the last equality,
which follows from the compatibility of the Cartan datum with the
involution.
Since the pairing $Y^\imath\times X^\imath\to\mathbb Z$
has right radical $\langle\bar\Lambda_0\rangle$, there
exists $\eta\in\{0,1\}$ such that
\begin{equation}\label{kl-minu}
\operatorname{wt}_K^\imath
-
\operatorname{wt}_L^\imath
+
\overline{\alpha}_i
=
\eta\overline{\Lambda}_0.
\end{equation}
It remains to determine the $\overline{\Lambda}_0$-component.
The apparent singularity of $R_i(u)$ at $u=0$ when $i=\pm1$
is removable, and a direct calculation gives
\begin{equation}\label{limi-123}
R_i(0)
=
\begin{cases}
-1,& |i|=1,\\
 1,& |i|>1.
\end{cases}
=
(-1)^{\delta_{|i|,1}}.
\end{equation}
On the other hand, by the definition of $\mathcal O_T(u)$ and \eqref{limi},
\[
\mathcal O_T(0)
=
-\frac{1}{2}(-1)^{\eps_0(T)}
\qquad
(T\in\{L,K\}).
\]
Evaluating \eqref{kl-equal} at $u=0$ after cancellation, and using
\eqref{limi-123}, therefore
yields
\begin{equation}\label{319}
\eps_0(K)-\eps_0(L)
\equiv
\delta_{|i|,1}
\pmod 2.
\end{equation}
The $\overline{\Lambda}_0$-component of
$\overline{\alpha}_i$ is $\delta_{|i|,1}$ modulo $2$. Indeed,
$
\alpha_i
=
2\Lambda_i-\Lambda_{i-1}-\Lambda_{i+1}$,
so this component is nonzero precisely when $i=\pm1$.
By Definition~\ref{def-wei}, the $\overline{\Lambda}_0$-component of
$\operatorname{wt}_T^\imath$ is $-\eps_0(T)$. Hence the
$\overline{\Lambda}_0$-component of the left-hand side of~\eqref{kl-minu} is
$$
-\eps_0(K)
+
\eps_0(L)
+
\delta_{|i|,1}
\equiv
0 \pmod {2},$$
where the equality follows from~\eqref{319}. Consequently $\eta=0$ in \eqref{kl-minu}. Therefore
\[
\operatorname{wt}_K^\imath
=
\operatorname{wt}_L^\imath-\overline{\alpha}_i
=
\operatorname{wt}_L^\imath+\overline{\alpha}_{-i}.
\]
This proves~(3).

For (4), let $V$ be finitely generated. The image $B_V$ of $B$ in
$\operatorname{End}_{\mathcal M}(V)$ is a finite-dimensional
commutative algebra. Its decomposition into local algebras gives a
finite direct sum decomposition of $V$. On each summand, every
$b\in B$ acts with a single generalized eigenvalue $\chi(b)$;
equivalently, the ideal $\ker\chi$ acts nilpotently on that summand.
All its simple subquotients have character $\chi$. Parts (2) and (6)
identify these characters with the folded weights that occur.
Morphisms preserve these summands, because they commute with the
$B$-action. Passing to directed unions therefore gives a canonical
decomposition $V=\bigoplus_\lambda V_\lambda$ for every object,
with $V_\lambda\in\mathcal M_\lambda$.
Here, for a general object, the generalized-character condition
means that $\ker\chi$ acts nilpotently on each finitely generated
subobject; the exponent need not be uniform.
To justify the product of categories in (4), realize $\mathcal M$
as locally finite-dimensional $H$-modules. Each finitely generated
projective $He_a$ has only finitely many of the weight summands
just constructed. Thus, for any family $V_\lambda\in\mathcal M_\lambda$,
$e_aV_\lambda\cong\operatorname{Hom}_H(He_a,V_\lambda)$ is nonzero
for only finitely many $\lambda$. The module
$\bigoplus_\lambda V_\lambda$ is therefore locally finite-dimensional.
The canonical central block decomposition gives the claimed equivalence
with $\prod_\lambda\mathcal M_\lambda$.
If all objects have finite length, each object has only finitely
many nonzero weight summands, giving the direct sum assertion.

For (5), let $P(L)$ be the indecomposable projective cover of a
simple object $L$. Part (4) places $P(L)$ in the same weight
subcategory as $L$. The functor $E_i$ preserves finitely generated
objects and has the exact right adjoint $E_{-i}$, so $E_iP(L)$
is a finite direct sum of indecomposable projectives.
If $P(K)$ occurs as a summand, then
\[
0\ne\operatorname{Hom}_{\mathcal M}(E_iP(L),K)
\cong\operatorname{Hom}_{\mathcal M}(P(L),E_{-i}K).
\]
Since $P(L)$ has simple head $L$ and a unique maximal subobject,
every nonzero quotient of $P(L)$ has $L$ as a simple quotient.
Hence $L$ is a simple subquotient of $E_{-i}K$, and (3) gives
$\operatorname{wt}_K^\imath=\operatorname{wt}_L^\imath-\bar\alpha_i$.
Every finitely generated object of $\mathcal M_\lambda$ is a
quotient of a finite direct sum of such projective covers of
weight $\lambda$. Exactness of $E_i$ proves (5) for these objects;
preservation of directed unions proves it for all objects.
\end{proof}

For $\lambda\in X^\imath$ with $\mathcal M_\lambda\ne0$, set
$c_\lambda:=\epsilon_0(\lambda)\in\{0,1\}$.
By Definition~\ref{def-wei}, we have
\begin{equation}\label{clambda}
c_\lambda\equiv\eps_0(L)\pmod2
\end{equation}
for every simple object $L\in\mathcal M_\lambda$.

\subsection{Central algebras and polynomial factors}
\label{subsec:polynomial-representatives}

Fix $\lambda\in X^\imath$ with $\mathcal M_\lambda\ne0$.
Let $\rho_\lambda:\operatorname{End}_{\AB}(0)\to
\operatorname{End}(\operatorname{Id}_{\mathcal M_\lambda})$
be induced by the $\AB$-action, and set
$Z_\lambda:=\operatorname{im}(\rho_\lambda)$.
For every simple object $L\in\mathcal M_\lambda$, evaluation at $L$
defines a $\Bbbk$-algebra homomorphism $\chi_L:Z_\lambda\to\Bbbk$
by $b_L=\chi_L(b)\operatorname{id}_L$ for $b\in Z_\lambda$.
By Lemma~\ref{weight-pro}(6), if $L$ and $K$ are simple objects of
$\mathcal M_\lambda$, then
$
\mathcal O_L(u)=\mathcal O_K(u)$.
Since the coefficients of $\mathcal O(u)$ are the dotted bubbles
generating $\operatorname{End}_{\AB}(0)$, it follows that
$
\chi_L=\chi_K$.
We therefore denote this common character by
$
\chi_\lambda\colon Z_\lambda\longrightarrow\Bbbk$
and define
$J_\lambda:=\ker(\chi_\lambda)$.
Equivalently,
\begin{equation}\label{jlambda}
J_\lambda
=
\left\{
b\in Z_\lambda\mid
b_L=0
\text{ for every simple object }L\in\mathcal M_\lambda
\right\}.
\end{equation}
The homomorphism $\chi_\lambda$ is surjective, since scalar
multiples of the identity belong to $Z_\lambda$. Hence
\[
Z_\lambda/J_\lambda\cong\Bbbk,
\]
and consequently $J_\lambda$ is a maximal ideal of $Z_\lambda$.

\begin{Lemma}\label{evalu-na}
Let $V$ be a nonzero finitely generated object of $\mathcal M_\lambda$,
and let
$\phi_V:\operatorname{End}(\operatorname{Id}_{\mathcal M_\lambda})
\to\operatorname{End}_{\mathcal M}(V)$ be evaluation at $V$.
Then $Z_{\lambda,V}:=\phi_V(Z_\lambda)$ is a finite-dimensional
local $\Bbbk$-algebra with maximal ideal
$J_{\lambda,V}:=\phi_V(J_\lambda)$.
\end{Lemma}
\begin{proof}
Finite-dimensionality follows from finite generation of $V$.
A simple quotient of $V$ shows that
$\ker(\phi_V|_{Z_\lambda})\subseteq J_\lambda$; hence
$Z_{\lambda,V}/J_{\lambda,V}\cong\Bbbk$.
For $b\in J_\lambda$, suppose that $b_V$ has a nonzero eigenvalue.
The corresponding generalized eigenspace is a nonzero finitely
generated direct summand of $V$. A simple quotient of this summand
would give a simple object on which $b$ acts by a nonzero scalar,
contrary to \eqref{jlambda}. Thus every element of $J_{\lambda,V}$
is nilpotent. Since $Z_{\lambda,V}$ is finite-dimensional and
commutative, $J_{\lambda,V}$ is a nilpotent ideal, and the quotient
calculation shows that it is the unique maximal ideal.
This is the local algebra used in \cite[Lemma~4.24]{GRS}.
\end{proof}

For the rest of this subsection, $V$ is a nonzero finitely generated
object of $\mathcal M_\lambda$. Its residue map is
\begin{equation}\label{red-map}
\rho_V:Z_{\lambda,V}\longrightarrow\Bbbk,
\qquad \rho_V(\phi_V(b))=\chi_\lambda(b)
\quad(b\in Z_\lambda).
\end{equation}
We apply $\rho_V$ coefficientwise to polynomials and formal Laurent
series. Set $\bar n_V(u):=\rho_V(n_V(u))$ and
$\varphi_i(V):=\operatorname{ord}_{u=i}\bar n_V(u)$ for $i\in\mathbb Z$.
Every coefficient of $\mathcal O_V(u)$ is a bubble endomorphism, so
\[
\rho_V(\mathcal O_V(u))=\mathcal O_L(u),
\qquad \bar n_V(u)=\mathcal O_L(u)m_V(u)
\]
for any simple object $L\in\mathcal M_\lambda$; compare
\cite[Lemma~4.25]{GRS}. Lemma~\ref{barn-dec} gives
\begin{equation}\label{barnv}
\bar n_V(u)=\left(u-\frac{(-1)^{a_L}}2\right)
\prod_{i\in\mathbb Z}(u-i)^{\varphi_i(V)}.
\end{equation}
In particular,
\[
\varphi_i(V)=\eps_i(V)+\eps_{-i}(L)-\eps_i(L)
\quad(i\in\mathbb Z),
\qquad \varphi_0(V)=\eps_0(V).
\]
The parity of $a_L$ depends only on $\lambda$, as shown in the proof
of Lemma~\ref{weight-pro}(6), so these formulas are independent of
$L$. Whenever a quotient of polynomials over $Z_{\lambda,V}$ appears
below, its denominator is a unit in $Z_{\lambda,V}((u^{-1}))$;
we interpret the quotient in that ring.

\begin{Defn}\label{mini-1}
Let $L$ be any simple object of $\mathcal M_\lambda$, and set
$\delta_L
:=1$ if $a_L$ is even, and
$0$ if $a_L$ is odd. The value of $\delta_L$ depends only on
$\lambda$. Recall that
$
m_V(u)
=
\prod_{i\in\mathbb Z}
(u-i)^{\eps_i(V)}$.
For every $i\in I=\mathbb Z$, define \begin{itemize}\item [(1)]
$m_{V,i}(u)
:=
\frac{m_V(u)}{(u-i)^{\eps_i(V)}}$ and $
\overline n_{V,i}(u)
:=
\frac{\overline n_V(u)}{(u-i)^{\varphi_i(V)}}$.
\item [(2)] $m_{V,0}^{+}(u)
:=
\prod_{i>0}(u-i)^{\eps_i(V)}$ and $
m_{V,0}^{-}(u)
:=
\prod_{i<0}(u-i)^{\eps_i(V)}$.
\item [(3)] $
\overline n_{V,0}^{+}(u)
:=
\left(
u-\frac{(-1)^{a_L}}{2}
\right)^{\delta_L}
\prod_{i>0}(u-i)^{\varphi_i(V)}
$, $
\overline n_{V,0}^{-}(u)
:=
\left(
u-\frac{(-1)^{a_L}}{2}
\right)^{1-\delta_L}
\prod_{i<0}(u-i)^{\varphi_i(V)}$, and
$
\overline t_{V,0}(u):=u^{\varphi_0(V)}$.
\end{itemize}\end{Defn}
By Lemma~\ref{evalu-na}, $J_{\lambda,V}$ is nilpotent.
For each $i\in\mathbb Z\setminus\{0\}$, set
$\overline t_{V,i}(u):=(u-i)^{\varphi_i(V)}$. Then
$\overline n_V(u)=\overline n_{V,i}(u)\overline t_{V,i}(u)$, and these two monic factors are coprime. By the monic Hensel lifting lemma \cite[Lemma~4.28]{GRS}, there are unique monic polynomials $n_{V,i}(u),t_{V,i}(u)\in Z_{\lambda,V}[u]$, reducing respectively to $\overline n_{V,i}(u)$ and $\overline t_{V,i}(u)$, such that
\[
 n_V(u)=n_{V,i}(u)t_{V,i}(u).
\]
For $i=0$, we have
$\overline n_V(u)
=
\overline n_{V,0}^{+}(u)\,
\overline n_{V,0}^{-}(u)\,
\overline t_{V,0}(u)$.
By construction, the three monic polynomials
$\overline n^+_{V,0}(u), \overline n^-_{V,0}(u)$ and $
  \overline t_{V,0}(u)$ are
 pairwise coprime in $\Bbbk[u]$. Hence, by
\cite[Lemma~4.28]{GRS}, applied twice, there exist unique monic
polynomials
$n^+_{V,0}(u),
  n^-_{V,0}(u)$ and $
  t_{V,0}(u)$ in $ Z_{\lambda,V}[u]$
whose images under the reduction map
\[
  Z_{\lambda,V}[u]\longrightarrow
  \bigl(Z_{\lambda,V}/J_{\lambda,V}\bigr)[u]
  \simeq \Bbbk[u]
\]
are $\overline n^+_{V,0}(u)$,
$\overline n^-_{V,0}(u)$, and
$\overline t_{V,0}(u)$, respectively, and which satisfy
\[
n_V(u)=n^+_{V,0}(u)n^-_{V,0}(u)t_{V,0}(u).
\]
Set $n_{V,0}(u):=n^+_{V,0}(u)n^-_{V,0}(u)$.
Then $n_V(u)=n_{V,0}(u)t_{V,0}(u)$, and $n_{V,0}(u)$
reduces to $\bar n_{V,0}(u)$. Moreover,
$\deg t_{V,0}=\varphi_0(V)=\eps_0(V)$.
For a simple object $L$, the residue map is the identity, and it follows from Definition~\ref{mini-1} that  \begin{equation}\label{zhfu}
\begin{aligned}
n_{L,0}^+(u)
&=\prod_{i>0}(u-i)^{\eps_{-i}(L)}(u-(-1)^{a_L}1/2)^{\delta_L},\\
n_{L,0}^-(u)
&=\prod_{i<0}(u-i)^{\eps_{-i}(L)}
(u-(-1)^{a_L}1/2)^{1-\delta_L},\\
{t}_{L,0}(u)&=u^{\eps_0(L)}.
\end{aligned}
\end{equation}

\begin{Lemma}\label{lem:fg-zero}
Let $V$ be a nonzero finitely generated object of $\mathcal M_\lambda$,
and let $L$ be a simple quotient of $V$. Then:
\begin{itemize}
\item[(1)] The rational function $\mathcal O_L(u)$ is regular at $0$, and
$\mathcal O_L(0)=(-1)^{\eps_0(L)+1}/2$.
\item[(2)] In $Z_{\lambda,V}((u^{-1}))$, we have
\begin{equation}\label{eq:fg-mixed-ratio}
\frac{n_{V,0}^+(u)n_{V,0}^-(-u)}
     {m_{V,0}^+(u)m_{V,0}^-(-u)}
=(-1)^{\eps_0(L)}\left(u-\frac12\right).
\end{equation}
\end{itemize}
\end{Lemma}
\begin{proof}
Part (1) follows from Lemma~\ref{barn-dec} and \eqref{limi}.
Put $d=\deg m_V$, $d_n^-=\deg n^-_{V,0}$, and
$d_m^-=\deg m^-_{V,0}$. Lemma~\ref{ou}(1) gives
\[
n_V(u)n_V(-u)
=\left(\frac14-u^2\right)m_V(u)m_V(-u).
\]
Since $\deg n_V=d+1$, multiplication by $(-1)^{d+1}$ makes both
sides monic. Group their reduced roots into
$\mathbb Z_{>0}\cup\{\frac12\}$,
$\mathbb Z_{<0}\cup\{-\frac12\}$, and $\{0\}$.
Uniqueness in monic Hensel lifting identifies the normalized
factors belonging to the first set:
\[
(-1)^{d_n^-}n^+_{V,0}(u)n^-_{V,0}(-u)
=\left(u-\frac12\right)(-1)^{d_m^-}
 m^+_{V,0}(u)m^-_{V,0}(-u).
\]
Consequently,
\begin{equation}\label{weight111}
\frac{n_{V,0}^+(u)n_{V,0}^-(-u)}
     {m_{V,0}^+(u)m_{V,0}^-(-u)}
=(-1)^{d_n^--d_m^-}\left(u-\frac12\right).
\end{equation}
The multiplicity formula following \eqref{barnv} now gives
\begin{equation}\label{weight222}
\begin{aligned}
d_n^--d_m^-
&=1-\delta_L+\sum_{i<0}\bigl(\varphi_i(V)-\eps_i(V)\bigr)\\
&=1-\delta_L+\sum_{i<0}\bigl(\eps_{-i}(L)-\eps_i(L)\bigr)\\
&\equiv a_L+\sum_{i\ne0}\eps_i(L)
\equiv\eps_0(L)\pmod2.
\end{aligned}
\end{equation}
This proves (2).
\end{proof}

We also record the following elementary consequence of the
nilpotence of $J_{\lambda,V}$; see the paragraph following
\cite[Lemma~4.28]{GRS}.
\begin{Lemma}~\label{key-coprim}
Let $f_1(u), f_2(u)\in Z_{\lambda,V}[u]$ be  monic polynomials.
For $r=1,2$, let $\bar f_r(u)$ denote the image of $f_r(u)$ under the map
$Z_{\lambda, V}[u]\rightarrow \Bbbk[u]$.
Assume that $\bar f_1(u)$ and   $\bar f_2(u)$ are coprime. Then the ideal generated by $f_1(u)$ and $f_2(u)$ is $Z_{\lambda, V}[u]$.  Thus,
there exist polynomials $g(u),h(u)\in Z_{\lambda,V}[u]$ such that
$
g(u)f_1(u)+h(u)f_2(u)=1$.
\end{Lemma}

\begin{proof}
Lift a B\'ezout identity for $\bar f_1$ and $\bar f_2$ to obtain
$g_0f_1+h_0f_2=1+r$ with $r\in J_{\lambda,V}[u]$.
Since $J_{\lambda,V}$ is nilpotent, so is $r$, and $1+r$ has a
polynomial inverse given by a finite geometric series. Multiplying
by this inverse gives the required identity.
\end{proof}

\section{Constructing the \texorpdfstring{$\imath$}{i}-Kac--Moody generators from an affine Brauer action}
\label{sec:forward-generators}

We retain Assumption~\ref{assumi} and the notation of
Section~\ref{sec:affine-spectral}. Unless stated otherwise, $V$ is
a nonzero finitely generated object of $\mathcal M_\lambda$.
We use the quiver orientation $i\to i+1$ and
$I_1=\mathbb Z_{>0}$ from Section~2.
Write $d_i(\lambda)=\lambda_i$ for $i\in I_1$ and $c=c_\lambda$.
We first define the components on finitely generated objects,
then extend them to all objects by naturality and directed unions.
Coefficients indexed by the zero object are omitted, and the
corresponding diagrams represent zero morphisms.

\subsection{The $2$-morphisms associated with type $A$ vertices}
\begin{Defn}\label{ration111}
Let $a,b\in \mathbb{Z}_{>0}$, $\sigma,\tau\in\{+,-\}$, and $i\in I_1$.
We identify the signs $+$ and $-$ with $+1$ and $-1$ in the
following expressions. For indeterminates $x,x_1,x_2,\ldots$, set
\begin{multicols}{2}
\begin{enumerate}[label=(\arabic*),ref=\arabic*,leftmargin=*,itemsep=.5ex,topsep=0pt]
\item $\ell^{ab}_{\sigma\tau}:=1+\sigma x_a+\tau x_b$,
\item $d^{ab}_{\sigma\tau}:=\sigma x_a+\tau x_b$,
\item $a_{V,i}(x):=\frac{n_{V,i}(x)}{x m_{V,i}(x)}$,
\item $b_{V,i}(x):=\frac{x m_{V,i}(x)}{n_{V,i}(x)}$,
\item $c^\sigma_V(x):= \frac{n^\sigma_{V,0}(x)}{m^\sigma_{V,0}(x)}$,
\item $d^\sigma_V(x):= \frac{m^\sigma_{V,0}(x)}{n^\sigma_{V,0}(x)}$.
\end{enumerate}
\end{multicols}
\end{Defn}

Black strands labelled $j$ represent the generalized eigenspace
functors $E_j$. Their dots are the restrictions of the affine Brauer
dot, so $x-j$ is nilpotent on $E_jV$. In a row of parallel strands,
$x_1,x_2,\ldots$ denote the dots from left to right at the height
where the coefficient is placed. A coefficient below a crossing
acts before that crossing, and a coefficient above it acts after
the crossing. Coefficients indexed by $V$ are transported by the
functors represented by the strands to their left. They commute
with the dots used for evaluation; we do not assume that they
commute with a crossing.

A rational expression is evaluated only when its denominator is
invertible on the specified generalized eigenspace. Removable
factors are cancelled before evaluation. A divided difference is
formed as a polynomial or formal series before its variables are
replaced by dots; no inverse of a nilpotent difference is taken.

For $i>0$, an upward (respectively, downward) red strand
labeled $i$ represents $E_i$ (respectively, $E_{-i}$).
Red strands labeled $0$ are left unoriented. The oriented diagrams
use the conventions of our half-integral construction in
\cite[Definitions~5.7 and~5.9]{GRS}, restricted to the nonexceptional
type $A$ vertices. These are auxiliary generators: their
normalization is different from the unoriented presentation in
Section~2. We give the precise conversion in
Subsection~\ref{subsec:forward-normalization}.

\begin{Defn}\label{natrual-basic}
Let $i\in I_1$.  For  a nonzero finitely generated
$V\in\mathcal M_\lambda$, define the following transformations:
\begin{itemize}
\item[(1)]
    $\begin{tikzpicture}[baseline=-0.5ex,scale=0.7 ]
        \draw[thick, darkred, ->] (0,-0.5) -- (0,0.5);
        \node at (.2,.05) {$\scriptstyle{\lambda}$};
\node[darkred] at (0, -0.7) {$\scriptstyle i$};
\end{tikzpicture}
    =\begin{tikzpicture}[baseline=-0.5ex ,scale=0.7]
        \draw[thick, black, ] (0,-0.5) -- (0,0.5);
    \node at (0, -0.7) {$\scriptstyle i$};
\end{tikzpicture}$,
$%
    \begin{tikzpicture}[baseline=-0.5ex,scale=0.7 ]
        \draw[thick, darkred, <-] (0,-0.5) -- (0,0.5);
\node[darkred] at (0, -0.7) {$\scriptstyle i$};
 \node at (.2,.05) {$\scriptstyle{\lambda}$};
\end{tikzpicture}
   = \begin{tikzpicture}[baseline=-0.5ex, scale=0.7  ]
        \draw[thick, black, ] (0,-0.5) -- (0,0.5);
    \node at (0, -0.7) {$\scriptstyle  {-i}$};
\end{tikzpicture}$,
$
    \begin{tikzpicture}[baseline=-0.5ex ,scale=0.7 ]
        \draw[thick, darkred, ->] (0,-0.5) -- (0,0.5);
        \fill[darkred]  (0, 0) circle (2pt);
        \node[darkred] at (0, -0.7) {$\scriptstyle i$};
 \node at (.2,.05) {$\scriptstyle{\lambda}$};
\end{tikzpicture}
    =
\begin{tikzpicture}[baseline=-0.5ex,scale=0.7 ] \draw[thick, black, ] (0,-0.5) -- (0,0.5);
        \node at (0, 0) {$\scriptstyle\bullet$};
        \node at (0.6,0) {$\scriptstyle{x-i}$};
        \node at (0, -0.7) {$\scriptstyle i$};
\end{tikzpicture}
        $,
        $
    \begin{tikzpicture}[baseline=-0.5ex ,scale=0.7 ]
        \draw[thick, darkred, <-] (0,-0.5) -- (0,0.5);
       \fill[darkred]  (0, 0) circle (2pt);
       \node[darkred] at (0, -0.7) {$\scriptstyle i$};
        \node at (.2,.05) {$\scriptstyle{\lambda}$};
       \end{tikzpicture}
   = \begin{tikzpicture}[baseline=-0.5ex ,scale=0.7  ]
        \draw[thick, black, ] (0,-0.5) -- (0,0.5);
         \node at (0, 0) {$\scriptstyle\bullet$};
          \node at (0.75,0) {$\scriptstyle{-x-i}$};
          \node at (0, -0.7) {$\scriptstyle  -i$};
\end{tikzpicture}$,

\item [(2)]
    $\begin{tikzpicture}[baseline=-0.5ex ,scale=0.7, ]
        \draw[thick, darkred, ->]
            (-0.4,0.5) .. controls (-0.4,-0.3) and (0.4,-0.3) .. (0.4,0.5);
            \node[darkred] at (-.4, .7) {$\scriptstyle{i}$};
 \node at (.5,.05) {$\scriptstyle{\lambda}$};
\end{tikzpicture}
            =\begin{tikzpicture}[baseline=-0.5ex, scale=0.7,  ]
             \draw[thick, black]
            (-0.4,0.5) .. controls (-0.4,-0.3) and (0.4,-0.3) .. (0.4,0.5);
            \node at (-.4, .7) {$\scriptstyle{-i}$};
             \node at (.4, .7) {$\scriptstyle{i}$};
\end{tikzpicture}$,
$
 \begin{tikzpicture}[baseline=-0.5ex,scale=0.7, , yshift=10pt]
        \draw[thick, darkred, ->]
            (-0.4,-0.5) .. controls (-0.4,0.3) and (0.4,0.3) .. (0.4,-0.5);\node[darkred] at (-0.4, -0.7) {$\scriptstyle i$};
 \node at (.5,.05) {$\scriptstyle{\lambda}$};
\end{tikzpicture}
= \begin{tikzpicture}[baseline=-0.5ex ,scale=0.7, ,   yshift=10pt]
        \draw[thick, ]
            (-0.4,-0.5) .. controls (-0.4,0.3) and (0.4,0.3) .. (0.4,-0.5);
                \node at (-0.4, -0.7) {$\scriptstyle i$};
  \node at (0.4, -0.7) {$\scriptstyle -i$};
\end{tikzpicture}$,

\item[(3)]
$\begin{tikzpicture}[baseline=-0.5ex,scale=0.7,  ]
        \draw[thick, darkred, <-]
            (-0.4,0.5) .. controls (-0.4,-0.3) and (0.4,-0.3) .. (0.4,0.5);
            \node[darkred] at (.4, .7) {$\scriptstyle{i}$};
\node at (.5,.05) {$\scriptstyle{\lambda}$};
\end{tikzpicture}=
\begin{tikzpicture}[baseline=-0.5ex,scale=0.7,]
        \draw[thick, black, ]
            (-0.4,0.5) .. controls (-0.4,-0.3) and (0.4,-0.3) .. (0.4,0.5);
             \node at (-.4, .7) {$\scriptstyle{i}$};
             \node at (.4, .7) {$\scriptstyle -i$};   \node at (-.38, .3) {$\scriptstyle\bullet$};
             \node at (-1.5, .4) {$\scriptstyle
 a_{E_{-i}V,i}(x)$};
\end{tikzpicture}$,
$    \begin{tikzpicture}[baseline=-0.5ex,scale=0.7,  yshift=10pt]
        \draw[thick, darkred, <-]
            (-0.4,-0.5) .. controls (-0.4,0.3) and (0.4,0.3) .. (0.4,-0.5);
             \node[darkred] at (.4, -.7) {$\scriptstyle{i}$};
\node at (.5,.05) {$\scriptstyle{\lambda}$};
\end{tikzpicture}$
=
$\begin{tikzpicture}[baseline=-0.5ex,scale=0.7, , yshift=10pt]
        \draw[thick]
            (-0.4,-0.5) .. controls (-0.4,0.3) and (0.4,0.3) .. (0.4,-0.5);
            \node at (-0.4, -0.7) {$\scriptstyle -i$};
  \node at (0.4, -0.73) {$\scriptstyle i$};
  \node at (.4, -.35) {$\scriptstyle\bullet$};
  \node at (1.2, -.4) {$\scriptstyle b_{V,i}(x)$};
\end{tikzpicture}$.
\end{itemize}
\end{Defn}
The coefficients in Definition~\ref{natrual-basic}(3) are
\(a_{E_{-i}V,i}(x)\) and \(b_{V,i}(x)\).
For each nonzero object \(W\) occurring in these
coefficients, evaluation takes place on \(E_iW\).
Since \(i\ne0\), and the complementary factors
\(m_{W,i}\) and \(n_{W,i}\) have nonzero residues at
\(u=i\), the required denominators act invertibly.

By the multiplicity formula following~\eqref{barnv},
$\varphi_i(V)=\eps_i(V)-d_i(\lambda)$ and
$\varphi_{-i}(V)=\eps_{-i}(V)+d_i(\lambda)$ for $i\in I_1$.

A diamond-marked black crossing denotes the local component
of the affine Brauer crossing obtained by pre-composing with
the source inclusions and post-composing with the target projections for the indicated
generalized eigenspaces; see \cite[(4.5)]{GRS}.

\begin{Defn} \label{natrual-basic1}
Let $i,j\in I_1$.  On a nonzero finitely generated
$V\in\mathcal M_\lambda$, define the following transformations:
\begin{itemize}
\item [(1)]
$
\begin{tikzpicture}[baseline = 0, darkred]
	\draw[->,thick] (0.28,-.3) to (-0.28,.4);
	\draw[->,thick] (-0.28,-.3) to (0.28,.4);
    \node at (-.3,-.45) {$\scriptstyle{j}$};
    \node at (.3,-.45) {$\scriptstyle{i}$};
    \node at (0.5,0.05) {$\color{gray}\scriptstyle{\lambda}$};
\end{tikzpicture}
=
\left\{
\begin{array}{ll}
\begin{tikzpicture}[baseline = 0, scale=0.8, ]
	\draw[-] (0.38,-.4) to (-0.38,.5);
	\draw[-] (-0.38,-.4) to (0.38,.5);
    \node at (0,0.05) {$\diamond$};
	\draw[-,densely dotted] (0.2,-.24) to (-0.19,-.24);
    \node at (-.4,-.55) {$\scriptstyle{i}$};
    \node at (.4,-.55) {$\scriptstyle{i}$};
    \node at (-.4,.65) {$\scriptstyle{i}$};
    \node at (.4,.65) {$\scriptstyle{i}$};
    \node at (1.3,-0.2) {$\scriptstyle{ (\ell_{+-}^{12})^{-1}}$};
\end{tikzpicture}
+
\begin{tikzpicture}[baseline = 0, scale=0.8, ]
	\draw[-,densely dotted] (0.1,.06) to (0.55,.06);
	\draw[-] (0.58,-.4) to (0.58,.5);
	\draw[-] (0.08,-.4) to (0.08,.5);
    \node at (0.08,-0.55) {$\scriptstyle{i}$};
    \node at (0.58,-0.55) {$\scriptstyle{i}$};
    \node at (1.6,0.1) {$\scriptstyle{(\ell_{+-}^{12})^{-1}} $};
\end{tikzpicture}
&\text{if $j=i$,}
\\
\begin{tikzpicture}[baseline = 0, scale=0.8, ]
	\draw[-] (0.38,-.4) to (-0.38,.5);
	\draw[-] (-0.38,-.4) to (0.38,.5);
	\draw[-,densely dotted] (0.2,-.24) to (-0.19,-.24);
    \node at (0,0.05) {$\diamond$};
    \node at (-.4,-.55) {$\scriptstyle{i+1}$};
    \node at (-.4,.69) {$\scriptstyle{i}$};
    \node at (.4,.69) {$\scriptstyle{i+1}$};
    \node at (.4,-.55) {$\scriptstyle{i}$};
    \node at (1.0,-0.2) {$\scriptstyle{d_{+-}^{12}}$};
\end{tikzpicture}
&\text{if $j=i+1$,}
\\
\begin{tikzpicture}[baseline = 0, scale=0.8, ]
	\draw[-] (0.38,-.4) to (-0.38,.5);
    \node at (0,0.05) {$\diamond$};
	\draw[-] (-0.38,-.4) to (0.38,.5);
	\draw[-,densely dotted] (0.18,-.24) to (-0.19,-.24);
    \node at (-.4,-.62) {$\scriptstyle{j}$};
    \node at (.4,-.6) {$\scriptstyle{i}$};
    \node at (.4,.75) {$\scriptstyle{j}$};
    \node at (-.4,.73) {$\scriptstyle{i}$};
    \node at (1.0,-0.2) {$\scriptstyle{\frac{d^{12}_{+-}} {\ell_{-+}^{12}}}$};
\end{tikzpicture}
&\text{if $j \neq i,i+1$.}
\end{array}\right.$

\item[(2)] $\mathord{
\begin{tikzpicture}[baseline = 0, scale=0.8, ]
	\draw[<-,thick,darkred] (0.28,-.3) to (-0.28,.4);
	\draw[<-,thick,darkred] (-0.28,-.3) to (0.28,.4);
   \node at (-0.28,.65) {$\scriptstyle{j}$};
   \node at (0.28,.65) {$\scriptstyle{i}$};
   \node at (.4,.05) {$\scriptstyle{\lambda}$};
\end{tikzpicture}
}
=t_{ji}^{-1}\mathord{
\begin{tikzpicture}[baseline = 0, scale=0.4]
\draw[->,thick,darkred] (-1.3,.4) to (-1.3,-1.2);
\draw[-,thick,darkred] (1.3,-.4) to (1.3,1.2);
\draw[-,thick,darkred] (-.5,1.1) to [out=180,in=90] (-1.3,.4);
\draw[-,thick,darkred] (.35,.4) to [out=90,in=0] (-.5,1.1);
\draw[-,thick,darkred] (.5,-1.1) to [out=0,in=-90] (1.3,-.4);
\draw[-,thick,darkred] (-.35,-.4) to [out=-90,in=180] (.5,-1.1);
\draw[-,thick,darkred] (-.35,-.4) to [out=90,in=-90] (.35,.4);
        \draw[-,thick,darkred] (0.35,-.5) to[out=180,in=0] (-0.35,.5);
        \draw[->,thick,darkred] (-0.35,.5) to[out=180,in=90] (-0.8,-1.2);
        \draw[-,thick,darkred] (0.35,-.5) to[out=0,in=-90] (0.8,1.2);
  \node at (0.8,1.9) {$\scriptstyle{j}$};
  \node at (1.3,1.8) {$\scriptstyle{i}$};
   \node at (1.55,0) {$\scriptstyle{\lambda}$};
\end{tikzpicture}
}$,  where  $t_{ij}=
-1$ (respectively, $1$) if $ j=i-1$ (respectively, otherwise).

\item [(3)] $\mathord{
\begin{tikzpicture}[baseline = 0]
	\draw[<-,thick,darkred] (0.28,-.3) to (-0.28,.4);
	\draw[->,thick,darkred] (-0.28,-.3) to (0.28,.4);
   \node at (-0.28,-.5) {$\scriptstyle{j}$};
   \node at (-0.28,.6) {$\scriptstyle{i}$};
   \node at (.4,.05) {$\scriptstyle{\lambda}$};
\end{tikzpicture}
}
=
\mathord{
\begin{tikzpicture}[baseline = 0, scale=0.6]
	\draw[->,thick,darkred] (0.3,-.5) to (-0.3,.5);
	\draw[-,thick,darkred] (-0.2,-.2) to (0.2,.3);
    \draw[-,thick,darkred] (0.2,.3) to[out=50,in=180] (0.5,.5);
    \draw[->,thick,darkred] (0.5,.5) to[out=0,in=90] (0.8,-.5);
    \draw[-,thick,darkred] (-0.2,-.2) to[out=230,in=0] (-0.5,-.5);
    \draw[-,thick,darkred] (-0.5,-.5) to[out=180,in=-90] (-0.8,.5);
   \node at (-0.8,.73) {$\scriptstyle{i}$};
   \node at (0.28,-.74) {$\scriptstyle{j}$};
   \node at (1.05,.05) {$\scriptstyle{\lambda}$};
\end{tikzpicture}
}$,  $\quad
\mathord{
\begin{tikzpicture}[baseline = 0]
	\draw[->,thick,darkred] (0.28,-.3) to (-0.28,.4);
	\draw[<-,thick,darkred] (-0.28,-.3) to (0.28,.4);
   \node at (-0.28,-.5) {$\scriptstyle{j}$};
   \node at (-0.28,.6) {$\scriptstyle{i}$};
   \node at (.4,.05) {$\scriptstyle{\lambda}$};
\end{tikzpicture}
}
=
\mathord{
\begin{tikzpicture}[baseline = 0, scale=0.6]
	\draw[<-,thick,darkred] (0.3,-.5) to (-0.3,.5);
	\draw[-,thick,darkred] (-0.2,-.2) to (0.2,.3);
        \draw[-,thick,darkred] (0.2,.3) to[out=50,in=180] (0.5,.5);
        \draw[-,thick,darkred] (0.5,.5) to[out=0,in=90] (0.8,-.5);
        \draw[-,thick,darkred] (-0.2,-.2) to[out=230,in=0] (-0.5,-.5);
        \draw[->,thick,darkred] (-0.5,-.5) to[out=180,in=-90] (-0.8,.5);
  \node at (-0.8,.68) {$\scriptstyle{i}$};
   \node at (0.28,-.69) {$\scriptstyle{j}$};
  \node at (1.05,.05) {$\scriptstyle{\lambda}$};
\end{tikzpicture}
}$,

\item [(4)]  $
   \sum\limits_{r\in \mathbb Z}
\begin{tikzpicture}[baseline = 0]
  \draw[->,thick,darkred] (0.2,0.2) to[out=90,in=0] (0,.4);
  \draw[-,thick,darkred] (0,0.4) to[out=180,in=90] (-.2,0.2);
\draw[-,thick,darkred] (-.2,0.2) to[out=-90,in=180] (0,0);
  \draw[-,thick,darkred] (0,0) to[out=0,in=-90] (0.2,0.2);
 \node at (0,-.13) {$\scriptstyle{i}$};
   \node at (-0.3,0.2) {$\scriptstyle{\lambda}$};
   \node at (0.2,0.2) {$\color{darkred}\scriptstyle\bullet$};
   \node at (0.4,0.2) {$\color{darkred}\scriptstyle{r}$};
\end{tikzpicture}
u^{-r-1}
       = \frac{t_{V,i}(u+i)}{u^{\eps_i(V)}}
$,
$\qquad
    \sum\limits_{r \in \mathbb{Z}}
\begin{tikzpicture}[baseline = 0]
  \draw[<-,thick,darkred] (0,0.4) to[out=180,in=90] (-.2,0.2);
  \draw[-,thick,darkred] (0.2,0.2) to[out=90,in=0] (0,.4);
 \draw[-,thick,darkred] (-.2,0.2) to[out=-90,in=180] (0,0);
  \draw[-,thick,darkred] (0,0) to[out=0,in=-90] (0.2,0.2);
 \node at (0,-.13) {$\scriptstyle{i}$};
   \node at (0.3,0.2) {$\scriptstyle{\lambda}$};
   \node at (-0.2,0.2) {$\color{darkred}\scriptstyle\bullet$};
   \node at (-0.4,0.2) {$\color{darkred}\scriptstyle{r}$};
\end{tikzpicture}
u^{-r-1}=
 \frac{ u^{\eps_i(V)}}{t_{V,i}(u+i)}$, where  $t_{V,i}(u)$
is the monic Hensel factor from
Subsection~\ref{subsec:polynomial-representatives}.
Both series are expanded at $u=\infty$ and use the auxiliary
oriented normalization. For negative dot numbers, these
coefficient identities specify the fake bubbles. For nonnegative
dot numbers, the bubbles are composites of the cup, cap, and dot
defined above. Their agreement with the displayed series follows
from the nonexceptional cases of \cite[(5.2) and~(5.4)]{GRS};
see also Lemma~\ref{km2c}.
\end{itemize}
\end{Defn}
The denominators in part (1) have residues $1$ when $j=i$ and
$1-j+i$ when $j\ne i,i+1$, so they are invertible.

\subsection{The $2$-morphisms associated with the type $0$ vertex}

\begin{Defn}\label{GammaLambda}
Let $V$ be a nonzero finitely generated object of
$\mathcal M_\lambda$.  Define
\begin{itemize}
\item[(1)] $
\Gamma_\lambda(u) := \prod_{i\in I_1}
\left[ \frac{1}{(u+i)(u-i)}\right]^{\sum\limits_{s\ge i}d_s(\lambda)}$.
\item [(2)] $
\Gamma_V(u) := \Gamma_\lambda(u) \cdot
\prod_{s\ge2}\left(\frac{\gamma_{V;s,1}(u)}{(u^2-1)^{d_s(\lambda)}}\right)^{s-1}
\prod_{s\ge i\ge 2}\frac{(u^2-i^2)^{d_s(\lambda)}}{\gamma_{V;s,i}(u)}
$, where, for $s\ge i\ge1$,
\[
\gamma_{V;s,i}(u):=
\frac{(-1)^{d_s(\lambda)}
t_{V,-s}(-u-s+i)t_{V,-s}(u-s+i)}
{(-u^2+i^2)^{\eps_{-s}(V)}}.
\]
\end{itemize}
\end{Defn}

\begin{Lemma}
\label{lem:zero-gamma-unit}
In the integral characteristic-zero setting,
$
\Gamma_V(u)\in Z_{\lambda,V}[[u]]^\times$ and $ \Gamma_V(-u)=\Gamma_V(u)$.
Consequently, both $\Gamma_V(x)$ and $\Gamma_V(x)^{-1}$
are defined on the relevant generalized $0$-eigenspaces.
\end{Lemma}
\begin{proof}
Only finitely many factors in Definition~\ref{GammaLambda} are nontrivial:
the $d_i(\lambda)$ vanish for $i\gg0$, and $V$ has finite dot
spectrum.
The factors $u^2-i^2$, for $i\ge1$, have unit constant terms.
The reduction of $t_{V,-s}(u)$ is a power of $u+s$, so its
evaluation at $-s+i$, for $i\ge1$, has nonzero residue.
Thus the numerator factors occurring in $\gamma_{V;s,i}(0)$
are also units in the local coefficient ring.
Every factor retained in the finite product is therefore
regular and invertible at $u=0$. This implies the first assertion.
The paired factors $t_{V,-s}(-u-s+i)t_{V,-s}(u-s+i)$
are even, which proves the second assertion.
\end{proof}

\begin{Defn}\label{notation-1}
Let $V$ be a nonzero finitely generated object of $\cM_\lambda$ and
put $c=c_\lambda$.  Define
\begin{itemize}
\item[(1)] $\eps=\max\bigl\{\eps_i(V),\eps_i(E_0V)\mid i\in I\bigr\}$,
\item[(2)] $S=\{i\in I_1\mid
\eps_i(V)+\eps_{-i}(V)+\eps_i(E_0V)+\eps_{-i}(E_0V)>0\}$.
\item[(3)] For $\eta\in\{+,-\}$, choose
$\varphi_V^\eta(u),h_V^\eta(u)\in Z_{\lambda,V}[u]$ such that
\[
f_V^\eta(u)n^\eta_{V,0}(u) +u^kh_V^\eta(u)=1,
\]
where $k\gg0$, and
$f_V^\eta(u):=\varphi_V^\eta(u)
\prod_{i\in S}(u-i)^{2\eps}(u+i)^{2\eps}$.

\item[(4)] $A_V(x):=m^+_{V,0}(x)f_V^+(x)\Gamma_V(x)$ and
$B_V(x):=\frac{(-1)^{c+1}}{2}
m^-_{V,0}(x)f_V^-(x)\Gamma_V(x)^{-1}$.

\item[(5)] $M_V(-x):=(1-2x)A_V(x)$ and
$W_V(-x):=(1+2x)B_V(x)$.

\item[(6)]
$r_V(x_1,x_2):=A_V\left(\frac{x_2-x_1}{2}\right)
B_V\left(\frac{x_1-x_2}{2}\right)$,

\item [(7)] $p_V(x_1,x_2):=\frac{A_V(x_2)-A_V(-x_1)}{x_1+x_2}$,

\item [(8)] $q_V(x_1,x_2):=\frac{B_V(x_1)-B_V(-x_2)}{x_1+x_2}$.
\end{itemize}
\end{Defn}

The set $S$ is finite. Let
$P_S(u)=\prod_{i\in S}(u-i)^{2\eps}(u+i)^{2\eps}$.
Both $P_S(0)$ and $n_{V,0}^\eta(0)$ are units in $Z_{\lambda,V}$.
Lemma~\ref{key-coprim}, applied to $P_Sn_{V,0}^\eta$ and $u^k$,
therefore gives the B\'ezout identity in part (3).
The rational functions $\Gamma_V^{\pm1}$ are expanded at zero.
For polynomial representatives of $A_V,B_V$, we retain the factors
$m^\pm_{V,0}f_V^\pm$ and truncate only $\Gamma_V^{\pm1}$ modulo
$u^k$. Thus the prescribed vanishing factors at nonzero spectral
values are retained. All divided differences are formed before
substituting dots. The truncation order $k$ is chosen as in the
next lemma.

Let $
\Delta_f(x_1,x_2):=\frac{f(x_1)-f(x_2)}{x_1-x_2}$. Then
$q_V(x_1,x_2)=\Delta_{B_V}(x_1,-x_2)$ and
$p_V(x_1,x_2)=\Delta_{A_V}(x_2,-x_1)$.

\begin{Lemma}\label{lem:finite-dd-evaluation}
Let $U\in\mathcal M$, let $A$ be a commutative unital subalgebra of
$\operatorname{End}_{\mathcal M}(U)$, and let $z_0,\ldots,z_s$ be
commuting nilpotent endomorphisms of $U$ that commute with $A$.
Set $B=A[z_0,\ldots,z_s]$ and $J=(z_0,\ldots,z_s)\subseteq B$.
Let $F(u)=P(u)Q(u)^{-1}$, where $P,Q\in A[u]$ and $Q(0)\in A^\times$.
Choose $n_i\ge1$ with $z_i^{n_i}=0$, and put
$N=1+\sum_{i=0}^s(n_i-1)$. Then $J^N=0$.
For every $q\ge0$ and $w_0,\ldots,w_q\in\{z_0,\ldots,z_s\}$,
with repetitions allowed, the divided difference $F[w_0,\ldots,w_q]$
is well-defined and depends only on $F$ modulo $u^{N+q}A[[u]]$.
Here $F[w_0]=F(w_0)$.
\end{Lemma}
\begin{proof}
Every monomial in the $z_i$ of total degree at least $N$ vanishes,
so $J^N=0$. Write $F(u)=\sum_{m\ge0}a_m u^m$ with $a_m\in A$.
The $q$-th divided difference of $u^m$ is zero for $m<q$ and is
the sum of all monomials of total degree $m-q$ in its $q+1$
variables otherwise. Hence
\[
F[w_0,\ldots,w_q]
=\sum_{m=q}^{N+q-1}a_m
\sum_{\substack{\alpha_0,\ldots,\alpha_q\ge0\\
                 \alpha_0+\cdots+\alpha_q=m-q}}
w_0^{\alpha_0}\cdots w_q^{\alpha_q},
\]
since terms with $m\ge N+q$ vanish. This gives the claimed bound.
For $q=0$, it agrees with $P(w_0)Q(w_0)^{-1}$ because
$Q(w_0)-Q(0)\in J$ and $Q(0)$ is a unit.
For a finite collection of divided differences of orders at most
$q_{\max}$, truncation modulo $u^{N+q_{\max}}$ suffices.
\end{proof}

On $E_0^rV$, the transported coefficients commute with the dots
by naturality, and dots on distinct strands commute by the
interchange law. These dots are nilpotent, so the lemma also
applies to their linear combinations, including $(x_2-x_1)/2$.

Set $\widehat A_V(u)=d_V^+(u)\Gamma_V(u)$ and
$\widehat B_V(u)=\frac{(-1)^{c+1}}2d_V^-(u)\Gamma_V(u)^{-1}$.
The polynomial representatives satisfy
$A_V-\widehat A_V,\ B_V-\widehat B_V\in u^kZ_{\lambda,V}[[u]]$.
By Lemma~\ref{lem:finite-dd-evaluation}, their evaluations and
divided differences at nilpotent dots depend only on
$\widehat A_V,\widehat B_V$ once $k$ satisfies the stated bound.
Identities involving these representatives are understood after
such evaluation unless stated as formal-series identities.

\begin{Lemma}\label{lem:MW-divided-difference}
Let $V$ be a nonzero finitely generated object of $\cM_\lambda$.
Under an evaluation as in Lemma~\ref{lem:finite-dd-evaluation},
let $x$ be nilpotent with $x^k=0$, where $k$ is the truncation
order in Definition~\ref{notation-1}. Then
$A_V(x)W_V(x)=1=B_V(x)M_V(x)$.
\end{Lemma}

\begin{proof}
Set $\kappa:=(-1)^{c+1}/2$.
After nilpotent evaluation, the B\'ezout identities give
$f_V^\pm(x)n_{V,0}^\pm(x)=1$. Hence
$
A_V(x)=d_V^+(x)\Gamma_V(x)$ and $
B_V(x)=\kappa\,d_V^-(x)\Gamma_V(x)^{-1}$.
Definition~\ref{ration111}(5),(6) and Lemma~\ref{lem:fg-zero}(2) give
$
c_V^+(x)c_V^-(-x)=\kappa(1-2x)$ and
$
\kappa\,d_V^+(x)d_V^-(-x)=(1-2x)^{-1}$.
Using $\Gamma_V(-x)=\Gamma_V(x)$, we obtain
$
A_V(x)B_V(-x)=(1-2x)^{-1}$ and $
B_V(x)A_V(-x)=(1+2x)^{-1}$.
Finally, by definition,
$
W_V(x)=(1-2x)B_V(-x)$ and $
M_V(x)=(1+2x)A_V(-x)$.
Thus $A_V(x)W_V(x)=1=B_V(x)M_V(x)$.
\end{proof}

Lemma~\ref{lem:MW-divided-difference} gives $A_V(0)B_V(0)=1$.
With $d=(x_2-x_1)/2$, define the divided difference
\begin{equation}\label{eq:gV-definition}
g_V(x_1,x_2):=\frac{1-A_V(d)B_V(d)}{x_1-x_2}
=\frac12\Delta_{A_VB_V}(d,0).
\end{equation}
The expression on the right explains the evaluation even when
$x_1-x_2$ is nilpotent.
For $\sigma\in\{+,-\}$, write $A_V^\sigma:=A_V(\sigma x)$
and $B_V^\sigma:=B_V(\sigma x)$.

\begin{Defn}\label{eq:red-cup}
Let $V$ be a nonzero finitely generated object of $\cM_\lambda$, and
write $d=(x_2-x_1)/2$ in the crossing formula. Define
\begin{multicols}{2}
\begin{enumerate}[label=(\arabic*),ref=\arabic*,leftmargin=*,itemsep=.5ex,topsep=0pt]
\item $
  \begin{tikzpicture}[baseline=-0.5ex ,scale=0.7, ]
        \draw[thick, darkred, -]
            (-0.4,0.5) .. controls (-0.4,-0.3) and (0.4,-0.3) .. (0.4,0.5);
            \node[darkred] at (-.4, .7) {$\scriptstyle{0}$};
            \node[darkred] at (.4, .7) {$\scriptstyle{0}$};
              \draw[-,darkg,thick] (.6,-0.5) to (.6,.5);
\node at (0.6,-.7) {$\darkg\scriptstyle{V}$};
\end{tikzpicture}
:=
\begin{tikzpicture}[baseline=-0.5ex,scale=0.7,]
        \draw[thick, black, ]
            (-0.4,0.5) .. controls (-0.4,-0.3) and (0.4,-0.3) .. (0.4,0.5);
             \node at (-.4, .7) {$\scriptstyle{0}$};
             \node at (.4, .7) {$\scriptstyle 0$};   \node at (-.38, .3) {$\scriptstyle\bullet$};
             \node at (-1.5, .4) {$\scriptstyle \frac{c^+_{E_{0}V}(x)}{\Gamma_{E_0 V}(x)}$};
\end{tikzpicture}$,

\item $\begin{tikzpicture}[baseline=-0.5ex,scale=0.7, , yshift=10pt]
        \draw[thick, darkred, -]
            (-0.4,-0.5) .. controls (-0.4,0.3) and (0.4,0.3) .. (0.4,-0.5);
            \node[darkred] at (-0.4, -0.7) {$\scriptstyle 0$};
        \node[darkred] at (0.4, -0.7) {$\scriptstyle 0$};
\draw[-,darkg,thick] (.7,-0.5) to (.7,.5);
\node at (0.7,-.7) {$\darkg\scriptstyle{V}$};
\end{tikzpicture}
 :=\begin{tikzpicture}[baseline=-0.5ex,scale=0.7, , yshift=10pt]
        \draw[thick]
            (-0.4,-0.5) .. controls (-0.4,0.3) and (0.4,0.3) .. (0.4,-0.5);
            \node at (-0.4, -0.7) {$\scriptstyle 0$};
  \node at (0.4, -0.73) {$\scriptstyle 0$};
  \node at (.4, -.35) {$\scriptstyle\bullet$};
  \node at (1.2, -.4) {$\scriptstyle A_V(x)$};
\end{tikzpicture}$,

\item $
\begin{tikzpicture}[baseline = 0, darkred]
	\draw[-,thick] (0.28,-.3) to (-0.28,.4);
	\draw[-,thick] (-0.28,-.3) to (0.28,.4);
   \node at (-.3,.55) {$\scriptstyle{0}$};
   \node at (.3,.55) {$\scriptstyle{0}$};
   \node at (-.3,-.55) {$\scriptstyle{0}$};
   \node at (.3,-.55) {$\scriptstyle{0}$};
   \node at (0.6,-.7) {$\darkg\scriptstyle{V}$};
   \draw[-,darkg,thick] (.6,-0.5) to (.6,.5);
\end{tikzpicture}
 :=\begin{tikzpicture}[baseline = 0, scale=0.8, ]
	\draw[-] (0.38,-.4) to (-0.38,.5);
	\draw[-] (-0.38,-.4) to (0.38,.5);
	\draw[-,densely dotted] (-0.19,-.2) to (0.19,-.2);
    \draw[-,densely dotted] (-0.2,.3) to (0.2,.3);
   \node at (-.4,-.55) {$\scriptstyle{0}$};
   \node at (.4,-.55) {$\scriptstyle{0}$};
   \node at (-.4,.65) {$\scriptstyle{0}$};
   \node at (.4,.65) {$\scriptstyle{0}$};
      \node at (1.0,-0.2) {$\scriptstyle A_V(d)$};
     \node at (1.0, 0.4) {$\scriptstyle B_V(d)$};
\end{tikzpicture}
 +\begin{tikzpicture}[baseline = 0, scale=0.8, ]
	\draw[-] (0.3,-.4) to (0.3,.5);
	\draw[-] (-0.3,-.4) to (-0.3,.5);
	\draw[-,densely dotted] (-0.29, .1) to (0.29, .1);
   \node at (-.3,-.55) {$\scriptstyle{0}$};
   \node at (.3,-.55) {$\scriptstyle{0}$};
   \node at (-.3,.65) {$\scriptstyle{0}$};
   \node at (.3,.65) {$\scriptstyle{0}$};
   \node at (1.0, 0.1) {$\scriptstyle{g_V(x_1,x_2)}$};
\end{tikzpicture}$.
\item $\begin{tikzpicture}[baseline=-0.5ex ,scale=0.7 ]
        \draw[thick, darkred, -] (0,-0.5) -- (0,0.5);
        \fill[darkred]  (0, 0) circle (2pt);
        \node[darkred] at (0, -0.7) {$\scriptstyle 0$};
        \draw[-,darkg,thick] (.4,-0.5) to (.4,.5);
\node at (.4,-.7) {$\darkg\scriptstyle{V}$};
\end{tikzpicture}:=\begin{tikzpicture}[baseline=-0.5ex,scale=0.7 ] \draw[thick, black, ] (0,-0.5) -- (0,0.5);
        \node at (0, 0) {$\scriptstyle\bullet$};
        \node at (0.6,0) {$\scriptstyle{x}$};
        \node at (0, -0.7) {$\scriptstyle 0$};
\end{tikzpicture}$
\end{enumerate}
\end{multicols}
\end{Defn}

We now check that these components are natural. The argument also
applies to the type $A$ coefficients above; compare
\cite[Lemma~4.30]{GRS}. Let $f:V\to W$ be a morphism between
finitely generated objects in the same central block, and put
$T=V\oplus W$. Since $m_T$ is a common annihilator of the two dot
endomorphisms, write
\[
 m_T=m_Vh,\qquad
 r_i:=\eps_i(T)-\eps_i(V)\geq0,\qquad
 h^\pm(u):=\prod_{\pm i>0}(u-i)^{r_i}.
\]
Restriction of bubble endomorphisms gives $n_T|_V=n_Vh$.
Uniqueness of monic Hensel lifting then gives
\[
 t_{T,i}|_V=t_{V,i}(u-i)^{r_i}\quad(i\ne0),\qquad
 n_{T,0}^\pm|_V=n_{V,0}^\pm h^\pm.
\]
The same factors $h^\pm$ occur in $m_{T,0}^\pm$. Consequently,
$c_T^\pm|_V=c_V^\pm$ and $d_T^\pm|_V=d_V^\pm$.
The complementary factor $h/(u-i)^{r_i}$ likewise cancels from
$a_{T,i}|_V$ and $b_{T,i}|_V$.
In $\gamma_{T;s,i}$, the extra factors coming from
$(u+s)^{r_{-s}}$ multiply to $(-u^2+i^2)^{r_{-s}}$ and cancel
against the denominator. Thus $\Gamma_T|_V=\Gamma_V$.
The same comparisons hold on $W$, and after applying $E_j$
to the direct sum $T$ for any $j\in I$.

All coefficients on $T$ commute with the endomorphism of $T$
whose only nonzero matrix entry is $f:V\to W$, since they are
built from central bubble endomorphisms. Their restrictions
therefore intertwine $f$. The same conclusion holds for the
normalized series $\widehat A,\widehat B$, and hence for
$A,B,p,q,g$ after sufficiently high nilpotent evaluation by
Lemma~\ref{lem:finite-dd-evaluation}. The affine Brauer generators
and the generalized eigenspace projections are natural, so the
composites in the definitions are natural as well. Morphisms
between different central blocks vanish. Finally, preservation of
the directed unions of finitely generated subobjects extends these
transformations to all objects of $\mathcal M$.

\begin{Lemma}\label{gamma-move1}
Let $V$ be a nonzero finitely generated object of $\mathcal M_\lambda$.
Then $E_0V\in\mathcal M_\lambda$, and, when $E_0V\ne0$, we have
$\begin{tikzpicture}[baseline=-0.5ex,scale=0.7,]
        \draw[thick, black, ]
            (-0.4,0.5) .. controls (-0.4,-0.3) and (0.4,-0.3) .. (0.4,0.5);
             \node at (-.4, .7) {$\scriptstyle{0}$};
             \node at (.4, .7) {$\scriptstyle 0$};
             \node at (-.38, .3) {$\scriptstyle\bullet$};
             \node at (-1.2, .4) {$\scriptstyle \Gamma^{-1}_{ E_0V}(x)$};
\end{tikzpicture}  = \begin{tikzpicture}[baseline=-0.5ex,scale=0.7,]
        \draw[thick, black, ]
            (-0.4,0.5) .. controls (-0.4,-0.3) and (0.4,-0.3) .. (0.4,0.5);
             \node at (-.4, .7) {$\scriptstyle{0}$};
             \node at (.4, .7) {$\scriptstyle 0$};
             \node at (.38, .3) {$\scriptstyle\bullet$};
             \node at (1.5, .3) {$\scriptstyle \Gamma^{-1}_V(x)$};
\end{tikzpicture}$.
\end{Lemma}

\begin{proof}
Since $\bar\alpha_0=0$, we have $E_0V\in\cM_\lambda$.
If $E_0V=0$, the assertion is immediate, so assume
$E_0V\ne0$. In the comparisons below, $E_0$ is applied
coefficientwise to rational functions.

Put $W=E_0V$ and let $x$ be its nilpotent dot. Work in the
finite-dimensional commutative algebra generated by $x$,
$E_0(Z_{\lambda,V})$, and $Z_{\lambda,W}$.
It is local: the elements $x$, $E_0(J_{\lambda,V})$, and
$J_{\lambda,W}$ generate a nilpotent ideal, and the quotient by
this ideal is $\Bbbk$. Thus monic Hensel lifting applies in this
common coefficient algebra. Set
\[
 N(u,x)=((u+x)^2-1)(u-x)^2,\qquad
 D(u,x)=((u-x)^2-1)(u+x)^2.
\]
Lemma~\ref{ou}(2), together with $n_X=O_Xm_X$, gives the polynomial
identity
\[
 n_W(u)E_0(m_V(u))D(u,x)
 =E_0(n_V(u))m_W(u)N(u,x).
\]
Modulo the maximal ideal, both $N$ and $D$ become $u^2(u^2-1)$.
Hence, for $s\ge2$, neither contributes a nontrivial $(u+s)$-primary
factor. Comparing the monic Hensel factors gives
\[
 \frac{t_{W,-s}(u)}{(u+s)^{\eps_{-s}(W)}}
 =E_0\!\left(\frac{t_{V,-s}(u)}{(u+s)^{\eps_{-s}(V)}}\right).
\]
Consequently,
$\gamma_{E_0V;s,i}(u)=E_0(\gamma_{V;s,i}(u))$
for $s\ge2$ and $1\le i\le s$.
Both sides are even functions of $u$ by Definition~\ref{GammaLambda}.
Moreover, we can rewrite $\Gamma_V(u)$ as
\[
 \Gamma_V(u)
 =
 \frac{1}{(u^2-1)^{d_1(\lambda)}}
 \prod_{s\ge2}
 \left(
   \frac{\gamma_{V;s,1}(u)}
        {(u^2-1)^{d_s(\lambda)}}
 \right)^s
 \prod_{\substack{s\ge i\ge1\\s\ge2}}
 \frac{1}{\gamma_{V;s,i}(u)}.
\]
The same formula holds for $\Gamma_{E_0V}(u)$, with
the same weight $\lambda$. It follows that
$\Gamma_{E_0V}(-u)=\Gamma_{E_0V}(u)=E_0(\Gamma_V(u))$.
Sliding the dot decoration through the cup replaces
$x$ by $-x$, so this equality gives the required identity.
\end{proof}

\begin{Lemma} \label{eq:red-cup-E0V}
When evaluated on a finitely generated object $V\in\mathcal M_\lambda$,
the following identities hold.
\begin{multicols}{2}
\begin{enumerate}[label=(\arabic*),ref=\arabic*,leftmargin=*,itemsep=.5ex,topsep=0pt]
\item
     $\begin{tikzpicture}[baseline=-0.5ex ,scale=0.7, ]
        \draw[thick, darkred, -]
            (-0.4,0.5) .. controls (-0.4,-0.3) and (0.4,-0.3) .. (0.4,0.5);
            \node[darkred] at (-.4, .7) {$\scriptstyle{0}$};
            \node[darkred] at (.4, .7) {$\scriptstyle{0}$};
              \draw[-,darkg,thick] (.6,-0.5) to (.6,.5);
\node at (0.6,-.7) {$\darkg\scriptstyle{V}$};
\end{tikzpicture}
 =\frac{1}{2}(-1)^{c+1} \begin{tikzpicture}[baseline=-0.5ex,scale=0.7,]
        \draw[thick, black, ]
            (-0.4,0.5) .. controls (-0.4,-0.3) and (0.4,-0.3) .. (0.4,0.5);
             \node at (-.4, .7) {$\scriptstyle{0}$};
             \node at (.4, .7) {$\scriptstyle 0$};
             \node at (.38, .3) {$\scriptstyle\bullet$};
             \node at (1.2, .4) {$\scriptstyle \frac{d^-_{V}(x)}{\Gamma_V(x)}$};
\end{tikzpicture}$,

\item
$\begin{tikzpicture}[baseline=-0.5ex,scale=0.7, , yshift=10pt]
        \draw[thick, darkred, -]
            (-0.4,-0.5) .. controls (-0.4,0.3) and (0.4,0.3) .. (0.4,-0.5);
        \node[darkred] at (-0.4, -0.7) {$\scriptstyle 0$};
        \node[darkred] at (0.4, -0.7) {$\scriptstyle 0$};
\draw[-,darkg,thick] (.7,-1.0) to (.7, 0);
\node at (0.7,-1.2) {$\darkg\scriptstyle{V}$};
\end{tikzpicture}
 =
 {2}(-1)^{c+1}\begin{tikzpicture}[baseline=-0.5ex,scale=0.7, , yshift=10pt]
        \draw[thick]
            (-0.4,-0.5) .. controls (-0.4,0.3) and (0.4,0.3) .. (0.4,-0.5);
            \node at (-0.4, -0.7) {$\scriptstyle 0$};
  \node at (0.4, -0.73) {$\scriptstyle 0$};
  \node at (-.4, -.35) {$\scriptstyle\bullet$};
  \node at (-2, -.4) {$\scriptstyle \Gamma_{E_0 V}(x)c^-_{E_0V}(x)$};
\end{tikzpicture}$.
\end{enumerate}
\end{multicols}
\end{Lemma}
\begin{proof}
Assume $E_0V\ne0$, put $W=E_0V$, and let $x$ be its dot.
In the common local coefficient algebra used in the proof of
Lemma~\ref{gamma-move1}, the same Hensel comparison gives
\[
 c_W^+(u)=E_0(c_V^+(u))\frac{u+x-1}{u-x-1},\qquad
 c_W^-(u)=E_0(c_V^-(u))\frac{u+x+1}{u-x+1}.
\]
Indeed, the positive nonzero factors of $N$ and $D$ in that proof
are $u+x-1$ and $u-x-1$, and the negative nonzero factors are
$u+x+1$ and $u-x+1$. Their remaining factors belong to the zero
primary component. The possible root $\frac12$ or $-\frac12$
from \eqref{barnv} stays in its prescribed factor because $V$ and
$W$ lie in the same central block.

Slide $n^+_{E_0V,0}(x)/m^+_{E_0V,0}(x)$ from the left leg through the cup
to the right; at the bend $x_1$ becomes $-x_2$, so $x_1+x_2=0$:
$$\begin{aligned}
\begin{tikzpicture}[baseline=-0.5ex,scale=0.7,]
        \draw[thick, black, ]
            (-0.4,0.5) .. controls (-0.4,-0.3) and (0.4,-0.3) .. (0.4,0.5);
             \node at (-.4, .7) {$\scriptstyle{0}$};
             \node at (.4, .7) {$\scriptstyle 0$};
             \node at (-.38, .3) {$\scriptstyle\bullet$};
             \node at (-2.2, .4) {$\scriptstyle \frac{n_{E_0V,0}^+(x)}
             {m_{E_0V,0}^+(x)} $};
\end{tikzpicture}
&=\left.
   \frac{n_{V,0}^+(x_1)}{m_{V,0}^+(x_1)}
   \frac{x_1+x_2-1}{x_1-x_2-1}
   \right|_{x_1+x_2=0} \begin{tikzpicture}[baseline=-0.5ex,scale=0.7,]
        \draw[thick, black, ]
            (-0.4,0.5) .. controls (-0.4,-0.3) and (0.4,-0.3) .. (0.4,0.5);
             \node at (-.4, .7) {$\scriptstyle{0}$};
             \node at (.4, .7) {$\scriptstyle 0$};
            \draw[-,densely dotted] (-.39, .3) to (.39,.3);
\end{tikzpicture}
\notag\\
&=\begin{tikzpicture}[baseline=-0.5ex,scale=0.7,]
        \draw[thick, black, ]
            (-0.4,0.5) .. controls (-0.4,-0.3) and (0.4,-0.3) .. (0.4,0.5);
             \node at (-.4, .7) {$\scriptstyle{0}$};
             \node at (.4, .7) {$\scriptstyle 0$};
             \node at (.38, .3) {$\scriptstyle\bullet$};
\end{tikzpicture}
             \frac{n_{V,0}^+(-x)}{m_{V,0}^+(-x)}\frac1{2x+1}
 =  \begin{tikzpicture}[baseline=-0.5ex,scale=0.7,]
        \draw[thick, black, ]
            (-0.4,0.5) .. controls (-0.4,-0.3) and (0.4,-0.3) .. (0.4,0.5);
             \node at (-.4, .7) {$\scriptstyle{0}$};
             \node at (.4, .7) {$\scriptstyle 0$};
             \node at (.38, .3) {$\scriptstyle\bullet$};
             \node at (2.1, .3) {$\scriptstyle (-1)^{c+1}\frac12\frac{m_{V,0}^-(x)}{n_{V,0}^-(x)}$};
\end{tikzpicture}.
\end{aligned}$$
The last equality uses Lemma~\ref{lem:fg-zero}(2).
Sliding the $\Gamma^{-1}$-factor through the cup by
Lemma~\ref{gamma-move1} proves (1).

For (2), use Lemma~\ref{lem:fg-zero}(2) with $u$ replaced
by $-u$, and invert the resulting unit identity. Apply the cap version
of the same primary-factor slide and transport the $\Gamma$-factor by
Lemma~\ref{gamma-move1}. The coefficient is $2(-1)^{c+1}$.
\end{proof}

\begin{Lemma}\label{lem:anothdefofcro}
On a nonzero finitely generated object $V\in\mathcal M_\lambda$, the
following two expressions for the zero crossing agree.
In each cup--cap term with a vertical dotted connector, $x_1$
denotes the dot on the right leg of the lower cap, and $x_2$
denotes the dot on the right leg of the upper cup. For the other
terms, the variables are read from left to right as usual.
\begin{equation} \begin{aligned}
\begin{tikzpicture}[baseline = 0, darkred]
	\draw[-,thick] (0.28,-.4) to (-0.28,.4);
	\draw[-,thick] (-0.28,-.4) to (0.28,.4);
    \node at (-.3,.55) {$\scriptstyle{0}$};
    \node at (.3,.55) {$\scriptstyle{0}$};
    \node at (-.3,-.55) {$\scriptstyle{0}$};
    \node at (.3,-.55) {$\scriptstyle{0}$};
    \node at (0.6,-.55) {$\darkg\scriptstyle{V}$};
    \draw[-,darkg,thick] (.6,-0.4) to (.6,.5);
\end{tikzpicture}
&=
\begin{tikzpicture}[baseline = 0, scale = 0.8]
	\draw[-] (0.38,-.4) to (-0.38,.5);
	\draw[-] (-0.38,-.4) to (0.38,.5);
	\draw[-,densely dotted] (-0.19,-.2) to (0.19,-.2);
    \node at (-.4,-.55) {$\scriptstyle {0}$};
    \node at (.4,-.55) {$\scriptstyle {0}$};
    \node at (-.4,.65) {$\scriptstyle {0}$};
    \node at (.4,.65) {$\scriptstyle {0}$};
    \node at (1,-0.2) {$\scriptstyle{r_V(x_1,x_2)}$};
\end{tikzpicture}
+
\begin{tikzpicture}[baseline = 0, scale = 0.8]
	\draw[-] (0.3,-.4) to (0.3,0.5);
	\draw[-] (-0.3,-.4) to (-0.3,0.5);
	\draw[-,densely dotted] (-0.29,0.1) to (0.29,0.1);
    \node at (-.3,-.55) {$\scriptstyle {0}$};
    \node at (.3,-.55) {$\scriptstyle {0}$};
    \node at (-.3,.65) {$\scriptstyle {0}$};
    \node at (.3,.65) {$\scriptstyle {0}$};
    \node at (1.2,0.1) {$\scriptstyle{\frac{1-r_V(x_1,x_2)}{x_1-x_2}}$};
\end{tikzpicture}
+
\begin{tikzpicture}[baseline = 0, scale = 0.8]
    \draw[-] (0,0.5) to[out=-90,in=-90] (0.8,0.5);
    \draw[-] (0,-0.4) to[out=90,in=90] (0.8,-0.4);
    \node at (0,-0.55) {$\scriptstyle {0}$};
    \node at (0.8,-0.55) {$\scriptstyle {0}$};
    \node at (0,0.6) {$\scriptstyle {0}$};
    \node at (0.8,0.6) {$\scriptstyle {0}$};
    \draw[densely dotted,-] (0.6,-0.2) -- (0.6,0.28);
    \node at (2.2,0) {$\scriptstyle q_V(x_2,x_1)A_V(x_1)$};
\end{tikzpicture}
\\
&=
\begin{tikzpicture}[baseline = 0, scale = 0.8]
    \draw[-] (0.38,-.4) to (-0.38,.5);
	\draw[-] (-0.38,-.4) to (0.38,.5);
	\draw[-,densely dotted] (-0.19,0.3) to (0.19,0.3);
    \node at (-.4,-.55) {$\scriptstyle {0}$};
    \node at (.4,-.55) {$\scriptstyle {0}$};
    \node at (-.4,.65) {$\scriptstyle {0}$};
    \node at (.4,.65) {$\scriptstyle {0}$};
    \node at (1,0.25) {$\scriptstyle{r_V(x_2,x_1)}$};
\end{tikzpicture}
+
\begin{tikzpicture}[baseline = 0, scale = 0.8]
    \draw[-] (0.3,-.4) to (0.3,0.5);
	\draw[-] (-0.3,-.4) to (-0.3,0.5);
	\draw[-,densely dotted] (-0.29,0.1) to (0.29,0.1);
    \node at (-.3,-.55) {$\scriptstyle {0}$};
    \node at (.3,-.55) {$\scriptstyle {0}$};
    \node at (-.3,.65) {$\scriptstyle {0}$};
    \node at (.3,.65) {$\scriptstyle {0}$};
    \node at (1.2,0.1) {$\scriptstyle{ \frac{1-r_V(x_2,x_1)}{x_1-x_2}}$};
\end{tikzpicture}
+
\begin{tikzpicture}[baseline = 0, scale = 0.8]
    \draw[-] (0,0.5) to[out=-90,in=-90] (0.8,0.5);
    \draw[-] (0,-0.4) to[out=90,in=90] (0.8,-0.4);
    \node at (0,-0.55) {$\scriptstyle {0}$};
    \node at (0.8,-0.55) {$\scriptstyle {0}$};
    \node at (0,0.6) {$\scriptstyle {0}$};
    \node at (0.8,0.6) {$\scriptstyle {0}$};
    \draw[densely dotted,-] (0.6,-0.2) -- (0.6,0.28);
    \node at (2.2,0) {$\scriptstyle B_V(x_2)p_V(x_2,x_1)$};
\end{tikzpicture}
\end{aligned}
\end{equation}
\end{Lemma}

\begin{proof}
For $\varphi(u)\in Z_{\lambda,V}[u]$, the affine Brauer relations give
the calculation below. Coefficient extraction is at infinity.
In the last two cup--cap diagrams, the circled minus on a leg
denotes $(t-x)^{-1}$ on that leg, where $t=u/2$.
\begingroup
\small
\allowdisplaybreaks[4]
\begin{align*}
\begin{tikzpicture}[x=1bp,y=1bp,scale=.92,
  baseline={([yshift=-.5ex]current bounding box.center)},
  line cap=round,line join=round,
  every node/.style={inner sep=.8pt,outer sep=0pt,font=\fontsize{8}{9}\selectfont}]
\draw[black,line width=.65bp] (2.217,4.319) -- (19.792,27.166);
\draw[black,line width=.65bp] (19.792,4.319) -- (2.217,27.166);
\draw[black,densely dotted,line width=.55bp] (5.462,22.948) -- (16.547,22.948);
\node[anchor=west,font=\fontsize{7.5}{7.5}\selectfont] at (23.292,22.948) {$\textstyle \varphi(x_2-x_1)$};
\node[anchor=north] at (1.302,3.13) {$0$};
\node[anchor=north] at (20.707,3.13) {$0$};
\node[anchor=south] at (1.302,28.355) {$0$};
\node[anchor=south] at (20.707,28.355) {$0$};
\end{tikzpicture}
-
\begin{tikzpicture}[x=1bp,y=1bp,scale=.92,
  baseline={([yshift=-.5ex]current bounding box.center)},
  line cap=round,line join=round,
  every node/.style={inner sep=.8pt,outer sep=0pt,font=\fontsize{8}{9}\selectfont}]
\draw[black,line width=.65bp] (2.217,4.319) -- (19.792,27.166);
\draw[black,line width=.65bp] (19.792,4.319) -- (2.217,27.166);
\draw[black,densely dotted,line width=.55bp] (5.462,8.537) -- (16.547,8.537);
\node[anchor=west,font=\fontsize{7.5}{7.5}\selectfont] at (23.292,8.537) {$\textstyle \varphi(x_1-x_2)$};
\node[anchor=north] at (1.302,3.13) {$0$};
\node[anchor=north] at (20.707,3.13) {$0$};
\node[anchor=south] at (1.302,28.355) {$0$};
\node[anchor=south] at (20.707,28.355) {$0$};
\end{tikzpicture}
&=
\left[
\begin{tikzpicture}[x=1bp,y=1bp,scale=.92,
  baseline={([yshift=-.5ex]current bounding box.center)},
  line cap=round,line join=round,
  every node/.style={inner sep=.8pt,outer sep=0pt,font=\fontsize{8}{9}\selectfont}]
\draw[black,line width=.65bp] (2.217,4.319) -- (19.792,27.166);
\draw[black,line width=.65bp] (19.792,4.319) -- (2.217,27.166);
\draw[black,densely dotted,line width=.55bp] (5.462,22.948) -- (16.547,22.948);
\node[anchor=west,font=\fontsize{7.5}{7.5}\selectfont] at (23.292,22.948) {$\textstyle \frac{\varphi(u)}{u-(x_2-x_1)}$};
\node[anchor=north] at (1.302,3.13) {$0$};
\node[anchor=north] at (20.707,3.13) {$0$};
\node[anchor=south] at (1.302,28.355) {$0$};
\node[anchor=south] at (20.707,28.355) {$0$};
\end{tikzpicture}
-
\begin{tikzpicture}[x=1bp,y=1bp,scale=.92,
  baseline={([yshift=-.5ex]current bounding box.center)},
  line cap=round,line join=round,
  every node/.style={inner sep=.8pt,outer sep=0pt,font=\fontsize{8}{9}\selectfont}]
\draw[black,line width=.65bp] (2.217,4.319) -- (19.792,27.166);
\draw[black,line width=.65bp] (19.792,4.319) -- (2.217,27.166);
\draw[black,densely dotted,line width=.55bp] (5.462,8.537) -- (16.547,8.537);
\node[anchor=west,font=\fontsize{7.5}{7.5}\selectfont] at (23.292,8.537) {$\textstyle \frac{\varphi(u)}{u-(x_1-x_2)}$};
\node[anchor=north] at (1.302,3.13) {$0$};
\node[anchor=north] at (20.707,3.13) {$0$};
\node[anchor=south] at (1.302,28.355) {$0$};
\node[anchor=south] at (20.707,28.355) {$0$};
\end{tikzpicture}\right]_{u^{-1}}
\\
&=
\left[
\begin{tikzpicture}[x=1bp,y=1bp,scale=.92,
  baseline={([yshift=-.5ex]current bounding box.center)},
  line cap=round,line join=round,
  every node/.style={inner sep=.8pt,outer sep=0pt,font=\fontsize{8}{9}\selectfont}]
\draw[black,line width=.65bp] (2.217,4.319) -- (19.792,27.166);
\draw[black,line width=.65bp] (19.792,4.319) -- (2.217,27.166);
\draw[black,densely dotted,line width=.55bp] (5.462,22.948) -- (16.547,22.948);
\draw[black,densely dotted,line width=.55bp] (5.462,8.537) -- (16.547,8.537);
\node[anchor=west,font=\fontsize{7.5}{7.5}\selectfont] at (23.292,22.948) {$\textstyle \frac{\varphi(u)}{u-(x_2-x_1)}$};
\node[anchor=west,font=\fontsize{7.5}{7.5}\selectfont] at (23.292,3.948) {$\textstyle \frac{u-(x_1-x_2)}{u-(x_1-x_2)}$};
\node[anchor=north] at (1.302,3.13) {$0$};
\node[anchor=north] at (20.707,3.13) {$0$};
\node[anchor=south] at (1.302,28.355) {$0$};
\node[anchor=south] at (20.707,28.355) {$0$};
\end{tikzpicture}
-
\begin{tikzpicture}[x=1bp,y=1bp,scale=.92,
  baseline={([yshift=-.5ex]current bounding box.center)},
  line cap=round,line join=round,
  every node/.style={inner sep=.8pt,outer sep=0pt,font=\fontsize{8}{9}\selectfont}]
\draw[black,line width=.65bp] (2.217,4.319) -- (19.792,27.166);
\draw[black,line width=.65bp] (19.792,4.319) -- (2.217,27.166);
\draw[black,densely dotted,line width=.55bp] (5.462,22.948) -- (16.547,22.948);
\draw[black,densely dotted,line width=.55bp] (5.462,8.537) -- (16.547,8.537);
\node[anchor=west,font=\fontsize{7.5}{7.5}\selectfont] at (23.292,22.948) {$\textstyle \frac{\varphi(u)(u-(x_2-x_1))}{u-(x_2-x_1)}$};
\node[anchor=west,font=\fontsize{7.5}{7.5}\selectfont] at (23.292,3.948) {$\textstyle \frac1{u-(x_1-x_2)}$};
\node[anchor=north] at (1.302,3.13) {$0$};
\node[anchor=north] at (20.707,3.13) {$0$};
\node[anchor=south] at (1.302,28.355) {$0$};
\node[anchor=south] at (20.707,28.355) {$0$};
\end{tikzpicture}\right]_{u^{-1}}
\\
&=
\left[2
\begin{tikzpicture}[x=1bp,y=1bp,scale=.92,
  baseline={([yshift=-.5ex]current bounding box.center)},
  line cap=round,line join=round,
  every node/.style={inner sep=.8pt,outer sep=0pt,font=\fontsize{8}{9}\selectfont}]
\draw[black,line width=.65bp] (2.217,34.197) .. controls (2.217,16.798) and (19.792,16.798) .. (19.792,34.197);
\draw[black,line width=.65bp] (2.217,4.32) .. controls (2.217,21.719) and (19.792,21.719) .. (19.792,4.32);
\draw[black,densely dotted,line width=.55bp] (18.88,27.988) -- (3.12,27.988);
\draw[black,densely dotted,line width=.55bp] (18.88,10.529) -- (3.12,10.529);
\node[anchor=north] at (2.217,2.82) {$0$};
\node[anchor=north] at (19.792,2.82) {$0$};
\node[anchor=south] at (2.217,35.697) {$0$};
\node[anchor=south] at (19.792,35.697) {$0$};
\node[anchor=west,font=\fontsize{7.5}{7.5}\selectfont] at (23.292,27.988) {$\textstyle \frac{\varphi(u)}{u-(x_2-x_1)}$};
\node[anchor=west,font=\fontsize{7.5}{7.5}\selectfont] at (23.292,8.988) {$\textstyle \frac1{u-(x_1-x_2)}$};
\end{tikzpicture}
-2
\begin{tikzpicture}[x=1bp,y=1bp,scale=.92,
  baseline={([yshift=-.5ex]current bounding box.center)},
  line cap=round,line join=round,
  every node/.style={inner sep=.8pt,outer sep=0pt,font=\fontsize{8}{9}\selectfont}]
\draw[black,line width=.65bp] (2.217,4.319) -- (2.217,27.166);
\draw[black,line width=.65bp] (19.792,4.319) -- (19.792,27.166);
\draw[black,densely dotted,line width=.55bp] (2.217,15.742) -- (19.792,15.742);
\node[anchor=west,font=\fontsize{7.5}{7.5}\selectfont] at (23.292,15.742) {$\textstyle \frac{\varphi(u)}{(u-(x_2-x_1))(u-(x_1-x_2))}$};
\node[anchor=north] at (2.217,2.819) {$0$};
\node[anchor=north] at (19.792,2.819) {$0$};
\node[anchor=south] at (2.217,28.666) {$0$};
\node[anchor=south] at (19.792,28.666) {$0$};
\end{tikzpicture}\right]_{u^{-1}}
\\
&=\left[2
\begin{tikzpicture}[x=1bp,y=1bp,scale=.92,
  baseline={([yshift=-.5ex]current bounding box.center)},
  line cap=round,line join=round,
  every node/.style={inner sep=.8pt,outer sep=0pt,font=\fontsize{8}{9}\selectfont}]
\draw[black,line width=.65bp] (2.217,34.197) .. controls (2.217,16.798) and (19.792,16.798) .. (19.792,34.197);
\draw[black,line width=.65bp] (2.217,4.32) .. controls (2.217,21.719) and (19.792,21.719) .. (19.792,4.32);
\fill[black] (18.441,26.832) circle[radius=1.15bp];
\fill[black] (18.441,11.685) circle[radius=1.15bp];
\node[anchor=north] at (2.217,2.82) {$0$};
\node[anchor=north] at (19.792,2.82) {$0$};
\node[anchor=south] at (2.217,35.697) {$0$};
\node[anchor=south] at (19.792,35.697) {$0$};
\node[anchor=west,font=\fontsize{7.5}{7.5}\selectfont] at (23.292,26.832) {$\textstyle \frac{\varphi(u)}{u-2x}$};
\node[anchor=west,font=\fontsize{7.5}{7.5}\selectfont] at (23.292,7.832) {$\textstyle \frac1{u+2x}$};
\end{tikzpicture}
-2
\begin{tikzpicture}[x=1bp,y=1bp,scale=.92,
  baseline={([yshift=-.5ex]current bounding box.center)},
  line cap=round,line join=round,
  every node/.style={inner sep=.8pt,outer sep=0pt,font=\fontsize{8}{9}\selectfont}]
\draw[black,line width=.65bp] (2.217,4.319) -- (2.217,27.166);
\draw[black,line width=.65bp] (19.792,4.319) -- (19.792,27.166);
\draw[black,densely dotted,line width=.55bp] (2.217,15.742) -- (19.792,15.742);
\node[anchor=west,font=\fontsize{7.5}{7.5}\selectfont] at (23.292,15.742) {$\textstyle \frac{\varphi(u)}{u^2-(x_2-x_1)^2}$};
\node[anchor=north] at (2.217,2.819) {$0$};
\node[anchor=north] at (19.792,2.819) {$0$};
\node[anchor=south] at (2.217,28.666) {$0$};
\node[anchor=south] at (19.792,28.666) {$0$};
\end{tikzpicture}\right]_{u^{-1}}
\\
&=\left[
\begin{tikzpicture}[x=1bp,y=1bp,scale=.92,
  baseline={([yshift=-.5ex]current bounding box.center)},
  line cap=round,line join=round,
  every node/.style={inner sep=.8pt,outer sep=0pt,font=\fontsize{8}{9}\selectfont}]
\draw[black,line width=.65bp] (0,54.75) .. controls (0,36.135) and (17.52,36.135) .. (17.52,54.75);
\draw[black,line width=.65bp] (0,0) .. controls (0,18.615) and (17.52,18.615) .. (17.52,0);
\draw[black,fill=white,line width=.55bp] (16.877,49.062) circle[radius=2.825bp];
\draw[black,line width=.55bp] (15.323,49.062) -- (18.431,49.062);
\draw[black,fill=white,line width=.55bp] (0.643,5.688) circle[radius=2.825bp];
\draw[black,line width=.55bp] (-0.911,5.688) -- (2.197,5.688);
\node[anchor=south] at (0,56.064) {$0$};
\node[anchor=north] at (0,-1.314) {$0$};
\node[anchor=south] at (17.52,56.064) {$0$};
\node[anchor=north] at (17.52,-1.314) {$0$};
\end{tikzpicture}\varphi(2t)\right]_{t^{-1}}
-
\begin{tikzpicture}[x=1bp,y=1bp,scale=.92,
  baseline={([yshift=-.5ex]current bounding box.center)},
  line cap=round,line join=round,
  every node/.style={inner sep=.8pt,outer sep=0pt,font=\fontsize{8}{9}\selectfont}]
\draw[black,line width=.65bp] (2.217,4.319) -- (2.217,27.166);
\draw[black,line width=.65bp] (19.792,4.319) -- (19.792,27.166);
\draw[black,densely dotted,line width=.55bp] (2.217,15.742) -- (19.792,15.742);
\node[anchor=west,font=\fontsize{7.5}{7.5}\selectfont] at (23.292,15.742) {$\textstyle \sum_{i\ge0}[2\varphi(u)]_{u^{2i+1}}(x_2-x_1)^{2i}$};
\node[anchor=north] at (2.217,2.819) {$0$};
\node[anchor=north] at (19.792,2.819) {$0$};
\node[anchor=south] at (2.217,28.666) {$0$};
\node[anchor=south] at (19.792,28.666) {$0$};
\end{tikzpicture}
\\
&=\left[
\begin{tikzpicture}[x=1bp,y=1bp,scale=.92,
  baseline={([yshift=-.5ex]current bounding box.center)},
  line cap=round,line join=round,
  every node/.style={inner sep=.8pt,outer sep=0pt,font=\fontsize{8}{9}\selectfont}]
\draw[black,line width=.65bp] (0,54.75) .. controls (0,36.135) and (17.52,36.135) .. (17.52,54.75);
\draw[black,line width=.65bp] (0,0) .. controls (0,18.615) and (17.52,18.615) .. (17.52,0);
\draw[black,fill=white,line width=.55bp] (16.877,49.062) circle[radius=2.825bp];
\draw[black,line width=.55bp] (15.323,49.062) -- (18.431,49.062);
\draw[black,fill=white,line width=.55bp] (0.643,5.688) circle[radius=2.825bp];
\draw[black,line width=.55bp] (-0.911,5.688) -- (2.197,5.688);
\node[anchor=south] at (0,56.064) {$0$};
\node[anchor=north] at (0,-1.314) {$0$};
\node[anchor=south] at (17.52,56.064) {$0$};
\node[anchor=north] at (17.52,-1.314) {$0$};
\end{tikzpicture}\varphi(2t)\right]_{t^{-1}}
-
\begin{tikzpicture}[x=1bp,y=1bp,scale=.92,
  baseline={([yshift=-.5ex]current bounding box.center)},
  line cap=round,line join=round,
  every node/.style={inner sep=.8pt,outer sep=0pt,font=\fontsize{8}{9}\selectfont}]
\draw[black,line width=.65bp] (2.217,4.319) -- (2.217,27.166);
\draw[black,line width=.65bp] (19.792,4.319) -- (19.792,27.166);
\draw[black,densely dotted,line width=.55bp] (2.217,15.742) -- (19.792,15.742);
\node[anchor=west,font=\fontsize{7.5}{7.5}\selectfont] at (23.292,15.742) {$\textstyle \frac{\varphi(x_2-x_1)-\varphi(x_1-x_2)}{x_2-x_1}$};
\node[anchor=north] at (2.217,2.819) {$0$};
\node[anchor=north] at (19.792,2.819) {$0$};
\node[anchor=south] at (2.217,28.666) {$0$};
\node[anchor=south] at (19.792,28.666) {$0$};
\end{tikzpicture}
\end{align*}
\endgroup
The last equality collects the odd-degree coefficients of $\varphi$.
By Lemma~\ref{lem:finite-dd-evaluation}, the same calculation applies
to sufficiently high polynomial representatives of $A_V$ and $B_V$.
Moving $B_V$ downward and $A_V$ upward gives the two displayed
expressions. Finally, for nilpotent $x_1,x_2$ and $d=(x_2-x_1)/2$,
\begin{equation}\label{rx1x2}\begin{aligned}
r_V(x_1,x_2)&= \frac{(-1)^{c+1}}2  m^+_{V,0}(d)f_V^+(d)m^-_{V,0}(-d)f_V^-(-d)\\
&=\frac{(-1)^{c+1}}2\frac{m^+_{V,0}(d)m^-_{V,0}(-d)}{n^+_{V,0}(d)n^-_{V,0}(-d)}
=\frac{(-1)^{c+1}}2\frac1{(-1)^{c+1}(\frac12-\frac{x_2-x_1}2)}
=\frac1{1+x_1-x_2}.
\end{aligned}\end{equation}
Here we use the evenness of $\Gamma_V$, the B\'ezout identities,
and Lemma~\ref{lem:fg-zero}(2). Henceforth, write
$r(x_1,x_2):=(1+x_1-x_2)^{-1}$ for the evaluated coefficient
$r_V(x_1,x_2)$.
\end{proof}

\subsection{The $2$-morphisms in mixed cases}

In this subsection, we omit arrows on red strands and use signed
colors: a strand labeled $j$ represents $E_j$.

\begin{Defn}\label{natrual-basic-0i}
Let either $i=0$ and $j\in I_1$, or $i\in I_1$ and $j=0$.
On a nonzero finitely generated object
$V\in\mathcal M_\lambda$, define
\[
\begin{tikzpicture}[baseline = 0]
	\draw[-,darkred,thick] (0.28,-.3) to (-0.28,.4);
	\draw[-,darkred,thick] (-0.28,-.3) to (0.28,.4);
    \node at (-.3,-.45) {$\scriptstyle{j}$};
    \node at (.3,-.45) {$\scriptstyle{i}$};
    \node at (0.5,0.05) {$\scriptstyle{\lambda}$};
\end{tikzpicture}
=
\left\{
\begin{array}{ll}
\begin{tikzpicture}[baseline = 0, scale=0.8, ]
	\draw[-] (0.38,-.4) to (-0.38,.5);
	\draw[-] (-0.38,-.4) to (0.38,.5);
	\draw[-,densely dotted] (0.2,-.24) to (-0.19,-.24);
    \node at (0,0.05) {$\diamond$};
    \node at (-.4,-.55) {$\scriptstyle{i+1}$};
    \node at (-.4,.69) {$\scriptstyle{i}$};
    \node at (.4,.69) {$\scriptstyle{i+1}$};
    \node at (.4,-.55) {$\scriptstyle{i}$};
    \node at (0.8,-0.2) {$\scriptstyle{d_{+-}^{12}}$};
\end{tikzpicture}
= \begin{tikzpicture}[baseline = 0, scale=0.8, ]
	\draw[-] (0.38,-.4) to (-0.38,.5);
	\draw[-] (-0.38,-.4) to (0.38,.5);
	\draw[-,densely dotted] (0.2,.3) to (-0.19,.3);
    \node at (0,0.05) {$\diamond$};
    \node at (-.4,-.55) {$\scriptstyle{i+1}$};
    \node at (-.4,.69) {$\scriptstyle{i}$};
    \node at (.4,.69) {$\scriptstyle{i+1}$};
    \node at (.4,-.55) {$\scriptstyle{i}$};
    \node at (0.8,0.3) {$\scriptstyle{d_{-+}^{12}}$};
\end{tikzpicture}
&\text{if $j=i+1$,}
\\
\begin{tikzpicture}[baseline = 0, scale=0.8, ]
	\draw[-] (0.38,-.4) to (-0.38,.5);
    \node at (0,0.05) {$\diamond$};
	\draw[-] (-0.38,-.4) to (0.38,.5);
	\draw[-,densely dotted] (0.18,-.24) to (-0.19,-.24);
    \node at (-.4,-.62) {$\scriptstyle{j}$};
    \node at (.4,-.6) {$\scriptstyle{i}$};
    \node at (.4,.75) {$\scriptstyle{j}$};
    \node at (-.4,.73) {$\scriptstyle{i}$};
    \node at (1.0,-0.2) {$\scriptstyle{\frac{d^{12}_{+-}} {\ell_{-+}^{12}}}$};
\end{tikzpicture}
=\begin{tikzpicture}[baseline = 0, scale=0.8, ]
	\draw[-] (0.38,-.4) to (-0.38,.5);
    \node at (0,0.05) {$\diamond$};
	\draw[-] (-0.38,-.4) to (0.38,.5);
	\draw[-,densely dotted] (0.18,.3) to (-0.19,.3);
    \node at (-.4,-.62) {$\scriptstyle{j}$};
    \node at (.4,-.6) {$\scriptstyle{i}$};
    \node at (.4,.75) {$\scriptstyle{j}$};
    \node at (-.4,.73) {$\scriptstyle{i}$};
    \node at (0.8,0.25) {$\scriptstyle{\frac{d^{12}_{-+}} {\ell_{+-}^{12}}}$};
\end{tikzpicture}
&\text{otherwise.}
\end{array}
\right.
\]
\end{Defn}
The denominators are invertible by the same residue calculation
as in Definition~\ref{natrual-basic1}(1).

\begin{Defn}\label{zero-neg_cr}
For $j\in I_1$, define
\begin{multicols}{2}
\begin{enumerate}[label=(\arabic*),ref=\arabic*,leftmargin=*,itemsep=.5ex,topsep=0pt]
\item
$\begin{tikzpicture}[scale=0.45, x=1cm,y=1cm,line width=.65pt,line cap=round,line join=round,baseline={([yshift=-.5ex]current bounding box.center)}]
    \draw[darkred] (0,0)--(1.2,1.8);
    \draw[darkred] (1.2,0)--(0,1.8);
    \node[font=\scriptsize,text=darkred,anchor=north,inner sep=.7pt,outer sep=0pt] at (0,-0.05) {$0$};
    \node[font=\scriptsize,text=darkred,anchor=north,inner sep=.7pt,outer sep=0pt] at (1.2,-0.05) {$-j$};
\end{tikzpicture}
:=
\begin{tikzpicture}[scale=0.45,x=1cm,y=1cm,line width=.65pt,line cap=round,line join=round,baseline={([yshift=-.5ex]current bounding box.center)}]
    \draw[draw=darkred] (0,1.8) .. controls (0,-0.65) and (0.6,0.1) .. (1.05,0.95) .. controls (1.5,1.8) and (2.1,2.45) .. (2.1,0);
    \draw[draw=darkred] (0.65,1.8) .. controls (0.65,1.15) and (1.45,0.65) .. (1.45,0);
    \node[font=\scriptsize,text=darkred,anchor=north,inner sep=.7pt,outer sep=0pt] at (1.45,-0.05) {$0$};
    \node[font=\scriptsize,text=darkred,anchor=north,inner sep=.7pt,outer sep=0pt] at (2.1,-0.05) {$-j$};
\end{tikzpicture}$

\item
$\begin{tikzpicture}[scale=0.45,x=1cm,y=1cm,line width=.65pt,line cap=round,line join=round,baseline={([yshift=-.5ex]current bounding box.center)}]
    \draw[darkred] (0,0)--(1.2,1.8);
    \draw[darkred] (1.2,0)--(0,1.8);
    \node[font=\scriptsize,text=darkred,anchor=north,inner sep=.7pt,outer sep=0pt] at (0,-0.05) {$-j$};
    \node[font=\scriptsize,text=darkred,anchor=north,inner sep=.7pt,outer sep=0pt] at (1.2,-0.05) {$0$};
\end{tikzpicture}
:=
\begin{tikzpicture}[scale=0.45,x=1cm,y=1cm,line width=.65pt,line cap=round,line join=round,baseline={([yshift=-.5ex]current bounding box.center)}]
    \draw[draw=darkred] (0,0) .. controls (0,2.45) and (0.6,1.7) .. (1.05,0.85) .. controls (1.5,0) and (2.1,-0.65) .. (2.1,1.8);
    \draw[draw=darkred] (0.65,0) .. controls (0.65,0.65) and (1.45,1.15) .. (1.45,1.8);
    \node[font=\scriptsize,text=darkred,anchor=north,inner sep=.7pt,outer sep=0pt] at (0,-0.05) {$-j$};
    \node[font=\scriptsize,text=darkred,anchor=north,inner sep=.7pt,outer sep=0pt] at (0.65,-0.05) {$0$};
\end{tikzpicture}$
\end{enumerate}
\end{multicols}
\end{Defn}

\subsection{Comparison with the unoriented generators}
\label{subsec:forward-normalization}
We now translate the auxiliary oriented formulas into the
normalization of Definition~\ref{def-2ikac}. The ordinary
Kac--Moody weight has simple-coroot coordinates $-\lambda_i$ for
$i>0$. This conversion requires scalar factors as well as a change
of notation. For $i>0$, put
\[
 \kappa_i(\lambda):=\zeta_i\gamma_i(\lambda).
\]
Write $y_i^+$ and $y_i^-$ for the upward and downward dots in
Definition~\ref{natrual-basic}(1), and write
\[
\begin{aligned}
 \eta_i^R&:\mathrm{Id}\longrightarrow E_{-i}E_i,
 &\varepsilon_i^R&:E_iE_{-i}\longrightarrow\mathrm{Id},\\
 \eta_i^L&:\mathrm{Id}\longrightarrow E_iE_{-i},
 &\varepsilon_i^L&:E_{-i}E_i\longrightarrow\mathrm{Id}
\end{aligned}
\]
for the cups and caps in parts (2) and (3), respectively.
All four transformations here are evaluated on
$\mathcal M_\lambda$. The unoriented dots of signed colors
$i$ and $-i$ act by
\begin{equation}\label{eq:forward-dot-normalization}
 X_i=y_i^+=x-i,\qquad X_{-i}=-y_i^-=x+i.
\end{equation}
The corresponding cups and caps are
\begin{equation}\label{eq:forward-adjunction-normalization}
\begin{aligned}
 \widetilde\eta_i^R&=\eta_i^R,
 &\widetilde\varepsilon_i^R&=\varepsilon_i^R,\\
 \widetilde\eta_i^L&=\kappa_i(\lambda)^{-1}\eta_i^L,
 &\widetilde\varepsilon_i^L&=\kappa_i(\lambda)\varepsilon_i^L.
\end{aligned}
\end{equation}
For signed nonzero colors $a,b$, let
$T^{\mathrm{or}}_{a,b}:E_aE_b\Rightarrow E_bE_a$ denote the
auxiliary crossing in Definition~\ref{natrual-basic1}; the
subscripts record the bottom colors from left to right.
The unoriented crossings are
\begin{equation}\label{eq:forward-crossing-normalization}
\begin{aligned}
 S_{i,j}&=T^{\mathrm{or}}_{i,j},
 &S_{-i,-j}&=t_{ij}T^{\mathrm{or}}_{-i,-j},\\
 S_{i,-j}&=T^{\mathrm{or}}_{i,-j},
 &S_{-i,j}&=T^{\mathrm{or}}_{-i,j}
 \qquad(i,j>0).
\end{aligned}
\end{equation}

For completeness, write $\mathcal B_a(u)$ for the unoriented
bubble series whose left leg has signed color $a$, and put
$F_i(u)=t_{V,i}(u+i)u^{-\eps_i(V)}$. The conversion of the two
series in Definition~\ref{natrual-basic1}(4) is
\begin{equation}\label{eq:forward-bubble-normalization}
 \mathcal B_{-i}(u)=\kappa_i(\lambda)F_i(u),\qquad
 \mathcal B_i(u)=-\kappa_i(\lambda)^{-1}F_i(-u)^{-1}.
\end{equation}
In particular, their leading coefficients and powers are those
of \eqref{ibubblerel}, and
$\mathcal B_{-i}(u)\mathcal B_i(-u)=-1$.

These formulas are the change of normalization for ordinary
Kac--Moody generators in \cite[Section~2.3]{BWW25}, followed by
the oriented--unoriented correspondence in
\cite[Theorem~3.7]{BWW25}. To see the crossing factor explicitly,
\eqref{ci1} gives, for $i,j>0$,
\[
 \frac{\kappa_i(\lambda-\bar\alpha_j)}{\kappa_i(\lambda)}
 =t_{ij}^{-1},\qquad t_{ii}=1.
\]
The double mate on the right of
Definition~\ref{natrual-basic1}(2) is multiplied by
$t_{ij}/t_{ji}$. Replacing its original coefficient
$t_{ji}^{-1}$ by $r_{ji}=1$ therefore multiplies the crossing
by $t_{ij}$, as in \eqref{eq:forward-crossing-normalization}.
The first mixed-orientation crossing uses only the unchanged
right cup and cap. In the second, the rescaled downward crossing
contributes $t_{ij}$, while the left cup and cap contribute
\[
 \kappa_i(\lambda-\bar\alpha_i+\bar\alpha_j)^{-1}
 \kappa_i(\lambda)=t_{ij}^{-1}.
\]

Thus both mixed-orientation crossings are unchanged. The
negative-color KLR polynomial is
$Q_{-i,-j}(X,Y)=Q_{i,j}(-X,-Y)$, as required by \eqref{Qij}.

The zero-color generators and the crossings involving zero are
unchanged. Indeed, $\bar\alpha_0=0$ and $t_{ii}=1$ imply
\[
 \kappa_i(\lambda-\bar\alpha_0)=\kappa_i(\lambda),\qquad
 \kappa_i(\lambda\pm\bar\alpha_i)=\kappa_i(\lambda).
\]
The first mate in Definition~\ref{zero-neg_cr} uses only the
right cup and cap. The second uses one left cup and one left cap,
whose scalar factors cancel by these equalities. The explicit
operators in that definition therefore remain the same.
In subsequent calculations, oriented diagrams continue to denote
the auxiliary generators; their images in the presentation of
Section~2 are given by
\eqref{eq:forward-dot-normalization}--\eqref{eq:forward-bubble-normalization}.

\section{Type \texorpdfstring{$A$}{A} and mixed relations}
\label{sec:forward-type-A-mixed}
We retain the assumptions and notation of
Section~\ref{sec:forward-generators}. Unless stated otherwise,
$j,k\in I_1$ and $V\in\mathcal M_\lambda$ is nonzero and finitely generated.
For nonzero $W\in\mathcal M_\mu$, write $c_W=c_\mu$.
Coefficients indexed by zero objects are omitted, and the corresponding
diagrams represent zero morphisms. We use the auxiliary oriented
generators of Section~4 and the conversion in
Subsection~\ref{subsec:forward-normalization}.

The variable $u$ is formal. On $E_iV$, the black dot $x$ has
nilpotent shift $x-i$. Functors act coefficientwise on rational
functions. Zero-color coefficients are simplified and expanded at
$u=0$ before dot evaluation, as in
Lemma~\ref{lem:finite-dd-evaluation}.
Hats denote the formal series $\widehat A_V,\widehat B_V$
defined in Section~4; identities for polynomial representatives
are understood after sufficiently high truncation and nilpotent evaluation.

\begin{Lemma}\label{km2c}
After the conversion in Subsection~\ref{subsec:forward-normalization},
the natural transformations in Definitions~\ref{natrual-basic}
and~\ref{natrual-basic1}, including the bubble coefficients, satisfy
the relations in Definition~\ref{def-2ikac} involving only nonzero
type $A$ vertices.
\end{Lemma}

\begin{proof}
By \cite[Remark~5.2]{GRS}, the transport formulas of
\cite[Lemma~5.1]{GRS} hold at the nonzero vertices in the integral
characteristic-zero case. The nonexceptional cases of
\cite[Lemmas~5.10--5.17]{GRS} therefore apply: their proofs use
these formulas and the identities among $O_V(u)$, $m_V(u)$ and $n_V(u)$.
They give the relations in the auxiliary oriented presentation.
The conversion in Subsection~\ref{subsec:forward-normalization}
gives the asserted relations.
\end{proof}

\begin{Lemma}\label{fw:Gamma-positive}
The following identities hold on $E_jV$ and $E_{-j}V$,
respectively, when these objects are nonzero:
\begin{multicols}{2}
\begin{enumerate}[label=(\arabic*),ref=\arabic*,leftmargin=*,itemsep=.5ex,topsep=0pt]
\item $\Gamma_{E_jV}(u)=E_j(\Gamma_V(u))
 \dfrac{u^2-x^2}{u^2-(x+1)^2}$,
\item $\Gamma_{E_{-j}V}(u)=E_{-j}(\Gamma_V(u))
 \dfrac{u^2-(x-1)^2}{u^2-x^2}$.
\end{enumerate}
\end{multicols}
\end{Lemma}
\begin{proof}
For $s>0$, put
$H_{V,s}(u)=t_{V,-s}(u)/(u+s)^{\eps_{-s}(V)}$.
The primary-factor comparison in Lemma~\ref{ou}(2) gives
\begin{equation} \label{eq:Gamma-negative-exact}
 \begin{aligned} \frac{H_{E_jV,s}(u)}{E_j(H_{V,s}(u))}
 & =\frac{(u+x+1)^{\delta_{s,j+1}}(u+x-1)^{\delta_{s,j-1}}}
 {(u+x)^{2\delta_{s,j}}}\\
\frac{H_{E_{-j}V,s}(u)}{E_{-j}(H_{V,s}(u))}
 & =\frac
 {(u-x)^{2\delta_{s,j}}} {(u-x+1)^{\delta_{s,j+1}}(u-x-1)^{\delta_{s,j-1}}}.  \\
 \end{aligned}
 \end{equation}
To justify these identities, first cross-multiply the bubble-slide identity to obtain a common
polynomial. Reduce $x$ to $j$ in the first case and to $-j$ in the second case, only to determine the primary supports of its factors.
Uniqueness of the monic primary factors of this common polynomial then identifies the exact
factors before reduction. The argument is the same as in the proof of
Lemma~\ref{gamma-move1}.

Put $a_s=\delta_{s,j-1}+\delta_{s,j+1}-2\delta_{s,j}$ for
$s\ge1$. The weight shift is
$d_s(\lambda-\sigma\bar\alpha_j)-d_s(\lambda)=\sigma a_s$
for $\sigma\in\{1,-1\}$, and $\sum_{s\ge1}s a_s=0$.
For the first identity, set $v=x-j$ and
\[
 P_r=(v+r)^2-u^2,\qquad
 K_s=\frac{P_1^{s-1}}{\prod_{r=2}^{s}P_r}\quad(s\ge1).
\]
Substitution in Definition~\ref{GammaLambda} gives
$\Gamma_{E_jV}/E_j(\Gamma_V)=\prod_{s\ge1}K_s^{a_s}
=P_j/P_{j+1}$, which is (1).
Here $K_1=1$, so $j=1$ is included.
For (2), take $v=-x-j$; the exponents are $-a_s$, and the
product is $P_{j+1}/P_j$.
\end{proof}

\begin{Lemma}\label{fw:parity-negative}
Let $j\neq0$. If $V\in\cM_\lambda$ and $E_jV\neq0$, then
$c_{E_jV}\equiv c_V+\delta_{|j|, 1}\pmod2$.
\end{Lemma}
\begin{proof}
By Definition~\ref{def-wei} and the definition of $c_\lambda$ in
\eqref{clambda}, $c_\lambda$ is the $\overline\Lambda_0$-coordinate of
$\lambda$ modulo two. Moreover, Lemma~\ref{weight-pro}(5) gives $E_jV\in\cM_{\lambda-\overline\alpha_j}$.
Since
$
 \alpha_j=2\Lambda_j-\Lambda_{j-1}-\Lambda_{j+1}$,
the $\overline\Lambda_0$-coordinate of $\overline\alpha_j$ is
$\delta_{|j|,1}$ modulo two. Substituting  $\bar\alpha_j$ proves the  claim.
\end{proof}

\begin{Lemma}\label{FGj}
Put $\delta=\delta_{j,1}$, $e=1-\delta$,
$F_j(u,x)=\frac{(u-x)^2}{(u-x-1)(u-x+1)^e}$ and
$G_j(u,x)=\frac{(u+x+1)(u+x-1)^e}{(u+x)^2}$. Then
\begin{enumerate}[label=(\arabic*),ref=\arabic*,leftmargin=*,itemsep=.5ex,topsep=.5ex]
\item $c^+_{E_jV}(u)=E_j(c^+_V(u))F_j(u,x)$;
\item $c^-_{E_jV}(u)=E_j(c^-_V(u))G_j(u,x)$;
\item $F_j(u,x)G_j(-u,x)=(-1)^\delta$.
\end{enumerate}
\end{Lemma}

\begin{proof}
Apply the monic primary-factor comparison from
Lemma~\ref{gamma-move1} to the slide identity in
Lemma~\ref{ou}\allowbreak(2), separating positive and negative nonzero roots.
When $j=1$, the factors with reduced root zero belong to the zero
primary component and are omitted; this accounts for $e=0$.
This gives (1)--(2). Substitution of $-u$ gives (3).
\end{proof}

\begin{Lemma} \label{fw:positive-pitchfork}
The two pitchfork relations \eqref{ipitchfork} hold for
$i=0$ and $j\in I_1$.

\end{Lemma}
\begin{proof}

Suppose first $j=1$. Expand the red cup and cap by
Lemma~\ref{eq:red-cup-E0V} and Definition~\ref{eq:red-cup}. Transport their
coefficients using Lemmas~\ref{fw:Gamma-positive} and
\ref{fw:parity-negative}, together with the exact primary-factor formulas
in Lemma~\ref{FGj}. After simplifying the complete rational coefficient, slide the
$0$-strand decoration through the mixed crossing and around the cup; the
zero variable changes sign at the bend. The affine Brauer relations then
give
\begin{align*}
\begin{tikzpicture}[scale=.65,transform shape,x=1cm,y=1cm,line width=.65pt,line cap=round,line join=round,baseline={([yshift=-.5ex]current bounding box.center)}]
\draw[draw=darkred] (0,2) .. controls (0,.8) and (.324,.28) .. (.9,.28) .. controls (1.476,.28) and (1.8,.8) .. (1.8,2);
\draw[draw=darkred] (.9,2) .. controls (.9,1.30) and (2.2,1.1) .. (2.2,0);
\node[font=\scriptsize,text=darkred,anchor=north,inner sep=.7pt,outer sep=0pt] at (2.2,-.05) {$1$};
\node[font=\scriptsize,text=darkred,anchor=south,inner sep=.7pt,outer sep=0pt] at (0,2.1) {$0$};
\node[font=\scriptsize,text=darkred,anchor=south,inner sep=.7pt,outer sep=0pt] at (0.9,2.1) {$1$};
\node[font=\scriptsize,text=darkred,anchor=south,inner sep=.7pt,outer sep=0pt] at (1.8,2.1) {$0$};
\draw[draw=darkg,line width=.6pt] ([xshift=10pt,yshift=22pt]2.2,-0.6)--([xshift=10pt]2.2,2.2);
\node[text=darkg,font=\scriptsize,anchor=north,inner sep=0pt,outer sep=0pt] at ([xshift=10pt,yshift=-1.2pt]2.2,0) {$\color{darkg}V$};
\end{tikzpicture}&=  %
\begin{tikzpicture}[scale=.55,transform shape,x=1cm,y=1cm,line width=.65pt,line cap=round,line join=round,baseline={([yshift=-.5ex]current bounding box.center)}]
\draw[draw=black] (0,2) .. controls (0,.8) and (.324,.28) .. (.9,.28) .. controls (1.476,.28) and (1.8,.8) .. (1.8,2);
\draw[draw=black] (.9,2) .. controls (.9,1.30) and (2.2,1.1) .. (2.2,0);
\node[font=\scriptsize,text=black,anchor=north,inner sep=.7pt,outer sep=0pt] at (2.2,-0.05) {$1$};
\node[font=\scriptsize,text=black,anchor=south,inner sep=.7pt,outer sep=0pt] at (0,2.1) {$0$};
\node[font=\scriptsize,text=black,anchor=south,inner sep=.7pt,outer sep=0pt] at (0.9,2.1) {$1$};
\node[font=\scriptsize,text=black,anchor=south,inner sep=.7pt,outer sep=0pt] at (1.8,2.1) {$0$};
\fill[black] (1.465632,0.53984) circle (1.65pt);
\node[font=\scriptsize,text=black,anchor=south east,inner sep=0pt,outer sep=0pt] at (1.485632,0.69984) {$B_{E_1V}(x)$};
\draw[densely dotted,line width=.4pt] (1.78634083,1.6)--(1.05758517,1.6);
\node[text=black,anchor=west,inner sep=.5pt,outer sep=0pt] at (1.92634083,1.6) {$\frac{x_2-x_1}{1-x_2+x_1}$};
\end{tikzpicture}
=
\begin{tikzpicture}[scale=.65,transform shape,x=1cm,y=1cm,line width=.65pt,line cap=round,line join=round,baseline={([yshift=-.5ex]current bounding box.center)}]
\draw[draw=black] (0,2) .. controls (0,.8) and (.324,.28) .. (.9,.28) .. controls (1.476,.28) and (1.8,.8) .. (1.8,2);
\draw[draw=black] (.9,2) .. controls (.9,1.30) and (2.2,1.1) .. (2.2,0);
\node[font=\scriptsize,text=black,anchor=north,inner sep=.7pt,outer sep=0pt] at (2.2,0) {$1$};
\node[font=\scriptsize,text=black,anchor=south,inner sep=.7pt,outer sep=0pt] at (0,2) {$0$};
\node[font=\scriptsize,text=black,anchor=south,inner sep=.7pt,outer sep=0pt] at (0.9,2) {$1$};
\node[font=\scriptsize,text=black,anchor=south,inner sep=.7pt,outer sep=0pt] at (1.8,2) {$0$};
\fill[black] (1.465632,0.53984) circle (1.65pt);
\node[font=\scriptsize,text=black,anchor=south east,inner sep=0pt,outer sep=0pt] at (1.485632,0.69984) {$B_V(x)$};
\draw[densely dotted,line width=.4pt] (1.78634083,1.6)--(1.05758517,1.6);
\node[text=black,anchor=west,inner sep=.5pt,outer sep=0pt] at (1.92634083,1.75) {$\frac{x_2-x_1}{1-x_2+x_1}$};
\draw[densely dotted,line width=.4pt] (1.59177803,0.74)--(1.95357656,0.74);
\node[text=black,anchor=west,inner sep=.5pt,outer sep=0pt] at (2.23357656,0.82) {$\frac{(x_1+x_2)^2}{x_1+x_2+1}\frac{x_1^2-(x_2+1)^2}{x_1^2-x_2^2}(-1)$};
\end{tikzpicture}
=
\begin{tikzpicture}[scale=.65,transform shape,x=1cm,y=1cm,line width=.65pt,line cap=round,line join=round,baseline={([yshift=-.5ex]current bounding box.center)}]
\draw[draw=black] (0,2) .. controls (0,.8) and (.324,.28) .. (.9,.28) .. controls (1.476,.28) and (1.8,.8) .. (1.8,2);
\draw[draw=black] (.9,2) .. controls (.9,1.30) and (2.2,1.1) .. (2.2,0);
\node[font=\scriptsize,text=black,anchor=north,inner sep=.7pt,outer sep=0pt] at (2.2,-.05) {$1$};
\node[font=\scriptsize,text=black,anchor=south,inner sep=.7pt,outer sep=0pt] at (0,2.1) {$0$};
\node[font=\scriptsize,text=black,anchor=south,inner sep=.7pt,outer sep=0pt] at (0.9,2.1) {$1$};
\node[font=\scriptsize,text=black,anchor=south,inner sep=.7pt,outer sep=0pt] at (1.8,2.1) {$0$};
\fill[black] (1.465632,0.53984) circle (1.65pt);
\node[font=\scriptsize,text=black,anchor=south east,inner sep=0pt,outer sep=0pt] at (1.485632,0.69984) {$B_V(x)$};
\draw[densely dotted,line width=.4pt] (1.78634083,1.6)--(1.05758517,1.6);
\node[font=\scriptsize,text=black,anchor=west,inner sep=.5pt,outer sep=0pt] at (1.92634083,1.75) {$\frac{x_2-x_1}{1-x_2+x_1}$};
\draw[densely dotted,line width=.4pt] (1.59177803,0.74)--(1.95357656,0.74);
\node[font=\scriptsize,text=black,anchor=west,inner sep=.5pt,outer sep=0pt] at (2.23357656,0.82) {$\frac{x_1+x_2}{x_1-x_2}(x_1-x_2-1)(-1)$};
\end{tikzpicture}
\\
&=
\begin{tikzpicture}[scale=.65,transform shape, x=1cm,y=1cm,line width=.65pt,line cap=round,line join=round,baseline={([yshift=-.5ex]current bounding box.center)}]
\draw[draw=black] (0,2) .. controls (0,.8) and (.324,.28) .. (.9,.28) .. controls (1.476,.28) and (1.8,.8) .. (1.8,2);
\draw[draw=black] (.9,2) .. controls (.9,1.30) and (2.2,1.1) .. (2.2,0);
\node[font=\scriptsize,text=black,anchor=north,inner sep=.7pt,outer sep=0pt] at (2.2,-0.05) {$1$};
\node[font=\scriptsize,text=black,anchor=south,inner sep=.7pt,outer sep=0pt] at (0,2.1) {$0$};
\node[font=\scriptsize,text=black,anchor=south,inner sep=.7pt,outer sep=0pt] at (0.9,2.1) {$1$};
\node[font=\scriptsize,text=black,anchor=south,inner sep=.7pt,outer sep=0pt] at (1.8,2.1) {$0$};
\fill[black] (1.761696,1.36032) circle (1.65pt);
\node[font=\scriptsize,text=black,anchor=west,inner sep=.3pt,outer sep=0pt] at (1.891696,1.38032) {$B_V(x)$};
\draw[densely dotted,line width=.4pt] (1.59177803,0.74)--(1.95357656,0.74);
\node[font=\scriptsize,text=black,anchor=west,inner sep=.5pt,outer sep=0pt] at (2.23357656,0.82) {$\frac{x_1-x_2}{1-x_1+x_2}\frac{x_1+x_2}{x_1-x_2}(1+x_2-x_1)$};
\end{tikzpicture}
=
\begin{tikzpicture}[scale=.65,transform shape,x=1cm,y=1cm,line width=.65pt,line cap=round,line join=round,baseline={([yshift=-.5ex]current bounding box.center)}]
\draw[draw=black] (0,2) .. controls (0,.8) and (.324,.28) .. (.9,.28) .. controls (1.476,.28) and (1.8,.8) .. (1.8,2);
\draw[draw=black] (.9,2) .. controls (.9,1.30) and (2.2,1.1) .. (2.2,0);
\node[font=\scriptsize,text=black,anchor=north,inner sep=.7pt,outer sep=0pt] at (2.2,0) {$1$};
\node[font=\scriptsize,text=black,anchor=south,inner sep=.7pt,outer sep=0pt] at (0,2) {$0$};
\node[font=\scriptsize,text=black,anchor=south,inner sep=.7pt,outer sep=0pt] at (0.9,2) {$1$};
\node[font=\scriptsize,text=black,anchor=south,inner sep=.7pt,outer sep=0pt] at (1.8,2) {$0$};
\fill[black] (1.761696,1.36032) circle (1.65pt);
\node[font=\scriptsize,text=black,anchor=west,inner sep=.3pt,outer sep=0pt] at (1.891696,1.38032) {$B_V(x)$};
\draw[densely dotted,line width=.4pt] (1.59177803,0.74)--(1.95357656,0.74);
\node[font=\scriptsize,text=black,anchor=west,inner sep=.5pt,outer sep=0pt] at (2.09357656,0.74) {$x_1+x_2$};
\end{tikzpicture}
=  %
\begin{tikzpicture}[scale=.65,transform shape,x=1cm,y=1cm,line width=.65pt,line cap=round,line join=round,baseline={([yshift=-.5ex]current bounding box.center)}]
\draw[draw=black] (0,2) .. controls (0,.8) and (.324,.28) .. (.9,.28) .. controls (1.476,.28) and (1.8,.8) .. (1.8,2);
\draw[draw=black] (.9,2) .. controls (.9,1.30) and (-0.4,1.1) .. (-0.4,0);
\node[font=\scriptsize,text=black,anchor=north,inner sep=.7pt,outer sep=0pt] at (-0.4,-0.05) {$1$};
\node[font=\scriptsize,text=black,anchor=south,inner sep=.7pt,outer sep=0pt] at (0,2.1) {$0$};
\node[font=\scriptsize,text=black,anchor=south,inner sep=.7pt,outer sep=0pt] at (0.9,2.1) {$1$};
\node[font=\scriptsize,text=black,anchor=south,inner sep=.7pt,outer sep=0pt] at (1.8,2.1) {$0$};
\fill[black] (1.465632,0.53984) circle (1.65pt);
\node[font=\scriptsize,text=black,anchor=south east,inner sep=0pt,outer sep=0pt] at (1.485632,0.69984) {$B_V(x)$};
\draw[densely dotted,line width=.4pt] (0.01365917,1.6)--(0.74241483,1.6);
\node[font=\scriptsize,text=black,anchor=east,inner sep=.5pt,outer sep=0pt] at (-0.10634083,1.6) {$x_2-x_1$};
\end{tikzpicture}
 =  %
\begin{tikzpicture}[scale=.65,transform shape,x=1cm,y=1cm,line width=.65pt,line cap=round,line join=round,baseline={([yshift=-.5ex]current bounding box.center)}]
\draw[draw=darkred] (0,2) .. controls (0,.8) and (.324,.28) .. (.9,.28) .. controls (1.476,.28) and (1.8,.8) .. (1.8,2);
\draw[draw=darkred] (.9,2) .. controls (.9,1.30) and (-0.4,1.1) .. (-0.4,0);
\node[font=\scriptsize,text=darkred,anchor=north,inner sep=.7pt,outer sep=0pt] at (-0.4,-0.05) {$1$};
\node[font=\scriptsize,text=darkred,anchor=south,inner sep=.7pt,outer sep=0pt] at (0,2.1) {$0$};
\node[font=\scriptsize,text=darkred,anchor=south,inner sep=.7pt,outer sep=0pt] at (0.9,2.1) {$1$};
\node[font=\scriptsize,text=darkred,anchor=south,inner sep=.7pt,outer sep=0pt] at (1.8,2.1) {$0$};
\draw[draw=darkg,line width=.6pt] ([xshift=10pt,yshift=22pt]1.8,-0.5)--([xshift=10pt]1.8,2.2);
\node[text=darkg,font=\scriptsize,anchor=north,inner sep=0pt,outer sep=0pt] at ([xshift=10pt,yshift=-1.2pt]1.8,0) {$\color{darkg}V$};
\end{tikzpicture}.
\end{align*}
Suppose $j\geq2$.  Move the decoration below the crossing, through the
bend, and back above on the other strand; the two factors cancel.  Thus
\begin{align*}
\begin{tikzpicture}[scale=.65,transform shape,x=1cm,y=1cm,line width=.65pt,line cap=round,line join=round,baseline={([yshift=-.5ex]current bounding box.center)}]
\draw[draw=darkred] (0,2) .. controls (0,.8) and (.324,.28) .. (.9,.28) .. controls (1.476,.28) and (1.8,.8) .. (1.8,2);
\draw[draw=darkred] (.9,2) .. controls (.9,1.30) and (2.2,1.1) .. (2.2,0);
\node[font=\scriptsize,text=darkred,anchor=north,inner sep=.7pt,outer sep=0pt] at (2.2,-.05) {$j$};
\node[font=\scriptsize,text=darkred,anchor=south,inner sep=.7pt,outer sep=0pt] at (0,2.1) {$0$};
\node[font=\scriptsize,text=darkred,anchor=south,inner sep=.7pt,outer sep=0pt] at (0.9,2.1) {$j$};
\node[font=\scriptsize,text=darkred,anchor=south,inner sep=.7pt,outer sep=0pt] at (1.8,2.1) {$0$};
\draw[draw=darkg,line width=.6pt] ([xshift=10pt,yshift=22pt]2.2,-0.5)--([xshift=10pt]2.2,2.2);
\node[text=darkg,font=\scriptsize,anchor=north,inner sep=0pt,outer sep=0pt] at ([xshift=10pt,yshift=-1.2pt]2.2,0) {$\color{darkg}V$};
\end{tikzpicture}&=  %
\begin{tikzpicture}[scale=.65,transform shape,x=1cm,y=1cm,line width=.65pt,line cap=round,line join=round,baseline={([yshift=-.5ex]current bounding box.center)}]
\draw[draw=black] (0,2) .. controls (0,.8) and (.324,.28) .. (.9,.28) .. controls (1.476,.28) and (1.8,.8) .. (1.8,2);
\draw[draw=black] (.9,2) .. controls (.9,1.30) and (2.2,1.1) .. (2.2,0);
\node[font=\scriptsize,text=black,anchor=north,inner sep=.7pt,outer sep=0pt] at (2.2,-.05) {$j$};
\node[font=\scriptsize,text=black,anchor=south,inner sep=.7pt,outer sep=0pt] at (0,2.1) {$0$};
\node[font=\scriptsize,text=black,anchor=south,inner sep=.7pt,outer sep=0pt] at (0.9,2.1) {$j$};
\node[font=\scriptsize,text=black,anchor=south,inner sep=.7pt,outer sep=0pt] at (1.8,2.1) {$0$};
\fill[black] (1.465632,0.53984) circle (1.65pt);
\node[font=\scriptsize,text=black,anchor=south east,inner sep=0pt,outer sep=0pt] at (1.485632,0.69984) {$B_{E_jV}^+$};
\draw[densely dotted,line width=.4pt] (1.78634083,1.6)--(1.05758517,1.6);
\node[font=\scriptsize,text=black,anchor=west,inner sep=.5pt,outer sep=0pt] at (1.92634083,1.6) {$\frac{x_2-x_1}{1-x_2+x_1}$};
\end{tikzpicture}
 =  %
\begin{tikzpicture}[scale=.650,transform shape,x=1cm,y=1cm,line width=.65pt,line cap=round,line join=round,baseline={([yshift=-.5ex]current bounding box.center)}]
\draw[draw=black] (0,2) .. controls (0,.8) and (.324,.28) .. (.9,.28) .. controls (1.476,.28) and (1.8,.8) .. (1.8,2);
\draw[draw=black] (.9,2) .. controls (.9,1.30) and (2.2,1.1) .. (2.2,0);
\node[font=\scriptsize,text=black,anchor=north,inner sep=.7pt,outer sep=0pt] at (2.2,-.05) {$j$};
\node[font=\scriptsize,text=black,anchor=south,inner sep=.7pt,outer sep=0pt] at (0,2.1) {$0$};
\node[font=\scriptsize,text=black,anchor=south,inner sep=.7pt,outer sep=0pt] at (0.9,2.1) {$j$};
\node[font=\scriptsize,text=black,anchor=south,inner sep=.7pt,outer sep=0pt] at (1.8,2.1) {$0$};
\fill[black] (1.465632,0.53984) circle (1.65pt);
\node[font=\scriptsize,text=black,anchor=south east,inner sep=0pt,outer sep=0pt] at (1.485632,0.69984) {$B_V^+$};
\draw[densely dotted,line width=.4pt] (1.78634083,1.6)--(1.05758517,1.6);
\node[font=\scriptsize,text=black,anchor=west,inner sep=.5pt,outer sep=0pt] at (1.92634083,1.75) {$\frac{x_2-x_1}{1-x_2+x_1}$};
\draw[densely dotted,line width=.4pt] (1.59177803,0.74)--(1.95357656,0.74);
\node[font=\scriptsize,text=black,anchor=west,inner sep=.5pt,outer sep=0pt] at (2.23357656,0.82) {$\frac{(x_1+x_2)^2}{(x_1+x_2)^2-1}\frac{x_1^2-(x_2+1)^2}{x_1^2-x_2^2}$};
\end{tikzpicture}
 =  %
\begin{tikzpicture}[scale=.65,transform shape,x=1cm,y=1cm,line width=.65pt,line cap=round,line join=round,baseline={([yshift=-.5ex]current bounding box.center)}]
\draw[draw=black] (0,2) .. controls (0,.8) and (.324,.28) .. (.9,.28) .. controls (1.476,.28) and (1.8,.8) .. (1.8,2);
\draw[draw=black] (.9,2) .. controls (.9,1.30) and (2.2,1.1) .. (2.2,0);
\node[font=\scriptsize,text=black,anchor=north,inner sep=.7pt,outer sep=0pt] at (2.2,-.05) {$j$};
\node[font=\scriptsize,text=black,anchor=south,inner sep=.7pt,outer sep=0pt] at (0,2.1) {$0$};
\node[font=\scriptsize,text=black,anchor=south,inner sep=.7pt,outer sep=0pt] at (0.9,2.1) {$j$};
\node[font=\scriptsize,text=black,anchor=south,inner sep=.7pt,outer sep=0pt] at (1.8,2.1) {$0$};
\fill[black] (1.465632,0.53984) circle (1.65pt);
\node[font=\scriptsize,text=black,anchor=south east,inner sep=0pt,outer sep=0pt] at (1.485632,0.69984) {$B_V(x)$};
\draw[densely dotted,line width=.4pt] (1.78634083,1.6)--(1.05758517,1.6);
\node[font=\scriptsize,text=black,anchor=west,inner sep=.5pt,outer sep=0pt] at (1.92634083,1.6) {$\frac{-(x_1+x_2)}{x_1+x_2-1}$};
\end{tikzpicture}
 =  %
\begin{tikzpicture}[scale=.650,transform shape,x=1cm,y=1cm,line width=.65pt,line cap=round,line join=round,baseline={([yshift=-.5ex]current bounding box.center)}]
\draw[draw=black] (0,2) .. controls (0,.8) and (.324,.28) .. (.9,.28) .. controls (1.476,.28) and (1.8,.8) .. (1.8,2);
\draw[draw=black] (.9,2) .. controls (.9,1.30) and (-0.4,1.1) .. (-0.4,0);
\node[font=\scriptsize,text=black,anchor=north,inner sep=.7pt,outer sep=0pt] at (-0.4,-.05) {$j$};
\node[font=\scriptsize,text=black,anchor=south,inner sep=.7pt,outer sep=0pt] at (0,2.1) {$0$};
\node[font=\scriptsize,text=black,anchor=south,inner sep=.7pt,outer sep=0pt] at (0.9,2.1) {$j$};
\node[font=\scriptsize,text=black,anchor=south,inner sep=.7pt,outer sep=0pt] at (1.8,2.1) {$0$};
\fill[black] (1.465632,0.53984) circle (1.65pt);
\node[font=\scriptsize,text=black,anchor=south east,inner sep=0pt,outer sep=0pt] at (1.485632,0.69984) {$B_V(x)$};
\draw[densely dotted,line width=.4pt] (0.01365917,1.6)--(0.74241483,1.6);
\node[font=\scriptsize,text=black,anchor=east,inner sep=.5pt,outer sep=0pt] at (-0.10634083,1.6) {$\frac{x_2-x_1}{1-x_2+x_1}$};
\end{tikzpicture}
 =  %
\begin{tikzpicture}[scale=.650,transform shape,x=1cm,y=1cm,line width=.65pt,line cap=round,line join=round,baseline={([yshift=-.5ex]current bounding box.center)}]
\draw[draw=darkred] (0,2) .. controls (0,.8) and (.324,.28) .. (.9,.28) .. controls (1.476,.28) and (1.8,.8) .. (1.8,2);
\draw[draw=darkred] (.9,2) .. controls (.9,1.30) and (-0.4,1.1) .. (-0.4,0);
\node[font=\scriptsize,text=darkred,anchor=north,inner sep=.7pt,outer sep=0pt] at (-0.4,-.05) {$j$};
\node[font=\scriptsize,text=darkred,anchor=south,inner sep=.7pt,outer sep=0pt] at (0,2.1) {$0$};
\node[font=\scriptsize,text=darkred,anchor=south,inner sep=.7pt,outer sep=0pt] at (0.9,2.1) {$j$};
\node[font=\scriptsize,text=darkred,anchor=south,inner sep=.7pt,outer sep=0pt] at (1.8,2.1) {$0$};
\draw[draw=darkg,line width=.6pt] ([xshift=10pt,yshift=22pt]1.8,-0.5)--([xshift=10pt]1.8,2);
\node[text=darkg,font=\scriptsize,anchor=north,inner sep=0pt,outer sep=0pt] at ([xshift=10pt,yshift=-1.2pt]1.8,0) {$\color{darkg}V$};
\end{tikzpicture}.
\end{align*}
This completes the proof of the first identity  in \eqref{ipitchfork}.

For the cap identity, put $e=1-\delta_{j,1}$. The transport formulas give
\[
 \frac{\widehat A_{E_jV}(u)}{E_j(\widehat A_V(u))}
 =\frac{(u-x+1)^e(u+x)}{(u-x)(u+x+1)}.
\]
Multiply by the mixed-crossing coefficient
$(x-u)/(1-x+u)^e$ and then replace $u$ by $-u$ at the cap.
The result is $(u-x)/(1-u+x)$, giving
$$\begin{aligned}
\begin{tikzpicture}[scale=.650,transform shape,x=1cm,y=1cm,line width=.65pt,line cap=round,line join=round,baseline={([yshift=-.5ex]current bounding box.center)}]
\draw[draw=darkred] (0,0) .. controls (0,1.2) and (.324,1.72) .. (.9,1.72) .. controls (1.476,1.72) and (1.8,1.2) .. (1.8,0);
\draw[draw=darkred] (.9,0) .. controls (.9,.70) and (2.2,0.9) .. (2.2,2);
\node[font=\scriptsize,text=darkred,anchor=south,inner sep=.7pt,outer sep=0pt] at (2.2,2.1) {$j$};
\node[font=\scriptsize,text=darkred,anchor=north,inner sep=.7pt,outer sep=0pt] at (0,-.05) {$0$};
\node[font=\scriptsize,text=darkred,anchor=north,inner sep=.7pt,outer sep=0pt] at (0.9,-.05) {$j$};
\node[font=\scriptsize,text=darkred,anchor=north,inner sep=.7pt,outer sep=0pt] at (1.8,-.05) {$0$};
\draw[draw=darkg,line width=.6pt] ([xshift=10pt,yshift=22pt]2.2,-0.5)--([xshift=10pt]2.2,2);
\node[text=darkg,font=\scriptsize,anchor=north,inner sep=0pt,outer sep=0pt] at ([xshift=10pt,yshift=-1.2pt]2.2,0) {$\color{darkg}V$};
\end{tikzpicture}&=  %
\begin{tikzpicture}[scale=.650,transform shape,x=1cm,y=1cm,line width=.65pt,line cap=round,line join=round,baseline={([yshift=-.5ex]current bounding box.center)}]
\draw[draw=black] (0,0) .. controls (0,1.2) and (.324,1.72) .. (.9,1.72) .. controls (1.476,1.72) and (1.8,1.2) .. (1.8,0);
\draw[draw=black] (.9,0) .. controls (.9,.70) and (2.2,0.9) .. (2.2,2);
\node[font=\scriptsize,text=black,anchor=south,inner sep=.7pt,outer sep=0pt] at (2.2,2.1) {$j$};
\node[font=\scriptsize,text=black,anchor=north,inner sep=.7pt,outer sep=0pt] at (0,-.05) {$0$};
\node[font=\scriptsize,text=black,anchor=north,inner sep=.7pt,outer sep=0pt] at (0.9,-.05) {$j$};
\node[font=\scriptsize,text=black,anchor=north,inner sep=.7pt,outer sep=0pt] at (1.8,-0.05){$0$};
\fill[black] (1.465632,1.46016) circle (1.65pt);
\node[font=\scriptsize,anchor=north east,inner sep=0pt,outer sep=0pt] at (1.485632,1.30016) {$A_{E_jV}(x)$};
\draw[densely dotted,line width=.4pt] (1.78634083,0.4)--(1.05758517,0.4);
\node[font=\scriptsize,anchor=west,inner sep=.5pt,outer sep=0pt] at (1.92634083,0.4) {$\frac{x_1-x_2}{(1-x_1+x_2)^{1-\delta_{j,1}}}$};
\end{tikzpicture}=\begin{tikzpicture}[scale=.550,transform shape,x=1cm,y=1cm,line width=.65pt,line cap=round,line join=round,baseline={([yshift=-.5ex]current bounding box.center)}]
\draw[draw=darkred] (0,0) .. controls (0,1.2) and (.324,1.72) .. (.9,1.72) .. controls (1.476,1.72) and (1.8,1.2) .. (1.8,0);
\draw[draw=darkred] (.9,0) .. controls (.9,.70) and (-0.4,0.9) .. (-0.4,2);
\node[font=\scriptsize,text=darkred,anchor=south,inner sep=.7pt,outer sep=0pt] at (-0.4,2.1) {$j$};
\node[font=\scriptsize,text=darkred,anchor=north,inner sep=.7pt,outer sep=0pt] at (0,-0.05) {$0$};
\node[font=\scriptsize,text=darkred,anchor=north,inner sep=.7pt,outer sep=0pt] at (0.9,-0.05) {$j$};
\node[font=\scriptsize,text=darkred,anchor=north,inner sep=.7pt,outer sep=0pt] at (1.8,-0.05) {$0$};
\draw[draw=darkg,line width=.6pt] ([xshift=10pt,yshift=22pt]1.8,-0.5)--([xshift=10pt]1.8,2);
\node[text=darkg,font=\scriptsize,anchor=north,inner sep=0pt,outer sep=0pt] at ([xshift=10pt,yshift=-1.2pt]1.8,0) {$\color{darkg}V$};
\end{tikzpicture}.
\end{aligned}
$$
This completes the proof of the second identity in \eqref{ipitchfork}. \end{proof}

For the next three lemmas and Proposition~\ref{prop:forward-mixed-braid},
fix $k>0$. Only displayed colors are prescribed: a black diagram
with boundary colors $\boldsymbol a,\boldsymbol b$ means its
component $\operatorname{pr}_{\boldsymbol b}D\operatorname{in}_{\boldsymbol a}$.
An internal color imposes an additional projection; otherwise,
the possible internal colors are summed over.

In a decoration of a marked pair, $x_1,x_2$ are its left and right
dots. If an annotation involves a third index, all its indices
refer to the three strands at that height. We write $x_L$ when
specifying the dot at the left endpoint of a connector explicitly.
Dots on distinct input and output cup--cap branches remain independent.

In each repeated cup--cap correction below, the upper
connector meets the left branch of the upper cup, consistently with
the dot variable used in the coefficient calculation.

\begin{Lemma}\label{mixed-3}

The following identity holds:

\begin{equation}\label{mixed-b1} 
\end{equation}
\end{Lemma}

\begin{proof}
Insert the generalized eigenspace projections at the crossings
and use the component formulas of \cite[Lemmas~4.11--4.12]{GRS}.
The components with unchanged or opposite colors contribute the
cup--cap correction below. Its upper factor is $(z-y)^{-1}$,
where $y$ is the dot on the $k$-strand and $z$ is the dot on the
left zero-colored leg of the upper cup. Sliding the remaining
coefficients gives
$$\begin{aligned} %
\end{aligned}$$
\end{proof}

\begin{Lemma}\label{mixed-1}

The following identity holds:
\begin{equation}\label{eq:mixed-1} \begin{aligned}
& %
.\\
\end{aligned}\end{equation}
\end{Lemma}

\begin{proof}
Let $Z$ be the left-hand side of \eqref{eq:mixed-1}.
Slide $(\ell_{-+}^{12})^{-1}$ through the zero-colored crossing
using $Sx_2-x_1S=e_1-1$, where $S$ is the black crossing
projected to colors $(0,0)$ and $e_1$ is the projected black cup--cap.
Retain both correction terms, and then use the projected
pitchfork and quadratic identities. Internal color labels in the
following calculation are essential.
We obtain
\begingroup
\small
\allowdisplaybreaks[4]
\begin{align*}
Z
&=%
\end{align*}
\endgroup
\end{proof}

\begin{Lemma}\label{mixed-2}
The following identity holds:

\begin{equation} \label{eq:mixed-2} \begin{aligned}
& %
\end{aligned}
\end{equation}
\end{Lemma}
\begin{proof}
Slide the middle decoration through the projected black crossing
of colors $(0,0)$,
retaining its identity and cup--cap terms as in
Lemma~\ref{mixed-1}. The internal zero in the first diagram on
the left fixes a color component, whereas the first diagram on
the right also contains the component with intermediate colors
\[
 (0,0,k)\longrightarrow(k,-k,k)
 \longrightarrow(k,0,0)\longrightarrow(k,0,0).
\]
By \cite[Lemmas~4.11--4.12]{GRS}, its three crossing coefficients
are $-(a+b)^{-1},(b-c)^{-1},(b-c)^{-1}$, with input dots
$(a,-a,b)$ and output dots $(b,c,-c)$ on the resulting cup--cap.
Subtracting this additional component with the displayed lower
decoration gives the first cup--cap correction on the right.
The other two follow from the quadratic relation
$1-(b-c)^{-2}$ for the projected mixed crossings.
\end{proof}

\begin{Prop}[Mixed braiding]\label{prop:forward-mixed-braid}
For $k>0$, relation \eqref{ibraid} on input colors $(0,0,k)$ is:
\begin{equation}\label{fw:mixed-braid}
.
\end{equation}
\end{Prop}
\begin{proof}
Expand the zero-color crossing in each braid by
Lemma~\ref{lem:anothdefofcro} and \eqref{rx1x2}, and expand the two
mixed crossings by Definition~\ref{natrual-basic-0i}. Use the cup/cap moves of
Lemma~\ref{fw:positive-pitchfork}. Move decorations only to the marked
boundary segments; in a disconnected cup--cap term, keep independent
variables on the input cap and output cup. This gives
\begin{equation}\label{mixed-red1} \begin{aligned} %
\end{aligned}
\end{equation}
Next move the remaining decorations through the cup--cap pair, changing
the sign of the variable at each bend.  Thus
\begin{equation}\label{eq:mixed-3}
%
\end{equation}

At every indicated two-strand level, $y_a$ denotes the red dot:
$y_a=x_a$ on a zero strand and $y_a=x_a-k$ on a positive $k$-strand.
The index refers to the marked level before or after the crossing. For
$\Delta=x_1-x_2$, the two projected off-diagonal black crossings have
product $1-\Delta^{-2}$. Hence the corresponding red double crossing is
$(y_1-y_2)^{\delta_{k,1}}$ on colors $(0,k)$ and
$(y_2-y_1)^{\delta_{k,1}}$ on colors $(k,0)$.

Let $Z_1$ (respectively, $Z_2$) be the third term on the right-hand side
of the first (respectively, second) identity in \eqref{mixed-red1}.  Move
$h_1=p_{E_kV}(x_3,x_1)$ and $h_2=q_V(x_3,x_1)$ below the final crossing;
the correction terms are exactly those displayed below.  Then
\begin{equation}\label{mixed-b2} \begin{aligned}
Z_1
&=

\end{align*}
\endgroup
In the last cup--cap term, let $x_1$ be the lower cap's right dot,
$x_2$ the through-strand dot, and $x_3$ the upper cup's left dot.
The boundary dot lists are
$(-x_1,x_1,x_2)$ and $(x_2,x_3,-x_3)$.
After moving the decorations around the bends, the two divided
differences are $h_1=p_{E_kV}(-x_3,x_1)$ and
$h_2=q_V(-x_3,x_1)$ in these boundary variables.
Lemma~\ref{compu-T} gives
$T=-\delta_{k,1}A_{E_kV}(x_1)B_V(-x_3)$; with
Lemma~\ref{eq:red-cup-E0V}, this is the red correction shown.
All generators in this relation are unchanged by
Subsection~\ref{subsec:forward-normalization}.
\end{proof}

Before treating the negative mixed pitchforks, we record both
zero-color triangle identities. Let $F=E_0$, and let
$\eta:\mathrm{Id}\Rightarrow F^2$ and $\varepsilon:F^2\Rightarrow\mathrm{Id}$
be the projected black cup and cap. By
Section~\ref{sec:forward-generators}, the maps
$a_V=A_V(x):FV\to FV$ are natural. They are invertible because
$A_V(0)B_V(0)=1$ and $x$ is nilpotent. Thus $a$ is a natural
automorphism of $F$. The red cup and cap
are
\[
 \eta_V^0=a_{FV}^{-1}\eta_V,\qquad
 \varepsilon_V^0=\varepsilon_VF(a_V).
\]
The first triangle identity follows by cancellation:
\[
 \varepsilon^0_{FV}F(\eta_V^0)
 =\varepsilon_{FV}F(a_{FV})F(a_{FV}^{-1})F(\eta_V)=\Id_{FV}.
\]
For the second, naturality of $a$ and of the black cup gives
\[
\begin{aligned}
 F(\varepsilon_V^0)\eta_{FV}^0
 &=F(\varepsilon_V)F^2(a_V)a_{F^2V}^{-1}\eta_{FV}\\
 &=a_V^{-1}F(\varepsilon_V)\eta_{FV}a_V
 =\Id_{FV}.
\end{aligned}
\]
Consequently, the zero-color mate operation is a bijection; its
inverse is obtained by bending the strand back using these cups
and caps.

\begin{Lemma}[Pitchfork]\label{fw:nonzero-pitchfork}

Let $j\in I_1$. The following pitchfork relations hold:
\begin{multicols}{2}
\begin{enumerate}[label=(\arabic*),ref=\arabic*,leftmargin=*,itemsep=1ex,topsep=0pt]
\item $\begin{tikzpicture}[scale=0.55,x=1cm,y=1cm,line width=.65pt,line cap=round,line join=round,baseline={([yshift=-.5ex]current bounding box.center)}]
\draw[draw=darkred] (0,2) .. controls (0,0.8) and (0.324,0.28) .. (0.9,0.28) .. controls (1.476,0.28) and (1.8,0.8) .. (1.8,2);
\draw[draw=darkred] (0.9,2) .. controls (0.9,1.3) and (2.2,1.1) .. (2.2,0);
\node[font=\scriptsize,text=darkred,anchor=north,inner sep=.7pt,outer sep=0pt] at (2.2,-0.05) {$0$};
\node[font=\scriptsize,text=darkred,anchor=south,inner sep=.7pt,outer sep=0pt] at (0,2.1) {$-j$};
\node[font=\scriptsize,text=darkred,anchor=south,inner sep=.7pt,outer sep=0pt] at (0.9,2.1) {$0$};
\node[font=\scriptsize,text=darkred,anchor=south,inner sep=.7pt,outer sep=0pt] at (1.8,2.1) {$j$};
\end{tikzpicture}=  %
\begin{tikzpicture}[scale=0.55,x=1cm,y=1cm,line width=.65pt,line cap=round,line join=round,baseline={([yshift=-.5ex]current bounding box.center)}]
\draw[draw=darkred] (0,2) .. controls (0,0.8) and (0.324,0.28) .. (0.9,0.28) .. controls (1.476,0.28) and (1.8,0.8) .. (1.8,2);
\draw[draw=darkred] (0.9,2) .. controls (0.9,1.3) and (-0.4,1.1) .. (-0.4,0);
\node[font=\scriptsize,text=darkred,anchor=north,inner sep=.7pt,outer sep=0pt] at (-0.4,-0.05) {$0$};
\node[font=\scriptsize,text=darkred,anchor=south,inner sep=.7pt,outer sep=0pt] at (0,2.1) {$-j$};
\node[font=\scriptsize,text=darkred,anchor=south,inner sep=.7pt,outer sep=0pt] at (0.9,2.1) {$0$};
\node[font=\scriptsize,text=darkred,anchor=south,inner sep=.7pt,outer sep=0pt] at (1.8,2.1) {$j$};
\end{tikzpicture}   $

\item $ \begin{tikzpicture}[scale=0.55,x=1cm,y=1cm,line width=.65pt,line cap=round,line join=round,baseline={([yshift=-.5ex]current bounding box.center)}]
\draw[draw=darkred] (0,0) .. controls (0,1.2) and (0.324,1.72) .. (0.9,1.72) .. controls (1.476,1.72) and (1.8,1.2) .. (1.8,0);
\draw[draw=darkred] (0.9,0) .. controls (0.9,0.7) and (-0.4,0.9) .. (-0.4,2);
\node[font=\scriptsize,text=darkred,anchor=south,inner sep=.7pt,outer sep=0pt] at (-0.4,2.1) {$0$};
\node[font=\scriptsize,text=darkred,anchor=north,inner sep=.7pt,outer sep=0pt] at (0,-0.05) {$j$};
\node[font=\scriptsize,text=darkred,anchor=north,inner sep=.7pt,outer sep=0pt] at (0.9,-0.05) {$0$};
\node[font=\scriptsize,text=darkred,anchor=north,inner sep=.7pt,outer sep=0pt] at (1.8,-0.05) {$-j$};
\end{tikzpicture}=  %
\begin{tikzpicture}[scale=0.62,x=1cm,y=1cm,line width=.65pt,line cap=round,line join=round,baseline={([yshift=-.5ex]current bounding box.center)}]
\draw[draw=darkred] (0,0) .. controls (0,1.2) and (0.324,1.72) .. (0.9,1.72) .. controls (1.476,1.72) and (1.8,1.2) .. (1.8,0);
\draw[draw=darkred] (0.9,0) .. controls (0.9,0.7) and (2.2,0.9) .. (2.2,2);
\node[font=\scriptsize,text=darkred,anchor=south,inner sep=.7pt,outer sep=0pt] at (2.2,2.1) {$0$};
\node[font=\scriptsize,text=darkred,anchor=north,inner sep=.7pt,outer sep=0pt] at (0,-0.05) {$j$};
\node[font=\scriptsize,text=darkred,anchor=north,inner sep=.7pt,outer sep=0pt] at (0.9,-0.05) {$0$};
\node[font=\scriptsize,text=darkred,anchor=north,inner sep=.7pt,outer sep=0pt] at (1.8,-0.05) {$-j$};
\end{tikzpicture}  $

\item $ \begin{tikzpicture}[scale=0.55,x=1cm,y=1cm,line width=.65pt,line cap=round,line join=round,baseline={([yshift=-.5ex]current bounding box.center)}]
\draw[draw=darkred] (0,2) .. controls (0,0.8) and (0.324,0.28) .. (0.9,0.28) .. controls (1.476,0.28) and (1.8,0.8) .. (1.8,2);
\draw[draw=darkred] (0.9,2) .. controls (0.9,1.3) and (2.2,1.1) .. (2.2,0);
\node[font=\scriptsize,text=darkred,anchor=north,inner sep=.7pt,outer sep=0pt] at (2.2,-0.05) {$0$};
\node[font=\scriptsize,text=darkred,anchor=south,inner sep=.7pt,outer sep=0pt] at (0,2.1) {$j$};
\node[font=\scriptsize,text=darkred,anchor=south,inner sep=.7pt,outer sep=0pt] at (0.9,2.1) {$0$};
\node[font=\scriptsize,text=darkred,anchor=south,inner sep=.7pt,outer sep=0pt] at (1.8,2.1) {$-j$};
\end{tikzpicture}=  %
\begin{tikzpicture}[scale=0.55,x=1cm,y=1cm,line width=.65pt,line cap=round,line join=round,baseline={([yshift=-.5ex]current bounding box.center)}]
\draw[draw=darkred] (0,2) .. controls (0,0.8) and (0.324,0.28) .. (0.9,0.28) .. controls (1.476,0.28) and (1.8,0.8) .. (1.8,2);
\draw[draw=darkred] (0.9,2) .. controls (0.9,1.3) and (-0.4,1.1) .. (-0.4,0);
\node[font=\scriptsize,text=darkred,anchor=north,inner sep=.7pt,outer sep=0pt] at (-0.4,-0.05) {$0$};
\node[font=\scriptsize,text=darkred,anchor=south,inner sep=.7pt,outer sep=0pt] at (0,2.1) {$j$};
\node[font=\scriptsize,text=darkred,anchor=south,inner sep=.7pt,outer sep=0pt] at (0.9,2.1) {$0$};
\node[font=\scriptsize,text=darkred,anchor=south,inner sep=.7pt,outer sep=0pt] at (1.8,2.1) {$-j$};
\end{tikzpicture} $

\item $\begin{tikzpicture}[scale=0.55,x=1cm,y=1cm,line width=.65pt,line cap=round,line join=round,baseline={([yshift=-.5ex]current bounding box.center)}]
\draw[draw=darkred] (0,0) .. controls (0,1.2) and (0.324,1.72) .. (0.9,1.72) .. controls (1.476,1.72) and (1.8,1.2) .. (1.8,0);
\draw[draw=darkred] (0.9,0) .. controls (0.9,0.7) and (-0.4,0.9) .. (-0.4,2);
\node[font=\scriptsize,text=darkred,anchor=south,inner sep=.7pt,outer sep=0pt] at (-0.4,2.1) {$0$};
\node[font=\scriptsize,text=darkred,anchor=north,inner sep=.7pt,outer sep=0pt] at (0,-0.05) {$-j$};
\node[font=\scriptsize,text=darkred,anchor=north,inner sep=.7pt,outer sep=0pt] at (0.9,-0.05) {$0$};
\node[font=\scriptsize,text=darkred,anchor=north,inner sep=.7pt,outer sep=0pt] at (1.8,-.05) {$j$};
\end{tikzpicture}=  %
\begin{tikzpicture}[scale=0.55,x=1cm,y=1cm,line width=.65pt,line cap=round,line join=round,baseline={([yshift=-.5ex]current bounding box.center)}]
\draw[draw=darkred] (0,0) .. controls (0,1.2) and (0.324,1.72) .. (0.9,1.72) .. controls (1.476,1.72) and (1.8,1.2) .. (1.8,0);
\draw[draw=darkred] (0.9,0) .. controls (0.9,0.7) and (2.2,0.9) .. (2.2,2);
\node[font=\scriptsize,text=darkred,anchor=south,inner sep=.7pt,outer sep=0pt] at (2.2,2.1) {$0$};
\node[font=\scriptsize,text=darkred,anchor=north,inner sep=.7pt,outer sep=0pt] at (0,-.05) {$-j$};
\node[font=\scriptsize,text=darkred,anchor=north,inner sep=.7pt,outer sep=0pt] at (0.9,-.05) {$0$};
\node[font=\scriptsize,text=darkred,anchor=north,inner sep=.7pt,outer sep=0pt] at (1.8,-.05) {$j$};
\end{tikzpicture}$

\item
$\begin{tikzpicture}[scale=0.55,x=1cm,y=1cm,line width=.65pt,line cap=round,line join=round,baseline={([yshift=-.5ex]current bounding box.center)}]
\draw[draw=darkred] (0,2) .. controls (0,0.8) and (0.324,0.28) .. (0.9,0.28) .. controls (1.476,0.28) and (1.8,0.8) .. (1.8,2);
\draw[draw=darkred] (0.9,2) .. controls (0.9,1.3) and (2.2,1.1) .. (2.2,0);
\node[font=\scriptsize,text=darkred,anchor=north,inner sep=.7pt,outer sep=0pt] at (2.2,-0.05) {$-j$};
\node[font=\scriptsize,text=darkred,anchor=south,inner sep=.7pt,outer sep=0pt] at (0,2.1) {$0$};
\node[font=\scriptsize,text=darkred,anchor=south,inner sep=.7pt,outer sep=0pt] at (0.9,2.1) {$-j$};
\node[font=\scriptsize,text=darkred,anchor=south,inner sep=.7pt,outer sep=0pt] at (1.8,2.1) {$0$};
\end{tikzpicture}=  %
\begin{tikzpicture}[scale=0.55,x=1cm,y=1cm,line width=.65pt,line cap=round,line join=round,baseline={([yshift=-.5ex]current bounding box.center)}]
\draw[draw=darkred] (0,2) .. controls (0,0.8) and (0.324,0.28) .. (0.9,0.28) .. controls (1.476,0.28) and (1.8,0.8) .. (1.8,2);
\draw[draw=darkred] (0.9,2) .. controls (0.9,1.3) and (-0.4,1.1) .. (-0.4,0);
\node[font=\scriptsize,text=darkred,anchor=north,inner sep=.7pt,outer sep=0pt] at (-0.4,-.05) {$-j$};
\node[font=\scriptsize,text=darkred,anchor=south,inner sep=.7pt,outer sep=0pt] at (0,2.1) {$0$};
\node[font=\scriptsize,text=darkred,anchor=south,inner sep=.7pt,outer sep=0pt] at (0.9,2.1) {$-j$};
\node[font=\scriptsize,text=darkred,anchor=south,inner sep=.7pt,outer sep=0pt] at (1.8,2.1) {$0$};
\end{tikzpicture}$
 \item $
\begin{tikzpicture}[scale=0.55,x=1cm,y=1cm,line width=.65pt,line cap=round,line join=round,baseline={([yshift=-.5ex]current bounding box.center)}]
\draw[draw=darkred] (0,0) .. controls (0,1.2) and (0.324,1.72) .. (0.9,1.72) .. controls (1.476,1.72) and (1.8,1.2) .. (1.8,0);
\draw[draw=darkred] (0.9,0) .. controls (0.9,0.7) and (-0.4,0.9) .. (-0.4,2);
\node[font=\scriptsize,text=darkred,anchor=south,inner sep=.7pt,outer sep=0pt] at (-0.4,2.1) {$-j$};
\node[font=\scriptsize,text=darkred,anchor=north,inner sep=.7pt,outer sep=0pt] at (0,-0.05) {$0$};
\node[font=\scriptsize,text=darkred,anchor=north,inner sep=.7pt,outer sep=0pt] at (0.9,-0.05) {$-j$};
\node[font=\scriptsize,text=darkred,anchor=north,inner sep=.7pt,outer sep=0pt] at (1.8,-0.05) {$0$};
\end{tikzpicture}=  %
\begin{tikzpicture}[scale=0.55,x=1cm,y=1cm,line width=.65pt,line cap=round,line join=round,baseline={([yshift=-.5ex]current bounding box.center)}]
\draw[draw=darkred] (0,0) .. controls (0,1.2) and (0.324,1.72) .. (0.9,1.72) .. controls (1.476,1.72) and (1.8,1.2) .. (1.8,0);
\draw[draw=darkred] (0.9,0) .. controls (0.9,0.7) and (2.2,0.9) .. (2.2,2);
\node[font=\scriptsize,text=darkred,anchor=south,inner sep=.7pt,outer sep=0pt] at (2.2,2.1) {$-j$};
\node[font=\scriptsize,text=darkred,anchor=north,inner sep=.7pt,outer sep=0pt] at (0,-0.05) {$0$};
\node[font=\scriptsize,text=darkred,anchor=north,inner sep=.7pt,outer sep=0pt] at (0.9,-0.05) {$-j$};
\node[font=\scriptsize,text=darkred,anchor=north,inner sep=.7pt,outer sep=0pt] at (1.8,-0.05) {$0$};
\end{tikzpicture}$.
\end{enumerate}
\end{multicols}
\end{Lemma}

\begin{proof}
Relations (1)--(4) follow from Definition~\ref{zero-neg_cr}
and the zig-zag identities for the $j$- and $(-j)$-colored strands.
For example, expanding the negative mixed crossing in the right-hand
side of (1) and cancelling the resulting cup--cap pair gives
    \[
\begin{tikzpicture}[scale=0.55,x=1cm,y=1cm,line width=.65pt,line cap=round,line join=round,baseline={([yshift=-.5ex]current bounding box.center)}]
\draw[draw=darkred] (0,2) .. controls (0,0.8) and (0.324,0.28) .. (0.9,0.28) .. controls (1.476,0.28) and (1.8,0.8) .. (1.8,2);
\draw[draw=darkred] (0.9,2) .. controls (0.9,1.3) and (-0.4,1.1) .. (-0.4,0);
\node[font=\scriptsize,text=darkred,anchor=north,inner sep=.7pt,outer sep=0pt] at (-0.4,-0.05) {$0$};
\node[font=\scriptsize,text=darkred,anchor=south,inner sep=.7pt,outer sep=0pt] at (0,2.1) {$-j$};
\node[font=\scriptsize,text=darkred,anchor=south,inner sep=.7pt,outer sep=0pt] at (0.9,2.1) {$0$};
\node[font=\scriptsize,text=darkred,anchor=south,inner sep=.7pt,outer sep=0pt] at (1.8,2.1) {$j$};
\end{tikzpicture}
=
\begin{tikzpicture}[scale=0.648,transform shape,x=1cm,y=1cm,line width=.65pt,line cap=round,line join=round,baseline={([yshift=-.5ex]current bounding box.center)}]
\draw[draw=darkred] (0,2) .. controls (0,-0.3) and (1,-0.3) .. (1,1.25) .. controls (1,2.55) and (2,2.55) .. (2,1.25) .. controls (2,-0.3) and (3,-0.3) .. (3,2);
\draw[draw=darkred] (0.75,2) .. controls (0.75,1.25) and (1.7,0.75) .. (1.7,0);
\node[font=\scriptsize,text=darkred,anchor=south,inner sep=.7pt,outer sep=0pt] at (0.75,2.1) {$0$};
\node[font=\scriptsize,text=darkred,anchor=north,inner sep=.7pt,outer sep=0pt] at (1.7,-0.05) {$0$};
\node[font=\scriptsize,text=darkred,anchor=south,inner sep=.7pt,outer sep=0pt] at (0,2.1) {$-j$};
\node[font=\scriptsize,text=darkred,anchor=south,inner sep=.7pt,outer sep=0pt] at (3,2.1) {$j$};
\end{tikzpicture}
=
\begin{tikzpicture}[scale=0.55,x=1cm,y=1cm,line width=.65pt,line cap=round,line join=round,baseline={([yshift=-.5ex]current bounding box.center)}]
\draw[draw=darkred] (0,2) .. controls (0,0.8) and (0.324,0.28) .. (0.9,0.28) .. controls (1.476,0.28) and (1.8,0.8) .. (1.8,2);
\draw[draw=darkred] (0.9,2) .. controls (0.9,1.3) and (2.2,1.1) .. (2.2,0);
\node[font=\scriptsize,text=darkred,anchor=north,inner sep=.7pt,outer sep=0pt] at (2.2,-0.05) {$0$};
\node[font=\scriptsize,text=darkred,anchor=south,inner sep=.7pt,outer sep=0pt] at (0,2.1) {$-j$};
\node[font=\scriptsize,text=darkred,anchor=south,inner sep=.7pt,outer sep=0pt] at (0.9,2.1) {$0$};
\node[font=\scriptsize,text=darkred,anchor=south,inner sep=.7pt,outer sep=0pt] at (1.8,2.1) {$j$};
\end{tikzpicture}
.
\]
This proves (1). For (2)--(4), expand the corresponding negative mixed
crossing and cancel the adjacent cup--cap pair in the same way.

It remains to prove (5)--(6). By the zero-color triangle identities,
bending the two zero strands converts (5) into (6), and the
operation is reversible. Explicitly, inserting a zero-color
zig-zag into each side of (5) gives
\begin{align*}
\begin{tikzpicture}[scale=0.65,transform shape,x=1cm,y=1cm,line width=.65pt,line cap=round,line join=round,baseline={([yshift=-.5ex]current bounding box.center)}]
\draw[draw=darkred] (0,2) .. controls (0,0.8) and (0.324,0.28) .. (0.9,0.28) .. controls (1.476,0.28) and (1.8,0.8) .. (1.8,2);
\draw[draw=darkred] (0.9,2) .. controls (0.9,1.3) and (2.2,1.1) .. (2.2,0);
\node[font=\scriptsize,text=darkred,anchor=north,inner sep=.7pt,outer sep=0pt] at (2.2,0) {$-j$};
\node[font=\scriptsize,text=darkred,anchor=south,inner sep=.7pt,outer sep=0pt] at (0,2) {$0$};
\node[font=\scriptsize,text=darkred,anchor=south,inner sep=.7pt,outer sep=0pt] at (0.9,2) {$-j$};
\node[font=\scriptsize,text=darkred,anchor=south,inner sep=.7pt,outer sep=0pt] at (1.8,2) {$0$};
\end{tikzpicture}&=\begin{tikzpicture}[scale=0.648,transform shape,x=1cm,y=1cm,line width=.65pt,line cap=round,line join=round,baseline={([yshift=-.5ex]current bounding box.center)}]
\draw[draw=darkred] (0,2) .. controls (0,-0.3) and (1,-0.3) .. (1,1.25) .. controls (1,2.55) and (2,2.55) .. (2,1.25) .. controls (2,-0.3) and (3,-0.3) .. (3,2);
\draw[draw=darkred] (0.6,2) .. controls (0.6,1.25) and (1.4,0.75) .. (1.4,0);
\node[font=\scriptsize,text=darkred,anchor=south,inner sep=.7pt,outer sep=0pt] at (0.6,2) {$-j$};
\node[font=\scriptsize,text=darkred,anchor=north,inner sep=.7pt,outer sep=0pt] at (1.4,0) {$-j$};
\node[font=\scriptsize,text=darkred,anchor=south,inner sep=.7pt,outer sep=0pt] at (0,2) {$0$};
\node[font=\scriptsize,text=darkred,anchor=south,inner sep=.7pt,outer sep=0pt] at (3,2) {$0$};
\end{tikzpicture},\quad
\begin{tikzpicture}[scale=0.65,transform shape,x=1cm,y=1cm,line width=.65pt,line cap=round,line join=round,baseline={([yshift=-.5ex]current bounding box.center)}]
\draw[draw=darkred] (0,2) .. controls (0,0.8) and (0.324,0.28) .. (0.9,0.28) .. controls (1.476,0.28) and (1.8,0.8) .. (1.8,2);
\draw[draw=darkred] (0.9,2) .. controls (0.9,1.3) and (-0.4,1.1) .. (-0.4,0);
\node[font=\scriptsize,text=darkred,anchor=north,inner sep=.7pt,outer sep=0pt] at (-0.4,0) {$-j$};
\node[font=\scriptsize,text=darkred,anchor=south,inner sep=.7pt,outer sep=0pt] at (0,2) {$0$};
\node[font=\scriptsize,text=darkred,anchor=south,inner sep=.7pt,outer sep=0pt] at (0.9,2) {$-j$};
\node[font=\scriptsize,text=darkred,anchor=south,inner sep=.7pt,outer sep=0pt] at (1.8,2) {$0$};
\end{tikzpicture}=\begin{tikzpicture}[scale=0.648,transform shape,x=1cm,y=1cm,line width=.65pt,line cap=round,line join=round,baseline={([yshift=-.5ex]current bounding box.center)}]
\draw[draw=darkred] (0,2) .. controls (0,-0.3) and (1,-0.3) .. (1,1.25) .. controls (1,2.55) and (2,2.55) .. (2,1.25) .. controls (2,-0.3) and (3,-0.3) .. (3,2);
\draw[draw=darkred] (2.4,2) .. controls (2.4,1.25) and (1.6,0.75) .. (1.6,0);
\node[font=\scriptsize,text=darkred,anchor=south,inner sep=.7pt,outer sep=0pt] at (2.4,2) {$-j$};
\node[font=\scriptsize,text=darkred,anchor=north,inner sep=.7pt,outer sep=0pt] at (1.6,0) {$-j$};
\node[font=\scriptsize,text=darkred,anchor=south,inner sep=.7pt,outer sep=0pt] at (0,2) {$0$};
\node[font=\scriptsize,text=darkred,anchor=south,inner sep=.7pt,outer sep=0pt] at (3,2) {$0$};
\end{tikzpicture}.
\end{align*}
The two diagrams on the right are equal once (6) is known.
We therefore prove (6).

Set $e=1-\delta_{j,1}$. We may assume $E_{-j}V\ne0$;
otherwise both sides of (6) are zero. Let $x$ be the black dot
on $E_{-j}V$, so that $x+j$ is nilpotent. The exact transport
identity for the formal series is
\begin{equation}\label{eq:negative-A-transport}
 \frac{\widehat A_{E_{-j}V}(u)}{E_{-j}(\widehat A_V(u))}
 =\frac{(u+x)(u-x+1)}{(u-x)(u+x+1)^e}.
\end{equation}
This follows from the positive primary factors in
Lemma~\ref{ou}(2) and Lemma~\ref{fw:Gamma-positive}(2), and hence
also holds for the polynomial representatives after nilpotent evaluation.
For the nonzero-color factors, put $N=E_{-j}V$, let $z$ be the black  dot
on $E_0N$ when this object is nonzero, and let $v$ be a formal
variable. To compare the factors away from $j$, put
$C_{N,j}(v)=n_{N,j}(v)/m_{N,j}(v)$. Primary-factor extraction from Lemma~\ref{ou}(2) gives
\[
 \frac{C_{E_0N,j}(v)}{E_0(C_{N,j}(v))}
 =\frac{(v+z+1)(v+z-1)^e(v-z)^2}
 {(v-z+1)(v-z-1)^e(v+z)^2}.
\]
These formulas specify the variables used in the following
calculation. We use the notation $\ell_{\sigma\tau}^{ab}$, $d_{\sigma\tau}^{ab}$,
$a_{V,i}$ and $b_{V,i}$ from Definition~\ref{ration111}.
First expand the negative crossing on the left-hand side of (6).
By Definitions~\ref{zero-neg_cr}, \ref{natrual-basic-0i}
and \ref{eq:red-cup}, the dot-sliding relations give

\begin{align*}
\begin{tikzpicture}[scale=1,x=1cm,y=1cm,
  baseline={([yshift=-0.5ex]current bounding box.center)},
  line cap=round,line join=round]
\begin{scope}
  \clip (-1.5,-0.4) rectangle (1.5,0.8);
  \draw[-,thick,darkred,name path=p150B] (-.35,-.5) to[out=0,in=180] (.35,.5);
  \draw[-,thick,darkred,name path=p150C] (.35,.5) to[out=0,in=90] (.8,-1.2);
  \draw[-,thick,darkred,name path=p150E] (-.35,.4) to[out=90,in=180] (.5,1.1);
  \draw[-,thick,darkred,name path=p150A] (.35,-.4) to[out=90,in=-90] (-.35,.4);
\end{scope}
\pgfresetboundingbox
\path (-0.36,-0.41) rectangle (0.84,0.8);
\node[text=darkred,anchor=south,inner sep=0.4pt,outer sep=0pt]
  at (-0.20,0.97) {$\scriptstyle -j$};
\node[text=darkred,anchor=north,inner sep=0.4pt,outer sep=0pt]
  at (-0.14,-0.57) {$\scriptstyle 0$};
\node[text=darkred,anchor=north,inner sep=0.4pt,outer sep=0pt]
  at (0.35,-0.57) {$\scriptstyle -j$};
\node[text=darkred,anchor=north,inner sep=0.4pt,outer sep=0pt]
  at (0.84,-0.57) {$\scriptstyle 0$};
\draw[darkg,line width=0.6pt] (1.49,-0.40)--(1.49,0.50);
\node[text=darkg,anchor=north,inner sep=0.4pt,outer sep=0pt]
  at (1.49,-0.57) {$\scriptstyle V$};
\end{tikzpicture}
&=\begin{tikzpicture}[scale=0.65,x=1cm,y=1cm,
  baseline={([yshift=-0.5ex]current bounding box.center)},
  line cap=round,line join=round]
\begin{scope}
  \clip (-1.5,-0.65) rectangle (1.5,1.2);
  \draw[-,thick,black,name path=p152F] (.35,-.4) to[out=-90,in=0] (-.5,-1.1);
  \draw[-,thick,black,name path=p152A] (.35,-.4) to[out=90,in=-90] (-.35,.4);
  \draw[-,thick,black,name path=p152E] (-.35,.4) to[out=90,in=180] (.5,1.1);
  \draw[-,thick,black,name path=p152O] (.5,1.1) to[out=0,in=90]  coordinate[pos=0.79] (dotA152) (1.3,.4);
  \draw[-,thick,black,name path=p152V] (1.3,.4) to (1.3,-1.2);
  \draw[-,thick,black,name path=p152D] (-.35,-.5) to[out=180,in=-90] (-.8,1.2);
  \draw[-,thick,black,name path=p152B] (-.35,-.5) to[out=0,in=180] (.35,.5);
  \draw[-,thick,black,name path=p152C] (.35,.5) to[out=0,in=90] (.8,-1.2);
  \path[name path=h152s] (-1.40,-0.22)--(1.40,-0.22);
  \path[name intersections={of=p152B and h152s,by=d152sL}];
  \path[name intersections={of=p152A and h152s,by=d152sR}];
  \fill[black] (dotA152) circle (1.65pt);
\end{scope}
\pgfresetboundingbox
\path (-0.84,-0.66) rectangle (1.31,1.2);
\draw[line width=0.4pt,line cap=round,
  dash pattern=on 0pt off 1.1pt] (d152sL)--(d152sR);
\node[text=black,anchor=south,inner sep=0.4pt,outer sep=0pt]
  at (-0.80,1.37) {$\scriptstyle -j$};
\node[text=black,anchor=north,inner sep=0.4pt,outer sep=0pt]
  at (0.30,-0.82) {$\scriptstyle 0$};
\node[text=black,anchor=north,inner sep=0.4pt,outer sep=0pt]
  at (0.83,-0.82) {$\scriptstyle -j$};
\node[text=black,anchor=north,inner sep=0.4pt,outer sep=0pt]
  at (1.30,-0.82) {$\scriptstyle 0$};
\node[text=black,anchor=west,inner sep=0.4pt,outer sep=0pt]
  at ($(dotA152)+(0.16,0.02)$) {$\scriptstyle A_V^+$};
\node[text=black,anchor=west,inner sep=0.4pt,outer sep=0pt]
  at ($(d152sR)+(1.17,0.00)$) {$\scriptstyle d_{+-}^{12}(\ell_{-+}^{12})^{\delta_{j,1}-1}$};
\end{tikzpicture}=\begin{tikzpicture}[scale=1,x=1cm,y=1cm,
  baseline={([yshift=-0.5ex]current bounding box.center)},
  line cap=round,line join=round]
\begin{scope}
  \clip (-1.5,-0.4) rectangle (1.5,0.8);
  \draw[-,thick,black,name path=p154C] (-.35,.5) to[out=180,in=90] (-.8,-1.2);
  \draw[-,thick,black,name path=p154B] (.35,-.5) to[out=180,in=0]  coordinate[pos=0.31] (dotA154) coordinate[pos=0.39] (dashLow154) (-.35,.5);
  \draw[-,thick,black,name path=p154E] (.35,.4) to[out=90,in=0] (-.5,1.1);
  \draw[-,thick,black,name path=p154A] (-.35,-.4) to[out=90,in=-90]  coordinate[pos=0.92] (dashHigh154) (.35,.4);
  \fill[black] (dotA154) circle (1.65pt);
\end{scope}
\pgfresetboundingbox
\path (-0.84,-0.41) rectangle (0.36,0.8);
\draw[line width=0.45pt,line cap=round,
  dash pattern=on 0pt off 1.1pt] (dashLow154)
  to[bend right=24] coordinate[pos=0.50] (dashMid154) (dashHigh154);
\node[text=black,anchor=south,inner sep=0.4pt,outer sep=0pt]
  at (0.20,0.97) {$\scriptstyle -j$};
\node[text=black,anchor=north,inner sep=0.4pt,outer sep=0pt]
  at (-0.84,-0.57) {$\scriptstyle 0$};
\node[text=black,anchor=north,inner sep=0.4pt,outer sep=0pt]
  at (-0.35,-0.57) {$\scriptstyle -j$};
\node[text=black,anchor=north,inner sep=0.4pt,outer sep=0pt]
  at (0.14,-0.57) {$\scriptstyle 0$};
\node[text=black,anchor=west,inner sep=0.4pt,outer sep=0pt]
  at ($(dotA154)+(0.17,-0.03)$) {$\scriptstyle A_V^+$};
\node[text=black,anchor=west,inner sep=0.4pt,outer sep=0pt]
  at ($(dashMid154)+(0.31,0.15)$) {$\scriptstyle d_{+-}^{12}(\ell_{-+}^{12})^{\delta_{j,1}-1}$};
\end{tikzpicture}.
\end{align*}
For the right-hand side of (6), expanding the negative crossing gives
\begin{align*}
\begin{tikzpicture}[scale=1,x=1cm,y=1cm,
  baseline={([yshift=-0.5ex]current bounding box.center)},
  line cap=round,line join=round]
\begin{scope}
  \clip (-1.5,-0.4) rectangle (1.5,0.8);
  \draw[-,thick,darkred,name path=p155C] (-.35,.5) to[out=180,in=90] (-.8,-1.2);
  \draw[-,thick,darkred,name path=p155B] (.35,-.5) to[out=180,in=0] (-.35,.5);
  \draw[-,thick,darkred,name path=p155E] (.35,.4) to[out=90,in=0] (-.5,1.1);
  \draw[-,thick,darkred,name path=p155A] (-.35,-.4) to[out=90,in=-90] (.35,.4);
\end{scope}
\pgfresetboundingbox
\path (-0.84,-0.41) rectangle (0.36,0.8);
\node[text=darkred,anchor=south,inner sep=0.4pt,outer sep=0pt]
  at (0.20,0.97) {$\scriptstyle -j$};
\node[text=darkred,anchor=north,inner sep=0.4pt,outer sep=0pt]
  at (-0.84,-0.57) {$\scriptstyle 0$};
\node[text=darkred,anchor=north,inner sep=0.4pt,outer sep=0pt]
  at (-0.35,-0.57) {$\scriptstyle -j$};
\node[text=darkred,anchor=north,inner sep=0.4pt,outer sep=0pt]
  at (0.14,-0.57) {$\scriptstyle 0$};
\draw[darkg,line width=0.6pt] (1.01,-0.40)--(1.01,0.50);
\node[text=darkg,anchor=north,inner sep=0.4pt,outer sep=0pt]
  at (1.01,-0.57) {$\scriptstyle V$};
\end{tikzpicture}
&=\begin{tikzpicture}[scale=0.65,x=1cm,y=1cm,
  baseline={([yshift=-0.5ex]current bounding box.center)},
  line cap=round,line join=round]
\begin{scope}
  \clip (-1.5,-0.65) rectangle (1.5,1.2);
  \draw[-,thick,darkred,name path=p156V] (-1.3,.4) to (-1.3,-1.2);
  \draw[-,thick,darkred,name path=p156O] (-.5,1.1) to[out=180,in=90] (-1.3,.4);
  \draw[-,thick,darkred,name path=p156E] (.35,.4) to[out=90,in=0] (-.5,1.1);
  \draw[-,thick,darkred,name path=p156A] (-.35,-.4) to[out=90,in=-90] (.35,.4);
  \draw[-,thick,darkred,name path=p156F] (-.35,-.4) to[out=-90,in=180] (.5,-1.1);
  \draw[-,thick,darkred,name path=p156D] (.35,-.5) to[out=0,in=-90] (.8,1.2);
  \draw[-,thick,darkred,name path=p156B] (.35,-.5) to[out=180,in=0] (-.35,.5);
  \draw[-,thick,darkred,name path=p156C] (-.35,.5) to[out=180,in=90] (-.8,-1.2);
\end{scope}
\pgfresetboundingbox
\path (-1.31,-0.66) rectangle (0.65,1.2);
\node[text=darkred,anchor=south,inner sep=0.4pt,outer sep=0pt]
  at (0.80,1.37) {$\scriptstyle -j$};
\node[text=darkred,anchor=north,inner sep=0.4pt,outer sep=0pt]
  at (-1.30,-0.82) {$\scriptstyle 0$};
\node[text=darkred,anchor=north,inner sep=0.4pt,outer sep=0pt]
  at (-0.83,-0.82) {$\scriptstyle -j$};
\node[text=darkred,anchor=north,inner sep=0.4pt,outer sep=0pt]
  at (-0.30,-0.82) {$\scriptstyle 0$};
\draw[darkg,line width=0.6pt] (1.49,-0.65)--(1.49,0.25);
\node[text=darkg,anchor=north,inner sep=0.4pt,outer sep=0pt]
  at (1.49,-0.82) {$\scriptstyle V$};
\end{tikzpicture}
=\begin{tikzpicture}[scale=0.85,x=1cm,y=1cm,
  baseline={([yshift=-0.5ex]current bounding box.center)},
  line cap=round,line join=round]
\begin{scope}
  \clip (-1.5,-0.65) rectangle (1.5,1.2);
  \draw[-,thick,black,name path=p157V] (-1.3,.4) to (-1.3,-1.2);
  \draw[-,thick,black,name path=p157O] (-.5,1.1) to[out=180,in=90] (-1.3,.4);
  \draw[-,thick,black,name path=p157E] (.35,.4) to[out=90,in=0]  coordinate[pos=0.57] (dotA157) (-.5,1.1);
  \draw[-,thick,black,name path=p157A] (-.35,-.4) to[out=90,in=-90] (.35,.4);
  \draw[-,thick,black,name path=p157F] (-.35,-.4) to[out=-90,in=180] (.5,-1.1);
  \draw[-,thick,black,name path=p157D] (.35,-.5) to[out=0,in=-90] (.8,1.2);
  \draw[-,thick,black,name path=p157B] (.35,-.5) to[out=180,in=0]  coordinate[pos=0.16] (dota157) coordinate[pos=0.90] (dotb157) (-.35,.5);
  \draw[-,thick,black,name path=p157C] (-.35,.5) to[out=180,in=90] (-.8,-1.2);
  \path[name path=h157s] (-1.40,0.20)--(1.40,0.20);
  \path[name intersections={of=p157B and h157s,by=d157sL}];
  \path[name intersections={of=p157A and h157s,by=d157sR}];
  \fill[black] (dotA157) circle (1.65pt);
  \fill[black] (dota157) circle (1.65pt);
\end{scope}
\pgfresetboundingbox
\path (-1.31,-0.66) rectangle (0.84,1.2);
\draw[line width=0.4pt,line cap=round,
  dash pattern=on 0pt off 1.1pt] (d157sL)--(d157sR);
\node[text=black,anchor=south,inner sep=0.4pt,outer sep=0pt]
  at (0.80,1.37) {$\scriptstyle -j$};
\node[text=black,anchor=north,inner sep=0.4pt,outer sep=0pt]
  at (-1.30,-0.82) {$\scriptstyle 0$};
\node[text=black,anchor=north,inner sep=0.4pt,outer sep=0pt]
  at (-0.83,-0.82) {$\scriptstyle -j$};
\node[text=black,anchor=north,inner sep=0.4pt,outer sep=0pt]
  at (-0.30,-0.82) {$\scriptstyle 0$};
\node[text=black,anchor=west,inner sep=0.4pt,outer sep=0pt]
  at (1.02,1.00) {$\scriptstyle A_{E_{-j}V}^+$};
\node[text=black,anchor=north west,inner sep=0.4pt,outer sep=0pt]
  at (1.02,-0.58) {$\scriptstyle a_{E_{-j}V,j}(x)$};
\node[text=black,anchor=east,inner sep=0.4pt,outer sep=0pt]
  at ($(d157sL)+(-1.40,0.00)$) {$\scriptstyle \frac{d_{+-}^{12}b_{E_0E_{-j}V,j}(x_1)}{\ell_{-+}^{12}}$};
\end{tikzpicture}.
\end{align*}
Here the factors $a_{E_{-j}V,j}$ and $b_{E_0E_{-j}V,j}$ are those of
Definition~\ref{ration111}(3)--(4). Moving the factor $A_{E_{-j}V}$ to
the $0$-colored strand below the crossing and  using Lemma~\ref{ou}(2) for the $m$- and $n$-factors
and Lemma~\ref{fw:Gamma-positive}(2) for the $\Gamma$-factor, we obtain
\begin{equation}\label{key111}
\begin{tikzpicture}[scale=0.65,x=1cm,y=1cm,
  baseline={([yshift=-0.5ex]current bounding box.center)},
  line cap=round,line join=round]
\begin{scope}
  \clip (-1.5,-0.65) rectangle (1.5,1.2);
  \draw[-,thick,black,name path=p159V] (-1.3,.4) to (-1.3,-1.2);
  \draw[-,thick,black,name path=p159O] (-.5,1.1) to[out=180,in=90] (-1.3,.4);
  \draw[-,thick,black,name path=p159E] (.35,.4) to[out=90,in=0] (-.5,1.1);
  \draw[-,thick,black,name path=p159A] (-.35,-.4) to[out=90,in=-90]  coordinate[pos=0.23] (dotA159) (.35,.4);
  \draw[-,thick,black,name path=p159F] (-.35,-.4) to[out=-90,in=180] (.5,-1.1);
  \draw[-,thick,black,name path=p159D] (.35,-.5) to[out=0,in=-90] (.8,1.2);
  \draw[-,thick,black,name path=p159B] (.35,-.5) to[out=180,in=0] (-.35,.5);
  \draw[-,thick,black,name path=p159C] (-.35,.5) to[out=180,in=90] (-.8,-1.2);
  \path[name path=h159l] (-1.40,0.46)--(1.40,0.46);
  \path[name intersections={of=p159B and h159l,by=d159lL}];
  \path[name intersections={of=p159D and h159l,by=d159lR}];
  \fill[black] (dotA159) circle (1.65pt);
\end{scope}
\pgfresetboundingbox
\path (-1.31,-0.66) rectangle (0.84,1.2);
\draw[line width=0.4pt,line cap=round,
  dash pattern=on 0pt off 1.1pt] (d159lL)--(d159lR);
\node[text=black,anchor=south,inner sep=0.4pt,outer sep=0pt]
  at (0.80,1.37) {$\scriptstyle -j$};
\node[text=black,anchor=north,inner sep=0.4pt,outer sep=0pt]
  at (-1.30,-0.82) {$\scriptstyle 0$};
\node[text=black,anchor=north,inner sep=0.4pt,outer sep=0pt]
  at (-0.83,-0.82) {$\scriptstyle -j$};
\node[text=black,anchor=north,inner sep=0.4pt,outer sep=0pt]
  at (-0.30,-0.82) {$\scriptstyle 0$};
\node[text=black,anchor=east,inner sep=0.4pt,outer sep=0pt]
  at ($(dotA159)+(-1.20,0.00)$) {$\scriptstyle A_V^+$};
\node[text=black,anchor=west,inner sep=0.4pt,outer sep=0pt]
  at ($(d159lR)+(0.23,0.00)$) {$\begin{gathered}
 \scriptstyle -\dfrac{(d_{++}^{12})^2}{d_{+-}^{12}\ell_{++}^{12}}
 \left(\dfrac{\ell_{-+}^{12}}{\ell_{--}^{12}}\right)^{1-\delta_{j,1}}\\[3pt]
 \scriptstyle {}\cdot\dfrac{\ell_{+-}^{23}}{d_{+-}^{23}}
 d_{++}^{23}(\ell_{++}^{23})^{\delta_{j,1}-1}
 \end{gathered}$};
\end{tikzpicture}
.\end{equation}

To simplify the coefficient, slide the dots through the
nonzero-colored cup. Its two legs carry $x_1$ and $x_3$, so
$x_1=-x_3$ on this cup, while $x_2$ is the zero-color dot. Hence
\[
\begin{gathered}
 d_{++}^{12}=d_{+-}^{23},\qquad
 d_{+-}^{12}=-d_{++}^{23},\qquad
 \ell_{++}^{12}=\ell_{+-}^{23},\quad
 \ell_{-+}^{12}=\ell_{++}^{23},\qquad
 \ell_{--}^{12}=\ell_{-+}^{23}.
\end{gathered}
\]
Substituting these identities into the displayed coefficient gives
\[
\begin{aligned}
 &\frac{(d_{+-}^{23})^2}{d_{++}^{23}\ell_{+-}^{23}}
 \left(\frac{\ell_{++}^{23}}{\ell_{-+}^{23}}\right)^{1-\delta_{j,1}}
 \frac{\ell_{+-}^{23}}{d_{+-}^{23}}\,
 d_{++}^{23}(\ell_{++}^{23})^{\delta_{j,1}-1}
 =d_{+-}^{23}(\ell_{-+}^{23})^{\delta_{j,1}-1}.
\end{aligned}
\]
Common factors are cancelled before evaluation on the generalized
eigenspaces; the same calculation covers both $j=1$ and $j>1$.
After relabelling the two remaining variables and cancelling the
cup--cap pair, the diagram in \eqref{key111} becomes

$\begin{tikzpicture}[scale=1,x=1cm,y=1cm,
  baseline={([yshift=-0.5ex]current bounding box.center)},
  line cap=round,line join=round]
\begin{scope}
  \clip (-1.5,-0.4) rectangle (1.5,0.8);
  \draw[-,thick,black,name path=p161C] (-.35,.5) to[out=180,in=90] (-.8,-1.2);
  \draw[-,thick,black,name path=p161B] (.35,-.5) to[out=180,in=0]  coordinate[pos=0.31] (dotA161) coordinate[pos=0.39] (dashLow161) (-.35,.5);
  \draw[-,thick,black,name path=p161E] (.35,.4) to[out=90,in=0] (-.5,1.1);
  \draw[-,thick,black,name path=p161A] (-.35,-.4) to[out=90,in=-90]  coordinate[pos=0.92] (dashHigh161) (.35,.4);
  \fill[black] (dotA161) circle (1.65pt);
\end{scope}
\pgfresetboundingbox
\path (-0.84,-0.41) rectangle (0.36,0.8);
\draw[line width=0.45pt,line cap=round,
  dash pattern=on 0pt off 1.1pt] (dashLow161)
  to[bend right=24] coordinate[pos=0.50] (dashMid161) (dashHigh161);
\node[text=black,anchor=south,inner sep=0.4pt,outer sep=0pt]
  at (0.20,0.97) {$\scriptstyle -j$};
\node[text=black,anchor=north,inner sep=0.4pt,outer sep=0pt]
  at (-0.84,-0.57) {$\scriptstyle 0$};
\node[text=black,anchor=north,inner sep=0.4pt,outer sep=0pt]
  at (-0.35,-0.57) {$\scriptstyle -j$};
\node[text=black,anchor=north,inner sep=0.4pt,outer sep=0pt]
  at (0.14,-0.57) {$\scriptstyle 0$};
\node[text=black,anchor=west,inner sep=0.4pt,outer sep=0pt]
  at ($(dotA161)+(0.17,-0.03)$) {$\scriptstyle A_V^+$};
\node[text=black,anchor=west,inner sep=0.4pt,outer sep=0pt]
  at ($(dashMid161)+(0.31,0.15)$) {$\scriptstyle d_{+-}^{12}(\ell_{-+}^{12})^{\delta_{j,1}-1}$};
\end{tikzpicture}$.
This is the same decorated diagram obtained from the left-hand side.
Hence (6) holds.
\end{proof}

These identities also hold in the normalization of Section~2.
The cups and caps in (1)--(2) are unchanged by
\eqref{eq:forward-adjunction-normalization}. In (3)--(4), the
left cup or cap has the same factor on both sides, since
$\kappa_j(\lambda-\bar\alpha_0)=\kappa_j(\lambda)$.
The zero-color cups and caps in (5)--(6), and all crossings
involving zero, are unchanged by the conversion in
Subsection~\ref{subsec:forward-normalization}.
\subsection{Completion of the forward construction}
\label{subsec:completion-forward}
We use the parameters of Section~2, in particular $r_{i,j}=1$.
The generators are converted as in
Subsection~\ref{subsec:forward-normalization}.
We now verify the remaining relations of Definition~\ref{def-2ikac}.

\begin{Lemma}\label{fw:basic-relations}
For $i,j\in I$, the relations \textnormal{\eqref{izigzag}--\eqref{ipitchfork}} hold.
\end{Lemma}

\begin{proof}
For $i\in I$, we have

\begin{equation*}
\begin{tikzpicture}[every path/.style={darkred, thick},scale=0.80,transform shape,baseline={([yshift=-.5ex]current bounding box.center)},x=.65cm,y=.65cm]
  \draw[line width=.65pt,line cap=round,line join=round] (0,1.8)..controls(0,-.2)and(1,-.2)..(1,.9)
  ..controls(1,2)and(2,2)..(2,0);
  \node[font=\scriptsize,text=black,inner sep=.8pt,outer sep=0pt,below]
  at (2,0) {$i$};
  \end{tikzpicture}=  %
\begin{tikzpicture}[every path/.style={darkred, thick},scale=0.80,transform shape,baseline={([yshift=-.5ex]current bounding box.center)},x=1cm,y=1cm]
  \draw[line width=.65pt,line cap=round,line join=round] (0,0)--(0,1.05);
  \node[font=\scriptsize,text=black,inner sep=.8pt,outer sep=0pt,below]
  at (0,0) {$i$};
  \end{tikzpicture}= %
\begin{tikzpicture}[every path/.style={darkred, thick},scale=0.80,transform shape,baseline={([yshift=-.5ex]current bounding box.center)},x=.65cm,y=.65cm]
  \begin{scope}[xscale=-1]
  \draw[line width=.65pt,line cap=round,line join=round] (0,1.8)..controls(0,-.2)and(1,-.2)..(1,.9)
  ..controls(1,2)and(2,2)..(2,0);
  \node[font=\scriptsize,text=black,inner sep=.8pt,outer sep=0pt,below]
  at (2,0) {$i$};
  \end{scope}
  \end{tikzpicture}.\qquad\text{\normalfont\eqref{izigzag}}
\end{equation*}

\begin{equation*}
\begin{tikzpicture}[every path/.style={darkred, thick},scale=0.80,transform shape,baseline={([yshift=-.5ex]current bounding box.center)},x=.65cm,y=.65cm]
  \draw[line width=.65pt,line cap=round,line join=round](0,1)..controls(0,-.15)and(1,-.15)..(1,1);
  \fill(.06075,.560125)circle(2pt);
  \node[font=\scriptsize,text=black,inner sep=.8pt,outer sep=0pt,above]
  at (0,1) {$i$};
  \end{tikzpicture}=- %
\begin{tikzpicture}[every path/.style={darkred, thick},scale=0.80,transform shape,baseline={([yshift=-.5ex]current bounding box.center)},x=.65cm,y=.65cm]
  \draw[line width=.65pt,line cap=round,line join=round](0,1)..controls(0,-.15)and(1,-.15)..(1,1);
  \fill(.93925,.560125)circle(2pt);
  \node[font=\scriptsize,text=black,inner sep=.8pt,outer sep=0pt,above]
  at (0,1) {$i$};
  \end{tikzpicture},\qquad
\begin{tikzpicture}[every path/.style={darkred, thick},scale=0.80,transform shape,baseline={([yshift=-.5ex]current bounding box.center)},x=.65cm,y=.65cm]
  \draw[line width=.65pt,line cap=round,line join=round](0,0)..controls(0,1.15)and(1,1.15)..(1,0);
  \fill(.06075,.439875)circle(2pt);
  \node[font=\scriptsize,text=black,inner sep=.8pt,outer sep=0pt,below]
  at (0,0) {$i$};
  \end{tikzpicture}=- %
\begin{tikzpicture}[every path/.style={darkred, thick},scale=0.80,transform shape,baseline={([yshift=-.5ex]current bounding box.center)},x=.65cm,y=.65cm]
  \draw[line width=.65pt,line cap=round,line join=round](0,0)..controls(0,1.15)and(1,1.15)..(1,0);
  \fill(.93925,.439875)circle(2pt);
  \node[font=\scriptsize,text=black,inner sep=.8pt,outer sep=0pt,below]
  at (0,0) {$i$};
  \end{tikzpicture}.\qquad\text{\normalfont\eqref{ipivotal}}
\end{equation*}

\begin{equation*}
\begin{tikzpicture}[every path/.style={darkred, thick},scale=0.80,transform shape,baseline={([yshift=-.5ex]current bounding box.center)},x=1cm,y=1cm]
  \draw[line width=.65pt,line cap=round,line join=round,] (0,1) .. controls (0,-.08) and (1,-.08) .. (1,1)
   (.5,1)--(1.22,0);
  \node[font=\scriptsize,text=black,inner sep=.8pt,outer sep=0pt,below]
  at (1.22,0) {$j$};
  \node[font=\scriptsize,text=black,inner sep=.8pt,outer sep=0pt,above]
  at (0,1) {$i$};
  \node[font=\scriptsize,text=black,inner sep=.8pt,outer sep=0pt,above]
  at (.5,1) {$j$};
  \node[font=\scriptsize,text=black,inner sep=.8pt,outer sep=0pt,above]
  at (1,1) {$\tau i$};
  \end{tikzpicture}=r_{i,j}\,  %
\begin{tikzpicture}[every path/.style={darkred, thick},scale=0.80,transform shape,baseline={([yshift=-.5ex]current bounding box.center)},x=1cm,y=1cm]
  \draw[line width=.65pt,line cap=round,line join=round] (0,1) .. controls (0,-.08) and (1,-.08) .. (1,1)
   (-.22,0)--(.5,1);
  \node[font=\scriptsize,text=black,inner sep=.8pt,outer sep=0pt,below]
  at (-.22,0) {$j$};
  \node[font=\scriptsize,text=black,inner sep=.8pt,outer sep=0pt,above]
  at (0,1) {$i$};
  \node[font=\scriptsize,text=black,inner sep=.8pt,outer sep=0pt,above]
  at (.5,1) {$j$};
  \node[font=\scriptsize,text=black,inner sep=.8pt,outer sep=0pt,above]
  at (1,1) {$\tau i$};
  \end{tikzpicture},
 \quad
\begin{tikzpicture}[every path/.style={darkred, thick},scale=0.80,transform shape,baseline={([yshift=-.5ex]current bounding box.center)},x=1cm,y=1cm]
  \draw[line width=.65pt,line cap=round,line join=round] (0,0) .. controls (0,1.08) and (1,1.08) .. (1,0)
   (-.22,1)--(.5,0);
  \node[font=\scriptsize,text=black,inner sep=.8pt,outer sep=0pt,above]
  at (-.22,1) {$j$};
  \node[font=\scriptsize,text=black,inner sep=.8pt,outer sep=0pt,below]
  at (0,0) {$\tau i$};
  \node[font=\scriptsize,text=black,inner sep=.8pt,outer sep=0pt,below]
  at (.5,0) {$j$};
  \node[font=\scriptsize,text=black,inner sep=.8pt,outer sep=0pt,below]
  at (1,0) {$i$};
  \end{tikzpicture}=r_{i,j}\,  %
\begin{tikzpicture}[every path/.style={darkred, thick},scale=0.80,transform shape,baseline={([yshift=-.5ex]current bounding box.center)},x=1cm,y=1cm]
  \draw[line width=.65pt,line cap=round,line join=round] (0,0) .. controls (0,1.08) and (1,1.08) .. (1,0)
   (.5,0)--(1.22,1);
  \node[font=\scriptsize,text=black,inner sep=.8pt,outer sep=0pt,above]
  at (1.22,1) {$j$};
  \node[font=\scriptsize,text=black,inner sep=.8pt,outer sep=0pt,below]
  at (0,0) {$\tau i$};
  \node[font=\scriptsize,text=black,inner sep=.8pt,outer sep=0pt,below]
  at (.5,0) {$j$};
  \node[font=\scriptsize,text=black,inner sep=.8pt,outer sep=0pt,below]
  at (1,0) {$i$};
  \end{tikzpicture}.
 \qquad\text{\normalfont\eqref{ipitchfork}}\end{equation*}
For nonzero colors, the relations follow from Lemma~\ref{km2c}.
At zero, the two triangle identities were proved above
Lemma~\ref{fw:nonzero-pitchfork}, the dot--bend identities follow
from Definition~\ref{eq:red-cup}, and the pitchfork relation is
relation~\textnormal{(4)} of Definition~\ref{nil}, verified in
Section~\ref{sec:forward-nil-brauer}.
The mixed pitchforks are Lemmas~\ref{fw:positive-pitchfork}
and~\ref{fw:nonzero-pitchfork} and their mates.
\end{proof}

\begin{Lemma}\label{fw:dot-crossing-relation}
For $i,j\in I$, relation \textnormal{\eqref{idotslide}} holds.
\end{Lemma}

\begin{proof}

\begin{align*}
\begin{tikzpicture}[every path/.style={darkred, thick},scale=0.80,transform shape,baseline={([yshift=-.5ex]current bounding box.center)},x=.75cm,y=.75cm]
  \draw[line width=.65pt,line cap=round,line join=round] (0,0)--(1,1) (0,1)--(1,0);
  \fill(.22,.22)circle(2pt);
  \node[font=\scriptsize,text=black,inner sep=.8pt,outer sep=0pt,below]
  at (0,0) {$i$};
  \node[font=\scriptsize,text=black,inner sep=.8pt,outer sep=0pt,below]
  at (1,0) {$j$};
  \end{tikzpicture}- %
\begin{tikzpicture}[every path/.style={darkred, thick},scale=0.80,transform shape,baseline={([yshift=-.5ex]current bounding box.center)},x=.75cm,y=.75cm]
  \draw[line width=.65pt,line cap=round,line join=round] (0,0)--(1,1) (0,1)--(1,0);
  \fill(.78,.78)circle(2pt);
  \node[font=\scriptsize,text=black,inner sep=.8pt,outer sep=0pt,below]
  at (0,0) {$i$};
  \node[font=\scriptsize,text=black,inner sep=.8pt,outer sep=0pt,below]
  at (1,0) {$j$};
  \end{tikzpicture}
 &=\delta_{i,j}\, %
\begin{tikzpicture}[every path/.style={darkred, thick},scale=0.80,transform shape,baseline={([yshift=-.5ex]current bounding box.center)},x=.7cm,y=.7cm]
  \draw[line width=.65pt,line cap=round,line join=round] (0,0)--(0,1.2) (1,0)--(1,1.2);
  \node[font=\scriptsize,text=black,inner sep=.8pt,outer sep=0pt,below]
  at (0,0) {$i$};
  \node[font=\scriptsize,text=black,inner sep=.8pt,outer sep=0pt,below]
  at (1,0) {$j$};
  \end{tikzpicture}-\delta_{i,\tau j}\, %
\begin{tikzpicture}[every path/.style={darkred, thick},scale=0.80,transform shape,baseline={([yshift=-.5ex]current bounding box.center)},x=.65cm,y=.65cm]
  \draw[line width=.65pt,line cap=round,line join=round] (0,1.8)..controls(0,.8)and(1,.8)..(1,1.8)
  (0,0)..controls(0,1)and(1,1)..(1,0);
  \node[font=\scriptsize,text=black,inner sep=.8pt,outer sep=0pt,above]
  at (0,1.8) {$j$};
  \node[font=\scriptsize,text=black,inner sep=.8pt,outer sep=0pt,below]
  at (0,0) {$i$};
  \end{tikzpicture}\\
 &= %
\begin{tikzpicture}[every path/.style={darkred, thick},scale=0.80,transform shape,baseline={([yshift=-.5ex]current bounding box.center)},x=.75cm,y=.75cm]
  \draw[line width=.65pt,line cap=round,line join=round] (0,0)--(1,1) (0,1)--(1,0);
  \fill(.22,.78)circle(2pt);
  \node[font=\scriptsize,text=black,inner sep=.8pt,outer sep=0pt,below]
  at (0,0) {$i$};
  \node[font=\scriptsize,text=black,inner sep=.8pt,outer sep=0pt,below]
  at (1,0) {$j$};
  \end{tikzpicture}- %
\begin{tikzpicture}[every path/.style={darkred, thick},scale=0.80,transform shape,baseline={([yshift=-.5ex]current bounding box.center)},x=.75cm,y=.75cm]
  \draw[line width=.65pt,line cap=round,line join=round] (0,0)--(1,1) (0,1)--(1,0);
  \fill(.78,.22)circle(2pt);
  \node[font=\scriptsize,text=black,inner sep=.8pt,outer sep=0pt,below]
  at (0,0) {$i$};
  \node[font=\scriptsize,text=black,inner sep=.8pt,outer sep=0pt,below]
  at (1,0) {$j$};
  \end{tikzpicture}.\qquad\text{\normalfont\eqref{idotslide}}
\end{align*}

If $i,j\ne0$, the relation follows from Lemma~\ref{km2c}.
If exactly one color is zero and the other is positive, the
projected black crossing interchanges the two dots. The coefficients
in Definition~\ref{natrual-basic-0i} commute with these dots, and the
shifts defining the red dots are interchanged as well.
This proves the relation in these cases.
The corresponding negative-color cases follow by taking mates
using Definition~\ref{zero-neg_cr} and the dot--bend identities.
For $i=j=0$, use the dot--crossing relation in
Section~\ref{sec:forward-nil-brauer}.
\end{proof}

\begin{Lemma}\label{fw:curl-relation}
For $i\in I$, relation \textnormal{\eqref{icurl}} holds.
\end{Lemma}

\begin{proof}

\begin{equation*}
\begin{tikzpicture}[every path/.style={darkred, thick},scale=0.80,transform shape,baseline={([yshift=-.5ex]current bounding box.center)},x=.65cm,y=.65cm]
  \draw[line width=.65pt,line cap=round,line join=round] (0,0)..controls(0,1.8)and(1,1.8)..(1,.9)
  ..controls(1,.2)and(0,.4)..(0,2);
  \node[font=\scriptsize,text=black,inner sep=.8pt,outer sep=0pt,below]
  at (0,0) {$i$};
  \end{tikzpicture}=\left[ %
\begin{tikzpicture}[anchorbase,scale=1.1,iQ,text=darkred]
\draw[-,draw=darkred] (-0.8,-0.5)\botlabel{i} -- (-0.8,0.5);
\node[cplus] at (-.8,0) {$+$};
\draw[fakebubble,draw=darkred] (-0.4,0) arc(180:-180:0.2);
\node at (-0.48,0) {\strandlabel{i}};
\bubblelabel{-.2,0}{u};
\end{tikzpicture}\right]_{u:-1}.\qquad\text{\normalfont\eqref{icurl}}
\end{equation*}

For $i\ne0$, use Lemma~\ref{km2c}.
For $i=0$, the fake bubble is zero by
\eqref{ifakebubblepolynomial}; the curl vanishes by bending
the cap-annihilation relation of
Section~\ref{sec:forward-nil-brauer} with the zero-color adjunction.
\end{proof}

\begin{Lemma}\label{fw:quadratic-relation}
For $i,j\in I$, relation \textnormal{\eqref{iquadratic}} holds.
\end{Lemma}

\begin{proof}
For $i,j\ne0$, use Lemma~\ref{km2c}.
On colors $(0,k)$ or $(k,0)$ with $k\in I_1$, the projected
black double crossing is $1-(x_1-x_2)^{-2}$.
Multiplying by the two coefficients in
Definition~\ref{natrual-basic-0i} gives
$(y_1-y_2)^{\delta_{k,1}}$ on $(0,k)$ and
$(y_2-y_1)^{\delta_{k,1}}$ on $(k,0)$.
These are the required $Q^\imath$-polynomials.
For $i=j=0$, the square vanishes by
Section~\ref{sec:forward-nil-brauer}, and both terms on the right
vanish by $Q^\imath_{0,0}=0$ and \eqref{ifakebubblepolynomial}.
The cases with one zero and one negative color follow by taking
mates using \eqref{izigzag}--\eqref{ipitchfork}.
For example,
\[ %
\begin{tikzpicture}[scale=0.80,transform shape,
  baseline={([yshift=-.5ex]current bounding box.center)},
  x=.6cm,y=.6cm]
  \begin{scope}[shift={(0,0.32)}]
    \draw[line width=.65pt,line cap=round,line join=round,darkred]
      (0,3.6)--(0,0.32);
    \draw[line width=.65pt,line cap=round,line join=round,darkred]
      (1,0.32) to[out=90,in=-90] (2,1.48);
    \draw[line width=.65pt,line cap=round,line join=round,darkred]
      (2,0) to[out=90,in=-90] (1,1.48);
  \end{scope}
  \begin{scope}[shift={(0,-0.32)}]
    \draw[line width=.65pt,line cap=round,line join=round,darkred]
      (3,3.28)--(3,0);
    \draw[line width=.65pt,line cap=round,line join=round,darkred]
      (2,2.12) to[out=90,in=-90] (1,3.6);
    \draw[line width=.65pt,line cap=round,line join=round,darkred]
      (1,2.12) to[out=90,in=-90] (2,3.28);
  \end{scope}
  \begin{scope}[shift={(0,-0.96)}]
    \draw[line width=.65pt,line cap=round,line join=round,darkred]
      (2,3.92)..controls(2,4.5)and(3,4.5)..(3,3.92);
  \end{scope}
  \begin{scope}[shift={(0,0.96)}]
    \draw[line width=.65pt,line cap=round,line join=round,darkred]
      (0,-0.32)..controls(0,-0.9)and(1,-0.9)..(1,-0.32);
  \end{scope}
  \node[font={\fontsize{7.5pt}{9pt}\selectfont},text=darkred,
    anchor=south,inner sep=.8pt,outer sep=0pt]
    at (0,4.12) {$0$};
  \node[font={\fontsize{7.5pt}{9pt}\selectfont},text=darkred,
    anchor=north,inner sep=.8pt,outer sep=0pt]
    at (3,-0.52) {$k$};
\end{tikzpicture}= %
\begin{tikzpicture}[scale=0.80,transform shape,
  baseline={([yshift=-.5ex]current bounding box.center)},
  x=.6cm,y=.6cm]
  \begin{scope}[shift={(0,-0.32)}]
    \draw[line width=.65pt,line cap=round,line join=round,darkred]
      (0,3.28)--(0,0);
    \draw[line width=.65pt,line cap=round,line join=round,darkred]
      (2,2.12) to[out=90,in=-90] (1,3.28);
    \draw[line width=.65pt,line cap=round,line join=round,darkred]
      (1,2.12) to[out=90,in=-90] (2,3.6);
  \end{scope}
  \begin{scope}[shift={(0,0.32)}]
    \draw[line width=.65pt,line cap=round,line join=round,darkred]
      (3,3.6)--(3,0.32);
    \draw[line width=.65pt,line cap=round,line join=round,darkred]
      (1,0) to[out=90,in=-90] (2,1.48);
    \draw[line width=.65pt,line cap=round,line join=round,darkred]
      (2,0.32) to[out=90,in=-90] (1,1.48);
  \end{scope}
  \begin{scope}[shift={(0,-0.96)}]
    \draw[line width=.65pt,line cap=round,line join=round,darkred]
      (0,3.92)..controls(0,4.5)and(1,4.5)..(1,3.92);
  \end{scope}
  \begin{scope}[shift={(0,0.96)}]
    \draw[line width=.65pt,line cap=round,line join=round,darkred]
      (2,-0.32)..controls(2,-0.9)and(3,-0.9)..(3,-0.32);
  \end{scope}
  \node[font={\fontsize{7.5pt}{9pt}\selectfont},text=darkred,
    anchor=south,inner sep=.8pt,outer sep=0pt]
    at (2,3.48) {$0$};
  \node[font={\fontsize{7.5pt}{9pt}\selectfont},text=darkred,
    anchor=north,inner sep=.8pt,outer sep=0pt]
    at (1,0.12) {$k$};
\end{tikzpicture}= %
\begin{tikzpicture}[scale=0.80,transform shape,
  baseline={([yshift=-.5ex]current bounding box.center)},
  x=.6cm,y=.6cm]
  \begin{scope}[shift={(0,0)}]
    \draw[line width=.65pt,line cap=round,line join=round,darkred]
      (0,3.28)--(0,0);
    \draw[line width=.65pt,line cap=round,line join=round,darkred]
      (3,3.6)--(3,0.32);
    \draw[line width=.65pt,line cap=round,line join=round,darkred]
      (1,3.28)--(1,0);
    \draw[line width=.65pt,line cap=round,line join=round,darkred]
      (2,3.6)--(2,0.32);
  \end{scope}
  \begin{scope}[shift={(0,-0.64)}]
    \draw[line width=.65pt,line cap=round,line join=round,darkred]
      (0,3.92)..controls(0,4.5)and(1,4.5)..(1,3.92);
  \end{scope}
  \begin{scope}[shift={(0,0.64)}]
    \draw[line width=.65pt,line cap=round,line join=round,darkred]
      (2,-0.32)..controls(2,-0.9)and(3,-0.9)..(3,-0.32);
  \end{scope}
  \draw[densely dotted,line cap=round]
    (17.07pt,30.73pt)--(34.14pt,30.73pt);
  \node[font={\fontsize{7.5pt}{9pt}\selectfont},
    anchor=west,inner sep=.8pt,outer sep=0pt]
    at (3.25,1.8) {$Q_{k,0}^{\imath}(x,y)$};
  \node[font={\fontsize{7.5pt}{9pt}\selectfont},text=darkred,
    anchor=south,inner sep=.8pt,outer sep=0pt]
    at (3,3.8) {$0$};
  \node[font={\fontsize{7.5pt}{9pt}\selectfont},text=darkred,
    anchor=north,inner sep=.8pt,outer sep=0pt]
    at (1,-0.2) {$k$};
\end{tikzpicture},\]
for $j=-k\in-I_1$. Thus,
\begin{equation}
\begin{tikzpicture}[baseline=8pt,scale=0.5,every path/.style={darkred, thick},]
		\draw[strand] (0,0) to[out=90,in=-90] (1,1);
		\draw[strand] (1,1) to[out=90,in=-90] (0,2);
		\draw[strand] (1,0) to[out=90,in=-90] (0,1);
		\draw[strand] (0,1) to[out=90,in=-90] (1,2);
	\node[font=\scriptsize,text=black,inner sep=.8pt,outer sep=0pt,below]
  at (0,0) {$0$};
  \node[font=\scriptsize,text=black,inner sep=.8pt,outer sep=0pt,below]
  at (1,0) {$-k$};
    \end{tikzpicture}
 =  %
\begin{tikzpicture}[scale=0.80,transform shape,
  baseline={([yshift=-.5ex]current bounding box.center)},
  x=1cm,y=1cm]
  \begin{scope}[shift={(0,0)}]
    \draw[line width=.65pt,line cap=round,line join=round,darkred]
      (0,0)--(0,1.25);
    \draw[line width=.65pt,line cap=round,line join=round,darkred]
      (.7,0)--(.7,1.25);
  \end{scope}
  \draw[densely dotted,line width=.5pt,darkred,line cap=round]
    (0pt,17.93pt)--(19.92pt,17.93pt);
  \node[font={\fontsize{7.5pt}{9pt}\selectfont},text=darkred,
    anchor=north,inner sep=.8pt,outer sep=0pt]
    at (0,-0.15) {$0$};
  \node[font={\fontsize{7.5pt}{9pt}\selectfont},text=darkred,
    anchor=north,inner sep=.8pt,outer sep=0pt]
    at (0.7,-0.15) {$-k$};
  \node[font={\fontsize{7.5pt}{9pt}\selectfont},text=darkred,
    anchor=west,inner sep=.8pt,outer sep=0pt]
    at (0.87,0.63) {$Q_{k,0}^{\imath}(-y,-x)$};
\end{tikzpicture}
 =  %
\begin{tikzpicture}[scale=0.80,transform shape,
  baseline={([yshift=-.5ex]current bounding box.center)},
  x=1cm,y=1cm]
  \begin{scope}[shift={(0,0)}]
    \draw[line width=.65pt,line cap=round,line join=round,darkred]
      (0,0)--(0,1.25);
    \draw[line width=.65pt,line cap=round,line join=round,darkred]
      (.7,0)--(.7,1.25);
  \end{scope}
  \draw[densely dotted,line width=.5pt,darkred,line cap=round]
    (0pt,17.93pt)--(19.92pt,17.93pt);
  \node[font={\fontsize{7.5pt}{9pt}\selectfont},text=darkred,
    anchor=north,inner sep=.8pt,outer sep=0pt]
    at (0,-0.15) {$0$};
  \node[font={\fontsize{7.5pt}{9pt}\selectfont},text=darkred,
    anchor=north,inner sep=.8pt,outer sep=0pt]
    at (0.7,-0.15) {$-k$};
  \node[font={\fontsize{7.5pt}{9pt}\selectfont},text=darkred,
    anchor=west,inner sep=.8pt,outer sep=0pt]
    at (0.87,0.63) {$Q_{0,-k}^{\imath}(x,y)$};
\end{tikzpicture},
\end{equation}
where
\begin{equation}\label{fw:Q-negative}
 Q_{0,-k}^{\imath}(x,y)=Q_{0,k}^{\imath}(-x,-y)
                     =Q_{k,0}^{\imath}(-y,-x).
\end{equation}
\end{proof}

\begin{Lemma}\label{fw:braid-relation}
For $i,j,k\in I$, the braid relation \eqref{ibraid} holds.
\end{Lemma}

\begin{proof}
Triples with no zero color are covered by Lemma~\ref{km2c},
and $(0,0,0)$ is proved in Section~\ref{sec:nil-brauer-braid}.
For exactly one zero and two positive colors, substitute the
crossings of Definitions~\ref{natrual-basic1}
and~\ref{natrual-basic-0i} into the projected black braid relation.
Taking cyclic mates using Lemma~\ref{fw:basic-relations} gives the
remaining signed cases; the rotation formulas are
\cite[(3.43)--(3.44)]{BWW25}, with $r_{i,j}=1$.
For $a\in I\setminus\{0\}$, the only possibly nonzero corrections are
\[
\begin{array}{c|c}
(a,0,a)&-\operatorname{sgn}(a)\delta_{|a|,1}\,\mathrm{Id}\\
(a,-a,0)&-\operatorname{sgn}(a)\delta_{|a|,1}\,C_L\\
(0,a,-a)& \operatorname{sgn}(a)\delta_{|a|,1}\,C_R,
\end{array}
\]
Here $C_L,C_R$ are the last two cup--cap diagrams in
\eqref{ibraid}, with their displayed scalar coefficients omitted.
For two zeros, Proposition~\ref{prop:forward-mixed-braid}
and its cyclic mates give
\[
\begin{array}{c|c}
(0,0,a)&-\operatorname{sgn}(a)\delta_{|a|,1}\,C_L\\
(0,a,0)& \operatorname{sgn}(a)\delta_{|a|,1}\,\mathrm{Id}\\
(a,0,0)&-\operatorname{sgn}(a)\delta_{|a|,1}\,C_R.
\end{array}
\]
These formulas agree with the divided-difference terms in \eqref{ibraid} and cover all color patterns.
For example, the triangle and cyclicity identities
\eqref{izigzag} and \eqref{ipitchfork}, verified in
Lemma~\ref{fw:basic-relations}, give
\begin{equation}
\begin{tikzpicture}[x=1cm,y=1cm,line cap=round,line join=round,
  baseline={([yshift=-0.5ex]current bounding box.center)}]
  \begin{scope}[scale=1.2]
  \draw[-,thick,draw=darkred] (0.45,.6) to (-0.45,-.6);
  \draw[-,thick,draw=darkred] (0.45,-.6) to (-0.45,.6);
  \draw[-,thick,draw=darkred] (0,-.6) to[out=90,in=-90] (-.45,0);
  \draw[-,thick,draw=darkred] (-0.45,0) to[out=90,in=-90] (0,0.6);
  \end{scope}
  \node[anchor=north,font={\fontsize{7pt}{8pt}\selectfont},text=darkred,inner sep=0.6pt,outer sep=0pt]
    at (-0.54,-0.85) {$0$};
  \node[anchor=north,font={\fontsize{7pt}{8pt}\selectfont},text=darkred,inner sep=0.6pt,outer sep=0pt]
    at (0,-0.85) {$0$};
  \node[anchor=north,font={\fontsize{7pt}{8pt}\selectfont},text=darkred,inner sep=0.6pt,outer sep=0pt]
    at (0.54,-0.85) {$s$};
\end{tikzpicture}
 -  %
\begin{tikzpicture}[x=1cm,y=1cm,line cap=round,line join=round,
  baseline={([yshift=-0.5ex]current bounding box.center)}]
  \begin{scope}[scale=1.2]
  \draw[-,thick,draw=darkred] (-0.45,-.6) to (0.45,.6);
  \draw[-,thick,draw=darkred] (-0.45,.6) to (0.45,-.6);
  \draw[-,thick,draw=darkred] (0,.6) to[out=-90,in=90] (.45,0);
  \draw[-,thick,draw=darkred] (0.45,0) to[out=-90,in=90] (0,-.6);
  \end{scope}

  \node[anchor=north,font={\fontsize{7pt}{8pt}\selectfont},text=darkred,inner sep=0.6pt,outer sep=0pt]
    at (-0.54,-0.85) {$0$};
  \node[anchor=north,font={\fontsize{7pt}{8pt}\selectfont},text=darkred,inner sep=0.6pt,outer sep=0pt]
    at (0,-0.85) {$0$};
  \node[anchor=north,font={\fontsize{7pt}{8pt}\selectfont},text=darkred,inner sep=0.6pt,outer sep=0pt]
    at (0.54,-0.85) {$s$};
\end{tikzpicture}
 =-\delta_{s,1}\,  %
\begin{tikzpicture}[x=1cm,y=1cm,line cap=round,line join=round,
  baseline={([yshift=-0.5ex]current bounding box.center)}]
  \begin{scope}[scale=1.2]
  \draw[-,thick,draw=darkred] (-0.45,.6) to (0.45,-.6);
  \draw[-,thick,draw=darkred] (0,.6) to[out=-90,in=180] (.25,0.2);
  \draw[-,thick,draw=darkred] (0.25,0.2) to[out=0,in=-90] (.45,0.6);
  \draw[-,thick,draw=darkred] (-.45,-.6) to[out=90,in=180] (-.2,-.2);
  \draw[-,thick,draw=darkred] (-.2,-.2) to[out=0,in=90] (0,-0.6);
  \end{scope}

  \node[anchor=north,font={\fontsize{7pt}{8pt}\selectfont},text=darkred,inner sep=0.6pt,outer sep=0pt]
    at (0.54,-0.85) {$s$};

  \node[anchor=north,font={\fontsize{7pt}{8pt}\selectfont},text=darkred,inner sep=0.6pt,outer sep=0pt]
    at (-0.54,-0.85) {$0$};
  \node[anchor=north,font={\fontsize{7pt}{8pt}\selectfont},text=darkred,inner sep=0.6pt,outer sep=0pt]
    at (0,-0.85) {$0$};
\end{tikzpicture}
 \qquad(s\in I_1).
\end{equation}
Taking its mate using the same identities gives
\begin{equation}\label{fw:negative-braid-last}
\begin{tikzpicture}[x=1cm,y=1cm,line cap=round,line join=round,
  baseline={([yshift=-0.5ex]current bounding box.center)}]
  \begin{scope}[scale=1.2]
  \draw[-,thick,draw=darkred] (0.45,.6) to (-0.45,-.6);
  \draw[-,thick,draw=darkred] (0.45,-.6) to (-0.45,.6);
  \draw[-,thick,draw=darkred] (0,-.6) to[out=90,in=-90] (-.45,0);
  \draw[-,thick,draw=darkred] (-0.45,0) to[out=90,in=-90] (0,0.6);
  \end{scope}

  \node[anchor=north,font={\fontsize{7pt}{8pt}\selectfont},text=darkred,inner sep=0.6pt,outer sep=0pt]
    at (-0.54,-0.85) {$0$};
  \node[anchor=north,font={\fontsize{7pt}{8pt}\selectfont},text=darkred,inner sep=0.6pt,outer sep=0pt]
    at (0,-0.85) {$-s$};
  \node[anchor=north,font={\fontsize{7pt}{8pt}\selectfont},text=darkred,inner sep=0.6pt,outer sep=0pt]
    at (0.54,-0.85) {$0$};
\end{tikzpicture}
 -  %
\begin{tikzpicture}[x=1cm,y=1cm,line cap=round,line join=round,
  baseline={([yshift=-0.5ex]current bounding box.center)}]
  \begin{scope}[scale=1.2]
  \draw[-,thick,draw=darkred] (-0.45,-.6) to (0.45,.6);
  \draw[-,thick,draw=darkred] (-0.45,.6) to (0.45,-.6);
  \draw[-,thick,draw=darkred] (0,.6) to[out=-90,in=90] (.45,0);
  \draw[-,thick,draw=darkred] (0.45,0) to[out=-90,in=90] (0,-.6);
  \end{scope}

  \node[anchor=north,font={\fontsize{7pt}{8pt}\selectfont},text=darkred,inner sep=0.6pt,outer sep=0pt]
    at (-0.54,-0.85) {$0$};
  \node[anchor=north,font={\fontsize{7pt}{8pt}\selectfont},text=darkred,inner sep=0.6pt,outer sep=0pt]
    at (0,-0.85) {$-s$};
  \node[anchor=north,font={\fontsize{7pt}{8pt}\selectfont},text=darkred,inner sep=0.6pt,outer sep=0pt]
    at (0.54,-0.85) {$0$};
\end{tikzpicture}
 =-\delta_{s,1}\,  %
\begin{tikzpicture}[x=1cm,y=1cm,line cap=round,line join=round,
  baseline={([yshift=-0.5ex]current bounding box.center)}]
  \begin{scope}[scale=1.2]
  \draw[-,thick,draw=darkred] (0.45,.6) to (0.45,-.6);
  \draw[-,thick,draw=darkred] (-0.45,.6) to (-0.45,-.6);
  \draw[-,thick,draw=darkred] (0,-.6) to (0,0.6);
  \end{scope}
  \node[anchor=north,font={\fontsize{7pt}{8pt}\selectfont},text=darkred,inner sep=0.6pt,outer sep=0pt]
    at (-0.54,-0.85) {$0$};

  \node[anchor=north,font={\fontsize{7pt}{8pt}\selectfont},text=darkred,inner sep=0.6pt,outer sep=0pt]
    at (0,-0.85) {$-s$};

  \node[anchor=north,font={\fontsize{7pt}{8pt}\selectfont},text=darkred,inner sep=0.6pt,outer sep=0pt]
    at (0.54,-0.85) {$0$};

\end{tikzpicture}.
\end{equation}
By our parameter convention,
\begin{equation}\label{fw:Q-last-source}
 Q_{0,-s}^{\imath}(x,y)
       =Q_{s,0}^{\imath}(-y,-x)
       =1-\delta_{s,1}+\delta_{s,1}(y-x)
       \qquad(s\in I_1),
\end{equation}
and hence
\begin{equation}\label{fw:Q-divdiff}
 \frac{Q_{0,-s}^{\imath}(x,y)-Q_{0,-s}^{\imath}(z,y)}{x-z}
 =-\delta_{s,1}.
\end{equation}
\end{proof}

The preceding calculations give the integral counterpart of
\cite[Theorem~6.16]{GRS}.

\begin{Theorem}
\label{thm:integral-forward-action}
\label{fw:red-to-black}
Suppose that a locally Schurian $\mathbb C$-linear category
$\mathcal M$ admits an affine Brauer categorification and that the
spectrum of the dot endomorphism on its generating endofunctor is
exactly $\mathbb Z$.  Then the following statements hold.
\begin{itemize}
\item[(1)] The weight decomposition of Lemma~\ref{weight-pro}, together
  with the $1$- and $2$-morphisms constructed in
  Definitions~\ref{natrual-basic}, \ref{natrual-basic1},
  \ref{eq:red-cup}, \ref{natrual-basic-0i}, and
  \ref{zero-neg_cr}, after the conversion in
  Subsection~\ref{subsec:forward-normalization}, makes $\mathcal M$ a generalized
  nilpotent $2$-representation of the integral quasi-split type
  $\operatorname{AIII}$ $2$-category $\mathfrak U^\imath$.
\item[(2)] The full subcategories
  $\mathcal M_\lambda^{\mathrm{fg}}$ of finitely generated objects form
  a nilpotent $2$-subrepresentation of $\mathfrak U^\imath$.
\item[(3)] If $\mathcal M$ is locally finite Abelian, then the resulting
  $2$-representation on $\mathcal M$ is nilpotent.
\end{itemize}
\end{Theorem}

\begin{proof}
Lemma~\ref{weight-pro} supplies the decomposition
$\mathcal M=\prod_{\lambda\in X^\imath}\mathcal M_\lambda$ and shows
that $E_i$ has the required weight shift. Lemma~\ref{km2c} gives the
nonzero-color relations. Lemmas~\ref{fw:positive-pitchfork} and
\ref{fw:nonzero-pitchfork} give the mixed cases of \eqref{ipitchfork}, while
the mixed braid relation~\eqref{fw:mixed-braid} gives the positive mixed braid.
Section~\ref{sec:forward-nil-brauer} gives the zero relations other than
the braid, and Section~\ref{sec:nil-brauer-braid} gives the all-zero
braid. Lemmas~\ref{fw:basic-relations}--\ref{fw:braid-relation} then
verify \eqref{izigzag}--\eqref{ibraid}.

It remains to check the bubble relations
\eqref{ibubblerel}--\eqref{ibubslide}. Write $H_{i,V}(u)$ for the
series in \eqref{ibubblegeneratingfunction}; thus $H_{0,V}(u)$ is the
ordinary zero-color bubble series minus $(2u)^{-1}\operatorname{Id}_V$.
Lemma~\ref{km2c} gives the nonzero-color bubble products and the
slides past nonzero strands.
For a nonzero bubble beside a zero strand with dot $y$,
the primary-factor comparison in Lemma~\ref{ou}(2) and
\eqref{eq:forward-bubble-normalization} give
\[
 \frac{H_{i,E_0V}(u)}{E_0H_{i,V}(u)}=
 \begin{cases}
  \dfrac{u+y}{u-y},& |i|=1,\\
  1,& |i|>1,
 \end{cases}
\]
which equals $Q^{\imath}_{-i,0}(-u,y)/Q^{\imath}_{i,0}(u,y)$.
For the zero-color series, put
$\epsilon=\varepsilon_0(V)$,
$\kappa=(-1)^{c_\lambda+1}/2$ and $C_V=A_VB_V$.
Use polynomial representatives with truncation order $K$.
The factorization in Lemma~\ref{nbone:lemma:R4} gives
\[
 C_V(v)\frac{O_V(v)}v
 =\kappa\frac{t_{V,0}(v)}{v^{\epsilon+1}}
   +v^{K-\epsilon-1}P_V(v)
\]
for a polynomial $P_V$. If $\beta_{r,V}$ is the red bubble with
$r\ge0$ dots, the calculation of Lemma~\ref{nbone:lemma:nil2} gives
\[
 \beta_{r,V}
 =\left[v^r C_V(v)
   \left(\frac{O_V(v)}v-1+\frac1{2v}\right)\right]_{v:-1}
 =\kappa\left[\frac{t_{V,0}(v)}{v^{\epsilon+1}}\right]_{v:-r-1}
   +\frac{\delta_{r,0}}2\operatorname{Id}_V.
\]
Here $K>\epsilon$ makes the tail polynomial, and $C_V(0)=1$.
Thus
\[
 H_{0,V}(u)=\kappa\frac{t_{V,0}(u)}{u^{\epsilon+1}}.
\]
Since $t_{V,0}$ is monic of degree $\epsilon$, this series has
leading term $\kappa u^{-1}$, as required by \eqref{ibubblerel}.
Uniqueness of the zero-primary factors in Lemma~\ref{ou}(1)
gives $t_{V,0}(u)t_{V,0}(-u)=(-1)^\epsilon u^{2\epsilon}$,
so
\[
 H_{0,V}(u)H_{0,V}(-u)
 =-\frac1{4u^2}\operatorname{Id}_V.
\]
Extracting the zero-primary factors in Lemma~\ref{ou}(2) on
$E_0V$, whose dot is $y$, also gives
\[
 \frac{H_{0,E_0V}(u)}{E_0H_{0,V}(u)}
 =\left(\frac{u-y}{u+y}\right)^2
 =\frac{R_{0,0}(-u,y)}{R_{0,0}(u,y)}.
\]
These are the zero-color bubble product and slide relations.
Let $j\ne0$ and $E_jV\ne0$. Extracting the factors at zero in
Lemma~\ref{ou}(2) and using Lemma~\ref{fw:parity-negative}, we obtain
\[
 \frac{H_{0,E_jV}(u)}{E_jH_{0,V}(u)}=
 \begin{cases}
  -\dfrac{u+y}{u-y},& |j|=1,\\
  1,& |j|>1,
 \end{cases}
 =\frac{R_{0,j}(-u,y)}{R_{0,j}(u,y)}.
\]
Here $y=x-j$ is the dot on the unoriented $j$-strand.
This proves the remaining bubble slide. Naturality and filtered
colimits extend the identities to all objects.

Hence the assignments above define the asserted $2$-representation.

Let $V$ be finitely generated.  By Lemma~\ref{sweet1}, $EV$ is
finitely generated and $\operatorname{End}_{\mathcal M}(EV)$ is
finite-dimensional.  Consequently the affine Brauer dot on $EV$ has
only finitely many generalized eigenvalues, so $E_iV=0$ for all but
finitely many $i\in\mathbb Z$.
After the conversion in
\eqref{eq:forward-dot-normalization}, the red dot on $E_iV$ is
$X_i=x-i$ for every $i\in I$, including $X_0=x$.
It is therefore nilpotent.
This proves generalized nilpotence.

By Lemma~\ref{sweet1}, $EV$ is finitely generated; its generalized
eigenspace summand $E_iV$ is therefore finitely generated as well.
The constructed natural transformations restrict to these summands.
Thus the families
$\mathcal M_\lambda^{\mathrm{fg}}$ form a nilpotent
$2$-subrepresentation, proving~\textup{(2)}.  If $\mathcal M$ is
locally finite Abelian, every object has finite length and therefore is
finitely generated.  The same argument then applies to every object,
which proves~\textup{(3)}.
\end{proof}

\section{Module category over the nil-Brauer category}
\label{sec:forward-nil-brauer}
Fix a weight $\lambda$. We verify relations~\textnormal{(1)}--\textnormal{(7)}
of Definition~\ref{nil} with $t=c_\lambda$; relation~\textnormal{(8)}
is proved in Section~\ref{sec:nil-brauer-braid}.
The two zero-color triangle identities proved immediately before
Lemma~\ref{fw:nonzero-pitchfork} give relation~\textnormal{(1)}.
The black dot--cap relation gives relation~\textnormal{(6)}, since
the normalizing factors commute with the dot.

By naturality and preservation of filtered colimits, it suffices to
check the remaining relations on finitely generated objects.
For a nonzero such object $V\in\mathcal M_\lambda$, write
$c=c_\lambda$ and $e=\varepsilon_0(V)$. We use $K$ for a sufficiently large
common B\'ezout truncation order.
For coefficient extractions at infinity, $A_V$ and $B_V$ denote
polynomial representatives: retain the factors $m^\pm_{V,0}f_V^\pm$
and replace $\Gamma_V^{\pm1}$ by their truncations modulo $u^K$.
All divided differences are formed before nilpotent evaluation, as in
Lemma~\ref{lem:finite-dd-evaluation}.

\begin{Lemma}\label{nbone:lemma:R3}
The generators satisfy relation~\textnormal{(5)} of Definition~\ref{nil}.
\end{Lemma}

\begin{proof}
Set $d=(x_2-x_1)/2$. Substituting the definition of the red crossing and applying the dot relation in the affine Brauer category gives
\begin{align*}
\mathrm{LHS}
&=
\begin{tikzpicture}[x=1bp,y=1bp,scale=.92,
  baseline={([yshift=-.5ex]current bounding box.center)},
  line cap=round,line join=round,
  every node/.style={inner sep=.8pt,outer sep=0pt,font=\fontsize{8}{9}\selectfont}]
\draw[black,line width=.65bp] (0,0) -- (17.52,26.28);
\draw[black,line width=.65bp] (17.52,0) -- (0,26.28);
\draw[black,densely dotted,line width=.55bp] (2.555,22.447) -- (14.965,22.447);
\draw[black,densely dotted,line width=.55bp] (2.555,3.832) -- (14.965,3.832);
\fill[black] (4.783,19.106) circle[radius=1.15bp];
\node[anchor=south] at (0,27.594) {$0$};
\node[anchor=north] at (0,-1.314) {$0$};
\node[anchor=south] at (17.52,27.594) {$0$};
\node[anchor=north] at (17.52,-1.314) {$0$};
\node[anchor=west,font=\fontsize{7.5}{7.5}\selectfont] at (21.02,22.447) {$\textstyle B_V(d)$};
\node[anchor=west,font=\fontsize{7.5}{7.5}\selectfont] at (21.02,3.832) {$\textstyle A_V(d)$};
\end{tikzpicture}
-
\begin{tikzpicture}[x=1bp,y=1bp,scale=.92,
  baseline={([yshift=-.5ex]current bounding box.center)},
  line cap=round,line join=round,
  every node/.style={inner sep=.8pt,outer sep=0pt,font=\fontsize{8}{9}\selectfont}]
\draw[black,line width=.65bp] (0,0) -- (17.52,26.28);
\draw[black,line width=.65bp] (17.52,0) -- (0,26.28);
\draw[black,densely dotted,line width=.55bp] (2.555,22.447) -- (14.965,22.447);
\draw[black,densely dotted,line width=.55bp] (2.555,3.832) -- (14.965,3.832);
\fill[black] (12.737,7.174) circle[radius=1.15bp];
\node[anchor=south] at (0,27.594) {$0$};
\node[anchor=north] at (0,-1.314) {$0$};
\node[anchor=south] at (17.52,27.594) {$0$};
\node[anchor=north] at (17.52,-1.314) {$0$};
\node[anchor=west,font=\fontsize{7.5}{7.5}\selectfont] at (21.02,22.447) {$\textstyle B_V(d)$};
\node[anchor=west,font=\fontsize{7.5}{7.5}\selectfont] at (21.02,3.832) {$\textstyle A_V(d)$};
\end{tikzpicture}
+
\begin{tikzpicture}[x=1bp,y=1bp,scale=.92,
  baseline={([yshift=-.5ex]current bounding box.center)},
  line cap=round,line join=round,
  every node/.style={inner sep=.8pt,outer sep=0pt,font=\fontsize{8}{9}\selectfont}]
\draw[black,line width=.65bp] (2.217,4.319) -- (2.217,27.166);
\draw[black,line width=.65bp] (19.792,4.319) -- (19.792,27.166);
\draw[black,densely dotted,line width=.55bp] (2.217,15.742) -- (19.792,15.742);
\node[anchor=west,font=\fontsize{7.5}{7.5}\selectfont] at (23.292,15.742) {$\textstyle g_V(x_1,x_2)(x_1-x_2)$};
\node[anchor=north] at (2.217,2.819) {$0$};
\node[anchor=north] at (19.792,2.819) {$0$};
\node[anchor=south] at (2.217,28.666) {$0$};
\node[anchor=south] at (19.792,28.666) {$0$};
\end{tikzpicture}\\
&=
\begin{tikzpicture}[x=1bp,y=1bp,scale=.92,
  baseline={([yshift=-.5ex]current bounding box.center)},
  line cap=round,line join=round,
  every node/.style={inner sep=.8pt,outer sep=0pt,font=\fontsize{8}{9}\selectfont}]
\draw[black,line width=.65bp] (2.217,4.319) -- (2.217,27.166);
\draw[black,line width=.65bp] (19.792,4.319) -- (19.792,27.166);
\draw[black,densely dotted,line width=.55bp] (2.217,15.742) -- (19.792,15.742);
\node[anchor=west,font=\fontsize{7.5}{7.5}\selectfont] at (23.292,15.742) {$\textstyle B_V(d)A_V(d)$};
\node[anchor=north] at (2.217,2.819) {$0$};
\node[anchor=north] at (19.792,2.819) {$0$};
\node[anchor=south] at (2.217,28.666) {$0$};
\node[anchor=south] at (19.792,28.666) {$0$};
\end{tikzpicture}
-
\begin{tikzpicture}[x=1bp,y=1bp,scale=.92,
  baseline={([yshift=-.5ex]current bounding box.center)},
  line cap=round,line join=round,
  every node/.style={inner sep=.8pt,outer sep=0pt,font=\fontsize{8}{9}\selectfont}]
\draw[black,line width=.65bp] (2.217,34.197) .. controls (2.217,16.798) and (19.792,16.798) .. (19.792,34.197);
\draw[black,line width=.65bp] (2.217,4.32) .. controls (2.217,21.719) and (19.792,21.719) .. (19.792,4.32);
\fill[black] (18.441,26.832) circle[radius=1.15bp];
\fill[black] (18.441,11.685) circle[radius=1.15bp];
\node[anchor=west,font=\fontsize{7.5}{7.5}\selectfont] at (23.197,26.832) {$\textstyle B_V(x)$};
\node[anchor=west,font=\fontsize{7.5}{7.5}\selectfont] at (23.197,11.685) {$\textstyle A_V(x)$};
\node[anchor=north] at (2.217,2.82) {$0$};
\node[anchor=north] at (19.792,2.82) {$0$};
\node[anchor=south] at (2.217,35.697) {$0$};
\node[anchor=south] at (19.792,35.697) {$0$};
\end{tikzpicture}
+
\begin{tikzpicture}[x=1bp,y=1bp,scale=.92,
  baseline={([yshift=-.5ex]current bounding box.center)},
  line cap=round,line join=round,
  every node/.style={inner sep=.8pt,outer sep=0pt,font=\fontsize{8}{9}\selectfont}]
\draw[black,line width=.65bp] (2.217,4.319) -- (2.217,27.166);
\draw[black,line width=.65bp] (19.792,4.319) -- (19.792,27.166);
\draw[black,densely dotted,line width=.55bp] (2.217,15.742) -- (19.792,15.742);
\node[anchor=west,font=\fontsize{7.5}{7.5}\selectfont] at (23.292,15.742) {$\textstyle g_V(x_1,x_2)(x_1-x_2)$};
\node[anchor=north] at (2.217,2.819) {$0$};
\node[anchor=north] at (19.792,2.819) {$0$};
\node[anchor=south] at (2.217,28.666) {$0$};
\node[anchor=south] at (19.792,28.666) {$0$};
\end{tikzpicture}.
\end{align*}

By the definition $g_V(x_1,x_2)=\frac{1-A_V(d)B_V(d)}{x_1-x_2}$, the coefficient of the straight-strand term is $1$. Moreover, {by Definition~\ref{eq:red-cup} and Lemma~\ref{eq:red-cup-E0V}}, the black cup decorated by $B_V$ and the black cap decorated by $A_V$ are precisely the red cup and cap. Hence the last expression equals
\[
\begin{tikzpicture}[x=1bp,y=1bp,scale=.92,
  baseline={([yshift=-.5ex]current bounding box.center)},
  line cap=round,line join=round,
  every node/.style={inner sep=.8pt,outer sep=0pt,font=\fontsize{8}{9}\selectfont}]
\draw[darkred,line width=.65bp] (2.217,4.319) -- (2.217,27.166);
\draw[darkred,line width=.65bp] (19.792,4.319) -- (19.792,27.166);
\node[anchor=north] at (2.217,2.819) {$0$};
\node[anchor=north] at (19.792,2.819) {$0$};
\node[anchor=south] at (2.217,28.666) {$0$};
\node[anchor=south] at (19.792,28.666) {$0$};
\end{tikzpicture}
-
\begin{tikzpicture}[x=1bp,y=1bp,scale=.92,
  baseline={([yshift=-.5ex]current bounding box.center)},
  line cap=round,line join=round,
  every node/.style={inner sep=.8pt,outer sep=0pt,font=\fontsize{8}{9}\selectfont}]
\draw[darkred,line width=.65bp] (0,32.85) .. controls (0,14.235) and (17.52,14.235) .. (17.52,32.85);
\draw[darkred,line width=.65bp] (0,0) .. controls (0,18.615) and (17.52,18.615) .. (17.52,0);
\node[anchor=south] at (0,34.35) {$0$};
\node[anchor=south] at (17.52,34.35) {$0$};
\node[anchor=north] at (0,-1.5) {$0$};
\node[anchor=north] at (17.52,-1.5) {$0$};
\end{tikzpicture}.
\]
\end{proof}

\begin{Lemma}\label{nbone:lemma:R4}
The generators satisfy relation~\textnormal{(3)} of Definition~\ref{nil}.
\end{Lemma}

\begin{proof}
Set \(d=(x_2-x_1)/2\), and write
\[
C_V(u):=A_V(u)B_V(u)
=\frac{(-1)^{c+1}}2m_{V,0}(u)f_V^+(u)f_V^-(u)
\bigl(1+u^K R_V(u)\bigr),
\]
where $R_V$ is a polynomial arising from the truncations of
$\Gamma_V$ and $\Gamma_V^{-1}$.
Since $C_V(0)=1$, set
\[
P_C(x):=\frac{C_V(x)-1}{2x}\in Z_{\lambda,V}[x].
\]

By the definitions of the red cap and the red crossing,
\begingroup
\setlength{\jot}{2pt}
\begin{align*}
\begin{tikzpicture}[x=1bp,y=1bp,scale=.80,
  baseline={([yshift=-.5ex]current bounding box.center)},
  line cap=round,line join=round,
  every node/.style={inner sep=.5pt,outer sep=0pt,font=\fontsize{7.2}{8}\selectfont}]
\draw[darkred,line width=.65bp] (2.217,4.319) .. controls (8.075,11.349) and (19.792,16.075) .. (19.792,25.409);
\draw[darkred,line width=.65bp] (19.792,4.319) .. controls (13.934,11.349) and (2.217,16.075) .. (2.217,25.409);
\draw[darkred,line width=.65bp] (2.217,25.409) .. controls (2.217,43.863) and (19.792,43.863) .. (19.792,25.409);
\node[text=darkred,anchor=north] at (2.217,2.319) {$0$};
\draw[darkg,line width=.65bp] (26.792,4.319) -- (26.792,26.116);
\node[text=darkg,anchor=north] at (26.792,2.519) {$V$};
\end{tikzpicture}
&=\begin{tikzpicture}[x=1bp,y=1bp,scale=.80,
  baseline={([yshift=-.5ex]current bounding box.center)},
  line cap=round,line join=round,
  every node/.style={inner sep=.5pt,outer sep=0pt,font=\fontsize{7.2}{8}\selectfont}]
\draw[black,line width=.65bp] (2.217,4.319) .. controls (8.075,11.349) and (19.792,16.075) .. (19.792,25.409);
\draw[black,line width=.65bp] (19.792,25.409) .. controls (19.792,42.984) and (2.217,42.984) .. (2.217,25.409);
\draw[black,line width=.65bp] (2.217,25.409) .. controls (2.217,16.075) and (13.934,11.349) .. (19.792,4.319);
\draw[black,densely dotted,line width=.55bp] (3.146,21.015) -- (18.863,21.015);
\draw[black,densely dotted,line width=.55bp] (6.765,8.713) -- (15.244,8.713);
\fill[black] (19.211,30.53) circle[radius=1.15bp];
\node[anchor=west,font=\fontsize{6.9}{7.4}\selectfont] at (23.292,30.53) {$\textstyle A_V(x)$};
\node[anchor=north] at (1.257,3.167) {$0$};
\node[anchor=north] at (20.752,3.167) {$0$};
\node[anchor=west,font=\fontsize{6.9}{7.4}\selectfont] at (23.231,19.03) {$\textstyle B_V(d)$};
\node[anchor=west,font=\fontsize{6.9}{7.4}\selectfont] at (23.292,7.53) {$\textstyle A_V(d)$};
\end{tikzpicture}\;+\;\begin{tikzpicture}[x=1bp,y=1bp,scale=.80,
  baseline={([yshift=-.5ex]current bounding box.center)},
  line cap=round,line join=round,
  every node/.style={inner sep=.5pt,outer sep=0pt,font=\fontsize{7.2}{8}\selectfont}]
\draw[black,line width=.65bp] (2.217,4.319) -- (2.217,10.47);
\draw[black,line width=.65bp] (2.217,10.47) .. controls (2.217,25.409) and (19.792,25.409) .. (19.792,10.47);
\draw[black,line width=.65bp] (19.792,10.47) -- (19.792,4.319);
\draw[black,densely dotted,line width=.55bp] (2.217,6.955) -- (19.792,6.955);
\fill[black] (19.706,12.232) circle[radius=1.15bp];
\node[anchor=west,font=\fontsize{6.9}{7.4}\selectfont] at (23.292,12.232) {$\textstyle A_V(x)$};
\node[anchor=north] at (2.217,2.819) {$0$};
\node[anchor=north] at (19.792,2.819) {$0$};
\node[anchor=west,font=\fontsize{6.9}{7.4}\selectfont] at (23.292,0.732) {$\textstyle g_V(x_1,x_2)$};
\end{tikzpicture}\\
&=\begin{tikzpicture}[x=1bp,y=1bp,scale=.80,
  baseline={([yshift=-.5ex]current bounding box.center)},
  line cap=round,line join=round,
  every node/.style={inner sep=.5pt,outer sep=0pt,font=\fontsize{7.2}{8}\selectfont}]
\draw[black,line width=.65bp] (2.217,4.319) .. controls (8.075,11.349) and (19.792,16.075) .. (19.792,25.409);
\draw[black,line width=.65bp] (19.792,25.409) .. controls (19.792,42.984) and (2.217,42.984) .. (2.217,25.409);
\draw[black,line width=.65bp] (2.217,25.409) .. controls (2.217,16.075) and (13.934,11.349) .. (19.792,4.319);
\draw[black,densely dotted,line width=.55bp] (6.765,8.713) -- (15.244,8.713);
\fill[black] (19.211,30.53) circle[radius=1.15bp];
\node[anchor=north] at (1.257,3.167) {$0$};
\node[anchor=north] at (20.752,3.167) {$0$};
\node[anchor=west,font=\fontsize{6.9}{7.4}\selectfont] at (23.292,30.53) {$\scriptstyle C_V(x)$};
\node[anchor=west,font=\fontsize{6.9}{7.4}\selectfont] at (23.292,8.713) {$\textstyle A_V(d)$};
\end{tikzpicture}\;+\;\begin{tikzpicture}[x=1bp,y=1bp,scale=.80,
  baseline={([yshift=-.5ex]current bounding box.center)},
  line cap=round,line join=round,
  every node/.style={inner sep=.5pt,outer sep=0pt,font=\fontsize{7.2}{8}\selectfont}]
\draw[black,line width=.65bp] (2.217,4.319) -- (2.217,10.47);
\draw[black,line width=.65bp] (2.217,10.47) .. controls (2.217,25.409) and (19.792,25.409) .. (19.792,10.47);
\draw[black,line width=.65bp] (19.792,10.47) -- (19.792,4.319);
\draw[black,densely dotted,line width=.55bp] (2.217,6.955) -- (19.792,6.955);
\fill[black] (19.706,12.232) circle[radius=1.15bp];
\node[anchor=west,font=\fontsize{6.9}{7.4}\selectfont] at (23.292,12.232) {$\textstyle A_V(x)$};
\node[anchor=north] at (2.217,2.819) {$0$};
\node[anchor=north] at (19.792,2.819) {$0$};
\node[anchor=west,font=\fontsize{6.9}{7.4}\selectfont] at (23.292,0.732) {$\textstyle g_V(x_1,x_2)$};
\end{tikzpicture}\\
&=\left[\begin{tikzpicture}[x=1bp,y=1bp,scale=.80,
  baseline={([yshift=-.5ex]current bounding box.center)},
  line cap=round,line join=round,
  every node/.style={inner sep=.5pt,outer sep=0pt,font=\fontsize{7.2}{8}\selectfont}]
\draw[black,line width=.65bp] (2.217,4.319) .. controls (8.075,11.349) and (19.792,16.075) .. (19.792,25.409);
\draw[black,line width=.65bp] (19.792,25.409) .. controls (19.792,42.984) and (2.217,42.984) .. (2.217,25.409);
\draw[black,line width=.65bp] (2.217,25.409) .. controls (2.217,16.075) and (13.934,11.349) .. (19.792,4.319);
\draw[black,densely dotted,line width=.55bp] (6.765,8.713) -- (15.244,8.713);
\fill[black] (19.211,30.53) circle[radius=1.15bp];
\node[anchor=north] at (1.257,3.167) {$0$};
\node[anchor=north] at (20.752,3.167) {$0$};
\node[anchor=west,font=\fontsize{6.9}{7.4}\selectfont] at (23.292,30.53) {$\textstyle \frac1{u-x}$};
\node[anchor=west,font=\fontsize{6.9}{7.4}\selectfont] at (23.292,8.713) {$\textstyle A_V(d)$};
\end{tikzpicture}C_V(u)\right]_{u:-1}
   \;+\;\begin{tikzpicture}[x=1bp,y=1bp,scale=.80,
  baseline={([yshift=-.5ex]current bounding box.center)},
  line cap=round,line join=round,
  every node/.style={inner sep=.5pt,outer sep=0pt,font=\fontsize{7.2}{8}\selectfont}]
\draw[black,line width=.65bp] (2.217,4.319) -- (2.217,10.47);
\draw[black,line width=.65bp] (2.217,10.47) .. controls (2.217,25.409) and (19.792,25.409) .. (19.792,10.47);
\draw[black,line width=.65bp] (19.792,10.47) -- (19.792,4.319);
\draw[black,densely dotted,line width=.55bp] (2.217,6.955) -- (19.792,6.955);
\fill[black] (19.706,12.232) circle[radius=1.15bp];
\node[anchor=west,font=\fontsize{6.9}{7.4}\selectfont] at (23.292,12.232) {$\textstyle A_V(x)$};
\node[anchor=north] at (2.217,2.819) {$0$};
\node[anchor=north] at (19.792,2.819) {$0$};
\node[anchor=west,font=\fontsize{6.9}{7.4}\selectfont] at (23.292,0.732) {$\textstyle g_V(x_1,x_2)$};
\end{tikzpicture}.
\end{align*}
\endgroup

The retained factors $m^\pm_{V,0}f_V^\pm$ in $C_V$ annihilate
the nonzero primary components, so the adjacent zero-color projection
may be omitted. Applying the
affine Brauer cap--resolvent identity gives
\begingroup
\setlength{\jot}{2pt}
\begin{align*}
\begin{tikzpicture}[x=1bp,y=1bp,scale=.80,
  baseline={([yshift=-.5ex]current bounding box.center)},
  line cap=round,line join=round,
  every node/.style={inner sep=.5pt,outer sep=0pt,font=\fontsize{7.2}{8}\selectfont}]
\draw[darkred,line width=.65bp] (2.217,4.319) .. controls (8.075,11.349) and (19.792,16.075) .. (19.792,25.409);
\draw[darkred,line width=.65bp] (19.792,4.319) .. controls (13.934,11.349) and (2.217,16.075) .. (2.217,25.409);
\draw[darkred,line width=.65bp] (2.217,25.409) .. controls (2.217,43.863) and (19.792,43.863) .. (19.792,25.409);
\node[text=darkred,anchor=north] at (2.217,2.319) {$0$};
\draw[darkg,line width=.65bp] (26.792,4.319) -- (26.792,26.116);
\node[text=darkg,anchor=north] at (26.792,2.519) {$V$};
\end{tikzpicture}
&=\left[
\begin{tikzpicture}[x=1bp,y=1bp,scale=.80,
  baseline={([yshift=-.5ex]current bounding box.center)},
  line cap=round,line join=round,
  every node/.style={inner sep=.5pt,outer sep=0pt,font=\fontsize{7.2}{8}\selectfont}]
\draw[black,line width=.65bp] (2.217,4.319) -- (2.217,9.592);
\draw[black,line width=.65bp] (2.217,9.592) .. controls (2.217,25.409) and (19.792,25.409) .. (19.792,9.592);
\draw[black,line width=.65bp] (19.792,9.592) -- (19.792,4.319);
\draw[black,fill=white,line width=.55bp] (2.222,10.069) circle[radius=2.834bp];
\draw[black,line width=.55bp] (0.664,10.069) -- (3.781,10.069);
\fill[black] (19.792,6.956) circle[radius=1.15bp];
\node[anchor=west,font=\fontsize{6.9}{7.4}\selectfont] at (23.292,6.956) {$\textstyle A_V(x)$};
\node[anchor=north] at (2.217,2.819) {$0$};
\node[anchor=north] at (19.792,2.819) {$0$};
\end{tikzpicture}
C_V(u)\frac{O_V(u)}u
\right]_{u:-1}
-
\left[\begin{tikzpicture}[x=1bp,y=1bp,scale=.80,
  baseline={([yshift=-.5ex]current bounding box.center)},
  line cap=round,line join=round,
  every node/.style={inner sep=.5pt,outer sep=0pt,font=\fontsize{7.2}{8}\selectfont}]
\draw[black,line width=.65bp] (2.217,4.319) -- (2.217,9.592);
\draw[black,line width=.65bp] (2.217,9.592) .. controls (2.217,25.409) and (19.792,25.409) .. (19.792,9.592);
\draw[black,line width=.65bp] (19.792,9.592) -- (19.792,4.319);
\draw[black,fill=white,line width=.55bp] (19.738,11.074) circle[radius=2.834bp];
\draw[black,line width=.55bp] (18.18,11.074) -- (21.297,11.074);
\fill[black] (19.792,5.59) circle[radius=1.15bp];
\node[anchor=west,font=\fontsize{6.9}{7.4}\selectfont] at (24.572,5.59) {$\textstyle A_V(x)$};
\node[anchor=north] at (2.217,2.819) {$0$};
\node[anchor=north] at (19.792,2.819) {$0$};
\end{tikzpicture}\frac{C_V(u)}{2u}\right]_{u:-1}
+
\;
\begin{tikzpicture}[x=1bp,y=1bp,scale=.80,
  baseline={([yshift=-.5ex]current bounding box.center)},
  line cap=round,line join=round,
  every node/.style={inner sep=.5pt,outer sep=0pt,font=\fontsize{7.2}{8}\selectfont}]
\draw[black,line width=.65bp] (2.217,4.319) -- (2.217,10.47);
\draw[black,line width=.65bp] (2.217,10.47) .. controls (2.217,25.409) and (19.792,25.409) .. (19.792,10.47);
\draw[black,line width=.65bp] (19.792,10.47) -- (19.792,4.319);
\draw[black,densely dotted,line width=.55bp] (2.217,6.955) -- (19.792,6.955);
\fill[black] (19.706,12.232) circle[radius=1.15bp];
\node[anchor=west,font=\fontsize{6.9}{7.4}\selectfont] at (23.292,12.232) {$\textstyle A_V(x)$};
\node[anchor=north] at (2.217,2.819) {$0$};
\node[anchor=north] at (19.792,2.819) {$0$};
\node[anchor=west,font=\fontsize{6.9}{7.4}\selectfont] at (23.292,0.732) {$\textstyle g_V(x_1,x_2)$};
\end{tikzpicture}\\
&=
\left[
\begin{tikzpicture}[x=1bp,y=1bp,scale=.80,
  baseline={([yshift=-.5ex]current bounding box.center)},
  line cap=round,line join=round,
  every node/.style={inner sep=.5pt,outer sep=0pt,font=\fontsize{7.2}{8}\selectfont}]
\draw[black,line width=.65bp] (2.217,4.319) -- (2.217,9.592);
\draw[black,line width=.65bp] (2.217,9.592) .. controls (2.217,25.409) and (19.792,25.409) .. (19.792,9.592);
\draw[black,line width=.65bp] (19.792,9.592) -- (19.792,4.319);
\draw[black,fill=white,line width=.55bp] (2.222,10.069) circle[radius=2.834bp];
\draw[black,line width=.55bp] (0.664,10.069) -- (3.781,10.069);
\fill[black] (19.792,6.956) circle[radius=1.15bp];
\node[anchor=west,font=\fontsize{6.9}{7.4}\selectfont] at (23.292,6.956) {$\textstyle A_V(x)$};
\node[anchor=north] at (2.217,2.819) {$0$};
\node[anchor=north] at (19.792,2.819) {$0$};
\end{tikzpicture}
\frac{(-1)^{c+1}}2
f_V^+(u)f_V^-(u)n^+_{V,0}(u)n^-_{V,0}(u)
\bigl(1+u^K R_V(u)\bigr)
\frac{t_{V,0}(u)}{u^{e+1}}
\right]_{u:-1}\\
&\quad{}+
\;\begin{tikzpicture}[x=1bp,y=1bp,scale=.80,
  baseline={([yshift=-.5ex]current bounding box.center)},
  line cap=round,line join=round,
  every node/.style={inner sep=.5pt,outer sep=0pt,font=\fontsize{7.2}{8}\selectfont}]
\draw[black,line width=.65bp] (2.217,4.319) -- (2.217,10.47);
\draw[black,line width=.65bp] (2.217,10.47) .. controls (2.217,25.409) and (19.792,25.409) .. (19.792,10.47);
\draw[black,line width=.65bp] (19.792,10.47) -- (19.792,4.319);
\draw[black,densely dotted,line width=.55bp] (2.217,6.955) -- (19.792,6.955);
\fill[black] (19.706,12.232) circle[radius=1.15bp];
\node[anchor=west,font=\fontsize{6.9}{7.4}\selectfont] at (23.292,12.232) {$\textstyle A_V(x)$};
\node[anchor=north] at (2.217,2.819) {$0$};
\node[anchor=north] at (19.792,2.819) {$0$};
\node[anchor=west,font=\fontsize{6.9}{7.4}\selectfont] at (23.292,0.732) {$\textstyle g_V(x_1,x_2)$};
\end{tikzpicture}
-
\left[\begin{tikzpicture}[x=1bp,y=1bp,scale=.80,
  baseline={([yshift=-.5ex]current bounding box.center)},
  line cap=round,line join=round,
  every node/.style={inner sep=.5pt,outer sep=0pt,font=\fontsize{7.2}{8}\selectfont}]
\draw[black,line width=.65bp] (2.217,4.319) -- (2.217,9.592);
\draw[black,line width=.65bp] (2.217,9.592) .. controls (2.217,25.409) and (19.792,25.409) .. (19.792,9.592);
\draw[black,line width=.65bp] (19.792,9.592) -- (19.792,4.319);
\draw[black,fill=white,line width=.55bp] (19.738,11.074) circle[radius=2.834bp];
\draw[black,line width=.55bp] (18.18,11.074) -- (21.297,11.074);
\fill[black] (19.792,5.59) circle[radius=1.15bp];
\node[anchor=west,font=\fontsize{6.9}{7.4}\selectfont] at (24.572,5.59) {$\textstyle A_V(x)$};
\node[anchor=north] at (2.217,2.819) {$0$};
\node[anchor=north] at (19.792,2.819) {$0$};
\end{tikzpicture}\frac{C_V(u)}{2u}\right]_{u:-1}
.
\end{align*}
\endgroup

Definition~\ref{notation-1} and the zero-primary factorization give
\[
m_V(u)=u^e m_{V,0}(u),\qquad
n_V(u)=n_{V,0}(u)t_{V,0}(u),
\]
with $t_{V,0}$ monic of degree $e$.  Thus
\[
\frac{O_V(u)}u=
\frac{n_{V,0}(u)}{m_{V,0}(u)}
\frac{t_{V,0}(u)}{u^{e+1}}.
\]
For the same finite representatives, the B\'ezout identity gives
\[
C_V(u)\frac{O_V(u)}u
=
\frac{(-1)^{c+1}}2\frac{t_{V,0}(u)}{u^{e+1}}
+u^{K-e-1}t_{V,0}(u)H_V(u)
\]
for some polynomial $H_V$.  The additional strand factor $(u-x)^{-1}$
makes the first summand $O(u^{-2})$ at infinity.  For the polynomial tail, coefficient extraction gives
\[
\left[\frac{u^{K-e-1}t_{V,0}(u)H_V(u)}{u-x}\right]_{u:-1}
=x^{K-e-1}t_{V,0}(x)H_V(x)=0
\]
for sufficiently large $K$, since $x$ is nilpotent; see
Lemma~\ref{lem:finite-dd-evaluation}. Hence
\begingroup
\setlength{\jot}{2pt}
\begin{align*}
\begin{tikzpicture}[x=1bp,y=1bp,scale=.80,
  baseline={([yshift=-.5ex]current bounding box.center)},
  line cap=round,line join=round,
  every node/.style={inner sep=.5pt,outer sep=0pt,font=\fontsize{7.2}{8}\selectfont}]
\draw[darkred,line width=.65bp] (2.217,4.319) .. controls (8.075,11.349) and (19.792,16.075) .. (19.792,25.409);
\draw[darkred,line width=.65bp] (19.792,4.319) .. controls (13.934,11.349) and (2.217,16.075) .. (2.217,25.409);
\draw[darkred,line width=.65bp] (2.217,25.409) .. controls (2.217,43.863) and (19.792,43.863) .. (19.792,25.409);
\node[text=darkred,anchor=north] at (2.217,2.319) {$0$};
\draw[darkg,line width=.65bp] (26.792,4.319) -- (26.792,26.116);
\node[text=darkg,anchor=north] at (26.792,2.519) {$V$};
\end{tikzpicture}
&=\begin{tikzpicture}[x=1bp,y=1bp,scale=.80,
  baseline={([yshift=-.5ex]current bounding box.center)},
  line cap=round,line join=round,
  every node/.style={inner sep=.5pt,outer sep=0pt,font=\fontsize{7.2}{8}\selectfont}]
\draw[black,line width=.65bp] (2.217,4.319) -- (2.217,10.47);
\draw[black,line width=.65bp] (2.217,10.47) .. controls (2.217,25.409) and (19.792,25.409) .. (19.792,10.47);
\draw[black,line width=.65bp] (19.792,10.47) -- (19.792,4.319);
\draw[black,densely dotted,line width=.55bp] (2.217,6.955) -- (19.792,6.955);
\fill[black] (19.706,12.232) circle[radius=1.15bp];
\node[anchor=west,font=\fontsize{6.9}{7.4}\selectfont] at (23.292,12.232) {$\textstyle A_V(x)$};
\node[anchor=north] at (2.217,2.819) {$0$};
\node[anchor=north] at (19.792,2.819) {$0$};
\node[anchor=west,font=\fontsize{6.9}{7.4}\selectfont] at (23.292,0.732) {$\textstyle g_V(x_1,x_2)$};
\end{tikzpicture}
-
\left[\begin{tikzpicture}[x=1bp,y=1bp,scale=.80,
  baseline={([yshift=-.5ex]current bounding box.center)},
  line cap=round,line join=round,
  every node/.style={inner sep=.5pt,outer sep=0pt,font=\fontsize{7.2}{8}\selectfont}]
\draw[black,line width=.65bp] (2.217,4.319) -- (2.217,9.592);
\draw[black,line width=.65bp] (2.217,9.592) .. controls (2.217,25.409) and (19.792,25.409) .. (19.792,9.592);
\draw[black,line width=.65bp] (19.792,9.592) -- (19.792,4.319);
\draw[black,fill=white,line width=.55bp] (19.738,11.074) circle[radius=2.834bp];
\draw[black,line width=.55bp] (18.18,11.074) -- (21.297,11.074);
\fill[black] (19.792,5.59) circle[radius=1.15bp];
\node[anchor=west,font=\fontsize{6.9}{7.4}\selectfont] at (24.572,5.59) {$\textstyle A_V(x)$};
\node[anchor=north] at (2.217,2.819) {$0$};
\node[anchor=north] at (19.792,2.819) {$0$};
\end{tikzpicture}\frac{C_V(u)}{2u}\right]_{u:-1}
\\
&=\begin{tikzpicture}[x=1bp,y=1bp,scale=.80,
  baseline={([yshift=-.5ex]current bounding box.center)},
  line cap=round,line join=round,
  every node/.style={inner sep=.5pt,outer sep=0pt,font=\fontsize{7.2}{8}\selectfont}]
\draw[black,line width=.65bp] (2.217,4.319) -- (2.217,10.47);
\draw[black,line width=.65bp] (2.217,10.47) .. controls (2.217,25.409) and (19.792,25.409) .. (19.792,10.47);
\draw[black,line width=.65bp] (19.792,10.47) -- (19.792,4.319);
\draw[black,densely dotted,line width=.55bp] (2.217,6.955) -- (19.792,6.955);
\fill[black] (19.706,12.232) circle[radius=1.15bp];
\node[anchor=west,font=\fontsize{6.9}{7.4}\selectfont] at (23.292,12.232) {$\textstyle A_V(x)$};
\node[anchor=north] at (2.217,2.819) {$0$};
\node[anchor=north] at (19.792,2.819) {$0$};
\node[anchor=west,font=\fontsize{6.9}{7.4}\selectfont] at (23.292,0.732) {$\textstyle g_V(x_1,x_2)$};
\end{tikzpicture}
-\;\begin{tikzpicture}[x=1bp,y=1bp,scale=.80,
  baseline={([yshift=-.5ex]current bounding box.center)},
  line cap=round,line join=round,
  every node/.style={inner sep=.5pt,outer sep=0pt,font=\fontsize{7.2}{8}\selectfont}]
\draw[black,line width=.65bp] (2.217,4.319) .. controls (2.217,26.639) and (19.792,26.639) .. (19.792,4.319);
\fill[black] (18.778,12.668) circle[radius=1.15bp];
\node[anchor=north] at (2.217,2.819) {$0$};
\node[anchor=north] at (19.792,2.819) {$0$};
\end{tikzpicture}\,
A_V(x){P_C(x)}.
\end{align*}
\endgroup

By the definition of $g_V$,
\[
g_V(-x,x)=\frac{1-A_V(x)B_V(x)}{-2x}=P_C(x).
\]
Thus the two diagrams in the last line cancel, proving
relation~\textnormal{(3)} of Definition~\ref{nil}.
\end{proof}

Apply the anti-isomorphism of Lemma~\ref{anti}, which interchanges
the black cup and cap and fixes the dot and crossing. Interchanging
$A_V$ and $B_V$ in the reflected calculation gives the corresponding
red cup identity:
\begin{equation}
\begin{tikzpicture}[x=1bp,y=1bp,scale=.80,
  baseline={([yshift=-.5ex]current bounding box.center)},
  line cap=round,line join=round,
  every node/.style={inner sep=.5pt,outer sep=0pt,font=\fontsize{7.2}{8}\selectfont}]
\draw[darkred,line width=.65bp] (2.217,13.84) .. controls (2.217,23.174) and (13.934,27.9) .. (19.792,34.93);
\draw[darkred,line width=.65bp] (19.792,13.84) .. controls (19.792,23.174) and (8.075,27.9) .. (2.217,34.93);
\draw[darkred,line width=.65bp] (2.217,13.84) .. controls (2.217,-4.613) and (19.792,-4.613) .. (19.792,13.84);
\node[text=darkred,anchor=south] at (2.217,36.93) {$0$};
\end{tikzpicture}=0.
\label{nbone:eq:52}
\end{equation}

\begin{Lemma}\label{nbone:lemma:nil2}
Relation~\textnormal{(2)} of Definition~\ref{nil} holds with $t=c_\lambda$.
\end{Lemma}

\begin{proof}
Let $\beta_V$ be the undotted red bubble on $V$, and use the
polynomial $C_V=A_VB_V$ from Lemma~\ref{nbone:lemma:R4}.
By Definition~\ref{eq:red-cup} and Lemma~\ref{eq:red-cup-E0V}, this is the black bubble decorated by $C_V$.
The support factors remove the nonzero primary components, so
\[
\beta_V=[C_V(u)b_V(u)]_{u:-1},
\qquad
b_V(u)=\frac{O_V(u)}u-1+\frac1{2u}.
\]
By the calculation in Lemma~\ref{nbone:lemma:R4}, the polynomial tail contributes no negative power,
while monicity of $t_{V,0}$ gives
\[
\left[C_V(u)\frac{O_V(u)}u\right]_{u:-1}
=\frac{(-1)^{c_\lambda+1}}2.
\]
Since $C_V(0)=1$,
\[
\beta_V=
\frac{(-1)^{c_\lambda+1}+1}{2}\,\operatorname{Id}_V
=c_\lambda\,\operatorname{Id}_V.
\]

\end{proof}

\begin{Lemma} \label{nbone:lemma:nil4}
Relation~\textnormal{(4)} of Definition~\ref{nil} holds.
\end{Lemma}
\begin{proof}
By relation~\textnormal{(1)} of Definition~\ref{nil}, it is enough to prove
the equivalent form
\[
\begin{tikzpicture}[x=1bp,y=1bp,scale=.74,
  baseline={([yshift=-.5ex]current bounding box.center)},
  line cap=round,line join=round,
  every node/.style={inner sep=.45pt,outer sep=0pt,font=\fontsize{7.0}{7.6}\selectfont}]
\draw[darkred,line width=.65bp] (2.217,4.319) .. controls (2.217,25.409) and (15.398,25.409) .. (22.428,13.107);
\draw[darkred,line width=.65bp] (22.428,13.107) .. controls (29.458,0.804) and (40.881,0.804) .. (40.881,22.773);
\draw[darkred,line width=.65bp] (13.64,4.319) -- (28.579,22.773);
\node[text=darkred,anchor=north] at (2.217,2.319) {$0$};
\node[text=darkred,anchor=north] at (13.64,2.319) {$0$};
\end{tikzpicture}
=
\begin{tikzpicture}[x=1bp,y=1bp,scale=.74,
  baseline={([yshift=-.5ex]current bounding box.center)},
  line cap=round,line join=round,
  every node/.style={inner sep=.45pt,outer sep=0pt,font=\fontsize{7.0}{7.6}\selectfont}]
\draw[darkred,line width=.65bp] (2.217,4.319) -- (19.792,27.167);
\draw[darkred,line width=.65bp] (19.792,4.319) -- (2.217,27.167);
\node[text=darkred,anchor=north] at (2.217,2.319) {$0$};
\node[text=darkred,anchor=north] at (19.792,2.319) {$0$};
\end{tikzpicture}.
\]
If $E_0V=0$, both sides vanish. Put $N=E_0V$ otherwise.
Substituting the red generators and straightening the black cup and cap
gives
\[
\begin{tikzpicture}[x=1bp,y=1bp,scale=.86,
baseline=17bp,
line cap=round,line join=round,
every node/.style={inner sep=1pt,outer sep=0pt,font=\scriptsize}]
\draw[black,line width=.65bp] (0,0)--(0,10)--(24,34)--(24,44);
\draw[black,line width=.65bp] (24,0)--(24,10)--(0,34)--(0,44);
\fill[black] (0,5) circle[radius=1.1bp];
\node[anchor=east] at (-5,5) {$ B_N(x)^{-1} $};
\fill[black] (24,39) circle[radius=1.1bp];
\node[anchor=west] at (29,39) {$ B_V(x) $};
\draw[black,densely dotted,line width=.55bp] (0,10)--(-8,10)--(-8,34)--(0,34);
\draw[black,densely dotted,line width=.55bp] (24,10)--(32,10)--(32,34)--(24,34);
\node[anchor=east] at (-12,22) {$B_N(\frac{x_1+x_2}{2})$};
\node[anchor=west] at (36,22) {$A_N(-\frac{x_1+x_2}{2})$};
\end{tikzpicture}\quad+\quad\begin{tikzpicture}[x=1bp,y=1bp,scale=.86,
baseline=17bp,
line cap=round,line join=round,
every node/.style={inner sep=1pt,outer sep=0pt,font=\scriptsize}]
\draw[black,line width=.65bp] (0,0)--(0,8) .. controls (0,24) and (24,24) .. (24,8)--(24,0);
\draw[black,line width=.65bp] (0,44)--(0,36) .. controls (0,20) and (24,20) .. (24,36)--(24,44);
\fill[black] (24,5) circle[radius=1.1bp];
\node[anchor=west] at (29,5) {$ A_V(x) $};
\fill[black] (24,39) circle[radius=1.1bp];
\node[anchor=west] at (29,39) {$ B_V(x) $};
\draw[black,densely dotted,line width=.55bp] (24,8)--(31,8)--(31,36)--(24,36);
\node[anchor=west] at (35,22) {$ g_N(x_1,-x_2) $};
\end{tikzpicture}.
\]
All boundary colors here are $0$. On each vertical dotted connector,
$x_1,x_2$ refer to its lower and upper ends; on a horizontal connector
they refer to the left and right strands. A single $x$ denotes the dot
at the marked point. In cup--cap coefficients below, $x_1$ is on the
right input leg and $x_2$ on the right output leg.

We use strand resolvents, with
\[
\left[\frac{f(u)}{u-x}\right]_{u:-1}=f(x).
\]
The circled signs denote $(u-x)^{-1}$ and $(u+x)^{-1}$, respectively.
The black dot relation gives the following two moves:
\[
\begin{tikzpicture}[x=1bp,y=1bp,scale=.86,
baseline=17bp,
line cap=round,line join=round,
every node/.style={inner sep=1pt,outer sep=0pt,font=\scriptsize}]
\draw[black,line width=.65bp] (0,0)--(0,10)--(24,34)--(24,44);
\draw[black,line width=.65bp] (24,0)--(24,10)--(0,34)--(0,44);
\draw[black,fill=white,line width=.55bp] (0,39) circle[radius=2.7bp];
\draw[black,line width=.55bp] (-1.5,39)--(1.5,39);
\end{tikzpicture}=\begin{tikzpicture}[x=1bp,y=1bp,scale=.86,
baseline=17bp,
line cap=round,line join=round,
every node/.style={inner sep=1pt,outer sep=0pt,font=\scriptsize}]
\draw[black,line width=.65bp] (0,0)--(0,10)--(24,34)--(24,44);
\draw[black,line width=.65bp] (24,0)--(24,10)--(0,34)--(0,44);
\draw[black,fill=white,line width=.55bp] (24,5) circle[radius=2.7bp];
\draw[black,line width=.55bp] (22.5,5)--(25.5,5);
\end{tikzpicture}+\begin{tikzpicture}[x=1bp,y=1bp,scale=.86,
baseline=17bp,
line cap=round,line join=round,
every node/.style={inner sep=1pt,outer sep=0pt,font=\scriptsize}]
\draw[black,line width=.65bp] (0,0)--(0,44) (24,0)--(24,44);
\draw[black,fill=white,line width=.55bp] (0,39) circle[radius=2.7bp];
\draw[black,line width=.55bp] (-1.5,39)--(1.5,39);
\draw[black,fill=white,line width=.55bp] (24,5) circle[radius=2.7bp];
\draw[black,line width=.55bp] (22.5,5)--(25.5,5);
\end{tikzpicture}-\begin{tikzpicture}[x=1bp,y=1bp,scale=.86,
baseline=17bp,
line cap=round,line join=round,
every node/.style={inner sep=1pt,outer sep=0pt,font=\scriptsize}]
\draw[black,line width=.65bp] (0,0)--(0,8) .. controls (0,24) and (24,24) .. (24,8)--(24,0);
\draw[black,line width=.65bp] (0,44)--(0,36) .. controls (0,20) and (24,20) .. (24,36)--(24,44);
\draw[black,fill=white,line width=.55bp] (24,39) circle[radius=2.7bp];
\draw[black,line width=.55bp] (22.5,39)--(25.5,39);
\draw[black,line width=.55bp] (24,37.5)--(24,40.5);
\draw[black,fill=white,line width=.55bp] (24,5) circle[radius=2.7bp];
\draw[black,line width=.55bp] (22.5,5)--(25.5,5);
\end{tikzpicture},
\]
\[
\begin{tikzpicture}[x=1bp,y=1bp,scale=.86,
baseline=17bp,
line cap=round,line join=round,
every node/.style={inner sep=1pt,outer sep=0pt,font=\scriptsize}]
\draw[black,line width=.65bp] (0,0)--(0,10)--(24,34)--(24,44);
\draw[black,line width=.65bp] (24,0)--(24,10)--(0,34)--(0,44);
\fill[black] (24,39) circle[radius=1.1bp];
\node[anchor=west] at (29,39) {$ f(x) $};
\end{tikzpicture}=\begin{tikzpicture}[x=1bp,y=1bp,scale=.86,
baseline=17bp,
line cap=round,line join=round,
every node/.style={inner sep=1pt,outer sep=0pt,font=\scriptsize}]
\draw[black,line width=.65bp] (0,0)--(0,10)--(24,34)--(24,44);
\draw[black,line width=.65bp] (24,0)--(24,10)--(0,34)--(0,44);
\fill[black] (0,5) circle[radius=1.1bp];
\node[anchor=east] at (-5,5) {$ f(x) $};
\end{tikzpicture}-\begin{tikzpicture}[x=1bp,y=1bp,scale=.86,
baseline=17bp,
line cap=round,line join=round,
every node/.style={inner sep=1pt,outer sep=0pt,font=\scriptsize}]
\draw[black,line width=.65bp] (0,0)--(0,44) (24,0)--(24,44);
\draw[black,densely dotted,line width=.55bp] (0,5)--(24,5);
\node[anchor=north] at (12,-5) {$ \frac{f(x_2)-f(x_1)}{x_2-x_1} $};
\end{tikzpicture}+\begin{tikzpicture}[x=1bp,y=1bp,scale=.86,
baseline=17bp,
line cap=round,line join=round,
every node/.style={inner sep=1pt,outer sep=0pt,font=\scriptsize}]
\draw[black,line width=.65bp] (0,0)--(0,8) .. controls (0,24) and (24,24) .. (24,8)--(24,0);
\draw[black,line width=.65bp] (0,44)--(0,36) .. controls (0,20) and (24,20) .. (24,36)--(24,44);
\draw[black,densely dotted,line width=.55bp] (24,8)--(31,8)--(31,36)--(24,36);
\node[anchor=west] at (35,22) {$ \frac{f(x_2)-f(-x_1)}{x_1+x_2} $};
\end{tikzpicture}.
\]
The first move is applied to the upper left leg with coefficient
$B_N((x+u)/2)$, where $x$ is on the lower left leg. The second moves
$A_N(-(x+u)/2)B_V(u)$ from the upper right leg, where $x$ is on the
lower right leg. Only dot powers are moved; the $N$-indexed coefficients
remain on their original strand.

For the straight-strand terms, put $h=(x_1+x_2)/2$. The two moves give,
in order, the following coefficients:
\[
\begin{aligned}
&\frac{A_N(-x_2)B_V(x_2)}{B_N(x_1)}
 \frac{B_N(x_1)-B_N(h)}{x_1-x_2},\\
&\frac{B_N(h)}{B_N(x_1)}
 \frac{A_N(-x_2)B_V(x_2)-A_N(-h)B_V(x_1)}{x_1-x_2}.
\end{aligned}
\]
For the cup--cap terms, with $h=(x_1+x_2)/2$ in the stated input--output
convention, they give
\[
\begin{aligned}
&-\frac{A_N(-h)B_V(x_2)}{B_N(-x_1)}
 \frac{B_N(0)-B_N(-h)}{x_1+x_2},\\
&\frac{B_N(0)}{B_N(-x_1)}
 \frac{A_N(-h)B_V(x_2)-A_N(0)B_V(-x_1)}{x_1+x_2}.
\end{aligned}
\]
These are regular divided differences: each numerator vanishes on the
corresponding diagonal. Division is performed before nilpotent
evaluation, as in Lemma~\ref{lem:finite-dd-evaluation}.

Lemma~\ref{gamma-move1} and the primary-factor calculation in the proof
of Lemma~\ref{eq:red-cup-E0V} give, with $y$ the dot on $N=E_0V$,
\[
B_N(u)^{-1}=E_0\bigl(B_V(u)^{-1}\bigr)
\frac{u+y+1}{u-y+1}.
\]
Together with $A_N(-u)B_N(u)=(1+2u)^{-1}$, this gives
\[
\frac{B_V(x_1)}{B_N(x_1)}=
\frac{1+x_1+x_2}{1+x_1-x_2},\qquad
A_N(-x_2)B_V(x_2)=1
\]
on straight strands. On the input cap, it gives
$B_N(-x_1)^{-1}=A_V(x_1)$. Thus the crossing contribution is
\[
\begin{tikzpicture}[x=1bp,y=1bp,scale=.86,
baseline=17bp,
line cap=round,line join=round,
every node/.style={inner sep=1pt,outer sep=0pt,font=\scriptsize}]
\draw[black,line width=.65bp] (0,0)--(0,10)--(24,34)--(24,44);
\draw[black,line width=.65bp] (24,0)--(24,10)--(0,34)--(0,44);
\draw[black,densely dotted,line width=.55bp] (0,5)--(24,5);
\node[anchor=north] at (12,-5) {$ \frac{A_N(-h)B_N(h)B_V(x_1)}{B_N(x_1)} $};
\end{tikzpicture}=\begin{tikzpicture}[x=1bp,y=1bp,scale=.86,
baseline=17bp,
line cap=round,line join=round,
every node/.style={inner sep=1pt,outer sep=0pt,font=\scriptsize}]
\draw[black,line width=.65bp] (0,0)--(0,10)--(24,34)--(24,44);
\draw[black,line width=.65bp] (24,0)--(24,10)--(0,34)--(0,44);
\draw[black,densely dotted,line width=.55bp] (0,5)--(24,5);
\node[anchor=north] at (12,-5) {$ \frac1{1+x_1-x_2} $};
\end{tikzpicture},
\]
where the dotted decoration lies entirely on the input side.
The sum of the two straight-strand coefficients is
\[
\begin{aligned}
&\frac{B_N(x_1)A_N(-x_2)B_V(x_2)
 -A_N(-h)B_N(h)B_V(x_1)}{B_N(x_1)(x_1-x_2)}\\
&\qquad=\frac{1-(1+x_1-x_2)^{-1}}{x_1-x_2}
=\frac1{1+x_1-x_2}.
\end{aligned}
\]
For the cup--cap terms, use $A_N(0)B_N(0)=1$ and add the $g_N$ term.
After factoring out $A_V(x_1)$, their coefficient is
\[
\begin{aligned}
&B_V(x_2)\frac{1-A_N(-h)B_N(-h)}{x_1+x_2}
 +\frac{B_V(-x_1)-B_V(x_2)A_N(-h)B_N(-h)}{-x_1-x_2}\\
&\qquad=\frac{B_V(x_2)-B_V(-x_1)}{x_1+x_2}
=q_V(x_2,x_1).
\end{aligned}
\]
Consequently the original expression is
\[
\begin{tikzpicture}[x=1bp,y=1bp,scale=.86,
baseline=17bp,
line cap=round,line join=round,
every node/.style={inner sep=1pt,outer sep=0pt,font=\scriptsize}]
\draw[black,line width=.65bp] (0,0)--(0,10)--(24,34)--(24,44);
\draw[black,line width=.65bp] (24,0)--(24,10)--(0,34)--(0,44);
\draw[black,densely dotted,line width=.55bp] (0,5)--(24,5);
\node[anchor=north] at (12,-5) {$ \frac1{1+x_1-x_2} $};
\end{tikzpicture}+\begin{tikzpicture}[x=1bp,y=1bp,scale=.86,
baseline=17bp,
line cap=round,line join=round,
every node/.style={inner sep=1pt,outer sep=0pt,font=\scriptsize}]
\draw[black,line width=.65bp] (0,0)--(0,44) (24,0)--(24,44);
\draw[black,densely dotted,line width=.55bp] (0,5)--(24,5);
\node[anchor=north] at (12,-5) {$ \frac1{1+x_1-x_2} $};
\end{tikzpicture}+\begin{tikzpicture}[x=1bp,y=1bp,scale=.86,
baseline=17bp,
line cap=round,line join=round,
every node/.style={inner sep=1pt,outer sep=0pt,font=\scriptsize}]
\draw[black,line width=.65bp] (0,0)--(0,8) .. controls (0,24) and (24,24) .. (24,8)--(24,0);
\draw[black,line width=.65bp] (0,44)--(0,36) .. controls (0,20) and (24,20) .. (24,36)--(24,44);
\fill[black] (24,5) circle[radius=1.1bp];
\node[anchor=west] at (29,5) {$ A_V(x_1) $};
\draw[black,densely dotted,line width=.55bp] (24,8)--(31,8)--(31,36)--(24,36);
\node[anchor=west] at (35,22) {$ q_V(x_2,x_1) $};
\end{tikzpicture}.
\]
This is the red crossing by Lemma~\ref{lem:anothdefofcro}.
\end{proof}

\begin{Lemma} \label{nilRel_7}
The generators satisfy relation~\textnormal{(7)} of Definition~\ref{nil}; equivalently,
\[
\begin{tikzpicture}[x=1bp,y=1bp,scale=.64,
  baseline={([yshift=-.5ex]current bounding box.center)},
  line cap=round,line join=round,
  every node/.style={inner sep=.4pt,outer sep=0pt,font=\fontsize{6.8}{7.3}\selectfont}]
\draw[darkred,line width=.65bp] (2.217,4.319) .. controls (8.075,11.554) and (19.792,16.528) .. (19.792,26.024);
\draw[darkred,line width=.65bp] (19.792,4.319) .. controls (13.934,11.554) and (2.217,16.528) .. (2.217,26.024);
\draw[darkred,line width=.65bp] (2.217,26.024) .. controls (2.217,35.52) and (13.934,40.494) .. (19.792,47.729);
\draw[darkred,line width=.65bp] (19.792,26.024) .. controls (19.792,35.52) and (8.075,40.494) .. (2.217,47.729);
\node[text=darkred,anchor=north] at (2.217,2.319) {$0$};
\node[text=darkred,anchor=north] at (19.792,2.319) {$0$};
\end{tikzpicture}
= 0 \ .
\]
\end{Lemma}

\begin{proof}
{Define $S_V$ by}
\[
S_V :=
\begin{tikzpicture}[x=1bp,y=1bp,scale=.56,
  baseline={([yshift=-.5ex]current bounding box.center)},
  line cap=round,line join=round,
  every node/.style={inner sep=.4pt,outer sep=0pt,font=\fontsize{6.8}{7.3}\selectfont}]
\draw[darkred,line width=.65bp] (0,6.57) .. controls (3.156,15.568) and (17.52,18.934) .. (17.52,28.47);
\draw[darkred,line width=.65bp] (17.52,28.47) .. controls (17.52,38.006) and (3.156,41.372) .. (0,50.37);
\draw[darkred,line width=.65bp] (17.52,6.57) .. controls (14.364,15.568) and (0,18.934) .. (0,28.47);
\draw[darkred,line width=.65bp] (0,28.47) .. controls (0,38.006) and (14.364,41.372) .. (17.52,50.37);
\draw[black,line width=.65bp] (0,0) .. controls (0,2.19) and (-0.739,4.462) .. (0,6.57);
\draw[black,line width=.65bp] (17.52,0) .. controls (17.52,2.19) and (18.259,4.462) .. (17.52,6.57);
\draw[black,line width=.65bp] (0,50.37) .. controls (-0.739,52.478) and (0,54.75) .. (0,56.94);
\draw[black,line width=.65bp] (17.52,50.37) .. controls (18.259,52.478) and (17.52,54.75) .. (17.52,56.94);
\draw[black,densely dotted,line width=.55bp] (-0.222,54.203) -- (17.742,54.203);
\draw[black,densely dotted,line width=.55bp] (-0.222,2.737) -- (17.742,2.737);
\node[anchor=south] at (0,58.254) {$0$};
\node[anchor=north] at (0,-1.314) {$0$};
\node[anchor=south] at (17.52,58.254) {$0$};
\node[anchor=north] at (17.52,-1.314) {$0$};
\node[anchor=west,font=\fontsize{6.3}{6.8}\selectfont] at (21.349,54.203) {$\textstyle 1-x_1+x_2$};
\node[anchor=west,font=\fontsize{6.3}{6.8}\selectfont] at (21.349,2.737) {$\textstyle 1+x_1-x_2$};
\draw[darkg,line width=.65bp] (100,0) -- (100,35.531);
\node[text=darkg,anchor=north] at (100,-1.8) {$V$};
\end{tikzpicture}
\]
Since the dot-polynomial factors \(1+x_1-x_2\) and \(1-x_1+x_2\) are invertible on $E_0^2V$, it suffices to prove $S_V=0$.

Put $\Delta=x_1-x_2$, $L=1+\Delta$, and $\bar L=1-\Delta$.
Let $X$ be the red crossing and $s$ the black crossing projected to the
zero--zero component.  Let $E_q$ and $E_p$ be the complete cup--cap terms
in the two formulas of Lemma~\ref{lem:anothdefofcro}; their decorations are
$q_V(x_2,x_1)A_V(x_1)$ and $B_V(x_2)p_V(x_2,x_1)$, respectively.  Set
\[
Q=E_qL,\qquad P=\bar L E_p.
\]
The two alternative crossing formulas in Lemma~\ref{lem:anothdefofcro} give
\[
XL=s+\operatorname{Id}+Q,\qquad
\bar L X=s-\operatorname{Id}+P.
\]
Hence
\[
S_V=\bar L X^2L
=(s-\operatorname{Id}+P)(s+\operatorname{Id}+Q)
=s^2-\operatorname{Id}+sQ-Q+Ps+P+PQ.
\]

Let $\sigma$ be the full affine Brauer crossing, let $\pi$ project
onto $E_0^2V$, and put $D=\pi\sigma(1-\pi)\sigma\pi$; thus $s=\pi\sigma\pi$.
{The preceding identities are on $E_0^2V$. From now on, extend
$X,s,P,Q,S_V$ by zero to $E^2V$. Under this convention
$\operatorname{Id}_{E_0^2V}$ becomes $\pi$, and the expansion reads
\[
S_V=s^2-\pi+sQ-Q+Ps+P+PQ.
\]
}
Choose compatible polynomial representatives as in
Definition~\ref{notation-1} and Lemma~\ref{lem:finite-dd-evaluation}.
Let $\widehat Q$ and $\widehat P$ be the following cup--cap operators,
with all unmarked colors summed:
\[
\widehat Q=\begin{tikzpicture}[x=1bp,y=1bp,scale=.56,
  baseline={([yshift=-.5ex]current bounding box.center)},
  line cap=round,line join=round,
  every node/.style={inner sep=.4pt,outer sep=0pt,font=\fontsize{6.8}{7.3}\selectfont}]
\draw[black,line width=.65bp] (2.217,57.044) -- (2.217,39.469);
\draw[black,line width=.65bp] (2.217,39.469) .. controls (2.217,28.045) and (19.792,28.045) .. (19.792,39.469);
\draw[black,line width=.65bp] (19.792,39.469) -- (19.792,57.044);
\draw[black,line width=.65bp] (2.217,4.319) -- (2.217,7.834);
\draw[black,line width=.65bp] (2.217,7.834) .. controls (2.217,19.258) and (19.792,19.258) .. (19.792,7.834);
\draw[black,line width=.65bp] (19.792,7.834) -- (19.792,4.319);
\draw[black,densely dotted,line width=.55bp] (17.046,33.043) -- (17.046,14.26);
\fill[black] (19.346,10.781) circle[radius=1.15bp];
\node[anchor=west,font=\fontsize{6.3}{6.8}\selectfont] at (23.292,10.781) {$\textstyle M_V(-x_1)$};
\node[anchor=west,font=\fontsize{6.3}{6.8}\selectfont] at (20.546,23.652) {$\textstyle q_V(x_2,x_1)$};
\end{tikzpicture},\qquad
\widehat P=\begin{tikzpicture}[x=1bp,y=1bp,scale=.56,
  baseline={([yshift=-.5ex]current bounding box.center)},
  line cap=round,line join=round,
  every node/.style={inner sep=.4pt,outer sep=0pt,font=\fontsize{6.8}{7.3}\selectfont}]
\draw[black,line width=.65bp] (2.217,57.044) -- (2.217,53.529);
\draw[black,line width=.65bp] (2.217,53.529) .. controls (2.217,42.105) and (19.792,42.105) .. (19.792,53.529);
\draw[black,line width=.65bp] (19.792,53.529) -- (19.792,57.044);
\draw[black,line width=.65bp] (2.217,4.319) -- (2.217,21.894);
\draw[black,line width=.65bp] (2.217,21.894) .. controls (2.217,33.317) and (19.792,33.317) .. (19.792,21.894);
\draw[black,line width=.65bp] (19.792,21.894) -- (19.792,4.319);
\draw[black,densely dotted,line width=.55bp] (17.046,47.103) -- (17.046,28.319);
\fill[black] (19.449,50.921) circle[radius=1.15bp];
\node[anchor=west,font=\fontsize{6.3}{6.8}\selectfont] at (23.292,50.921) {$\textstyle W_V(-x_2)$};
\node[anchor=west,font=\fontsize{6.3}{6.8}\selectfont] at (20.546,37.711) {$\textstyle p_V(x_2,x_1)$};
\end{tikzpicture}.
\]
In these cup--cap coefficients, $x_1$ and $x_2$ are the independent dots
on the right legs of the input cap and output cup, respectively;
horizontal decorations use the left-to-right strand order.
Moving the normalization factors
through the cup and cap gives $Q=\pi\widehat Q\pi$ and $P=\pi\widehat P\pi$.
{Only opposite-color pairs contribute to a cup--cap operator.
On a block with zero input and nonzero output, the retained factors in
$B_V$ give $B_V(x_2)=0$ and hence $W_V(-x_2)=0$. Since
$B_V(-x_1)M_V(-x_1)=1$ on the input, the two coefficients become
\[
q_V(x_2,x_1)M_V(-x_1)=-\frac1{x_1+x_2},\qquad
W_V(-x_2)p_V(x_2,x_1)=0.
\]
With nonzero input and zero output, the retained factors instead give
$A_V(x_1)=M_V(-x_1)=0$. The unit identity on the output gives
$A_V(-x_2)W_V(-x_2)=1$, hence
\[
q_V(x_2,x_1)M_V(-x_1)=0,\qquad
W_V(-x_2)p_V(x_2,x_1)=-\frac1{x_1+x_2}.
\]
Here $x_1+x_2$ is invertible, since its eigenvalue is nonzero.
The black dot relation gives crossing coefficient $1/(x_1+x_2)$,
and hence}
\[
\begin{gathered}
(1-\pi)\widehat Q\pi=-(1-\pi)\sigma\pi,\qquad
\pi\widehat P(1-\pi)=-\pi\sigma(1-\pi),\\
\pi\widehat Q(1-\pi)=0=(1-\pi)\widehat P\pi.
\end{gathered}
\]
Consequently,
\[
\begin{aligned}
\pi(\sigma\widehat Q-\widehat Q)\pi&=sQ-Q-D,\\
\pi(\widehat P\sigma+\widehat P)\pi&=Ps+P-D,\\
\pi\widehat P\widehat Q\pi&=PQ+D.
\end{aligned}
\]
We calculate these three products with all intermediate colors summed.

{In the coefficients below, form the regular divided differences
at zero from $\widehat A_V,\widehat B_V$, and then take sufficiently
high polynomial representatives, as in
Lemma~\ref{lem:finite-dd-evaluation}. Coefficient extractions at infinity
use the finite representatives retaining $m^\pm_{V,0}f_V^\pm$.
}

\noindent\emph{The \(q_V\) cup--cap terms.}
The first product, $\pi(\sigma\widehat Q-\widehat Q)\pi$, is
\[
\begin{tikzpicture}[x=1bp,y=1bp,scale=.56,
  baseline={([yshift=-.5ex]current bounding box.center)},
  line cap=round,line join=round,
  every node/.style={inner sep=.4pt,outer sep=0pt,font=\fontsize{6.8}{7.3}\selectfont}]
\draw[black,line width=.65bp] (2.217,57.044) .. controls (8.075,50.6) and (19.792,46.594) .. (19.792,37.711);
\draw[black,line width=.65bp] (19.792,37.711) .. controls (19.792,28.045) and (2.217,28.045) .. (2.217,37.711);
\draw[black,line width=.65bp] (2.217,37.711) .. controls (2.217,46.594) and (13.934,50.6) .. (19.792,57.044);
\draw[black,line width=.65bp] (2.217,4.319) -- (2.217,7.834);
\draw[black,line width=.65bp] (2.217,7.834) .. controls (2.217,19.258) and (19.792,19.258) .. (19.792,7.834);
\draw[black,line width=.65bp] (19.792,7.834) -- (19.792,4.319);
\draw[black,densely dotted,line width=.55bp] (17.046,32.274) -- (17.046,14.26);
\draw[black,densely dotted,line width=.55bp] (2.217,6.78) -- (19.792,6.78);
\fill[black] (19.346,10.781) circle[radius=1.15bp];
\node[anchor=west,font=\fontsize{6.3}{6.8}\selectfont] at (23.292,10.781) {$\textstyle A_V(x)$};
\node[anchor=north] at (2.217,2.819) {$0$};
\node[anchor=north] at (19.792,2.819) {$0$};
\node[anchor=south] at (1.208,58.154) {$0$};
\node[anchor=south] at (20.801,58.154) {$0$};
\node[anchor=west,font=\fontsize{6.3}{6.8}\selectfont] at (20.546,23.267) {$\textstyle q_V(x_2,x_1)$};
\node[anchor=west,font=\fontsize{6.3}{6.8}\selectfont] at (23.292,-0.719) {$\textstyle 1+x_1-x_2$};
\end{tikzpicture}
-
\begin{tikzpicture}[x=1bp,y=1bp,scale=.56,
  baseline={([yshift=-.5ex]current bounding box.center)},
  line cap=round,line join=round,
  every node/.style={inner sep=.4pt,outer sep=0pt,font=\fontsize{6.8}{7.3}\selectfont}]
\draw[black,line width=.65bp] (2.217,57.044) -- (2.217,39.469);
\draw[black,line width=.65bp] (2.217,39.469) .. controls (2.217,28.045) and (19.792,28.045) .. (19.792,39.469);
\draw[black,line width=.65bp] (19.792,39.469) -- (19.792,57.044);
\draw[black,line width=.65bp] (2.217,4.319) -- (2.217,7.834);
\draw[black,line width=.65bp] (2.217,7.834) .. controls (2.217,19.258) and (19.792,19.258) .. (19.792,7.834);
\draw[black,line width=.65bp] (19.792,7.834) -- (19.792,4.319);
\draw[black,densely dotted,line width=.55bp] (17.046,33.043) -- (17.046,14.26);
\draw[black,densely dotted,line width=.55bp] (2.217,6.78) -- (19.792,6.78);
\fill[black] (19.346,10.781) circle[radius=1.15bp];
\node[anchor=west,font=\fontsize{6.3}{6.8}\selectfont] at (23.292,10.781) {$\textstyle A_V(x)$};
\node[anchor=north] at (2.217,2.819) {$0$};
\node[anchor=north] at (19.792,2.819) {$0$};
\node[anchor=south] at (2.217,58.544) {$0$};
\node[anchor=south] at (19.792,58.544) {$0$};
\node[anchor=west,font=\fontsize{6.3}{6.8}\selectfont] at (20.546,23.652) {$\textstyle q_V(x_2,x_1)$};
\node[anchor=west,font=\fontsize{6.3}{6.8}\selectfont] at (23.292,-0.719) {$\textstyle 1+x_1-x_2$};
\end{tikzpicture}
=
\left[\begin{tikzpicture}[x=1bp,y=1bp,scale=.56,
  baseline={([yshift=-.5ex]current bounding box.center)},
  line cap=round,line join=round,
  every node/.style={inner sep=.4pt,outer sep=0pt,font=\fontsize{6.8}{7.3}\selectfont}]
\draw[black,line width=.65bp] (2.217,57.044) .. controls (8.075,50.6) and (19.792,46.594) .. (19.792,37.711);
\draw[black,line width=.65bp] (19.792,37.711) .. controls (19.792,28.045) and (2.217,28.045) .. (2.217,37.711);
\draw[black,line width=.65bp] (2.217,37.711) .. controls (2.217,46.594) and (13.934,50.6) .. (19.792,57.044);
\draw[black,line width=.65bp] (2.217,4.319) -- (2.217,7.834);
\draw[black,line width=.65bp] (2.217,7.834) .. controls (2.217,19.258) and (19.792,19.258) .. (19.792,7.834);
\draw[black,line width=.65bp] (19.792,7.834) -- (19.792,4.319);
\draw[black,fill=white,line width=.55bp] (19.576,39.753) circle[radius=2.834bp];
\draw[black,line width=.55bp] (18.018,39.753) -- (21.135,39.753);
\fill[black] (19.465,10.385) circle[radius=1.15bp];
\node[anchor=north] at (2.217,2.819) {$0$};
\node[anchor=north] at (19.792,2.819) {$0$};
\node[anchor=south] at (1.208,58.154) {$0$};
\node[anchor=south] at (20.801,58.154) {$0$};
\node[anchor=west,font=\fontsize{6.3}{6.8}\selectfont] at (23.292,10.385) {$\textstyle q_V(u,x)M_V(-x)$};
\end{tikzpicture}
-
\begin{tikzpicture}[x=1bp,y=1bp,scale=.56,
  baseline={([yshift=-.5ex]current bounding box.center)},
  line cap=round,line join=round,
  every node/.style={inner sep=.4pt,outer sep=0pt,font=\fontsize{6.8}{7.3}\selectfont}]
\draw[black,line width=.65bp] (2.217,57.045) -- (2.217,39.47);
\draw[black,line width=.65bp] (2.217,39.47) .. controls (2.217,28.046) and (19.792,28.046) .. (19.792,39.47);
\draw[black,line width=.65bp] (19.792,39.47) -- (19.792,57.045);
\draw[black,line width=.65bp] (2.217,4.32) -- (2.217,7.835);
\draw[black,line width=.65bp] (2.217,7.835) .. controls (2.217,19.259) and (19.792,19.259) .. (19.792,7.835);
\draw[black,line width=.65bp] (19.792,7.835) -- (19.792,4.32);
\draw[black,fill=white,line width=.55bp] (19.792,39.47) circle[radius=2.834bp];
\draw[black,line width=.55bp] (18.233,39.47) -- (21.351,39.47);
\fill[black] (19.465,10.386) circle[radius=1.15bp];
\node[anchor=north] at (2.217,2.82) {$0$};
\node[anchor=north] at (19.792,2.82) {$0$};
\node[anchor=south] at (2.217,58.545) {$0$};
\node[anchor=south] at (19.792,58.545) {$0$};
\node[anchor=west,font=\fontsize{6.3}{6.8}\selectfont] at (23.292,10.386) {$\textstyle q_V(u,x)M_V(-x)$};
\end{tikzpicture}\right]_{u:-1}
\]
We retain the bubble contribution and the cup--cap term with coefficient
\[
\frac{1-M_V(-x_1)B_V(-x_2)}{x_2-x_1}
\]
for the final sum. The remaining cup--cap term is
\[
\begin{tikzpicture}[x=1bp,y=1bp,scale=.64,
  baseline={([yshift=-.5ex]current bounding box.center)},
  line cap=round,line join=round,
  every node/.style={inner sep=.4pt,outer sep=0pt,font=\fontsize{6.8}{7.3}\selectfont}]
\draw[black,line width=.65bp] (2.217,34.197) .. controls (2.217,16.798) and (19.792,16.798) .. (19.792,34.197);
\draw[black,line width=.65bp] (2.217,4.32) .. controls (2.217,21.719) and (19.792,21.719) .. (19.792,4.32);
\draw[black,densely dotted,line width=.55bp] (17.046,24.41) -- (17.046,14.107);
\node[anchor=north] at (2.217,2.82) {$0$};
\node[anchor=north] at (19.792,2.82) {$0$};
\node[anchor=south] at (2.217,35.697) {$0$};
\node[anchor=south] at (19.792,35.697) {$0$};
\node[anchor=west,font=\fontsize{6.3}{6.8}\selectfont] at (20.605,19.259) {$\textstyle T_B(x_1,x_2)$};
\end{tikzpicture}
\]
where
{\[
\begin{aligned}
T_B(x_1,x_2)=\frac1{x_1+x_2}\Bigg[&1-B_V(x_2)M_V(-x_1)
 +\frac{1-B_V(-x_2)M_V(-x_1)}{-x_1+x_2}\\
&-\frac{(B_V(x_2)-B_V(-x_2))M_V(-x_1)}{2x_2}\Bigg].
\end{aligned}
\]}

\noindent\emph{The \(p_V\) cup--cap terms.}
Similarly, the second product, $\pi(\widehat P\sigma+\widehat P)\pi$, is
\[
\begin{tikzpicture}[x=1bp,y=1bp,scale=.56,
  baseline={([yshift=-.5ex]current bounding box.center)},
  line cap=round,line join=round,
  every node/.style={inner sep=.4pt,outer sep=0pt,font=\fontsize{6.8}{7.3}\selectfont}]
\draw[black,line width=.65bp] (2.217,57.044) -- (2.217,53.529);
\draw[black,line width=.65bp] (2.217,53.529) .. controls (2.217,42.105) and (19.792,42.105) .. (19.792,53.529);
\draw[black,line width=.65bp] (19.792,53.529) -- (19.792,57.044);
\draw[black,line width=.65bp] (2.217,4.319) .. controls (8.075,10.763) and (19.792,14.769) .. (19.792,23.652);
\draw[black,line width=.65bp] (19.792,23.652) .. controls (19.792,33.317) and (2.217,33.317) .. (2.217,23.652);
\draw[black,line width=.65bp] (2.217,23.652) .. controls (2.217,14.769) and (13.934,10.763) .. (19.792,4.319);
\draw[black,densely dotted,line width=.55bp] (17.046,47.103) -- (17.046,29.088);
\draw[black,densely dotted,line width=.55bp] (2.217,54.935) -- (19.792,54.935);
\fill[black] (19.449,50.921) circle[radius=1.15bp];
\node[anchor=west,font=\fontsize{6.3}{6.8}\selectfont] at (23.292,50.921) {$\textstyle B_V(x)$};
\node[anchor=north] at (1.208,3.209) {$0$};
\node[anchor=north] at (20.801,3.209) {$0$};
\node[anchor=south] at (2.217,58.544) {$0$};
\node[anchor=south] at (19.792,58.544) {$0$};
\node[anchor=west,font=\fontsize{6.3}{6.8}\selectfont] at (20.546,38.096) {$\textstyle p_V(x_2,x_1)$};
\node[anchor=west,font=\fontsize{6.3}{6.8}\selectfont] at (23.292,62.421) {$\textstyle 1-x_1+x_2$};
\end{tikzpicture}
+
\begin{tikzpicture}[x=1bp,y=1bp,scale=.56,
  baseline={([yshift=-.5ex]current bounding box.center)},
  line cap=round,line join=round,
  every node/.style={inner sep=.4pt,outer sep=0pt,font=\fontsize{6.8}{7.3}\selectfont}]
\draw[black,line width=.65bp] (2.217,57.044) -- (2.217,53.529);
\draw[black,line width=.65bp] (2.217,53.529) .. controls (2.217,42.105) and (19.792,42.105) .. (19.792,53.529);
\draw[black,line width=.65bp] (19.792,53.529) -- (19.792,57.044);
\draw[black,line width=.65bp] (2.217,4.319) -- (2.217,21.894);
\draw[black,line width=.65bp] (2.217,21.894) .. controls (2.217,33.317) and (19.792,33.317) .. (19.792,21.894);
\draw[black,line width=.65bp] (19.792,21.894) -- (19.792,4.319);
\draw[black,densely dotted,line width=.55bp] (17.046,47.103) -- (17.046,28.319);
\draw[black,densely dotted,line width=.55bp] (2.217,54.935) -- (19.792,54.935);
\fill[black] (19.449,50.921) circle[radius=1.15bp];
\node[anchor=west,font=\fontsize{6.3}{6.8}\selectfont] at (23.292,50.921) {$\textstyle B_V(x)$};
\node[anchor=north] at (2.217,2.819) {$0$};
\node[anchor=north] at (19.792,2.819) {$0$};
\node[anchor=south] at (2.217,58.544) {$0$};
\node[anchor=south] at (19.792,58.544) {$0$};
\node[anchor=west,font=\fontsize{6.3}{6.8}\selectfont] at (20.546,37.711) {$\textstyle p_V(x_2,x_1)$};
\node[anchor=west,font=\fontsize{6.3}{6.8}\selectfont] at (23.292,62.421) {$\textstyle 1-x_1+x_2$};
\end{tikzpicture}
=
\left[\begin{tikzpicture}[x=1bp,y=1bp,scale=.56,
  baseline={([yshift=-.5ex]current bounding box.center)},
  line cap=round,line join=round,
  every node/.style={inner sep=.4pt,outer sep=0pt,font=\fontsize{6.8}{7.3}\selectfont}]
\draw[black,line width=.65bp] (2.217,57.044) -- (2.217,53.529);
\draw[black,line width=.65bp] (2.217,53.529) .. controls (2.217,42.105) and (19.792,42.105) .. (19.792,53.529);
\draw[black,line width=.65bp] (19.792,53.529) -- (19.792,57.044);
\draw[black,line width=.65bp] (2.217,4.319) .. controls (8.075,10.763) and (19.792,14.769) .. (19.792,23.652);
\draw[black,line width=.65bp] (19.792,23.652) .. controls (19.792,33.318) and (2.217,33.318) .. (2.217,23.652);
\draw[black,line width=.65bp] (2.217,23.652) .. controls (2.217,14.769) and (13.934,10.763) .. (19.792,4.319);
\draw[black,fill=white,line width=.55bp] (2.433,21.61) circle[radius=2.834bp];
\draw[black,line width=.55bp] (0.874,21.61) -- (3.991,21.61);
\fill[black] (19.465,50.978) circle[radius=1.15bp];
\node[anchor=north] at (1.208,3.209) {$0$};
\node[anchor=north] at (20.801,3.209) {$0$};
\node[anchor=south] at (2.217,58.544) {$0$};
\node[anchor=south] at (19.792,58.544) {$0$};
\node[anchor=west,font=\fontsize{6.3}{6.8}\selectfont] at (23.292,50.978) {$\textstyle W_V(-x)p_V(x,u)$};
\end{tikzpicture}
+
\begin{tikzpicture}[x=1bp,y=1bp,scale=.56,
  baseline={([yshift=-.5ex]current bounding box.center)},
  line cap=round,line join=round,
  every node/.style={inner sep=.4pt,outer sep=0pt,font=\fontsize{6.8}{7.3}\selectfont}]
\draw[black,line width=.65bp] (2.217,57.045) -- (2.217,53.53);
\draw[black,line width=.65bp] (2.217,53.53) .. controls (2.217,42.106) and (19.792,42.106) .. (19.792,53.53);
\draw[black,line width=.65bp] (19.792,53.53) -- (19.792,57.045);
\draw[black,line width=.65bp] (2.217,4.32) -- (2.217,21.895);
\draw[black,line width=.65bp] (2.217,21.895) .. controls (2.217,33.319) and (19.792,33.319) .. (19.792,21.895);
\draw[black,line width=.65bp] (19.792,21.895) -- (19.792,4.32);
\draw[black,fill=white,line width=.55bp] (19.792,21.895) circle[radius=2.834bp];
\draw[black,line width=.55bp] (18.233,21.895) -- (21.351,21.895);
\fill[black] (19.465,50.979) circle[radius=1.15bp];
\node[anchor=north] at (2.217,2.82) {$0$};
\node[anchor=north] at (19.792,2.82) {$0$};
\node[anchor=south] at (2.217,58.545) {$0$};
\node[anchor=south] at (19.792,58.545) {$0$};
\node[anchor=west,font=\fontsize{6.3}{6.8}\selectfont] at (23.292,50.979) {$\textstyle W_V(-x)p_V(x,u)$};
\end{tikzpicture}\right]_{u:-1}
\]
We retain the bubble contribution and the cup--cap term with coefficient
\[
\frac{W_V(-x_2)A_V(-x_1)-1}{x_2-x_1}
\]
for the final sum. The remaining cup--cap term is
\[
\begin{tikzpicture}[x=1bp,y=1bp,scale=.64,
  baseline={([yshift=-.5ex]current bounding box.center)},
  line cap=round,line join=round,
  every node/.style={inner sep=.4pt,outer sep=0pt,font=\fontsize{6.8}{7.3}\selectfont}]
\draw[black,line width=.65bp] (2.217,34.197) .. controls (2.217,16.798) and (19.792,16.798) .. (19.792,34.197);
\draw[black,line width=.65bp] (2.217,4.32) .. controls (2.217,21.719) and (19.792,21.719) .. (19.792,4.32);
\draw[black,densely dotted,line width=.55bp] (17.046,24.41) -- (17.046,14.107);
\node[anchor=north] at (2.217,2.82) {$0$};
\node[anchor=north] at (19.792,2.82) {$0$};
\node[anchor=south] at (2.217,35.697) {$0$};
\node[anchor=south] at (19.792,35.697) {$0$};
\node[anchor=west,font=\fontsize{6.3}{6.8}\selectfont] at (20.605,19.259) {$\textstyle T_A(x_1,x_2)$};
\end{tikzpicture}
\]
where
{\[
\begin{aligned}
T_A(x_1,x_2)=\frac1{x_1+x_2}\Bigg[&A_V(x_1)W_V(-x_2)-1
 -\frac{(A_V(x_1)-A_V(-x_1))W_V(-x_2)}{2x_1}\\
&+\frac{1-A_V(-x_1)W_V(-x_2)}{x_1-x_2}\Bigg].
\end{aligned}
\]}

\noindent\emph{The bubble terms.}
{Use the full black bubble, with all colors summed, and its series
$b_V(u)$ from Lemma~\ref{nbone:lemma:nil2}:}
\[
b_V(u)
{=\frac1u\sum_{r\ge0}\begin{tikzpicture}[x=1bp,y=1bp,scale=.86,
baseline={([yshift=-.5ex]current bounding box.center)},
every node/.style={inner sep=1pt,outer sep=0pt,font=\scriptsize}]
\draw[black,line width=.65bp] (8,11) ellipse[x radius=8bp,y radius=11bp];
\fill[black] (16,11) circle[radius=1.2bp];
\node[anchor=west] at (20,11) {$\frac{x^r}{u^r}$};
\draw[darkg,line width=.65bp] (48,0)--(48,22);
\node[text=darkg,anchor=north] at (48,-2) {$V$};
\end{tikzpicture}}
=
\frac{n_{V,0}(u)}{m_{V,0}(u)}
 \frac{t_{V,0}(u)}{u^{e+1}}-1+\frac1{2u}.
\]
The product $\widehat P\widehat Q$ and the two retained bubble
contributions have the common factor
$b_V(u)/((u+x_1)(u+x_2))$. Their remaining factors sum to
\[
\begin{aligned}
&\bigl(W_V(-x_2)A_V(u)-1\bigr)
 \bigl(M_V(-x_1)B_V(u)-1\bigr)\\
&\qquad{}+M_V(-x_1)B_V(u)+W_V(-x_2)A_V(u)\\
&=M_V(-x_1)W_V(-x_2)A_V(u)B_V(u)+1.
\end{aligned}
\]
The constant term contributes no $u^{-1}$-coefficient, since
$b_V(u)=O(u^{-1})$ at infinity.
For the support-preserving polynomial representatives used here,
\[
\frac{n_{V,0}(u)}{m_{V,0}(u)}
\frac{t_{V,0}(u)}{u^{e+1}}A_V(u)B_V(u)
=
\frac{(-1)^{c+1}}2\frac{t_{V,0}(u)}{u^{e+1}}
+u^{K-e-1}t_{V,0}(u)H_V(u).
\]
The leading term has order $u^{-1}$; after the two strand resolvents in the
bubble terms it starts at order $u^{-3}$ and therefore has no
$u^{-1}$-coefficient.  {For the polynomial tail, let $I$ be the ideal generated by the commuting nilpotent dots. If $I^N=0$, choosing} $K\ge e+N+2$ makes every resulting first
divided difference vanish after nilpotent evaluation.

The $1/(2u)$ part of the bubble series gives
\[
\left[\frac{1}{2u}A_V(u)B_V(u)\begin{tikzpicture}[x=1bp,y=1bp,scale=.56,
  baseline={([yshift=-.5ex]current bounding box.center)},
  line cap=round,line join=round,
  every node/.style={inner sep=.4pt,outer sep=0pt,font=\fontsize{6.8}{7.3}\selectfont}]
\draw[black,line width=.65bp] (2.217,57.044) -- (2.217,53.529);
\draw[black,line width=.65bp] (2.217,53.529) .. controls (2.217,42.105) and (19.792,42.105) .. (19.792,53.529);
\draw[black,line width=.65bp] (19.792,53.529) -- (19.792,57.044);
\draw[black,line width=.65bp] (2.217,4.319) -- (2.217,7.834);
\draw[black,line width=.65bp] (2.217,7.834) .. controls (2.217,19.258) and (19.792,19.258) .. (19.792,7.834);
\draw[black,line width=.65bp] (19.792,7.834) -- (19.792,4.319);
\draw[black,fill=white,line width=.55bp] (2.459,51.315) circle[radius=2.834bp];
\draw[black,line width=.55bp] (0.901,51.315) -- (4.018,51.315);
\draw[black,fill=white,line width=.55bp] (2.459,10.048) circle[radius=2.834bp];
\draw[black,line width=.55bp] (0.901,10.048) -- (4.018,10.048);
\fill[black] (19.465,50.978) circle[radius=1.15bp];
\fill[black] (19.465,10.385) circle[radius=1.15bp];
\node[anchor=west,font=\fontsize{6.3}{6.8}\selectfont] at (23.292,50.978) {$\textstyle W_V(-x)$};
\node[anchor=west,font=\fontsize{6.3}{6.8}\selectfont] at (23.292,10.385) {$\textstyle M_V(-x)$};
\node[anchor=north] at (2.217,2.819) {$0$};
\node[anchor=north] at (19.792,2.819) {$0$};
\node[anchor=south] at (2.217,58.544) {$0$};
\node[anchor=south] at (19.792,58.544) {$0$};
\end{tikzpicture}\right]_{u:-1}
=
\begin{tikzpicture}[x=1bp,y=1bp,scale=.64,
  baseline={([yshift=-.5ex]current bounding box.center)},
  line cap=round,line join=round,
  every node/.style={inner sep=.4pt,outer sep=0pt,font=\fontsize{6.8}{7.3}\selectfont}]
\draw[black,line width=.65bp] (2.217,34.197) .. controls (2.217,16.798) and (19.792,16.798) .. (19.792,34.197);
\draw[black,line width=.65bp] (2.217,4.32) .. controls (2.217,21.719) and (19.792,21.719) .. (19.792,4.32);
\draw[black,densely dotted,line width=.55bp] (17.046,24.41) -- (17.046,14.107);
\node[anchor=north] at (2.217,2.82) {$0$};
\node[anchor=north] at (19.792,2.82) {$0$};
\node[anchor=south] at (2.217,35.697) {$0$};
\node[anchor=south] at (19.792,35.697) {$0$};
\node[anchor=west,font=\fontsize{6.3}{6.8}\selectfont] at (20.605,19.259) {$\textstyle T_0(x_1,x_2)$};
\end{tikzpicture}
\]
where $T_0(x_1,x_2)
=\frac{M_V(-x_1)W_V(-x_2)}{-x_1+x_2}
\left(
\frac{A_V(-x_1)B_V(-x_1)-1}{-2x_1}
-\frac{A_V(-x_2)B_V(-x_2)-1}{-2x_2}
\right)$.

The $-1$ part of the bubble series and the terms retained from the
$q_V$ and $p_V$ calculations give, respectively,
\[
\left[-A_V(u)B_V(u)\begin{tikzpicture}[x=1bp,y=1bp,scale=.56,
  baseline={([yshift=-.5ex]current bounding box.center)},
  line cap=round,line join=round,
  every node/.style={inner sep=.4pt,outer sep=0pt,font=\fontsize{6.8}{7.3}\selectfont}]
\draw[black,line width=.65bp] (2.217,57.044) -- (2.217,53.529);
\draw[black,line width=.65bp] (2.217,53.529) .. controls (2.217,42.105) and (19.792,42.105) .. (19.792,53.529);
\draw[black,line width=.65bp] (19.792,53.529) -- (19.792,57.044);
\draw[black,line width=.65bp] (2.217,4.319) -- (2.217,7.834);
\draw[black,line width=.65bp] (2.217,7.834) .. controls (2.217,19.258) and (19.792,19.258) .. (19.792,7.834);
\draw[black,line width=.65bp] (19.792,7.834) -- (19.792,4.319);
\draw[black,fill=white,line width=.55bp] (2.459,51.315) circle[radius=2.834bp];
\draw[black,line width=.55bp] (0.901,51.315) -- (4.018,51.315);
\draw[black,fill=white,line width=.55bp] (2.459,10.048) circle[radius=2.834bp];
\draw[black,line width=.55bp] (0.901,10.048) -- (4.018,10.048);
\fill[black] (19.465,50.978) circle[radius=1.15bp];
\fill[black] (19.465,10.385) circle[radius=1.15bp];
\node[anchor=west,font=\fontsize{6.3}{6.8}\selectfont] at (23.292,50.978) {$\textstyle W_V(-x)$};
\node[anchor=west,font=\fontsize{6.3}{6.8}\selectfont] at (23.292,10.385) {$\textstyle M_V(-x)$};
\node[anchor=north] at (2.217,2.819) {$0$};
\node[anchor=north] at (19.792,2.819) {$0$};
\node[anchor=south] at (2.217,58.544) {$0$};
\node[anchor=south] at (19.792,58.544) {$0$};
\end{tikzpicture}\right]_{u:-1}
+
\left[B_V(u)\begin{tikzpicture}[x=1bp,y=1bp,scale=.56,
  baseline={([yshift=-.5ex]current bounding box.center)},
  line cap=round,line join=round,
  every node/.style={inner sep=.4pt,outer sep=0pt,font=\fontsize{6.8}{7.3}\selectfont}]
\draw[black,line width=.65bp] (2.217,57.045) -- (2.217,39.47);
\draw[black,line width=.65bp] (2.217,39.47) .. controls (2.217,28.046) and (19.792,28.046) .. (19.792,39.47);
\draw[black,line width=.65bp] (19.792,39.47) -- (19.792,57.045);
\draw[black,line width=.65bp] (2.217,4.32) -- (2.217,7.835);
\draw[black,line width=.65bp] (2.217,7.835) .. controls (2.217,19.259) and (19.792,19.259) .. (19.792,7.835);
\draw[black,line width=.65bp] (19.792,7.835) -- (19.792,4.32);
\draw[black,fill=white,line width=.55bp] (2.217,39.47) circle[radius=2.834bp];
\draw[black,line width=.55bp] (0.658,39.47) -- (3.776,39.47);
\draw[black,fill=white,line width=.55bp] (2.459,10.049) circle[radius=2.834bp];
\draw[black,line width=.55bp] (0.901,10.049) -- (4.018,10.049);
\fill[black] (19.465,10.386) circle[radius=1.15bp];
\node[anchor=west,font=\fontsize{6.3}{6.8}\selectfont] at (23.292,10.386) {$\textstyle M_V(-x)$};
\node[anchor=north] at (2.217,2.82) {$0$};
\node[anchor=north] at (19.792,2.82) {$0$};
\node[anchor=south] at (2.217,58.545) {$0$};
\node[anchor=south] at (19.792,58.545) {$0$};
\end{tikzpicture}\right]_{u:-1}
+
\left[A_V(u)\begin{tikzpicture}[x=1bp,y=1bp,scale=.56,
  baseline={([yshift=-.5ex]current bounding box.center)},
  line cap=round,line join=round,
  every node/.style={inner sep=.4pt,outer sep=0pt,font=\fontsize{6.8}{7.3}\selectfont}]
\draw[black,line width=.65bp] (2.217,57.044) -- (2.217,53.529);
\draw[black,line width=.65bp] (2.217,53.529) .. controls (2.217,42.105) and (19.792,42.105) .. (19.792,53.529);
\draw[black,line width=.65bp] (19.792,53.529) -- (19.792,57.044);
\draw[black,line width=.65bp] (2.217,4.319) -- (2.217,21.894);
\draw[black,line width=.65bp] (2.217,21.894) .. controls (2.217,33.317) and (19.792,33.317) .. (19.792,21.894);
\draw[black,line width=.65bp] (19.792,21.894) -- (19.792,4.319);
\draw[black,fill=white,line width=.55bp] (2.217,21.894) circle[radius=2.834bp];
\draw[black,line width=.55bp] (0.658,21.894) -- (3.776,21.894);
\draw[black,fill=white,line width=.55bp] (2.459,51.315) circle[radius=2.834bp];
\draw[black,line width=.55bp] (0.901,51.315) -- (4.018,51.315);
\fill[black] (19.465,50.978) circle[radius=1.15bp];
\node[anchor=west,font=\fontsize{6.3}{6.8}\selectfont] at (23.292,50.978) {$\textstyle W_V(-x)$};
\node[anchor=north] at (2.217,2.819) {$0$};
\node[anchor=north] at (19.792,2.819) {$0$};
\node[anchor=south] at (2.217,58.544) {$0$};
\node[anchor=south] at (19.792,58.544) {$0$};
\end{tikzpicture}\right]_{u:-1}
\]
Their sum vanishes, since the numerator of their common coefficient is
\[
\begin{aligned}
&-A_V(-x_1)W_V(-x_2)+B_V(-x_2)M_V(-x_1)\\
&\qquad{}+1-B_V(-x_2)M_V(-x_1)+A_V(-x_1)W_V(-x_2)-1=0.
\end{aligned}
\]

Thus the sum of the three projected products is
\[
sQ-Q+Ps+P+PQ-D=\begin{tikzpicture}[x=1bp,y=1bp,scale=.64,
  baseline={([yshift=-.5ex]current bounding box.center)},
  line cap=round,line join=round,
  every node/.style={inner sep=.4pt,outer sep=0pt,font=\fontsize{6.8}{7.3}\selectfont}]
\draw[black,line width=.65bp] (2.217,34.196) .. controls (2.217,16.797) and (19.792,16.797) .. (19.792,34.196);
\draw[black,line width=.65bp] (2.217,4.319) .. controls (2.217,21.718) and (19.792,21.718) .. (19.792,4.319);
\draw[black,densely dotted,line width=.55bp] (17.046,24.409) -- (17.046,14.106);
\node[anchor=north] at (2.217,2.819) {$0$};
\node[anchor=north] at (19.792,2.819) {$0$};
\node[anchor=south] at (2.217,35.696) {$0$};
\node[anchor=south] at (19.792,35.696) {$0$};
\end{tikzpicture}\,
\left(T_0(x_1,x_2)+T_B(x_1,x_2)+T_A(x_1,x_2)\right).
\]
Using the relations $A_V(-x_2)W_V(-x_2)=1=B_V(-x_1)M_V(-x_1)$ and $A_V(-u)B_V(u)=\frac1{1+2u}$,
we get
\[
T_0(x_1,x_2) + T_B(x_1,x_2) + T_A(x_1,x_2) = 0.
\]
It follows that
\begin{equation}\label{eq:nil-square-defect}
S_V=(s^2-\pi)+D=\pi(\sigma^2-\operatorname{Id})\pi.
\end{equation}
No quadratic relation for $\sigma$ was used in this calculation.
Now $\sigma^2=\operatorname{Id}$ gives $S_V=0$, proving
relation~\textnormal{(7)} of Definition~\ref{nil}.
\end{proof}

\section{The nil-Brauer braid relation}
\label{sec:nil-brauer-braid}

{We prove relation~\textnormal{(8)} of Definition~\ref{nil}.} Let
\(\mathbf L\) and \(\mathbf R\) denote the left- and right-hand sides of the
red braid relation, evaluated on a finitely
generated object \(V\in\cM_\lambda\).
If $E_0V=0$, the relation is immediate; hence assume $E_0V\ne0$.

{Use the three-term expansions in Lemma~\ref{lem:anothdefofcro},
labelling the black-crossing, parallel-strand, and cup--cap terms by
$1,2,3$. Use the lower formula for the top crossing and the upper
formula for the other two. As explained at the beginning of
Appendix~\ref{app:graphical-certificates}, the complete expansions
allow projection only at the outer boundary. Thus
\[
\mathbf L=\sum_{a,b,c\in\{1,2,3\}}L_{abc},\qquad
\mathbf R=\sum_{a,b,c\in\{1,2,3\}}R_{abc},
\]
with source diagrams as in Appendix~\ref{27terms-app}.
These terms are evaluated using compatible finite polynomial
representatives. The middle crossing and parallel terms are paired
before rational simplification; individual internal factors
$1+x_a-x_b$ need not be invertible.}
 Set $l_{a,b}=1+x_a-x_b$ and
\[
 \widetilde L_{abc}
 =
 l_{2,1}l_{3,2}\,
 L_{abc}\,
 l_{1,2}l_{1,3}l_{2,3},
 \qquad
 \widetilde R_{abc}
 =
 l_{2,1}l_{3,2}\,
 R_{abc}\,
 l_{1,2}l_{1,3}l_{2,3}.
\]
We write \(\widetilde{\mathbf L}\) and \(\widetilde{\mathbf R}\) for the
corresponding sums. These factors act on the outer copy of $E_0^3V$,
where all three dots are nilpotent. Each factor is therefore invertible,
so it is enough to prove
\(\widetilde{\mathbf L}=\widetilde{\mathbf R}\).

We organize the calculation by the number of black crossings. Here $s_i$
is a full black crossing and $e_i$ a cup--cap on positions $i,i+1$.
Words in these symbols denote full diagrams projected only at their
outer boundary. We use the following fifteen representatives:
\[
\begin{array}{c|l}
3 & s_1s_2s_1\\
2 & s_1s_2,\ s_2s_1,\ s_1e_2s_1=s_2e_1s_2\\
1 & s_1,\ s_2,\ s_1e_2,\ s_2e_1,\ e_1s_2,\ e_2s_1\\
0 & 1,\ e_1,\ e_2,\ e_1e_2,\ e_2e_1.
\end{array}
\]
Let \(\mathcal F_{\le d}\) be the span of the terms with at most \(d\)
black crossings.  Moving a dot-polynomial through a black crossing gives
the same-crossing leading term plus a divided-difference correction with
one fewer crossing; moving through a cup or cap does not increase the
crossing number.  Hence the comparison is triangular:
\[
3\longrightarrow2\longrightarrow1\longrightarrow0.
\]
All lower-order corrections produced at one stage are retained for the
next stage. We compare the coefficient functions in these explicit
reductions; linear independence of the operators on $V$ is not needed.
{The graphical reductions and one-crossing comparisons are in
Appendix~\ref{app:graphical-certificates}; the zero-crossing
coefficient calculations are in Appendix~\ref{app:coefficient-certificates}.}

{For each representative $\tau$, write
\begin{equation}\label{eq:normal-form-symbol}
[\tau;F]
\end{equation}
for decoration by $F(x_1,x_2,x_3)$ at the marks specified in that
comparison. The positions, including those below a crossing, determine
the variables; they are not relabelled under reflection.
Let $\Pi$ project onto $E_0^3V$. For example,
\[
[s_1s_2s_1;1]=\Pi s_1s_2s_1\Pi,
\]
where the crossings on the right are full black crossings.
}
\subsubsection*{Three crossings}

{Only $\widetilde L_{111}$ and $\widetilde R_{111}$ contribute
modulo $\mathcal F_{\le2}$. By \eqref{rx1x2}, moving the decorations
to the upper marked points gives crossing factor
$(l_{2,1}l_{3,1}l_{3,2})^{-1}$ on either side. The lower normalizing
factor becomes $l_{3,2}l_{3,1}l_{2,1}$, leaving $l_{2,1}l_{3,2}$.
Thus the full black braid relation gives
\[
\widetilde L_{111}\equiv
\begin{tikzpicture}[x=1bp,y=1bp,scale=.72,baseline=24bp,
line cap=round,line join=round,draw=black,
every node/.style={inner sep=1pt,outer sep=0pt,font=\scriptsize}]
\draw[black,line width=.75bp] (0,0)--(0,4)--(20,16)--(20,24)--(40,36)--(40,44)--(40,56)--(40,76);
\draw[black,line width=.75bp] (20,0)--(20,4)--(0,16)--(0,24)--(0,36)--(0,44)--(20,56)--(20,76);
\draw[black,line width=.75bp] (40,0)--(40,4)--(40,16)--(40,24)--(20,36)--(20,44)--(0,56)--(0,76);
\fill[black] (0,69) circle[radius=1.2bp];
\node[anchor=south] at (0,78) {$ 1 $};
\fill[black] (20,69) circle[radius=1.2bp];
\node[anchor=south] at (20,78) {$ 2 $};
\fill[black] (40,69) circle[radius=1.2bp];
\node[anchor=south] at (40,78) {$ 3 $};
\draw[black,densely dotted,line width=.55bp] (0,69)--(40,69);
\node[anchor=west] at (46,69) {$l_{2,1}l_{3,2}$};
\end{tikzpicture}\;=\;\begin{tikzpicture}[x=1bp,y=1bp,scale=.72,baseline=24bp,
line cap=round,line join=round,draw=black,
every node/.style={inner sep=1pt,outer sep=0pt,font=\scriptsize}]
\draw[black,line width=.75bp] (0,0)--(0,4)--(0,16)--(0,24)--(20,36)--(20,44)--(40,56)--(40,76);
\draw[black,line width=.75bp] (20,0)--(20,4)--(40,16)--(40,24)--(40,36)--(40,44)--(20,56)--(20,76);
\draw[black,line width=.75bp] (40,0)--(40,4)--(20,16)--(20,24)--(0,36)--(0,44)--(0,56)--(0,76);
\fill[black] (0,69) circle[radius=1.2bp];
\node[anchor=south] at (0,78) {$ 1 $};
\fill[black] (20,69) circle[radius=1.2bp];
\node[anchor=south] at (20,78) {$ 2 $};
\fill[black] (40,69) circle[radius=1.2bp];
\node[anchor=south] at (40,78) {$ 3 $};
\draw[black,densely dotted,line width=.55bp] (0,69)--(40,69);
\node[anchor=west] at (46,69) {$l_{2,1}l_{3,2}$};
\end{tikzpicture}
\equiv\widetilde R_{111}
\pmod{\mathcal F_{\le2}}.
\]
The labels $1,2,3$ specify the upper dot positions. All terms created
by moving the decorations are retained in the seven-term expansions
of Appendix~\ref{app:three-crossing-full}. Hence
\[
\widetilde{\mathbf L}\equiv\widetilde{\mathbf R}
\pmod{\mathcal F_{\le2}}.
\]
}
\subsubsection*{Two crossings}
For each representative $d$ above, denote its coefficient functions
in the reductions of $\widetilde{\mathbf L}$ and $\widetilde{\mathbf R}$ by
$C_{d,L}$ and $C_{d,R}$, respectively, using the indicated dot positions.
Modulo \(\mathcal F_{\le1}\), there are three Brauer diagrams.
For \(d=s_2s_1\),
\[
 C_{d,L}
 =2l_{3,2}-l_{1,2}l_{3,2}
 =l_{2,1}l_{3,2}
 =C_{d,R},
\]
and for \(d=s_1s_2\),
\[
 C_{d,L}
 =l_{2,1}l_{3,2}
 =2l_{2,1}-l_{2,1}l_{2,3}
 =C_{d,R}.
\]
The corresponding graphical reductions are in
Appendix~\ref{app:two-crossing-full}.

It remains to compare
\[
 d=s_1e_2s_1=s_2e_1s_2.
\]
{The marked positions are as follows: $x_1$ is on the through
strand, $x_2$ on the right leg of the upper cup, and $x_3$ on the
right leg of the lower cap.
\[
\begin{tikzpicture}[scale=.8,baseline=0pt,line cap=round,line join=round,
every node/.style={inner sep=1pt,outer sep=0pt,font=\scriptsize}]
\draw[black,thick] (-.6,.9) .. controls (-.6,.55) and (-.3,.2) .. (0,.2)
 .. controls (.3,.2) and (.6,.55) .. (.6,.9);
\draw[black,thick] (-.6,-.9) .. controls (-.6,-.55) and (-.3,-.2) .. (0,-.2)
 .. controls (.3,-.2) and (.6,-.55) .. (.6,-.9);
\draw[black,thick] (0,-.9)--(0,.9);
\fill[black] (0,.65) circle[radius=1.2pt];
\fill[black] (.4125,.41875) circle[radius=1.2pt];
\fill[black] (.4125,-.41875) circle[radius=1.2pt];
\node[anchor=east] at (-.08,.65) {$1$};
\node[anchor=west] at (.5,.41875) {$2$};
\node[anchor=west] at (.5,-.41875) {$3$};
\end{tikzpicture}
\]
}
For this comparison write $N=E_0V$. Coefficients indexed by $N$ are evaluated on that object,
and those indexed by $V$ are transported coefficientwise. The two coefficients are
\begin{align}
 C_{d,L}
 &=
 -(1-x_1+x_2)
 +(1-x_1+x_2)(1+x_1-x_3)
 q_V(x_2,x_3)M_V(-x_3),\label{eq:braid-degree2-L}\\
 C_{d,R}
 &=
 -(1+x_1+x_2)
 +(1+x_1+x_2)(1-x_1-x_3)
 q_N(x_2,x_3)M_N(-x_3).
 \label{eq:braid-degree2-R}
\end{align}
{The primary-factor calculation in the proof of
Lemma~\ref{eq:red-cup-E0V}, together with Lemma~\ref{gamma-move1},
gives the following identity for the zero-color series:}
\[
B_N(u)=\frac{1+u-x_1}{1+u+x_1}\,E_0(B_V(u)).
\]
Substituting this into
\[
q_X(x_2,x_3)M_X(-x_3)
=\frac{B_X(x_2)/B_X(-x_3)-1}{x_2+x_3},\qquad X\in\{V,N\},
\]
{the difference between \eqref{eq:braid-degree2-R} and
\eqref{eq:braid-degree2-L} is
\[
-2x_1+
\frac{(1-x_1+x_2)(1+x_1-x_3)-(1+x_1+x_2)(1-x_1-x_3)}{x_2+x_3}
=0.
\]
The numerator is $2x_1(x_2+x_3)$. Thus no inverse of the nilpotent
operator $x_2+x_3$ is used; divide before evaluation, as in
Lemma~\ref{lem:finite-dd-evaluation}. Hence}
\[
 \widetilde{\mathbf L}\equiv
 \widetilde{\mathbf R}\pmod{\mathcal F_{\le1}}.
\]

\subsubsection*{One crossing}

Modulo \(\mathcal F_{\le0}\), we compare
\[
 s_1,\quad s_2,\quad s_1e_2,\quad s_2e_1,\quad e_1s_2,\quad e_2s_1.
\]
{For $s_1$ and $s_2$, use the marked positions in
Lemma~\ref{lem:appA-one-crossing-ss}, retaining their labels under reflection.}
For \(d=s_2\),
\[
 C_{d,L}
 =2l_{3,2}-l_{3,2}l_{1,2}
 =l_{2,1}l_{3,2}
 =l_{2,1}(1-l_{1,3}+l_{1,2})
 =C_{d,R},
\]
while for \(d=s_1\),
\[
 C_{d,L}=l_{1,2}l_{2,3}=C_{d,R}.
\]
For the four diagrams
\[
 e_1s_2,\qquad e_2s_1,\qquad s_1e_2,\qquad s_2e_1,
\]
the source reductions and their higher-crossing corrections are given in
Appendices~\ref{app:one-crossing-sources} and~\ref{app:mixed-one-crossing-graphics}.
The coefficient comparisons are given in Appendix~\ref{app:one-crossing-comparison},
Lemmas~\ref{lem:appA-one-crossing-es} and~\ref{lem:appA-one-crossing-se}. Thus
\[
 \widetilde{\mathbf L}\equiv
 \widetilde{\mathbf R}\pmod{\mathcal F_{\le0}}.
\]

\subsubsection*{Zero crossings}

It remains to compare
\[
 1,\qquad e_1,\qquad e_2,\qquad e_1e_2,\qquad e_2e_1.
\]
For the identity diagram,
\[
 C_{1,L}=l_{2,1}l_{3,2}=C_{1,R}.
\]

For the remaining four Brauer diagrams, we first collect the direct
zero-crossing terms and all lower-order corrections retained above.
{The graphical reductions are in
Appendix~\ref{app:zero-crossing-graphics}. The coefficients for
$e_1e_2$ and $e_2e_1$ agree by Lemmas~\ref{CRCLE} and
\ref{lem:e2e1-coefficients}. For $e_1$ and $e_2$, include the bubbles
before comparing coefficients: Lemmas~\ref{lem:e1-coefficients} and
\ref{lem:e2-coefficients} give, respectively,
\[
\begin{aligned}
C_{e_1,L}-C_{e_1,R}&=L_3+B_t-R_3=0,\\
C_{e_2,L}-C_{e_2,R}&=L_4-R_4-B_r=0.
\end{aligned}
\]
Here $L_3,R_3,L_4,R_4$ are the terms without bubbles, and $B_t,B_r$
are the bubble contributions defined there. Finite polynomial
representatives are retained until coefficient extraction is complete.
}
Consequently
\[
 \widetilde{\mathbf L}=\widetilde{\mathbf R},
\]
and hence \(\mathbf L=\mathbf R\). This proves relation~\textnormal{(8)}.

\section{Module categories over the affine Brauer category}
\label{sec:integral-reverse-complete}
Assume throughout that $I=\mathbb Z$, $I^\imath=\mathbb Z_{>0}$,
and $\mathcal M=(\mathcal M_\lambda)_{\lambda\in X^\imath}$ is a
locally Schurian nilpotent $2$-representation of the integral
$\imath$-Kac--Moody $2$-category. We write $E_{-i}=F_i$ for $i>0$;
the self-adjoint functor $E_0$ occurs once in $E=\bigoplus_{i\in I}E_i$.
We work over $\mathbb C$, with
$\lambda_i=\langle\alpha_i^{\theta,\vee},\lambda\rangle$ for $i>0$
as in \eqref{lambdacoe}; thus
$E_i\colon\mathcal M_\lambda\to\mathcal M_{\lambda-\bar\alpha_i}$.
{Each $E_i$ preserves finitely generated objects and filtered colimits
by biadjunction. On $\mathcal M=\prod_\lambda\mathcal M_\lambda$, set
\[
 (EV)_\mu=\bigoplus_{i\in I}E_iV_{\mu+\bar\alpha_i}.
\]
For a finitely generated projective $P\in\mathcal M_\mu$, adjunction
gives $\Hom(P,E_iV_{\mu+\bar\alpha_i})\cong
\Hom(E_{-i}P,V_{\mu+\bar\alpha_i})$. These spaces are finite-dimensional
and vanish for all but finitely many $i$ by
Definition~\ref{def:integral-nilpotence}. Thus $E$ is well defined.
On finitely generated objects, $E$ and each fixed iterate are finite
sums. As in \cite[Section~7.1]{GRS}, we define the transformations
there and extend by filtered colimits.
}\newcommand{\IRstrand}[2]{\begin{tikzpicture}[baseline=-1mm,scale=.7]
 \draw[#1] (0,-.4)--(0,.5);\node at (0,-.6){$\scriptstyle #2$};\end{tikzpicture}}
\newcommand{\IRcross}[3]{\begin{tikzpicture}[baseline=0,scale=.7]
 \draw[#1] (-.38,-.4)--(.38,.5) (.38,-.4)--(-.38,.5);
 \node at (-.4,-.6){$\scriptstyle #2$};\node at (.4,-.6){$\scriptstyle #3$};\end{tikzpicture}}
\newcommand{\IRparallel}[3]{\begin{tikzpicture}[baseline=0,scale=.7]
 \draw[#1] (-.3,-.4)--(-.3,.5) (.3,-.4)--(.3,.5);
 \node at (-.3,-.6){$\scriptstyle #2$};\node at (.3,-.6){$\scriptstyle #3$};
 \node at (-.3,.7){$\scriptstyle #2$};\node at (.3,.7){$\scriptstyle #3$};\end{tikzpicture}}

\newcommand{\IRcap}[3]{\begin{tikzpicture}[baseline=0,scale=.7]
 \draw[#1] (-.4,-.4)..controls(-.4,.5)and(.4,.5)..(.4,-.4);
 \node at (-.4,-.6){$\scriptstyle #2$};\node at (.4,-.6){$\scriptstyle #3$};\end{tikzpicture}}
\newcommand{\IRcup}[3]{\begin{tikzpicture}[baseline=0,scale=.7]
 \draw[#1] (-.4,.4)..controls(-.4,-.5)and(.4,-.5)..(.4,.4);
 \node at (-.4,.6){$\scriptstyle #2$};\node at (.4,.6){$\scriptstyle #3$};\end{tikzpicture}}
\newcommand{\IRsnake}[3]{\begin{tikzpicture}[baseline=0,scale=.65,xscale=#3]
 \draw[#1] (-.65,.5)..controls(-.65,-.7)and(0,-.7)..(0,0)
 ..controls(0,.7)and(.65,.7)..(.65,-.5);
 \node at (.65,-.7){$\scriptstyle #2$};\end{tikzpicture}}
\newcommand{\IRcapcross}[4]{\begin{tikzpicture}[baseline=0,scale=.7]
 \begin{scope}[xscale=#4]
  \draw[#1,thick] (-.58,-.48)
   ..controls(-.58,.50)and(.58,.50)..(.58,-.48);
  \draw[#1,thick] (0,-.48) to[out=90,in=-90] (-.82,.64);
  \node at (-.82,.84){$\scriptstyle #3$};
 \end{scope}
 \node at (-.58,-.70){$\scriptstyle #2$};
 \node at (0,-.70){$\scriptstyle #3$};
 \node at (.58,-.70){$\scriptstyle #2$};\end{tikzpicture}}
\newcommand{\IRdouble}[3]{\begin{tikzpicture}[baseline=0,scale=.7]
 \draw[#1,thick] (-.40,-.62) to[out=90,in=-90] (.40,0)
  to[out=90,in=-90] (-.40,.62);
 \draw[#1,thick] (.40,-.62) to[out=90,in=-90] (-.40,0)
  to[out=90,in=-90] (.40,.62);
 \node at (-.40,-.82){$\scriptstyle #2$};\node at (.40,-.82){$\scriptstyle #3$};
 \node at (-.40,.82){$\scriptstyle #2$};\node at (.40,.82){$\scriptstyle #3$};\end{tikzpicture}}

\newcommand{\IRbraid}[5]{\begin{tikzpicture}[baseline=0,scale=.62]
 \begin{scope}[xscale=#5]
  \draw[#1,thick] (-.72,-.78) to[out=90,in=-90] (.72,.78);
  \draw[#1,thick] (.72,-.78) to[out=90,in=-90] (-.72,.78);
  \draw[#1,thick] (0,-.78) to[out=90,in=-90] (-.72,0)
   to[out=90,in=-90] (0,.78);
 \end{scope}
 \node at (-.72,-1.02){$\scriptstyle #2$};\node at (0,-1.02){$\scriptstyle #3$};
 \node at (.72,-1.02){$\scriptstyle #4$};
 \node at (-.72,1.02){$\scriptstyle #4$};\node at (0,1.02){$\scriptstyle #3$};
 \node at (.72,1.02){$\scriptstyle #2$};\end{tikzpicture}}

\newcommand{\IRcrossdot}[4]{\begin{tikzpicture}[baseline=0,scale=.7]
 \draw[#1] (-.4,-.5)--(.4,.5) (.4,-.5)--(-.4,.5);
 \fill[#1] #4 circle (1.7pt);
 \node at (-.4,-.7){$\scriptstyle #2$};\node at (.4,-.7){$\scriptstyle #3$};
\end{tikzpicture}}

\newcommand{\IRdotcap}[4]{\begin{tikzpicture}[baseline=0,scale=.7]
 \draw[#1] (-.4,-.4)..controls(-.4,.5)and(.4,.5)..(.4,-.4);
 \fill[#1] #4 circle (1.7pt);
 \node at (-.4,-.6){$\scriptstyle #2$};\node at (.4,-.6){$\scriptstyle #3$};
\end{tikzpicture}}
\newcommand{\IRdotcup}[4]{\begin{tikzpicture}[baseline=0,scale=.7]
 \draw[#1] (-.4,.4)..controls(-.4,-.5)and(.4,-.5)..(.4,.4);
 \fill[#1] #4 circle (1.7pt);
 \node at (-.4,.6){$\scriptstyle #2$};\node at (.4,.6){$\scriptstyle #3$};
\end{tikzpicture}}
\newcommand{\IRcupcap}[3]{\begin{tikzpicture}[baseline=0,scale=.7]
 \draw[#1] (-.4,-.45)..controls(-.4,.05)and(.4,.05)..(.4,-.45);
 \draw[#1] (-.4,.55)..controls(-.4,.05)and(.4,.05)..(.4,.55);
 \node at (-.4,-.65){$\scriptstyle #2$};\node at (.4,-.65){$\scriptstyle #3$};
\end{tikzpicture}}
\newcommand{\IRcrosscap}[3]{\begin{tikzpicture}[baseline=0,scale=.7]
 \draw[#1] (-.4,-.5)--(.4,.12)
 ..controls(.4,.68)and(-.4,.68)..(-.4,.12)--(.4,-.5);
 \node at (-.4,-.7){$\scriptstyle #2$};\node at (.4,-.7){$\scriptstyle #3$};
\end{tikzpicture}}
\newcommand{\IRcupcross}[3]{\begin{tikzpicture}[baseline=0,scale=.7]
 \draw[#1] (-.4,.5)--(.4,-.12)
 ..controls(.4,-.68)and(-.4,-.68)..(-.4,-.12)--(.4,.5);
 \node at (-.4,.7){$\scriptstyle #2$};\node at (.4,.7){$\scriptstyle #3$};
\end{tikzpicture}}
\newcommand{\IRcupcrossside}[4]{\begin{tikzpicture}[baseline=0,scale=.65]
 \begin{scope}[xscale=#4]
  \draw[#1,thick] (-.58,.48)
   ..controls(-.58,-.50)and(.58,-.50)..(.58,.48);
  \draw[#1,thick] (0,.48) to[out=-90,in=90] (.82,-.64);
  \node at (.82,-.84){$\scriptstyle #3$};
 \end{scope}
 \node at (-.58,.70){$\scriptstyle #2$};
 \node at (0,.70){$\scriptstyle #3$};
 \node at (.58,.70){$\scriptstyle #2$};
\end{tikzpicture}}
\newcommand{\IRcapcrosscomp}[6]{\begin{tikzpicture}[baseline=0,scale=.65]
 \begin{scope}[xscale=#6]
  \draw[#1,thick] (-.58,-.48)
   ..controls(-.58,.50)and(.58,.50)..(.58,-.48);
  \draw[#1,thick] (0,-.48) to[out=90,in=-90] (-.82,.64);
  \node at (-.82,.84){$\scriptstyle #5$};
 \end{scope}
 \node at (-.58,-.70){$\scriptstyle #2$};
 \node at (0,-.70){$\scriptstyle #3$};
 \node at (.58,-.70){$\scriptstyle #4$};
\end{tikzpicture}}
\newcommand{\IRbubble}[2]{\begin{tikzpicture}[baseline=0,scale=.72]
 \draw[#1] (0,.38)..controls(-.5,.38)and(-.5,-.38)..(0,-.38)
 ..controls(.5,-.38)and(.5,.38)..(0,.38);
 \node at (0,-.62){$\scriptstyle #2$};
\end{tikzpicture}}

\newcommand{\IRcupcapdeco}[3]{\begin{tikzpicture}[baseline=0,scale=.72]
 \draw[#1,thick] (-.42,-.58)..controls(-.42,-.05)and(.42,-.05)..(.42,-.58);
 \draw[#1,thick] (-.42,.58)..controls(-.42,.05)and(.42,.05)..(.42,.58);
 \draw[#1,densely dotted] (-.42,-.42)--(.42,-.42);
 \draw[#1,densely dotted] (-.42,.42)--(.42,.42);
 \node[anchor=west] at (.52,-.42){$\scriptstyle #2$};
 \node[anchor=west] at (.52,.42){$\scriptstyle #3$};
 \node at (-.42,-.78){$\scriptstyle0$};\node at (.42,-.78){$\scriptstyle0$};
\end{tikzpicture}}

For $i>0$, give $E_i$ an upward arrow and $E_{-i}$ a downward
arrow; the downward dot is the negative of the unoriented dot.
We use the following cups and caps. In each picture the exterior
region has weight $\lambda$:
\begin{align*}
 \begin{tikzpicture}[baseline=0,scale=.75]
  \draw[darkred,thick,->] (-.38,.4)..controls(-.38,-.4)and(.38,-.4)..(.38,.4);
  \node at (-.38,.61){$\scriptstyle i$};\node at (.38,.61){$\scriptstyle i$};
 \end{tikzpicture}
 &=
 \begin{tikzpicture}[baseline=0,scale=.75]
  \draw[darkred,thick] (-.38,.4)..controls(-.38,-.4)and(.38,-.4)..(.38,.4);
  \node at (-.38,.61){$\scriptstyle-i$};\node at (.38,.61){$\scriptstyle i$};
 \end{tikzpicture},
 &
 \begin{tikzpicture}[baseline=0,scale=.75]
  \draw[darkred,thick,<-] (-.38,-.4)..controls(-.38,.4)and(.38,.4)..(.38,-.4);
  \node at (-.38,-.61){$\scriptstyle i$};\node at (.38,-.61){$\scriptstyle i$};
 \end{tikzpicture}
 &=\frac1{\zeta_i\gamma_i(\lambda)}
 \begin{tikzpicture}[baseline=0,scale=.75]
  \draw[darkred,thick] (-.38,-.4)..controls(-.38,.4)and(.38,.4)..(.38,-.4);
  \node at (-.38,-.61){$\scriptstyle-i$};\node at (.38,-.61){$\scriptstyle i$};
 \end{tikzpicture},\\
 \begin{tikzpicture}[baseline=0,scale=.75]
  \draw[darkred,thick,<-] (-.38,.4)..controls(-.38,-.4)and(.38,-.4)..(.38,.4);
  \node at (-.38,.61){$\scriptstyle i$};\node at (.38,.61){$\scriptstyle i$};
 \end{tikzpicture}
 &=\zeta_i\gamma_i(\lambda)
 \begin{tikzpicture}[baseline=0,scale=.75]
  \draw[darkred,thick] (-.38,.4)..controls(-.38,-.4)and(.38,-.4)..(.38,.4);
  \node at (-.38,.61){$\scriptstyle i$};\node at (.38,.61){$\scriptstyle-i$};
 \end{tikzpicture},
 &
 \begin{tikzpicture}[baseline=0,scale=.75]
  \draw[darkred,thick,->] (-.38,-.4)..controls(-.38,.4)and(.38,.4)..(.38,-.4);
  \node at (-.38,-.61){$\scriptstyle i$};\node at (.38,-.61){$\scriptstyle i$};
 \end{tikzpicture}
 &=
 \begin{tikzpicture}[baseline=0,scale=.75]
  \draw[darkred,thick] (-.38,-.4)..controls(-.38,.4)and(.38,.4)..(.38,-.4);
  \node at (-.38,-.61){$\scriptstyle i$};\node at (.38,-.61){$\scriptstyle-i$};
 \end{tikzpicture}.
\end{align*}
Since $i$ and $-i$ are nonadjacent, \eqref{ci1} gives
$\gamma_i(\lambda\pm\bar\alpha_i)=\gamma_i(\lambda)$.
Thus the scalar factors cancel in both zigzags. The downward-dot
convention turns \eqref{ipivotal} into the ordinary oriented
dot adjunction. Oriented crossings involving downward strands are
given by the mate diagrams in Definition~\ref{natrual-basic1}(2)--(3),
using these cups and caps and the indicated factor $t_{ji}^{-1}$.
Upward crossings are the original positive-color crossings. {The scaling factors satisfy, by \eqref{ci1},}
\[
 \frac{\gamma_i(\lambda-\bar\alpha_j)}{\gamma_i(\lambda)}
 =(-1)^{\delta_{i,j+1}+\delta_{i,1-j}}.
\]
In particular, the factor is $1$ for $j=0$.

Write $\circlearrowleft_{\lambda,i}(u)$ and
$\circlearrowright_{\lambda,i}(u)$ for the resulting counterclockwise
and clockwise bubble series, with the dot on the upward branch.
The fake-bubble coefficients use the same normalization.
On an object $V\in\mathcal M_\lambda$ we write
$\circlearrowleft_{V,i}=\circlearrowleft_{\lambda,i}|_V$ and
$\circlearrowright_{V,i}=\circlearrowright_{\lambda,i}|_V$.
In terms of the unoriented series of
\eqref{ibubblegeneratingfunction}, these conventions mean
\begin{align*}
 \circlearrowleft_{\lambda,i}(u)
 &=\frac1{\zeta_i\gamma_i(\lambda)}
 \begin{tikzpicture}[baseline=0,scale=.75]
  \draw[darkred,thick] (0,.36)..controls(-.48,.36)and(-.48,-.36)..(0,-.36)
   ..controls(.48,-.36)and(.48,.36)..(0,.36);
  \fill[darkred] (.36,0) circle(1.6pt);
  \node[anchor=west] at (.43,0){$\scriptstyle(u)$};
  \node at (-.59,0){$\scriptstyle-i$};
 \end{tikzpicture},
 &
 \circlearrowright_{\lambda,i}(u)
 &=-\zeta_i\gamma_i(\lambda)
 \begin{tikzpicture}[baseline=0,scale=.75]
  \draw[darkred,thick] (0,.36)..controls(-.48,.36)and(-.48,-.36)..(0,-.36)
   ..controls(.48,-.36)and(.48,.36)..(0,.36);
  \fill[darkred] (.36,0) circle(1.6pt);
  \node[anchor=west] at (.43,0){$\scriptstyle(-u)$};
  \node at (-.59,0){$\scriptstyle i$};
 \end{tikzpicture}.
\end{align*}
The minus sign in the second formula is the identity
$(u+y)^{-1}=-(-u-y)^{-1}$ after moving the dot to the right branch
by \eqref{ipivotal}. Relations \eqref{ci1} and \eqref{gammai} give
\[
 \gamma_i(\lambda)\gamma_{-i}(\lambda)=(-1)^{\lambda_i},
 \qquad \zeta_i\zeta_{-i}=-1.
\]
Consequently, \eqref{ibubblerel} and \eqref{iinfgrass} imply
\[
\begin{aligned}
 \circlearrowleft_{\lambda,i}(u)&\in u^{-\lambda_i}\one_\lambda
   +u^{-\lambda_i-1}\End(\one_\lambda)[[u^{-1}]],\\
 \circlearrowright_{\lambda,i}(u)&\in u^{\lambda_i}\one_\lambda
   +u^{\lambda_i-1}\End(\one_\lambda)[[u^{-1}]],\\
 \circlearrowleft_{\lambda,i}(u)\circlearrowright_{\lambda,i}(u)
 &=\one_\lambda.
\end{aligned}
\]

\begin{Defn}\label{gen-red1}
Fix $\lambda\in X^\imath$ and a finitely generated object
$V\in\mathcal M_\lambda$. We use the following notation.
\begin{itemize}
\item[(1)] $E=\bigoplus_{i\in I}E_i$.
\item [(2)] {Let $\varepsilon_i(V)\in\mathbb N$ be determined by requiring that $u^{\varepsilon_i(V)}$ is the minimal polynomial of}
\(\begin{tikzpicture}[baseline = -1mm,darkred, scale=0.7]
	\draw[-,thick] (0.08,-.4) to (0.08,.4);
      \node at (.08,0) {$\bullet$};
       \node at (.08,-.6) {$\scriptstyle{i}$};
      \node at (0.4,-.3) {$\color{black}\scriptstyle{\lambda}$};
\end{tikzpicture}
\begin{tikzpicture}[baseline = 0, scale=0.7]
  \draw[-,darkg,thick] (0.65,-0.25) to (0.65,.4);
   \node at (0.65,-.45) {$\darkg\scriptstyle{V}$};
\end{tikzpicture}\), where $i\in I$. We set $\varepsilon_i(V)=0$ if $E_iV=0$.
\item [(3)] $
 \begin{tikzpicture}[baseline = 0,scale=0.8]
  \draw[-,thick,darkred] (0.2,0.2) to[out=90,in=0] (0,.4);
  \draw[-,thick,darkred] (0,0.4) to[out=180,in=90] (-.2,0.2);
\draw[-,thick,darkred] (-.2,0.2) to[out=-90,in=180] (0,0);
  \draw[-,thick,darkred] (0,0) to[out=0,in=-90] (0.2,0.2);
   \node at (0.85,0.2) {$\scriptstyle{\lambda}$};
    \node at (-0.46,0.2) {$\scriptstyle{j}$};
   \node at (0.54,0.2) {$\color{black}\scriptstyle{(u)}$};
        \end{tikzpicture} =   \sum_{r\in \mathbb Z}  \begin{tikzpicture}[baseline = 0, scale=0.8, ]
  \draw[-,thick,darkred] (0.2,0.2) to[out=90,in=0] (0,.4);
  \draw[-,thick,darkred] (0,0.4) to[out=180,in=90] (-.2,0.2);
\draw[-,thick,darkred] (-.2,0.2) to[out=-90,in=180] (0,0);
  \draw[-,thick,darkred] (0,0) to[out=0,in=-90] (0.2,0.2);
 \node at (-0.46,0.2) {$\scriptstyle{j}$};
   \node at (0.2,0.2) {$\color{darkred}\scriptstyle\bullet$};
   \node at (0.4,0.2) {$\color{darkred}\scriptstyle{r}$};
    \node at (0.4,0.45) {$\scriptstyle{\lambda}$};
    \end{tikzpicture}
u^{-r-1}$ and $ \begin{tikzpicture}[baseline = 0,scale=0.8]
  \draw[-,thick,darkred] (0.2,0.2) to[out=90,in=0] (0,.4);
  \draw[-,thick,darkred] (0,0.4) to[out=180,in=90] (-.2,0.2);
\draw[-,thick,darkred] (-.2,0.2) to[out=-90,in=180] (0,0);
  \draw[-,thick,darkred] (0,0) to[out=0,in=-90] (0.2,0.2);
   \node at (0.85,0.2) {$\scriptstyle{\lambda}$};
    \node at (-0.42,0.2) {$\scriptstyle{i}$};
   \node at (0.54,0.2) {$\color{black}\scriptstyle{(u)}$};
        \end{tikzpicture} =   \sum_{r\in \mathbb Z}  \begin{tikzpicture}[baseline = 0, scale=0.8, ]
  \draw[-,thick,darkred] (0.2,0.2) to[out=90,in=0] (0,.4);
  \draw[-,thick,darkred] (0,0.4) to[out=180,in=90] (-.2,0.2);
\draw[-,thick,darkred] (-.2,0.2) to[out=-90,in=180] (0,0);
  \draw[-,thick,darkred] (0,0) to[out=0,in=-90] (0.2,0.2);
 \node at (-0.46,0.4) {$\scriptstyle{i}$};
   \node at (-0.2,0.2) {$\color{darkred}\scriptstyle\bullet$};
   \node at (-0.4,0.2) {$\color{darkred}\scriptstyle{r}$};
    \node at (0.4,0.4) {$\scriptstyle{\lambda}$};
    \end{tikzpicture}
u^{-r-1}$, where   $i\in I_1$ and $j\in -I_1\bigcup\{0\}$.
\end{itemize}
At color zero, (3) is the ordinary dotted-bubble series. The zero-color
series in \eqref{ibubblegeneratingfunction} is obtained by subtracting
$\one_\lambda/(2u)$; it is used below only when explicitly cited.
\end{Defn}

\subsection{The zero-color bubble series}

{Define $c(\lambda)\in\{0,1\}$ by}
\[
 c(\lambda)=\begin{cases}
 0,&\lambda=\sum_{i>0}\lambda_i\bar\Lambda_i,\\
 1,&\lambda=\sum_{i>0}\lambda_i\bar\Lambda_i+\bar\Lambda_0.
 \end{cases}
\]
{Then $\gamma_0(\lambda)=(-1)^{c(\lambda)}$.} This is the parameter $t$ in the corresponding nil-Brauer category (see \cite[\S~3.3, first version]{BWW25}).  In
particular, the undotted red $0$-bubble is
\[
 \begin{tikzpicture}[baseline = 0, scale=0.8, ]
  \draw[-,thick,darkred] (0.2,0.2) to[out=90,in=0] (0,.4);
  \draw[-,thick,darkred] (0,0.4) to[out=180,in=90] (-.2,0.2);
\draw[-,thick,darkred] (-.2,0.2) to[out=-90,in=180] (0,0);
  \draw[-,thick,darkred] (0,0) to[out=0,in=-90] (0.2,0.2);
 \node at (0,-.19) {$\scriptstyle{0}$};
   \node at (-0.36,0.2) {$\scriptstyle{\lambda}$};
    \end{tikzpicture} =c(\lambda)\one_\lambda.
\]

Define
\begin{equation}\label{eq:O-lambda-zero}
 O_{\lambda,0}(u)
 :=(-1)^{c(\lambda)}
 \bigl(\one_\lambda-2u \begin{tikzpicture}[baseline = 0,scale=0.8]
  \draw[-,thick,darkred] (0.2,0.2) to[out=90,in=0] (0,.4);
  \draw[-,thick,darkred] (0,0.4) to[out=180,in=90] (-.2,0.2);
\draw[-,thick,darkred] (-.2,0.2) to[out=-90,in=180] (0,0);
  \draw[-,thick,darkred] (0,0) to[out=0,in=-90] (0.2,0.2);
 \node at (0,-.19) {$\scriptstyle{0}$};
   \node at (-0.36,0.2) {$\scriptstyle{\lambda}$};
   \node at (0.54,0.2) {$\color{black}\scriptstyle{(u)}$};
        \end{tikzpicture}\bigr).
\end{equation}
For $V\in\M_\lambda$, put
\[
 O_{V,0}(u):=O_{\lambda,0}(u)|_V.
\]
The nil-Brauer relations \cite[(2.22)--(2.23), Theorem~2.5]{BWW24} give
\begin{equation}\label{eq:Ozero-coefficients}
 \Coeff{O_{\lambda,0}(u)}{u:0}=\one_\lambda,
 \quad \Coeff{O_{\lambda,0}(u)}{u:r}=0,\quad
 \Coeff{O_{\lambda,0}(u)}{u:-r}
 =-2(-1)^{c(\lambda)}\begin{tikzpicture}[baseline = 0, scale=0.8, ]
  \draw[-,thick,darkred] (0.2,0.2) to[out=90,in=0] (0,.4);
  \draw[-,thick,darkred] (0,0.4) to[out=180,in=90] (-.2,0.2);
\draw[-,thick,darkred] (-.2,0.2) to[out=-90,in=180] (0,0);
  \draw[-,thick,darkred] (0,0) to[out=0,in=-90] (0.2,0.2);
 \node at (0,-.19) {$\scriptstyle{0}$};
   \node at (-0.36,0.2) {$\scriptstyle{\lambda}$};
   \node at (0.2,0.2) {$\color{darkred}\scriptstyle\bullet$};
   \node at (0.4,0.2) {$\color{darkred}\scriptstyle{r}$};
    \end{tikzpicture}
 \quad(r\geq1),
\end{equation}
and
\begin{equation}\label{eq:Ozero-inverse}
 O_{\lambda,0}(u)O_{\lambda,0}(-u)=\one_\lambda.
\end{equation}

Define
\begin{equation}\label{eq:tildeO-lambda}
 \widetilde O_\lambda(u):=O_{\lambda,0}(u)
 \prod_{i>0}(-1)^{\lambda_i}
 \circlearrowleft_{\lambda,i}(u-i)
 \circlearrowright_{\lambda,i}(-u-i).
\end{equation}
The product is interpreted after evaluation on finitely generated
objects, on which only finitely many factors differ from $1$.

\begin{Lemma}\label{lem:tildeO-inverse}
One has
\[
 \widetilde O_\lambda(u)\widetilde O_\lambda(-u)=\one_\lambda.
\]
\end{Lemma}

\begin{proof}
{In $\widetilde O_\lambda(u)\widetilde O_\lambda(-u)$, the two
oriented series cancel at each nonzero color. The zero-color factors
cancel by \eqref{eq:Ozero-inverse}.}
\end{proof}

From now on, fix $\lambda\in X^\imath $ and a finitely generated object
$V\in\M_\lambda$.  We may assume $V\ne0$, since identities on the zero object are
automatic. Let $Z_V\subseteq Z(\End V)$ be the commutative algebra
generated by evaluated red bubbles. We show below that $Z_V$ is local.

\newcommand{\IRzerobubbledeco}[1]{\begin{tikzpicture}[baseline=0,scale=.75]
 \draw[darkred,thick] (0,.42) ..controls(-.56,.42)and(-.56,-.42)..(0,-.42)
  ..controls(.56,-.42)and(.56,.42)..(0,.42);
 \fill[darkred] (.42,0) circle(1.6pt);
 \node[anchor=west] at (.50,0){$\scriptstyle #1$};
 \node at (0,-.62){$\scriptstyle0$};
\end{tikzpicture}}
\newcommand{\IRzerostranddeco}[1]{\begin{tikzpicture}[baseline=0,scale=.75]
 \draw[darkred,thick] (0,-.48)--(0,.48);
 \fill[darkred] (0,0) circle(1.6pt);
 \node[anchor=west] at (.12,0){$\scriptstyle #1$};
 \node at (0,-.68){$\scriptstyle0$};
\end{tikzpicture}}
\newcommand{\IRzerocapdeco}[1]{\begin{tikzpicture}[baseline=0,scale=.75]
 \draw[darkred,thick] (-.42,-.48) --(-.42,0)
  ..controls(-.42,.60)and(.42,.60)..(.42,0)--(.42,-.48);
 \fill[darkred] (.42,-.16) circle(1.6pt);
 \node[anchor=west] at (.53,-.16){$\scriptstyle #1$};
 \node at (-.42,-.68){$\scriptstyle0$};\node at (.42,-.68){$\scriptstyle0$};
\end{tikzpicture}}
\newcommand{\IRzerocapcrossdeco}[1]{\begin{tikzpicture}[baseline=0,scale=.75]
 \draw[darkred,thick] (-.42,-.72)--(.42,-.20)--(.42,.12)
  ..controls(.42,.72)and(-.42,.72)..(-.42,.12)--(-.42,-.20)--(.42,-.72);
 \fill[darkred] (.42,0) circle(1.6pt);
 \node[anchor=west] at (.53,0){$\scriptstyle #1$};
 \node at (-.42,-.92){$\scriptstyle0$};\node at (.42,-.92){$\scriptstyle0$};
\end{tikzpicture}}

\begin{Lemma}\label{lem:Bzero-polynomial}
Suppose $f(u)\in Z_V[u]$ and
$\IRzerostranddeco{f(y)}|_V=0$. Then
\[
 g(u):=\IRzerobubbledeco{(u-y)^{-1}}\,f(u)\Big|_V
 \in Z_V[u].
\]
If $f(u)\in\mathbb C[u]$, one also has
\begin{equation}\label{eq:zero-difference-quotient}
 \IRzerostranddeco{g(-y)}\Big|_V
 =\IRzerostranddeco{\frac{f(y)-f(-y)}{2y}}\Big|_V.
\end{equation}
The quotient is formed in $\mathbb C[y]$ before evaluation, and
coefficients from $Z_V$ act through $E_0$.
\end{Lemma}
\begin{proof}
All Laurent expansions are at infinity. For $r\ge0$,
\[
 [g(u)]_{u:-r-1}
 =\IRzerobubbledeco{y^rf(y)}\Big|_V=0,
\]
which proves polynomiality. The dot-crossing relation and the vanishing
red curl give the following generating-function identity, equivalently
\cite[(2.24)]{BWW24}:
\[
 \IRzerocapcrossdeco{(u-y)^{-1}}
 =\IRzerocapdeco{(u+y)^{-1}}\,
   \IRzerobubbledeco{(u-y)^{-1}}
  -\IRzerocapdeco{((u-y)(u+y))^{-1}}.
\]
Multiply by $f(u)$ and extract $u^{-1}$. The left-hand side vanishes,
because $f(y)=0$ on $E_0V$, and the marked right leg above the
crossing is obtained by applying $E_0$ to this endomorphism.
The right-hand side becomes
\[
 0=\IRzerocapdeco{g(-y)}\Big|_V
   -\IRzerocapdeco{\frac{f(y)-f(-y)}{2y}}\Big|_V,
\]
since
$[f(u)/(u^2-y^2)]_{u:-1}=(f(y)-f(-y))/(2y)$.
Taking the adjoint mate gives \eqref{eq:zero-difference-quotient}.
\end{proof}

\begin{Lemma}\label{lem:Ozero-annihilator}
Let $f(u)\in\mathbb C[u]$ be monic with
$\IRzerostranddeco{f(y)}|_V=0$, and set
$g(u):=O_{V,0}(u)f(u)$. Then $g\in Z_V[u]$ is monic of degree
$\deg f$ and
\[
 \IRzerostranddeco{g(-y)}\Big|_V=0.
\]
\end{Lemma}
\begin{proof}
Polynomiality follows from Lemma~\ref{lem:Bzero-polynomial} and
\eqref{eq:O-lambda-zero}; the leading term of $O_{V,0}$ is $1$.
By \eqref{eq:zero-difference-quotient},
\[
 \IRzerostranddeco{g(-y)}\Big|_V
 =(-1)^{c(\lambda)}
 \IRzerostranddeco{f(-y)+2y\frac{f(y)-f(-y)}{2y}}\Big|_V
 =(-1)^{c(\lambda)}\IRzerostranddeco{f(y)}\Big|_V=0.
\]
Only polynomials are evaluated at the nilpotent dot.
\end{proof}

\begin{Cor}\label{cor:Ozero-simple}
If $L\in\M_\lambda$ is simple, then $O_{L,0}(u)=1$.
\end{Cor}
\begin{proof}
By Lemma~\ref{lem:Ozero-annihilator},
$O_{L,0}(u)u^{\eps_0(L)}\in\mathbb C[u]$ is monic of degree
$\eps_0(L)$ and annihilates $-y$ on $E_0L$. Its minimal polynomial is
$u^{\eps_0(L)}$, so the two monic polynomials coincide.
\end{proof}

\begin{Lemma}\label{lem:zero-bubble-parity}
Let $L\in\M_\lambda$ be simple. Then
$c(\lambda)\equiv\eps_0(L)\pmod2$.
\end{Lemma}
\begin{proof}
Corollary~\ref{cor:Ozero-simple} and \eqref{eq:O-lambda-zero} give
\[
 \IRzerobubbledeco{(u-y)^{-1}}\Big|_L=\frac{c(\lambda)}u.
\]
If $E_0L=0$, the bubble vanishes and $c(\lambda)=0$.
Otherwise apply \eqref{eq:zero-difference-quotient} to
$f(u)=u^{\eps_0(L)}$:
\[
 c(\lambda)(-1)^{\eps_0(L)-1}
 \IRzerostranddeco{y^{\eps_0(L)-1}}\Big|_L
 =\frac{1-(-1)^{\eps_0(L)}}2
 \IRzerostranddeco{y^{\eps_0(L)-1}}\Big|_L.
\]
The displayed strand is nonzero by minimality, proving the assertion.
\end{proof}

\subsection{The total red series and the black bubble series}

\begin{Defn} \label{evalu-bub} Suppose $i\in I^\imath$. We use the
normalized oriented bubbles of the type $A$ convention: the two
series below are inverse, with leading terms $u^{-\lambda_i}$ and
$u^{\lambda_i}$, respectively. For any $\lambda\in X^\imath$ and any
finitely generated object \(V\in\mathcal M_\lambda\), define
\begin{itemize} \item [(1)] $\circlearrowleft_{V,i}(u)=\begin{tikzpicture}[baseline = 0,scale=0.8]
  \draw[->,thick,darkred] (0.2,0.2) to[out=90,in=0] (0,.4);
  \draw[-,thick,darkred] (0,0.4) to[out=180,in=90] (-.2,0.2);
\draw[-,thick,darkred] (-.2,0.2) to[out=-90,in=180] (0,0);
  \draw[-,thick,darkred] (0,0) to[out=0,in=-90] (0.2,0.2);
 \node at (0,-.19) {$\scriptstyle{i}$};
   \node at (-0.36,0.2) {$\scriptstyle{\lambda}$};
   \node at (0.2,0.2) {$\color{darkred}\scriptstyle\bullet$};
   \node at (0.54,0.2) {$\color{black}\scriptstyle{(u)}$};
        \end{tikzpicture}  \Big|_V$,
        \item [(2)]  $\circlearrowright_{V,i}(u)=\begin{tikzpicture}[baseline = 0,scale=0.8]
  \draw[<-,thick,darkred] (0,0.4) to[out=180,in=90] (-.2,0.2);
  \draw[-,thick,darkred] (0.2,0.2) to[out=90,in=0] (0,.4);
 \draw[-,thick,darkred] (-.2,0.2) to[out=-90,in=180] (0,0);
  \draw[-,thick,darkred] (0,0) to[out=0,in=-90] (0.2,0.2);
 \node at (0,-.19) {$\scriptstyle{i}$};
   \node at (0.3,0.45) {$\scriptstyle{\lambda}$};
   \node at (-0.2,0.2) {$\color{darkred}\scriptstyle\bullet$};
   \node at (0.45,0.15) {$\color{black}\scriptstyle{(u)}$};
\end{tikzpicture} \Big|_V$,
\item [(3)]$m_V(u)=\prod_{j\in I}(u-j)^{\varepsilon_j(V)}$, and  $a_V= \text{deg } m_V(u)$,
\item [(4)] $
\widetilde O_V(u):=\widetilde O_\lambda(u)\Big|_V$,
\end{itemize}
\end{Defn}

\begin{Lemma}\label{lem:red-series-simple}
Let $i\in I^\imath$ and let $L\in\M_\lambda$ be simple. Then
\begin{align}
\circlearrowleft_{L,i}(u)&=u^{\eps_{-i}(L)-\eps_i(L)},\label{eq:Bplus-simple}\\
 \eps_{-i}(L)-\eps_i(L)
 &=-\lambda_i,\label{eq:weight-difference}\\
 \circlearrowright_{L,i}(u)&=u^{\eps_i(L)-\eps_{-i}(L)}.\label{eq:Bminus-simple}
\end{align}
\end{Lemma}

\begin{proof}
The nonexceptional calculation of \cite[Lemmas~7.4 and~7.6]{GRS}
applies, since $i$ and $-i$ are nonadjacent. The two monic polynomials
\[
 u^{\eps_i(L)}\circlearrowleft_{L,i}(u),\qquad
 u^{\eps_{-i}(L)}\circlearrowright_{L,i}(u)
\]
{annihilate the dots on $E_{-i}L$ and $E_iL$, respectively. Their
degrees give
$\eps_{-i}(L)\le\eps_i(L)-\lambda_i$ and
$\eps_i(L)\le\eps_{-i}(L)+\lambda_i$.
Both inequalities are equalities, so each polynomial is the
corresponding minimal polynomial.}
\end{proof}

For every $\lambda\in X^\imath$, including those with
$\mathcal M_\lambda=0$, define $a_\lambda\in\{0,1\}$ by
\[
 a_\lambda\equiv c(\lambda)+\sum_{i>0}\lambda_i\pmod2.
\]
\begin{rem}\label{rem:a-lambda-well-defined}
If $L\in\M_\lambda$ is simple, then
\[
 a_L=\eps_0(L)+\sum_{i>0}(\eps_i(L)+\eps_{-i}(L))
 \equiv c(\lambda)+\sum_{i>0}\lambda_i
 \equiv a_\lambda\pmod2
\]
by Lemmas~\ref{lem:zero-bubble-parity} and
\ref{lem:red-series-simple}.
\end{rem}

By Corollary~\ref{cor:Ozero-simple} and
Lemma~\ref{lem:red-series-simple}, every bubble has the same scalar
value on the simple objects of weight $\lambda$. The proof of
\cite[Lemma~7.8]{GRS} therefore shows that $Z_V$ is a finite-dimensional
local algebra with residue field $\mathbb C$.

Define
\begin{equation}\label{eq:nV-definition}
 n_V(u):=\left(u-\frac12(-1)^{a_\lambda}\right)
 \widetilde O_V(u)m_V(u).
\end{equation}
In the evaluated notation of Definition~\ref{evalu-bub},
\[
 \widetilde O_V(u)=O_{V,0}(u)
 \prod_{i>0}(-1)^{\lambda_i}
 \circlearrowleft_{V,i}(u-i)\circlearrowright_{V,i}(-u-i).
\]
Hence
\begin{align*}
 n_V(u)={}&\left(u-\frac12(-1)^{a_\lambda}\right)
 \bigl(u^{\eps_0(V)}O_{V,0}(u)\bigr)\\
 &\quad\cdot\prod_{i>0}
 \bigl((u-i)^{\eps_i(V)}\circlearrowleft_{V,i}(u-i)\bigr)
 \bigl((-1)^{\lambda_i}(u+i)^{\eps_{-i}(V)}
       \circlearrowright_{V,i}(-u-i)\bigr).
\end{align*}
For arbitrary finitely generated $V$, the nonexceptional calculation
of \cite[Lemma~7.4]{GRS}, with the normalization fixed above, shows that
$u^{\eps_i(V)}\circlearrowleft_{V,i}(u)$ and
$u^{\eps_{-i}(V)}\circlearrowright_{V,i}(u)$ are monic polynomials in $Z_V[u]$ of degrees
$\eps_i(V)-\lambda_i$ and $\eps_{-i}(V)+\lambda_i$, respectively.
Together with Lemma~\ref{lem:Ozero-annihilator}, this shows that each
parenthesized factor is monic and that their degrees sum to $\deg m_V+1$.
Thus $n_V$ is monic of that degree. All but finitely many factors are $1$.

For $z\in Z_V$, write $\bar z$ for its image in $\mathbb C$.
As in \cite[Remark~7.9]{GRS}, let $t_{V,i}$ be the unique monic factor
of $n_V$ whose reduction is a power of $u-i$ and whose complement has
coprime reduction. Put
\[
 m_{V,i}(u):=\frac{m_V(u)}{(u-i)^{\eps_i(V)}},
 \qquad n_{V,i}(u):=\frac{n_V(u)}{t_{V,i}(u)}.
\]

\begin{Lemma}\label{lem:tVi-integral}
For every $i\in I$,
\[
 t_{V,i}(u)=
 \begin{cases}
(u-i)^{\varepsilon_i(V)}\circlearrowleft_{V,i}(u-i),&i\in I^\imath,\\[2mm]
 (-1)^{\lambda_{-i}}(u-i)^{\varepsilon_{i}(V)}\circlearrowright_{V,-i}(-u+i),&i\in-I^\imath,\\[2mm]
 O_{V,0}(u)u^{\eps_0(V)},&i=0.
 \end{cases}
\]
Moreover, for every simple object $L\in\M_\lambda$,
\begin{equation}\label{eq:tildeO-simple}
 \widetilde O_L(u)=(-1)^{a_L}\frac{m_L(-u)}{m_L(u)}.
\end{equation}
\end{Lemma}

\begin{proof}
The preceding factorization identifies $t_{V,i}$ by uniqueness,
as in \cite[Lemma~7.10]{GRS}. Finally,
Corollary~\ref{cor:Ozero-simple} and
Lemma~\ref{lem:red-series-simple} give
\[
 \widetilde O_L(u)=\prod_{i>0}
 \left(\frac{u+i}{u-i}\right)^{\lambda_i}
 =(-1)^{a_L}\frac{m_L(-u)}{m_L(u)},
\]
as required.
\end{proof}

We will also use the complementary factors
\[
\begin{aligned}
 m^+_{V,0}(u)&:=\prod_{i>0}(u-i)^{\eps_i(V)},&
 m^-_{V,0}(u)&:=\prod_{i<0}(u-i)^{\eps_i(V)},\\
 n^+_{V,0}(u)&:=\left(u-\tfrac12\right)^{1-a_\lambda}
                 \prod_{i>0}t_{V,i}(u),&
 n^-_{V,0}(u)&:=\left(u+\tfrac12\right)^{a_\lambda}
                 \prod_{i<0}t_{V,i}(u).
\end{aligned}
\]
They satisfy $m_{V,0}=m^+_{V,0}m^-_{V,0}$ and
$n_{V,0}=n^+_{V,0}n^-_{V,0}$, and their constant terms are units.
These polynomials are constructed entirely from the given red action.

\begin{Lemma}\label{lem:zero-ratio-residue}
For every finitely generated $V\in\M_\lambda$,
\begin{equation}\label{eq:zero-ratio-residue}
 {\frac{n_{V,0}(0)}{m_{V,0}(0)}}
 =(-1)^{c(\lambda)+1}\frac12.
\end{equation}
\end{Lemma}

\begin{proof}
By Lemma~\ref{lem:tVi-integral} and the inverse relation for the
oriented series,
\begin{align*}
 \frac{n^+_{V,0}(u)}{m^+_{V,0}(u)}
 \frac{n^-_{V,0}(-u)}{m^-_{V,0}(-u)}
 &=(u-\tfrac12)^{1-a_\lambda}(-u+\tfrac12)^{a_\lambda}
   \prod_{i>0}(-1)^{\lambda_i}\\
 &=(-1)^{c(\lambda)}(u-\tfrac12).
\end{align*}
For later use, \eqref{eq:nV-definition} and
Lemma~\ref{lem:tildeO-inverse} also give
\begin{equation}\label{eq:nV-product}
 n_V(u)n_V(-u)
 =\left(\frac12-u\right)\left(\frac12+u\right)m_V(u)m_V(-u).
\end{equation}
Evaluation of the first identity at zero gives the asserted exact
equality, equivalently
\begin{equation}\label{eq:zero-ratio-exact}
 \frac{n_{V,0}(0)}{m_{V,0}(0)}
 =\begin{cases}-\dfrac12,&c(\lambda)=0,\\[1mm]
                \dfrac12,&c(\lambda)=1.
   \end{cases}\qedhere
\end{equation}
\end{proof}

\subsection{Bubble slides}
On $E_i$ and $E_{-i}$, for $i>0$, let $y$ be the oriented red dot
and put $x=y+i$ and $x=-y-i$, respectively. On $E_0$ put $x=y$.
A black strand in the following slides denotes this shifted dot;
no black crossing has yet been defined.

\begin{Lemma}\label{lem:integral-parity-slide}
For every $i\in I$, $a_{\lambda-\bar\alpha_i}=a_\lambda$, and
\begin{equation}\label{eq:integral-O-slide-diagram}
  \mathord{
\begin{tikzpicture}[baseline = -1mm]
	\draw[-] (0.08,-.3) to (0.08,.3);
    \node at (.08, -.4) {$\scriptstyle{i}$};
      \node at (-0.74,0) {$\widetilde O_{\lambda-\bar\alpha_i}\scriptstyle (u)$};
\end{tikzpicture}
  }=
\mathord{
\begin{tikzpicture}[baseline = -1mm]
	\draw[-] (0.08,-.3) to (0.08,.3);
    \node at (.08, -.4) {$\scriptstyle{i}$};
      \node at (0.6,0) {$\widetilde O_\lambda (u) $};
      \node at (.08,0) {$\scriptstyle\bullet$};
      \node at (-1.2,0) {$\scriptstyle{\frac{(u+x)^2-1}{(u-x)^2-1}\frac{(u-x)^2}{(u+x)^2}}$};
\end{tikzpicture}}.
\end{equation}
Here the black dot is $x$ on the $i$-colored summand.
\end{Lemma}
\begin{proof}
For $j=0$, $\bar\alpha_0=0$. For $j\ne0$,
\[
 \sum_{i>0}\langle\alpha_i^{\theta,\vee},\bar\alpha_j\rangle
 \equiv\delta_{j,1}+\delta_{j,-1}\pmod2,
 \qquad
 c(\lambda-\bar\alpha_j)-c(\lambda)
 \equiv\delta_{j,1}+\delta_{j,-1}\pmod2.
\]
The definition of $a_\lambda$ therefore gives
$a_{\lambda-\bar\alpha_j}=a_\lambda$ for every $j$.

In each nonzero bubble pair in \eqref{eq:tildeO-lambda}, the
normalization factors $\zeta_i\gamma_i(\lambda)$ and their inverses
cancel. Writing both raw bubbles with dots on their right branches
gives arguments $u-i,u+i$ and the sign $(-1)^{\lambda_i+1}$.
The normalized zero series is $-2u(-1)^{c(\lambda)}$ times the
unoriented series in \eqref{ibubblegeneratingfunction}.
The parity just proved shows that the ratio of all these scalar
factors across a $j$-strand is $1$.

Let $y$ be its unoriented dot, so that the black dot is $x=j+y$.
Applying \eqref{ibubslide} to the raw factors and reindexing the
numerator gives the transport factor
\begin{align*}
 \prod_{k\in\mathbb Z}
 \frac{R_{k,j}(-u-k,y)}{R_{k,j}(u-k,y)}
 &=\frac{(u-x)^2}{(u+x)^2}
   \frac{(u+x-1)(u+x+1)}{(u-x-1)(u-x+1)}\\
 &=\frac{(u+x)^2-1}{(u-x)^2-1}
   \frac{(u-x)^2}{(u+x)^2}.
\end{align*}
By \eqref{Qij}, only $k=j,j\pm1$ contribute: the equal-color factor
gives the square ratio, and the two adjacent factors give the other
ratio. This calculation also covers $j=0,\pm1$ and proves
\eqref{eq:integral-O-slide-diagram}.
\end{proof}

Set $\widehat O_\lambda(u)=
(u-\tfrac12(-1)^{a_\lambda})\widetilde O_\lambda(u)$.
Lemmas~\ref{lem:tildeO-inverse} and \ref{lem:integral-parity-slide} give
\begin{align}
 \widehat O_\lambda(u)\widehat O_\lambda(-u)&=\tfrac14-u^2,
 \label{eq:integral-widehat-inverse}\\
  \mathord{
\begin{tikzpicture}[baseline = -1mm]
	\draw[-] (0.08,-.3) to (0.08,.3);
    \node at (.08, -.4) {$\scriptstyle{i}$};
      \node at (-0.74,0) {$\scriptstyle\widehat O_{\lambda-\bar\alpha_i} (u)$};
\end{tikzpicture}
  }=
\mathord{
\begin{tikzpicture}[baseline = -1mm]
	\draw[-] (0.08,-.3) to (0.08,.3);
    \node at (.08, -.4) {$\scriptstyle{i}$};
      \node at (0.6,0) {$\widehat O_\lambda (u) $};
      \node at (.08,0) {$\scriptstyle\bullet$};
      \node at (-1.2,0) {$\scriptstyle{\frac{(u+x)^2-1}{(u-x)^2-1}\frac{(u-x)^2}{(u+x)^2}}$};
\end{tikzpicture}}
 \label{eq:integral-widehat-slide}\\
 \widehat O_\lambda(u)|_V&=\frac{n_V(u)}{m_V(u)}.
 \label{eq:integral-widehat-black}
\end{align}
These identities follow from the red relations.

Assume $E_jV\ne0$, and let $x$ be its shifted dot. Work in the
commutative subalgebra of $\End(E_jV)$ generated by $x-j$,
$E_j(Z_V)$ and $Z_{E_jV}$. It is finite-dimensional and local:
$x-j$ and the images of the two maximal ideals generate a nilpotent
ideal with quotient $\mathbb C$.
After clearing denominators in \eqref{eq:integral-widehat-slide},
uniqueness of monic primary factors in this algebra gives, for $i\in\mathbb Z$,
\begin{equation}\label{eq:ir-primary-transport}
\begin{split}
 &\frac{t_{E_jV,i}(u)}{(u-i)^{\varepsilon_i(E_jV)}}\,
 E_j\!\left(\frac{(u-i)^{\varepsilon_i(V)}}{t_{V,i}(u)}\right)\\
 &\qquad=
 \frac{(u+x-1)^{\delta_{i,1-j}}(u+x+1)^{\delta_{i,-1-j}}
       (u-x)^{2\delta_{i,j}}}
 {(u-x-1)^{\delta_{i,j+1}}(u-x+1)^{\delta_{i,j-1}}
       (u+x)^{2\delta_{i,-j}}}.
\end{split}
\end{equation}
Group factors with the same reduction before applying uniqueness;
this includes $j=0,\pm1$. The factor at $(-1)^{a_\lambda}/2$ is unchanged.

Use the complementary factors constructed above to define
$\Gamma_\lambda,\gamma_{V;s,i}$ and $\Gamma_V$ by the formulas of
Definition~\ref{GammaLambda}, with $d_s(\lambda)=\lambda_s$:
\begin{align*}
 \gamma_{V;s,i}(u)&=
 \frac{(-1)^{d_s(\lambda)}
 t_{V,-s}(-u-s+i)t_{V,-s}(u-s+i)}
 {(i^2-u^2)^{\varepsilon_{-s}(V)}},\qquad 1\le i\le s,\\
 \Gamma_V(u)&=
 \frac{1}{(u^2-1)^{d_1(\lambda)}}
 \prod_{s\ge2}
 \left(\frac{\gamma_{V;s,1}(u)}
 {(u^2-1)^{d_s(\lambda)}}\right)^s
 \prod_{\substack{s\ge2\\1\le i\le s}}
 \gamma_{V;s,i}(u)^{-1}.
\end{align*}
Only finitely many factors differ from $1$. The reductions of the
numerator factors of $\gamma_{V;s,i}(0)$ are powers of $i\ne0$.
Thus $\Gamma_V$ is an even unit in $Z_V[[u]]$.
Set
\begin{equation}\label{eq:ir-WM-definitions}
 W_V(u)=\frac{n^+_{V,0}(u)}{m^+_{V,0}(u)\Gamma_V(u)},\qquad
 M_V(u)=2(-1)^{c(\lambda)+1}
 \frac{n^-_{V,0}(u)\Gamma_V(u)}{m^-_{V,0}(u)}.
\end{equation}
{Apply the B\'ezout construction of Definition~\ref{notation-1} to
these red data, retaining its vanishing factors at nonzero colors.
At nilpotent dots, the resulting finite representatives satisfy
$A_V=W_V^{-1}$ and $B_V=M_V^{-1}$ once sufficiently many terms
have been retained.}

\begin{Lemma}\label{lem:ir-normalization}
The functions in \eqref{eq:ir-WM-definitions} satisfy
\begin{equation}\label{eq:ir-WM-products}
 W_V(u)M_V(-u)=1-2u,\qquad
 M_V(u)W_V(-u)=1+2u.
\end{equation}
Their nilpotent evaluations, including all divided differences used
below, are natural and independent of the polynomial representatives
once sufficiently many terms have been retained.
\end{Lemma}
\begin{proof}
The inverse relation for the oriented bubbles gives
\[
 \frac{n^+_{V,0}(u)n^-_{V,0}(-u)}
 {m^+_{V,0}(u)m^-_{V,0}(-u)}
 =(-1)^{c(\lambda)}(u-\tfrac12).
\]
Since $\Gamma_V$ is even, this proves \eqref{eq:ir-WM-products}.
For naturality, use the common formulas on $V\oplus V'$ for a morphism
$V\to V'$. On restriction to $V$, the extra factors in $m_{V\oplus V'}$
also occur in $n_{V\oplus V'}$ and cancel in every ratio defining
$\Gamma$, $W$ and $M$; likewise for $V'$. Finally, the dots are
nilpotent, and a divided difference lowers degree by at most its order.
Thus sufficiently many terms of each series determine every
calculation, independently of the chosen polynomial representative.
\end{proof}

\begin{Lemma}\label{lem:ir-Gamma-transport}
With $x$ the dot on the indicated summand,
\begin{align}
 \Gamma_{E_0V}(u)&=E_0(\Gamma_V(u)),
 \label{eq:ir-Gamma-zero}\\
 \Gamma_{E_jV}(u)&=E_j(\Gamma_V(u))
 \frac{u^2-x^2}{u^2-(x+1)^2}, &&j>0,
 \label{eq:ir-Gamma-positive}\\
 \Gamma_{E_{-j}V}(u)&=E_{-j}(\Gamma_V(u))
 \frac{u^2-(x-1)^2}{u^2-x^2}, &&j>0.
 \label{eq:ir-Gamma-negative}
\end{align}
Moreover, on $E_0V$,
\begin{align}
 W_{E_0V}(u)&=E_0(W_V(u))\frac{u+x-1}{u-x-1},
 \label{eq:ir-W-zero-transport}\\
 M_{E_0V}(u)&=E_0(M_V(u))\frac{u+x+1}{u-x+1}.
 \label{eq:ir-M-zero-transport}
\end{align}
\end{Lemma}
\begin{proof}
Substitute \eqref{eq:ir-primary-transport} into the product defining
$\Gamma_V$, as in the proof of Lemma~\ref{fw:Gamma-positive}.
For $E_0$, the factors indexed by $s\ge2$ are unchanged. For $E_j$
with $j>0$, only the factors at $-j,-j\pm1$ change; cancellation
leaves $(u^2-x^2)/(u^2-(x+1)^2)$. For $E_{-j}$ the same computation
gives its displayed counterpart. Factors indexed by zero are absent,
so this also covers $j=1$.
The last two identities follow by multiplying, respectively, the
positive and negative factors in \eqref{eq:ir-primary-transport}
and using \eqref{eq:ir-Gamma-zero}.
\end{proof}

\begin{Lemma}\label{lem:ir-W-positive-transport}
For $j>0$, with $x$ the dot on $E_jV$, one has
\begin{equation}\label{eq:ir-W-positive-transport}
 W_{E_jV}(u)=E_j(W_V(u))
 \frac{(u-x)(u+x+1)}{(u+x)(u-x+1)^{1-\delta_{j,1}}}.
\end{equation}
\end{Lemma}
\begin{proof}
Multiplying the factors at positive colors in
\eqref{eq:ir-primary-transport} gives
\[
 \frac{n^+_{E_jV,0}(u)}{m^+_{E_jV,0}(u)}
 =E_j\!\left(\frac{n^+_{V,0}(u)}{m^+_{V,0}(u)}\right)
 \frac{(u-x)^2}{(u-x-1)(u-x+1)^{1-\delta_{j,1}}}.
\]
Divide by \eqref{eq:ir-Gamma-positive} and cancel the common factors
before nilpotent evaluation. This proves the formula, including $j=1$.
\end{proof}

\subsection{Reverse dots, cups, and caps}

\begin{Defn}
\label{natrual-basic-red}
For $\lambda\in X^\imath$, $i\in I^\imath$, $j,k\in I$,
and finitely generated $V\in\mathcal M_\lambda$, define the following
maps. At color zero the black strand is $E_0$ and its dot is the red dot.
The decorations in (10)--(11) are on the marked right legs.
\begin{itemize}
\item[(1)]
    $\begin{tikzpicture}[baseline=-0.5ex,scale=0.6]
        \draw[thick, black] (0,-0.5) -- (0,0.5);
        \node at (0, -0.82) {$i$};
        \node at (0.35,0.05) {$\scriptstyle{\lambda}$};
    \end{tikzpicture}=
    \begin{tikzpicture}[baseline=-0.5ex, scale=0.6]
        \draw[thick, darkred, ->] (0,-0.5) -- (0,0.5);
        \node at (0.35,0.05) {$\scriptstyle{\lambda}$};
        \node[darkred] at (0, -0.82) {$i$};
    \end{tikzpicture}$, \hskip.5cm
    $\begin{tikzpicture}[baseline=-0.5ex,scale=0.6]
        \draw[thick, black] (0,-0.5) -- (0,0.5);
        \node at (0, -0.78) {$-i$};
        \node at (0.35,0.05) {$\scriptstyle{\lambda}$};
    \end{tikzpicture}
    =    \begin{tikzpicture}[baseline=-0.5ex, scale=0.6]
        \draw[thick, darkred, <-] (0,-0.5) -- (0,0.5);
        \node[darkred] at (0, -0.8) {$i$};
        \node at (0.35,0.05) {$\scriptstyle{\lambda}$};
    \end{tikzpicture}
   $,
   \item [(2)]  $ \begin{tikzpicture}[baseline=-0.5ex, scale=0.6]
        \draw[thick, black] (0,-0.5) -- (0,0.5);
        \node at (0, 0) {$\scriptstyle\bullet$};
        \node at (0.27,0.4) {$\scriptstyle{x}$};
        \node at (0.55,-0.2) {$\scriptstyle{\lambda}$};
        \node at (0, -0.79) {$i$};
    \end{tikzpicture}
    =
    \begin{tikzpicture}[baseline=-0.5ex, scale=0.6]
        \draw[thick, darkred, ->] (0,-0.5) -- (0,0.5);
        \fill[darkred]  (0, 0) circle (4pt);
        \node at (0.6,0) {$\scriptstyle{y+i}$};
        \node[darkred] at (0, -0.78) {$i$};
        \node at (1.7,-0.15) {$\scriptstyle{\lambda}$};
    \end{tikzpicture}$,
    \hskip.5cm $\begin{tikzpicture}[baseline=-0.5ex,scale=0.6]
        \draw[thick, black] (0,-0.5) -- (0,0.5);
        \node at (0, 0) {$\scriptstyle\bullet$};
        \node at (0.28,0) {$\scriptstyle{x}$};
        \node at (0.55,-0.22) {$\scriptstyle{\lambda}$};
        \node at (0, -0.7) {$-i$};
    \end{tikzpicture}
    =\begin{tikzpicture}[baseline=-0.5ex, scale=0.6]
        \draw[thick, darkred, <-] (0,-0.5) -- (0,0.5);
        \fill[darkred]  (0, 0) circle (4pt);
        \node at (0.74,0) {$\scriptstyle{-y-i}$};
        \node[darkred] at (0, -0.78) {$i$};
        \node at (1.45,-0.15) {$\scriptstyle{\lambda}$};
    \end{tikzpicture}$,
\item [(3)]
    $\begin{tikzpicture}[baseline=-0.5ex, scale=0.8]
        \draw[thick, black]
            (-0.4,0.5) .. controls (-0.4,-0.3) and (0.4,-0.3) .. (0.4,0.5);
        \node at (-0.4, 0.7) {$\scriptstyle{-i}$};
        \node at (0.4, 0.7) {$\scriptstyle{i}$};
        \node at (0.65,0.1) {$\scriptstyle{\lambda}$};
    \end{tikzpicture}
    =
    \begin{tikzpicture}[baseline=-0.5ex,scale=0.8]
        \draw[thick, darkred, ->]
            (-0.4,0.5) .. controls (-0.4,-0.3) and (0.4,-0.3) .. (0.4,0.5);
        \node[darkred] at (-0.4, 0.7) {$\scriptstyle{i}$};
        \node at (0.8,0) {$\scriptstyle{\lambda}$};
    \end{tikzpicture}
$,\item [(4)]
$ \begin{tikzpicture}[baseline=-0.5ex, scale=0.8, yshift=10pt]
        \draw[thick]
            (-0.4,-0.5) .. controls (-0.4,0.3) and (0.4,0.3) .. (0.4,-0.5);
        \node at (-0.4, -0.7) {$i$};
        \node at (0.4, -0.7) {$-i$};
        \node at (0.65,0.1) {$\scriptstyle{\lambda}$};
    \end{tikzpicture}=
    \begin{tikzpicture}[baseline=-0.5ex,scale=0.8, yshift=10pt]
        \draw[thick, darkred, ->]
            (-0.4,-0.5) .. controls (-0.4,0.3) and (0.4,0.3) .. (0.4,-0.5);
        \node[darkred] at (-0.4, -0.7) {$i$};
        \node at (0.65,0.05) {$\scriptstyle{\lambda}$};
    \end{tikzpicture}$,

\item[(5)]
$\begin{tikzpicture}[baseline=-0.5ex,scale=0.8]
    \draw[thick, black]
        (-0.4,0.5) .. controls (-0.4,-0.3) and (0.4,-0.3) .. (0.4,0.5);
    \node at (-0.4, 0.7) {$\scriptstyle{i}$};
    \node at (0.4, 0.7) {$\scriptstyle -i$};
    \node at (0.65,0.1) {$\scriptstyle{\lambda}$};
\end{tikzpicture}\begin{tikzpicture}[baseline = 0]
   \draw[-,darkg,thick] (0.65,-0.4) to (0.65,.4);
    \node at (0.65,-.6) {$\darkg\scriptstyle{V}$};
 \end{tikzpicture}= \begin{tikzpicture}[baseline=-0.5ex, scale=0.8]
    \draw[thick, darkred, <-]
        (-0.4,0.5) .. controls (-0.4,-0.3) and (0.4,-0.3) .. (0.4,0.5);
    \node[darkred] at (0.4, 0.7) {$\scriptstyle{i}$};
    \node at (-0.38, 0.28) {$\scriptstyle\bullet$};
    \node at (-1.6, 0.4) {$\frac{xm_{E_{-i}V, i}(x)}{n_{E_{-i}V, i}(x)} $};
    \node at (0.65,0.05) {$\scriptstyle{\lambda}$};
\end{tikzpicture}$,
 \item [(6)] $\begin{tikzpicture}[baseline=-0.5ex,scale=0.8, yshift=10pt]
        \draw[thick]
            (-0.4,-0.5) .. controls (-0.4,0.3) and (0.4,0.3) .. (0.4,-0.5);
        \node at (-0.4, -0.7) {$-i$};
        \node at (0.4, -0.70) {$i$};
        \node at (0.65,0.1) {$\scriptstyle{\lambda}$};
    \end{tikzpicture}\begin{tikzpicture}[baseline = 0]
   \draw[-,darkg,thick] (0.65,-0.4) to (0.65,.4);
    \node at (0.65,-.6) {$\darkg\scriptstyle{V}$};
 \end{tikzpicture}= \quad
    \begin{tikzpicture}[baseline=-0.5ex, scale=0.8, yshift=10pt]
        \draw[thick, darkred, <-]
            (-0.4,-0.5) .. controls (-0.4,0.3) and (0.4,0.3) .. (0.4,-0.5);
        \node[darkred] at (0.4, -0.7) {$\scriptstyle{i}$};
        \node at (0.4, -0.4) {$\scriptstyle\bullet$};
        \node at (1.5, -0.4) {$\scriptstyle\frac{n_{V, i}(x)}{xm_{V, i}(x)}$};
        \node at (0.65,0.05) {$\scriptstyle{\lambda}$};
    \end{tikzpicture}$.
\item [(7)]
  $  \begin{tikzpicture}[baseline=0,scale=0.8]
        \draw[thick, black] (-0.4,0.5) .. controls (-0.4,-0.3) and (0.4,-0.3) .. (0.4,0.5);
        \node at (-.4, .7) {$\scriptstyle{j}$};
        \node at (.4, .7) {$\scriptstyle{k}$};
        \node at (.5, .1){$\scriptstyle{\lambda}$};
    \end{tikzpicture}=0$, if $k\neq -j$, \item [(8)]
    $\begin{tikzpicture}[baseline=-0.5ex,scale=0.8,yshift=10pt]
        \draw[thick] (-0.4,-0.5) .. controls (-0.4,0.3) and (0.4,0.3) .. (0.4,-0.5);
        \node at (-0.4, -0.7) {$\scriptstyle k$};
        \node at (0.4, -0.7) {$\scriptstyle j$};
        \node at (.6, .1){$\scriptstyle{\lambda}$};
    \end{tikzpicture}=0$, if  $k\neq -j$, \item [(9)]
   $ \begin{tikzpicture}[baseline=0]
        \draw[-,thick] (0,-0.4) to (0,.4);
        \node at (0,-.57) {$\scriptstyle{k}$};
        \node at (0,.57) {$\scriptstyle{j}$};
        \node at (0.2,0) {$\scriptstyle{\lambda}$};
    \end{tikzpicture}=0 $  if $ k\neq j$.
    \item[(10)]
$\begin{tikzpicture}[baseline=-0.5ex,scale=0.8]
    \draw[thick, black]
        (-0.4,0.5) .. controls (-0.4,-0.3) and (0.4,-0.3) .. (0.4,0.5);
    \node at (-0.4, 0.7) {$\scriptstyle{0}$};
    \node at (0.4, 0.7) {$\scriptstyle 0$};
    \node at (0.65,0.1) {$\scriptstyle{\lambda}$};
\end{tikzpicture}\begin{tikzpicture}[baseline = 0]
   \draw[-,darkg,thick] (0.65,-0.4) to (0.65,.4);
    \node at (0.65,-.6) {$\darkg\scriptstyle{V}$};
 \end{tikzpicture}=
 \begin{tikzpicture}[baseline=-0.5ex, scale=0.8]
    \draw[thick, darkred, -]
        (-0.4,0.5) .. controls (-0.4,-0.3) and (0.4,-0.3) .. (0.4,0.5);
    \node[darkred] at (0.4, 0.7) {$\scriptstyle{0}$};
    \node at (0.38, 0.28) {$\scriptstyle\bullet$};
    \node at (1.9, 0.4) {$\scriptstyle{M_V(x)} $};
    \node at (0.65,0.05) {$\scriptstyle{\lambda}$};
\end{tikzpicture}$,
 \item [(11)] $\begin{tikzpicture}[baseline=-0.5ex,scale=0.8, yshift=10pt]
        \draw[thick]
            (-0.4,-0.5) .. controls (-0.4,0.3) and (0.4,0.3) .. (0.4,-0.5);
        \node at (-0.4, -0.7) {$0$};
        \node at (0.4, -0.70) {$0$};
        \node at (0.65,0.1) {$\scriptstyle{\lambda}$};
    \end{tikzpicture}
    \begin{tikzpicture}[baseline = 0]
   \draw[-,darkg,thick] (0.65,-0.4) to (0.65,.4);
    \node at (0.65,-.6) {$\darkg\scriptstyle{V}$};
 \end{tikzpicture}= \quad
    \begin{tikzpicture}[baseline=-0.5ex, scale=0.8, yshift=10pt]
        \draw[thick, darkred, -]
            (-0.4,-0.5) .. controls (-0.4,0.3) and (0.4,0.3) .. (0.4,-0.5);
        \node[darkred] at (0.4, -0.7) {$\scriptstyle{0}$};
        \node at (0.4, -0.4) {$\scriptstyle\bullet$};
        \node at (1.5, -0.4) {$\scriptstyle W_V(x)$};
        \node at (0.65,0.05) {$\scriptstyle{\lambda}$};
    \end{tikzpicture}$.
\end{itemize}
\end{Defn}

{The fractions in (5)--(6) are regular at the indicated nonzero
color. Their evaluations are natural by the argument following
\cite[Definition~7.11]{GRS}, using uniqueness of the monic factors
$t_{V,i}$. At zero, naturality follows from
Lemma~\ref{lem:ir-normalization}.

For the black bubble with $r$ dots, write $\omega_{V,r}$ for its}
value on $V$, and put
\[
 O_V(u)=u-\frac12+\sum_{r\ge0}\omega_{V,r}u^{-r}.
\]
\begin{Lemma}\label{lem:red-black-series}
For every finitely generated object $V\in\M_\lambda$,
\begin{equation}\label{eq:red-black-main}
 \widetilde O_V(u)
 \left(u-\frac12(-1)^{a_\lambda}\right)=O_V(u).
\end{equation}
Equivalently,
\[
 O_V(u)=\frac{n_V(u)}{m_V(u)}.
\]
\end{Lemma}

\begin{proof}
Follow \cite[Proposition~7.20]{GRS}, taking
$S=\{0\}\cup\{i\ne0\mid E_iV\ne0\}$ and a sufficiently large
$k>\max_i\varepsilon_i(V)$. For $i\in S\setminus\{0\}$ choose
\begin{equation}\label{eq:bezout-nonzero}
 u\,m_{V,i}(u)d_i(u)+(u-i)^kh_i(u)=1,
\end{equation}
and at zero choose
\begin{equation}\label{eq:bezout-zero}
 m_{V,0}(u)d_0(u)+u^kh_0(u)=1.
\end{equation}
The B\'ezout identities imply
\[
 \frac1{u m_V(u)}
 -\sum_{i\in S\setminus\{0\}}\frac{d_i(u)}{(u-i)^{\eps_i(V)}}
 -\frac{d_0(u)}{u^{\eps_0(V)+1}}\in Z_V[u].
\]
{At zero the difference is regular because $k\ge\eps_0(V)+1$;
at $i\ne0$ it is regular because $k>\eps_i(V)$. Thus it has no
principal part. Multiplication by $n_V$ adds only a polynomial,
which contributes no negative Laurent coefficient.}
For $i\ne0$, the bubble calculations \cite[(7.14)--(7.16)]{GRS},
using Lemma~\ref{lem:tVi-integral} and
\eqref{eq:ir-primary-transport}, give
\begin{equation}\label{eq:black-bubble-nonzero}
\begin{tikzpicture}[baseline = 1.25mm]
	\draw[-] (0,0.4) to[out=180,in=90] (-.2,0.2);
	\draw[-] (0.2,0.2) to[out=90,in=0] (0,.4);
	\draw[-] (-.2,0.2) to[out=-90,in=180] (0,0);
	\draw[-] (0,0) to[out=0,in=-90] (0.2,0.2);
	\node at (0.2,0.2) {$\scriptstyle\bullet$};
    \node at (0.28,-0.05) {$\scriptstyle i$};
    \node at (-0.45,0.2) {$\scriptstyle -i$};
	\node at (0.45,0.2) {$\scriptstyle{x^r}$};
\end{tikzpicture}
=\Coeff{\dfrac{n_V(u)d_i(u)}{(u-i)^{\eps_i(V)}}}{u:-r-1}
 \qquad(r\geq0).
\end{equation}
At zero, the two $\Gamma_V$ factors cancel. By
\eqref{eq:O-lambda-zero} and \eqref{eq:zero-ratio-exact},
\begin{equation}\label{eq:black-bubble-zero}
\begin{tikzpicture}[baseline = 1.25mm]
	\draw[-] (0,0.4) to[out=180,in=90] (-.2,0.2);
	\draw[-] (0.2,0.2) to[out=90,in=0] (0,.4);
	\draw[-] (-.2,0.2) to[out=-90,in=180] (0,0);
	\draw[-] (0,0) to[out=0,in=-90] (0.2,0.2);
	\node at (0.2,0.2) {$\scriptstyle\bullet$};
    \node at (0.28,-0.05) {$\scriptstyle 0$};
	\node at (0.45,0.2) {$\scriptstyle{x^r}$};
\end{tikzpicture}
 =\Coeff{\dfrac{n_V(u)d_0(u)}{u^{\eps_0(V)}}}{u:-r}
 +\delta_{r,0}\frac12.
\end{equation}
The additional constant is
$-(-1)^{c(\lambda)}n_{V,0}(0)/m_{V,0}(0)=1/2$.
Summing these identities identifies the negative and constant
coefficients, and both sides have leading term $u$.
\end{proof}

\subsection{Reverse crossings}
\label{subsec:integral-reverse-crossings}
\begin{Defn}\label{cross-red123} Suppose   $\lambda\in X^\imath$, %
and $i,j\in\mathbb Z_{\ge0}$. In (1), exclude $(i,j)=(0,0)$.
In a cup--cap decoration joining the two components, $x_1$ is the
dot on the left outgoing leg and $x_2$ the dot on the left incoming
leg; thus $x_1+x_2$ has residue $k+l$ in (3) and $2k$ in (4).
All sums are evaluated componentwise. Define
\begin{itemize}
\item [(1)] %
$
\begin{tikzpicture}[baseline = 0, scale=0.7]
	\draw[-] (0.38,-.4) to (-0.38,.5);
	\draw[-] (-0.38,-.4) to (0.38,.5);
      \node at (0,0.05) {$\diamond$};
   \node at (-.4,-.62) {$\scriptstyle{j}$};
   \node at (.4,-.62) {$\scriptstyle{i}$};
   \node at (-.4,.72) {$\scriptstyle{i}$};
   \node at (.4,.72) {$\scriptstyle{j}$};
    \node at (0.7,0) {$\scriptstyle{\lambda}$};
\end{tikzpicture}
=
\left\{
\begin{array}{ll}
\begin{tikzpicture}[baseline = 0]
	\draw[-,thick, darkred] (0.28,-.3) to (-0.28,.4);
	\draw[-,thick, darkred] (-0.28,-.3) to (0.28,.4);
   \node at (-.3,-.45) {$\scriptstyle{j}$};
   \node at (.3,-.45) {$\scriptstyle{i}$};
   \node at (0.5,0.05) {$\scriptstyle{\lambda}$};
      \draw[-,densely dotted] (0.2,-.24) to (-0.19,-.24);
      \node at (1.3,-0.2) {$\scriptstyle{x_1-x_2+1}$};
\end{tikzpicture}
-\begin{tikzpicture}[baseline = 0mm, scale=0.8]
	\draw[-,darkred] (0.58,-.4) to (0.58,.5);
	\draw[-,darkred] (0.08,-.4) to (0.08,.5);
      \node at (0.08,-0.6) {$\scriptstyle{i}$};
      \node at (0.58,-0.6) {$\scriptstyle{i}$};
 \node at (0.65,0) {$\scriptstyle{\lambda}$};
\end{tikzpicture}
&\text{if $j=i\ne0$,}\\
\begin{tikzpicture}[baseline = 0,scale=0.86]
	\draw[-,darkred] (0.38,-.4) to (-0.38,.5);
	\draw[-,darkred] (-0.38,-.4) to (0.38,.5);
   \node at (-.4,-.55) {$\scriptstyle{i+1}$};
   \node at (-.4,.65) {$\scriptstyle{i}$};
   \node at (.4,.65) {$\scriptstyle{i+1}$};
    \node at (0.5,0.09) {$\scriptstyle{\lambda}$};
   \node at (.4,-.55) {$\scriptstyle{i}$};
      \draw[-,densely dotted] (0.2,-.24) to (-0.19,-.24);
      \node at (1.3,-0.2) {$\scriptstyle{(x_1-x_2)^{-1}}$};
\end{tikzpicture}
&\text{if $j=i+1$,}\\
-\begin{tikzpicture}[baseline = 0, scale=0.8]
	\draw[-,darkred] (0.38,-.4) to (-0.38,.5);
	\draw[-,darkred] (-0.38,-.4) to (0.38,.5);
	\draw[-,densely dotted] (0.18,-.24) to (-0.19,-.24);
   \node at (-.4,-.62) {$\scriptstyle{j}$};
   \node at (.4,-.62) {$\scriptstyle{i}$};
   \node at (.4,.72) {$\scriptstyle{j}$};
   \node at (-.4,.72) {$\scriptstyle{i}$};
      \node at (1.6,-0.2) {$\scriptstyle{\frac{x_1-x_2-1} {x_1-x_2}}$};
        \node at (0.5,0.09) {$\scriptstyle{\lambda}$};
\end{tikzpicture}
&\text{if $j \neq i,i+1$.}
\end{array}\right.
$
\item[(2)] The three formulas below apply for $i>0$, $j>0$, and $i,j>0$, respectively.
\[
\begin{aligned}
\begin{tikzpicture}[baseline = 0, scale=0.8 ]
	\draw[-] (0.38,-.4) to (-0.38,.5);
	\draw[-] (-0.38,-.4) to (0.38,.5);
      \node at (0,0.05) {$\diamond$};
   \node at (-.4,-.62) {$\scriptstyle{j}$};
   \node at (.4,-.6) {$\scriptstyle{-i}$};
   \node at (-.4,.68) {$\scriptstyle{-i}$};
   \node at (.4,.72) {$\scriptstyle{j}$};
     \node at (0.5,0.09) {$\scriptstyle{\lambda}$};
\end{tikzpicture}&=
\mathord{
\begin{tikzpicture}[baseline = 0,scale=0.7,  ]
	\draw[-,thick] (0.3,-.5) to (-0.3,.5);
	\draw[-,thick] (-0.2,-.2) to (0.2,.3);
        \draw[-,thick] (0.2,.3) to[out=50,in=180] (0.5,.5);
        \draw[-,thick] (0.5,.5) to[out=0,in=90] (0.8,-.5);
        \draw[-,thick] (-0.2,-.2) to[out=230,in=0] (-0.5,-.5);
        \draw[-,thick] (-0.5,-.5) to[out=180,in=-90] (-0.8,.5);
        \node at (0,0.05) {$\diamond$};
  \node at (-0.8,.7) {$\scriptstyle{-i}$};
   \node at (0,.7) {$\scriptstyle{j}$};
   \node at (0.28,-.69) {$\scriptstyle{j}$};
     \node at (0.68,-.69) {$\scriptstyle{-i}$};
   \node at (1.05,.05) {$\scriptstyle{\lambda}$};
\end{tikzpicture}
},\\[1.5ex]
\begin{tikzpicture}[baseline = 0, scale=0.8 ]
	\draw[-] (0.38,-.4) to (-0.38,.5);
	\draw[-] (-0.38,-.4) to (0.38,.5);
      \node at (0,0.05) {$\diamond$};
   \node at (-.4,-.62) {$\scriptstyle{-j}$};
   \node at (.4,-.6) {$\scriptstyle{i}$};
   \node at (-.4,.7) {$\scriptstyle{i}$};
   \node at (.4,.72) {$\scriptstyle{-j}$};
    \node at (0.5,0.09) {$\scriptstyle{\lambda}$};
\end{tikzpicture}&=
\mathord{
\begin{tikzpicture}[baseline = 0, yscale=-1, scale=0.8]
	\draw[-,thick] (0.3,-.5) to (-0.3,.5);
	\draw[-,thick] (-0.2,-.2) to (0.2,.3);
        \draw[-,thick] (0.2,.3) to[out=50,in=180] (0.5,.5);
        \draw[-,thick] (0.5,.5) to[out=0,in=90] (0.8,-.5);
        \draw[-,thick] (-0.2,-.2) to[out=230,in=0] (-0.5,-.5);
        \draw[-,thick] (-0.5,-.5) to[out=180,in=-90] (-0.8,.5);
        \node at (0,0.05) {$\diamond$};
  \node at (-0.8,.65) {$\scriptstyle{-j}$};
   \node at (0.28,-.62) {$\scriptstyle{i}$};
    \node at (-0.1,.65) {$\scriptstyle{i}$};
   \node at (0.78,-.62) {$\scriptstyle{-j}$};
   \node at (1.05,.05) {$\scriptstyle{\lambda}$};
\end{tikzpicture}
},\\[1.5ex]
\begin{tikzpicture}[baseline = 0,scale=0.8 ]
	\draw[-] (0.38,-.4) to (-0.38,.5);
	\draw[-] (-0.38,-.4) to (0.38,.5);
      \node at (0,0.05) {$\diamond$};
   \node at (-.4,-.55) {$\scriptstyle{-j}$};
   \node at (.4,-.55) {$\scriptstyle{-i}$};
   \node at (-.4,.65) {$\scriptstyle{-i}$};
   \node at (.4,.65) {$\scriptstyle{-j}$};
  \node at (0.5,0.09) {$\scriptstyle{\lambda}$};
\end{tikzpicture}&=\mathord{
\begin{tikzpicture}[baseline = 0,scale=0.4 ]
\draw[-,thick] (-1.3,.4) to (-1.3,-1.2);
\draw[-,thick] (1.3,-.4) to (1.3,1.2);
\draw[-,thick] (-.5,1.1) to [out=180,in=90] (-1.3,.4);
\draw[-,thick] (.35,.4) to [out=90,in=0] (-.5,1.1);
\draw[-,thick] (.5,-1.1) to [out=0,in=-90] (1.3,-.4);
\draw[-,thick] (-.35,-.4) to [out=-90,in=180] (.5,-1.1);
\draw[-,thick] (-.35,-.4) to [out=90,in=-90] (.35,.4);
        \draw[-,thick] (0.35,-.5) to[out=180,in=0] (-0.35,.5);
        \draw[-,thick] (-0.35,.5) to[out=180,in=90] (-0.8,-1.2);
        \draw[-,thick] (0.35,-.5) to[out=0,in=-90] (0.8,1.2);
        \node at (0,0) {$\diamond$};
          \node at (-0.65,-1.6) {$\scriptstyle{-i}$};
  \node at (-1.55,-1.6) {$\scriptstyle{-j}$};
  \node at (0.65,1.6) {$\scriptstyle{-i}$};
  \node at (1.55,1.6) {$\scriptstyle{-j}$};
   \node at (1.62,0) {$\scriptstyle{\lambda}$};
\end{tikzpicture}
}.
\end{aligned}
\]
\item[(3)]  $ \begin{tikzpicture}[baseline = 0, scale=0.8]
	\draw[-] (0.38,-.4) to (-0.38,.5);
	\draw[-] (-0.38,-.4) to (0.38,.5);
      \node at (0,0.05) {$\diamond$};
   \node at (-.4,-.55) {$\scriptstyle{k}$};
   \node at (.4,-.55) {$\scriptstyle{l}$};
   \node at (-.4,.65) {$\scriptstyle{k}$};
   \node at (.4,.65) {$\scriptstyle{l}$};
  \node at (0.5,0.09) {$\scriptstyle{\lambda}$};
\end{tikzpicture}=
\begin{tikzpicture}[baseline = 0mm, scale=0.8 ]
	\draw[-,densely dotted] (0.14,.06) to (0.5,.06);
	\draw[-] (0.58,-.4) to (0.58,.5);
	\draw[-] (0.08,-.4) to (0.08,.5);
      \node at (0.08,-0.55) {$\scriptstyle{k}$};
      \node at (0.58,-0.55) {$\scriptstyle{l}$};
      \node at (1.6,0.1) {$\scriptstyle{(x_1-x_2)^{-1}}$};
        \node at (1.,-0.22) {$\scriptstyle{\lambda}$};
\end{tikzpicture}$,
$ \begin{tikzpicture}[baseline = 0, scale=0.8, ]
	\draw[-] (0.38,-.4) to (-0.38,.5);
	\draw[-] (-0.38,-.4) to (0.38,.5);
      \node at (0,0.05) {$\diamond$};
   \node at (-.4,-.55) {$\scriptstyle{l}$};
   \node at (.4,-.55) {$\scriptstyle{-l}$};
   \node at (-.4,.65) {$\scriptstyle{k}$};
   \node at (.4,.65) {$\scriptstyle{-k}$};
  \node at (0.5,0.09) {$\scriptstyle{\lambda}$};
\end{tikzpicture}=-~
\begin{tikzpicture}[baseline = -1mm]
    \draw[-] (0,0.35) to[out=-90,in=-90] (0.6,0.35);
    \draw[-] (0,-0.35) to[out=90,in=90] (0.6,-0.35);
\node at (-.18,-0.45) {$\scriptstyle l$};
    \node at (0.68,-0.45) {$\scriptstyle -l$};
    \node at (-.18,0.45) {$\scriptstyle k$};
    \node at (0.68, 0.45) {$\scriptstyle -k$};
    \draw[densely dotted,-] (0.2,-0.20) -- (0.2,0.18);
       \node at (1.,0.1) {$\scriptstyle{\lambda}$};
    \node at (1.45,-0.1) {$\scriptstyle (x_1+x_2)^{-1}$};
\end{tikzpicture}$  for all $k, l\in I$ such that  $k\neq \pm l$.
\item[(4)] %
$\begin{tikzpicture}[baseline = 0,scale=0.8, ]
	\draw[-] (0.38,-.4) to (-0.38,.5);
	\draw[-] (-0.38,-.4) to (0.38,.5);
      \node at (0,0.05) {$\diamond$};
   \node at (-.4,-.55) {$\scriptstyle{k}$};
   \node at (.4,-.55) {$\scriptstyle{-k}$};
   \node at (-.4,.65) {$\scriptstyle{k}$};
   \node at (.4,.65) {$\scriptstyle{-k}$};
     \node at (0.5,0.09) {$\scriptstyle{\lambda}$};
  \node at (0.5,0.09) {$\scriptstyle{\lambda}$};
\end{tikzpicture}= \begin{tikzpicture}[baseline = 0mm]
	\draw[-,densely dotted] (0.14,.06) to (0.5,.06);
	\draw[-] (0.58,-.4) to (0.58,.5);
	\draw[-] (0.08,-.4) to (0.08,.5);
      \node at (0.08,-0.55) {$\scriptstyle{k}$};
      \node at (0.58,-0.55) {$\scriptstyle{-k}$};
      \node at (1.6,0.1) {$\scriptstyle{(x_1-x_2)^{-1}}$};
        \node at (1.,-0.2) {$\scriptstyle{\lambda}$};
\end{tikzpicture}-~
\begin{tikzpicture}[baseline = -1mm]
    \draw[-] (0,0.35) to[out=-90,in=-90] (0.6,0.35);
    \draw[-] (0,-0.35) to[out=90,in=90] (0.6,-0.35);
\node at (-.18,-0.45) {$\scriptstyle k$};
    \node at (0.68,-0.45) {$\scriptstyle -k$};
    \node at (-.18,0.45) {$\scriptstyle k$};
    \node at (0.68, 0.45) {$\scriptstyle -k$};
    \draw[densely dotted,-] (0.2,-0.20) -- (0.2,0.18);
       \node at (1.,0.2) {$\scriptstyle{\lambda}$};
    \node at (1.45,-0.1) {$\scriptstyle (x_1+x_2)^{-1}$};
\end{tikzpicture}$ for all $k\in I\setminus\{0\}$.
\item[(5)]
$\begin{tikzpicture}[baseline = 0, scale=0.8,  ]
	\draw[-] (0.38,-.4) to (-0.38,.5);
	\draw[-] (-0.38,-.4) to (0.38,.5);
      \node at (0,0.05) {$\diamond$};
   \node at (-.4,-.55) {$\scriptstyle{t}$};
   \node at (.4,-.55) {$\scriptstyle{c}$};
   \node at (-.4,.65) {$\scriptstyle{k}$};
   \node at (.4,.65) {$\scriptstyle{l}$};
   \node at (0.7,0) {$\scriptstyle{\lambda}$};
\end{tikzpicture}=0$, when    $\{c,t\}\neq \{k,l\}$ and either  $c\neq -t $ or $k\neq -l$,   where  $c,t,k,l\in I$.
\item[(6)]
   $\begin{tikzpicture}[baseline=-0.5ex,scale=0.8 ]
        \draw[thick, black, ] (0,-0.5) -- (0,0.5);
\end{tikzpicture}=  \sum_{\mu\in X^\imath}
    \begin{tikzpicture}[baseline=-0.5ex,scale=0.8,  ]
        \draw[thick, black, ] (0,-0.5) -- (0,0.5);
     \node at (.2,.05) {$\scriptstyle{\mu}$};
\end{tikzpicture}$, where  $\begin{tikzpicture}[baseline=-0.5ex,scale=0.8,  ]
        \draw[thick, black, ] (0,-0.5) -- (0,0.5);
     \node at (.2,.05) {$\scriptstyle{\lambda}$};
\end{tikzpicture}=  \sum_{i\in I}
    \begin{tikzpicture}[baseline=-0.5ex,scale=0.8,  ]
        \draw[thick, black, ] (0,-0.5) -- (0,0.5);
    \node at (0, -0.7) {$i$};
     \node at (.2,.05) {$\scriptstyle{\lambda}$};
\end{tikzpicture}$.
 \item [(7)]
$\begin{tikzpicture}[baseline=-0.5ex, scale=0.8 ]
        \draw[thick, black, ] (0,-0.5) -- (0,0.5);
        \node at (0, 0) {$\scriptstyle\bullet$};
        \node at (0.2,0) {$\scriptstyle{x}$};
        \end{tikzpicture} =\sum_{\mu\in X^\imath}
        \begin{tikzpicture}[baseline=-0.5ex,scale=0.8,  ]
        \draw[thick, black, ] (0,-0.5) -- (0,0.5);
        \node at (0, 0) {$\scriptstyle\bullet$};
        \node at (0.2,0) {$\scriptstyle{x}$};
          \node at (.3,-.3) {$\scriptstyle{\mu}$};
        \end{tikzpicture}
 $, where
 $\begin{tikzpicture}[baseline={(0, -0.1)},scale=0.8 ]
        \draw[thick, black, ] (0,-0.5) -- (0,0.5);
        \node at (0, 0) {$\scriptstyle\bullet$};
        \node at (0.2,0) {$\scriptstyle{x}$};
          \node at (.2,-.3) {$\scriptstyle{\lambda}$};
        \end{tikzpicture}= \sum_{i\in I}  \begin{tikzpicture}[baseline=-0.5ex,scale=0.8,  ]
        \draw[thick, black, ] (0,-0.5) -- (0,0.5);
        \node at (0, 0) {$\scriptstyle\bullet$};
        \node at (0.2,0) {$\scriptstyle{x}$};
          \node at (.2,-.3) {$\scriptstyle{\lambda}$};
        \node at (0, -0.7) {$i$};
        \end{tikzpicture}$.
\item [(8)] $\begin{tikzpicture}[baseline=0,scale=0.8 ]
             \draw[thick, black]
            (-0.4,0.5) .. controls (-0.4,-0.3) and (0.4,-0.3) .. (0.4,0.5);
\end{tikzpicture}=\sum_{\mu\in X^\imath}  \begin{tikzpicture}[baseline=0,scale=0.8, ]
             \draw[thick, black]
            (-0.4,0.5) .. controls (-0.4,-0.3) and (0.4,-0.3) .. (0.4,0.5);
              \node at (.5, .1){$\scriptstyle{\mu}$};
\end{tikzpicture}
            $, where     %
    $\begin{tikzpicture}[baseline=0,scale=0.8 ]
             \draw[thick, black]
            (-0.4,0.5) .. controls (-0.4,-0.3) and (0.4,-0.3) .. (0.4,0.5);
              \node at (.5, .1){$\scriptstyle{\lambda}$};
\end{tikzpicture}=\sum_{i\in I}\begin{tikzpicture}[baseline=0,scale=0.8, ]
             \draw[thick, black]
            (-0.4,0.5) .. controls (-0.4,-0.3) and (0.4,-0.3) .. (0.4,0.5);
            \node at (-.4, .7) {$\scriptstyle{-i}$};
             \node at (.4, .7) {$\scriptstyle{i}$};   \node at (.5, .1){$\scriptstyle{\lambda}$};
\end{tikzpicture}     $
\item [(9)] $ \begin{tikzpicture}[baseline=-0.5ex,scale=0.8,  yshift=10pt]
        \draw[thick, ]
            (-0.4,-0.5) .. controls (-0.4,0.3) and (0.4,0.3) .. (0.4,-0.5);
\end{tikzpicture}=\sum_{\mu\in X^\imath}  \begin{tikzpicture}[baseline=-0.5ex,scale=0.8, yshift=10pt]
        \draw[thick, ]
            (-0.4,-0.5) .. controls (-0.4,0.3) and (0.4,0.3) .. (0.4,-0.5);
               \node at (.6, .1){$\scriptstyle{\mu}$};
\end{tikzpicture}$, where     %
   $ \begin{tikzpicture}[baseline=-0.5ex,scale=0.8,  yshift=10pt]
        \draw[thick, ]
            (-0.4,-0.5) .. controls (-0.4,0.3) and (0.4,0.3) .. (0.4,-0.5);
               \node at (.6, .1){$\scriptstyle{\lambda}$};
\end{tikzpicture}=\sum_{i\in I}\begin{tikzpicture}[baseline=-0.5ex,scale=0.8,  yshift=10pt]
        \draw[thick, ]
            (-0.4,-0.5) .. controls (-0.4,0.3) and (0.4,0.3) .. (0.4,-0.5);
                \node at (-0.4, -0.7) {$\scriptstyle i$};
  \node at (0.4, -0.7) {$\scriptstyle -i$}; \node at (.6, .1){$\scriptstyle{\lambda}$};
\end{tikzpicture}$.

\item[(10)] $\begin{tikzpicture}[baseline = 0, scale=0.8 ]
	\draw[-] (0.38,-.4) to (-0.38,.5);
	\draw[-] (-0.38,-.4) to (0.38,.5);
\end{tikzpicture}=\sum_{\mu\in X^\imath}\begin{tikzpicture}[baseline = 0,scale=0.8 ]
	\draw[-] (0.38,-.4) to (-0.38,.5);
	\draw[-] (-0.38,-.4) to (0.38,.5);
   \node at (0.7,0) {$\scriptstyle{\mu}$};
\end{tikzpicture}$, where $\begin{tikzpicture}[baseline = 0,scale=0.8, ]
	\draw[-] (0.38,-.4) to (-0.38,.5);
	\draw[-] (-0.38,-.4) to (0.38,.5);
   \node at (0.7,0) {$\scriptstyle{\lambda}$};
\end{tikzpicture}=\sum_{i, j, k, l\in I}\begin{tikzpicture}[baseline = 0,scale=0.8, ]
	\draw[-] (0.38,-.4) to (-0.38,.5);
	\draw[-] (-0.38,-.4) to (0.38,.5);
      \node at (0,0.05) {$\diamond$};
   \node at (-.4,-.55) {$\scriptstyle{i}$};
   \node at (.4,-.55) {$\scriptstyle{j}$};
   \node at (-.4,.65) {$\scriptstyle{k}$};
   \node at (.4,.65) {$\scriptstyle{l}$};
   \node at (0.7,0) {$\scriptstyle{\lambda}$};
\end{tikzpicture}$.
\end{itemize}
\end{Defn}

\begin{Defn}\label{rem:integral-zero-crossing}
{On each horizontal pair of marked points, read $x_1,x_2$ from
left to right. Put $z=(x_2-x_1)/2$ on $E_0^2V$, and define}
\begin{equation}\label{eq:ir-zero-crossing-correction}
 h_V(x_1,x_2)=\frac{1-M_V(z)W_V(z)}{x_1-x_2}.
\end{equation}
Since $M_V(0)W_V(0)=1$, this quotient is a formal power series:
the numerator is divisible by $z$, and $x_1-x_2=-2z$.
Its nilpotent evaluation is understood as in
Lemma~\ref{lem:ir-normalization}. The zero crossing is
\begin{equation}\label{eq:ir-zero-crossing-definition}
 \IRcross{black}{0}{0}
 =
 \begin{tikzpicture}[baseline=0,scale=.72]
  \draw[darkred,thick](-.42,-.60)--(.42,.60)
   (.42,-.60)--(-.42,.60);
  \draw[darkred,densely dotted](-.43,-.42)--(.43,-.42)
   (-.43,.42)--(.43,.42);
  \node[anchor=west] at(.52,-.42){$\scriptstyle W_V(z)$};
  \node[anchor=west] at(.52,.42){$\scriptstyle M_V(z)$};
  \node at(-.42,-.8){$\scriptstyle0$};
  \node at(.42,-.8){$\scriptstyle0$};
 \end{tikzpicture}
 \ +\
 \begin{tikzpicture}[baseline=0,scale=.72]
  \draw[darkred,thick](-.34,-.60)--(-.34,.60)
   (.34,-.60)--(.34,.60);
  \draw[darkred,densely dotted](-.38,0)--(.38,0);
  \node[anchor=west] at(.48,0){$\scriptstyle h_V(x_1,x_2)$};
  \node at(-.34,-.8){$\scriptstyle0$};
  \node at(.34,-.8){$\scriptstyle0$};
 \end{tikzpicture}.
\end{equation}
The factor $M_V$ lies above the crossing and $W_V$ below it.
\end{Defn}

{Multiply Definition~\ref{eq:red-cup}(3) on the left by $M_V(z)$
and on the right by $W_V(z)$. The parallel term becomes
$-M_V(z)g_V(x_1,x_2)W_V(z)=h_V(x_1,x_2)$ by
\eqref{eq:gV-definition}, giving \eqref{eq:ir-zero-crossing-definition}.}

\newcommand{\IRDcupdeco}[2]{%
 \begin{tikzpicture}[baseline=0,scale=.72]
  \draw[darkred,thick](-.42,.64)--(-.42,.12)
   ..controls(-.42,-.43)and(.42,-.43)..(.42,.12)--(.42,.64);
  \fill[darkred]({#1*.42},.38)circle(1.6pt);
  \ifnum#1=1
   \node[anchor=west]at(.56,.38){$\scriptstyle #2$};
  \else
   \node[anchor=east]at(-.56,.38){$\scriptstyle #2$};
  \fi
  \node at(-.42,.86){$\scriptstyle0$};
  \node at(.42,.86){$\scriptstyle0$};
 \end{tikzpicture}}
\newcommand{\IRDcapdeco}[1]{%
 \begin{tikzpicture}[baseline=0,scale=.72]
  \draw[darkred,thick](-.42,-.64)--(-.42,-.12)
   ..controls(-.42,.43)and(.42,.43)..(.42,-.12)--(.42,-.64);
  \fill[darkred](.42,-.38)circle(1.6pt);
  \node[anchor=west]at(.56,-.38){$\scriptstyle #1$};
  \node at(-.42,-.86){$\scriptstyle0$};
  \node at(.42,-.86){$\scriptstyle0$};
 \end{tikzpicture}}

\newcommand{\IRDparallelcoef}[1]{%
 \begin{tikzpicture}[baseline=0,scale=.72]
  \draw[darkred,thick](-.36,-.58)--(-.36,.58)
   (.36,-.58)--(.36,.58);
  \draw[darkred,densely dotted](-.36,0)--(.36,0);
  \node[anchor=west]at(.5,0){$\scriptstyle #1$};
  \node at(-.36,-.8){$\scriptstyle0$};
  \node at(.36,-.8){$\scriptstyle0$};
  \node at(-.36,.8){$\scriptstyle0$};
  \node at(.36,.8){$\scriptstyle0$};
 \end{tikzpicture}}
\subsection{Relations on zero-color strands}
\label{subsec:integral-zero-relations}

Put $c=c(\lambda)$; all regions have weight $\lambda$ since
$\bar\alpha_0=0$. Red and black diagrams denote the given and
reconstructed generators, respectively, evaluated on $V$.

For reference, the relations from Definition~\ref{nil} used below are
repeated in the same diagrammatic notation:
\begin{align}
 \IRsnake{darkred}{0}{1}
 &=\IRstrand{darkred}{0}
 =\IRsnake{darkred}{0}{-1},
 &&\label{eq:ir-red-snake}\\[3pt]
 \IRdotcup{darkred}{0}{0}{(-.3854592,.20128)}
 &=-\IRdotcup{darkred}{0}{0}{(.3854592,.20128)},
 &
 \IRdotcap{darkred}{0}{0}{(-.3854592,-.20128)}
 &=-\IRdotcap{darkred}{0}{0}{(.3854592,-.20128)},
 \label{eq:ir-red-dots}\\[3pt]
 \IRcrossdot{darkred}{0}{0}{(-.34,.42)}
 -\IRcrossdot{darkred}{0}{0}{(.34,-.42)}
 &=\IRparallel{darkred}{0}{0}-\IRcupcap{darkred}{0}{0},
 &&\label{eq:ir-red-affine-left}\\[-1pt]
 \IRcrossdot{darkred}{0}{0}{(-.34,-.42)}
 -\IRcrossdot{darkred}{0}{0}{(.34,.42)}
 &=\IRparallel{darkred}{0}{0}-\IRcupcap{darkred}{0}{0},
 &&\label{eq:ir-red-affine-right}\\[3pt]
 \IRcrosscap{darkred}{0}{0}&=0,
 &\IRbubble{darkred}{0}&=c\,1_\lambda,
 \label{eq:ir-red-curl-bubble}\\[3pt]
 \IRcapcross{darkred}{0}{0}{1}
 &=\IRcapcross{darkred}{0}{0}{-1},
 &
 \IRcupcrossside{darkred}{0}{0}{1}
 &=\IRcupcrossside{darkred}{0}{0}{-1},
 \label{eq:ir-red-mates}\\[3pt]
 \IRdouble{darkred}{0}{0}&=0,
 &
 \IRbraid{darkred}{0}{0}{0}{1}
 &=\IRbraid{darkred}{0}{0}{0}{-1}.
 \label{eq:ir-red-quadratic-braid}
\end{align}
The second dot-crossing identity is the mate of the first:
bending both strands uses \eqref{eq:ir-red-snake},
\eqref{eq:ir-red-dots} and \eqref{eq:ir-red-mates}, with one minus
sign for each bent dot.

\begin{Lemma}
\label{lem:ir-zero-adjunction}
The reverse cup, cap and dot satisfy
\begin{align}
 \IRsnake{black}{0}{1}
 &=\IRstrand{black}{0}
 =\IRsnake{black}{0}{-1},
 \label{eq:ir-zero-snake-diagram}\\[3pt]
 \IRdotcup{black}{0}{0}{(-.3854592,.20128)}
 &=-\IRdotcup{black}{0}{0}{(.3854592,.20128)},
 &
 \IRdotcap{black}{0}{0}{(-.3854592,-.20128)}
 &=-\IRdotcap{black}{0}{0}{(.3854592,-.20128)}.
 \label{eq:ir-zero-dot-diagram}
\end{align}
\end{Lemma}
\begin{proof}
The red dot bend, \eqref{eq:ir-W-zero-transport}, and
\eqref{eq:ir-WM-products} give
\begin{equation}\label{eq:ir-cup-unit-slide}
 \IRDcupdeco{-1}{W_{E_0V}(x)}
 =
 \IRDcupdeco{1}{\dfrac{W_V(-x)}{1+2x}}
 =
 \IRDcupdeco{1}{M_V(x)^{-1}}.
\end{equation}
Thus the black cup can also be written
\[
 \IRcup{black}{0}{0}
 =\IRDcupdeco{1}{M_V(x)}
 =\IRDcupdeco{-1}{W_{E_0V}(x)^{-1}}.
\]
In each triangle identity, naturality and the dot bend move the cup
and cap decorations to the same strand, where they are inverse.
The red triangle identities then give
\eqref{eq:ir-zero-snake-diagram}. The dot identities follow from
\eqref{eq:ir-red-dots}, since the decorations commute with the dots.
\end{proof}

\begin{Lemma}
\label{lem:ir-zero-affine}
One has
\begin{align}
 \IRcrossdot{black}{0}{0}{(-.34,.42)}
 -\IRcrossdot{black}{0}{0}{(.34,-.42)}
 &=\IRparallel{black}{0}{0}-\IRcupcap{black}{0}{0},
 \label{eq:ir-zero-affine-left}\\[-1pt]
 \IRcrossdot{black}{0}{0}{(-.34,-.42)}
 -\IRcrossdot{black}{0}{0}{(.34,.42)}
 &=\IRparallel{black}{0}{0}-\IRcupcap{black}{0}{0}.
 \label{eq:ir-zero-affine-right}
\end{align}
\end{Lemma}
\begin{proof}
On each bend, $z=(x_2-x_1)/2$ becomes the right-leg dot. Hence
\[
 \IRcupcapdeco{darkred}{W_V(z)}{M_V(z)}
 =\IRcupcap{black}{0}{0}.
\]
Substitute \eqref{eq:ir-zero-crossing-definition} and use
\eqref{eq:ir-red-affine-left}:
\[
\begin{aligned}
 &\IRcrossdot{black}{0}{0}{(-.34,.42)}
  -\IRcrossdot{black}{0}{0}{(.34,-.42)}
 \\
 &\quad=\IRDparallelcoef{M_V(z)W_V(z)}
  -\IRcupcapdeco{darkred}{W_V(z)}{M_V(z)}
  +\IRDparallelcoef{1-M_V(z)W_V(z)}
 \\
 &\quad=\IRparallel{black}{0}{0}-\IRcupcap{black}{0}{0}.
\end{aligned}
\]
The second identity is the same calculation with
\eqref{eq:ir-red-affine-right}.
\end{proof}

\newcommand{\IRFcross}[6]{\begin{tikzpicture}[baseline=0,scale=.70,line width=.65pt]
 \draw[#1](-.4,-.5)--(.4,.5)(.4,-.5)--(-.4,.5);
 \if\relax\detokenize{#2#3#4#5}\relax\else
  \node[fill=white,inner sep=.3pt]at(0,0){$\scriptstyle\diamond$};\fi
 \if\relax\detokenize{#6}\relax\else\fill[#1](#6)circle(1.65pt);\fi
 \node at(-.4,-.71){$\scriptstyle#2$};\node at(.4,-.71){$\scriptstyle#3$};
 \node at(-.4,.71){$\scriptstyle#4$};\node at(.4,.71){$\scriptstyle#5$};
\end{tikzpicture}}
\newcommand{\IRFparallel}[4]{\begin{tikzpicture}[baseline=0,scale=.70,line width=.65pt]
 \draw[#1](-.4,-.5)--(-.4,.5)(.4,-.5)--(.4,.5);
 \if\relax\detokenize{#4}\relax\else
  \draw[#1,densely dotted](-.4,0)--(.4,0);
  \node[anchor=west]at(.53,0){$\scriptstyle#4$};\fi
 \node at(-.4,-.71){$\scriptstyle#2$};\node at(.4,-.71){$\scriptstyle#3$};
 \node at(-.4,.71){$\scriptstyle#2$};\node at(.4,.71){$\scriptstyle#3$};
\end{tikzpicture}}
\newcommand{\IRFcupcap}[6]{\begin{tikzpicture}[baseline=0,scale=.70,line width=.65pt]
 \draw[#1](-.4,-.65)--(-.4,-.3)
  ..controls(-.4,.0)and(.4,.0)..(.4,-.3)--(.4,-.65);
 \draw[#1](-.4,.65)--(-.4,.3)
  ..controls(-.4,.0)and(.4,.0)..(.4,.3)--(.4,.65);
 \if\relax\detokenize{#6}\relax\else
  \draw[#1,densely dotted](-.4,-.42)..controls(-.8,-.25)and(-.8,.25)..(-.4,.42);
  \node[anchor=east]at(-.83,0){$\scriptstyle#6$};\fi
 \node at(-.4,-.86){$\scriptstyle#2$};\node at(.4,-.86){$\scriptstyle#3$};
 \node at(-.4,.86){$\scriptstyle#4$};\node at(.4,.86){$\scriptstyle#5$};
\end{tikzpicture}}

\newcommand{\IRFmixed}[2]{\begin{tikzpicture}[baseline=0,scale=.64,line width=.65pt]
 \draw[darkred](-.6,-.85)--(-.6,-.35)
  ..controls(-.6,.70)and(.6,.70)..(.6,-.35)--(.6,-.85);
 \draw[darkred,->](0,-.85)--(0,-.35)to[out=90,in=-90]({#1*.9},.8);
 \draw[darkred,densely dotted](0,-.5)--(.6,-.5);
 \fill[darkred](.6,-.78)circle(1.65pt);
 \node[anchor=west]at(.8,-.80){$\scriptstyle W_V(x_3)$};
 \node[anchor=south west,inner sep=1pt]at(1.02,-.27){$\scriptstyle#2$};
 \node at(-.6,-1.06){$\scriptstyle0$};\node at(0,-1.06){$\scriptstyle j$};
 \node at(.6,-1.06){$\scriptstyle0$};\node at({#1*.9},1.02){$\scriptstyle j$};
\end{tikzpicture}}

\subsection{Verification of the affine Brauer relations}
\label{subsec:ir-global-relations}
{Unmarked intermediate colors are summed over. The sums are finite
on finitely generated objects. For nonzero colors we use
\cite[Sections~7.2--7.3]{GRS}, with the oriented normalization fixed
above. For zero colors we use the finite representatives
$A_V=W_V^{-1}$ and $B_V=M_V^{-1}$ as specified in
Lemma~\ref{lem:ir-normalization}.}

\begin{Lemma}\label{lem:ir-global-adjunction}
The reverse cup, cap and dot satisfy
\begin{align}
 \IRsnake{black}{}{1}&=\IRstrand{black}{}
                    =\IRsnake{black}{}{-1},
 \label{eq:ir-global-snake-diagram}\\[3pt]
 \IRdotcup{black}{}{}{(-.3854592,.20128)}
 &=-\IRdotcup{black}{}{}{(.3854592,.20128)},
 &\IRdotcap{black}{}{}{(-.3854592,-.20128)}
 &=-\IRdotcap{black}{}{}{(.3854592,-.20128)}.
 \label{eq:ir-global-dot-diagram}
\end{align}
\end{Lemma}
\begin{proof}
On a colored summand the zero case is
Lemma~\ref{lem:ir-zero-adjunction}. For nonzero colors the defining
marks cancel in both zigzags, as in \cite[Corollary~7.15]{GRS};
the dot changes sign at a bend. Sum these component identities.
\end{proof}

\begin{Lemma}\label{lem:ir-global-affine}
The full crossing satisfies
\begin{align}
 \IRFcross{black}{}{}{}{}{-.28,.35}
 -\IRFcross{black}{}{}{}{}{.28,-.35}
 &=\IRFparallel{black}{}{}{}-\IRFcupcap{black}{}{}{}{}{},
 \label{eq:ir-global-affine-left}\\[2pt]
 \IRFcross{black}{}{}{}{}{-.28,-.35}
 -\IRFcross{black}{}{}{}{}{.28,.35}
 &=\IRFparallel{black}{}{}{}-\IRFcupcap{black}{}{}{}{}{}.
 \label{eq:ir-global-affine-right}
\end{align}
\end{Lemma}
\begin{proof}
The proof of \cite[Lemma~7.14]{GRS} applies to every component except
the one with four zero boundary colors: substitute
Definition~\ref{cross-red123}(1)--(5), use the type $A$ dot relations,
and use adjunction for negative colors. The parallel and cup--cap entries
give the two terms on the right, with invertible denominators wherever
they occur. The four-zero component is Lemma~\ref{lem:ir-zero-affine}.
\end{proof}

\newcommand{\IRFOrigRedBent}{%
\begin{tikzpicture}[x=1bp,y=1bp,scale=.68,
  baseline={([yshift=-.5ex]current bounding box.center)},
  line cap=round,line join=round,
  every node/.style={inner sep=.45pt,outer sep=0pt,font=\fontsize{7.0}{7.6}\selectfont}]
\draw[darkred,line width=.65bp] (2.217,4.319) .. controls (2.217,25.409) and (15.398,25.409) .. (22.428,13.107);
\draw[darkred,line width=.65bp] (22.428,13.107) .. controls (29.458,0.804) and (40.881,0.804) .. (40.881,22.773);
\draw[darkred,line width=.65bp] (13.64,4.319) -- (28.579,22.773);
\node[text=darkred,anchor=north] at (2.217,2.319) {$0$};
\node[text=darkred,anchor=north] at (13.64,2.319) {$0$};
\draw[darkg,line width=.65bp] (47.881,4.319) -- (47.881,21.479);
\node[text=darkg,anchor=north] at (47.881,2.519) {$V$};
\end{tikzpicture}%
}

\newcommand{\IRFOrigRedStraight}{%
\begin{tikzpicture}[x=1bp,y=1bp,scale=.74,
  baseline={([yshift=-.5ex]current bounding box.center)},
  line cap=round,line join=round,
  every node/.style={inner sep=.45pt,outer sep=0pt,font=\fontsize{7.0}{7.6}\selectfont}]
\draw[darkred,line width=.65bp] (2.217,4.319) -- (19.792,27.167);
\draw[darkred,line width=.65bp] (19.792,4.319) -- (2.217,27.167);
\node[text=darkred,anchor=north] at (2.217,2.319) {$0$};
\node[text=darkred,anchor=north] at (19.792,2.319) {$0$};
\draw[darkg,line width=.65bp] (26.792,4.319) -- (26.792,21.479);
\node[text=darkg,anchor=north] at (26.792,2.519) {$V$};
\end{tikzpicture}%
}

\newcommand{\IRFOrigBlackBent}{%
\begin{tikzpicture}[x=1bp,y=1bp,scale=.82,baseline=20bp,
 line cap=round,line join=round,
 every node/.style={inner sep=1pt,outer sep=0pt,font=\scriptsize}]
\draw[black,line width=.65bp] (0,0)--(0,36)
 ..controls (0,54) and (32,54)..(32,36)--(32,8)
 ..controls (32,-10) and (64,-10)..(64,8)--(64,44);
\draw[black,line width=.65bp] (16,0)--(48,44);
\draw[black,densely dotted,line width=.55bp] (32,32)--(39.272727,32);
\draw[black,densely dotted,line width=.55bp] (24.727273,12)--(32,12);
\fill[black] (27,46.125) circle[radius=1.1bp];
\fill[black] (64,16) circle[radius=1.1bp];
\node[anchor=south] at (27,51) {$W_{E_0N}(x)^{-1}$};
\node[anchor=west] at (70,16) {$M_V(x)^{-1}$};
\node[anchor=west] at (70,35) {$M_N(\frac{x_2-x_1}{2})^{-1}$};
\node[anchor=west] at (70,0) {$W_N(\frac{x_2-x_1}{2})^{-1}$};
\node[anchor=north] at (0,-5) {$0$};
\node[anchor=north] at (16,-5) {$0$};
\node[anchor=south] at (48,47) {$0$};
\node[anchor=south] at (64,47) {$0$};
\end{tikzpicture}%
}
\newcommand{\IRFOrigBlackStraight}{%
\begin{tikzpicture}[x=1bp,y=1bp,scale=.86,baseline=17bp,
 line cap=round,line join=round,
 every node/.style={inner sep=1pt,outer sep=0pt,font=\scriptsize}]
\draw[black,line width=.65bp] (0,0)--(0,10)--(24,34)--(24,44);
\draw[black,line width=.65bp] (24,0)--(24,10)--(0,34)--(0,44);
\fill[black] (0,5) circle[radius=1.1bp];
\fill[black] (24,39) circle[radius=1.1bp];
\draw[black,densely dotted,line width=.55bp]
 (0,10)--(-8,10)--(-8,34)--(0,34)
 (24,10)--(32,10)--(32,34)--(24,34);
\node[anchor=east] at (-5,5) {$M_N(x)$};
\node[anchor=west] at (29,39) {$M_V(x)^{-1}$};
\node[anchor=east] at (-12,22) {$M_N(\frac{x_1+x_2}{2})^{-1}$};
\node[anchor=west] at (36,22) {$W_N(-\frac{x_1+x_2}{2})^{-1}$};
\node[anchor=north] at (0,-4) {$0$};
\node[anchor=north] at (24,-4) {$0$};
\node[anchor=south] at (0,48) {$0$};
\node[anchor=south] at (24,48) {$0$};
\end{tikzpicture}%
}

\begin{Lemma}\label{lem:ir-zero-full-mate}\label{lem:ir-zero-mates}
One has
{\begin{equation}\label{eq:ir-zero-mates-diagram}
 \IRcapcross{black}{0}{0}{1}=\IRcapcross{black}{0}{0}{-1},
 \qquad
 \IRcupcrossside{black}{0}{0}{1}=\IRcupcrossside{black}{0}{0}{-1}.
\end{equation}}
\end{Lemma}
\begin{proof}
By adjunction it suffices to compare the bent and straight crossings.
{Put $N=E_0V$. On horizontal connectors below, read $x_1,x_2$
from left to right; on those beside the straight crossing, read them
from bottom to top. A single $x$ is the dot at its marked point.
The calculation in Lemma~\ref{nbone:lemma:nil4}, before identifying
the bent and straight black crossings, gives}
{\begin{equation}\label{eq:ir-zero-mate-defect}
\begin{aligned}
 \IRFOrigRedBent-\IRFOrigRedStraight
 &=\IRFOrigBlackBent\\[8pt]
 &\quad-\IRFOrigBlackStraight.
\end{aligned}
\end{equation}}
This calculation uses only dot slides, adjunction,
\eqref{eq:ir-WM-products}, and
\eqref{eq:ir-W-zero-transport}--\eqref{eq:ir-M-zero-transport}.
The left-hand side is zero by \eqref{eq:ir-red-mates}.
On the right, move the decorations through the cup and cap using
\eqref{eq:ir-W-zero-transport}--\eqref{eq:ir-M-zero-transport}.
Both terms then have the same invertible decorations on their exterior
legs. Canceling them proves the bent-crossing identity; adjunction
gives the two displayed relations. The cup and cap force every
intermediate color here to be zero.
\end{proof}

\begin{Lemma}\label{lem:ir-full-pitchfork}\label{lem:ir-global-mates}
The full crossing satisfies
\begin{equation}\label{eq:ir-global-mates-diagram}
 \IRcapcross{black}{}{}{1}=\IRcapcross{black}{}{}{-1},
 \qquad
 \IRcupcrossside{black}{}{}{1}=\IRcupcrossside{black}{}{}{-1}.
\end{equation}
\end{Lemma}
\begin{proof}
By adjunction it suffices to prove the cap relation. Fix bottom
colors $(i,j,l)$ and top color $k$. If $(l,k)\ne(-i,j)$, substitute
Definition~\ref{cross-red123}(3)--(5) and use the proof of
\cite[Lemma~7.18]{GRS}; dot bends and the triangle identities give
the same straight-strand terms on both sides.

Let $(l,k)=(-i,j)$. For $i,j\ne0$, use
\cite[Lemma~7.16, Corollary~7.17 and Lemma~7.18]{GRS}, with
\eqref{eq:ir-primary-transport}. For $i\ne0,j=0$, the identity is
the mate formula of Definition~\ref{cross-red123}(2).

For $i=0,j>0$, \eqref{eq:ir-W-positive-transport} gives
{\[
\begin{aligned}
 \IRcapcrosscomp{black}{0}{j}{0}{j}{-1}
 &=\IRFmixed{1}{\begin{gathered}
  \scriptstyle\frac{(x_3-x_2)(x_3+x_2+1)}
  {(x_3+x_2)(x_3-x_2+1)^{1-\delta_{j,1}}}\\[-1pt]
  \scriptstyle{}\times\frac{(1+x_3-x_2)^{1-\delta_{j,1}}}{x_2-x_3}
 \end{gathered}}\\[5pt]
 &=\IRFmixed{1}{-\frac{1+x_2+x_3}{x_2+x_3}}
 =\IRFmixed{-1}{-\frac{1+x_2+x_3}{x_2+x_3}}
 =\IRcapcrosscomp{black}{0}{j}{0}{j}{1}.
\end{aligned}
\]}
The middle equality is \eqref{ipitchfork}; the last uses $x_1=-x_3$
at the cap. Cancel the factors with exponent $1-\delta_{j,1}$ before
evaluation, so this includes $j=1$. Taking a mate gives $j<0$.
The case $i=j=0$ is Lemma~\ref{lem:ir-zero-full-mate}.
\end{proof}

For $i>0$, the normalization factors
$\zeta_i\gamma_i(\lambda)$ cancel in the mixed and quadratic relations
and give the normalized bubbles in the two curl relations. Since
$Q^\imath_{i,-i}=-1$, these are the nonexceptional relations used in
\cite[Lemmas~7.21--7.22]{GRS}. Their calculations apply with
$O_V=n_V/m_V$: the factors $u$ and $4u^2-1$ are units at $\pm i$,
and the full bubble includes the color zero.

\newcommand{\IRCcap}[6]{%
\begin{tikzpicture}[baseline=0,scale=.72,line cap=round,line join=round]
 \begin{scope}[yscale=#2]
  \draw[#1,thick](-.42,-.70)--(-.42,0)
    ..controls(-.42,.70)and(.42,.70)..(.42,0)--(.42,-.70);
  \if\relax\detokenize{#5}\relax\else
   \fill[#1](-.42,-.36)circle(1.6pt);
   \node[anchor=east]at(-.54,-.36){$\scriptstyle #5$};\fi
  \if\relax\detokenize{#6}\relax\else
   \fill[#1](.42,-.36)circle(1.6pt);
   \node[anchor=west]at(.54,-.36){$\scriptstyle #6$};\fi
  \node at(-.42,-.91){$\scriptstyle #3$};
  \node at(.42,-.91){$\scriptstyle #4$};
 \end{scope}
\end{tikzpicture}}

\newcommand{\IRCcurl}[8]{%
\begin{tikzpicture}[baseline=0,scale=.72,line cap=round,line join=round]
 \begin{scope}[yscale=#2]
  \draw[#1,thick](-.42,-1)--(-.42,-.50)--(.42,.30)--(.42,.62)
   ..controls(.42,1.28)and(-.42,1.28)..(-.42,.62)
   --(-.42,.30)--(.42,-.50)--(.42,-1);
  \if\relax\detokenize{#5}\relax\else
   \node[anchor=east]at(-.50,.35){$\scriptstyle #5$};
   \node[anchor=west]at(.50,.35){$\scriptstyle-#5$};\fi
  \if\relax\detokenize{#6}\relax\else
   \fill[#1](.42,.62)circle(1.6pt);
   \node[anchor=west]at(.54,.74){$\scriptstyle #6$};\fi
  \if\relax\detokenize{#7}\relax\else
   \fill[#1](-.42,-.72)circle(1.6pt);
   \node[anchor=east]at(-.54,-.72){$\scriptstyle #7$};\fi
  \if\relax\detokenize{#8}\relax\else
   \fill[#1](.42,-.72)circle(1.6pt);
   \node[anchor=west]at(.54,-.72){$\scriptstyle #8$};\fi
  \node at(-.42,-1.21){$\scriptstyle #3$};
  \node at(.42,-1.21){$\scriptstyle #4$};
 \end{scope}
\end{tikzpicture}}
\newcommand{\IRCstrand}[1]{%
\begin{tikzpicture}[baseline=0,scale=.72]
 \draw[thick](0,-.55)--(0,.55);
 \fill(0,0)circle(1.6pt);
 \node[anchor=west]at(.12,0){$\scriptstyle #1$};
 \node at(0,-.77){$\scriptstyle0$};
\end{tikzpicture}}
\newcommand{\IRCmoment}{%
\begin{tikzpicture}[baseline=0,scale=.72]
 \draw[thick](.72,.45)..controls(.12,.45)and(.12,-.45)..(.72,-.45)
   ..controls(1.32,-.45)and(1.32,.45)..(.72,.45);
 \fill(1.17,0)circle(1.6pt);
 \node[anchor=west]at(1.29,0){$\scriptstyle(u)$};
 \draw[thick](-.75,-.60)--(-.75,.60);
 \fill(-.75,0)circle(1.6pt);
 \node[anchor=east]at(-.87,0){$\scriptstyle\dfrac{C_V(u)}{u-x}$};
 \node at(-.75,-.81){$\scriptstyle0$};
\end{tikzpicture}}

\begin{Lemma}\label{lem:ir-zero-full-cap}
The components of crossing invariance with zero-colored boundary strands are
\[
 \sum_{i\in I}\IRCcurl{black}{1}{0}{0}{i}{}{}{}
 =\IRcap{black}{0}{0},
 \qquad
 \sum_{i\in I}\IRCcurl{black}{-1}{0}{0}{i}{}{}{}
 =\IRcup{black}{0}{0}.
\]
All intermediate colors are summed over.
\end{Lemma}
\begin{proof}
Use the calculation of Lemma~\ref{nbone:lemma:R4}, with the actual
polynomial representatives $A_V,B_V$ from Definition~\ref{notation-1}.
Put
\begin{equation}\label{eq:ir-C-definition}
 C_V(u):=A_V(u)B_V(u).
\end{equation}
For any prescribed $K$, choose these representatives so that
\[
 C_V(u)=\frac{(-1)^{c+1}}2\,
 m_{V,0}(u)f_V^+(u)f_V^-(u)(1+u^KR_V(u)),
 \qquad R_V(u)\in Z_V[u],
\]
retaining the vanishing factors at every nonzero color occurring in
the calculation. Thus $C_V$ agrees with $(W_VM_V)^{-1}$ to the
required order at zero and $C_V(0)=1$. In particular,
\[
 O_V(u)C_V(u)-
 \frac{(-1)^{c+1}t_{V,0}(u)}{2u^{\varepsilon_0(V)}}
 \in u^{K-\varepsilon_0(V)}Z_V[u]
\]
for arbitrarily large $K$. Choose $K$ so large that
$x^{K-\varepsilon_0(V)-1}=0$ at each marked zero dot.
The polynomial remainder then contributes zero, and extracting the
full bubble coefficient gives
\begin{equation}\label{eq:ir-curl-moment}
 \left[\IRCmoment\right]_{u:-1}
 =-\IRCstrand{C_V(x)}
  +\IRCstrand{\dfrac{C_V(x)-1}{2x}}.
\end{equation}
The last quotient is formed in $Z_V[u]$ before evaluation.
The bubble marked $(u)$ is $O_V(u)/u-1+1/(2u)$. Substitute this
identity in the dot-sliding calculation of Lemma~\ref{nbone:lemma:R4},
using \eqref{eq:ir-zero-crossing-definition}. The red curl relation
then gives
\begin{equation}\label{eq:ir-zero-cap-defect}
 0=\sum_{i\in I}\IRCcurl{black}{1}{0}{0}{i}{}{}
       {\dfrac{C_V(-x)}{W_V(x)}}
   -\IRCcap{black}{1}{0}{0}{}{\dfrac{C_V(-x)}{W_V(x)}}.
\end{equation}
The reflected calculation gives
\begin{equation}\label{eq:ir-zero-cup-defect}
 0=\IRcupcross{darkred}{0}{0}
  =\sum_{i\in I}\IRCcurl{black}{-1}{0}{0}{i}{}{}
       {\dfrac{C_V(-x)}{M_V(x)}}
   -\IRCcap{black}{-1}{0}{0}{}{\dfrac{C_V(-x)}{M_V(x)}}.
\end{equation}
All intermediate colors occur in these sums, through the full bubble
in \eqref{eq:ir-curl-moment}. The common decorations are units at zero;
canceling them proves both identities.
\end{proof}

\begin{Lemma}\label{lem:ir-global-cap}
The full cup and cap are crossing invariant. Equivalently, for
every $i\in I$,
\begin{equation}\label{eq:ir-zero-curl-diagram}
 \sum_{j\in I}\IRCcurl{black}{1}{i}{-i}{j}{}{}{}
 =\IRcap{black}{i}{-i},
 \qquad
 \sum_{j\in I}\IRCcurl{black}{-1}{i}{-i}{j}{}{}{}
 =\IRcup{black}{i}{-i}.
\end{equation}
\end{Lemma}
\begin{proof}
For nonopposite boundary colors both sides vanish. For $i\ne0$,
apply the nonexceptional calculation of \cite[Lemma~7.21]{GRS}, using
Lemma~\ref{lem:tVi-integral} and $O_V=n_V/m_V$. The normalized red
curl relation \eqref{icurl} makes the difference zero, and its
decorations are units at $\pm i$. The case $i=0$ is
Lemma~\ref{lem:ir-zero-full-cap}.
\end{proof}

\newcommand{\IRTdouble}[8]{%
 \begin{tikzpicture}[baseline={([yshift=-.5ex]current bounding box.center)},
  scale=.66,line cap=round,line join=round]
  \draw[#1,thick](-.4,-1)--(-.4,-.7)--(.4,0)--(-.4,.7)--(-.4,1);
  \draw[#1,thick](.4,-1)--(.4,-.7)--(-.4,0)--(.4,.7)--(.4,1);
  \node at(-.4,-1.20){$\scriptstyle #2$};
  \node at(.4,-1.20){$\scriptstyle #3$};
  \node at(-.4,1.20){$\scriptstyle #4$};
  \node at(.4,1.20){$\scriptstyle #5$};
  \if\relax\detokenize{#6}\relax\else
   \node[anchor=east,fill=white,inner sep=.3pt]at(-.43,0)
    {$\scriptstyle #6$};\fi
  \if\relax\detokenize{#7}\relax\else
   \node[anchor=west,fill=white,inner sep=.3pt]at(.43,0)
    {$\scriptstyle #7$};\fi
  #8
 \end{tikzpicture}}
\newcommand{\IRTnormdouble}{%
 \IRTdouble{darkred}{0}{0}{0}{0}{}{}{%
  \draw[black,densely dotted](-.43,-.85)--(.43,-.85)
   (-.43,.85)--(.43,.85);
  \node[anchor=west]at(.55,-.85){$\scriptstyle1+x_1-x_2$};
  \node[anchor=west]at(.55,.85){$\scriptstyle1-x_1+x_2$};}}

\begin{Lemma}\label{lem:ir-zero-full-square}\label{lem:ir-zero-quadratic-braid}
The zero component of the full quadratic relation is
\begin{equation}\label{eq:ir-zero-square-diagram}
 \sum_{k\in I}\IRTdouble{black}{0}{0}{0}{0}{k}{-k}{}
 =\IRparallel{black}{0}{0}.
\end{equation}
\end{Lemma}
\begin{proof}
{Apply \eqref{eq:nil-square-defect}, obtained in the proof of
Lemma~\ref{nilRel_7} before imposing the black quadratic relation,
with $A_V=W_V^{-1}$ and $B_V=M_V^{-1}$. That identity compares the
normalized red square with the full black square minus the identity.
Thus the intermediate sum includes every opposite-color pair,
and gives}
\begin{equation}\label{eq:ir-zero-square-defect}
 \IRTnormdouble
 =\sum_{k\in I}\IRTdouble{black}{0}{0}{0}{0}{k}{-k}{}
  -\IRparallel{black}{0}{0}.
\end{equation}
The red double crossing is zero, proving the assertion.
\end{proof}

\begin{Lemma}\label{lem:ir-global-square}
The full quadratic relation holds:
\begin{equation}\label{eq:ir-global-square-diagram}
 \IRTdouble{black}{i}{j}{k}{l}{}{}{}
 =\delta_{k,i}\delta_{l,j}\,\IRparallel{black}{i}{j}.
\end{equation}
The unmarked intermediate colors on the left are summed.
\end{Lemma}
\begin{proof}
As in \cite[Lemma~7.22]{GRS}, the dot relations and
Lemma~\ref{lem:ir-global-cap} move a dot through the double crossing.
Thus a component vanishes unless $(k,l)=(i,j)$, by nilpotence of the
shifted input dots. If $i\ne-j$, the ordinary calculation in that
lemma applies, including a single zero color. If $j=-i\ne0$, use
its nonexceptional calculation with the normalized red relation
\eqref{iquadratic}; the required decorations are units at $\pm i$.
The remaining case is Lemma~\ref{lem:ir-zero-full-square}.
\end{proof}

\newcommand{\IRTMbraid}[6]{%
 \begin{tikzpicture}[baseline={([yshift=-.5ex]current bounding box.center)},
  scale=.75,line cap=round,line join=round]
  \begin{scope}[xscale=#1]
   \draw[thick](.45,.60)--(-.45,-.60);
   \draw[thick](.45,-.60)--(-.45,.60);
   \draw[thick](0,-.60)..controls(0,-.30667)and(-.45,-.29334)..(-.45,0)
    ..controls(-.45,.29334)and(0,.30667)..(0,.60);
  \end{scope}
  \foreach \IRTMpos in {-.45,0,.45}{%
   \draw[thick](\IRTMpos,-.60)--(\IRTMpos,-.96)
    (\IRTMpos,.60)--(\IRTMpos,1.35);}
  \if\relax\detokenize{#2}\relax\else
   \foreach \IRTMcol [count=\IRTMn from 0] in {#2}{%
    \node at({.45*(\IRTMn-1)},-1.14){$\scriptstyle\IRTMcol$};}\fi
  \if\relax\detokenize{#3}\relax\else
   \foreach \IRTMcol [count=\IRTMn from 0] in {#3}{%
    \node at({.45*(\IRTMn-1)},1.55){$\scriptstyle\IRTMcol$};}\fi
  \if\relax\detokenize{#4}\relax\else
   \draw[densely dotted](-.47,-.80)--(.47,-.80);
   \node[anchor=west]at(.65,-.80){$\scriptstyle #4$};\fi
  \ifnum#5=1
   \draw[densely dotted](-.45,.76)--(0,.76) (0,1.16)--(.45,1.16);
   \node[anchor=west]at(.65,.76){$\scriptstyle1+x_2-x_1$};
   \node[anchor=west]at(.65,1.16){$\scriptstyle1+x_3-x_2$};
  \fi
  #6
 \end{tikzpicture}}

\begin{Lemma}\label{lem:ir-global-braid}
The reverse crossings satisfy
\begin{equation}\label{eq:ir-global-braid-diagram}
 \IRTMbraid{1}{}{}{}{0}{}=\IRTMbraid{-1}{}{}{}{0}{}.
\end{equation}
\end{Lemma}
\begin{proof}
The dot, cup/cap and quadratic relations move each dot through the
difference of the braid diagrams, reversing its position. The argument
of \cite[Lemma~7.23]{GRS} therefore reduces to input $(i,j,k)$ and
output $(k,j,i)$.

If all three colors are nonzero, apply
\cite[Lemma~7.24, Corollary~7.25, and Lemmas~7.26,~7.28]{GRS}.
The required bubble and mate identities have been proved above;
the denominators at a repeated opposite pair have nonzero residues
$i$ and $4i^2-1$. If exactly one color is zero, adjunction and
\eqref{eq:ir-global-mates-diagram} reduce to $(i,j,0)$ or $(i,0,k)$
with the other colors positive. The type $A$ calculation of
\cite[Lemma~7.24]{GRS} then applies, also at the pair $0,1$.

If two colors are zero, use adjunction to reduce to $(0,0,k)$ with
$k>0$. Repeat the calculation of \eqref{fw:mixed-braid}, keeping
its two leading black diagrams separate. Equations
\eqref{mixed-red1}--\eqref{mixed-b2} and Lemma~\ref{compu-T},
with the notation of Definition~\ref{notation-1}, cancel the other terms,
{including those through $(k,-k,k)$. Reading $x_1,x_2,x_3$ from
left to right at the lower boundary, the red relation gives}
\begin{equation}\label{eq:ir-mixed-braid-defect}
 \begin{aligned}
 &\IRTMbraid{1}{0,0,k}{k,0,0}
  {\dfrac{(x_1-x_3)(x_2-x_3)}
  {(1+x_1-x_2)(1-x_1+x_3)(1-x_2+x_3)}}{0}{}
 \\
 &\quad-
 \IRTMbraid{-1}{0,0,k}{k,0,0}
  {\dfrac{(x_1-x_3)(x_2-x_3)}
  {(1+x_1-x_2)(1-x_1+x_3)(1-x_2+x_3)}}{0}{}=0.
 \end{aligned}
\end{equation}
On the $(0,0,k)$ input, the common factor has residue
$k^2/(k+1)^2\ne0$ for every $k>0$, and is therefore invertible.

{For three zero colors, use Section~\ref{sec:nil-brauer-braid} with
the lower crossing formula at the top and the upper formula at the
other two crossings. Choose supported finite polynomial
representatives as in Definition~\ref{notation-1}. The projection
and pairing argument at the beginning of
Appendix~\ref{app:graphical-certificates} gives the source terms in
Appendix~\ref{27terms-app}, with only the external colors fixed.
Retain the polynomials until coefficient extraction. In the
subsequent reductions, keep the two leading black braid diagrams
separate; the remaining comparisons use the dot, cup/cap and
quadratic relations already proved. The red braid then gives
the following identity, with $x_1,x_2,x_3$ read from left to right
at the upper boundary:
}\begin{equation}\label{eq:ir-allzero-braid-defect}
 \IRTMbraid{1}{0,0,0}{0,0,0}{}{1}{}
 -\IRTMbraid{-1}{0,0,0}{0,0,0}{}{1}{}=0.
\end{equation}
The two exterior factors are units at the zero-colored dots.
Canceling them proves this case.
\end{proof}

\subsection{The reconstructed action}
\label{subsec:ir-action}

\begin{Theorem}\label{prop:ir-action}
The reverse assignments define a $\mathbb C$-linear monoidal functor
\[
 \mathcal{AB}\longrightarrow\End_{\mathbb C}(\mathcal M).
\]
Its generating functor is $E=\bigoplus_{i\in\mathbb Z}E_i$, with
dot spectrum contained in $\mathbb Z$. On every nonzero finitely
generated object $V\in\mathcal M_\lambda$, the dot has minimal polynomial
\[
 m_V(u)=\prod_{i\in\mathbb Z}(u-i)^{\varepsilon_i(V)}.
\]
\end{Theorem}
\begin{proof}
Lemmas~\ref{lem:ir-global-adjunction}--\ref{lem:ir-global-affine},
\ref{lem:ir-full-pitchfork}, \ref{lem:ir-global-cap},
\ref{lem:ir-global-square}, and \ref{lem:ir-global-braid}
verify all defining relations of Definition~\ref{C defn}.
{The assignments are natural by Lemma~\ref{lem:ir-normalization}
and the argument following Definition~\ref{natrual-basic-red}.
The extension by filtered colimits described at the start of this
section gives the stated functor.}
The dot is $y+i$ on $E_i$ for $i>0$, $-y-i$ on $E_{-i}$,
and $y$ on $E_0$. Nilpotence and finite color support give the
displayed polynomial and identify the generalized $i$-eigenspace functor with $E_i$.
\end{proof}

\begin{Cor}\label{cor:ir-relative}
Let $\mathcal R$ be an ambient locally Schurian
$2$-representation. For each weight $\lambda$, let
$\mathcal N_\lambda\subseteq\mathcal R_\lambda$ be a full additive
subcategory of finitely generated ambient objects, closed under
isomorphisms. Set $\mathcal N=\bigoplus_\lambda\mathcal N_\lambda$ and
suppose that every $E_i$ preserves $\mathcal N$.
Suppose every object of $\mathcal N$ has finite
color support and nilpotent red dots. Then the reverse assignments
define an affine Brauer action on
$\bigoplus_\lambda\mathcal N_\lambda$.
\end{Cor}
\begin{proof}
Apply the preceding proof to the ambient finitely generated objects
of $\mathcal N$. Stability and finite support ensure that all functors
and diagrams lie in $\mathcal N$. For a morphism, the common formulas
on the direct sum of its source and target give naturality, as in
Lemma~\ref{lem:ir-normalization}; fullness then restricts every
transformation. Exactness of the $E_i$ makes finite color support and
dot nilpotence descend to ambient simple subquotients, so the preceding
simple-object arguments apply.
\end{proof}

\section{The cyclotomic Brauer--\texorpdfstring{$\imath$}{i}-Kac--Moody isomorphism}
\label{sec:integral-cyclotomic}

{Following \cite[Section~8]{GRS}, we apply
Theorem~\ref{thm:integral-forward-action} and
Corollary~\ref{cor:ir-relative} to the cyclotomic quotients.
We work over $\mathbb C$, with $I=\mathbb Z$,
$I^\imath=\mathbb Z_{>0}$ and the weight convention of
Definition~\ref{def-wei}. Here $\mathfrak U^\imath$ is the integral
$\imath$-Kac--Moody $2$-category.} Its weight lattice is
\[
 X^\imath=X/(1+\tau)X
 =\bigoplus_{i>0}\mathbb Z\bar\Lambda_i\oplus
   (\mathbb Z/2\mathbb Z)\bar\Lambda_0.
\]
Thus $\bar\alpha_0=0$, $\bar\alpha_{-i}=-\bar\alpha_i$, and
$E_{-i}=F_i$ for $i>0$, whereas $E_0$ is self-adjoint.

{The cyclotomic isomorphism constructed below transfers the grading
of the $2$-representation quotient to the cyclotomic Brauer category.}

\subsection{Integral cyclotomic parameters}
\label{subsec:ic-parameters}

Fix a positive integer $a$ and $\mathbf u=(u_1,\ldots,u_a)\in\mathbb Z^a$. Set
\begin{equation}\label{eq:ic-m}
 m(u)=\prod_{r=1}^a(u-u_r)
     =\prod_{i\in\mathbb Z}(u-i)^{\epsilon_i},
 \qquad \epsilon_i=\#\{r\mid u_r=i\}.
\end{equation}
Fix the $\mathbf u$-admissible bubble parameters $\omega_r\in\mathbb C$
by the expansion at infinity
\begin{equation}\label{eq:ic-admissibility}
 O_{\mathbf u}(u):=u-\frac12+\sum_{r\ge0}\omega_ru^{-r}
 =\left(u-\frac12(-1)^a\right)
   \prod_{r=1}^a\frac{u+u_r}{u-u_r}.
\end{equation}

\begin{Defn}\label{def:ic-brauer}
The cyclotomic Brauer category $\mathcal{CB}_{\mathbf u}$ is the
quotient of $\mathcal{AB}$ by the left tensor ideal generated by
\begin{equation}\label{eq:ic-brauer-relations}
 m(x):E\longrightarrow E,
 \qquad \bigcirc_r-\omega_r1_{\mathbf1}\quad(r\ge0),
\end{equation}
where $\bigcirc_r$ is the black bubble with $r$ dots.
Its locally unital path algebra is
\begin{equation}\label{eq:ic-B}
 B_{\mathbf u}
 =\bigoplus_{d,e\ge0}
   \operatorname{Hom}_{\mathcal{CB}_{\mathbf u}}(E^d,E^e),
\end{equation}
with multiplication given by composition whenever the objects match,
and zero otherwise. Write $1_d=1_{E^d}$.
\end{Defn}

The admissible basis theorem \cite[Corollary~4.27 and (4.45)]{RS-cyc},
transported to the left tensor convention of Definition~\ref{C defn},
gives local finite-dimensionality and
\[
 \operatorname{End}_{\mathcal{CB}_{\mathbf u}}(\mathbf1)=\mathbb C,
 \qquad
 \operatorname{End}_{\mathcal{CB}_{\mathbf u}}({E})
 =\mathbb C[x]/(m(x)).
\]
By \cite[Proposition~8.3]{GRS}, $\operatorname{lfdmod}\text{-}B_{\mathbf u}$
carries the induced affine Brauer action.

\begin{Lemma}\label{lem:ic-exact-spectrum}
The dot spectrum of this action is exactly $\mathbb Z$.
\end{Lemma}
\begin{proof}
First take the $u_r$ to be independent parameters. At rank $d$,
the generic seminormal formulas \cite[Theorems~4.13 and~5.3]{AMR-cyc}
give, for every dot $x$, the identity
\[
 \prod_{r=1}^a\prod_{c=-d}^{d}
       (x-u_r-c)(x+u_r+c)=0.
\]
{The admissible basis \cite[Theorem~5.5]{AMR-cyc} makes the algebra
free over the parameter ring, so the identity descends from its
fraction field and survives every integral specialization, including
coincident parameters. Hence all dot eigenvalues are integers.} Conversely, fix a
root $u_r$ of $m$. The quotient of the
rank-$d$ cyclotomic Brauer algebra by its cup--cap generators is the
cyclotomic degenerate Hecke algebra \cite[Proposition~7.2]{AMR-cyc}.
For $\eta=\pm1$, it has a one-dimensional module on which the crossings
act by $\eta$ and the successive dots, numbered as strands are added
on the left, act by
\[
 u_r,\ u_r+\eta,\ \ldots,\ u_r+\eta(d-1).
\]
{Hence $u_r+\eta(d-1)$ is an eigenvalue of the dot in the regular
rank-$d$ algebra, and therefore on
$E(1_{d-1}B_{\mathbf u})\cong1_dB_{\mathbf u}$.
Varying $d\geq1$ and $\eta$ gives every integer; apply
Lemma~\ref{endoe-dec}.}
\end{proof}

Theorem~\ref{thm:integral-forward-action} therefore applies to
$\operatorname{lfdmod}\text{-}B_{\mathbf u}$. In particular,
$E=\bigoplus_iE_i$, and only finitely many colors occur on each
finitely generated projective.

Put
\begin{equation}\label{eq:ic-n}
 s=\frac12(-1)^a,\qquad
 n(u)=O_{\mathbf u}(u)m(u)
     =(u-s)\prod_{i\in\mathbb Z}(u-i)^{\phi_i},
 \qquad \phi_i=\epsilon_{-i}.
\end{equation}
Since $s\notin\mathbb Z$, we have
\begin{equation}\label{eq:ic-zero-exponent}
 \phi_0=\epsilon_0.
\end{equation}

Define the initial weight and its zero-color parameter by
\begin{equation}\label{eq:ic-kappa}
 \kappa=\sum_{i>0}(\epsilon_i-\epsilon_{-i})\bar\Lambda_i
              +c_0\bar\Lambda_0,
 \qquad c_0\in\{0,1\},\quad c_0\equiv\epsilon_0\pmod2.
\end{equation}
Then
\begin{equation}\label{eq:ic-weight-parameters}
 h_i:=-\langle\alpha_i^{\theta,\vee},\kappa\rangle
       =\phi_i-\epsilon_i\quad(i>0),\qquad
 c(\kappa)=c_0,
 \qquad a_\kappa\equiv a\pmod2.
\end{equation}
The last equality follows from
$a\equiv\epsilon_0+\sum_{i>0}(\epsilon_{-i}-\epsilon_i)\pmod2$.

Let $P=1_0B_{\mathbf u}$. Its bubble series is $O_{\mathbf u}$ and
$m_P=m$. Every simple subquotient $L$ of $P$ has the same bubble
series. Comparing its orders at $\pm i$ and its value
$O_{\mathbf u}(0)=(-1)^{\epsilon_0+1}/2$ in Lemma~\ref{barn-dec}
gives $\operatorname{wt}_L^\imath=\kappa$. Thus
Lemma~\ref{weight-pro}(4) places $P$ in the central block
$(\operatorname{lfdmod}\text{-}B_{\mathbf u})_\kappa$.
Moreover, $E^dP\cong1_dB_{\mathbf u}$.
For
$\mathbf i=(i_1,\ldots,i_d)\in I^d$, put
\[
 E_{\mathbf i}=E_{i_d}\cdots E_{i_1},\qquad
 \operatorname{wt}(\mathbf i)=\sum_{r=1}^d\bar\alpha_{i_r},\qquad
 E_{\emptyset}=1.
\]
The joint generalized-eigenspace decomposition gives orthogonal
idempotents $e_{\mathbf i}\in1_dB_{\mathbf u}1_d$ such that
\begin{equation}\label{eq:ic-B-idempotents}
 1_d=\sum_{\mathbf i\in I^d}e_{\mathbf i},\qquad
 e_{\mathbf i}B_{\mathbf u}\cong E_{\mathbf i}P,
 \qquad E_{\mathbf i}P\in
 (\operatorname{lfdmod}\text{-}B_{\mathbf u})_
 {\kappa-\operatorname{wt}(\mathbf i)}.
\end{equation}
For fixed $d$ the sum is finite after zero terms are removed.

\subsection{The integral cyclotomic 2-representation}
\label{subsec:ic-quotient}

For $i>0$, write $C_{\kappa,i}(u)$ and $D_{\kappa,i}(u)$ for the
counterclockwise and clockwise series
$\circlearrowleft_{\kappa,i}(u)$ and $\circlearrowright_{\kappa,i}(u)$
of Definition~\ref{evalu-bub}, including fake bubbles. The prescribed scalar
series at the initial object are
\begin{equation}\label{eq:ic-local-bubbles}
 C_i(u)=u^{h_i},\qquad D_i(u)=u^{-h_i}.
\end{equation}
Equivalently, writing these series as
$\sum_{r\in\mathbb Z}\omega_i^{(r)}u^{-r-1}$ and
$\sum_{r\in\mathbb Z}\widetilde\omega_i^{(r)}u^{-r-1}$, one has
\begin{equation}\label{eq:ic-local-coefficients}
 \omega_i^{(r)}=\delta_{r,-h_i-1},\qquad
 \widetilde\omega_i^{(r)}=\delta_{r,h_i-1}.
\end{equation}
At color zero, let $\beta_{\kappa,r}$ denote the red bubble with
$r$ dots. The prescribed values are
\begin{equation}\label{eq:ic-zero-bubbles}
 \beta_{\kappa,r}=c_0\delta_{r,0}\,1_\kappa\quad(r\ge0).
\end{equation}
Thus the zero bubble series and its normalized version are
\begin{equation}\label{eq:ic-zero-series}
 b_{\kappa,0}(u)=c_0u^{-1},\qquad
 O_{\kappa,0}(u)=(-1)^{c_0}(1-2u b_{\kappa,0}(u))=1.
\end{equation}
These replace the bubble parameters of \cite[Definition~8.4]{GRS}.

Consider the principal $2$-representation
\[
 \mathcal R(\kappa)_\lambda
 =\operatorname{Hom}_{\mathfrak U^\imath}(\kappa,\lambda),
 \qquad \lambda\in X^\imath,
\]
with action by left horizontal composition. We use invariant ideals as in \cite[Section~8.2]{GRS}.

\begin{Defn}\label{def:ic-quotient}
Let $\mathcal I(m)$ be the smallest $\mathfrak U^\imath$-invariant
ideal of $\mathcal R(\kappa)$ containing the following morphisms:
\begin{itemize}
\item For each $i>0$, the dotted-strand morphisms
      \begin{equation}\label{eq:ic-nonzero-dots}
      y_i^{\epsilon_i}1_{E_i1_\kappa},\qquad
      y_{-i}^{\phi_i}1_{E_{-i}1_\kappa}.
      \end{equation}
\item The zero-color dotted-strand morphism
      \begin{equation}\label{eq:ic-zero-dot}
      y_0^{\epsilon_0}1_{E_01_\kappa}.
      \end{equation}
\item For each $i>0$, the coefficients of
      \begin{equation}\label{eq:ic-bubble-generators}
      C_{\kappa,i}(u)-u^{h_i}1_\kappa,\qquad
      D_{\kappa,i}(u)-u^{-h_i}1_\kappa.
      \end{equation}
\item For each $r\ge0$, the zero-color bubble morphism
      \begin{equation}\label{eq:ic-zero-bubble-generators}
      \beta_{\kappa,r}-c_0\delta_{r,0}1_\kappa.
      \end{equation}
\end{itemize}
Here $y_i$ denotes the original red dot. An exponent zero means
the identity morphism of the corresponding object. Define
\begin{equation}\label{eq:ic-Q-category}
 \mathcal Q(m)_\lambda
   =\mathcal R(\kappa)_\lambda/\mathcal I(m)_\lambda,
 \qquad
 \mathcal Q(m)=\bigoplus_{\lambda\in X^\imath}\mathcal Q(m)_\lambda.
\end{equation}
Objects have only finitely many nonzero weight components, and
morphisms between different weight components are zero.
\end{Defn}

The relations are imposed at $1_\kappa$ and extended by left horizontal
composition. In particular, $\epsilon_0=0$ means $E_01_\kappa=0$.
They vanish on $P$: the dot bounds follow from $m_P=m$, while
Definition~\ref{natrual-basic1}(4) gives the nonzero bubbles from
$t_{P,i}(u)=(u-i)^{\epsilon_{-i}}$. At zero, the bubble calculations
in the proofs of Lemmas~\ref{nbone:lemma:R4} and~\ref{nbone:lemma:nil2} give
\[
 \frac1{2u}+\frac{(-1)^{c_0+1}}2
       \frac{u^{\epsilon_0}}{u^{\epsilon_0+1}}
 =\frac{c_0}{u}.
\]

\begin{rem}\label{rem:ic-full-data}
The polynomial $m$, and not only $\kappa$, is part of the quotient data:
$\kappa$ records $\epsilon_i-\epsilon_{-i}$ and
$\epsilon_0\pmod2$, whereas the relations use the individual exponents.
\end{rem}

{The full subcategory on the word objects $E_{\mathbf i}1_\kappa$ has path algebra}
\begin{equation}\label{eq:ic-Q-algebra}
 Q_m=
 \bigoplus_{\substack{\mathbf i,\mathbf j\\
       \operatorname{wt}(\mathbf i)=\operatorname{wt}(\mathbf j)}}
 \operatorname{Hom}_{\mathcal Q(m)}
     (E_{\mathbf j}1_\kappa,E_{\mathbf i}1_\kappa).
\end{equation}
Write $1_{\mathbf i}=1_{E_{\mathbf i}1_\kappa}$ and
$P'=1_\emptyset Q_m$. After additive idempotent completion, Yoneda
identifies these objects with the corresponding projective modules;
in particular,
\begin{equation}\label{eq:ic-Q-projectives}
 E_{\mathbf i}P'\cong1_{\mathbf i}Q_m.
\end{equation}

\begin{Lemma}\label{lem:ic-support-grading}
For each sequence $\mathbf i$, only finitely many colors $j$ satisfy
$E_jE_{\mathbf i}1_\kappa\ne0$. For each $d\geq0$, only finitely many
$\mathbf i\in I^d$ satisfy $E_{\mathbf i}1_\kappa\ne0$.
The ideal $\mathcal I(m)$ and the algebra $Q_m$ are graded.
\end{Lemma}
\begin{proof}
Let $W=E_{i_d}\cdots E_{i_1}1_\kappa$. If $\epsilon_j=0$ and
$j\notin\{i_r-1,i_r,i_r+1,-i_r\}$ for every $r$, the red quadratic
relations make each crossing of $j$ with $i_r$ invertible. Hence
\[
 E_jW\cong E_{i_d}\cdots E_{i_1}E_j1_\kappa=0.
\]
There are only finitely many exceptions. Induction on sequence length
proves the second assertion; both assertions extend to finite sums
and summands.

{In the grading of Section~2, dots have degree $2$, and}
\[
\begin{aligned}
 \deg\bigl([C_{\kappa,i}(u)]_{u:-r-1}\bigr)&=2(r+1+h_i),\\
 \deg\bigl([D_{\kappa,i}(u)]_{u:-r-1}\bigr)&=2(r+1-h_i),
 &\deg\beta_{\kappa,r}&=2r.
\end{aligned}
\]
Every nonzero prescribed scalar occurs in degree zero. Thus all
generators of $\mathcal I(m)$ are homogeneous, including the
identity relations with exponent zero. Horizontal composition
preserves the grading.
\end{proof}

\newcommand{\ICLFcross}[4]{%
\begin{tikzpicture}[iQ,centerzero,scale=.9,darkred]
\draw[-] (-.3,-.3)\botlabel{#1} -- (.3,.3);
\draw[-] (.3,-.3)\botlabel{#2} -- (-.3,.3);
\if\relax\detokenize{#3}\relax\else
\draw[densely dotted] (-.21,-.21)--(.21,-.21);
\node[anchor=north,font=\scriptsize,inner sep=1pt]at(0,-.78){$#3$};
\fi
\if\relax\detokenize{#4}\relax\else
\draw[densely dotted] (-.21,.21)--(.21,.21);
\node[anchor=south,font=\scriptsize,inner sep=1pt]at(0,.35){$#4$};
\fi
\end{tikzpicture}%
}
\newcommand{\ICLFparallel}[3]{%
\begin{tikzpicture}[iQ,centerzero,scale=.9,darkred]
\draw[-] (-.2,-.3)\botlabel{#1} -- (-.2,.3);
\draw[-] (.2,-.3)\botlabel{#2} -- (.2,.3);
\if\relax\detokenize{#3}\relax\else
\draw[densely dotted] (-.2,0)--(.2,0);
\node[anchor=south,font=\scriptsize,inner sep=1pt]at(0,.35){$#3$};
\fi
\end{tikzpicture}%
}
\newcommand{\ICLFcap}[3]{%
\begin{tikzpicture}[iQ,centerzero,scale=.9,darkred]
\draw[-] (-.3,-.25)\botlabel{#1}
  to[out=90,in=90,looseness=3](.3,-.25)\botlabel{#2};
\draw[densely dotted] (-.3,-.25)--(.3,-.25);
\node[anchor=north,font=\scriptsize,inner sep=1pt]at(0,-.72){$#3$};
\end{tikzpicture}%
}
\newcommand{\ICLFcupcap}[1]{%
\begin{tikzpicture}[iQ,centerzero,scale=.9,darkred]
\draw[-] (-.2,-.3)\botlabel{0}
  to[looseness=2,out=90,in=90](.2,-.3)\botlabel{0};
\draw[-] (.2,.3)
  to[looseness=2,out=-90,in=-90](-.2,.3);
\draw[densely dotted] (-.2,.3)--(.2,.3);
\node[anchor=south,font=\scriptsize,inner sep=1pt]at(0,.4){$#1$};
\end{tikzpicture}%
}

\begin{Lemma}\label{lem:ic-local-finiteness}
Every dot on $E_{\mathbf i}1_\kappa$ is nilpotent, and
$\operatorname{Hom}_{\mathcal Q(m)}(E_{\mathbf i}1_\kappa,E_{\mathbf j}1_\kappa)$
is finite-dimensional.
\end{Lemma}
\begin{proof}
{We prove nilpotence by induction on word length, starting with
\eqref{eq:ic-nonzero-dots}--\eqref{eq:ic-zero-dot}.
On $E_jE_iU$, only the leftmost dot remains to be considered.
In the diagrams below, $x_1,x_2$ are the left and right dots,
and all diagrams act on $U$. By induction, choose $N,L\geq1$
annihilating the right dots on $E_jE_iU$ and $E_iE_jU$, respectively.}

Suppose first that $i\ne j$. If $j=-i$, the cap dot bend gives
\begin{equation}\label{eq:ic-cap-dot-bound}
 \ICLFcap{j}{i}{x_1^N}
 =(-1)^N\ICLFcap{j}{i}{x_2^N}=0.
\end{equation}
Thus right multiplication by $x_1^N$ kills all cup--cap terms in
\eqref{idotslide} and \eqref{iquadratic}; such terms are absent for
$j\ne-i$. Iterating \eqref{idotslide} therefore gives
\[
 \ICLFcross{j}{i}{x_1^{N+L}}{}
 =\ICLFcross{j}{i}{x_1^N}{x_2^L}=0.
\]
Composing with the opposite crossing and applying \eqref{iquadratic} gives
\[
 \ICLFparallel{j}{i}
 {Q^\imath_{j,i}(x_1,x_2)\,x_1^{N+L}}=0.
\]
The coefficient in \eqref{Qij}--\eqref{iQ} is a nonzero constant or
a nonzero multiple of $x_1-x_2$.
In the latter case, the displayed equality gives
$x_1^{N+L+1}=x_2x_1^{N+L}$; iterating $N$ times and using $x_2^N=0$
gives $x_1^{2N+L}=0$. The constant case is immediate.

It remains to consider $i=j$. The nilHecke relations, and at
$i=0$ the relations \eqref{eq:ir-red-affine-left},
\eqref{eq:ir-red-affine-right},
\eqref{eq:ir-red-curl-bubble} and
\eqref{eq:ir-red-quadratic-braid}, imply
\begin{equation}\label{eq:ic-dot-interchanger-square}
 \left(
   \ICLFcross{i}{i}{}{x_1-x_2}
   -\ICLFparallel{i}{i}{}
 \right)^2
 =\ICLFparallel{i}{i}{}.
\end{equation}
Products here mean vertical composition. For $i\ne0$, the
parenthesized difference interchanges the two dots, and hence
$x_1^N=0$. At $i=0$, the corresponding calculation is
\[
\begin{aligned}
 &\left(
   \ICLFcross{0}{0}{}{x_1-x_2}-\ICLFparallel{0}{0}{}
  \right)\ICLFparallel{0}{0}{x_1}\\
 &\quad-\ICLFparallel{0}{0}{x_2}\left(
   \ICLFcross{0}{0}{}{x_1-x_2}-\ICLFparallel{0}{0}{}
  \right)
 =-\ICLFcupcap{x_1-x_2}.
\end{aligned}
\]
The cap bend \eqref{eq:ic-cap-dot-bound} also applies to two zero strands.
Multiplying the displayed commutator by $x_1^N$ therefore kills its
cup--cap term. Repeating the resulting interchange $N$ times moves the
remaining $N$ dots to the right strand, where $x_2^N=0$; the square
identity \eqref{eq:ic-dot-interchanger-square} then gives $x_1^{2N}=0$.

{As in \cite[Proposition~8.5]{GRS}, use the spanning diagrams of
\cite[Definition~5.8 and Lemma~5.9]{BWW25}, which include color zero.
For fixed boundaries there are finitely many reduced matching diagrams.
At each marked point, a horizontal slice gives a word object, so the
nilpotence just proved bounds its dot powers. The bubbles in the
rightmost $\kappa$-region are scalar by Definition~\ref{def:ic-quotient}.
Thus finitely many dotted diagrams span the morphism space.}
\end{proof}

{As in \cite[Section~8.2]{GRS}, biadjunction induces a locally
Schurian $2$-representation on $\operatorname{lfdmod}\text{-}Q_m$
and preserves $\operatorname{pmod}\text{-}Q_m$.
Lemmas~\ref{lem:ic-support-grading}--\ref{lem:ic-local-finiteness}
make the latter a nilpotent $2$-subrepresentation: finite color support
and dot nilpotence pass to finite sums and summands of the word
projectives. Corollary~\ref{cor:ir-relative} therefore gives the
reverse affine Brauer action.}

The scalar series in the quotient have the required normalization:
\begin{equation}\label{eq:ic-series-match}
 \widetilde O_\kappa(u)
 =\prod_{i>0}(-1)^{h_i}(u-i)^{h_i}(-u-i)^{-h_i}
 =\prod_{r=1}^a\frac{u+u_r}{u-u_r}.
\end{equation}
The zero roots contribute $1$ to this product, and their parity is
recorded in $c_0$. Hence
\begin{equation}\label{eq:ic-black-series-match}
 \left(u-\frac12(-1)^{a_\kappa}\right)
   \widetilde O_\kappa(u)=O_{\mathbf u}(u).
\end{equation}
The dot on $EP'$ is annihilated by $m$, and
\eqref{eq:ic-black-series-match} gives the prescribed black bubbles.
Thus evaluation at $P'$ factors through $\mathcal{CB}_{\mathbf u}$.

\subsection{The isomorphism}
\label{subsec:ic-isomorphism}

\begin{Lemma}\label{lem:ic-inverse-generators}
On the projectives generated by $P$, the forward assignment followed
by the reverse assignment fixes the dots, cups, caps, crossings, and
bubbles. The reverse assignment followed by the forward assignment does
the same on the projectives generated by $P'$.
\end{Lemma}
\begin{proof}
{We compare the formulas on the word projectives.
The dot substitutions are inverse by Definitions~\ref{natrual-basic}(1),
\ref{eq:red-cup}(4), and \ref{natrual-basic-red}(2).
For nonzero colors, use \cite[proof of Theorem~8.6]{GRS}.
Lemmas~\ref{lem:tVi-integral} and~\ref{lem:red-black-series}
recover the primary factors. Replacing $m_V$ by a multiple and
multiplying $n_V$ by the same factor leaves all the ratios defining
the coefficients unchanged, by the argument following
Definition~\ref{eq:red-cup} and Lemma~\ref{lem:ir-normalization}.
Thus these comparisons commute with evaluation using common
annihilating polynomials and sufficiently high polynomial representatives.

}Fix such a projective $V$ of weight $\lambda$. At zero,
Lemmas~\ref{nbone:lemma:R4} and~\ref{nbone:lemma:nil2},
with $r$ dots inserted, give
\begin{equation}\label{eq:ic-zero-factor-recovery}
 \left.\IRzerobubbledeco{x^r}\right|_V
 =\frac{\delta_{r,0}}2\,1_V
  +\frac{(-1)^{c(\lambda)+1}}2
    \left[\frac{t_{V,0}(u)}{u^{\varepsilon_0(V)}}\right]_{u:-r}
 \qquad(r\geq0).
\end{equation}
{Use $C_V=A_VB_V$ from \eqref{eq:ir-C-definition}, retaining its
nonzero-color vanishing factors and choosing $K>\varepsilon_0(V)$
sufficiently large. Then
$O_VC_V-(-1)^{c(\lambda)+1}t_{V,0}/(2u^{\varepsilon_0(V)})\in uZ_V[u]$.
In
$[u^rC_V(u)(O_V(u)/u-1+1/(2u))]_{u:-1}$,
this polynomial remainder contributes zero, while $C_V(0)=1$
gives $\delta_{r,0}/2$. This proves the formula and recovers
$O_{V,0}(u)=t_{V,0}(u)/u^{\varepsilon_0(V)}$.
Hence both assignments use the same bubble series and weight.}

The nonzero cups and caps have reciprocal decorations.
At zero, \eqref{eq:ir-cup-unit-slide} puts the inverse cup
decoration on its original left leg. The two cancellations are
\[
 \IRDcupdeco{-1}{W_{E_0V}(x)^{-1}W_{E_0V}(x)}
 =\IRcup{darkred}{0}{0},
 \qquad
 \IRDcapdeco{W_V(x)A_V(x)}
 =\IRcap{darkred}{0}{0}.
\]
The reverse composites cancel the same factors. The nonzero crossings
are recovered by \cite[proof of Theorem~8.6]{GRS}; for the crossings
between zero and a positive color, Definitions~\ref{natrual-basic-0i}
and~\ref{cross-red123}(1) are inverse substitutions.

{For the zero crossing, use $z=(x_2-x_1)/2$ as in
\eqref{eq:ir-zero-crossing-definition}. Keeping each decoration on
its displayed side of the crossing, the two composites are}
\begin{align}
 &\begin{tikzpicture}[baseline=0,scale=.67]
  \draw[black,thick](-.42,-.60)--(.42,.60)
   (.42,-.60)--(-.42,.60);
  \draw[black,densely dotted](-.294,-.42)--(.294,-.42)
   (-.294,.42)--(.294,.42);
  \node[anchor=west]at(.52,-.42){$\scriptstyle A_V(z)W_V(z)$};
  \node[anchor=west]at(.52,.42){$\scriptstyle M_V(z)B_V(z)$};
  \node at(-.42,-.80){$\scriptstyle0$};
  \node at(.42,-.80){$\scriptstyle0$};
 \end{tikzpicture}
 +\IRFparallel{black}{0}{0}
   {M_V(z)g_V(x_1,x_2)W_V(z)+h_V(x_1,x_2)}
 =\IRcross{black}{0}{0},
 \label{eq:ic-zero-crossing-composite-black}\\[3pt]
 &\begin{tikzpicture}[baseline=0,scale=.67]
  \draw[darkred,thick](-.42,-.60)--(.42,.60)
   (.42,-.60)--(-.42,.60);
  \draw[darkred,densely dotted](-.294,-.42)--(.294,-.42)
   (-.294,.42)--(.294,.42);
  \node[anchor=west]at(.52,-.42){$\scriptstyle W_V(z)A_V(z)$};
  \node[anchor=west]at(.52,.42){$\scriptstyle B_V(z)M_V(z)$};
  \node at(-.42,-.80){$\scriptstyle0$};
  \node at(.42,-.80){$\scriptstyle0$};
 \end{tikzpicture}
 +\IRFparallel{darkred}{0}{0}
   {B_V(z)h_V(x_1,x_2)A_V(z)+g_V(x_1,x_2)}
 =\IRcross{darkred}{0}{0}.
 \label{eq:ic-zero-crossing-composite-red}
\end{align}
The parallel terms vanish by \eqref{eq:gV-definition} and
\eqref{eq:ir-zero-crossing-correction}, and the remaining decorations
cancel. Crossings with a negative color follow from the cup--cap
formulas in Definition~\ref{cross-red123}(2) and
\eqref{eq:ir-global-mates-diagram}, including the prescribed scalars.

The original and reconstructed black crossings satisfy the dot
relations of Lemma~\ref{lem:ir-global-affine} with the same dots,
cups and caps. Unless the output colors are the input colors in
reverse order, their difference is killed by an invertible difference
of dot actions, and is therefore zero.
In particular,
\begin{equation}\label{eq:ic-opposite-component-recovery}
\begin{aligned}
 \IRFcross{black}{i}{-i}{k}{-k}{}
 &=\delta_{k,i}\,
   \IRFparallel{black}{i}{-i}{(x_1-x_2)^{-1}}\\
 &\quad-\IRFcupcap{black}{i}{-i}{k}{-k}{(x_1+x_2)^{-1}}
 \qquad(k\ne-i).
\end{aligned}
\end{equation}
The straight term is omitted for $k\ne i$. In the cup--cap term,
$x_1,x_2$ are the dots at the lower and upper ends of the connector;
its denominator has nonzero residue $i+k$. This includes
$00\leftrightarrow(k,-k)$.
All other components are those of Definition~\ref{cross-red123}(3)--(5).
Finally, the cup, cap and dot comparison recovers the real bubbles,
and the inverse-series relation recovers the fake bubbles.
\end{proof}

\begin{Theorem}\label{thm:ic-isomorphism}
For the integral $\mathbf u$-admissible parameters fixed above,
there is an isomorphism of locally unital $\mathbb C$-algebras
\begin{equation}\label{eq:ic-isomorphism}
 \Phi:Q_m\xrightarrow{\ \sim\ }B_{\mathbf u},\qquad
 \Phi(1_{\mathbf i})=e_{\mathbf i}.
\end{equation}
In particular, it identifies
\[
 \operatorname{Hom}_{\mathcal Q(m)}
 (E_{\mathbf j}1_\kappa,E_{\mathbf i}1_\kappa)
 \cong e_{\mathbf i}B_{\mathbf u}e_{\mathbf j}.
\]
\end{Theorem}

\begin{proof}
Apply the argument of \cite[Theorem~8.6]{GRS}. Evaluation of the
forward action at $P$ factors through $\mathcal Q(m)$ and gives
\[
 \Phi:Q_m\longrightarrow B_{\mathbf u},\qquad
 \Phi(1_{\mathbf i})=e_{\mathbf i},
\]
using
$\operatorname{Hom}_{B_{\mathbf u}}(e_{\mathbf j}B_{\mathbf u},
e_{\mathbf i}B_{\mathbf u})=e_{\mathbf i}B_{\mathbf u}e_{\mathbf j}$.
Conversely, the reverse action at $P'$ factors through
$\mathcal{CB}_{\mathbf u}$. Since
\[
 E^dP'\cong\bigoplus_{\mathbf i\in I^d}1_{\mathbf i}Q_m
 \cong q_dQ_m,\qquad q_d=\sum_{\mathbf i\in I^d}1_{\mathbf i},
\]
with finitely many nonzero summands, it defines
$\Psi:B_{\mathbf u}\to Q_m$. On $E_{\mathbf i}P'$, the reconstructed
black dots have the prescribed eigenvalues and nilpotent shifts.
Their joint generalized eigenspace projectors therefore select exactly
these summands, giving $\Psi(e_{\mathbf i})=1_{\mathbf i}$.
Lemma~\ref{lem:ic-inverse-generators} shows that the two composites
fix all generators and distinguished idempotents. Thus
$\Psi\Phi=\operatorname{Id}_{Q_m}$ and
$\Phi\Psi=\operatorname{Id}_{B_{\mathbf u}}$.
\end{proof}

We now transport the grading of $Q_m$ from
Lemma~\ref{lem:ic-support-grading} through the ungraded isomorphism
$\Phi$. For $r\in\mathbb Z$, define
\begin{equation}\label{eq:ic-category-grading}
 \operatorname{Hom}_{\mathcal{CB}_{\mathbf u}}^{r}(E^d,E^e)
 :=\Phi\bigl(q_e(Q_m)_rq_d\bigr).
\end{equation}
{Since $q_d$ has degree zero and $\Phi(q_d)=1_d$, these spaces
give a direct-sum grading with identities in degree zero.
Composition adds degrees because $\Phi$ is an algebra homomorphism.
Thus $\mathcal{CB}_{\mathbf u}$ is a graded $\mathbb C$-linear
category and $\Phi$ is an isomorphism of graded path algebras.

}\begin{Cor}\label{cor:ic-fixed-rank}
Let $\mathcal B_{a,d}(\mathbf u)$ denote the $\mathbf u$-admissible
cyclotomic Brauer (Nazarov--Wenzl) algebra of level $a$ on $d$ strands.
Put
\[
 X_d=\bigoplus_{\mathbf i\in I^d}E_{\mathbf i}1_\kappa,
 \qquad q_d=\sum_{\mathbf i\in I^d}1_{\mathbf i},
\]
{omitting zero summands. By Lemma~\ref{lem:ic-support-grading},
these sums are finite and $q_d$ has degree zero. We have}
\begin{equation}\label{eq:ic-fixed-rank}
 \mathcal B_{a,d}(\mathbf u)
 \cong\operatorname{End}_{\mathcal{CB}_{\mathbf u}}(E^d)
 \cong q_dQ_mq_d
 \cong\operatorname{End}_{\mathcal Q(m)}(X_d).
\end{equation}
The grading of $\mathcal{CB}_{\mathbf u}$ restricts to a
$\mathbb Z$-grading on $\mathcal B_{a,d}(\mathbf u)$, and the
isomorphisms in \eqref{eq:ic-fixed-rank} are graded.
\end{Cor}
\begin{proof}
Use \cite[(4.45)]{RS-cyc} and restrict the isomorphism in
Theorem~\ref{thm:ic-isomorphism} to the subalgebras
$1_dB_{\mathbf u}1_d$ and $q_dQ_mq_d$, as in
\cite[Corollary~8.7]{GRS}.
The grading follows by taking $e=d$ in
\eqref{eq:ic-category-grading}.
\end{proof}

\subsection{Parabolic category $\mathcal O$ in types $B,C,D$}
\label{subsec:parabolic-category-O}

Let $\mathfrak g$ be a complex simple Lie algebra of type $B_n$, $C_n$, or
$D_n$, and let $V$ be its natural representation.  Write
$\varepsilon_{\mathfrak g}=1$ in orthogonal type and
$\varepsilon_{\mathfrak g}=-1$ in symplectic type.  Normalize the Casimir
$C$ by its eigenvalue
$(\lambda,\lambda+2\rho)$ on a highest-weight module of highest
weight $\lambda$, where $(\varepsilon_i,\varepsilon_j)=\delta_{i,j}$, and put
\[
 \Omega=\tfrac12\bigl(\Delta(C)-C\otimes1-1\otimes C\bigr),
 \qquad c_V=\dim V-\varepsilon_{\mathfrak g}.
\]
For the invariant form $(\ ,\ )$ on $V$, choose dual bases satisfying
$(v_a^*,v_b)=\delta_{a,b}$.  The affine Brauer action on $E=V\otimes-$ is
given by
\begin{align*}
 \lcup_M(m)&=\sum_a v_a\otimes v_a^*\otimes m,&
 \lcap_M(v\otimes w\otimes m)&=(v,w)m,\\
 \swap_M(v\otimes w\otimes m)
   &=\varepsilon_{\mathfrak g}w\otimes v\otimes m,&
 \left(\begin{tikzpicture}[baseline=-.5ex,scale=.5]
   \draw[thick] (0,-.5)--(0,.5);
   \fill (0,0) circle (.1);
 \end{tikzpicture}\right)_M
   &=\varepsilon_{\mathfrak g}
       \bigl(\Omega_{V,M}+\tfrac12c_V\operatorname{Id}_{V\otimes M}\bigr).
\end{align*}
These are the maps of \cite[Theorem~3.8 and (3.14)--(3.16)]{RS-cyc},
after horizontal reflection, the tensor flip, and simultaneous rescaling
of cup and cap by $\varepsilon_{\mathfrak g}$. The affine relation is
\[
 \Omega_{V,V}\circ\operatorname{flip}_{V,V}
 =\operatorname{Id}_{V\otimes V}-(\lcup\circ\lcap)_{\mathbb C}.
\]
The dot--cup and dot--cap relations follow from invariance of the form.

Let $\mathcal O_{\mathbb Z}$ be the full subcategory of ordinary BGG
category $\mathcal O$ whose simple composition factors have highest weights
$\lambda=\sum_{j=1}^n\lambda_j\varepsilon_j$,
$\lambda_j\in\mathbb Z$.  For a standard parabolic subalgebra
$\mathfrak p$, put
$\mathcal O_{\mathbb Z}^{\mathfrak p}
 =\mathcal O_{\mathbb Z}\cap\mathcal O^{\mathfrak p}$.
Both categories are stable under $E$, so the affine Brauer action restricts
to them.

\begin{Theorem}\label{thm:parabolic-O-i-action}
The dot spectrum on $\mathcal O_{\mathbb Z}^{\mathfrak p}$ is contained in
$\tfrac12+\mathbb Z$ in types $C,D$, and in $\mathbb Z$ in type $B$.
It is the full indicated lattice on ordinary $\mathcal O_{\mathbb Z}$.
Moreover, $\mathcal O_{\mathbb Z}^{\mathfrak p}$ carries a nilpotent
$2$-representation of the even component $\mathfrak U^\imath_+$ of the
half-integral quasi-split type $\operatorname{AIII}$ $2$-category in types
$C,D$, and of the integral $\mathfrak U^\imath$ in type $B$.
\end{Theorem}

\begin{proof}
The tensor identity and a Borel-stable flag of $V$ give a Verma filtration
of $V\otimes M(\lambda)$ with factors $M(\lambda+\nu)$, where $\nu$ runs
through the weights of $V$.  The displayed Casimir formula acts on this
factor by
\[
 \frac{\varepsilon_{\mathfrak g}}2
       \bigl((\lambda+\nu,\lambda+\nu+2\rho)
                 -(\lambda,\lambda+2\rho)\bigr).
\]
For $\nu=\pm\varepsilon_j$ this is
\[
 \varepsilon_{\mathfrak g}
       \bigl(\pm(\lambda_j+\rho_j)+\tfrac12\bigr),
 \qquad
 \rho_j=
 \begin{cases}
  n-j+1,&\mathfrak g=\mathfrak{sp}_{2n},\\
  n-j,&\mathfrak g=\mathfrak{so}_{2n},\\
  n-j+\tfrac12,&\mathfrak g=\mathfrak{so}_{2n+1}.
 \end{cases}
\]
The weight $\nu=0$ in type $B$ gives eigenvalue $0$.
{Exactness of tensoring, the Verma surjections onto simple objects,
and finite length give containment for both categories. Varying
$\lambda_1\in\mathbb Z$ in the $\nu=\varepsilon_1$ factor gives every
point of the stated lattice. Exactness of the generalized-eigenspace
functors then gives the same spectrum on simple objects of
ordinary $\mathcal O_{\mathbb Z}$.}

Choose a finite-dimensional basic algebra $A_\chi$ for each block of
$\mathcal O_{\mathbb Z}$ and set $A=\bigoplus_\chi A_\chi$.
The self-biadjoint functor $E$ preserves finitely generated projectives
and extends, with its natural transformations, to
$\operatorname{lfdmod}\text{-}A$. This locally Schurian category has the
same simple objects and spectrum. Apply
Theorem~\ref{thm:integral-forward-action} in type $B$ and
\cite[Theorem~6.16]{GRS} in types $C,D$.
Restrict to the full subcategory $\mathcal O_{\mathbb Z}^{\mathfrak p}$,
{which is closed under $E$ and direct summands.
Finite-dimensional endomorphism spaces give finite color support
and nilpotent shifted dots on each object.}
\end{proof}

Containment can be strict. For $\mathfrak g=\mathfrak{sp}_{2n}$ and
$\mathfrak p=\mathfrak g$, all objects are finite-dimensional. If
$\lambda$ is dominant, the eigenvalue on $L(\lambda+\varepsilon_j)$ is
$-(\lambda_j+n-j+\tfrac32)\leq-\tfrac32$. If
$L(\lambda-\varepsilon_j)$ occurs, then $\lambda_j\geq1$, and its
eigenvalue is $\lambda_j+n-j+\tfrac12\geq\tfrac32$.
Thus neither $\tfrac12$ nor $-\tfrac12$ occurs.

\section{The mixed affine case in odd characteristic}
\label{sec:mixed-affine}

Let $\Bbbk$ be algebraically closed of characteristic $p=2r+1>3$.
{The images of $\mathbb Z$ and $\frac12+\mathbb Z$ in $\Bbbk$ coincide,
so both $0$ and $\pm\frac12$ occur. We combine the zero-color
calculations above with the nonzero-color calculations of \cite{GRS},
using the changes of parameters and normalization below.}

\subsection{The mixed datum and spectral weights}
\label{subsec:mixed-affine-datum}

Put $I=\mathbb F_p\subset\Bbbk$, with arrows $i\to i-1$ and
$\tau(i)=-i$. All color indices are taken modulo $p$. Thus
\[
 -\tfrac12=r,\qquad\tfrac12=-r,\qquad
 I_0=\{0\},\qquad I_1=\{1,\ldots,r-1,-r\}.
\]
The fixed vertex is $0$ and the adjacent two-element orbit is
$\{r,-r\}$. This orientation is opposite to that of Section~4.
Take the ambient realization
\[
 X=\bigoplus_{i\in I}\mathbb Z\Lambda_i\oplus\mathbb ZD,
 \quad Y=\operatorname{Hom}_{\mathbb Z}(X,\mathbb Z),\quad
 \tau(\Lambda_i)=\Lambda_{-i},\quad\tau(D)=D,
\]
where $\langle\alpha_i^\vee,\Lambda_j\rangle=\delta_{ij}$,
$\langle\alpha_i^\vee,D\rangle=0$, and
\[
 \alpha_0=2\Lambda_0-\Lambda_1-\Lambda_{-1}+2D,
 \qquad \alpha_j=2\Lambda_j-\Lambda_{j-1}-\Lambda_{j+1}\quad(j\ne0).
\]
{The roots are independent and sum to $2D$. For $i\in I_1\cup\{0\}$,
extend $\alpha_j\mapsto(-1)^{\delta_{j,i+1}}$ to a character
$c_i:X\to\Bbbk^\times$; for $i\in I_1$, define $c_{-i}$ by
\eqref{ci1}. These extensions exist because $\Bbbk^\times$ is
divisible, and $c_i(D)^2=-1$.}

Use the parameters
\[
 Q_{ij}(u,v)=\begin{cases}
 0,&i=j,\\ u-v,&j=i-1,\\ v-u,&j=i+1,\\1,&\text{otherwise},
 \end{cases}
 \qquad Q^\imath_{ij}=(-1)^{\delta_{i,-j}}Q_{ij},
\]
\begin{gather*}
 \varsigma_0=-1,\quad\varsigma_{-r}=1,\quad
 \varsigma_i=0\ (i\ne0,-r),\\
 \zeta_0=-\tfrac12,\quad\zeta_r=1,\quad\zeta_{-r}=-2,
 \quad\zeta_i=1,\ \zeta_{-i}=-1\ (1\le i<r).
\end{gather*}
These choices have $r_{i,j}=1$ and satisfy
\cite[(3.12)--(3.18)]{BWW25}.
Let $\mathfrak U^\imath_p$ be the $2$-category of
\cite[Definition~3.4]{BWW25}, using the corrected
$Q^\imath$ sign of \cite[(3.24), footnote~2]{BWW25v2}
and the diagrammatic convention of Definition~\ref{def-2ikac}.

In $X^\imath=X/(1+\tau)X$, set
\[
 X^\imath_{p,\mathrm{sp}}=
 \left\{\sum_{i=1}^r\lambda_i\bar\Lambda_i+c_0\bar\Lambda_0
 \ \middle|\ \lambda_i\in\mathbb Z,\ c_0\in\mathbb Z/2\mathbb Z,
 \ c_0+\sum_i\lambda_i\equiv0\pmod2\right\}.
\]
The $\bar D$-coordinate is zero.
For a lift $\hat\lambda\in X$, the integer coordinates satisfy
\[
 \lambda_i=\langle\alpha_i^\vee-\alpha_{-i}^\vee,\hat\lambda\rangle,
 \qquad \lambda_{-i}=-\lambda_i,\qquad \lambda_0=0.
\]
The parameter $c_0$ records the parity at the zero vertex.
The projected roots are
\begin{align*}
 \bar\alpha_0&=0,\\
 \bar\alpha_1&=2\bar\Lambda_1-\bar\Lambda_2-\bar\Lambda_0,\\
 \bar\alpha_j&=2\bar\Lambda_j-\bar\Lambda_{j-1}-\bar\Lambda_{j+1}
 \quad(1<j<r),\\
 \bar\alpha_r&=3\bar\Lambda_r-\bar\Lambda_{r-1},\\
 \bar\alpha_{-j}&=-\bar\alpha_j\quad(1\le j\le r).
\end{align*}
Thus this subgroup is stable under root shifts; denote the associated
full $2$-subcategory by $\mathfrak U^\imath_{p,\mathrm{sp}}$.
{Its projected root lattice has index $p$: in the basis
$\bar\Lambda_i+\bar\Lambda_0$, the root matrix has diagonal
$(2,\ldots,2,3)$, adjacent entries $-1$, and determinant $2r+1=p$.}
Since $\Bbbk^\times$ has no nontrivial $p$-torsion, the characters
\eqref{gammai} on this subgroup are determined by their root values.
In particular, $\gamma_0(\lambda)=(-1)^{c_0}$ and
$\gamma_i(\lambda)\gamma_{-i}(\lambda)=(-1)^{\lambda_i}$ for $1\le i\le r$.

\begin{Prop}\label{prop:mixed-spectral-weights}
Let $\mathcal M$ be a locally Schurian affine Brauer module category
with dot spectrum contained in $I$. For a simple object $L$, write
\[
 m_L(u)=\prod_{i\in I}(u-i)^{\epsilon_i(L)},\qquad
 a_L=\deg m_L,\qquad s_L=\frac{(-1)^{a_L}}2.
\]
Its spectral weight lies in $X^\imath_{p,\mathrm{sp}}$, with
\[
 \lambda_i(L)=\epsilon_i(L)-\epsilon_{-i}(L)-\delta_{i,s_L}
 \quad(1\le i\le r),\qquad c_0(L)=\epsilon_0(L)\pmod2.
\]
Moreover,
\begin{equation}\label{eq:mixed-scalar-bubble}
 O_L(u)=(u+r)\prod_{i=1}^r
 \left(\frac{u+i}{u-i}\right)^{\lambda_i(L)},
 \qquad O_L(0)=\frac{(-1)^{c_0(L)+1}}2.
\end{equation}
The central block decomposition for the dotted-bubble action is
\[
 \mathcal M\simeq\prod_{\lambda\in X^\imath_{p,\mathrm{sp}}}
 \mathcal M_\lambda,
 \qquad E_j\mathcal M_\lambda\subseteq\mathcal M_{\lambda-\bar\alpha_j}.
\]
Only finitely many weights occur on a finite-length object.
\end{Prop}
\begin{proof}
The identity of \cite[Lemma~4.16]{GRS}, valid in characteristic
different from $2$, is
\[
 O_L(u)=\left(u-\frac{(-1)^{a_L}}2\right)
 \prod_{i\in I}\left(\frac{u+i}{u-i}\right)^{\epsilon_i(L)}.
\]
Cancel the zero factors and, when $a_L$ is odd, use
$u-r=(u+r)((u+r)/(u-r))^{-1}$. This gives
\eqref{eq:mixed-scalar-bubble} and
$c_0+\sum_i\lambda_i\equiv a_L-\mathbf1_{a_L\text{ odd}}\equiv0$.
{For $1\le i\le r$, the orders
$-\operatorname{ord}_{u=i}O_L(u)=\lambda_i(L)$ are integers,
not residues modulo $p$. They determine the weight, including
$c_0$ by the parity condition, so the usual generalized-character
decomposition applies.}

For a simple subquotient $K$ of $E_jL$, \cite[Lemma~4.17]{GRS} gives
\begin{equation}\label{eq:mixed-forward-slide}
 \frac{O_K(u)}{O_L(u)}
 =\frac{((u+j)^2-1)(u-j)^2}{((u-j)^2-1)(u+j)^2}.
\end{equation}
Orders at the nonzero residues give $\lambda\mapsto\lambda-\bar\alpha_j$.
At $j=r$, the equality $r+1=-r$ gives the coefficient $3$;
at $j=\pm1$, evaluation at zero changes the sign and hence $c_0$.
The same calculation in the finite-dimensional coefficient algebra
gives the shift of generalized characters; passage to filtered
colimits gives the assertion on all objects.
\end{proof}

\subsection{The mixed bubble normalization}
\label{subsec:mixed-normalization}

Let $\mathcal M$ be a locally Schurian or locally finite
Abelian nilpotent $2$-representation of
$\mathfrak U^\imath_{p,\mathrm{sp}}$, as in
Definition~\ref{def:integral-nilpotence}. For $i\in I_1$, orient $E_i$
upwards and $E_{-i}$ downwards, with the negative dot on a downward
strand. Use the oriented notation of Definition~\ref{evalu-bub},
with cup and cap coefficients
\[
 \begin{array}{c|cc}
 \text{endpoints}&\text{cup}&\text{cap}\\ \hline
 (-i,i)&1&(\zeta_i\gamma_i(\lambda))^{-1}\\
 (i,-i)&\zeta_i\gamma_i(\lambda)&1
 \end{array}
\]
in exterior region $\lambda$; use the same factors for fake bubbles.
Here $\gamma_i(\lambda\pm\bar\alpha_i)=\gamma_i(\lambda)$, so
both adjunctions hold. Explicitly,
\[
\begin{aligned}
 \circlearrowleft_{\lambda,i}(u)
 &=\frac1{\zeta_i\gamma_i(\lambda)}
 \begin{tikzpicture}[baseline=-1mm,iQ]
  \draw[-](-.25,0) arc(180:-180:.25);
  \node at(-.42,0){\strandlabel{-i}};\node at(.55,0){$(u)$};
  \region{.95,0}{\lambda};
 \end{tikzpicture},&
 \circlearrowright_{\lambda,i}(u)
 &=-\zeta_i\gamma_i(\lambda)
 \begin{tikzpicture}[baseline=-1mm,iQ]
  \draw[-](-.25,0) arc(180:-180:.25);
  \node at(-.42,0){\strandlabel{i}};\node at(.63,0){$(-u)$};
  \region{1.05,0}{\lambda};
 \end{tikzpicture},\\
 O_{\lambda,0}(u)&=-2u(-1)^{c_0}
 \begin{tikzpicture}[baseline=-1mm,iQ]
  \draw[-](-.25,0) arc(180:-180:.25);
  \node at(-.42,0){\strandlabel{0}};\node at(.55,0){$(u)$};
  \region{.95,0}{\lambda};
 \end{tikzpicture}.
\end{aligned}
\]

\begin{Lemma}\label{lem:mixed-normalization}
For nonzero finitely generated $V\in\mathcal M_\lambda$, let
$\varepsilon_i(V)$ be the nilpotence exponent on $E_iV$, or $0$
if $E_iV=0$. For $j\in I_1$,
\begin{align*}
 \circlearrowleft_{V,j}(u)&=u^{-\lambda_j+\delta_{j,-r}}(1+O(u^{-1})),&
 \circlearrowright_{V,j}(u)&=2^{\delta_{j,-r}}u^{\lambda_j}(1+O(u^{-1})),\\
 \circlearrowleft_{V,j}(u)\circlearrowright_{V,j}(u)
 &=2^{\delta_{j,-r}}u^{\delta_{j,-r}},&
 O_{V,0}(u)O_{V,0}(-u)&=1.
\end{align*}
Also $O_{V,0}(u)=1+O(u^{-1})$, and
$u^{\varepsilon_0(V)}O_{V,0}(u)$,
$u^{\varepsilon_j(V)}\circlearrowleft_{V,j}(u)$ and
$u^{\varepsilon_{-j}(V)}\circlearrowright_{V,j}(u)$ are polynomials.
For the affine dot $i$ minus the original colored dot, put
\begin{equation}\label{eq:mixed-total-normalization}
 \widehat O_V(u)=O_{V,0}(-u)
 \prod_{j\in I_1}(-1)^{\lambda_j+\delta_{j,-r}}2^{-\delta_{j,-r}}
 \circlearrowleft_{V,j}(j-u)\circlearrowright_{V,j}(u+j).
\end{equation}
This series has leading term $u$, and
\begin{equation}\label{eq:mixed-formal-inverse}
 \widehat O_V(u)\widehat O_V(-u)=\frac14-u^2.
\end{equation}
For $m_V(u)=\prod_{i\in I}(u-i)^{\varepsilon_i(V)}$,
{$n_V(u)=m_V(u)\widehat O_V(u)$ is a monic polynomial of degree
$\deg m_V+1$.}
If $V$ is simple, then
\begin{equation}\label{eq:mixed-normalized-scalar}
 \widehat O_V(u)=(u+r)\prod_{i=1}^r
 \left(\frac{u+i}{u-i}\right)^{\lambda_i}.
\end{equation}
\end{Lemma}
\begin{proof}
The leading and product identities are \eqref{ibubblerel} and
\eqref{iinfgrass}; the exceptional product is $2u$.
Polynomiality is the coefficient argument of \cite[Lemma~7.4]{GRS}:
a sufficiently dotted real bubble factors through an annihilating power
of the dot, and the fake part is a polynomial. At zero, use the
same argument with the factor $-2u(-1)^{c_0}$, as in
Definition~\ref{evalu-bub}. Thus the only denominators in
\eqref{eq:mixed-total-normalization} are removed by $m_V$.
Pairing the local products at $u$ and $-u$ gives
$r^2-u^2=1/4-u^2$; their leading terms give degree $\deg m_V+1$.
On a simple object, the local polynomial products are scalar
multiples of powers of $u$, so unique factorization gives
\eqref{eq:mixed-normalized-scalar}, as in \cite[Lemma~7.6]{GRS}.
\end{proof}

\subsection{Compatibility of the forward and reverse constructions}
\label{subsec:mixed-comparison}

First pass to arrows $i\to i+1$ by setting
\[
 Q'_{ij}(u,v)=Q_{ij}(-u,-v),\quad
 c'_i(\mu)=(-1)^{\langle\alpha_i^\vee,\mu\rangle}c_i(\mu),\quad
 \zeta'_i=(-1)^{\varsigma_i-1}\zeta_i.
\]
The isomorphism from this presentation to the original one fixes
cups and caps and negates dots and crossings. It sends each bubble
series, including its fake part, to minus that series at $-u$.
Substitution in \eqref{ibubblegeneratingfunction}--\eqref{ibraid}
verifies the isomorphism: a curl, quadratic relation and braid
relation acquire, respectively, the common signs $-1,1,-1$.
In particular $\gamma'_i=(-1)^{\lambda_i}\gamma_i$ for $i\ne0$,
while $\gamma'_0=\gamma_0$ and $\zeta'_0=\zeta_0$.
All red generators below are in this primed presentation.

Set $K=\{-r,\ldots,-1\}$ and $J=I_1$.
We use $K$ to orient the nonzero strands and $J$ in the products
defining $W$ and $M$. The nonzero residues lift to the half-integral path
by $a\mapsto a-p/2$ for $1\le a\le2r$; the lifts of $K$ are
$1/2,3/2,\ldots,r-1/2$. {Their sums and differences have absolute value at most $p-2$.
Thus reduction preserves equality, opposition and adjacency on the
nonzero path, as required for the primary-factor formulas in
\cite[Remark~5.2]{GRS}.}

For $i,j\in K$, put $t_{ij}=-1$ if $j=i-1$ and $t_{ij}=1$ otherwise.
Multiply cups, caps and crossings by the following coefficients,
where $\lambda$ is the exterior region:
\[
 \begin{array}{c|cc}
 \text{endpoints}&\text{cup}&\text{cap}\\ \hline
 (-i,i)&1&(\zeta'_i\gamma'_i(\lambda))^{-1}\\
 (i,-i)&(-1)^{\delta_{i,-r}}\zeta'_i\gamma'_i(\lambda)&1
 \end{array}
\]
\[
 \begin{array}{c|cccc}
 \text{source}&(i,j)&(i,-j)&(-i,j)&(-i,-j)\\ \hline
 \text{crossing}&1&1&(-1)^{\delta_{i,j}\delta_{i,-r}}&t_{ij}
 \end{array}
\]
The bubble coefficient is the product of its cup and cap coefficients.
Indeed,
\begin{equation}\label{eq:mixed-character-table}
 \gamma'_i(\bar\alpha_j)=t_{ij}(-1)^{\delta_{i,j}\delta_{i,-r}},
 \qquad\gamma'_i(\lambda)\gamma'_{-i}(\lambda)=(-1)^{\lambda_i}.
\end{equation}
The resulting parameters are those of
\cite[(3.56)]{BWW25}. The same substitution in
\eqref{ibraid} multiplies every term by the product of the three
crossing coefficients, including its two cup--cap corrections.
Thus \cite[Theorem~3.8]{BWW25}, with weight $-\lambda$,
gives exactly the oriented convention of \cite[Definitions~3.1--3.2]{GRS}.
{At $i=-r$, \eqref{eq:mixed-character-table} gives
$\gamma'_{-r}(\bar\alpha_{-r})=-1$, which accounts for the
extra minus sign in the cup coefficient.}

Use \cite[Definitions~7.11 and~7.13]{GRS} for the nonzero black
generators, with these $K$-adjunctions and the polynomial $n_V$.
The black dot on $E_i$ is $i$ plus its primed dot, hence $i$ minus
its original dot. In the sums defining the crossings, include all
colors of $I$, also zero.
For each residue $i$, write
\[
 m_{V,i}(u)=\frac{m_V(u)}{(u-i)^{\varepsilon_i(V)}},\qquad
 n_{V,i}(u)=\frac{n_V(u)}{t_{V,i}(u)},
\]
{where $t_{V,i}$ is the monic factor of $n_V$ whose reduction is
a power of $u-i$. The bubble coefficients generate a finite-dimensional
local algebra, by Lemma~\ref{lem:mixed-normalization} and the argument
of \cite[Lemma~7.8]{GRS}. The factors therefore exist uniquely by
\cite[Lemma~4.28 and Remark~7.9]{GRS}.}

The bubble slide \eqref{ibubslide}, multiplied over $I$, gives
\begin{equation}\label{eq:mixed-total-primary-slide}
 \frac{\widehat O_{E_jV}(u)}{E_j(\widehat O_V(u))}
 =\frac{((u+x)^2-1)(u-x)^2}{((u-x)^2-1)(u+x)^2},
\end{equation}
where $x$ is the black dot of residue $j$.
Here the product of the factors $R'_{i,j}(u-i,x-j)$ is
$-((u-x)^2-1)/(u-x)^2$; this includes $0,r,-r$.
Define
\begin{align*}
 c_V^+(u)&=\prod_{i\in J}\frac{t_{V,i}(u)}{(u-i)^{\varepsilon_i(V)}},&
 c_V^-(u)&=\prod_{i\in-J}\frac{t_{V,i}(u)}{(u-i)^{\varepsilon_i(V)}},\\
 W_V(u)&=c_V^+(u),&M_V(u)&=2(-1)^{c_0+1}c_V^-(u).
\end{align*}
These are units at $u=0$. On an ordinary orbit,
$\circlearrowleft'_{V,-j}(v)=(-1)^{\lambda_j}
\circlearrowright'_{V,j}(-v)$ and the analogous reversed identity
show that the factors $t_{V,i}$ do not depend on the chosen
representative. The actual nonzero cups and caps continue to use $K$.
The local bubble products and \eqref{eq:mixed-total-primary-slide} give
\begin{align*}
 W_V(u)M_V(-u)&=1-2u,&M_V(u)W_V(-u)&=1+2u,\\
 \frac{W_{E_0V}(u)}{E_0(W_V(u))}&=\frac{u+x-1}{u-x-1},&
 \frac{M_{E_0V}(u)}{E_0(M_V(u))}&=\frac{u+x+1}{u-x+1}.
\end{align*}
Use these $W,M$ in Definition~\ref{natrual-basic-red}(10)--(11)
and \eqref{eq:ir-zero-crossing-correction}--\eqref{eq:ir-zero-crossing-definition}.

For $k\in K$, let $u$ be the zero dot and $x$ the black dot of
residue $k$. {Fix a square root of $-1$ in $\Bbbk$, and let
$f_k(u,x)\in\Bbbk[[u,x-k]]$ be even in $u$ and satisfy}
\begin{equation}\label{eq:mixed-K-factors}
 f_k(u,x)^2=\begin{cases}
 \dfrac1{(u^2-x^2)(u^2-(x-1)^2)},&k=-r,\\[4pt]
 u^2-x^2,&k=1-r,\\[2pt]
 \dfrac{u^2-x^2}{u^2-(x-1)^2},&1-r<k\le-1,
 \end{cases}
\end{equation}
with $f_{-r}(0,x)=1/(x(x-1))$,
$f_{1-r}(0,x)=\sqrt{-1}\,x$, and
$f_k(0,x)=x/(x-1)$ in the last range, which is empty for $p=5$.
{The denominators are units at $(u,x)=(0,k)$.
The stated initial values determine the square roots by coefficient
comparison, dividing only by twice a unit. Nilpotence makes their
evaluations finite. Decorate the input strands of the off-diagonal
red crossings by the following coefficients to obtain the black crossings:}
\begin{equation}\label{eq:mixed-K-crossing-factors}
 \begin{array}{c|cc}
 \text{source}&(k,0)&(0,k)\\ \hline
 \text{coefficient}&
 f_k(u,x)\dfrac{1+u-x}{x-u}&
 f_k(u,x)^{-1}\dfrac{(1+x-u)^{1-\delta_{k,-1}}}{u-x}.
 \end{array}
\end{equation}
In particular the nonunit factor $1+x-u$ is omitted at $k=-1$.
The same $K$-adjunctions define the crossings with $-k$.
All straight and cup--cap entries are those of
\cite[Definition~7.13(3)--(5)]{GRS}, with the black cups and caps
already defined. {Extract the $J$-factors from \eqref{eq:mixed-total-primary-slide}
and multiply by \eqref{eq:mixed-K-factors}. In each of its three cases
this gives}
\begin{equation}\label{eq:mixed-K-transport-factor}
 \frac{W_{E_kV}(u)}{E_k(W_V(u))}f_k(u,x)^2
 =\frac{(u-x)(u+x+1)^{1-\delta_{k,-1}}}{(u+x)(u-x+1)}.
\end{equation}
This is the mixed transport of Section~\ref{sec:integral-reverse-complete},
with the negative color obtained by adjunction. Multiplying it by
$(1+u-x)/(x-u)$ gives $-(u+x+1)^{1-\delta_{k,-1}}/(u+x)$,
the other mixed coefficient after bending the zero dot.

To apply the previous calculations, use the same polynomial algebra
in the independent coefficients of an even series $H(u)$, with its
constant coefficient inverted, for all objects. Tensor each morphism
space with this algebra over $\Bbbk$ and extend the functors coefficientwise by
\[
 E_k(H(u))=H(u)f_k(u,x),\qquad
 E_{-k}(H(u))=H(u)f_k(u,-x)^{-1},\qquad E_0(H(u))=H(u).
\]
{The extension is faithful. The functors are $\Bbbk$-linear and act
on the added coefficients as displayed. Naturality follows from the
dot relations: equal-color products are symmetric, and opposite-color
products equal $1$ at opposite dots. Their differences from $1$
are therefore divisible by the sum of the dots, which kills the
cup--cap terms. This also applies to the black generators:
straight terms commute with the factors, and cup--cap terms restrict
to $1$ at both ends.
}
Take $\Gamma=H^2$ and replace $W,M$ by $W/H^2,MH^2$.
Equation~\eqref{eq:mixed-K-transport-factor} is then precisely
the transport used in the integral calculations.
Conjugating every zero strand by $H(x)^{-1}$ gives the displayed
generators. For the zero crossing this follows directly from
\eqref{eq:ir-zero-crossing-definition}: the quotient
\[
 \frac{H(x_1)H(x_2)}{H((x_2-x_1)/2)^2}
\]
{is symmetric and equals $1$ at $x_2=-x_1$. Its difference from
$1$ is the sum of the dots times a series in their sum and product.
By \eqref{eq:ir-red-affine-left}--\eqref{eq:ir-red-affine-right}
and \eqref{eq:ir-red-dots}, the sum commutes with the red crossing
and kills cup--cap terms, while the commutator with the product
consists of dotted cup--caps. Hence the displayed quotient commutes
with the red crossing and the reconstructed zero crossing.
The conjugating factors cancel on each intermediate color word.
Since the resulting coefficients lie in the original morphism spaces,
the identities descend by faithfulness. This extension is used only
for finite diagram calculations; no locally Schurian recognition
theorem is applied to it.
}
Among triples containing zero, the only new adjacency is in the
six permutations of $(0,r,-r)$. Their relation \eqref{ibraid} has
zero correction, since their entries are distinct and
$Q^\imath_{r,0}=Q^\imath_{-r,0}=1$. For example,
\begin{equation}\label{eq:mixed-three-colour-target}
\begin{tikzpicture}[iQ,centerzero,scale=1.1]
 \draw[-](-.4,-.4)\botlabel{0}--(.4,.4)\toplabel{0};
 \draw[-](0,-.4)\botlabel{r} to[out=135,in=down](-.32,0)
  to[out=up,in=225](0,.4)\toplabel{r};
 \draw[-](.4,-.4)\botlabel{-r}--(-.4,.4)\toplabel{-r};
\end{tikzpicture}
\;=\;
\begin{tikzpicture}[iQ,centerzero,scale=1.1]
 \draw[-](-.4,-.4)\botlabel{0}--(.4,.4)\toplabel{0};
 \draw[-](0,-.4)\botlabel{r} to[out=45,in=down](.32,0)
  to[out=up,in=-45](0,.4)\toplabel{r};
 \draw[-](.4,-.4)\botlabel{-r}--(-.4,.4)\toplabel{-r};
\end{tikzpicture}.
\end{equation}

\subsection{The mixed constructions}
\label{subsec:mixed-recognition}

Since $I=\mathbb F_p$ is finite, finite color support is
automatic. Thus generalized nilpotence and nilpotence coincide for
this datum under Definition~\ref{def:integral-nilpotence}; dot
nilpotence is required on finitely generated objects.

\begin{Cor}\label{cor:mixed-recognition}
An affine Brauer action on a locally Schurian $\Bbbk$-linear category
with dot spectrum contained in $I$ gives a nilpotent
$2$-representation of $\mathfrak U^\imath_{p,\mathrm{sp}}$.
Conversely, a nilpotent $2$-representation of this $2$-category on
a locally Schurian category gives an affine Brauer action with
spectrum contained in $I$.
Color $i$ acts by $E_i$, and its colored dot is $i$ minus the affine dot.
\end{Cor}
\begin{proof}
We use the constructions above on finitely generated objects and
their iterated images. The polynomial evaluations are justified
by Lemma~\ref{lem:finite-dd-evaluation}.

In the reverse direction, the nonzero adjunction and dot
calculations are \cite[Lemma~7.14 and Corollary~7.15]{GRS};
at zero use Lemma~\ref{lem:ir-global-adjunction}.
{The bubble calculation of \cite[Proposition~7.20]{GRS},
with the zero terms of Lemma~\ref{lem:red-black-series}, gives
the full series $n_V/m_V$, including its constant coefficient:
$c_V^+(0)c_V^-(0)=(-1)^{c_0+1}/2$.}
The cup--crossing and cap relations are
\cite[Lemmas~7.18 and~7.21]{GRS} and
Lemmas~\ref{lem:ir-full-pitchfork} and~\ref{lem:ir-global-cap}.
The square follows from \cite[Lemma~7.22]{GRS} and
Lemma~\ref{lem:ir-zero-full-square}, specifically
\eqref{eq:ir-zero-square-defect},
summing over all intermediate colors, including $(r,-r)$ and $(-r,r)$.
{The required bubble and transport identities were proved above.}

For the braid, the nonzero cases are the calculations of
\cite[Lemma~7.23, Lemma~7.24, Corollary~7.25 and Lemmas~7.26--7.28]{GRS}.
{With one zero, use the upward calculation and
\eqref{eq:ir-global-mates-diagram}. For $(0,-r,-r)$, the two
$f_{-r}$ factors are symmetric in the equal-color dots, so their
divided difference is zero. Bending a nonzero strand and moving
the zero strand with the invertible mixed crossings gives the six
permutations of $(0,r,-r)$, without inverting the $r,-r$ crossing.}

With two zeros, apply the same adjunction to the entire calculation
of \eqref{fw:mixed-braid}, keeping its two black braid diagrams
separate as in \eqref{eq:ir-mixed-braid-defect}:
\[
 \operatorname{Nat}(E_0^2E_{-k},E_{-k}E_0^2)
 \longrightarrow\operatorname{Nat}(E_kE_0^2,E_0^2E_k),\qquad k\in K.
\]
Every cup--cap and divided-difference term is transported by this
same adjunction. For the coefficients, use
\eqref{eq:mixed-K-transport-factor} and Lemma~\ref{compu-T}, with
Definition~\ref{notation-1} and the shift \eqref{appC:shift}.
{For $j=-k\in\{1,\ldots,r\}$, the original leading coefficient has
nonzero residue $j^2/(j+1)^2$, so it can be cancelled before
applying adjunction. At $k=-1$ use $j=1$; at $j=r$ the residue is
$1$, since $r+1=-r$. No inverse of $4j^2-1$ is needed.}
For three zeros, use \eqref{eq:ir-allzero-braid-defect} with the
same $W,M$ and the complete bubble series.
These calculations prove the black braid and complete the verification
of the affine Brauer relations.

In the forward direction, use Proposition~\ref{prop:mixed-spectral-weights}
and the inverse generator formulas. Their nonzero relations are
\cite[Lemmas~5.10--5.17 and Theorem~6.15]{GRS}; those involving
zero are Lemmas~\ref{fw:basic-relations}--\ref{fw:braid-relation},
with the same four braid cases just described. These are uses of
the cited local calculations: their denominators are units on the
specified summands, and all sums include every intermediate color.
The normalization and change of coefficients above therefore apply
in both directions and give identities in the original category.
Naturality follows from common annihilating polynomials on
$V\oplus V'$, as in Lemma~\ref{lem:ir-normalization}.

{Finally, biadjunction preserves finite generation and filtered
colimits. As in the proof of Theorem~\ref{prop:ir-action}, the
transformations and relations extend uniquely from finitely generated
objects. On finitely generated objects the shifted dots are nilpotent;
finite color support is automatic because $I$ is finite.}
\end{proof}

\begin{rem}
All factors are formed at the distinct residues of $I$ before
evaluation. No infinite product is reduced modulo $p$ and no
division by $p$ is used. For $p>3$, neither $r$ nor $-r$ is adjacent
to $0$. At $p=3$, both are adjacent to $0$, this case is treated in Subsection~10.5. 
\end{rem}
\subsection{Characteristic three}
\label{subsec:mixed-characteristic-three}

Let $\Bbbk$ be algebraically closed of characteristic $3$.
Use the ambient realization and parameters of
Subsection~\ref{subsec:mixed-affine-datum} with $r=1$.
Thus
\[
\begin{gathered}
 I=\{0,1,-1\},\qquad I_1=J=K=\{-1\},\\
 \varsigma_0=-1,\quad\varsigma_{-1}=1,\quad\varsigma_1=0,
 \qquad
 \zeta_0=-\tfrac12,\quad\zeta_{-1}=-2,\quad\zeta_1=1.
\end{gathered}
\]
The spectral subgroup is
\[
 X^\imath_{3,\mathrm{sp}}
 =\{\lambda_1\bar\Lambda_1+c_0\bar\Lambda_0:
   \lambda_1\in\mathbb Z,\ c_0\in\mathbb Z/2\mathbb Z,\
   \lambda_1+c_0\equiv0\pmod2\},\qquad\bar D=0.
\]
The projected roots are
\[
 \bar\alpha_0=0,\qquad
 \bar\alpha_1=3\bar\Lambda_1-\bar\Lambda_0
             =3(\bar\Lambda_1+\bar\Lambda_0),\qquad
 \bar\alpha_{-1}=-\bar\alpha_1.
\]
Thus the root lattice has index $3$ in this subgroup.
Since $\Bbbk^\times$ has no nontrivial $3$-torsion, the character
argument of Subsection~\ref{subsec:mixed-affine-datum} applies.
The parameter identities hold with $r=1$.
Write $\mathfrak U^\imath_{3,\mathrm{sp}}$ for the resulting
full $2$-subcategory.

\begin{Cor}\label{cor:mixed-characteristic-three}
An affine Brauer action on a locally Schurian $\Bbbk$-linear
category with dot spectrum contained in $\mathbb F_3$ gives a
nilpotent $2$-representation of
$\mathfrak U^\imath_{3,\mathrm{sp}}$.
Conversely, a nilpotent $2$-representation of this $2$-category
on a locally Schurian category gives an affine Brauer action
with spectrum contained in $\mathbb F_3$.
Color $i$ acts by $E_i$, and its colored dot is $i$ minus the
affine dot.
\end{Cor}

\begin{proof}
We use the local calculations in the proof of
Corollary~\ref{cor:mixed-recognition} and the results cited there.
We check the coefficients in characteristic $3$ and then the
relations in both directions.
All coefficients are evaluated on nonzero finitely generated
objects and their nonzero iterated images.

The proof of Proposition~\ref{prop:mixed-spectral-weights}
applies with $r=1$. In particular,
\[
 O_L(u)=(u+1)\left(\frac{u+1}{u-1}\right)^{\lambda_1(L)},
 \qquad c_0(L)\equiv\lambda_1(L)\pmod2,
\]
and the orders at $1,-1$ give the root shifts above.
The proof of Lemma~\ref{lem:mixed-normalization} also applies:
the three primary supports are distinct, the exceptional local
product is $2u$, and the zero-color product is unchanged.
Hence the monic polynomials $m_V,n_V$, their primary factors,
and \eqref{eq:mixed-total-primary-slide} are available before
using the braid relation.
In the forward direction this identity is the affine Brauer
bubble slide; in the reverse direction it follows from the
red bubble slides.

We use the primed presentation of
Subsection~\ref{subsec:mixed-comparison}.
The nonzero normalization still has
$\gamma'_{-1}(\bar\alpha_{-1})=-1$, with the exceptional minus
sign in the $K$-cup coefficient.
Put
\[
 W_{*,V}(u)=\frac{t_{V,-1}(u)}{(u+1)^{\varepsilon_{-1}(V)}},
 \qquad
 M_{*,V}(u)=2(-1)^{c_0(V)+1}
           \frac{t_{V,1}(u)}{(u-1)^{\varepsilon_1(V)}}.
\]
By \eqref{eq:mixed-formal-inverse} and
$c_0\equiv\lambda_1\pmod2$,
$W_{*,V}(u)M_{*,V}(-u)=1-2u$.
Set
\begin{equation}\label{eq:char3-WM}
 W_V^{(3)}(u)=(1-2u)M_{*,V}(u),\qquad
 M_V^{(3)}(u)=\frac{W_{*,V}(u)}{1-2u}.
\end{equation}
These are units at zero, with
\begin{gather}
 W_V^{(3)}(u)M_V^{(3)}(-u)=1-2u,\qquad
 M_V^{(3)}(u)W_V^{(3)}(-u)=1+2u,\label{eq:char3-products}\\
 W_V^{(3)}(u)M_V^{(3)}(u)
 =W_{*,V}(u)M_{*,V}(u).\notag
\end{gather}
Thus the zero-color bubble calculation, including its constant
coefficient, is unchanged.

Put
\[
 R(u,x)=
 \frac{(u+x-1)(u+x+1)(u-x)^2}
 {(u-x-1)(u-x+1)(u+x)^2}.
\]
The factors of $R(u,x)$ with reduced root $1$, for $x$ of
residue $0,-1,1$, respectively, are
\[
 \frac{u+x-1}{u-x-1},\qquad
 \frac{1}{(u-x+1)(u+x)^2},\qquad
 (u+x+1)(u-x)^2.
\]
The parity $c_0$ is unchanged by $E_0$ and changes by $1$
under $E_{\pm1}$. Primary-factor uniqueness in
\eqref{eq:mixed-total-primary-slide} gives
\begin{align}
 \frac{W^{(3)}_{E_0V}(u)}{E_0(W_V^{(3)}(u))}
 &=\frac{u+x-1}{u-x-1},&
 \frac{M^{(3)}_{E_0V}(u)}{E_0(M_V^{(3)}(u))}
 &=\frac{u+x+1}{u-x+1},\label{eq:char3-zero-slides}\\
 \frac{W^{(3)}_{E_{-1}V}(u)}{E_{-1}(W_V^{(3)}(u))}
 &=-\frac{1}{(u-x+1)(u+x)^2},&
 \frac{W^{(3)}_{E_1V}(u)}{E_1(W_V^{(3)}(u))}
 &=-(u+x+1)(u-x)^2.\label{eq:char3-raw-slides}
\end{align}

There is a unique series $f(u,x)\in\Bbbk[[u,x+1]]$,
even in $u$, such that
\[
 f(u,x)^2=x^2-u^2,\qquad f(0,x)=x.
\]
Its coefficients are determined successively by division by
$2x$, which is a unit.
Use $f$ for $k=-1$ in
\eqref{eq:mixed-K-crossing-factors}, and the same
$K$-adjunctions for crossings with color $1$.
Together with $W^{(3)},M^{(3)}$, these give the reverse
generators of Subsection~\ref{subsec:mixed-comparison}.
The forward generators are given by the inverse formulas,
using $(W^{(3)})^{-1},(M^{(3)})^{-1}$ for zero-color coefficients.
Both mixed coefficients in
\eqref{eq:mixed-K-crossing-factors} are units at $(u,x)=(0,-1)$.
Moreover, \eqref{eq:char3-raw-slides} gives
\begin{equation}\label{eq:char3-mixed-transport}
 \frac{W^{(3)}_{E_{-1}V}(u)}{E_{-1}(W_V^{(3)}(u))}\,f(u,x)^2
 =\frac{u-x}{(u+x)(u-x+1)}.
\end{equation}

The coefficient-extension argument of
Subsection~\ref{subsec:mixed-comparison} now applies.
Adjoin the coefficients of an even series $H(u)$ with invertible
constant term and set
\[
 E_{-1}(H(u))=H(u)f(u,x),\qquad
 E_1(H(u))=H(u)f(u,-x)^{-1},\qquad E_0(H(u))=H(u).
\]
Equal-color factors are symmetric; for an opposite pair their
product is $1$ when the dots sum to zero.
The dot, cup and cap relations give naturality as before,
and the extension is faithful.
For $\overline W=W^{(3)}/H^2$ and $\overline M=M^{(3)}H^2$,
the zero slides and products are unchanged, while
\begin{equation}\label{eq:char3-normalized-slides}
 \frac{\overline W_{E_1V}(u)}{E_1(\overline W_V(u))}
 =\frac{(u-x)(u+x+1)}{u+x},\qquad
 \frac{\overline W_{E_{-1}V}(u)}{E_{-1}(\overline W_V(u))}
 =\frac{u-x}{(u+x)(u-x+1)}.
\end{equation}
These are the positive $j=1$ transport of
Lemma~\ref{lem:ir-W-positive-transport} and its negative-color
form. The transports of $\overline M$ follow from
\eqref{eq:char3-products}.
Conjugation by $H(x)^{-1}$ on each zero strand gives the original
generators, as in Subsection~\ref{subsec:mixed-comparison}.

In the calculations below, take
$A=\overline W^{-1}$ and $B=\overline M^{-1}$.
Use finite polynomial representatives with sufficiently high
vanishing orders at $1,-1$, as in
Definition~\ref{notation-1} and Lemma~\ref{lem:finite-dd-evaluation}.
In particular,
\[
 A_V(u)B_V(-u)=\frac{1}{1-2u},\qquad
 A_{E_1V}(u)=E_1(A_V(u))
            \frac{u+x}{(u-x)(u+x+1)}.
\]
These are the two coefficient identities used in
Lemma~\ref{compu-T} for $\delta=1$.
The denominators involving a zero and a nonzero dot have
nonzero residues. Divided differences at coincident zero dots
are formed before evaluation.
Finite representatives retain all intermediate-color terms in
the bubble calculations.

We check the forward construction.
The nonzero-color generators, adjunctions and local relations
are those of the half-integral calculations in
\cite[Sections~5--6]{GRS}, with $-1=\tfrac12$ and
$1=-\tfrac12$, in the normalization of
Subsection~\ref{subsec:mixed-comparison}.
For zero and mixed colors, the calculations of
Lemmas~\ref{fw:basic-relations}--\ref{fw:quadratic-relation}
and Section~\ref{sec:forward-nil-brauer} use the products and
slides above. These prove all relations except the braid.

The bubble relations also follow from the primary factors,
as in Theorem~\ref{thm:integral-forward-action}.
Write $H_{i,V}$ for the primed red bubble series and
$\kappa_V=(-1)^{c_0(V)+1}/2$.
For the zero color,
\[
 H_{0,V}(u)=\kappa_V
     \frac{t_{V,0}(u)}{u^{\varepsilon_0(V)+1}},\qquad
 H_{0,V}(u)H_{0,V}(-u)=-\frac{1}{4u^2}.
\]
If $y=x-j$ is the primed colored dot on $E_j$, then
\[
 \frac{H_{0,E_jV}(u)}{E_j(H_{0,V}(u))}
 =\begin{cases}
 ((u-y)/(u+y))^2,&j=0,\\
 -(u+y)/(u-y),&j=\pm1,
 \end{cases}
 \qquad
 \frac{H_{i,E_0V}(u)}{E_0(H_{i,V}(u))}
 =\frac{u+y}{u-y}\quad(i=\pm1).
\]
These are the primed bubble slides.
The nonzero bubble leading terms and products are the
half-integral ones.

For the red braid, the rotation argument of
Lemma~\ref{fw:braid-relation}, using the established
adjunction and cyclicity identities, reduces the $27$ triples
to six representatives:
\[
\begin{array}{c|c|l}
\text{representative}&\text{orbit size}&\text{calculation}\\ \hline
(-1,-1,-1)&6&\text{ordinary half-integral}\\
(1,-1,1)&2&\text{exceptional half-integral}\\
(0,-1,-1)&6&\text{ordinary primary crossing}\\
(-1,0,-1)&6&\text{ordinary primary crossing}\\
(0,0,1)&6&\text{Proposition~\ref{prop:forward-mixed-braid}
 and Lemma~\ref{compu-T}}\\
(0,0,0)&1&\text{Section~\ref{sec:nil-brauer-braid}}.
\end{array}
\]
The rotation $(i,j,k)\mapsto(j,k,-i)$ is given by
\cite[(3.43)--(3.44)]{BWW25}; it uses adjunction and cyclicity,
not an inverse mixed red crossing.
For the two one-zero representatives, use
\[
 Q'{}^\imath_{-1,0}(u,v)=u-v.
\]
Since no pair of distinct strands has opposite colors,
the ordinary calculation in Lemma~\ref{fw:braid-relation}
gives zero correction for $(0,-1,-1)$ and the identity
correction for $(-1,0,-1)$.
Their rotated triples include the three-distinct-color cases:
\[
\begin{array}{c|rrrrrr}
(i,j,k)&(0,1,-1)&(0,-1,1)&(1,-1,0)&(-1,1,0)&(1,0,-1)&(-1,0,1)\\ \hline
D(i,j,k)&C_R&-C_R&-C_L&C_L&0&0.
\end{array}
\]
Here $D$ is the left-hand side of \eqref{ibraid}, and $C_L,C_R$
are its last two cup--cap diagrams without scalar coefficients.
The signs are the divided differences of $Q'{}^\imath$;
\eqref{eq:mixed-three-colour-target} is not used.
For $(0,0,1)$, the coefficient identity above is
\eqref{appC:shift} with $\delta=1$, so
Proposition~\ref{prop:forward-mixed-braid} applies.
For $(0,0,0)$, use the same products, zero slides and
finite polynomial representatives.
This proves the forward defining relations.

For the reverse construction, the black adjunction, dot and
pitchfork relations follow from
Lemmas~\ref{lem:ir-global-adjunction},
\ref{lem:ir-global-affine}, and \ref{lem:ir-full-pitchfork},
using the half-integral nonzero-color formulas.
The nonzero principal-part calculation of
\cite[Proposition~7.20]{GRS} and the zero calculation of
Lemma~\ref{lem:red-black-series} give the full black bubble
$n_V/m_V$, including its constant coefficient.
The cap and quadratic relations follow from
\cite[Lemmas~7.21--7.22]{GRS} and
Lemmas~\ref{lem:ir-zero-full-cap},
\ref{lem:ir-zero-full-square}, and \ref{lem:ir-global-square}.
For the opposite pair $1,-1$, use the exceptional
half-integral case, not the nonexceptional integral argument
that inverts $4i^2-1$.

It remains to verify the black braid.
The dot, cap and quadratic relations reduce it to reversed
output colors as in \cite[Lemma~7.23]{GRS}, using
\cite[Corollary~2.7]{GRS} independently of the black braid.
The established adjunction, pitchfork and quadratic relations
allow the same rotations, so the six representatives suffice.
For $(-1,-1,-1)$, $(0,-1,-1)$ and $(-1,0,-1)$, use
\cite[Lemma~7.24]{GRS}, with the mixed coefficients above when
zero occurs. For $(1,-1,1)$, use the exceptional calculation
of \cite[Lemma~7.27]{GRS}, keeping the two black braid diagrams
separate until applying the assumed red relation.
The canceled factors are units: in the notation of
\cite[(7.33)]{GRS}, $h_1,h',h'_3$ all have residue $1$.
For $(0,0,1)$, repeat the calculation leading to
\eqref{eq:ir-mixed-braid-defect}, keeping the two leading
black braid diagrams. The remaining terms cancel by the
identities above and Lemma~\ref{compu-T}.
The common leading factor has residue $1/4=1$, hence is a unit.
For $(0,0,0)$, use \eqref{eq:ir-allzero-braid-defect}, whose
proof uses only the lower relations and keeps the two leading
black braid diagrams. Its exterior factors have residue $1$.
Thus the reverse verification does not assume the black braid.

Conjugation and faithfulness give these identities in the
original coefficient category.
Naturality follows by choosing common annihilating polynomials
on direct sums; evaluations are finite on finitely generated
objects. Biadjunction preserves finite generation and filtered
colimits, so the relations extend as in
Theorem~\ref{prop:ir-action}.
The shifted dots are nilpotent on finitely generated objects,
and finite color support is automatic.
Undoing the priming isomorphism gives the stated dot convention.
\end{proof}

\appendix
\section{Graphical reductions}
\label{app:graphical-certificates}

\label{app:braid-graphical-calculations}

This appendix records the reductions of the normalized
terms $\widetilde L_{ijk}$ and $\widetilde R_{ijk}$ used in the braid
proof. Fix a finitely generated object $V\in\cM_\lambda$ and put
$N=E_0V$.
We proceed by decreasing crossing number, retaining all lower-order
corrections. The one-crossing coefficients are compared in
Subsection~\ref{app:one-crossing-coefficients}; the zero-crossing
reductions in Subsection~\ref{app:zero-crossing-graphics} supply the
coefficient calculations in Appendix~B. The unnormalized source
diagrams are collected in Subsection~\ref{27terms-app}.

Let $\Pi$ project onto $E_0^3V$, and let $X_i$ be the red crossing on
positions $i,i+1$. Write $U_i$ and $D_i$ for the full three-term upper
and lower expressions in Lemma~\ref{lem:anothdefofcro}, respectively.
Use compatible polynomial representatives of sufficiently high order;
extensions of coefficients on $E_0V$ preserve the primary decomposition
of the remaining strand. The off-diagonal cancellations in the proof
of Lemma~\ref{nilRel_7} give
\[
{ \Pi D_i=\Pi D_i\Pi,\qquad U_i\Pi=\Pi U_i\Pi,\qquad
 \Pi U_i\Pi=\Pi D_i\Pi=X_i.}
\]
Indeed, on the projected input cap the normalization uses
$M(-x)=(1-2x)A(x)$, and on the projected output cup it uses
$W(-y)=(1+2y)B(y)$. These identities are needed only on the
corresponding zero primary components. Consequently,
\[
{ \Pi D_1U_2U_1\Pi=X_1X_2X_1,\qquad
 \Pi D_2U_1U_2\Pi=X_2X_1X_2.}
\]
{Thus the top crossing uses the lower expression (the $p$ term),
and the other two use the upper expression (the $q$ term), as in
the source diagrams. Expanding each product gives the twenty-seven
source terms.}
This removes the intermediate projections from the complete sums.

For the middle crossing, put $\Delta=x_1-x_2$ and choose a scalar
polynomial $R(t)\equiv(1+t)^{-1}\pmod{t^K}$ with $R(0)=1$.
Its parallel term uses the polynomial $H(t)=(1-R(t))/t$.
Take $K$ sufficiently large for all divided differences in the calculation.
Fix the upper and lower source terms, and combine the middle crossing
and parallel terms before rational simplification. On the left, the
lower term can produce only colors $(j,-j,0)$. If $j\ne0$, the middle
cup--cap term vanishes, and the final projection kills the exchanged
color sequence. On the remaining diagonal component,
Definition~\ref{C defn}(4) gives $s_{\mathrm{diag}}=\Delta^{-1}$, so
\[
 s_{\mathrm{diag}}R(\Delta)+H(\Delta)=\Delta^{-1}.
\]
Here $\Delta$ has nonzero eigenvalue $-j$; no inverse of $1+\Delta$
is used. For $j=0$, every surviving middle component is zero-colored,
where the chosen truncation applies. The right-hand calculation uses
$(0,j,-j)$ and is identical. Thus each pair
$L_{a1c}+L_{a2c}$, and likewise $R_{a1c}+R_{a2c}$, depends only on
the prescribed finite Taylor expansion of $R$ at zero.
The rational notation below is interpreted after this pairing, with
all dot-moving corrections retained.

We use the marked-position convention of the main proof:
$[D;f]$ denotes the representative $D$ decorated by $f$ at its marked
positions. Variables on disconnected cup and cap components remain
independent. Write $l_{a,b}=1+x_a-x_b$ and $l_{a,b}^{-1}=(1+x_a-x_b)^{-1}$.
We also retain the notation of Definition~\ref{ration111}:
$\ell_{\epsilon\eta}^{ij}=1+\epsilon x_i+\eta x_j$ and
$d_{\epsilon\eta}^{ij}=\epsilon x_i+\eta x_j$.
The diagram labels abbreviate
\[
q_{Xij}^{\epsilon\eta}=q_X(\epsilon x_i,\eta x_j),
\quad
p_{Xij}^{\epsilon\eta}=p_X(\epsilon x_i,\eta x_j),
\quad
F_{Xi}^{\epsilon}=F_X(\epsilon x_i),
\]
where $X\in\{V,N\}$, $\epsilon,\eta\in\{+,-\}$,
and $F\in\{A,B,M,W\}$.
A label without a numerical subscript uses the dot on that strand.
Unless the marked positions are explicitly numbered or specified
in the text, the arguments of a polynomial attached to a dashed
connector are ordered by increasing vertical coordinate and,
at equal height, from left to right. Thus a horizontal connector
is read from left to right. This convention applies to straight
and curved connectors alike. All contacts with solid strands,
including interior contacts, are counted. The ordering is local
to each connector; explicitly assigned variables are retained.

Transport through a crossing uses
Definition~\ref{C defn}(4), including its straight and cup--cap
corrections. Transport around a bend uses
Definition~\ref{C defn}(1),(5) and reverses the dot's sign.
Definition~\ref{C defn}(2),(7) reduces the undecorated diagrams to
the fifteen fixed representatives. The notation
$[\,\cdot\,]_{{u:-1}}$ denotes coefficient extraction in the
prescribed expansion at infinity; $u$ is a formal parameter.
In functional-calculus identities, rational decorations are first replaced
by polynomial representatives of sufficiently high order, as in the main
proof. Coefficient extraction is then applied with the resolvents
$(u-x)^{-1}$ expanded at infinity. The displayed rational divided
differences abbreviate the resulting regular expressions. A circled minus on a strand denotes $(u-x)^{-1}$, where $x$ is
the dot at that position; a circled plus denotes $(u+x)^{-1}$.
An undecorated point denotes the dot itself. A point accompanied by a
function label denotes that function of the dot. These pointwise
conventions do not identify variables on different components.

\subsection{The three-crossing term}\label{app:three-crossing-full}
\label{app:three-crossing}

\paragraph{Left-hand reduction.}
\begin{align}
&{}\tilde L_{111} = \begin{tikzpicture}[
baseline=(current bounding box.center),scale=0.8,>=stealth,
  line cap=round,line join=round,font={\fontsize{10pt}{12pt}\selectfont}
]
  \draw[grsw strand]
    (0.45,0.6)
      -- (-0.45,-0.6);
  \draw[grsw strand]
    (0.45,-0.6)
      -- (-0.45,0.6);
  \draw[grsw strand]
    (0,-0.6)
      .. controls (0,-0.31) and (-0.45,-0.29) .. (-0.45,0);
  \draw[grsw strand]
    (-0.45,0)
      .. controls (-0.45,0.29) and (0,0.31) .. (0,0.6);
  \draw[grsw morphism,dash pattern=on 0.978678pt off 0.978678pt,line width=0.350000pt]
    (-0.06,0.43)
      -- (0.33,0.43);
  \draw[grsw morphism,dash pattern=on 0.993018pt off 0.993018pt,line width=0.350000pt]
    (-0.07,-0.09)
      -- (0.07,-0.09);
  \draw[grsw morphism,dash pattern=on 0.869375pt off 0.869375pt,line width=0.350000pt]
    (-0.36,-0.48)
      -- (0.36,-0.48);
  \node[anchor=west,inner sep=1pt,outer sep=0pt] at (0.61,0.47) {$\scriptstyle l_{2,1}$};
  \node[anchor=east,inner sep=1pt,outer sep=0pt] at (-0.61,-0.05) {$\scriptstyle r(x_1,x_2)$};
  \node[anchor=west,inner sep=1pt,outer sep=0pt] at (0.61,-0.44) {$\scriptstyle l_{1,3}l_{2,3}$};
\end{tikzpicture} = \begin{tikzpicture}[
baseline=(current bounding box.center),scale=0.8,>=stealth,
  line cap=round,line join=round,font={\fontsize{10pt}{12pt}\selectfont}
]
  \draw[grsw strand]
    (0.45,0.6)
      -- (-0.45,-0.6);
  \draw[grsw strand]
    (0.45,-0.6)
      -- (-0.45,0.6);
  \draw[grsw strand]
    (0,-0.6)
      .. controls (0,-0.31) and (-0.45,-0.29) .. (-0.45,0);
  \draw[grsw strand]
    (-0.45,0)
      .. controls (-0.45,0.29) and (0,0.31) .. (0,0.6);
  \draw[grsw morphism,dash pattern=on 0.978678pt off 0.978678pt,line width=0.350000pt]
    (-0.06,0.43)
      -- (0.33,0.43);
  \draw[grsw morphism,dash pattern=on 0.894003pt off 0.894003pt,line width=0.350000pt]
    (0,-0.57)
      -- (0.43,-0.57);
  \node[anchor=west,inner sep=1pt,outer sep=0pt] at (0.61,0.47) {$\scriptstyle l_{2,1}$};
  \node[anchor=west,inner sep=1pt,outer sep=0pt] at (0.61,-0.53) {$\scriptstyle l_{1,2}$};
\end{tikzpicture} + \begin{tikzpicture}[
baseline=(current bounding box.center),scale=0.48,>=stealth,
  line cap=round,line join=round,font={\fontsize{10pt}{12pt}\selectfont}
]
  \draw[grsw strand]
    (-0.9,1)
      -- (1.1,-1);
  \draw[grsw strand]
    (0.1,1)
      -- (-0.9,-1);
  \draw[grsw strand]
    (-0.1,-1)
      -- (0.9,1);
  \draw[grsw morphism,dash pattern=on 0.910501pt off 0.910501pt,line width=0.350000pt]
    (-0.01,0.78)
      -- (0.79,0.78);
  \node[anchor=west,inner sep=1pt,outer sep=0pt] at (1.32,0.82) {$\scriptstyle l_{2,1}$};
\end{tikzpicture} + \widetilde L_{1,7}\notag\\[4pt]
&{}= \begin{tikzpicture}[
baseline=(current bounding box.center),scale=0.8,>=stealth,
  line cap=round,line join=round,font={\fontsize{10pt}{12pt}\selectfont}
]
  \draw[grsw strand]
    (0.45,0.6)
      -- (-0.45,-0.6);
  \draw[grsw strand]
    (0.45,-0.6)
      -- (-0.45,0.6);
  \draw[grsw strand]
    (0,-0.6)
      .. controls (0,-0.31) and (-0.45,-0.29) .. (-0.45,0);
  \draw[grsw strand]
    (-0.45,0)
      .. controls (-0.45,0.29) and (0,0.31) .. (0,0.6);
  \draw[grsw morphism,dash pattern=on 0.978678pt off 0.978678pt,line width=0.350000pt]
    (-0.06,0.43)
      -- (0.33,0.43);
  \draw[grsw morphism,dash pattern=on 0.783626pt off 0.783626pt,line width=0.350000pt]
    (-0.38,0.18)
      -- (-0.14,0.18);
  \node[anchor=west,inner sep=1pt,outer sep=0pt] at (0.61,0.47) {$\scriptstyle l_{2,1}$};
  \node[anchor=east,inner sep=1pt,outer sep=0pt] at (-0.61,0.22) {$\scriptstyle l_{1,2}$};
\end{tikzpicture} + \begin{tikzpicture}[
baseline=(current bounding box.center),scale=0.8,>=stealth,
  line cap=round,line join=round,font={\fontsize{10pt}{12pt}\selectfont}
]
  \draw[grsw strand]
    (-0.45,0.6)
      -- (0.45,-0.6);
  \draw[grsw strand]
    (0,-0.6)
      .. controls (0,-0.42) and (-0.07,-0.2) .. (-0.25,-0.2);
  \draw[grsw strand]
    (-0.25,-0.2)
      .. controls (-0.43,-0.2) and (-0.45,-0.42) .. (-0.45,-0.6);
  \draw[grsw strand]
    (0.45,0.6)
      .. controls (0.45,0.42) and (0.38,0.2) .. (0.2,0.2);
  \draw[grsw strand]
    (0.2,0.2)
      .. controls (0.02,0.2) and (0,0.42) .. (0,0.6);
  \draw[grsw morphism,dash pattern=on 0.839478pt off 0.839478pt,line width=0.350000pt]
    (0.02,0.4)
      -- (0.42,0.4);
  \node[anchor=west,inner sep=1pt,outer sep=0pt] at (0.61,0.44) {$\scriptstyle l_{2,1}$};
\end{tikzpicture} + \begin{tikzpicture}[
baseline=(current bounding box.center),scale=0.8,>=stealth,
  line cap=round,line join=round,font={\fontsize{10pt}{12pt}\selectfont}
]
  \draw[grsw strand]
    (0.45,0.6)
      -- (0.45,-0.6);
  \draw[grsw strand]
    (-0.45,0.6)
      -- (-0.45,-0.6);
  \draw[grsw strand]
    (0,-0.6)
      -- (0,0.6);
  \draw[grsw morphism,dash pattern=on 0.931184pt off 0.931184pt,line width=0.350000pt]
    (0,-0.04)
      -- (0.45,-0.04);
  \node[anchor=west,inner sep=1pt,outer sep=0pt] at (0.61,0) {$\scriptstyle l_{2,1}$};
\end{tikzpicture} - \begin{tikzpicture}[
baseline=(current bounding box.center),scale=0.8,>=stealth,
  line cap=round,line join=round,font={\fontsize{10pt}{12pt}\selectfont}
]
  \draw[grsw strand]
    (0.3,0.6)
      .. controls (0.3,0.37) and (0.23,0.1) .. (0,0.1);
  \draw[grsw strand]
    (0,0.1)
      .. controls (-0.23,0.1) and (-0.3,0.37) .. (-0.3,0.6);
  \draw[grsw strand]
    (0.3,-0.6)
      .. controls (0.3,-0.37) and (0.23,-0.1) .. (0,-0.1);
  \draw[grsw strand]
    (0,-0.1)
      .. controls (-0.23,-0.1) and (-0.3,-0.37) .. (-0.3,-0.6);
  \draw[grsw strand]
    (0,-0.6)
      -- (0,0.6);
  \draw[grsw morphism,dash pattern=on 0.945884pt off 0.945884pt,line width=0.350000pt]
    (0,0.46)
      -- (0.29,0.46);
  \node[anchor=west,inner sep=1pt,outer sep=0pt] at (0.46,0.5) {$\scriptstyle l_{2,1}$};
\end{tikzpicture} + \widetilde L_{1,7}.\notag\\[4pt]
&{}=\widetilde L_{1,1} +2\widetilde L_{1,2} -2\widetilde L_{1,3} +\widetilde L_{1,4} +\widetilde L_{1,5} -\widetilde L_{1,6} +\widetilde L_{1,7}
\end{align}
where
\[
\begin{alignedat}[c]{3}
\widetilde L_{1,1}
  &=[s_2s_1s_2;l_{2,1}l_{3,2}],
&\hspace{0.4em}\widetilde L_{1,2}
  &=[s_2s_1;l_{3,2}],
&\hspace{0.4em}\widetilde L_{1,3}
  &=[e_1e_2;l_{3,2}],\\[6pt]
\widetilde L_{1,4}
  &=[e_2e_1;l_{3,2}],
&\hspace{0.4em}\widetilde L_{1,5}
  &=[1;l_{3,2}],
&\hspace{0.4em}\widetilde L_{1,6}
  &=[s_2e_1s_2; l_{3,2}],
\end{alignedat}
\hspace{0.6em}
\widetilde L_{1,7}
 =\left[\begin{tikzpicture}[
baseline=(current bounding box.center),scale=0.8,>=stealth,
  line cap=round,line join=round,font={\fontsize{10pt}{12pt}\selectfont}
]
  \draw[grsw strand]
    (-0.3,0.6)
      -- (0.6,-0.6);
  \draw[grsw strand]
    (-0.05,-0.6)
      .. controls (-0.05,-0.46) and (-0.16,-0.35) .. (-0.3,-0.35);
  \draw[grsw strand]
    (-0.3,-0.35)
      .. controls (-0.43,-0.35) and (-0.5,-0.47) .. (-0.5,-0.6);
  \draw[grsw strand]
    (0.5,0.6)
      .. controls (0.5,0.18) and (0.27,-0.25) .. (-0.15,-0.25);
  \draw[grsw strand]
    (-0.15,-0.25)
      .. controls (-0.28,-0.25) and (-0.35,-0.13) .. (-0.35,0);
  \draw[grsw strand]
    (-0.35,0)
      .. controls (-0.35,0.28) and (0.05,0.32) .. (0.05,0.6);
  \draw[grsw morphism,dash pattern=on 0.909027pt off 0.909027pt,line width=0.350000pt]
    (0.01,0.46)
      -- (0.49,0.46);
  \draw[grsw morphism,dash pattern=on 0.918032pt off 0.918032pt,line width=0.350000pt]
    (-0.48,-0.49)
      -- (0.52,-0.49);
  \node[anchor=west,inner sep=1pt,outer sep=0pt] at (0.73,0.5) {$\scriptstyle l_{2,1}$};
  \node[anchor=center,inner sep=0pt,outer sep=0pt] at (0.34,-0.26) {$\scriptscriptstyle\bullet$};
  \node[anchor=west,inner sep=1pt,outer sep=0pt] at (0.73,-0.26) {$\scriptstyle r(u,x)$};
  \node[anchor=east,inner sep=1pt,outer sep=0pt] at (-0.63,-0.45) {$\scriptstyle l_{1,3}l_{2,3}$};
  \node[cplus,scale=0.7,anchor=center,inner sep=0pt,outer sep=0pt] at (0.01,-0.23) {$-$};
  \node[cplus,scale=0.7,anchor=center,inner sep=0pt,outer sep=0pt] at (-0.49,-0.52) {$-$};
\end{tikzpicture}\right]_{{u:-1}}.\]

\paragraph{Right-hand reduction.}
\begin{align}
&{}\tilde R_{111} = \begin{tikzpicture}[
baseline=(current bounding box.center),scale=0.8,>={Stealth[reversed]},
  line cap=round,line join=round,font={\fontsize{10pt}{12pt}\selectfont}
]
  \draw[grsw strand]
    (-0.45,0.6)
      -- (0.45,-0.6);
  \draw[grsw strand]
    (-0.45,-0.6)
      -- (0.45,0.6);
  \draw[grsw strand]
    (0,-0.6)
      .. controls (0,-0.31) and (0.45,-0.29) .. (0.45,0);
  \draw[grsw strand]
    (0.45,0)
      .. controls (0.45,0.29) and (0,0.31) .. (0,0.6);
  \draw[grsw morphism,dash pattern=on 0.978678pt off 0.978678pt,line width=0.350000pt] (-0.33,0.43) -- (0.06,0.43);
  \draw[grsw morphism,dash pattern=on 0.993018pt off 0.993018pt,line width=0.350000pt] (-0.07,-0.09) -- (0.07,-0.09);
  \draw[grsw morphism,dash pattern=on 0.869375pt off 0.869375pt,line width=0.350000pt] (-0.36,-0.48) -- (0.36,-0.48);
  \node[anchor=east,inner sep=1pt,outer sep=0pt] at (-0.61,0.47) {$\scriptstyle l_{2,1}$};
  \node[anchor=east,inner sep=1pt,outer sep=0pt] at (-0.61,-0.05) {$\scriptstyle r(x_1,x_2)$};
  \node[anchor=west,inner sep=1pt,outer sep=0pt] at (0.61,-0.44) {$\scriptstyle l_{1,3}l_{1,2}$};
\end{tikzpicture} = \begin{tikzpicture}[
baseline=(current bounding box.center),scale=0.8,>=stealth,
  line cap=round,line join=round,font={\fontsize{10pt}{12pt}\selectfont}
]
  \draw[grsw strand]
    (-0.45,0.6)
      -- (0.45,-0.6);
  \draw[grsw strand]
    (-0.45,-0.6)
      -- (0.45,0.6);
  \draw[grsw strand]
    (0,-0.6)
      .. controls (0,-0.31) and (0.45,-0.29) .. (0.45,0);
  \draw[grsw strand]
    (0.45,0)
      .. controls (0.45,0.29) and (0,0.31) .. (0,0.6);
  \draw[grsw morphism,dash pattern=on 0.978678pt off 0.978678pt,line width=0.350000pt] (-0.33,0.43) -- (0.06,0.43);
  \draw[grsw morphism,dash pattern=on 0.894003pt off 0.894003pt,line width=0.350000pt] (-0.43,-0.57) -- (0,-0.57);
  \node[anchor=east,inner sep=1pt,outer sep=0pt] at (-0.61,0.47) {$\scriptstyle l_{2,1}$};
  \node[anchor=east,inner sep=1pt,outer sep=0pt] at (-0.61,-0.53) {$\scriptstyle l_{1,2}$};
\end{tikzpicture} + \begin{tikzpicture}[
baseline=(current bounding box.center),scale=0.48,>=stealth,
  line cap=round,line join=round,font={\fontsize{10pt}{12pt}\selectfont}
]
  \draw[grsw strand]
    (0.9,1)
      -- (-1.1,-1);
  \draw[grsw strand]
    (-0.1,1)
      -- (0.9,-1);
  \draw[grsw strand]
    (0.1,-1)
      -- (-0.9,1);
  \draw[grsw morphism,dash pattern=on 0.910501pt off 0.910501pt,line width=0.350000pt] (-0.79,0.78) -- (0.01,0.78);
  \node[anchor=east,inner sep=1pt,outer sep=0pt] at (-1.32,0.82) {$\scriptstyle l_{2,1}$};
\end{tikzpicture} - \widetilde R_{1,7}\notag\\[4pt]
&{}= \begin{tikzpicture}[
baseline=(current bounding box.center),scale=0.8,>=stealth,
  line cap=round,line join=round,font={\fontsize{10pt}{12pt}\selectfont}
]
  \draw[grsw strand]
    (-0.45,0.6)
      -- (0.45,-0.6);
  \draw[grsw strand]
    (-0.45,-0.6)
      -- (0.45,0.6);
  \draw[grsw strand]
    (0,-0.6)
      .. controls (0,-0.31) and (0.45,-0.29) .. (0.45,0);
  \draw[grsw strand]
    (0.45,0)
      .. controls (0.45,0.29) and (0,0.31) .. (0,0.6);
  \draw[grsw morphism,dash pattern=on 0.978678pt off 0.978678pt,line width=0.350000pt] (-0.33,0.43) -- (0.06,0.43);
  \draw[grsw morphism,dash pattern=on 0.783626pt off 0.783626pt,line width=0.350000pt] (0.14,0.18) -- (0.38,0.18);
  \node[anchor=east,inner sep=1pt,outer sep=0pt] at (-0.61,0.47) {$\scriptstyle l_{2,1}$};
  \node[anchor=west,inner sep=1pt,outer sep=0pt] at (0.61,0.22) {$\scriptstyle l_{1,2}$};
\end{tikzpicture} + \begin{tikzpicture}[
baseline=(current bounding box.center),scale=0.8,>=stealth,
  line cap=round,line join=round,font={\fontsize{10pt}{12pt}\selectfont}
]
  \draw[grsw strand]
    (0.45,0.6)
      -- (-0.45,-0.6);
  \draw[grsw strand]
    (0,-0.6)
      .. controls (0,-0.42) and (0.07,-0.2) .. (0.25,-0.2);
  \draw[grsw strand]
    (0.25,-0.2)
      .. controls (0.43,-0.2) and (0.45,-0.42) .. (0.45,-0.6);
  \draw[grsw strand]
    (-0.45,0.6)
      .. controls (-0.45,0.42) and (-0.38,0.2) .. (-0.2,0.2);
  \draw[grsw strand]
    (-0.2,0.2)
      .. controls (-0.02,0.2) and (0,0.42) .. (0,0.6);
  \draw[grsw morphism,dash pattern=on 0.839478pt off 0.839478pt,line width=0.350000pt] (-0.42,0.4) -- (-0.02,0.4);
  \node[anchor=east,inner sep=1pt,outer sep=0pt] at (-0.61,0.44) {$\scriptstyle l_{2,1}$};
\end{tikzpicture} + \begin{tikzpicture}[
baseline=(current bounding box.center),scale=0.8,>=stealth,
  line cap=round,line join=round,font={\fontsize{10pt}{12pt}\selectfont}
]
  \draw[grsw strand]
    (-0.45,0.6)
      -- (-0.45,-0.6);
  \draw[grsw strand]
    (0.45,0.6)
      -- (0.45,-0.6);
  \draw[grsw strand]
    (0,-0.6)
      -- (0,0.6);
  \draw[grsw morphism,dash pattern=on 0.931184pt off 0.931184pt,line width=0.350000pt] (-0.45,-0.04) -- (0,-0.04);
  \node[anchor=east,inner sep=1pt,outer sep=0pt] at (-0.61,0) {$\scriptstyle l_{2,1}$};
\end{tikzpicture} - \begin{tikzpicture}[
baseline=(current bounding box.center),scale=0.8,>=stealth,
  line cap=round,line join=round,font={\fontsize{10pt}{12pt}\selectfont}
]
  \draw[grsw strand]
    (-0.3,0.6)
      .. controls (-0.3,0.37) and (-0.23,0.1) .. (0,0.1);
  \draw[grsw strand]
    (0,0.1)
      .. controls (0.23,0.1) and (0.3,0.37) .. (0.3,0.6);
  \draw[grsw strand]
    (-0.3,-0.6)
      .. controls (-0.3,-0.37) and (-0.23,-0.1) .. (0,-0.1);
  \draw[grsw strand]
    (0,-0.1)
      .. controls (0.23,-0.1) and (0.3,-0.37) .. (0.3,-0.6);
  \draw[grsw strand]
    (0,-0.6)
      -- (0,0.6);
  \draw[grsw morphism,dash pattern=on 0.945884pt off 0.945884pt,line width=0.350000pt] (-0.29,0.46) -- (0,0.46);
  \node[anchor=east,inner sep=1pt,outer sep=0pt] at (-0.46,0.5) {$\scriptstyle l_{2,1}$};
\end{tikzpicture} - \widetilde R_{1,7}.\notag\\[4pt]
&{}= \widetilde R_{1,1} +2\widetilde R_{1,2} -2\widetilde R_{1,3} +\widetilde R_{1,4} +\widetilde R_{1,5} -\widetilde R_{1,6} -\widetilde R_{1,7} ,
\end{align}
 where
\[
\begin{alignedat}[c]{3}
\widetilde R_{1,1}
  &=[s_2s_1s_2;l_{2,1}l_{3,2}],
&\hspace{0.4em}\widetilde R_{1,2}
  &=[s_1s_2;l_{2,1}],
&\hspace{0.4em}\widetilde R_{1,3}
  &=[e_2e_1;l_{2,1}],\\[6pt]
\widetilde R_{1,4}
  &=[e_1e_2;l_{2,1}],
&\hspace{0.4em}\widetilde R_{1,5}
  &=[1;l_{2,1}],
&\hspace{0.4em}\widetilde R_{1,6}
  &=[s_2e_1s_2;l_{2,1}],
\end{alignedat}
\hspace{0.6em}
\widetilde R_{1,7} =\left[
\right]_{{u:-1}}.
\]

\subsection{The two-crossing terms}\label{app:two-crossing-full}
\label{app:two-crossing-through}
\subsubsection{The $s_2s_1$-term}
The $s_2s_1$-coefficient receives contributions from $\widetilde L_{211}$, $\widetilde R_{112}$, and $2\widetilde L_{1,2}$. The factor $2$ is the multiplicity in the three-crossing expansion. We have
\[\begin{aligned}
\tilde L_{211} & = %
\right]_{{u:-1}}
\]

\subsubsection{The $s_1s_2$-term}
The $s_1s_2$-coefficient receives contributions from $\widetilde L_{112}$, $\widetilde R_{211}$, and $2\widetilde R_{1,2}$. Again, the factor $2$ comes from the three-crossing expansion.
We have:
\[\begin{aligned}
\tilde L_{112} & = %
\right]_{{u:-1}}
\]

\subsubsection{The cup--cap term}
\label{app:two-crossing-cupcap}

The source term on the left is
\[
\tilde L_{131} = %
\]
and the reflected source term is
\[
\tilde R_{131} = %
\]
These are the source diagrams used in the comparison of
\(s_1e_2s_1=s_2e_1s_2\) in the main proof.

\subsection{One-crossing source terms}
\label{app:one-crossing-sources}
There are twelve one-crossing summands on each side. We group the
left- and right-hand reductions by their source indices.
{The one-crossing corrections from higher-crossing terms are
recorded in Subsection~\ref{app:mixed-one-crossing-graphics}.}

\subsubsection{Direct reductions}

\noindent\textit{Left-hand terms.}\par\nopagebreak[4]

\[
\begin{alignedat}[c]{2}
\tilde L_{312} &=
%
.
\end{alignedat}
\]

\subsubsection{The sources $321$, $231$, and $221$}

\noindent\textit{Left-hand terms.}\par\nopagebreak[4]
\begin{align*}
\tilde L_{321} & = %
\end{align*}

\begin{proof} By definition, we have
 \[
\begin{aligned}
\tilde L_{321} & =%
\\
\end{aligned}\]
{For the sources $321$ and $231$, move $l_{1,2}$ to the marked
upper segments; for $221$, move $l_{1,3}$. Retain the corrections
from Definition~\ref{C defn}(4), then apply ~\eqref{AB relations1}.} \end{proof}

\noindent\textit{Right-hand terms.}\par\nopagebreak[4]
\begin{align*}
\tilde R_{321} & = %
\end{align*}

\begin{proof}
By definition, we have
\begin{align*}
&{}\tilde R_{321} = %
\end{align*} {For the sources $321$ and $231$, move $l_{1,2}$ to the marked
upper segments; for $221$, move $l_{1,3}$. Retain the corrections
from Definition~\ref{C defn}(4), then apply ~\eqref{AB relations1}.} \end{proof}

\subsubsection{The sources $213$, $132$, and $123$}

\noindent\textit{Left-hand terms.}\par\nopagebreak[4]
 \begin{align*}
\tilde L_{213} & = -%
\end{align*}

\begin{proof}
    By definition, we have $\tilde L_{213} =%
$. {Move $l_{2,1}^{-1}$ to the marked upper segments, retaining the
corrections from Definition~\ref{C defn}(4).  ~\eqref{AB relations1}
gives the $213$ reduction; the $132$ and $123$ reductions follow
from the definitions.}
\end{proof}

\noindent\textit{Right-hand terms.}\par\nopagebreak[4]
 \begin{align*}
\tilde R_{213} & = -%
\end{align*}
\begin{proof}
    By definition, we have $$\tilde R_{213} = %
.$$
{Move $l_{2,1}^{-1}$ to the marked upper segments, retaining the
corrections from Definition~\ref{C defn}(4). ~\eqref{AB relations1}
gives the $213$ reduction; the $132$ and $123$ reductions follow
from the definitions.}
\end{proof}

\subsection{One-crossing corrections from higher-crossing terms}
\label{app:mixed-one-crossing-graphics}

In this subsection, $z_i$ denotes a decorated diagram without coefficient
extraction. We write $y_i(u)$ for a diagram depending on $u$ and set
\[
  Y_i=[y_i(u)]_{{u:-1}}.
\]
The labels are local to this subsection. External signs and numerical
coefficients are written explicitly in the identities below.

\subsubsection*{Catalogue of decorated diagrams}
\par\medskip
\begingroup

\providecommand{\appAZYentry}[3]{%
	\hbox{$\thickmuskip=3mu\displaystyle #1_{#2}=#3$}%
}

\setlength{\tabcolsep}{0pt}
\renewcommand{\arraystretch}{1}
\setlength{\LTleft}{\fill}
\setlength{\LTright}{\fill}
\setlength{\LTpre}{0pt}
\setlength{\LTpost}{0pt}

\begin{longtable}{@{}l@{\hspace{3pt}}l@{\hspace{3pt}}l@{\hspace{3pt}}l@{}}

	\appAZYentry{z}{1}{
		\begin{tikzpicture}[baseline=(current bounding box.center),
			scale=0.8,line cap=round,line join=round,
			font={\fontsize{10pt}{12pt}\selectfont}]
			\begin{scope}[xscale=-1,yscale=-1]
				\draw[-,thick] (-0.4,0.6) to[out=270, in=180] (0.2,0.1);
				\draw[-,thick] (0.2,0.1) to[out = 0, in = 270] (0.4,0.6);
				\draw[-,thick] (-0.4,-0.6) to (0,0.6);
				\draw[-,thick] (0,-0.6) to[out=90, in=180] (0.2,-0.3);
				\draw[-,thick] (0.2,-0.3) to[out = 0, in = 90] (0.4,-0.6);
			\end{scope}
		\end{tikzpicture}
	}%
	&
	\appAZYentry{z}{2}{
		\begin{tikzpicture}[baseline=(current bounding box.center),
			scale=0.8,line cap=round,line join=round,
			font={\fontsize{10pt}{12pt}\selectfont}]
			\begin{scope}[xscale=1,yscale=-1]
				\draw[-,thick] (-0.4,0.6) to[out=270, in=180] (0.2,0.1);
				\draw[-,thick] (0.2,0.1) to[out = 0, in = 270] (0.4,0.6);
				\draw[-,thick] (-0.4,-0.6) to (0,0.6);
				\draw[-,thick] (0,-0.6) to[out=90, in=180] (0.2,-0.3);
				\draw[-,thick] (0.2,-0.3) to[out = 0, in = 90] (0.4,-0.6);
			\end{scope}
		\end{tikzpicture}
	}%
	&
	\appAZYentry{z}{3}{
		\begin{tikzpicture}[baseline=(current bounding box.center),
			scale=0.8,line cap=round,line join=round,
			font={\fontsize{10pt}{12pt}\selectfont}]
			\begin{scope}[xscale=-1,yscale=1]
				\draw[-,thick] (0.6,-0.6) to[out=90, in=0]  (0.3,-0.1);
				\draw[-,thick] (0.3,-0.1) to[out = 180, in = 90](0,-0.6);
				\draw[-,thick] (0.6,0.6) to[out=-90, in=0]  (0.3,0.1);
				\draw[-,thick] (0.3,0.1) to[out = -180, in = -90] (0,0.6);
				\draw[-,thick] (0.3,0.6) to   (0.3,-0.6);
			\end{scope}
			\draw[densely dashed,line width=.35pt,-,line width=0.350000pt,rounded corners=3pt] (-0.47,0.3)--(-0.62,0.3)--(-0.62,0.03)--(-0.26,0.03);
			\node[anchor=east,inner sep=1pt,outer sep=0pt] at (-0.76,0.3) {$\scriptstyle \ell_{--}^{12}$};
		\end{tikzpicture}
	}%
	&
	\appAZYentry{z}{4}{
		\begin{tikzpicture}[baseline=(current bounding box.center),
			scale=0.8,line cap=round,line join=round,
			font={\fontsize{10pt}{12pt}\selectfont}]
			\begin{scope}[xscale=-1,yscale=1]
				\draw[-,thick] (0.6,-0.6) to[out=90, in=0]  (0.3,-0.1);
				\draw[-,thick] (0.3,-0.1) to[out = 180, in = 90](0,-0.6);
				\draw[-,thick] (0.6,0.6) to[out=-90, in=0]  (0.3,0.1);
				\draw[-,thick] (0.3,0.1) to[out = -180, in = -90] (0,0.6);
				\draw[-,thick] (0.3,0.6) to   (0.3,-0.6);
			\end{scope}
			\draw[densely dashed,line width=.35pt,-,line width=0.350000pt,rounded corners=3pt] (-0.47,0.3)--(-0.62,0.3)--(-0.62,0.03)--(-0.26,0.03);
			\node[anchor=east,inner sep=1pt,outer sep=0pt] at (-0.76,0.3) {$\scriptstyle \ell_{--}^{12}$};
			\draw[densely dashed,line width=.35pt,-,line width=0.350000pt,rounded corners=3pt] (-0.47,-0.3)--(-0.62,-0.3)--(-0.62,-0.03)--(-0.26,-0.03);
			\node[anchor=east,inner sep=1pt,outer sep=0pt] at (-0.76,-0.3) {$\scriptstyle \ell_{++}^{12}$};
			\draw[densely dashed,line width=.35pt,-,line width=0.350000pt] (-0.02,-0.4) -- (-0.02,0.4);
			\node[anchor=west,inner sep=1pt,outer sep=0pt] at (0.16,0) {$\scriptstyle q_{V21}^{++}M_{V1}^{-}$};
		\end{tikzpicture}
	}\\
\end{longtable}

\begin{longtable}{@{}l@{\hspace{3pt}}l@{\hspace{3pt}}l@{}}

	\appAZYentry{y}{31}{
		\begin{tikzpicture}[baseline=(current bounding box.center),
			scale=0.9,line cap=round,line join=round,
			font={\fontsize{10pt}{12pt}\selectfont}]
			\begin{scope}[xscale=-1,yscale=-1]
				\draw[-,thick] (0,0.6) to (0.45,-0.6);
				\draw[-,thick] (-0.45,0.6) to[out=-90,in=180] (0,0.2);
				\draw[-,thick] (0,0.2) to[out=0,in=-90] (0.45,0.6);
				\draw[-,thick] (-0.45,-0.6) to[out=90,in=180] (-0.2,-0.2);
				\draw[-,thick] (-0.2,-0.2) to[out=0,in=90] (0,-0.6);
			\end{scope}
			\draw[densely dashed,line width=.35pt,-,line width=0.350000pt] (-0.43,0.54) -- (0.45,0.54);
			\node[anchor=south,inner sep=1pt,outer sep=0pt] at (0,0.74) {$\scriptstyle \frac{\ell_{-+}^{12}p_V(-x_2,-u)\, W_{V2}^{+}}{u+x_1}$};
			\draw[densely dashed,line width=.35pt,-,line width=0.350000pt] (-0.45,-0.54) -- (0.45,-0.54);
			\node[anchor=north,inner sep=1pt,outer sep=0pt] at (0,-0.74) {$\scriptstyle \frac{\ell_{+-}^{12}}{u-x_1}$};
		\end{tikzpicture}
	}%
	&
	\appAZYentry{y}{9}{
		\begin{tikzpicture}[baseline=(current bounding box.center),
			scale=0.8,line cap=round,line join=round,
			font={\fontsize{10pt}{12pt}\selectfont}]
			\begin{scope}[xscale=-1,yscale=1]
				\draw[-,thick] (0.6,-0.6) to[out=90, in=0]  (0.3,-0.1);
				\draw[-,thick] (0.3,-0.1) to[out = 180, in = 90](0,-0.6);
				\draw[-,thick] (0.6,0.6) to[out=-90, in=0]  (0.3,0.1);
				\draw[-,thick] (0.3,0.1) to[out = -180, in = -90] (0,0.6);
				\draw[-,thick] (0.3,0.6) to   (0.3,-0.6);
			\end{scope}
			\draw[densely dashed,line width=.35pt,-,line width=0.350000pt,rounded corners=3pt] (-0.47,0.3)--(-0.62,0.3)--(-0.62,0.03)--(-0.26,0.03);
			\node[anchor=east,inner sep=1pt,outer sep=0pt] at (-0.76,0.3) {$\scriptstyle \ell_{--}^{12}$};
			\draw[densely dashed,line width=.35pt,-,line width=0.350000pt,rounded corners=3pt] (-0.47,-0.3)--(-0.62,-0.3)--(-0.62,-0.03)--(-0.26,-0.03);
			\node[anchor=east,inner sep=1pt,outer sep=0pt] at (-0.76,-0.3) {$\scriptstyle \ell_{++}^{12}$};
			\draw[densely dashed,line width=.35pt,-,line width=0.350000pt] (-0.02,-0.4) -- (-0.02,0.4);
			\node[anchor=west,inner sep=1pt,outer sep=0pt] at (0.16,0) {$\scriptstyle \frac{q_V(x_2,u)M_V(-u)}{u-x_1}$};
		\end{tikzpicture}
	}%
	&
	\appAZYentry{y}{10}{
		\begin{tikzpicture}[baseline=(current bounding box.center),
			scale=0.8,line cap=round,line join=round,
			font={\fontsize{10pt}{12pt}\selectfont}]
			\begin{scope}[xscale=-1,yscale=1]
				\draw[-,thick] (0.6,-0.6) to[out=90, in=0]  (0.3,-0.1);
				\draw[-,thick] (0.3,-0.1) to[out = 180, in = 90](0,-0.6);
				\draw[-,thick] (0.6,0.6) to[out=-90, in=0]  (0.3,0.1);
				\draw[-,thick] (0.3,0.1) to[out = -180, in = -90] (0,0.6);
				\draw[-,thick] (0.3,0.6) to   (0.3,-0.6);
			\end{scope}
			\draw[densely dashed,line width=.35pt,-,line width=0.350000pt,rounded corners=3pt] (-0.47,0.3)--(-0.62,0.3)--(-0.62,0.03)--(-0.26,0.03);
			\node[anchor=east,inner sep=1pt,outer sep=0pt] at (-0.76,0.3) {$\scriptstyle \ell_{--}^{12}$};
			\draw[densely dashed,line width=.35pt,-,line width=0.350000pt,rounded corners=3pt] (-0.47,-0.3)--(-0.62,-0.3)--(-0.62,-0.03)--(-0.26,-0.03);
			\node[anchor=east,inner sep=1pt,outer sep=0pt] at (-0.76,-0.3) {$\scriptstyle \frac{\ell_{++}^{12}}{u+x_1}$};
				\fill (-0.01,0.5) circle[radius=1.25pt];
			\node[anchor=west,inner sep=1pt,outer sep=0pt] at (0.16,0.52) {$\scriptstyle q_V(x,u)M_V(-u)$};
		\end{tikzpicture}
	}\\
\end{longtable}
\endgroup
\par\medskip

\begingroup
\allowdisplaybreaks[4]

\definecolor{diagramred}{RGB}{174,48,48}

\subsubsection*{The \(e_1s_2,e_2s_1\) family}
With the notation above, the required local reductions are
\begin{align*}
- %
  & = Y_{28} = Y_{29} - Y_{30} + Y_{31}
\end{align*}
Move each marked decoration through the crossing by
Definition~\ref{C defn}(4), then use the cup--cap relations in
Definition~\ref{C defn}(1),(5). For bracketed terms, take the
$u^{-1}$-coefficient after the graphical reduction.

\subsubsection*{The \(s_1e_2,s_2e_1\) family}
The corresponding reductions are
\begin{align*}
%
\right]_{{u:-1}}
  & = Y_{44} =  Y_{45} -Y_{46} + Y_{47}
\end{align*}

\endgroup

\subsection{Comparison of one-crossing coefficients}
\label{app:one-crossing-coefficients}\label{app:one-crossing-comparison}
For each representative $d$, collect its contributions from Subsections~\ref{app:one-crossing-sources} and \ref{app:mixed-one-crossing-graphics}, together with the one-crossing terms already displayed in the preceding reductions. The dot positions for each comparison are marked below.
\begin{Lemma}\label{lem:appA-one-crossing-ss}
For $d\in \{   s_1,s_2\}$,
we have $C_{d,L}=C_{d,R}$.
\end{Lemma}
\begin{proof}
For $d=s_2$, the positions of the $x_i$ are indicated by the following diagram $
\begin{tikzpicture}
    \draw[] (0,0) to (0,1);
    \draw (0.5,0) to (1,1);
    \draw (1,0) to (0.5,1);
    \node at (0.65,0.3) {$\scriptstyle\bullet$};
    \node at (0.6,0.8) {$\scriptstyle\bullet$};\node at (0,0.5) {$\scriptstyle\bullet$};
        \node at (0.45,0.3) {$\scriptstyle3$};
    \node at (0.4,0.8) {$\scriptstyle2$};\node at (-0.2,0.5) {$\scriptstyle1$};
\end{tikzpicture}$.
Then we have
$$
C_{d,L} =2l_{3,2} -l_{3,2}l_{1,2} = l_{2,1}l_{3,2} = l_{2,1}(1-l_{1,3}+l_{1,2}) =C_{d,R}.$$
For $d=s_1$, the positions of the $x_i$ are indicated by the following diagram
$
\begin{tikzpicture}[xscale=-1]
    \draw[] (0,0) to (0,1);
    \draw (0.5,0) to (1,1);
    \draw (1,0) to (0.5,1);
    \node at (0.65,0.3) {$\scriptstyle\bullet$};
    \node at (0.6,0.8) {$\scriptstyle\bullet$};\node at (0,0.5) {$\scriptstyle\bullet$};
        \node at (0.45,0.3) {$\scriptstyle3$};
    \node at (0.4,0.8) {$\scriptstyle2$};\node at (-0.2,0.5) {$\scriptstyle1$};
\end{tikzpicture}$.
Then we have
\[
C_{d,L} =l_{1,2}(1-l_{3,1}+l_{2,1}) = l_{1,2}l_{2,3} =2l_{2,3} -l_{2,3}l_{2,1} =C_{d,R}.
\]
\end{proof}
\begin{Lemma}\label{lem:appA-one-crossing-es}
For $d\in \{e_1s_2,e_2s_1\}$,
we have $C_{d,L}=C_{d,R}$.
\end{Lemma}
\begin{proof}
   For $d=e_1s_2$, the positions of the $x_i$
are indicated by the following diagram $
\begin{tikzpicture}[baseline=(current bounding box.center),
  scale=0.8,line cap=round,line join=round,
  font={\fontsize{10pt}{12pt}\selectfont}]
  \begin{scope}[xscale=-1,yscale=-1]
  \draw[-,thick] (-0.4,0.6) to[out=270, in=180] (0.2,0.1);
  \draw[-,thick] (0.2,0.1) to[out = 0, in = 270] (0.4,0.6);
  \draw[-,thick] (-0.4,-0.6) to (0,0.6);
  \draw[-,thick] (0,-0.6) to[out=90, in=180] (0.2,-0.3);
  \draw[-,thick] (0.2,-0.3) to[out = 0, in = 90] (0.4,-0.6);
  \path (-0.4,-0.6) to  coordinate[pos=0.225] (aSix111p1) (0,0.6);
\path (0,-0.6) to[out=90, in=180]  coordinate[pos=0.4] (aSix111p2) (0.2,-0.3);
\path (-0.4,0.6) to[out=270, in=180]  coordinate[pos=0.72] (aSix111p3) (0.2,0.1);
\end{scope}
  \fill (aSix111p1) circle[radius=1.25pt];
  \node[inner sep=1pt] at (0.5,0.33) {$\scriptstyle 1$};
  \fill (aSix111p2) circle[radius=1.25pt];
  \node[inner sep=1pt] at (0.04,0.75) {$\scriptstyle 2$};
  \fill (aSix111p3) circle[radius=1.25pt];
  \node[inner sep=1pt] at (0.36,-0.09) {$\scriptstyle 3$};
\end{tikzpicture}$.
Then
\begin{align*}
 C_{d,L}
=& 2d_{-+}^{12}
 +\ell_{+-}^{12}W_N(-x_2)p_N(x_2,x_3)
  \bigl(2-d_{++}^{33}\bigr)
+\ell_{+-}^{12}W_N(-x_2)p_N(x_2,-x_1)
  q_V(x_1,x_3)M_V(-x_3)\ell_{+-}^{13}\\[4pt]
&{}-2q_V(x_1,x_3)M_V(-x_3)\ell_{+-}^{13}
-\ell_{-+}^{12}\ell_{+-}^{13}M_V(-x_3)
  \frac{q_V(x_1,x_3)-q_V(x_2,x_3)}{d_{+-}^{12}}\\[4pt]
&{}-\ell_{+-}^{12}\ell_{+-}^{13}W_N(-x_2)
  \frac{p_N(x_2,-x_1)-p_N(x_2,x_3)}{d_{--}^{13}}.
\end{align*}
On the right,
\[
C_{d,R}=2x_2-2x_2q_N(x_2,x_3)M_N(-x_3)(1-x_3-x_1).
\]
Using
$B_N(u)=\allowbreak\frac{1+u-x_1}{1+u+x_1}B_V(u)$ and
$A_X(-u)B_X(u)=\allowbreak\frac{1}{1+2u}$, together with the definitions of
$p_X,\allowbreak q_X,\allowbreak M_X,\allowbreak W_X$, the displayed expression for $C_{d,L}$ simplifies
to the displayed $C_{d,R}$.

 For $d=e_2s_1$, the positions of the $x_i$
are indicated by the following diagram
$
\begin{tikzpicture}[baseline=(current bounding box.center),
  scale=0.8,line cap=round,line join=round,
  font={\fontsize{10pt}{12pt}\selectfont}]
  \begin{scope}[xscale=1,yscale=-1]
  \draw[-,thick] (-0.4,0.6) to[out=270, in=180] (0.2,0.1);
  \draw[-,thick] (0.2,0.1) to[out = 0, in = 270] (0.4,0.6);
  \draw[-,thick] (-0.4,-0.6) to (0,0.6);
  \draw[-,thick] (0,-0.6) to[out=90, in=180] (0.2,-0.3);
  \draw[-,thick] (0.2,-0.3) to[out = 0, in = 90] (0.4,-0.6);
  \path (-0.4,-0.6) to  coordinate[pos=0.225] (aSix112p1) (0,0.6);
\path (0.2,-0.3) to[out = 0, in = 90]  coordinate[pos=0.65] (aSix112p2) (0.4,-0.6);
\path (0.2,0.1) to[out = 0, in = 270]  coordinate[pos=0.2] (aSix112p3) (0.4,0.6);
\end{scope}
  \fill (aSix112p1) circle[radius=1.25pt];
  \node[inner sep=1pt] at (-0.5,0.33) {$\scriptstyle 1$};
  \fill (aSix112p2) circle[radius=1.25pt];
  \node[inner sep=1pt] at (0.54,0.49) {$\scriptstyle 2$};
  \fill (aSix112p3) circle[radius=1.25pt];
  \node[inner sep=1pt] at (0.51,-0.13) {$\scriptstyle 3$};
\end{tikzpicture}$. We have
\begin{align*}
C_{d,L}={}&2x_2-1
 -(2x_2-1)q_V(x_2,x_3)M_V(-x_3)(1+x_1-x_3)\\[4pt]
 &+(1-2x_1)(1+x_1-x_3)M_V(-x_3)
 \frac{q_V(-x_1,x_3)-q_V(x_2,x_3)}{-x_1-x_2}.
\end{align*}
On the right,
\begin{align*}
C_{d,R}
={}&\ell_{--}^{12}W_V(-x_2)p_V(x_2,x_3)
  \bigl(2-d_{++}^{33}\bigr)
+\ell_{--}^{12}W_V(-x_2)p_V(x_2,x_1)
  q_N(-x_1,x_3)M_N(-x_3)\ell_{--}^{13}\\[4pt]
&-q_N(-x_1,x_3)M_N(-x_3)\ell_{--}^{13}
+\ell_{--}^{12}W_V(-x_2)\ell_{--}^{13}
  \frac{p_V(x_2,x_1)-p_V(x_2,x_3)}{d_{-+}^{13}}
  +2d_{++}^{12}-1.
\end{align*}
Substitution of the transition identity and the definitions of
$p_X,q_X,M_X,W_X$ gives $C_{d,L}=C_{d,R}$.
{The scalar calculation is described in Subsection~\ref{app:maple-checks}.}

\end{proof}

\begin{Lemma}\label{lem:appA-one-crossing-se} For $d\in \{ s_1e_2, s_2e_1\}$, we have
$C_{d,L}=C_{d,R}$.
\end{Lemma}

\begin{proof}
  For $d=s_1e_2$,  the positions of the $x_i$
are indicated by the following diagram
\[
\begin{tikzpicture}[baseline=(current bounding box.center),
  scale=0.8,line cap=round,line join=round,
  font={\fontsize{10pt}{12pt}\selectfont}]
  \begin{scope}[xscale=1,yscale=1]
  \draw[-,thick] (-0.4,0.6) to[out=270, in=180] (0.2,0.1);
  \draw[-,thick] (0.2,0.1) to[out = 0, in = 270] (0.4,0.6);
  \draw[-,thick] (-0.4,-0.6) to (0,0.6);
  \draw[-,thick] (0,-0.6) to[out=90, in=180] (0.2,-0.3);
  \draw[-,thick] (0.2,-0.3) to[out = 0, in = 90] (0.4,-0.6);
  \path (-0.4,-0.6) to  coordinate[pos=0.225] (aSix113p1) (0,0.6);
\path (0.2,-0.3) to[out = 0, in = 90]  coordinate[pos=0.65] (aSix113p2) (0.4,-0.6);
\path (0.2,0.1) to[out = 0, in = 270]  coordinate[pos=0.2] (aSix113p3) (0.4,0.6);
\end{scope}
  \fill (aSix113p1) circle[radius=1.25pt];
  \node[inner sep=1pt] at (-0.5,-0.33) {$\scriptstyle 1$};
  \fill (aSix113p2) circle[radius=1.25pt];
  \node[inner sep=1pt] at (0.54,-0.49) {$\scriptstyle 2$};
  \fill (aSix113p3) circle[radius=1.25pt];
  \node[inner sep=1pt] at (0.51,0.13) {$\scriptstyle 3$};
\end{tikzpicture}
\]
Then we have
 \begin{align*}
C_{d,R}
={}&\ell_{++}^{13}\ell_{++}^{12}\ell_{+-}^{12}
  \frac{r(x_1,x_2)-r(x_1,-x_3)}{d_{++}^{23}}-\ell_{++}^{13}
 -\ell_{++}^{13}\ell_{++}^{12}
  A_N(-x_1)q_N(x_3,-x_1)\\[4pt]
&-\ell_{++}^{12}\ell_{+-}^{12}M_V(-x_2)q_V(x_3,x_2)
 +\ell_{++}^{13}M_N(x_1)
  q_N(x_3,-x_1)\\[4pt]
&+
  \ell_{++}^{13}\ell_{++}^{12}\ell_{+-}^{12}
  A_N(-x_1)M_V(-x_2)q_V(x_1,x_2)q_N(x_3,-x_1)\\[4pt]
&+
  \ell_{++}^{13}\ell_{++}^{12}\ell_{+-}^{12}M_V(-x_2)
  \frac{r(x_1,-x_1)q_V(x_1,x_2)
       -r(x_1,-x_3)q_V(x_3,x_2)}
       {d_{-+}^{13}}\\[4pt]
&-
  \frac{\ell_{++}^{13}\ell_{++}^{12}\ell_{+-}^{12}}
       {d_{-+}^{13}}
  \left[
    \frac{r(x_1,-x_1)-r(x_1,x_2)}{d_{--}^{12}}
   -\frac{r(x_1,-x_3)-r(x_1,x_2)}
         {d_{--}^{23}}
  \right].
\end{align*}
and
\begin{align*}
&{}C_{d,L} =-l_{3,1} + l_{3,1}(x_1-x_2)M_V(-x_2)q_V(x_3,x_2)\\[4pt]
&{}+l_{3,1}(1+2x_1)\frac{M_V(x_1)q_V(x_3,-x_1)-M_V(-x_2)q_V(x_3,x_2)}{-x_1-x_2}
\end{align*}

At the marked through strand,
$B_N(u)=(1+u-x_1)(1+u+x_1)^{-1}B_V(u)$.
Direct substitution of this identity,
$A_X(-u)B_X(u)=(1+2u)^{-1}$, and the definitions of
$p_X,q_X,M_X,W_X$ and $r$ gives $C_{d,L}=C_{d,R}$.

For $d=s_2e_1$,
the positions of the $x_i$
are indicated by the following diagram
$
\begin{tikzpicture}[baseline=(current bounding box.center),
  scale=0.8,line cap=round,line join=round,
  font={\fontsize{10pt}{12pt}\selectfont}]
  \begin{scope}[xscale=-1,yscale=1]
  \draw[-,thick] (-0.4,0.6) to[out=270, in=180] (0.2,0.1);
  \draw[-,thick] (0.2,0.1) to[out = 0, in = 270] (0.4,0.6);
  \draw[-,thick] (-0.4,-0.6) to (0,0.6);
  \draw[-,thick] (0,-0.6) to[out=90, in=180] (0.2,-0.3);
  \draw[-,thick] (0.2,-0.3) to[out = 0, in = 90] (0.4,-0.6);
  \path (-0.4,-0.6) to  coordinate[pos=0.225] (aSix114p3) (0,0.6);
\path (0.2,-0.3) to[out = 0, in = 90]  coordinate[pos=0.65] (aSix114p1) (0.4,-0.6);
\path (0.2,0.1) to[out = 0, in = 270]  coordinate[pos=0.2] (aSix114p2) (0.4,0.6);
\end{scope}
  \fill (aSix114p3) circle[radius=1.25pt];
  \node[inner sep=1pt] at (0.5,-0.33) {$\scriptstyle 3$};
  \fill (aSix114p1) circle[radius=1.25pt];
  \node[inner sep=1pt] at (-0.55,-0.49) {$\scriptstyle 1$};
  \fill (aSix114p2) circle[radius=1.25pt];
  \node[inner sep=1pt] at (-0.52,0.13) {$\scriptstyle 2$};
\end{tikzpicture}$. We have
\begin{align*}
C_{d,L}
={}&-\ell_{--}^{23}\ell_{--}^{13}\ell_{+-}^{13}
  \frac{r(x_1,x_3)-r(-x_2,x_3)}{d_{++}^{12}}-\ell_{--}^{23}
 -\ell_{--}^{23}\ell_{--}^{13}
  A_V(x_3)q_V(-x_2,x_3)\\[4pt]
&-\ell_{--}^{13}\ell_{+-}^{13}M_N(x_1)q_N(-x_2,-x_1)
 +2\ell_{--}^{23}M_V(-x_3)
  q_V(-x_2,x_3)\\[4pt]
&-
  \ell_{--}^{23}\ell_{++}^{13}
  \frac{M_V(x_1)q_V(-x_2,-x_1)
       -M_V(-x_3)q_V(-x_2,x_3)}
       {d_{--}^{13}}\\[4pt]
&
  +\ell_{--}^{23}\ell_{--}^{13}\ell_{+-}^{13}
   A_V(x_3)M_N(x_1)q_V(-x_2,x_3)q_N(-x_3,-x_1)
  \\[4pt]
&-\ell_{--}^{23}\ell_{--}^{13}\ell_{+-}^{13}M_N(x_1)
  \frac{r(-x_2,x_3)q_N(-x_2,-x_1)
       -r(-x_3,x_3)q_N(-x_3,-x_1)}
       {d_{-+}^{23}}\\[4pt]
&-\frac{\ell_{--}^{23}\ell_{--}^{13}\ell_{+-}^{13}}
       {d_{-+}^{23}}
  \left[
    \frac{r(-x_2,x_3)-r(x_1,x_3)}{d_{--}^{12}}
   -\frac{r(-x_3,x_3)-r(x_1,x_3)}{d_{--}^{13}}
  \right].
\end{align*}

 and
 \[
C_{d,R} = -l_{3,2} +l_{3,2}l_{1,3}M_N(x_1)q_N(-x_2,-x_1).
\]

Here the marked through-strand variable is $x_3$.
Use $B_N(u)=(1+u-x_3)(1+u+x_3)^{-1}B_V(u)$ and
$A_X(-u)B_X(u)=(1+2u)^{-1}$.
Direct substitution of the definitions into the displayed coefficients
gives $C_{d,L}=C_{d,R}$.
\end{proof}

\subsection{Zero-crossing graphical calculations}
\label{app:zero-crossing-graphics}

We first record the direct source terms, then the zero-crossing corrections and bubble evaluations. The corresponding scalar coefficient comparisons are given in Appendix~B; the marked dot positions in each calculation specify its variables.

\subsubsection{Direct zero-crossing source terms}
\label{app:zero-crossing-sources}

The following identities concern the original source terms with no crossings. Corrections inherited from earlier reductions are collected in Subsection~\ref{app:zero-crossing-corrections}.

\newsavebox{\appAZeroPairBox}
\newlength{\appAZeroPairWidth}
\setlength{\appAZeroPairWidth}{\dimexpr(\linewidth-12pt)/2\relax}
\newcommand{\appAZeroPairCell}[1]{%
  \sbox{\appAZeroPairBox}{$\displaystyle #1$}%
  \makebox[\appAZeroPairWidth][l]{%
    \ifdim\wd\appAZeroPairBox>\appAZeroPairWidth
      \resizebox{\appAZeroPairWidth}{!}{\usebox{\appAZeroPairBox}}%
    \else
      \usebox{\appAZeroPairBox}%
    \fi
  }%
}

\begin{Lemma}
    When evaluated on a finitely generated object $V$ in $\cM_\lambda$, we have
\[
\begin{gathered}
\appAZeroPairCell{\tilde L_{222} ={} -
}
\end{gathered}
\]
\end{Lemma}

\begin{samepage}
\begin{Lemma}
    When evaluated on a finitely generated object $V$ in $\cM_\lambda$, we have
\[
\begin{gathered}
\appAZeroPairCell{\tilde R_{222} ={} -%
}
\end{gathered}
\]
\end{Lemma}
\end{samepage}

The subscripts select the parallel-strand ($2$) and cup--cap ($3$)
summands of the source diagrams in Subsection~\ref{27terms-app}.
No crossing reduction is needed for these terms.

\subsubsection{Zero-crossing corrections}
\label{app:zero-crossing-corrections}
{We retain all zero-crossing corrections from the preceding reductions.}

\paragraph{Initial transport.}
In this paragraph, $z_i$ is defined by the corresponding expression in the table below, and
$Y_i=[y_i(u)]_{{u:-1}}$ for $1\leq i\leq6$.
\par\smallskip
\begingroup

\pgfkeysifdefined{/tikz/grsw morphism/.@cmd}{}{%
	\tikzset{grsw morphism/.style={-,densely dashed,line width=.62pt}}%
}

\newsavebox{\grswInitialBox}
\newlength{\grswInitialCellWidth}
\setlength{\grswInitialCellWidth}{\dimexpr(\linewidth-20pt)/3\relax}
\newcommand{\grswInitialDiagram}[1]{%
  \sbox{\grswInitialBox}{$\thickmuskip=3mu\displaystyle #1$}%
  \raisebox{0pt}[34pt][24pt]{%
    \makebox[\grswInitialCellWidth][c]{%
      \ifdim\wd\grswInitialBox>\grswInitialCellWidth
        \resizebox{\grswInitialCellWidth}{!}{\usebox{\grswInitialBox}}%
      \else
        \usebox{\grswInitialBox}%
      \fi
    }%
  }%
}

\setlength{\tabcolsep}{0pt}
\renewcommand{\arraystretch}{1}
\setlength{\LTleft}{\fill}
\setlength{\LTright}{\fill}
\setlength{\LTpre}{0pt}
\setlength{\LTpost}{0pt}


  & = z_{7} +z_{8} -z_{9}
\end{align*}
Move the marked decoration through each crossing by
Definition~\ref{C defn}(4), and around a cup or cap by
Definition~\ref{C defn}(1),(5). The terms produced are the diagrams
named in the table; repeated contributions give the displayed integer
coefficients. In the fourth line, take the $u^{-1}$-coefficient after
the diagrammatic reduction.

\paragraph{Polynomial corrections.}
The labels in the following table are local to this paragraph.
Here again $Y_i=[y_i(u)]_{{u:-1}}$ for $1\leq i\leq6$.
\par\smallskip
\begingroup

\pgfkeysifdefined{/tikz/grsw morphism/.@cmd}{}{%
  \tikzset{grsw morphism/.style={-,densely dashed,line width=.62pt}}%
}
\pgfkeysifdefined{/tikz/cplus/.@cmd}{}{%
  \tikzset{cplus/.style={circle,draw=black,fill=white,
    line width=0.2pt,inner sep=0.025pt,minimum size=.2pt,
    font=\fontsize{5}{6}\selectfont}}%
}

\newsavebox{\grswPolyBox}
\newlength{\grswPolyCellWidth}
\setlength{\grswPolyCellWidth}{\dimexpr(\linewidth-40pt)/3\relax}
\newcommand{\grswAlignedDiagram}[2]{%
  \sbox{\grswPolyBox}{$\thickmuskip=3mu\displaystyle #1={}#2$}%
  \raisebox{0pt}[34pt][24pt]{%
    \makebox[\grswPolyCellWidth][c]{%
      \ifdim\wd\grswPolyBox>\grswPolyCellWidth
        \resizebox{\grswPolyCellWidth}{!}{\usebox{\grswPolyBox}}%
      \else
        \usebox{\grswPolyBox}%
      \fi
    }%
  }%
}

\setlength{\tabcolsep}{0pt}
\renewcommand{\arraystretch}{1}
\setlength{\LTleft}{\fill}
\setlength{\LTright}{\fill}
\setlength{\LTpre}{0pt}
\setlength{\LTpost}{0pt}

\right]_{{u:-1}}
  & = Y_{4} - Y_{5} + Y_{6}
\end{align*}
Apply the same crossing and bend relations to each source,
retaining the exterior factors $-1$ and $-2$ on the left. The first
four lines give the displayed $z_i$-terms. In the last two lines,
take the $u^{-1}$-coefficient after reduction to obtain the $Y_i$-terms.

\paragraph{Reflected corrections and divided differences.}

\par\medskip
The diagrams $z_i$ and $y_i(u)$ in this paragraph are defined below.
\par\smallskip

\par\medskip
\begingroup
\allowdisplaybreaks[1]

\ifdefined\tthreeOriginalWidth\else
\newlength{\tthreeOriginalWidth}
\fi
\ifdefined\tthreeOriginalBox\else
\newsavebox{\tthreeOriginalBox}
\fi
\ifdefined\tthreeOriginalLabelWidth\else
\newlength{\tthreeOriginalLabelWidth}
\fi
\setlength{\tthreeOriginalWidth}{\dimexpr(\linewidth-10pt)/2\relax}
\setlength{\tthreeOriginalLabelWidth}{3.5em}

\newcommand{\tthreeOriginalCell}[2]{%
	\begin{minipage}[c]{\tthreeOriginalWidth}%
		\sbox{\tthreeOriginalBox}{$\displaystyle #2$}%
		\makebox[\tthreeOriginalLabelWidth][r]{$\displaystyle#1={}$}%
		\ifdim\wd\tthreeOriginalBox>\dimexpr\linewidth-\tthreeOriginalLabelWidth\relax
			\resizebox{\dimexpr\linewidth-\tthreeOriginalLabelWidth\relax}{!}{\usebox{\tthreeOriginalBox}}%
		\else
			\usebox{\tthreeOriginalBox}%
		\fi
	\end{minipage}%
}

\setlength{\tabcolsep}{0pt}
\setlength{\LTleft}{\fill}
\setlength{\LTright}{\fill}
\setlength{\LTpre}{0pt}
\setlength{\LTpost}{0pt}
\renewcommand{\arraystretch}{1}

\endgroup
\par\medskip
\noindent where
\begin{align*}
	\Lambda_1 &:= l_{1,2}l_{1,3}
	\frac{r(x_1,-x_1)q_{V13}^{++}-r(x_1,u)q_V(-u,x_3)}{-x_1-u}M_{V3}^{-},\\
	\Lambda_2 &:= l_{1,2}l_{1,3}
\frac{r(x_1,u)-r(x_1,-x_2)}{u+x_2},\\
\Lambda_3 &:= l_{1,2}l_{1,3}\frac{1}{-x_1-u}
\left(
\frac{r(x_1,-x_1)-r(x_1,-x_2)}{-x_1+x_2}
-\frac{r(x_1,u)-r(x_1,-x_2)}{u+x_2}
\right).
\end{align*}
Set $Y_i=[y_i(u)]_{{u:-1}}$ for $1\leq i\leq12$.
\par\medskip
Moving the marked decorations through the reflected crossings gives

\begin{equation}\label{eq:appA-zero-12}
-\begin{tikzpicture}[
 baseline=(current bounding box.center),scale=0.8,
 line cap=round,line join=round,font={\fontsize{10pt}{12pt}\selectfont}]
 \draw[-,thick] (-0.4,0.6) to[out=270,in=180] (0.2,0.1);
 \draw[-,thick] (0.2,0.1) to[out=0,in=270] (0.4,0.6);
 \draw[-,thick] (-0.4,-0.6) -- (0,0.6);
 \draw[-,thick] (0,-0.6) to[out=90,in=180] (0.2,-0.3);
 \draw[-,thick] (0.2,-0.3) to[out=0,in=90] (0.4,-0.6);
 \path[overlay,name path=sourceCup]
 (-0.4,0.6) to[out=270,in=180] (0.2,0.1);
 \path[overlay,name path=sourceThrough] (-0.4,-0.6) -- (0,0.6);
 \path[overlay,name path=sourceCut] (-0.6,0.48) -- (0.2,0.48);
 \path[overlay,name intersections={of=sourceCup and sourceCut,by=sourceLeft}];
 \path[overlay,name intersections={of=sourceThrough and sourceCut,by=sourceRight}]; \draw[densely dashed,line width=.35pt] (sourceLeft) -- (sourceRight);
 \node[anchor=east,inner sep=1pt] at (-0.52,0.48)
 {$\scriptstyle l_{2,1}$};
\end{tikzpicture}
=-z_1+2z_2-2z_3.
\end{equation}

\begin{equation}\label{eq:appA-zero-13}
\begin{tikzpicture}
    [baseline=(current bounding box.center),
  scale=0.9,line cap=round,line join=round,
  font={\fontsize{10pt}{12pt}\selectfont}]
  \begin{scope}[xscale=-1,yscale=1]
  \draw[-,thick] (0,0.6) to (0.45,-0.6);
  \draw[-,thick] (-0.45,0.6) to[out=-90,in=180] (0,0.2);
  \draw[-,thick] (0,0.2) to[out=0,in=-90] (0.45,0.6);
  \draw[-,thick] (-0.45,-0.6) to[out=90,in=180] (-0.2,-0.2);
  \draw[-,thick] (-0.2,-0.2) to[out=0,in=90] (0,-0.6);
  \end{scope}
\begin{scope}[overlay,xscale=-1,yscale=1]
\path[overlay,name path global=clear275through] (0,0.6) to (0.45,-0.6);
\path[overlay,name path global=clear275othercup] (-0.45,0.6) to[out=-90,in=180] (0,0.2);
\end{scope}
\path[overlay,name path=clear275low] (-1,0.02)--(1,0.02);
\path[overlay,name path=clear275high] (-1,0.46)--(1,0.46);
\path[overlay,name intersections={of=clear275through and clear275low,by=clear275L}];
\path[overlay,name intersections={of=clear275othercup and clear275high,by=clear275R}];
\draw[grsw morphism,rounded corners=1pt,line width=0.350000pt] (clear275L)--(0.54,0.02)--(0.54,0.46)--(clear275R);
\node[anchor=west,inner sep=1pt] at (0.62,0.2) {$\scriptstyle -q_{N21}^{+-}A_{N1}^{-}$};
\begin{scope}[overlay,xscale=-1,yscale=1]
\path[overlay,name path global=tail275p3s] (0,0.2) to[out=0,in=-90] (0.45,0.6);
\end{scope}
\path[overlay,name path=tail275p3g] (-0.48,0.49)--(-0.38,0.49);
\path[overlay,name intersections={of=tail275p3s and tail275p3g,by=tail275p3}];
\begin{scope}[overlay,xscale=-1,yscale=1]
\path[overlay,name path global=tail275p4s] (0,0.6) to (0.45,-0.6);
\end{scope}
\path[overlay,name path=tail275p4g] (-0.09,0.49)--(0.01,0.49);
\path[overlay,name intersections={of=tail275p4s and tail275p4g,by=tail275p4}];
\draw[grsw morphism,dash pattern=on 1.3pt off 1.1pt,line width=0.350000pt] (tail275p3)--(tail275p4);
\node[anchor=east,inner sep=1pt,outer sep=0pt] at (-0.59,0.49) {$\scriptstyle l_{2,1}$};
\begin{scope}[overlay,xscale=-1,yscale=1]
\path[overlay,name path global=tail275p5s] (0,0.6) to (0.45,-0.6);
\end{scope}
\path[overlay,name path=tail275p5g] (-0.46,-0.5)--(-0.36,-0.5);
\path[overlay,name intersections={of=tail275p5s and tail275p5g,by=tail275p5}];
\begin{scope}[overlay,xscale=-1,yscale=1]
\path[overlay,name path global=tail275p6s] (-0.2,-0.2) to[out=0,in=90] (0,-0.6);
\end{scope}
\path[overlay,name path=tail275p6g] (-0.05,-0.5)--(0.05,-0.5);
\path[overlay,name intersections={of=tail275p6s and tail275p6g,by=tail275p6}];
\draw[grsw morphism,dash pattern=on 1.3pt off 1.1pt,line width=0.350000pt] (tail275p5)--(tail275p6);
\node[anchor=west,inner sep=1pt,outer sep=0pt] at (0.59,-0.5) {$\scriptstyle l_{1,2}$};
\end{tikzpicture}
  = z_{4} +Y_{1} -Y_{2}
\end{equation}

\begin{equation}\label{eq:appA-zero-14}
\begin{tikzpicture}
[baseline=(current bounding box.center),
  scale=0.8,line cap=round,line join=round,
  font={\fontsize{10pt}{12pt}\selectfont}]
  \begin{scope}[xscale=1,yscale=1]
  \draw[-,thick] (-0.4,0.6) to[out=270, in=180] (0.2,0.1);
  \draw[-,thick] (0.2,0.1) to[out = 0, in = 270] (0.4,0.6);
  \draw[-,thick] (-0.4,-0.6) to (0,0.6);
  \draw[-,thick] (0,-0.6) to[out=90, in=180] (0.2,-0.3);
  \draw[-,thick] (0.2,-0.3) to[out = 0, in = 90] (0.4,-0.6);
  \end{scope}
\begin{scope}[overlay,xscale=1,yscale=1]
\path[overlay,name path global=tail279p1s] (0.2,-0.3) to[out = 0, in = 90] (0.4,-0.6);
\end{scope}
\path[overlay,name path=tail279p1g] (0.36,-0.46)--(0.36,-0.36);
\path[overlay,name intersections={of=tail279p1s and tail279p1g,by=tail279p1}];
\begin{scope}[overlay,xscale=1,yscale=1]
\path[overlay,name path global=tail279p2s] (0.2,0.1) to[out = 0, in = 270] (0.4,0.6);
\end{scope}
\path[overlay,name path=tail279p2g] (0.36,0.16)--(0.36,0.26);
\path[overlay,name intersections={of=tail279p2s and tail279p2g,by=tail279p2}];
\draw[grsw morphism,dash pattern=on 1.3pt off 1.1pt,line width=0.350000pt] (tail279p1)--(tail279p2);
\node[anchor=west,inner sep=1pt,outer sep=0pt] at (0.56,-0.1) {$\scriptstyle q_{V21}^{++}M_{V1}^{-}$};
\begin{scope}[overlay,xscale=1,yscale=1]
\path[overlay,name path global=tail279p3s] (-0.4,-0.6) to (0,0.6);
\end{scope}
\path[overlay,name path=tail279p3g] (-0.42,-0.5)--(-0.32,-0.5);
\path[overlay,name intersections={of=tail279p3s and tail279p3g,by=tail279p3}];
\begin{scope}[overlay,xscale=1,yscale=1]
\path[overlay,name path global=tail279p4s] (0.2,-0.3) to[out = 0, in = 90] (0.4,-0.6);
\end{scope}
\path[overlay,name path=tail279p4g] (0.34,-0.5)--(0.44,-0.5);
\path[overlay,name intersections={of=tail279p4s and tail279p4g,by=tail279p4}];
\draw[grsw morphism,dash pattern=on 1.3pt off 1.1pt,line width=0.350000pt] (tail279p3)--(tail279p4);
\node[anchor=east,inner sep=1pt,outer sep=0pt] at (-0.56,-0.5) {$\scriptstyle l_{1,2}l_{1,3}$};
\begin{scope}[overlay,xscale=1,yscale=1]
\path[overlay,name path global=tail279p5s] (-0.4,0.6) to[out=270, in=180] (0.2,0.1);
\end{scope}
\path[overlay,name path=tail279p5g] (-0.43,0.46)--(-0.33,0.46);
\path[overlay,name intersections={of=tail279p5s and tail279p5g,by=tail279p5}];
\path[overlay,name path=fixA5Through] (-0.4,-0.6)--(0,0.6);
\path[overlay,name path=fixA5Cup]
(-0.4,0.6) to[out=270,in=180] (0.2,0.1);
\path[overlay,name path=fixA5Level] (-0.5,0.13)--(0.3,0.13);
\path[name intersections={of=fixA5Level and fixA5Through,by=fixA5L}];
\path[name intersections={of=fixA5Level and fixA5Cup,by=fixA5R}];
\begin{scope}[overlay,xscale=1,yscale=1]
\path[overlay,name path global=tail279p7s] (-0.4,-0.6) to (0,0.6);
\end{scope}
\path[overlay,name path=tail279p7g] (-0.18,0.01)--(-0.18,0.11);
\path[overlay,name intersections={of=tail279p7s and tail279p7g,by=tail279p7}];
\begin{scope}[overlay,xscale=1,yscale=1]
\path[overlay,name path global=tail279p8s] (-0.4,0.6) to[out=270, in=180] (0.2,0.1);
\end{scope}
\path[overlay,name path=tail279p8g] (-0.18,0.16)--(-0.18,0.26);
\path[overlay,name intersections={of=tail279p8s and tail279p8g,by=tail279p8}];
\path[overlay,name path=A5Cup]
  (-0.40,0.60) to[out=270,in=180] (0.20,0.10);
\path[overlay,name path=A5Through]
  (-0.40,-0.60)--(0.00,0.60);
\path[overlay,name path=A5Top]
  (-1,0.46)--(1,0.46);
\path[name intersections={of=A5Top and A5Cup,by=A5TopL}];
\path[name intersections={of=A5Top and A5Through,by=A5TopR}];
\draw[grsw morphism,line width=0.35pt] (A5TopL)--(A5TopR);
\node[anchor=east,inner sep=1pt] at (-0.56,0.49)
  {$\scriptstyle l_{2,1}$};
\path[overlay,name path=polyY4low]
  (-1,-0.04)--(1,-0.04);
\path[overlay,name path=polyY4high]
  (-1,0.14)--(1,0.14);\path[overlay,name path=polyY4through]
  (-0.4,-0.6)--(0,0.6);\path[overlay,name path=polyY4cup]
  (-0.4,0.6) to[out=270,in=180] (0.2,0.1);
\path[overlay,name intersections={
  of=polyY4low and polyY4through,by=polyY4left}];
\path[overlay,name intersections={
  of=polyY4high and polyY4cup,by=polyY4right}];

\draw[grsw morphism,line width=0.45pt,
  dash pattern=on 0.8pt off 0.7pt]
  (polyY4left) .. controls +(0.06,-0.14)
  and +(0,-0.14) .. (polyY4right);

\node[anchor=east,inner sep=1pt,outer sep=0pt]
  at (-0.56,0.04) {$\scriptstyle l_{1,2}^{-1}$};
\end{tikzpicture}
  = z_{5} -2z_{6} +2z_{7}
\end{equation}

\begin{equation}\label{eq:appA-zero-15}
\begin{tikzpicture}
    [baseline=(current bounding box.center),
  scale=0.9,line cap=round,line join=round,
  font={\fontsize{10pt}{12pt}\selectfont}]
  \begin{scope}[xscale=-1,yscale=1]
  \draw[-,thick] (0,0.6) to (0.45,-0.6);
  \draw[-,thick] (-0.45,0.6) to[out=-90,in=180] (0,0.2);
  \draw[-,thick] (0,0.2) to[out=0,in=-90] (0.45,0.6);
  \draw[-,thick] (-0.45,-0.6) to[out=90,in=180] (-0.2,-0.2);
  \draw[-,thick] (-0.2,-0.2) to[out=0,in=90] (0,-0.6);
  \end{scope}
\begin{scope}[overlay,xscale=-1,yscale=1]
\path[overlay,name path global=clear283through] (0,0.6) to (0.45,-0.6);
\path[overlay,name path global=clear283othercup] (-0.45,0.6) to[out=-90,in=180] (0,0.2);
\end{scope}
\path[overlay,name path=clear283low] (-1,0.02)--(1,0.02);
\path[overlay,name path=clear283high] (-1,0.46)--(1,0.46);
\path[overlay,name intersections={of=clear283through and clear283low,by=clear283L}];
\path[overlay,name intersections={of=clear283othercup and clear283high,by=clear283R}];
\draw[grsw morphism,rounded corners=1pt,line width=0.350000pt] (clear283L)--(0.54,0.02)--(0.54,0.46)--(clear283R);
\node[anchor=west,inner sep=1pt] at (0.62,0.2) {$\scriptstyle q_{N21}^{+-}A_{N1}^{-}$};
\begin{scope}[overlay,xscale=-1,yscale=1]
\path[overlay,name path global=tail283p3s] (0,0.2) to[out=0,in=-90] (0.45,0.6);
\end{scope}
\path[overlay,name path=tail283p3g] (-0.48,0.49)--(-0.38,0.49);
\path[overlay,name intersections={of=tail283p3s and tail283p3g,by=tail283p3}];
\begin{scope}[overlay,xscale=-1,yscale=1]
\path[overlay,name path global=tail283p4s] (0,0.6) to (0.45,-0.6);
\end{scope}
\path[overlay,name path=tail283p4g] (-0.09,0.49)--(0.01,0.49);
\path[overlay,name intersections={of=tail283p4s and tail283p4g,by=tail283p4}];
\draw[grsw morphism,dash pattern=on 1.3pt off 1.1pt,line width=0.350000pt] (tail283p3)--(tail283p4);
\node[anchor=east,inner sep=1pt,outer sep=0pt] at (-0.59,0.49) {$\scriptstyle l_{2,1}$};
\begin{scope}[overlay,xscale=-1,yscale=1]
\path[overlay,name path global=tail283p5s] (0,0.6) to (0.45,-0.6);
\end{scope}
\path[overlay,name path=tail283p5g] (-0.46,-0.5)--(-0.36,-0.5);
\path[overlay,name intersections={of=tail283p5s and tail283p5g,by=tail283p5}];
\begin{scope}[overlay,xscale=-1,yscale=1]
\path[overlay,name path global=tail283p6s] (-0.45,-0.6) to[out=90,in=180] (-0.2,-0.2);
\end{scope}
\path[overlay,name path=tail283p6g] (0.39,-0.5)--(0.49,-0.5);
\path[overlay,name intersections={of=tail283p6s and tail283p6g,by=tail283p6}];
\draw[grsw morphism,dash pattern=on 1.3pt off 1.1pt,line width=0.350000pt] (tail283p5)--(tail283p6);
\node[anchor=west,inner sep=1pt,outer sep=0pt] at (0.59,-0.5) {$\scriptstyle l_{1,2}l_{1,3}q_{V13}^{++}M_{V3}^{-}$};
\end{tikzpicture}
  = z_{8} + Y_{3} -Y_{4}
\end{equation}

\begin{equation}\label{eq:appA-zero-16}
\begin{tikzpicture}
    [baseline=(current bounding box.center),
  scale=0.9,line cap=round,line join=round,
  font={\fontsize{10pt}{12pt}\selectfont}]
  \begin{scope}[xscale=-1,yscale=1]
  \draw[-,thick] (0,0.6) to (0.45,-0.6);
  \draw[-,thick] (-0.45,0.6) to[out=-90,in=180] (0,0.2);
  \draw[-,thick] (0,0.2) to[out=0,in=-90] (0.45,0.6);
  \draw[-,thick] (-0.45,-0.6) to[out=90,in=180] (-0.2,-0.2);
  \draw[-,thick] (-0.2,-0.2) to[out=0,in=90] (0,-0.6);
  \end{scope}
\begin{scope}[overlay,xscale=-1,yscale=1]
\path[overlay,name path global=clear287through] (0,0.6) to (0.45,-0.6);
\path[overlay,name path global=clear287othercup] (-0.45,0.6) to[out=-90,in=180] (0,0.2);
\end{scope}
\path[overlay,name path=clear287low] (-1,0.02)--(1,0.02);
\path[overlay,name path=clear287high] (-1,0.46)--(1,0.46);
\path[overlay,name intersections={of=clear287through and clear287low,by=clear287L}];
\path[overlay,name intersections={of=clear287othercup and clear287high,by=clear287R}];
\draw[grsw morphism,rounded corners=1pt,line width=0.350000pt] (clear287L)--(0.54,0.02)--(0.54,0.46)--(clear287R);
\node[anchor=west,inner sep=1pt] at (0.62,0.2) {$\scriptstyle (1+x_1+x_2)q_{N21}^{+-}M_{N1}^{+}$};
\end{tikzpicture}
  = z_{9} + Y_{5} -Y_{6}
\end{equation}

\begin{equation}\label{eq:appA-zero-17}
\left[ \begin{tikzpicture}
    [baseline=(current bounding box.center),
  scale=0.9,line cap=round,line join=round,
  font={\fontsize{10pt}{12pt}\selectfont}]
  \begin{scope}[xscale=-1,yscale=1]
  \draw[-,thick] (0,0.6) to (0.45,-0.6);
  \draw[-,thick] (-0.45,0.6) to[out=-90,in=180] (0,0.2);
  \draw[-,thick] (0,0.2) to[out=0,in=-90] (0.45,0.6);
  \draw[-,thick] (-0.45,-0.6) to[out=90,in=180] (-0.2,-0.2);
  \draw[-,thick] (-0.2,-0.2) to[out=0,in=90] (0,-0.6);
  \end{scope}
\begin{scope}[overlay,xscale=-1,yscale=1]
\path[overlay,name path global=tail291p1s] (0,0.2) to[out=0,in=-90] (0.45,0.6);
\end{scope}
\path[overlay,name path=tail291p1g] (-0.48,0.49)--(-0.38,0.49);
\path[overlay,name intersections={of=tail291p1s and tail291p1g,by=tail291p1}];
\begin{scope}[overlay,xscale=-1,yscale=1]
\path[overlay,name path global=tail291p2s] (0,0.6) to (0.45,-0.6);
\end{scope}
\path[overlay,name path=tail291p2g] (-0.09,0.49)--(0.01,0.49);
\path[overlay,name intersections={of=tail291p2s and tail291p2g,by=tail291p2}];
\draw[grsw morphism,dash pattern=on 1.3pt off 1.1pt,line width=0.350000pt] (tail291p1)--(tail291p2);
\node[anchor=east,inner sep=1pt,outer sep=0pt] at (-0.59,0.49) {$\scriptstyle l_{2,1}$};
\begin{scope}[overlay,xscale=-1,yscale=1]
\path[overlay,name path global=tail291p3s] (0,0.2) to[out=0,in=-90] (0.45,0.6);
\end{scope}
\path[overlay,name path=tail291p3g] (-0.43,0.38)--(-0.33,0.38);
\path[overlay,name intersections={of=tail291p3s and tail291p3g,by=tail291p3}];
\node[cplus,scale=.7,inner sep=0pt] at (tail291p3) {$-$};
\begin{scope}[overlay,xscale=-1,yscale=1]
\path[overlay,name path global=tail291p4s] (0,0.6) to (0.45,-0.6);
\end{scope}
\path[overlay,name path=tail291p4g] (-0.32,-0.13)--(-0.22,-0.13);
\path[overlay,name intersections={of=tail291p4s and tail291p4g,by=tail291p4}];
\node[cplus,scale=.7,inner sep=0pt] at (tail291p4) {$+$};
\begin{scope}[overlay,xscale=-1,yscale=1]
\path[overlay,name path global=tail291p5s] (0,0.6) to (0.45,-0.6);
\end{scope}
\path[overlay,name path=tail291p5g] (-0.46,-0.5)--(-0.36,-0.5);
\path[overlay,name intersections={of=tail291p5s and tail291p5g,by=tail291p5}];
\begin{scope}[overlay,xscale=-1,yscale=1]
\path[overlay,name path global=tail291p6s] (-0.45,-0.6) to[out=90,in=180] (-0.2,-0.2);
\end{scope}
\path[overlay,name path=tail291p6g] (0.39,-0.5)--(0.49,-0.5);
\path[overlay,name intersections={of=tail291p6s and tail291p6g,by=tail291p6}];
\draw[grsw morphism,dash pattern=on 1.3pt off 1.1pt,line width=0.350000pt] (tail291p5)--(tail291p6);
\node[anchor=west,inner sep=1pt,outer sep=0pt] at (0.59,-0.5) {$\scriptstyle l_{1,2}l_{1,3}r(x_1,u)q_V(-u,x_3)M_{V3}^{-}$};
\end{tikzpicture}\right]_{{u:-1}}
  = z_{10} + Y_{7} -Y_{8}
\end{equation}

\begin{equation}\label{eq:appA-zero-18}
\left[ \begin{tikzpicture}
    [baseline=(current bounding box.center),
  scale=0.9,line cap=round,line join=round,
  font={\fontsize{10pt}{12pt}\selectfont}]
  \begin{scope}[xscale=-1,yscale=1]
  \draw[-,thick] (0,0.6) to (0.45,-0.6);
  \draw[-,thick] (-0.45,0.6) to[out=-90,in=180] (0,0.2);
  \draw[-,thick] (0,0.2) to[out=0,in=-90] (0.45,0.6);
  \draw[-,thick] (-0.45,-0.6) to[out=90,in=180] (-0.2,-0.2);
  \draw[-,thick] (-0.2,-0.2) to[out=0,in=90] (0,-0.6);
  \end{scope}
\begin{scope}[overlay,xscale=-1,yscale=1]
\path[overlay,name path global=tail295p1s] (0,0.2) to[out=0,in=-90] (0.45,0.6);
\end{scope}
\path[overlay,name path=tail295p1g] (-0.48,0.49)--(-0.38,0.49);
\path[overlay,name intersections={of=tail295p1s and tail295p1g,by=tail295p1}];
\begin{scope}[overlay,xscale=-1,yscale=1]
\path[overlay,name path global=tail295p2s] (0,0.6) to (0.45,-0.6);
\end{scope}
\path[overlay,name path=tail295p2g] (-0.09,0.49)--(0.01,0.49);
\path[overlay,name intersections={of=tail295p2s and tail295p2g,by=tail295p2}];
\draw[grsw morphism,dash pattern=on 1.3pt off 1.1pt,line width=0.350000pt] (tail295p1)--(tail295p2);
\node[anchor=east,inner sep=1pt,outer sep=0pt] at (-0.59,0.49) {$\scriptstyle l_{2,1}$};
\begin{scope}[overlay,xscale=-1,yscale=1]
\path[overlay,name path global=tail295p3s] (0,0.2) to[out=0,in=-90] (0.45,0.6);
\end{scope}
\path[overlay,name path=tail295p3g] (-0.43,0.38)--(-0.33,0.38);
\path[overlay,name intersections={of=tail295p3s and tail295p3g,by=tail295p3}];
\node[cplus,scale=.7,inner sep=0pt] at (tail295p3) {$-$};
\begin{scope}[overlay,xscale=-1,yscale=1]
\path[overlay,name path global=tail295p4s] (0,0.6) to (0.45,-0.6);
\end{scope}
\path[overlay,name path=tail295p4g] (-0.32,-0.13)--(-0.22,-0.13);
\path[overlay,name intersections={of=tail295p4s and tail295p4g,by=tail295p4}];
\node[cplus,scale=.7,inner sep=0pt] at (tail295p4) {$+$};
\begin{scope}[overlay,xscale=-1,yscale=1]
\path[overlay,name path global=tail295p5s] (0,0.6) to (0.45,-0.6);
\end{scope}
\path[overlay,name path=tail295p5g] (-0.46,-0.5)--(-0.36,-0.5);
\path[overlay,name intersections={of=tail295p5s and tail295p5g,by=tail295p5}];
\begin{scope}[overlay,xscale=-1,yscale=1]
\path[overlay,name path global=tail295p6s] (-0.45,-0.6) to[out=90,in=180] (-0.2,-0.2);
\end{scope}
\path[overlay,name path=tail295p6g] (0.39,-0.5)--(0.49,-0.5);
\path[overlay,name intersections={of=tail295p6s and tail295p6g,by=tail295p6}];
\draw[grsw morphism,dash pattern=on 1.3pt off 1.1pt,line width=0.350000pt] (tail295p5)--(tail295p6);
\node[anchor=west,inner sep=1pt,outer sep=0pt] at (0.59,-0.5) {$\scriptstyle -l_{1,2}l_{1,3}r(x_1,u)r(x_1,-x_2)$};
\end{tikzpicture}\right]_{{u:-1}}
  = -z_{11} -2Y_{9} -Y_{10} +2Y_{11} +Y_{12}
\end{equation}
In each line, use the crossing relation at the marked connector
and the bend relations at the adjacent cup or cap. The divided
differences on the diagrams are formed before extracting any
$u^{-1}$-coefficient.

\paragraph{Remaining one-crossing sources.}

The following diagrams belong to the remaining one-crossing sources.
\par\medskip
\begingroup
\ifdefined\tfourColWidth\else
\newlength{\tfourColWidth}
\fi
\ifdefined\tfourNameWidth\else
\newlength{\tfourNameWidth}
\fi
\ifdefined\tfourPictureWidth\else
\newlength{\tfourPictureWidth}
\fi
\ifdefined\tfourPictureBox\else
\newsavebox{\tfourPictureBox}
\fi
\setlength{\tfourColWidth}{\dimexpr(\linewidth-12pt)/2\relax}
\setlength{\tfourNameWidth}{3em}
\setlength{\tfourPictureWidth}{\dimexpr\tfourColWidth-\tfourNameWidth\relax}

\newcommand{\tfourCell}[2]{%
	\begingroup
	\sbox{\tfourPictureBox}{$\displaystyle #2$}%
	\makebox[\tfourColWidth][l]{%
		\makebox[\tfourNameWidth][r]{$\displaystyle #1={}$}%
		\ifdim\wd\tfourPictureBox>\tfourPictureWidth
		\resizebox{\tfourPictureWidth}{!}{\usebox{\tfourPictureBox}}%
		\else
		\usebox{\tfourPictureBox}%
		\fi
	}%
	\endgroup
}

\setlength{\tabcolsep}{0pt}
\setlength{\LTleft}{\fill}
\setlength{\LTright}{\fill}
\setlength{\LTpre}{0pt}
\setlength{\LTpost}{0pt}
\renewcommand{\arraystretch}{1}
\begin{longtable}{@{}l@{\hspace{12pt}}l@{}}
\multicolumn{2}{@{}l}{Define}\\
\end{longtable}
\endgroup
\par\medskip

Set $Y_i=[y_i(u)]_{{u:-1}}$ for $1\leq i\leq8$.
\par\medskip
Crossing transport and the cup--cap relations give

\begin{equation}\label{eq:appA-zero-19}
%
\right]_{{u:-1}}
  = z_{14} +Y_{7} -Y_{8}
\end{equation}
The terms on the right arise from the indicated crossing relation
and its cup--cap corrections. For a bracketed diagram, take its
$u^{-1}$-coefficient after reduction.

\paragraph{Corrections from higher-crossing sources.}

\par\medskip
The diagrams below are the terms left by the higher-crossing
reductions.
\par\smallskip
\begingroup
\ifdefined\higherCrossingWidth\else
\newlength{\higherCrossingWidth}
\fi
\ifdefined\higherCrossingLabelWidth\else
\newlength{\higherCrossingLabelWidth}
\fi
\ifdefined\higherCrossingBox\else
\newsavebox{\higherCrossingBox}
\fi
\setlength{\higherCrossingWidth}{\dimexpr(\linewidth-10pt)/2\relax}
\setlength{\higherCrossingLabelWidth}{3.5em}
\newcommand{\higherCrossingCell}[2]{%
	\begin{minipage}[c]{\higherCrossingWidth}%
		\sbox{\higherCrossingBox}{$\displaystyle #2$}%
		\makebox[\higherCrossingLabelWidth][r]{$\displaystyle#1={}$}%
		\ifdim\wd\higherCrossingBox>\dimexpr\linewidth-\higherCrossingLabelWidth\relax
			\resizebox{\dimexpr\linewidth-\higherCrossingLabelWidth\relax}{!}{\usebox{\higherCrossingBox}}%
		\else
			\usebox{\higherCrossingBox}%
		\fi
	\end{minipage}%
}
\setlength{\tabcolsep}{0pt}
\setlength{\LTleft}{\fill}
\setlength{\LTright}{\fill}
\setlength{\LTpre}{0pt}
\setlength{\LTpost}{0pt}
\renewcommand{\arraystretch}{1}

\endgroup
\par\medskip
Set $Y_i=[y_i(u)]_{{u:-1}}$ for $1\leq i\leq3$ and $8\leq i\leq15$.
\par\medskip
Applying the crossing and bend relations to each source gives

\begin{equation}\label{eq:appA-zero-26}
\left[ %
  = z_{11} + z_{12} - z_{13}
\end{equation}
These are the lower-crossing terms retained from the preceding
reductions. 

\paragraph{{Combining the terms.}}

\par\medskip
Define the diagrams in the following combinations.
\par\smallskip
\noindent\makebox[\linewidth][c]{%
$\displaystyle y_{1}={}%
$%
}\par\medskip
Set $Y_i=[y_i(u)]_{{u:-1}}$ for $1\leq i\leq15$.
\par\medskip
The bend relations and cancellation of opposite terms give

\begin{equation}\label{eq:appA-zero-36}
 %
\right]_{{u:-1}}
  = Y_{7} -Y_{8} +Y_{9}
\end{equation}

\begin{equation}\label{eq:appA-zero-41}
z_1
  = Y_{10} - Y_{11} + Y_{12}
\end{equation}
\begin{equation}\label{eq:appA-zero-42}
Y_3=Y_{13}-Y_{14}+Y_{15}.
\end{equation}
The four bubble terms are the last terms of
\eqref{eq:appA-zero-39} and \eqref{eq:appA-zero-40} and the first terms
on the right of \eqref{eq:appA-zero-41} and \eqref{eq:appA-zero-42}.
We evaluate them next.

\paragraph{Bubble evaluations.}

\par\medskip
The following diagrams record the required bubble terms.
\par\smallskip
\noindent\makebox[\linewidth][c]{%
$\displaystyle y_{1}={}
$%
}\par\medskip
Set $Y_i=[y_i(u)]_{{u:-1}}$ for $1\leq i\leq5$.
\par\medskip
Closing the indicated cups and caps and extracting the
$u^{-1}$-coefficients gives

\begin{equation}\label{eq:appA-zero-44}
\left[ %
\right]_{{u:-1}}
  = Y_{5}
\end{equation}

Close the marked cup--cap component, apply the bubble relations,
and then extract the $u^{-1}$-coefficient in each equality above.
\subsubsection{{The four bubble terms}}
\label{app:nodal-bubble-sum}
The symbols $y_1(u),\ldots,y_{24}(u)$ are local to this subsection;
throughout, $Y_i=[y_i(u)]_{{u:-1}}$. The diagrams, including their
marked decorations, are defined immediately before they are used.

{Write $\mathcal B_1,\ldots,\mathcal B_4$ for the four bubble terms
below. We compute their sum in Lemma~\ref{lem:nb-bubble-sum}.}

\par\medskip
\noindent Define\par\smallskip
\noindent\makebox[\linewidth][c]{%
$\displaystyle y_{1}={}\begin{tikzpicture}[baseline=(current bounding box.center),
  scale=0.8,line cap=round,line join=round,
  font={\fontsize{10pt}{12pt}\selectfont}]
  \begin{scope}[xscale=-1,yscale=1]
  \draw[-,thick] (0,-0.6) to[out=90, in=0] (-0.2,-0.3);
 \path[overlay,name path global=ndS1] (0,-0.6) to[out=90, in=0] (-0.2,-0.3);
  \draw[-,thick] (-0.2,-0.3) to[out = 180, in = 90] (-0.4,-0.6);
 \path[overlay,name path global=ndS2] (-0.2,-0.3) to[out = 180, in = 90] (-0.4,-0.6);
  \draw[-,thick] (-0.2,0.2) to[out=180,in=90] (-0.4,0);
 \path[overlay,name path global=ndS3] (-0.2,0.2) to[out=180,in=90] (-0.4,0);
  \draw[-,thick] (0,0) to[out=90,in=0] (-0.2,0.2);
 \path[overlay,name path global=ndS4] (0,0) to[out=90,in=0] (-0.2,0.2);
  \draw[-,thick] (-0.4,0) to[out=-90,in=180] (-0.2,-0.2);
 \path[overlay,name path global=ndS5] (-0.4,0) to[out=-90,in=180] (-0.2,-0.2);
  \draw[-,thick] (-0.2,-0.2) to[out=0,in=-90] (0,0);
 \path[overlay,name path global=ndS6] (-0.2,-0.2) to[out=0,in=-90] (0,0);
  \draw[-,thick] (0.4,0.6) to (0.4,-0.6);
 \path[overlay,name path global=ndS7] (0.4,0.6) to (0.4,-0.6);
  \draw[-,thick] (0,0.6) to[out=-90, in=0] (-0.2,0.3);
 \path[overlay,name path global=ndS8] (0,0.6) to[out=-90, in=0] (-0.2,0.3);
  \draw[-,thick] (-0.2,0.3) to[out = -180, in = -90] (-0.4,0.6);
 \path[overlay,name path global=ndS9] (-0.2,0.3) to[out = -180, in = -90] (-0.4,0.6);
  \end{scope}
 \path[overlay,name path=ndG1] (-2,0.54)--(2,0.54);
 \path[overlay,name intersections={of=ndG1 and ndS7,by=ndP1}];
 \path[overlay,name path=ndG2] (-2,0.54)--(2,0.54);
 \path[overlay,name intersections={of=ndG2 and ndS9,by=ndP2}];
 \draw[grsw morphism,line width=0.350000pt] (ndP1)--(ndP2);
 \node[anchor=south,inner sep=1pt] at (0,0.76) {$\scriptstyle l_{2,1}$};
 \path[overlay,name path=ndG3] (-2,-0.54)--(2,-0.54);
 \path[overlay,name intersections={of=ndG3 and ndS7,by=ndP3}];
 \path[overlay,name path=ndG4] (-2,-0.54)--(2,-0.54);
 \path[overlay,name intersections={of=ndG4 and ndS2,by=ndP4}];
 \draw[grsw morphism,line width=0.350000pt] (ndP3)--(ndP4);
 \node[anchor=north,inner sep=1pt] at (0,-0.76) {$\scriptstyle l_{1,2}l_{1,3}M_{V3}^{-}\frac{q_V(-u,x_3)-q_{V13}^{-+}}{-u+x_1}$};
 \path[overlay,name path=ndG5] (-2,0.04)--(2,0.04);
 \path[overlay,name intersections={of=ndG5 and ndS4,by=ndP5}];
 \node[cplus,scale=.7,inner sep=0pt] at (ndP5) {$-$};
 \path[overlay,name path=ndG6] (-2,0.43)--(2,0.43);
 \path[overlay,name intersections={of=ndG6 and ndS9,by=ndP6}];
 \node[cplus,scale=.7,inner sep=0pt] at (ndP6) {$-$};
\end{tikzpicture}$%
\hspace{8pt}%
$\displaystyle y_{2}={}\begin{tikzpicture}[baseline=(current bounding box.center),
  scale=0.8,line cap=round,line join=round,
  font={\fontsize{10pt}{12pt}\selectfont}]
  \begin{scope}[xscale=-1,yscale=1]
  \draw[-,thick] (0,-0.6) to[out=90, in=0] (-0.2,-0.3);
 \path[overlay,name path global=ndS1] (0,-0.6) to[out=90, in=0] (-0.2,-0.3);
  \draw[-,thick] (-0.2,-0.3) to[out = 180, in = 90] (-0.4,-0.6);
 \path[overlay,name path global=ndS2] (-0.2,-0.3) to[out = 180, in = 90] (-0.4,-0.6);
  \draw[-,thick] (-0.2,0.2) to[out=180,in=90] (-0.4,0);
 \path[overlay,name path global=ndS3] (-0.2,0.2) to[out=180,in=90] (-0.4,0);
  \draw[-,thick] (0,0) to[out=90,in=0] (-0.2,0.2);
 \path[overlay,name path global=ndS4] (0,0) to[out=90,in=0] (-0.2,0.2);
  \draw[-,thick] (-0.4,0) to[out=-90,in=180] (-0.2,-0.2);
 \path[overlay,name path global=ndS5] (-0.4,0) to[out=-90,in=180] (-0.2,-0.2);
  \draw[-,thick] (-0.2,-0.2) to[out=0,in=-90] (0,0);
 \path[overlay,name path global=ndS6] (-0.2,-0.2) to[out=0,in=-90] (0,0);
  \draw[-,thick] (0.4,0.6) to (0.4,-0.6);
 \path[overlay,name path global=ndS7] (0.4,0.6) to (0.4,-0.6);
  \draw[-,thick] (0,0.6) to[out=-90, in=0] (-0.2,0.3);
 \path[overlay,name path global=ndS8] (0,0.6) to[out=-90, in=0] (-0.2,0.3);
  \draw[-,thick] (-0.2,0.3) to[out = -180, in = -90] (-0.4,0.6);
 \path[overlay,name path global=ndS9] (-0.2,0.3) to[out = -180, in = -90] (-0.4,0.6);
  \end{scope}
 \path[overlay,name path=ndG1] (-2,0.54)--(2,0.54);
 \path[overlay,name intersections={of=ndG1 and ndS7,by=ndP1}];
 \path[overlay,name path=ndG2] (-2,0.54)--(2,0.54);
 \path[overlay,name intersections={of=ndG2 and ndS9,by=ndP2}];
 \draw[grsw morphism,line width=0.350000pt] (ndP1)--(ndP2);
 \node[anchor=south,inner sep=1pt] at (0,0.76) {$\scriptstyle l_{2,1}W_{V3}^{-}(p_V(x_3,-u)+\frac{p_{V31}^{+-}-p_V(x_3,-u)}{u-x_1})$};
 \path[overlay,name path=ndG3] (-2,-0.54)--(2,-0.54);
 \path[overlay,name intersections={of=ndG3 and ndS7,by=ndP3}];
 \path[overlay,name path=ndG4] (-2,-0.54)--(2,-0.54);
 \path[overlay,name intersections={of=ndG4 and ndS2,by=ndP4}];
 \draw[grsw morphism,line width=0.350000pt] (ndP3)--(ndP4);
 \node[anchor=north,inner sep=1pt] at (0,-0.76) {$\scriptstyle l_{1,2}$};
 \path[overlay,name path=ndG5] (-2,0.04)--(2,0.04);
 \path[overlay,name intersections={of=ndG5 and ndS4,by=ndP5}];
 \node[cplus,scale=.7,inner sep=0pt] at (ndP5) {$-$};
 \path[overlay,name path=ndG6] (-2,-0.43)--(2,-0.43);
 \path[overlay,name intersections={of=ndG6 and ndS2,by=ndP6}];
 \node[cplus,scale=.7,inner sep=0pt] at (ndP6) {$-$};
\end{tikzpicture}$%
\hspace{8pt}%
$\displaystyle y_{3}={}\begin{tikzpicture}[baseline=(current bounding box.center),
  scale=0.8,line cap=round,line join=round,
  font={\fontsize{10pt}{12pt}\selectfont}]
  \begin{scope}[xscale=-1,yscale=1]
  \draw[-,thick] (0,-0.6) to[out=90, in=0] (-0.2,-0.3);
 \path[overlay,name path global=ndS1] (0,-0.6) to[out=90, in=0] (-0.2,-0.3);
  \draw[-,thick] (-0.2,-0.3) to[out = 180, in = 90] (-0.4,-0.6);
 \path[overlay,name path global=ndS2] (-0.2,-0.3) to[out = 180, in = 90] (-0.4,-0.6);
  \draw[-,thick] (-0.2,0.2) to[out=180,in=90] (-0.4,0);
 \path[overlay,name path global=ndS3] (-0.2,0.2) to[out=180,in=90] (-0.4,0);
  \draw[-,thick] (0,0) to[out=90,in=0] (-0.2,0.2);
 \path[overlay,name path global=ndS4] (0,0) to[out=90,in=0] (-0.2,0.2);
  \draw[-,thick] (-0.4,0) to[out=-90,in=180] (-0.2,-0.2);
 \path[overlay,name path global=ndS5] (-0.4,0) to[out=-90,in=180] (-0.2,-0.2);
  \draw[-,thick] (-0.2,-0.2) to[out=0,in=-90] (0,0);
 \path[overlay,name path global=ndS6] (-0.2,-0.2) to[out=0,in=-90] (0,0);
  \draw[-,thick] (0.4,0.6) to (0.4,-0.6);
 \path[overlay,name path global=ndS7] (0.4,0.6) to (0.4,-0.6);
  \draw[-,thick] (0,0.6) to[out=-90, in=0] (-0.2,0.3);
 \path[overlay,name path global=ndS8] (0,0.6) to[out=-90, in=0] (-0.2,0.3);
  \draw[-,thick] (-0.2,0.3) to[out = -180, in = -90] (-0.4,0.6);
 \path[overlay,name path global=ndS9] (-0.2,0.3) to[out = -180, in = -90] (-0.4,0.6);
  \end{scope}
 \path[overlay,name path=ndG1] (-2,0.54)--(2,0.54);
 \path[overlay,name intersections={of=ndG1 and ndS7,by=ndP1}];
 \path[overlay,name path=ndG2] (-2,0.54)--(2,0.54);
 \path[overlay,name intersections={of=ndG2 and ndS9,by=ndP2}];
 \draw[grsw morphism,line width=0.350000pt] (ndP1)--(ndP2);
 \node[anchor=south,inner sep=1pt] at (0,0.76) {$\scriptstyle l_{2,1}W_{V3}^{-}p_V(x_3,-u)$};
 \path[overlay,name path=ndG3] (-2,-0.54)--(2,-0.54);
 \path[overlay,name intersections={of=ndG3 and ndS7,by=ndP3}];
 \path[overlay,name path=ndG4] (-2,-0.54)--(2,-0.54);
 \path[overlay,name intersections={of=ndG4 and ndS2,by=ndP4}];
 \draw[grsw morphism,line width=0.350000pt] (ndP3)--(ndP4);
 \node[anchor=north,inner sep=1pt] at (0,-0.76) {$\scriptstyle l_{1,2}l_{1,3}M_{V3}^{-}q_V(-u,x_3)$};
 \path[overlay,name path=ndG5] (-2,0.04)--(2,0.04);
 \path[overlay,name intersections={of=ndG5 and ndS4,by=ndP5}];
 \node[cplus,scale=.7,inner sep=0pt] at (ndP5) {$-$};
 \path[overlay,name path=ndG6] (-2,0)--(2,0);
 \path[overlay,name intersections={of=ndG6 and ndS7,by=ndP6}];
 \node[cplus,scale=.7,inner sep=0pt] at (ndP6) {$-$};
\end{tikzpicture}$%
}\par\medskip
\noindent\makebox[\linewidth][c]{%
$\displaystyle y_{4}={}\begin{tikzpicture}[baseline=(current bounding box.center),
  scale=0.8,line cap=round,line join=round,
  font={\fontsize{10pt}{12pt}\selectfont}]
  \begin{scope}[xscale=-1,yscale=1]
  \draw[-,thick] (0,-0.6) to[out=90, in=0] (-0.2,-0.3);
 \path[overlay,name path global=ndS1] (0,-0.6) to[out=90, in=0] (-0.2,-0.3);
  \draw[-,thick] (-0.2,-0.3) to[out = 180, in = 90] (-0.4,-0.6);
 \path[overlay,name path global=ndS2] (-0.2,-0.3) to[out = 180, in = 90] (-0.4,-0.6);
  \draw[-,thick] (-0.2,0.2) to[out=180,in=90] (-0.4,0);
 \path[overlay,name path global=ndS3] (-0.2,0.2) to[out=180,in=90] (-0.4,0);
  \draw[-,thick] (0,0) to[out=90,in=0] (-0.2,0.2);
 \path[overlay,name path global=ndS4] (0,0) to[out=90,in=0] (-0.2,0.2);
  \draw[-,thick] (-0.4,0) to[out=-90,in=180] (-0.2,-0.2);
 \path[overlay,name path global=ndS5] (-0.4,0) to[out=-90,in=180] (-0.2,-0.2);
  \draw[-,thick] (-0.2,-0.2) to[out=0,in=-90] (0,0);
 \path[overlay,name path global=ndS6] (-0.2,-0.2) to[out=0,in=-90] (0,0);
  \draw[-,thick] (0.4,0.6) to (0.4,-0.6);
 \path[overlay,name path global=ndS7] (0.4,0.6) to (0.4,-0.6);
  \draw[-,thick] (0,0.6) to[out=-90, in=0] (-0.2,0.3);
 \path[overlay,name path global=ndS8] (0,0.6) to[out=-90, in=0] (-0.2,0.3);
  \draw[-,thick] (-0.2,0.3) to[out = -180, in = -90] (-0.4,0.6);
 \path[overlay,name path global=ndS9] (-0.2,0.3) to[out = -180, in = -90] (-0.4,0.6);
  \end{scope}
 \path[overlay,name path=ndG1] (-2,0.54)--(2,0.54);
 \path[overlay,name intersections={of=ndG1 and ndS7,by=ndP1}];
 \path[overlay,name path=ndG2] (-2,0.54)--(2,0.54);
 \path[overlay,name intersections={of=ndG2 and ndS9,by=ndP2}];
 \draw[grsw morphism,line width=0.350000pt] (ndP1)--(ndP2);
 \node[anchor=south,inner sep=1pt] at (0,0.76) {$\scriptstyle l_{2,1}W_{V3}^{-}p_V(x_3,-u)$};
 \path[overlay,name path=ndG3] (-2,-0.54)--(2,-0.54);
 \path[overlay,name intersections={of=ndG3 and ndS7,by=ndP3}];
 \path[overlay,name path=ndG4] (-2,-0.54)--(2,-0.54);
 \path[overlay,name intersections={of=ndG4 and ndS2,by=ndP4}];
 \draw[grsw morphism,line width=0.350000pt] (ndP3)--(ndP4);
 \node[anchor=north,inner sep=1pt] at (0,-0.76) {$\scriptstyle l_{1,2}l_{1,3}\frac{r(x_1,x_3)-1}{x_1-x_3}$};
 \path[overlay,name path=ndG5] (-2,0.04)--(2,0.04);
 \path[overlay,name intersections={of=ndG5 and ndS4,by=ndP5}];
 \node[cplus,scale=.7,inner sep=0pt] at (ndP5) {$-$};
 \path[overlay,name path=ndG6] (-2,0)--(2,0);
 \path[overlay,name intersections={of=ndG6 and ndS7,by=ndP6}];
 \node[cplus,scale=.7,inner sep=0pt] at (ndP6) {$-$};
\end{tikzpicture}$%
}\par\medskip
\par\medskip
\begin{equation}\label{eq:nb-bubble-1}
\mathcal B_1
  = -Y_{1}
\end{equation}

\begin{equation}\label{eq:nb-bubble-2}
\mathcal B_2
  = -Y_{2}
\end{equation}

\begin{equation}\label{eq:nb-bubble-3}
\mathcal B_3
  = -Y_{3}
\end{equation}

\begin{equation}\label{eq:nb-bubble-4}
\mathcal B_4
  = -Y_{4}
\end{equation}

\par\medskip
\noindent Define\par\smallskip
\noindent\makebox[\linewidth][c]{%
$\displaystyle y_{5}={}\begin{tikzpicture}[baseline=(current bounding box.center),
  scale=0.8,line cap=round,line join=round,
  font={\fontsize{10pt}{12pt}\selectfont}]
  \begin{scope}[xscale=-1,yscale=1]
  \draw[-,thick] (0,-0.6) to[out=90, in=0] (-0.2,-0.3);
 \path[overlay,name path global=ndS1] (0,-0.6) to[out=90, in=0] (-0.2,-0.3);
  \draw[-,thick] (-0.2,-0.3) to[out = 180, in = 90] (-0.4,-0.6);
 \path[overlay,name path global=ndS2] (-0.2,-0.3) to[out = 180, in = 90] (-0.4,-0.6);
  \draw[-,thick] (-0.2,0.2) to[out=180,in=90] (-0.4,0);
 \path[overlay,name path global=ndS3] (-0.2,0.2) to[out=180,in=90] (-0.4,0);
  \draw[-,thick] (0,0) to[out=90,in=0] (-0.2,0.2);
 \path[overlay,name path global=ndS4] (0,0) to[out=90,in=0] (-0.2,0.2);
  \draw[-,thick] (-0.4,0) to[out=-90,in=180] (-0.2,-0.2);
 \path[overlay,name path global=ndS5] (-0.4,0) to[out=-90,in=180] (-0.2,-0.2);
  \draw[-,thick] (-0.2,-0.2) to[out=0,in=-90] (0,0);
 \path[overlay,name path global=ndS6] (-0.2,-0.2) to[out=0,in=-90] (0,0);
  \draw[-,thick] (0.4,0.6) to (0.4,-0.6);
 \path[overlay,name path global=ndS7] (0.4,0.6) to (0.4,-0.6);
  \draw[-,thick] (0,0.6) to[out=-90, in=0] (-0.2,0.3);
 \path[overlay,name path global=ndS8] (0,0.6) to[out=-90, in=0] (-0.2,0.3);
  \draw[-,thick] (-0.2,0.3) to[out = -180, in = -90] (-0.4,0.6);
 \path[overlay,name path global=ndS9] (-0.2,0.3) to[out = -180, in = -90] (-0.4,0.6);
  \end{scope}
 \path[overlay,name path=ndG1] (-2,0.54)--(2,0.54);
 \path[overlay,name intersections={of=ndG1 and ndS7,by=ndP1}];
 \path[overlay,name path=ndG2] (-2,0.54)--(2,0.54);
 \path[overlay,name intersections={of=ndG2 and ndS9,by=ndP2}];
 \draw[grsw morphism,line width=0.350000pt] (ndP1)--(ndP2);
 \node[anchor=south,inner sep=1pt] at (0,0.76) {$\scriptstyle l_{2,1}W_{V3}^{-}$};
 \path[overlay,name path=ndG3] (-2,-0.54)--(2,-0.54);
 \path[overlay,name intersections={of=ndG3 and ndS7,by=ndP3}];
 \path[overlay,name path=ndG4] (-2,-0.54)--(2,-0.54);
 \path[overlay,name intersections={of=ndG4 and ndS2,by=ndP4}];
 \draw[grsw morphism,line width=0.350000pt] (ndP3)--(ndP4);
 \node[anchor=north,inner sep=1pt] at (0,-0.76) {$\scriptstyle l_{1,2}l_{1,3}M_{V3}^{-}$};
 \path[overlay,name path=ndG5] (-2,0.43)--(2,0.43);
 \path[overlay,name intersections={of=ndG5 and ndS9,by=ndP5}];
 \node[cplus,scale=.7,inner sep=0pt] at (ndP5) {$-$};
 \path[overlay,name path=ndG6] (-2,0.04)--(2,0.04);
 \path[overlay,name intersections={of=ndG6 and ndS4,by=ndP6}];
 \node[cplus,scale=.7,inner sep=0pt] at (ndP6) {$-$};
 \path[overlay,name path=ndG7] (-2,-0.43)--(2,-0.43);
 \path[overlay,name intersections={of=ndG7 and ndS2,by=ndP7}];
 \node[cplus,scale=.7,inner sep=0pt] at (ndP7) {$-$};
 \path[overlay,name path=ndG8] (-2,0)--(2,0);
 \path[overlay,name intersections={of=ndG8 and ndS7,by=ndP8}];
 \node[cplus,scale=.7,inner sep=0pt] at (ndP8) {$-$};
 \node[anchor=east,inner sep=1pt] at (-0.56,0.1) {$\scriptstyle -A_V(-u)B_V(-u)$};
\end{tikzpicture}$%
\hspace{8pt}%
$\displaystyle y_{6}={}\begin{tikzpicture}[baseline=(current bounding box.center),
  scale=0.8,line cap=round,line join=round,
  font={\fontsize{10pt}{12pt}\selectfont}]
  \begin{scope}[xscale=-1,yscale=1]
  \draw[-,thick] (0,-0.6) to[out=90, in=0] (-0.2,-0.3);
 \path[overlay,name path global=ndS1] (0,-0.6) to[out=90, in=0] (-0.2,-0.3);
  \draw[-,thick] (-0.2,-0.3) to[out = 180, in = 90] (-0.4,-0.6);
 \path[overlay,name path global=ndS2] (-0.2,-0.3) to[out = 180, in = 90] (-0.4,-0.6);
  \draw[-,thick] (0.4,0.6) to (0.4,-0.6);
 \path[overlay,name path global=ndS3] (0.4,0.6) to (0.4,-0.6);
  \draw[-,thick] (0,0.6) to[out=-90, in=0] (-0.2,0.3);
 \path[overlay,name path global=ndS4] (0,0.6) to[out=-90, in=0] (-0.2,0.3);
  \draw[-,thick] (-0.2,0.3) to[out = -180, in = -90] (-0.4,0.6);
 \path[overlay,name path global=ndS5] (-0.2,0.3) to[out = -180, in = -90] (-0.4,0.6);
  \end{scope}
 \path[overlay,name path=ndG1] (-2,0.54)--(2,0.54);
 \path[overlay,name intersections={of=ndG1 and ndS3,by=ndP1}];
 \path[overlay,name path=ndG2] (-2,0.54)--(2,0.54);
 \path[overlay,name intersections={of=ndG2 and ndS5,by=ndP2}];
 \draw[grsw morphism,line width=0.350000pt] (ndP1)--(ndP2);
 \node[anchor=south,inner sep=1pt] at (0,0.76) {$\scriptstyle l_{2,1}W_{V3}^{-}$};
 \path[overlay,name path=ndG3] (-2,-0.54)--(2,-0.54);
 \path[overlay,name intersections={of=ndG3 and ndS3,by=ndP3}];
 \path[overlay,name path=ndG4] (-2,-0.54)--(2,-0.54);
 \path[overlay,name intersections={of=ndG4 and ndS2,by=ndP4}];
 \draw[grsw morphism,line width=0.350000pt] (ndP3)--(ndP4);
 \node[anchor=north,inner sep=1pt] at (0,-0.76) {$\scriptstyle l_{1,2}l_{1,3}M_{V3}^{-}$};
 \path[overlay,name path=ndG5] (-2,0.43)--(2,0.43);
 \path[overlay,name intersections={of=ndG5 and ndS5,by=ndP5}];
 \node[cplus,scale=.7,inner sep=0pt] at (ndP5) {$-$};
 \path[overlay,name path=ndG6] (-2,-0.43)--(2,-0.43);
 \path[overlay,name intersections={of=ndG6 and ndS2,by=ndP6}];
 \node[cplus,scale=.7,inner sep=0pt] at (ndP6) {$-$};
 \path[overlay,name path=ndG7] (-2,0)--(2,0);
 \path[overlay,name intersections={of=ndG7 and ndS3,by=ndP7}];
 \node[cplus,scale=.7,inner sep=0pt] at (ndP7) {$-$};
 \node[anchor=east,inner sep=1pt] at (-0.56,0.1) {$\scriptstyle \frac{(-1)^{c+1}}2m_{V,0}^+(-u)f_V^+(-u)m_{V,0}^-(-u)f_V^-(-u)\frac1{-u}(\frac{n_V(-u)}{m_V(-u)}+u+\frac12)$};
\end{tikzpicture}$%
}\par\medskip
\par\medskip
\begin{Lemma}\label{lem:nb-bubble-sum}
The sum of these four terms is

\begin{equation}\label{eq:nb-bubble-sum}
\mathcal B_1 +\mathcal B_2 +\mathcal B_3 +\mathcal B_4
  = Y_{5}
\end{equation}

Equivalently,

\begin{equation}\label{eq:nb-bubble-expanded}
\mathcal B_1 +\mathcal B_2 +\mathcal B_3 +\mathcal B_4
  = Y_{6}
\end{equation}

\end{Lemma}

\begin{proof}
{We first derive the four bubble terms.}
The following pair cancels:
\newcommand{\Auppercurl}[2]{%
\begin{tikzpicture}[baseline=(current bounding box.center),scale=0.8,
  line cap=round,line join=round,font={\fontsize{10pt}{12pt}\selectfont}]
  \draw[thick] (-0.45,0.62)--(-0.45,-0.62);
  \draw[thick] (-0.05,-0.62) to[out=90,in=180] (0.20,-0.42)
                    to[out=0,in=90] (0.45,-0.62);
  \draw[thick] (0,0.62)
    .. controls (0,0.34) and (0.42,0.38) .. (0.31,0.05)
    .. controls (0.24,-0.17) and (-0.18,-0.17) .. (-0.24,0.05)
    .. controls (-0.33,0.38) and (0.40,0.34) .. (0.40,0.62);
  \draw[densely dashed,line width=.35pt] (-0.45,0.54)--(0.40,0.54);
  \node[anchor=south,inner sep=1pt] at (0,0.76) {$\scriptstyle #1$};
  \draw[densely dashed,line width=.35pt] (-0.45,-0.54)--(0.45,-0.54);
  \node[anchor=north,inner sep=1pt] at (0,-0.76) {$\scriptstyle #2$};
  \node[cplus,scale=.7,inner sep=0pt] at (-0.24,0.05) {$-$};
\end{tikzpicture}}

\newcommand{\Athroughcurl}[6]{%
\begin{tikzpicture}[baseline=(current bounding box.center),scale=0.8,
  line cap=round,line join=round,font={\fontsize{10pt}{12pt}\selectfont}]
  \draw[thick] (-0.45,0.62)--(-0.45,0.28)--(0.17,-0.12);
  \draw[thick] (-0.45,-0.62)--(-0.45,-0.28)--(0.17,0.12);
  \draw[thick] (0.17,0.12)
    .. controls (0.64,0.46) and (0.64,-0.46) .. (0.17,-0.12);
  \draw[thick] (-0.05,0.62) to[out=-90,in=180] (0.20,0.42)
                    to[out=0,in=-90] (0.45,0.62);
  \draw[thick] (-0.05,-0.62) to[out=90,in=180] (0.20,-0.42)
                    to[out=0,in=90] (0.45,-0.62);
  \draw[densely dashed,line width=.35pt] (-0.45,0.54)--(0.45,0.54);
  \node[anchor=south,inner sep=1pt] at (0,0.76) {$\scriptstyle #1$};
  \draw[densely dashed,line width=.35pt] (-0.45,-0.54)--(0.45,-0.54);
  \node[anchor=north,inner sep=1pt] at (0,-0.76) {$\scriptstyle #2$};
  \ifx\relax#3\relax\else
    \node[anchor=west,inner sep=1pt] at (0.56,0.26) {$\scriptstyle #3$};
  \fi
  \ifx\relax#4\relax\else
    \node[anchor=west,inner sep=1pt] at (0.56,-0.26) {$\scriptstyle #4$};
  \fi
  \ifx\relax#5\relax\else
    \node[anchor=east,inner sep=1pt] at (-0.56,0) {$\scriptstyle #5$};
  \fi
  \def\Auppermark{upper}%
  \def\Amarkarg{#6}%
  \ifx\Amarkarg\Auppermark
    \node[cplus,scale=.7,inner sep=0pt] at (-0.45,0.40) {$-$};
  \else
    \node[cplus,scale=.7,inner sep=0pt] at (0.3,-0.19) {$-$};
  \fi
\end{tikzpicture}} \newcommand{\AthroughQurl}[6]{%
\begin{tikzpicture}[baseline=(current bounding box.center),scale=0.8,
  line cap=round,line join=round,font={\fontsize{10pt}{12pt}\selectfont}]
  \draw[thick] (-0.45,0.62)--(-0.45,0.28)--(0.17,-0.12);
  \draw[thick] (-0.45,-0.62)--(-0.45,-0.28)--(0.17,0.12);
  \draw[thick] (0.17,0.12)
    .. controls (0.64,0.46) and (0.64,-0.46) .. (0.17,-0.12);
  \draw[thick] (-0.05,0.62) to[out=-90,in=180] (0.20,0.42)
                    to[out=0,in=-90] (0.45,0.62);
  \draw[thick] (-0.05,-0.62) to[out=90,in=180] (0.20,-0.42)
                    to[out=0,in=90] (0.45,-0.62);
  \draw[densely dashed,line width=.35pt] (-0.45,0.54)--(0.45,0.54);
  \node[anchor=south,inner sep=1pt] at (0,0.76) {$\scriptstyle #1$};
  \draw[densely dashed,line width=.35pt] (-0.45,-0.54)--(0.45,-0.54);
  \node[anchor=north,inner sep=1pt] at (0,-0.76) {$\scriptstyle #2$};

  \path[overlay,name path=acUpperLoop] (0.17,0.12)
    .. controls (0.405,0.29) and (0.5225,0.145) .. (0.5225,0);
  \path[overlay,name path=acLowerLoop] (0.5225,0)
    .. controls (0.5225,-0.145) and (0.405,-0.29) .. (0.17,-0.12);
  \path[overlay,name path=acUpperCap] (0.20,0.42)
    to[out=0,in=-90] (0.45,0.62);
  \path[overlay,name path=acLowerCup] (0.20,-0.42)
    to[out=0,in=90] (0.45,-0.62);
  \path[overlay,name path=acVertical] (0.34,-0.65)--(0.34,0.65);
  \ifx\relax#3\relax\else
    \path[name intersections={of=acVertical and acUpperLoop,by=acU1}];
    \path[name intersections={of=acVertical and acUpperCap,by=acU2}];
    \draw[grsw morphism,line width=.35pt,dash pattern=on .8pt off .8pt] (acU1)--(acU2);

    \node[anchor=west,inner sep=1pt] at (0.56,0.26) {$\scriptstyle #3$};
  \fi
  \ifx\relax#4\relax\else
    \path[name intersections={of=acVertical and acLowerCup,by=acL1}];
    \path[name intersections={of=acVertical and acLowerLoop,by=acL2}];

    \node[anchor=west,inner sep=1pt] at (0.56,-0.26) {$\scriptstyle #4$};
  \fi
  \ifx\relax#5\relax\else
    \node[anchor=east,inner sep=1pt] at (-0.56,0) {$\scriptstyle #5$};
  \fi
  \def\Auppermark{upper}%
  \def\Amarkarg{#6}%
  \ifx\Amarkarg\Auppermark
    \node[cplus,scale=.7,inner sep=0pt] at (-0.45,0.40) {$-$};
  \else
    \node[cplus,scale=.7,inner sep=0pt] at (0.3,-0.19) {$-$};
  \fi
\end{tikzpicture}} 
\newcommand{\AthroughPurl}[6]{%
\begin{tikzpicture}[baseline=(current bounding box.center),scale=0.8,
  line cap=round,line join=round,font={\fontsize{10pt}{12pt}\selectfont}]
  \draw[thick] (-0.45,0.62)--(-0.45,0.28)--(0.17,-0.12);
  \draw[thick] (-0.45,-0.62)--(-0.45,-0.28)--(0.17,0.12);
  \draw[thick] (0.17,0.12)
    .. controls (0.64,0.46) and (0.64,-0.46) .. (0.17,-0.12);
  \draw[thick] (-0.05,0.62) to[out=-90,in=180] (0.20,0.42)
                    to[out=0,in=-90] (0.45,0.62);
  \draw[thick] (-0.05,-0.62) to[out=90,in=180] (0.20,-0.42)
                    to[out=0,in=90] (0.45,-0.62);
  \draw[densely dashed,line width=.35pt] (-0.45,0.54)--(0.45,0.54);
  \node[anchor=south,inner sep=1pt] at (0,0.76) {$\scriptstyle #1$};
  \draw[densely dashed,line width=.35pt] (-0.45,-0.54)--(0.45,-0.54);
  \node[anchor=north,inner sep=1pt] at (0,-0.76) {$\scriptstyle #2$};

  \path[overlay,name path=acUpperLoop] (0.17,0.12)
    .. controls (0.405,0.29) and (0.5225,0.145) .. (0.5225,0);
  \path[overlay,name path=acLowerLoop] (0.5225,0)
    .. controls (0.5225,-0.145) and (0.405,-0.29) .. (0.17,-0.12);
  \path[overlay,name path=acUpperCap] (0.20,0.42)
    to[out=0,in=-90] (0.45,0.62);
  \path[overlay,name path=acLowerCup] (0.20,-0.42)
    to[out=0,in=90] (0.45,-0.62);
  \path[overlay,name path=acVertical] (0.34,-0.65)--(0.34,0.65);
  \ifx\relax#3\relax\else
    \path[name intersections={of=acVertical and acUpperLoop,by=acU1}];
    \path[name intersections={of=acVertical and acUpperCap,by=acU2}];
    \draw[grsw morphism,line width=.35pt,dash pattern=on .8pt off .8pt] (acU1)--(acU2);

    \node[anchor=west,inner sep=1pt] at (0.56,0.26) {$\scriptstyle #3$};
  \fi
  \ifx\relax#4\relax\else
    \path[name intersections={of=acVertical and acLowerCup,by=acL1}];
    \path[name intersections={of=acVertical and acLowerLoop,by=acL2}];
    \draw[grsw morphism,line width=.35pt,dash pattern=on .8pt off .8pt] (acL1)--(acL2);

    \node[anchor=west,inner sep=1pt] at (0.56,-0.26) {$\scriptstyle #4$};
  \fi
  \ifx\relax#5\relax\else
    \node[anchor=east,inner sep=1pt] at (-0.56,0) {$\scriptstyle #5$};
  \fi
  \def\Auppermark{upper}%
  \def\Amarkarg{#6}%
  \ifx\Amarkarg\Auppermark
    \node[cplus,scale=.7,inner sep=0pt] at (-0.45,0.40) {$-$};
  \else
    \node[cplus,scale=.7,inner sep=0pt] at (0.3,-0.19) {$-$};
  \fi
\end{tikzpicture}
}
\begin{equation}\label{eq:nb-cancel}
-\left[
\Auppercurl
 {l_{2,1}}
 {l_{1,2}l_{1,3}\dfrac{1-r(x_1,u)}{(x_1-u)(u-x_3)}}
\right]_{{u:-1}}
+\left[
\Auppercurl
 {l_{2,1}}
 {l_{1,2}l_{1,3}\dfrac{r(x_1,u)-r(x_1,x_3)}{(u-x_3)(u-x_1)}}
\right]_{{u:-1}}
=0.
\end{equation}
\par\medskip
\noindent Define\par\smallskip
\newcommand{\cleanPCurl}[2]{%
\begin{tikzpicture}[baseline=(current bounding box.center),scale=0.8,
  line cap=round,line join=round,draw=black,text=black,
  font={\fontsize{10pt}{12pt}\selectfont}]
  \draw[thick] (-0.45,0.62)--(-0.45,-0.62);
  \draw[thick] (-0.05,0.62) to[out=-90,in=180] (0.20,0.42)
                    to[out=0,in=-90] (0.45,0.62);
  \draw[thick] (0,-0.62)
    .. controls (0,-0.34) and (0.42,-0.38) .. (0.31,-0.05)
    .. controls (0.24,0.17) and (-0.18,0.17) .. (-0.24,-0.05)
    .. controls (-0.33,-0.38) and (0.40,-0.34) .. (0.40,-0.62);
  \draw[densely dashed,line width=.35pt] (-0.45,0.54)--(0.45,0.54);
  \node[anchor=south,inner sep=1pt] at (0,0.76) {$\scriptstyle #1$};
  \draw[densely dashed,line width=.35pt] (-0.45,-0.54)--(0.40,-0.54);
  \node[anchor=north,inner sep=1pt] at (0,-0.76) {$\scriptstyle #2$};
  \node[cplus,draw=black,text=black,scale=.7,inner sep=0pt]
    at (-0.24,-0.05) {$-$};
\end{tikzpicture}}
\noindent\resizebox{\linewidth}{!}{%
$\displaystyle y_{7}={}\cleanPCurl
 {l_{2,1}W_{V3}^{-}(p_V(x_3,-u)-\frac{p_V(x_3,-u)-p_{V31}^{+-}}{u-x_1})}
 {{l_{1,2}}}%
$%
\hspace{8pt}%
$\displaystyle y_{8}={}\begin{tikzpicture}[baseline=(current bounding box.center),
  scale=0.8,line cap=round,line join=round,
  font={\fontsize{10pt}{12pt}\selectfont}]
  \begin{scope}[xscale=-1,yscale=1]
  \draw[-,thick] (0,-0.6) to[out=90, in=0] (-0.2,-0.3);
 \path[overlay,name path global=ndS1] (0,-0.6) to[out=90, in=0] (-0.2,-0.3);
  \draw[-,thick] (-0.2,-0.3) to[out = 180, in = 90] (-0.4,-0.6);
 \path[overlay,name path global=ndS2] (-0.2,-0.3) to[out = 180, in = 90] (-0.4,-0.6);
  \draw[-,thick] (0.4,0.6) to (0.4,-0.6);
 \path[overlay,name path global=ndS3] (0.4,0.6) to (0.4,-0.6);
  \draw[-,thick] (0,0.6) to[out=-90, in=0] (-0.2,0.3);
 \path[overlay,name path global=ndS4] (0,0.6) to[out=-90, in=0] (-0.2,0.3);
  \draw[-,thick] (-0.2,0.3) to[out = -180, in = -90] (-0.4,0.6);
 \path[overlay,name path global=ndS5] (-0.2,0.3) to[out = -180, in = -90] (-0.4,0.6);
  \end{scope}
 \path[overlay,name path=ndG1] (-2,0.54)--(2,0.54);
 \path[overlay,name intersections={of=ndG1 and ndS3,by=ndP1}];
 \path[overlay,name path=ndG2] (-2,0.54)--(2,0.54);
 \path[overlay,name intersections={of=ndG2 and ndS5,by=ndP2}];
 \draw[grsw morphism,line width=0.350000pt] (ndP1)--(ndP2);
 \node[anchor=south,inner sep=1pt] at (0,0.76) {$\scriptstyle l_{2,1}W_{V3}^{-}(p_V(x_3,-u)+\frac{p_{V31}^{+-}-p_V(x_3,-u)}{u-x_1})$};
 \path[overlay,name path=ndG3] (-2,-0.54)--(2,-0.54);
 \path[overlay,name intersections={of=ndG3 and ndS3,by=ndP3}];
 \path[overlay,name path=ndG4] (-2,-0.54)--(2,-0.54);
 \path[overlay,name intersections={of=ndG4 and ndS2,by=ndP4}];
 \draw[grsw morphism,line width=0.350000pt] (ndP3)--(ndP4);
 \node[anchor=north,inner sep=1pt] at (0,-0.76) {$\scriptstyle l_{1,2}$};
 \path[overlay,name path=ndG5] (-2,-0.43)--(2,-0.43);
 \path[overlay,name intersections={of=ndG5 and ndS2,by=ndP5}];
 \node[cplus,scale=.7,inner sep=0pt] at (ndP5) {$-$};
\end{tikzpicture}$%
\hspace{8pt}%
$\displaystyle y_{9}={}\begin{tikzpicture}[baseline=(current bounding box.center),
  scale=0.8,line cap=round,line join=round,
  font={\fontsize{10pt}{12pt}\selectfont}]
  \begin{scope}[xscale=-1,yscale=1]
  \draw[-,thick] (0,-0.6) to[out=90, in=0] (-0.2,-0.3);
 \path[overlay,name path global=ndS1] (0,-0.6) to[out=90, in=0] (-0.2,-0.3);
  \draw[-,thick] (-0.2,-0.3) to[out = 180, in = 90] (-0.4,-0.6);
 \path[overlay,name path global=ndS2] (-0.2,-0.3) to[out = 180, in = 90] (-0.4,-0.6);
  \draw[-,thick] (0.4,0.6) to (0.4,-0.6);
 \path[overlay,name path global=ndS3] (0.4,0.6) to (0.4,-0.6);
  \draw[-,thick] (0,0.6) to[out=-90, in=0] (-0.2,0.3);
 \path[overlay,name path global=ndS4] (0,0.6) to[out=-90, in=0] (-0.2,0.3);
  \draw[-,thick] (-0.2,0.3) to[out = -180, in = -90] (-0.4,0.6);
 \path[overlay,name path global=ndS5] (-0.2,0.3) to[out = -180, in = -90] (-0.4,0.6);
  \end{scope}
 \path[overlay,name path=ndG1] (-2,0.54)--(2,0.54);
 \path[overlay,name intersections={of=ndG1 and ndS3,by=ndP1}];
 \path[overlay,name path=ndG2] (-2,0.54)--(2,0.54);
 \path[overlay,name intersections={of=ndG2 and ndS5,by=ndP2}];
 \draw[grsw morphism,line width=0.350000pt] (ndP1)--(ndP2);
 \node[anchor=south,inner sep=1pt] at (0,0.76) {$\scriptstyle l_{2,1}W_{V3}^{-}(p_V(x_3,-u)+\frac{p_{V31}^{+-}-p_V(x_3,-u)}{u-x_1})$};
 \path[overlay,name path=ndG3] (-2,-0.54)--(2,-0.54);
 \path[overlay,name intersections={of=ndG3 and ndS3,by=ndP3}];
 \path[overlay,name path=ndG4] (-2,-0.54)--(2,-0.54);
 \path[overlay,name intersections={of=ndG4 and ndS2,by=ndP4}];
 \draw[grsw morphism,line width=0.350000pt] (ndP3)--(ndP4);
 \node[anchor=north,inner sep=1pt] at (0,-0.76) {$\scriptstyle l_{1,2}$};
 \path[overlay,name path=ndG5] (-2,-0.43)--(2,-0.43);
 \path[overlay,name intersections={of=ndG5 and ndS2,by=ndP5}];
 \node[cplus,scale=.7,inner sep=0pt] at (ndP5) {$-$};
 \path[overlay,name path=ndG6] (-2,-0.43)--(2,-0.43);
 \path[overlay,name intersections={of=ndG6 and ndS1,by=ndP6}];
 \node[cplus,scale=.7,inner sep=0pt] at (ndP6) {$-$};
\end{tikzpicture}$%
}\par\medskip
\par\medskip

For the terms involving $p_V$, we have

\begin{equation}\label{eq:nb-reduce-p}
\begin{aligned}
\left[\cleanPCurl
 {l_{2,1}W_{V3}^{-}p_V(x_3,-u)}
 {{l_{1,2}}}\right]_{{u:-1}}
{}-\left[\cleanPCurl
 {l_{2,1}W_{V3}^{-}\frac{p_V(x_3,-u)-p_{V31}^{+-}}{u-x_1}}
 {{l_{1,2}}}\right]_{{u:-1}}
&=Y_7\\[4pt]
&=Y_8+Y_9-Y_2.
\end{aligned}
\end{equation}

\par\medskip
\noindent Define\par\smallskip
\noindent\makebox[\linewidth][c]{%
$\displaystyle y_{10}={}\begin{tikzpicture}[baseline=(current bounding box.center),
  scale=0.8,line cap=round,line join=round,
  font={\fontsize{10pt}{12pt}\selectfont}]
  \begin{scope}[xscale=-1,yscale=1]
  \draw[-,thick] (0,-0.6) to[out=90, in=0] (-0.2,-0.3);
 \path[overlay,name path global=ndS1] (0,-0.6) to[out=90, in=0] (-0.2,-0.3);
  \draw[-,thick] (-0.2,-0.3) to[out = 180, in = 90] (-0.4,-0.6);
 \path[overlay,name path global=ndS2] (-0.2,-0.3) to[out = 180, in = 90] (-0.4,-0.6);
  \draw[-,thick] (-0.2,0.2) to[out=180,in=90] (-0.4,0);
 \path[overlay,name path global=ndS3] (-0.2,0.2) to[out=180,in=90] (-0.4,0);
  \draw[-,thick] (0,0) to[out=90,in=0] (-0.2,0.2);
 \path[overlay,name path global=ndS4] (0,0) to[out=90,in=0] (-0.2,0.2);
  \draw[-,thick] (-0.4,0) to[out=-90,in=180] (-0.2,-0.2);
 \path[overlay,name path global=ndS5] (-0.4,0) to[out=-90,in=180] (-0.2,-0.2);
  \draw[-,thick] (-0.2,-0.2) to[out=0,in=-90] (0,0);
 \path[overlay,name path global=ndS6] (-0.2,-0.2) to[out=0,in=-90] (0,0);
  \draw[-,thick] (0.4,0.6) to (0.4,-0.6);
 \path[overlay,name path global=ndS7] (0.4,0.6) to (0.4,-0.6);
  \draw[-,thick] (0,0.6) to[out=-90, in=0] (-0.2,0.3);
 \path[overlay,name path global=ndS8] (0,0.6) to[out=-90, in=0] (-0.2,0.3);
  \draw[-,thick] (-0.2,0.3) to[out = -180, in = -90] (-0.4,0.6);
 \path[overlay,name path global=ndS9] (-0.2,0.3) to[out = -180, in = -90] (-0.4,0.6);
  \end{scope}
 \path[overlay,name path=ndG1] (-2,0.54)--(2,0.54);
 \path[overlay,name intersections={of=ndG1 and ndS7,by=ndP1}];
 \path[overlay,name path=ndG2] (-2,0.54)--(2,0.54);
 \path[overlay,name intersections={of=ndG2 and ndS9,by=ndP2}];
 \draw[grsw morphism,line width=0.350000pt] (ndP1)--(ndP2);
 \node[anchor=south,inner sep=1pt] at (0,0.76) {$\scriptstyle l_{2,1}$};
 \path[overlay,name path=ndG3] (-2,-0.54)--(2,-0.54);
 \path[overlay,name intersections={of=ndG3 and ndS7,by=ndP3}];
 \path[overlay,name path=ndG4] (-2,-0.54)--(2,-0.54);
 \path[overlay,name intersections={of=ndG4 and ndS2,by=ndP4}];
 \draw[grsw morphism,line width=0.350000pt] (ndP3)--(ndP4);
 \node[anchor=north,inner sep=1pt] at (0,-0.76) {$\scriptstyle l_{1,2}l_{1,3}$};
 \path[overlay,name path=ndG5] (0.34,-2)--(0.34,2);
 \path[overlay,name intersections={of=ndG5 and ndS3,by=ndL5}];
 \path[overlay,name intersections={of=ndG5 and ndS9,by=ndR5}];
 \draw[grsw morphism,line width=0.350000pt] (ndL5)--(ndR5);
 \node[anchor=west,inner sep=1pt] at (0.56,0.31) {$\scriptstyle W_{V2}^{-}p_{V21}^{++}$};
 \path[overlay,name path=ndG6] (0.34,-2)--(0.34,2);
 \path[overlay,name intersections={of=ndG6 and ndS2,by=ndL6}];
 \path[overlay,name intersections={of=ndG6 and ndS5,by=ndR6}];
 \draw[grsw morphism,line width=0.350000pt] (ndL6)--(ndR6);
 \node[anchor=west,inner sep=1pt] at (0.56,-0.31) {$\scriptstyle M_{V1}^{-}q_{V21}^{++}$};
 \path[overlay,name path=ndG7] (-2,0.04)--(2,0.04);
 \path[overlay,name intersections={of=ndG7 and ndS4,by=ndP7}];
 \node[cplus,scale=.7,inner sep=0pt] at (ndP7) {$-$};
 \path[overlay,name path=ndG8] (-2,0)--(2,0);
 \path[overlay,name intersections={of=ndG8 and ndS7,by=ndP8}];
 \node[cplus,scale=.7,inner sep=0pt] at (ndP8) {$-$};
 \path[overlay,name path=ndG9] (-2,0.31)--(2,0.31);
 \path[overlay,name intersections={of=ndG9 and ndS7,by=ndP9}];
 \fill (ndP9) circle (1.25pt);
 \node[anchor=east,inner sep=1pt] at (-0.56,0.31) {$\scriptstyle r(x,u)$};
\end{tikzpicture}$%
\hspace{8pt}%
$\displaystyle y_{11}={}\AthroughPurl
 {l_{2,1}}
 {l_{1,2}l_{1,3}r(x_1,u)}
 {W_{V2}^{-}p_{V21}^{++}}
 {M_{V1}^{-}q_{V21}^{++}}
 {}
 {inner}
$%
\hspace{8pt}%
$\displaystyle y_{12}={}\begin{tikzpicture}[baseline=(current bounding box.center),
  scale=0.8,line cap=round,line join=round,
  font={\fontsize{10pt}{12pt}\selectfont}]
  \begin{scope}[xscale=-1,yscale=1]
  \draw[-,thick] (0,-0.6) to[out=90, in=0] (-0.2,-0.3);
 \path[overlay,name path global=ndS1] (0,-0.6) to[out=90, in=0] (-0.2,-0.3);
  \draw[-,thick] (-0.2,-0.3) to[out = 180, in = 90] (-0.4,-0.6);
 \path[overlay,name path global=ndS2] (-0.2,-0.3) to[out = 180, in = 90] (-0.4,-0.6);
  \draw[-,thick] (0.4,0.6) to (0.4,-0.6);
 \path[overlay,name path global=ndS3] (0.4,0.6) to (0.4,-0.6);
  \draw[-,thick] (0,0.6) to[out=-90, in=0] (-0.2,0.3);
 \path[overlay,name path global=ndS4] (0,0.6) to[out=-90, in=0] (-0.2,0.3);
  \draw[-,thick] (-0.2,0.3) to[out = -180, in = -90] (-0.4,0.6);
 \path[overlay,name path global=ndS5] (-0.2,0.3) to[out = -180, in = -90] (-0.4,0.6);
  \end{scope}
 \path[overlay,name path=ndG1] (-2,0.54)--(2,0.54);
 \path[overlay,name intersections={of=ndG1 and ndS3,by=ndP1}];
 \path[overlay,name path=ndG2] (-2,0.54)--(2,0.54);
 \path[overlay,name intersections={of=ndG2 and ndS5,by=ndP2}];
 \draw[grsw morphism,line width=0.350000pt] (ndP1)--(ndP2);
 \node[anchor=south,inner sep=1pt] at (0,0.76) {$\scriptstyle l_{2,1}W_{V3}^{-}p_{V31}^{++}$};
\node[anchor=north,inner sep=1pt] at (0,-0.76) {$\scriptstyle l_{1,2}l_{1,3}M_{V3}^{-}q_{V13}^{++}$};
 \path[overlay,name path=ndG3] (-2,-0.54)--(2,-0.54);
 \path[overlay,name intersections={of=ndG3 and ndS3,by=ndP3}];
 \path[overlay,name path=ndG4] (-2,-0.54)--(2,-0.54);
 \path[overlay,name intersections={of=ndG4 and ndS2,by=ndP4}];
 \draw[grsw morphism,line width=0.350000pt] (ndP3)--(ndP4);
 \path[overlay,name path=ndG5] (-2,0)--(2,0);
 \path[overlay,name intersections={of=ndG5 and ndS3,by=ndP5}];
 \fill (ndP5) circle (1.25pt);
 \node[anchor=east,inner sep=1pt] at (-0.56,0) {$\scriptstyle \frac{r(x,u)}{(u-x)(u+x)}$};
\end{tikzpicture}$%
\hspace{8pt}%
$\displaystyle y_{13}={}\Athroughcurl
 {l_{2,1}W_{V3}^{-}p_V(x_3,-u)}
 {l_{1,2}l_{1,3}M_{V3}^{-}q_V(-u,x_3)}
 {}{}{}
 {inner}
$%
}\par\medskip
\par\medskip

For the mixed terms, we first use\nopagebreak[4]

\begin{equation}\label{eq:nb-mixed-preparation}
\begin{aligned}
\begin{tikzpicture}[baseline=(current bounding box.center),
  scale=0.8,line cap=round,line join=round,
  font={\fontsize{10pt}{12pt}\selectfont}]
  \begin{scope}[xscale=-1,yscale=1]
  \draw[-,thick] (0,-0.6) to[out=90, in=0] (-0.2,-0.3);
 \path[overlay,name path global=ndS1] (0,-0.6) to[out=90, in=0] (-0.2,-0.3);
  \draw[-,thick] (-0.2,-0.3) to[out = 180, in = 90] (-0.4,-0.6);
 \path[overlay,name path global=ndS2] (-0.2,-0.3) to[out = 180, in = 90] (-0.4,-0.6);
  \draw[-,thick] (-0.2,0.2) to[out=180,in=90] (-0.4,0);
 \path[overlay,name path global=ndS3] (-0.2,0.2) to[out=180,in=90] (-0.4,0);
  \draw[-,thick] (0,0) to[out=90,in=0] (-0.2,0.2);
 \path[overlay,name path global=ndS4] (0,0) to[out=90,in=0] (-0.2,0.2);
  \draw[-,thick] (-0.4,0) to[out=-90,in=180] (-0.2,-0.2);
 \path[overlay,name path global=ndS5] (-0.4,0) to[out=-90,in=180] (-0.2,-0.2);
  \draw[-,thick] (-0.2,-0.2) to[out=0,in=-90] (0,0);
 \path[overlay,name path global=ndS6] (-0.2,-0.2) to[out=0,in=-90] (0,0);
  \draw[-,thick] (0.4,0.6) to (0.4,-0.6);
 \path[overlay,name path global=ndS7] (0.4,0.6) to (0.4,-0.6);
  \draw[-,thick] (0,0.6) to[out=-90, in=0] (-0.2,0.3);
 \path[overlay,name path global=ndS8] (0,0.6) to[out=-90, in=0] (-0.2,0.3);
  \draw[-,thick] (-0.2,0.3) to[out = -180, in = -90] (-0.4,0.6);
 \path[overlay,name path global=ndS9] (-0.2,0.3) to[out = -180, in = -90] (-0.4,0.6);
  \end{scope}
 \path[overlay,name path=ndG1] (-2,0.54)--(2,0.54);
 \path[overlay,name intersections={of=ndG1 and ndS7,by=ndP1}];
 \path[overlay,name path=ndG2] (-2,0.54)--(2,0.54);
 \path[overlay,name intersections={of=ndG2 and ndS9,by=ndP2}];
 \draw[grsw morphism,line width=0.350000pt] (ndP1)--(ndP2);
 \node[anchor=south,inner sep=1pt] at (0,0.76) {$\scriptstyle l_{2,1}$};
 \path[overlay,name path=ndG3] (-2,-0.54)--(2,-0.54);
 \path[overlay,name intersections={of=ndG3 and ndS7,by=ndP3}];
 \path[overlay,name path=ndG4] (-2,-0.54)--(2,-0.54);
 \path[overlay,name intersections={of=ndG4 and ndS2,by=ndP4}];
 \draw[grsw morphism,line width=0.350000pt] (ndP3)--(ndP4);
 \node[anchor=north,inner sep=1pt] at (0,-0.76) {$\scriptstyle l_{1,2}l_{1,3}$};
 \path[overlay,name path=ndG5] (0.34,-2)--(0.34,2);
 \path[overlay,name intersections={of=ndG5 and ndS3,by=ndL5}];
 \path[overlay,name intersections={of=ndG5 and ndS9,by=ndR5}];
 \draw[grsw morphism,line width=0.350000pt] (ndL5)--(ndR5);
 \node[anchor=west,inner sep=1pt] at (0.56,0.31) {$\scriptstyle W_{V2}^{-}p_{V21}^{++}$};
 \path[overlay,name path=ndG6] (0.34,-2)--(0.34,2);
 \path[overlay,name intersections={of=ndG6 and ndS2,by=ndL6}];
 \path[overlay,name intersections={of=ndG6 and ndS5,by=ndR6}];
 \draw[grsw morphism,line width=0.350000pt] (ndL6)--(ndR6);
 \node[anchor=west,inner sep=1pt] at (0.56,-0.31) {$\scriptstyle M_{V1}^{-}q_{V21}^{++}$};
 \path[overlay,name path=ndG7] (-2,0.04)--(2,0.04);
 \path[overlay,name intersections={of=ndG7 and ndS7,by=ndP7}];
 \path[overlay,name path=ndG8] (-2,0.04)--(2,0.04);
 \path[overlay,name intersections={of=ndG8 and ndS4,by=ndP8}];
 \draw[grsw morphism,line width=0.350000pt] (ndP7)--(ndP8);
 \node[anchor=east,inner sep=1pt] at (-0.56,0.04) {$\scriptstyle \frac{1-r(x_1,x_2)}{x_1-x_2}$};
\end{tikzpicture}  +Y_{11}
&= Y_{10} + Y_{11}\\[4pt]
  &= Y_{12} + Y_{13}
\end{aligned}
\end{equation}
\par\medskip
\noindent Define\par\smallskip
\newcommand{\cleanQCurl}[2]{%
\begin{tikzpicture}[baseline=(current bounding box.center),scale=0.8,
  line cap=round,line join=round,font={\fontsize{10pt}{12pt}\selectfont}]
  \draw[thick] (-0.45,0.62)--(-0.45,-0.62);
  \draw[thick] (-0.05,-0.62) to[out=90,in=180] (0.20,-0.42)
                    to[out=0,in=90] (0.45,-0.62);
  \draw[thick] (0,0.62)
    .. controls (0,0.34) and (0.42,0.38) .. (0.31,0.05)
    .. controls (0.24,-0.17) and (-0.18,-0.17) .. (-0.24,0.05)
    .. controls (-0.33,0.38) and (0.40,0.34) .. (0.40,0.62);
  \draw[densely dashed,line width=.35pt] (-0.45,0.54)--(0.40,0.54);
  \node[anchor=south,inner sep=1pt] at (0,0.76) {$\scriptstyle #1$};
  \draw[densely dashed,line width=.35pt] (-0.45,-0.54)--(0.45,-0.54);
  \node[anchor=north,inner sep=1pt] at (0,-0.76) {$\scriptstyle #2$};
  \node[cplus,scale=.7,inner sep=0pt] at (-0.24,0.05) {$-$};
\end{tikzpicture}}
\noindent\makebox[\linewidth][c]{%
$\displaystyle y_{14}={}{\cleanQCurl
 {{l_{2,1}}}
 {l_{1,2}l_{1,3}M_{V3}^{-}\frac{q_V(-u,x_3)-q_{V13}^{-+}}{x_1-u}}}$%
\hspace{8pt}%
$\displaystyle y_{15}={}\begin{tikzpicture}[baseline=(current bounding box.center),
  scale=0.8,line cap=round,line join=round,
  font={\fontsize{10pt}{12pt}\selectfont}]
  \begin{scope}[xscale=-1,yscale=1]
  \draw[-,thick] (0,-0.6) to[out=90, in=0] (-0.2,-0.3);
 \path[overlay,name path global=ndS1] (0,-0.6) to[out=90, in=0] (-0.2,-0.3);
  \draw[-,thick] (-0.2,-0.3) to[out = 180, in = 90] (-0.4,-0.6);
 \path[overlay,name path global=ndS2] (-0.2,-0.3) to[out = 180, in = 90] (-0.4,-0.6);
  \draw[-,thick] (0.4,0.6) to (0.4,-0.6);
 \path[overlay,name path global=ndS3] (0.4,0.6) to (0.4,-0.6);
  \draw[-,thick] (0,0.6) to[out=-90, in=0] (-0.2,0.3);
 \path[overlay,name path global=ndS4] (0,0.6) to[out=-90, in=0] (-0.2,0.3);
  \draw[-,thick] (-0.2,0.3) to[out = -180, in = -90] (-0.4,0.6);
 \path[overlay,name path global=ndS5] (-0.2,0.3) to[out = -180, in = -90] (-0.4,0.6);
  \end{scope}
 \path[overlay,name path=ndG1] (-2,0.54)--(2,0.54);
 \path[overlay,name intersections={of=ndG1 and ndS3,by=ndP1}];
 \path[overlay,name path=ndG2] (-2,0.54)--(2,0.54);
 \path[overlay,name intersections={of=ndG2 and ndS5,by=ndP2}];
 \draw[grsw morphism,line width=0.350000pt] (ndP1)--(ndP2);
 \node[anchor=south,inner sep=1pt] at (0,0.76) {$\scriptstyle l_{2,1}$};
 \path[overlay,name path=ndG3] (-2,-0.54)--(2,-0.54);
 \path[overlay,name intersections={of=ndG3 and ndS3,by=ndP3}];
 \path[overlay,name path=ndG4] (-2,-0.54)--(2,-0.54);
 \path[overlay,name intersections={of=ndG4 and ndS2,by=ndP4}];
 \draw[grsw morphism,line width=0.350000pt] (ndP3)--(ndP4);
 \node[anchor=north,inner sep=1pt] at (0,-0.76) {$\scriptstyle l_{1,2}l_{1,3}M_{V3}^{-}\frac{q_V(-u,x_3)-q_{V13}^{-+}}{-u+x_1}$};
 \path[overlay,name path=ndG5] (-2,0.43)--(2,0.43);
 \path[overlay,name intersections={of=ndG5 and ndS5,by=ndP5}];
 \node[cplus,scale=.7,inner sep=0pt] at (ndP5) {$-$};
\end{tikzpicture}$%
\hspace{8pt}%
$\displaystyle y_{16}={}\begin{tikzpicture}[baseline=(current bounding box.center),
  scale=0.8,line cap=round,line join=round,
  font={\fontsize{10pt}{12pt}\selectfont}]
  \begin{scope}[xscale=-1,yscale=1]
  \draw[-,thick] (0,-0.6) to[out=90, in=0] (-0.2,-0.3);
 \path[overlay,name path global=ndS1] (0,-0.6) to[out=90, in=0] (-0.2,-0.3);
  \draw[-,thick] (-0.2,-0.3) to[out = 180, in = 90] (-0.4,-0.6);
 \path[overlay,name path global=ndS2] (-0.2,-0.3) to[out = 180, in = 90] (-0.4,-0.6);
  \draw[-,thick] (0.4,0.6) to (0.4,-0.6);
 \path[overlay,name path global=ndS3] (0.4,0.6) to (0.4,-0.6);
  \draw[-,thick] (0,0.6) to[out=-90, in=0] (-0.2,0.3);
 \path[overlay,name path global=ndS4] (0,0.6) to[out=-90, in=0] (-0.2,0.3);
  \draw[-,thick] (-0.2,0.3) to[out = -180, in = -90] (-0.4,0.6);
 \path[overlay,name path global=ndS5] (-0.2,0.3) to[out = -180, in = -90] (-0.4,0.6);
  \end{scope}
 \path[overlay,name path=ndG1] (-2,0.54)--(2,0.54);
 \path[overlay,name intersections={of=ndG1 and ndS3,by=ndP1}];
 \path[overlay,name path=ndG2] (-2,0.54)--(2,0.54);
 \path[overlay,name intersections={of=ndG2 and ndS5,by=ndP2}];
 \draw[grsw morphism,line width=0.350000pt] (ndP1)--(ndP2);
 \node[anchor=south,inner sep=1pt] at (0,0.76) {$\scriptstyle l_{2,1}$};
 \path[overlay,name path=ndG3] (-2,-0.54)--(2,-0.54);
 \path[overlay,name intersections={of=ndG3 and ndS3,by=ndP3}];
 \path[overlay,name path=ndG4] (-2,-0.54)--(2,-0.54);
 \path[overlay,name intersections={of=ndG4 and ndS2,by=ndP4}];
 \draw[grsw morphism,line width=0.350000pt] (ndP3)--(ndP4);
 \node[anchor=north,inner sep=1pt] at (0,-0.76) {$\scriptstyle l_{1,2}l_{1,3}M_{V3}^{-}\frac{q_V(-u,x_3)-q_{V13}^{-+}}{-u+x_1}$};
 \path[overlay,name path=ndG5] (-2,0.43)--(2,0.43);
 \path[overlay,name intersections={of=ndG5 and ndS5,by=ndP5}];
 \node[cplus,scale=.7,inner sep=0pt] at (ndP5) {$-$};
 \path[overlay,name path=ndG6] (-2,0.43)--(2,0.43);
 \path[overlay,name intersections={of=ndG6 and ndS4,by=ndP6}];
 \node[cplus,scale=.7,inner sep=0pt] at (ndP6) {$-$};
\end{tikzpicture}$%
}\par\medskip
\par\medskip

{For $q_V$, use $r(x_1,u)=(1+x_1-u)^{-1}$ to combine the
two decorations:
\[
 q_V(-u,x_3)r(x_1,u)
 -\frac{q_V(-u,x_3)r(x_1,u)-q_{V13}^{-+}}{u-x_1}
 =\frac{q_V(-u,x_3)-q_{V13}^{-+}}{x_1-u}.
\]
The upper curl is retained until the second equality below:}

\begin{equation}\label{eq:nb-reduce-q}
\begin{aligned}
\left[ {\cleanQCurl
 {{l_{2,1}}}
 {l_{1,2}l_{1,3}M_{V3}^{-}q_V(-u,x_3)r(x_1,u)}}\right]_{{u:-1}}- \left[ {\cleanQCurl
 {{l_{2,1}}}
 {l_{1,2}l_{1,3}M_{V3}^{-}\frac{q_V(-u,x_3)r(x_1,u)-q_{V13}^{-+}}{u-x_1}}}\right]_{{u:-1}}
  &= Y_{14}\\[4pt]
  &= Y_{15} + Y_{16} -Y_{1}
\end{aligned}
\end{equation}
\par\medskip
\noindent Define\par\smallskip
\noindent\makebox[\linewidth][c]{%
$\displaystyle y_{17}={}\begin{tikzpicture}[baseline=(current bounding box.center),
  scale=0.8,line cap=round,line join=round,
  font={\fontsize{10pt}{12pt}\selectfont}]
  \begin{scope}[xscale=-1,yscale=1]
  \draw[-,thick] (0,-0.6) to[out=90, in=0] (-0.2,-0.3);
 \path[overlay,name path global=ndS1] (0,-0.6) to[out=90, in=0] (-0.2,-0.3);
  \draw[-,thick] (-0.2,-0.3) to[out = 180, in = 90] (-0.4,-0.6);
 \path[overlay,name path global=ndS2] (-0.2,-0.3) to[out = 180, in = 90] (-0.4,-0.6);
  \draw[-,thick] (-0.2,0.2) to[out=180,in=90] (-0.4,0);
 \path[overlay,name path global=ndS3] (-0.2,0.2) to[out=180,in=90] (-0.4,0);
  \draw[-,thick] (0,0) to[out=90,in=0] (-0.2,0.2);
 \path[overlay,name path global=ndS4] (0,0) to[out=90,in=0] (-0.2,0.2);
  \draw[-,thick] (-0.4,0) to[out=-90,in=180] (-0.2,-0.2);
 \path[overlay,name path global=ndS5] (-0.4,0) to[out=-90,in=180] (-0.2,-0.2);
  \draw[-,thick] (-0.2,-0.2) to[out=0,in=-90] (0,0);
 \path[overlay,name path global=ndS6] (-0.2,-0.2) to[out=0,in=-90] (0,0);
  \draw[-,thick] (0.4,0.6) to (0.4,-0.6);
 \path[overlay,name path global=ndS7] (0.4,0.6) to (0.4,-0.6);
  \draw[-,thick] (0,0.6) to[out=-90, in=0] (-0.2,0.3);
 \path[overlay,name path global=ndS8] (0,0.6) to[out=-90, in=0] (-0.2,0.3);
  \draw[-,thick] (-0.2,0.3) to[out = -180, in = -90] (-0.4,0.6);
 \path[overlay,name path global=ndS9] (-0.2,0.3) to[out = -180, in = -90] (-0.4,0.6);
  \end{scope}
 \path[overlay,name path=ndG1] (-2,0.54)--(2,0.54);
 \path[overlay,name intersections={of=ndG1 and ndS7,by=ndP1}];
 \path[overlay,name path=ndG2] (-2,0.54)--(2,0.54);
 \path[overlay,name intersections={of=ndG2 and ndS9,by=ndP2}];
 \draw[grsw morphism,line width=0.350000pt] (ndP1)--(ndP2);
 \node[anchor=south,inner sep=1pt] at (0,0.76) {$\scriptstyle l_{2,1}$};
 \path[overlay,name path=ndG3] (-2,-0.54)--(2,-0.54);
 \path[overlay,name intersections={of=ndG3 and ndS7,by=ndP3}];
 \path[overlay,name path=ndG4] (-2,-0.54)--(2,-0.54);
 \path[overlay,name intersections={of=ndG4 and ndS2,by=ndP4}];
 \draw[grsw morphism,line width=0.350000pt] (ndP3)--(ndP4);
 \node[anchor=north,inner sep=1pt] at (0,-0.76) {$\scriptstyle l_{1,2}l_{1,3}\frac{1}{u-x_3}(\frac{1-r(x_1,u)}{x_1-u}-\frac{1-r(x_1,x_3)}{x_1-x_3})$};
 \path[overlay,name path=ndG5] (0.34,-2)--(0.34,2);
 \path[overlay,name intersections={of=ndG5 and ndS3,by=ndL5}];
 \path[overlay,name intersections={of=ndG5 and ndS9,by=ndR5}];
 \draw[grsw morphism,line width=0.350000pt] (ndL5)--(ndR5);
 \node[anchor=west,inner sep=1pt] at (0.56,0.26) {$\scriptstyle W_{V2}^{-}p_{V21}^{++}$};
 \path[overlay,name path=ndG6] (-2,0.04)--(2,0.04);
 \path[overlay,name intersections={of=ndG6 and ndS4,by=ndP6}];
 \node[cplus,scale=.7,inner sep=0pt] at (ndP6) {$-$};
\end{tikzpicture}$%
\hspace{8pt}%
$\displaystyle y_{18}={}\AthroughQurl
 {l_{2,1}}
 {l_{1,2}l_{1,3}\dfrac{r(x_1,u)-r(x_1,x_3)}{x_3-u}}
 {W_{V2}^{-}p_{V21}^{++}}
 {}
 {}
 {inner}
$%
\hspace{8pt}%
$\displaystyle y_{19}={}\begin{tikzpicture}[baseline=(current bounding box.center),
  scale=0.8,line cap=round,line join=round,
  font={\fontsize{10pt}{12pt}\selectfont}]
  \begin{scope}[xscale=-1,yscale=1]
  \draw[-,thick] (0,-0.6) to[out=90, in=0] (-0.2,-0.3);
 \path[overlay,name path global=ndS1] (0,-0.6) to[out=90, in=0] (-0.2,-0.3);
  \draw[-,thick] (-0.2,-0.3) to[out = 180, in = 90] (-0.4,-0.6);
 \path[overlay,name path global=ndS2] (-0.2,-0.3) to[out = 180, in = 90] (-0.4,-0.6);
  \draw[-,thick] (0.4,0.6) to (0.4,-0.6);
 \path[overlay,name path global=ndS3] (0.4,0.6) to (0.4,-0.6);
  \draw[-,thick] (0,0.6) to[out=-90, in=0] (-0.2,0.3);
 \path[overlay,name path global=ndS4] (0,0.6) to[out=-90, in=0] (-0.2,0.3);
  \draw[-,thick] (-0.2,0.3) to[out = -180, in = -90] (-0.4,0.6);
 \path[overlay,name path global=ndS5] (-0.2,0.3) to[out = -180, in = -90] (-0.4,0.6);
  \end{scope}
 \path[overlay,name path=ndG1] (-2,0.54)--(2,0.54);
 \path[overlay,name intersections={of=ndG1 and ndS3,by=ndP1}];
 \path[overlay,name path=ndG2] (-2,0.54)--(2,0.54);
 \path[overlay,name intersections={of=ndG2 and ndS5,by=ndP2}];
 \draw[grsw morphism,line width=0.350000pt] (ndP1)--(ndP2);\node[anchor=south,inner sep=1pt] at (0,0.76) {$\scriptstyle l_{2,1}W_{V3}^{-}p_{V31}^{++}$};
 \path[overlay,name path=ndG3] (-2,-0.54)--(2,-0.54);
 \path[overlay,name intersections={of=ndG3 and ndS3,by=ndP3}];
 \path[overlay,name path=ndG4] (-2,-0.54)--(2,-0.54);
 \path[overlay,name intersections={of=ndG4 and ndS2,by=ndP4}];
 \draw[grsw morphism,line width=0.350000pt] (ndP3)--(ndP4);
 \node[anchor=north,inner sep=1pt] at (0,-0.76) {$\scriptstyle l_{1,2}l_{1,3}\frac{r(x_1,u)-r(x_1,x_3)}{x_3-u}$};
 \path[overlay,name path=ndG5] (-2,0)--(2,0);
 \path[overlay,name intersections={of=ndG5 and ndS3,by=ndP5}];
 \fill (ndP5) circle (1.25pt);
 \node[anchor=east,inner sep=1pt] at (-0.56,0) {$\scriptstyle \frac1{(u-x)(u+x)}$};
\end{tikzpicture}$%
}\par\medskip
\noindent\makebox[\linewidth][c]{%
$\displaystyle y_{20}={}
\AthroughQurl
 {l_{2,1}}
 {l_{1,2}l_{1,3}\dfrac{r(x_1,u)-r(x_1,x_3)}{x_3-u}}
 {W_{V2}^{-}p_{V21}^{++}}
 {}
 {}
 {upper}
$%
}\par\medskip
\par\medskip

The remaining difference terms satisfy

\begin{equation}\label{eq:nb-difference-preparation}
\begin{aligned}
- \left[ \begin{tikzpicture}[baseline=(current bounding box.center),
  scale=0.8,line cap=round,line join=round,
  font={\fontsize{10pt}{12pt}\selectfont}]
  \begin{scope}[xscale=-1,yscale=1]
  \draw[-,thick] (0,-0.6) to[out=90, in=0] (-0.2,-0.3);
 \path[overlay,name path global=ndS1] (0,-0.6) to[out=90, in=0] (-0.2,-0.3);
  \draw[-,thick] (-0.2,-0.3) to[out = 180, in = 90] (-0.4,-0.6);
 \path[overlay,name path global=ndS2] (-0.2,-0.3) to[out = 180, in = 90] (-0.4,-0.6);
  \draw[-,thick] (-0.2,0.2) to[out=180,in=90] (-0.4,0);
 \path[overlay,name path global=ndS3] (-0.2,0.2) to[out=180,in=90] (-0.4,0);
  \draw[-,thick] (0,0) to[out=90,in=0] (-0.2,0.2);
 \path[overlay,name path global=ndS4] (0,0) to[out=90,in=0] (-0.2,0.2);
  \draw[-,thick] (-0.4,0) to[out=-90,in=180] (-0.2,-0.2);
 \path[overlay,name path global=ndS5] (-0.4,0) to[out=-90,in=180] (-0.2,-0.2);
  \draw[-,thick] (-0.2,-0.2) to[out=0,in=-90] (0,0);
 \path[overlay,name path global=ndS6] (-0.2,-0.2) to[out=0,in=-90] (0,0);
  \draw[-,thick] (0.4,0.6) to (0.4,-0.6);
 \path[overlay,name path global=ndS7] (0.4,0.6) to (0.4,-0.6);
  \draw[-,thick] (0,0.6) to[out=-90, in=0] (-0.2,0.3);
 \path[overlay,name path global=ndS8] (0,0.6) to[out=-90, in=0] (-0.2,0.3);
  \draw[-,thick] (-0.2,0.3) to[out = -180, in = -90] (-0.4,0.6);
 \path[overlay,name path global=ndS9] (-0.2,0.3) to[out = -180, in = -90] (-0.4,0.6);
  \end{scope}
 \path[overlay,name path=ndG1] (-2,0.54)--(2,0.54);
 \path[overlay,name intersections={of=ndG1 and ndS7,by=ndP1}];
 \path[overlay,name path=ndG2] (-2,0.54)--(2,0.54);
 \path[overlay,name intersections={of=ndG2 and ndS9,by=ndP2}];
 \draw[grsw morphism,line width=0.350000pt] (ndP1)--(ndP2);
 \node[anchor=south,inner sep=1pt] at (0,0.76) {$\scriptstyle l_{2,1}$};
 \path[overlay,name path=ndG3] (-2,-0.54)--(2,-0.54);
 \path[overlay,name intersections={of=ndG3 and ndS7,by=ndP3}];
 \path[overlay,name path=ndG4] (-2,-0.54)--(2,-0.54);
 \path[overlay,name intersections={of=ndG4 and ndS2,by=ndP4}];
 \draw[grsw morphism,line width=0.350000pt] (ndP3)--(ndP4);
 \node[anchor=north,inner sep=1pt] at (0,-0.76) {$\scriptstyle l_{1,2}l_{1,3}$};
 \path[overlay,name path=ndG5] (0.34,-2)--(0.34,2);
 \path[overlay,name intersections={of=ndG5 and ndS3,by=ndL5}];
 \path[overlay,name intersections={of=ndG5 and ndS9,by=ndR5}];
 \draw[grsw morphism,line width=0.350000pt] (ndL5)--(ndR5);
 \node[anchor=west,inner sep=1pt] at (0.56,0.26) {$\scriptstyle W_{V2}^{-}p_{V21}^{++}$};
 \path[overlay,name path=ndG6] (-2,0.04)--(2,0.04);
 \path[overlay,name intersections={of=ndG6 and ndS4,by=ndP6}];
 \node[cplus,scale=.7,inner sep=0pt] at (ndP6) {$-$};
 \path[overlay,name path=ndG7] (-2,-0.43)--(2,-0.43);
 \path[overlay,name intersections={of=ndG7 and ndS2,by=ndP7}];
 \node[cplus,scale=.7,inner sep=0pt] at (ndP7) {$-$};
 \path[overlay,name path=ndG8] (-2,-0.06)--(2,-0.06);
 \path[overlay,name intersections={of=ndG8 and ndS7,by=ndP8}];
 \fill (ndP8) circle (1.25pt);
 \node[anchor=east,inner sep=1pt] at (-0.56,-0.06) {$\scriptstyle \frac{1-r(x,u)}{x-u}$};
\end{tikzpicture}\right]_{{u:-1}} -\left[
\AthroughQurl
 {l_{2,1}}
 {l_{1,2}l_{1,3}r(x_1,u)(u-x_3)^{-1}}
 {W_{V2}^{-}p_{V21}^{++}}
 {}
 {}
 {inner}
\right]_{{u:-1}}
  &= -Y_{17}+Y_{18}\\[4pt]
  &= Y_{19}+Y_{20}
\end{aligned}
\end{equation}
\par\medskip
\noindent Define\par\smallskip
\noindent\makebox[\linewidth][c]{%
$\displaystyle y_{21}={}\begin{tikzpicture}[baseline=(current bounding box.center),
  scale=0.8,line cap=round,line join=round,
  font={\fontsize{10pt}{12pt}\selectfont}]
  \begin{scope}[xscale=-1,yscale=1]
  \draw[-,thick] (0,-0.6) to[out=90, in=0] (-0.2,-0.3);
 \path[overlay,name path global=ndS1] (0,-0.6) to[out=90, in=0] (-0.2,-0.3);
  \draw[-,thick] (-0.2,-0.3) to[out = 180, in = 90] (-0.4,-0.6);
 \path[overlay,name path global=ndS2] (-0.2,-0.3) to[out = 180, in = 90] (-0.4,-0.6);
  \draw[-,thick] (0.4,0.6) to (0.4,-0.6);
 \path[overlay,name path global=ndS3] (0.4,0.6) to (0.4,-0.6);
  \draw[-,thick] (0,0.6) to[out=-90, in=0] (-0.2,0.3);
 \path[overlay,name path global=ndS4] (0,0.6) to[out=-90, in=0] (-0.2,0.3);
  \draw[-,thick] (-0.2,0.3) to[out = -180, in = -90] (-0.4,0.6);
 \path[overlay,name path global=ndS5] (-0.2,0.3) to[out = -180, in = -90] (-0.4,0.6);
  \end{scope}
 \path[overlay,name path=ndG1] (-2,0.54)--(2,0.54);
 \path[overlay,name intersections={of=ndG1 and ndS3,by=ndP1}];
 \path[overlay,name path=ndG2] (-2,0.54)--(2,0.54);
 \path[overlay,name intersections={of=ndG2 and ndS5,by=ndP2}];
 \draw[grsw morphism,line width=0.350000pt] (ndP1)--(ndP2);
 \node[anchor=south,inner sep=1pt] at (0,0.76) {$\scriptstyle l_{2,1}W_{V3}^{-}p_V(x_3,-u)$};
 \path[overlay,name path=ndG3] (-2,-0.54)--(2,-0.54);
 \path[overlay,name intersections={of=ndG3 and ndS3,by=ndP3}];
 \path[overlay,name path=ndG4] (-2,-0.54)--(2,-0.54);
 \path[overlay,name intersections={of=ndG4 and ndS2,by=ndP4}];
 \draw[grsw morphism,line width=0.350000pt] (ndP3)--(ndP4);
 \node[anchor=north,inner sep=1pt] at (0,-0.76) {$\scriptstyle l_{1,2}l_{1,3}M_{V3}^{-}q_V(-u,x_3)$};
 \path[overlay,name path=ndG5] (-2,0)--(2,0);
 \path[overlay,name intersections={of=ndG5 and ndS3,by=ndP5}];
 \fill (ndP5) circle (1.25pt);
 \node[anchor=east,inner sep=1pt] at (-0.56,0) {$\scriptstyle \frac1{(u-x)(u+x)}$};
\end{tikzpicture}$%
\hspace{8pt}%
$\displaystyle y_{22}={}\begin{tikzpicture}[baseline=(current bounding box.center),
  scale=0.8,line cap=round,line join=round,
  font={\fontsize{10pt}{12pt}\selectfont}]
  \begin{scope}[xscale=-1,yscale=1]
  \draw[-,thick] (0,-0.6) to[out=90, in=0] (-0.2,-0.3);
 \path[overlay,name path global=ndS1] (0,-0.6) to[out=90, in=0] (-0.2,-0.3);
  \draw[-,thick] (-0.2,-0.3) to[out = 180, in = 90] (-0.4,-0.6);
 \path[overlay,name path global=ndS2] (-0.2,-0.3) to[out = 180, in = 90] (-0.4,-0.6);
  \draw[-,thick] (0.4,0.6) to (0.4,-0.6);
 \path[overlay,name path global=ndS3] (0.4,0.6) to (0.4,-0.6);
  \draw[-,thick] (0,0.6) to[out=-90, in=0] (-0.2,0.3);
 \path[overlay,name path global=ndS4] (0,0.6) to[out=-90, in=0] (-0.2,0.3);
  \draw[-,thick] (-0.2,0.3) to[out = -180, in = -90] (-0.4,0.6);
 \path[overlay,name path global=ndS5] (-0.2,0.3) to[out = -180, in = -90] (-0.4,0.6);
  \end{scope}
 \path[overlay,name path=ndG1] (-2,0.54)--(2,0.54);
 \path[overlay,name intersections={of=ndG1 and ndS3,by=ndP1}];
 \path[overlay,name path=ndG2] (-2,0.54)--(2,0.54);
 \path[overlay,name intersections={of=ndG2 and ndS5,by=ndP2}];
 \draw[grsw morphism,line width=0.350000pt] (ndP1)--(ndP2);
 \node[anchor=south,inner sep=1pt] at (0,0.76) {$\scriptstyle l_{2,1}W_{V3}^{-}p_V(x_3,-u)$};
 \path[overlay,name path=ndG3] (-2,-0.54)--(2,-0.54);
 \path[overlay,name intersections={of=ndG3 and ndS3,by=ndP3}];
 \path[overlay,name path=ndG4] (-2,-0.54)--(2,-0.54);
 \path[overlay,name intersections={of=ndG4 and ndS2,by=ndP4}];
 \draw[grsw morphism,line width=0.350000pt] (ndP3)--(ndP4);
 \node[anchor=north,inner sep=1pt] at (0,-0.76) {$\scriptstyle l_{1,2}l_{1,3}M_{V3}^{-}q_V(-u,x_3)$};
 \path[overlay,name path=ndG5] (-2,0)--(2,0);
 \path[overlay,name intersections={of=ndG5 and ndS3,by=ndP5}];
 \node[cplus,scale=.7,inner sep=0pt] at (ndP5) {$-$};
\end{tikzpicture}$%
}\par\medskip
\par\medskip

The last term in \eqref{eq:nb-mixed-preparation} reduces as follows:

\begin{equation}\label{eq:nb-reduce-pq}
\begin{aligned}
\left[
\Athroughcurl
 {l_{2,1}W_{V3}^{-}p_V(x_3,-u)}
 {l_{1,2}l_{1,3}M_{V3}^{-}q_V(-u,x_3)}
 {}{}{}
 {inner}
\right]_{{u:-1}}
  &=Y_{13}\\[4pt]
  &=-Y_3+Y_{21}+Y_{22}.
\end{aligned}
\end{equation}

\par\medskip
\noindent Define\par\smallskip
\noindent\makebox[\linewidth][c]{%
$\displaystyle y_{23}={}\begin{tikzpicture}[baseline=(current bounding box.center),
  scale=0.8,line cap=round,line join=round,
  font={\fontsize{10pt}{12pt}\selectfont}]
  \begin{scope}[xscale=-1,yscale=1]
  \draw[-,thick] (0,-0.6) to[out=90, in=0] (-0.2,-0.3);
 \path[overlay,name path global=ndS1] (0,-0.6) to[out=90, in=0] (-0.2,-0.3);
  \draw[-,thick] (-0.2,-0.3) to[out = 180, in = 90] (-0.4,-0.6);
 \path[overlay,name path global=ndS2] (-0.2,-0.3) to[out = 180, in = 90] (-0.4,-0.6);
  \draw[-,thick] (0.4,0.6) to (0.4,-0.6);
 \path[overlay,name path global=ndS3] (0.4,0.6) to (0.4,-0.6);
  \draw[-,thick] (0,0.6) to[out=-90, in=0] (-0.2,0.3);
 \path[overlay,name path global=ndS4] (0,0.6) to[out=-90, in=0] (-0.2,0.3);
  \draw[-,thick] (-0.2,0.3) to[out = -180, in = -90] (-0.4,0.6);
 \path[overlay,name path global=ndS5] (-0.2,0.3) to[out = -180, in = -90] (-0.4,0.6);
  \end{scope}
 \path[overlay,name path=ndG1] (-2,0.54)--(2,0.54);
 \path[overlay,name intersections={of=ndG1 and ndS3,by=ndP1}];
 \path[overlay,name path=ndG2] (-2,0.54)--(2,0.54);
 \path[overlay,name intersections={of=ndG2 and ndS5,by=ndP2}];
 \draw[grsw morphism,line width=0.350000pt] (ndP1)--(ndP2);
 \node[anchor=south,inner sep=1pt] at (0,0.76) {$\scriptstyle l_{2,1}W_{V3}^{-}p_V(x_3,-u)$};
 \path[overlay,name path=ndG3] (-2,-0.54)--(2,-0.54);
 \path[overlay,name intersections={of=ndG3 and ndS3,by=ndP3}];
 \path[overlay,name path=ndG4] (-2,-0.54)--(2,-0.54);
 \path[overlay,name intersections={of=ndG4 and ndS2,by=ndP4}];
 \draw[grsw morphism,line width=0.350000pt] (ndP3)--(ndP4);
 \node[anchor=north,inner sep=1pt] at (0,-0.76) {$\scriptstyle l_{1,2}l_{1,3}\frac{r(x_1,x_3)-1}{x_1-x_3}$};
 \path[overlay,name path=ndG5] (-2,0)--(2,0);
 \path[overlay,name intersections={of=ndG5 and ndS3,by=ndP5}];
 \fill (ndP5) circle (1.25pt);
 \node[anchor=east,inner sep=1pt] at (-0.56,0) {$\scriptstyle \frac1{(u-x)(u+x)}$};
\end{tikzpicture}$%
\hspace{8pt}%
$\displaystyle y_{24}={}\begin{tikzpicture}[baseline=(current bounding box.center),
  scale=0.8,line cap=round,line join=round,
  font={\fontsize{10pt}{12pt}\selectfont}]
  \begin{scope}[xscale=-1,yscale=1]
  \draw[-,thick] (0,-0.6) to[out=90, in=0] (-0.2,-0.3);
 \path[overlay,name path global=ndS1] (0,-0.6) to[out=90, in=0] (-0.2,-0.3);
  \draw[-,thick] (-0.2,-0.3) to[out = 180, in = 90] (-0.4,-0.6);
 \path[overlay,name path global=ndS2] (-0.2,-0.3) to[out = 180, in = 90] (-0.4,-0.6);
  \draw[-,thick] (0.4,0.6) to (0.4,-0.6);
 \path[overlay,name path global=ndS3] (0.4,0.6) to (0.4,-0.6);
  \draw[-,thick] (0,0.6) to[out=-90, in=0] (-0.2,0.3);
 \path[overlay,name path global=ndS4] (0,0.6) to[out=-90, in=0] (-0.2,0.3);
  \draw[-,thick] (-0.2,0.3) to[out = -180, in = -90] (-0.4,0.6);
 \path[overlay,name path global=ndS5] (-0.2,0.3) to[out = -180, in = -90] (-0.4,0.6);
  \end{scope}
 \path[overlay,name path=ndG1] (-2,0.54)--(2,0.54);
 \path[overlay,name intersections={of=ndG1 and ndS3,by=ndP1}];
 \path[overlay,name path=ndG2] (-2,0.54)--(2,0.54);
 \path[overlay,name intersections={of=ndG2 and ndS5,by=ndP2}];
 \draw[grsw morphism,line width=0.350000pt] (ndP1)--(ndP2);
 \node[anchor=south,inner sep=1pt] at (0,0.76) {$\scriptstyle l_{2,1}W_{V3}^{-}p_V(x_3,-u)$};
 \path[overlay,name path=ndG3] (-2,-0.54)--(2,-0.54);
 \path[overlay,name intersections={of=ndG3 and ndS3,by=ndP3}];
 \path[overlay,name path=ndG4] (-2,-0.54)--(2,-0.54);
 \path[overlay,name intersections={of=ndG4 and ndS2,by=ndP4}];
 \draw[grsw morphism,line width=0.350000pt] (ndP3)--(ndP4);
 \node[anchor=north,inner sep=1pt] at (0,-0.76) {$\scriptstyle l_{1,2}l_{1,3}\frac{r(x_1,x_3)-1}{x_1-x_3}$};
 \path[overlay,name path=ndG5] (-2,0)--(2,0);
 \path[overlay,name intersections={of=ndG5 and ndS3,by=ndP5}];
 \node[cplus,scale=.7,inner sep=0pt] at (ndP5) {$-$};
\end{tikzpicture}$%
}\par\medskip
\par\medskip

The last term in \eqref{eq:nb-difference-preparation} similarly gives

\begin{equation}\label{eq:nb-reduce-r}
\left[
\Athroughcurl
 {l_{2,1}W_{V3}^{-}p_V(x_3,-u)}
 {l_{1,2}l_{1,3}\dfrac{r(x_1,x_3)-1}{x_1-x_3}}
 {}{}{}
 {inner}
\right]_{{u:-1}}
=-Y_4+Y_{23}+Y_{24}.
\end{equation}

\renewcommand{\qedsymbol}{}
The bubble terms in \eqref{eq:nb-reduce-q}, \eqref{eq:nb-reduce-p},
\eqref{eq:nb-reduce-pq}, and \eqref{eq:nb-reduce-r} are the four terms
listed above. We now evaluate them.

For the $q_V$ difference quotient, we obtain

\begin{equation}\label{eq:nb-eval-q}
\mathcal B_1 =\left[ 
\right]_{{u:-1}}
\end{equation}

In $\mathcal B_2$, separate the $p_V(x_3,-u)$ term from the difference
quotient. The latter gives\nopagebreak[4]

\begin{equation}\label{eq:nb-eval-p}
-\left[ %
\right]_{{u:-1}}
\end{equation}

The product term $\mathcal B_3$ gives

\begin{equation}\label{eq:nb-eval-pq}
\mathcal B_3 =\left[ %
\right]_{{u:-1}}
\end{equation}

To combine the remaining part of $\mathcal B_2$ with
$\mathcal B_4$, note that the sign of the second term below follows from
$l_{1,3}(r(x_1,x_3)-1)/(x_1-x_3)=-1$ and the minus sign in
$\mathcal B_4$. We then use
$-(u-x_3)^{-1}+(u-x_1)^{-1}
=(x_1-x_3)((u-x_1)(u-x_3))^{-1}$ to obtain

\begin{align}
&{}-\left[ %
\right]_{{u:-1}}\label{eq:nb-eval-rest}
\end{align}

Adding \eqref{eq:nb-eval-q}--\eqref{eq:nb-eval-rest} gives
\eqref{eq:nb-bubble-sum}. Substituting the expressions for $A_V$ and $B_V$
yields \eqref{eq:nb-bubble-expanded}.
\end{proof}

\subsection{Source diagrams}
\label{27terms-app}
The 27 unnormalized summands on each side are paired by their common source index and grouped by crossing number. As in the main proof, $1$, $2$, and $3$ designate the black-crossing, parallel-strand, and cup--cap choices, respectively. Thus the four groups contain $1$, $6$, $12$, and $8$ terms on each side. The normalized reductions above use the corresponding $\widetilde L_{ijk}$ and $\widetilde R_{ijk}$.

\subsubsection{Three crossings}

\leavevmode\par\nobreak
\noindent\begin{minipage}[c]{0.49\linewidth}\raggedright
\(L_{111} =
\)
\end{minipage}\par\smallskip

\subsubsection{Two crossings}
\leavevmode\par\nobreak
\noindent\begin{minipage}[c]{0.49\linewidth}\raggedright
\(L_{211} =%
\)
\end{minipage}\par\smallskip

\subsubsection{One crossing}
\leavevmode\par\nobreak
\noindent\begin{minipage}[c]{0.49\linewidth}\raggedright
\(L_{221} =%
\)
\end{minipage}\par\smallskip

\subsubsection{No crossings}
\leavevmode\par\nobreak
\noindent\begin{minipage}[c]{0.49\linewidth}\raggedright
\(L_{222} =%
\)
\end{minipage}\par\smallskip

\subsection{Rational-function calculations}
\label{app:maple-checks}
The supplementary Maple file \texttt{appendix\_A\_checks.mpl}
contains the scalar checks for
Lemmas~\ref{lem:appA-one-crossing-ss}--\ref{lem:appA-one-crossing-se}.
It substitutes
$M_X=B_X^{-1}$, $W_X=A_X^{-1}$,
$A_X(-u)B_X(u)=(1+2u)^{-1}$ and
$B_N(u)=(1+u-x_i)(1+u+x_i)^{-1}B_V(u)$,
where $x_i$ is the marked through-strand variable,
and simplifies the coefficient differences to zero.
The identities are formed in the normalized formal series;
Lemma~\ref{lem:finite-dd-evaluation} gives their dot evaluations.

\section{Coefficients of zero-crossing terms}
\label{app:coefficient-certificates}

We record the coefficients obtained from the graphical
reductions of Appendix~\ref{app:graphical-certificates} and the
identities needed to compare them. The routine algebraic
simplifications of the displayed expressions are omitted.
Throughout this appendix, $C_L(d)$ and $C_R(d)$ are the full coefficients
of the marked diagram $d$ in $\tilde L$ and $\tilde R$, respectively,
including all bubble contributions.

We retain the notation of Definitions~\ref{ration111} and
\ref{notation-1}. All sums include the correction terms from
Appendix~\ref{app:graphical-certificates}.
Divided differences are formed before nilpotent dot evaluation, as in
Lemma~\ref{lem:finite-dd-evaluation}. A dot in an argument indicates
the variable in which the divided difference is taken; all other
arguments are held fixed, even when they coincide with its nodes. Thus
\begin{align*}
 p_N(x_2,\cdot)[x_1,x_3]
 &=\frac{p_N(x_2,x_1)-p_N(x_2,x_3)}{x_1-x_3},\\
 f[x_1,x_2,x_3]
 &=\frac{f[x_1,x_2]-f[x_2,x_3]}{x_1-x_3}.
\end{align*}
{Repeated nodes are interpreted by regular extension; in particular,
$f[a,a]=f'(a)$ and $f[a,a,a]=f''(a)/2$.
No difference of nilpotent dots is inverted.} Auxiliary functions are local to their subsections.
For these evaluations we use the normalized formal series preceding
Lemma~\ref{lem:MW-divided-difference} and a common sufficiently high truncation.
Finite representatives are retained until coefficient extraction at
infinity is complete.

\subsection{Zero-crossing coefficients}
\label{app:zero-crossing-coefficients}
\subsubsection{The $e_1e_2$-term}
\label{app:first-zero-coefficient}
For $d=e_1e_2$,
we choose the positions of dots as
$\begin{tikzpicture}[baseline=(current bounding box.center),
  scale=0.8,line cap=round,line join=round,
  font={\fontsize{10pt}{12pt}\selectfont}]
  \begin{scope}[yscale=1]
  \draw[-,thick] (0.45,.6) to (-0.45,-.6);
  \node at (.5,.5){$\scriptstyle 2$};
  \node at (.35,.4){$\scriptstyle\bullet$};
  \node at (.1,.5){$\scriptstyle 1$};
  \node at (0,.4){$\scriptstyle\bullet$};
  \node at (.6,-.5){$\scriptstyle 3$};
  \node at (.4,-.4){$\scriptstyle\bullet$};
  \draw[-,thick] (0,-.6) to[out=90,in=180] (.25,-.2);
  \draw[-,thick] (0.25,-.2) to[out=0,in=90] (.45,-0.6);
  \draw[-,thick] (-.45,.6) to[out=-90,in=180] (-.2,0.2);
  \draw[-,thick] (-.2,0.2) to[out=0,in=-90] (0,0.6);
  \end{scope}
\end{tikzpicture}$.
The lower and upper boundary variables are
\((x_2,-x_3,x_3)\) and \((-x_1,x_1,x_2)\), respectively.
Set
\[
\begin{aligned}
 \Phi(u)&=M_V(-u)q_V(x_2,u),
 &F(u)&=r(x_2,u)q_V(-u,x_3),\\
 h(u)&=q_V(x_2,u)A_V(u)(1+x_2-u).
\end{aligned}
\]

\begin{Lemma}\label{CRCLE}
{The coefficients satisfy $C_L(e_1e_2)=C_R(e_1e_2)$, where}
\begin{align*}
C_L(e_1e_2)={}&1-2\ell_{-+}^{12}+(1-\ell_{+-}^{23})\Phi(x_3)-\ell_{++}^{22}\Phi[-x_2,x_3]\\
 &+\ell_{-+}^{12}W_N(-x_1)\bigl\{(1-2x_3)p_N(x_1,x_3)
       \\ &\qquad +p_N(x_1,-x_2)\bigl(\Phi(x_3)\ell_{+-}^{23}-1-2h(x_3)+\ell_{++}^{22}h[-x_2,x_3]\bigr)\\
 &\qquad+2p_N(x_1,\cdot)[-x_2,x_3]-\ell_{++}^{22}p_N(x_1,\cdot)[-x_2,-x_2,x_3]\bigr\} ,
\end{align*}
and
\begin{align*}
C_R(e_1e_2)={}&\ell_{++}^{11}-2\\ &+(\ell_{++}^{11}-2)A_N(-x_2)q_N(x_1,-x_2)\ell_{++}^{23}(1-\Phi(x_3)\ell_{+-}^{23})\\ &+(1-2x_1)M_N(x_2)q_N(x_1,-x_2)\\
 &+\ell_{++}^{22}q_N(\cdot,-x_2)[x_1,x_2]\bigl(M_N(x_2)-A_N(-x_2)\ell_{++}^{23}(1-\Phi(x_3)\ell_{+-}^{23})\bigr)\\
 &+\ell_{++}^{23}\ell_{+-}^{23}\bigl\{M_V(-x_3)\bigl(\ell_{++}^{11}F(-x_1)+2F[-x_1,-x_2]-\ell_{++}^{22}F[-x_1,-x_2,-x_2]\bigr)\\
 &\qquad-\ell_{++}^{11}r(x_2,\cdot)[-x_1,x_3]+2r(x_2,\cdot)[-x_2,x_3]\\ &\qquad -2r(x_2,\cdot)[-x_1,-x_2,x_3]
             +\ell_{++}^{22}r(x_2,\cdot)[-x_1,-x_2,-x_2,x_3]-2\Phi(x_3)/\ell_{++}^{22}\bigr\}.
\end{align*}

\end{Lemma}

\begin{proof}
The nested differences in the graphical reduction give the
repeated-node brackets displayed above. The normalized series satisfy
\[
 \widehat A_X(u)\widehat B_X(-u)=\frac{1}{1-2u},\qquad
 \widehat B_N(u)=\frac{1+u-x_2}{1+u+x_2}\widehat B_V(u),
\]
by Lemmas~\ref{lem:MW-divided-difference}, \ref{ou}(2) and
\ref{gamma-move1}; the coefficientwise action of $E_0$ is understood.
Substituting Definition~\ref{notation-1} into the displayed coefficients
and using these identities, together with
$M_X=\widehat B_X^{-1}$ and $W_X=\widehat A_X^{-1}$, gives
$C_L(e_1e_2)-C_R(e_1e_2)=0$.
Lemma~\ref{lem:finite-dd-evaluation} then permits the common nilpotent
evaluation, including repeated nodes.
\end{proof}

\subsubsection{ The $e_2e_1$-term}
\label{app:zero-coeff-2} For $d=e_2e_1$,
we choose the positions of dots as follows
$\begin{tikzpicture}[baseline=(current bounding box.center),
  scale=0.8,line cap=round,line join=round,
  font={\fontsize{10pt}{12pt}\selectfont}]
  \begin{scope}[yscale=-1]
  \draw[-,thick] (0.45,.6) to (-0.45,-.6);
  \node at (.5,.5){$\scriptstyle 2$};
  \node at (.3,.4){$\scriptstyle\bullet$};
  \node at (-.3,.5){$\scriptstyle 1$};
  \node at (-.4,.4){$\scriptstyle\bullet$};
  \node at (.3,-.5){$\scriptstyle 3$};
  \node at (0,-.4){$\scriptstyle\bullet$};
  \draw[-,thick] (0,-.6) to[out=90,in=180] (.25,-.2);
  \draw[-,thick] (0.25,-.2) to[out=0,in=90] (.45,-0.6);
  \draw[-,thick] (-.45,.6) to[out=-90,in=180] (-.2,0.2);
  \draw[-,thick] (-.2,0.2) to[out=0,in=-90] (0,0.6);
  \end{scope}
\end{tikzpicture}$.
The lower and upper boundary variables are
\((x_1,-x_1,x_2)\) and \((x_2,x_3,-x_3)\), respectively.
Set
\[
 J(u)=r(u,x_2)q_N(u,-x_1),\qquad
 \Phi(u)=M_N(-u)q_N(-x_2,u).
\]

\begin{Lemma}\label{lem:e2e1-coefficients}
{The coefficients satisfy $C_L(e_2e_1)=C_R(e_2e_1)$, where}
\begin{align*}
C_L(e_2e_1)={}&(\ell_{--}^{33}-2)(1+\ell_{--}^{12}A_V(x_2)q_V(-x_3,x_2)-M_V(-x_2)q_V(-x_3,x_2))
       \\ &+\ell_{--}^{12}\ell_{+-}^{12}M_N(x_1)q_N(-x_2,-x_1)\bigl((2-\ell_{--}^{33})A_V(x_2)q_V(-x_3,x_2)-2/\ell_{--}^{22}\bigr)\\ &+\ell_{--}^{33}\ell_{--}^{12}\ell_{+-}^{12}M_N(x_1)J(-x_3)\\
 &+\ell_{--}^{33}\ell_{--}^{12}\ell_{+-}^{12}r(\cdot,x_2)[-x_3,x_1]-2\ell_{--}^{12}\ell_{+-}^{12}r(\cdot,x_2)[x_1,-x_2]
       \\ &+\ell_{--}^{22}q_V(\cdot,x_2)[-x_2,-x_3]\bigl(M_V(-x_2)-\ell_{--}^{12}A_V(x_2)\\ &\qquad +\ell_{--}^{12}\ell_{+-}^{12}A_V(x_2)M_N(x_1)q_N(-x_2,-x_1)\bigr)\\
 &-2\ell_{--}^{12}\ell_{+-}^{12}M_N(x_1)J[-x_2,-x_3]-\ell_{--}^{12}\ell_{+-}^{12}\ell_{--}^{22}M_N(x_1)J[-x_2,-x_2,-x_3]\\
 &-2\ell_{--}^{12}\ell_{+-}^{12}r(\cdot,x_2)[-x_2,-x_3,x_1]-\ell_{--}^{12}\ell_{+-}^{12}\ell_{--}^{22}r(\cdot,x_2)[-x_2,-x_2,-x_3,x_1] ,
\end{align*}
and
\begin{align*}
C_R(e_2e_1)={}&1-2\ell_{-+}^{23}\\ &+
       \bigl((x_1-x_2)\Phi(-x_1)+\ell_{--}^{22}\Phi[-x_1,x_2]\bigr)\bigl(\ell_{-+}^{23}W_V(x_3)p_V(-x_3,x_2)-1\bigr)\\
 &+\ell_{-+}^{23}W_V(x_3)\bigl\{(1+2x_1)p_V(-x_3,-x_1)-p_V(-x_3,x_2)
               \\ &\qquad +2p_V(-x_3,\cdot)[-x_1,x_2]-\ell_{--}^{22}p_V(-x_3,\cdot)[-x_1,x_2,x_2]\bigr\}.
\end{align*}

\end{Lemma}

\begin{proof}
Substitution of Definition~\ref{notation-1} and the normalized
identities in the proof of Lemma~\ref{CRCLE} gives
$C_L(e_2e_1)=C_R(e_2e_1)$.
The repeated-node divided differences are evaluated using
Lemma~\ref{lem:finite-dd-evaluation}.
\end{proof}

\subsubsection{The $e_1$-term}
\label{app:zero-coeff-3}
For $d=e_1$, we choose the positions of dots as
$
\begin{tikzpicture}
[baseline=(current bounding box.center),
  scale=0.8,line cap=round,line join=round,
  font={\fontsize{10pt}{12pt}\selectfont}]
  \begin{scope}[xscale=1,yscale=1]
  \draw[-,thick] (0,-0.6) to[out=90, in=0] (-.2,-0.3);
  \draw[-,thick] (-.2,-0.3) to[out = 180, in = 90] (-0.4,-.6);
  \node at (.1,.5){$\scriptstyle 2$};
  \node at (0,.4){$\scriptstyle\bullet$};
  \node at (.1,-.5){$\scriptstyle 1$};
  \node at (-.1,-.4){$\scriptstyle\bullet$};
  \node at (.6,-.5){$\scriptstyle 3$};
  \node at (.4,-.4){$\scriptstyle\bullet$};
  \draw[-,thick] (.4,.6) to (.4,-.6);
  \draw[-,thick] (0,.6) to[out=-90, in=0] (-.2,0.3);
  \draw[-,thick] (-.2,0.3) to[out = -180, in = -90] (-0.4,.6);
  \end{scope}
\end{tikzpicture}$.
The variables \(x_1,x_2\) on the two arcs are independent.
The strand variables at the upper and lower boundaries are
\((-x_2,x_2,x_3)\) and \((-x_1,x_1,x_3)\), respectively.
Define
\begin{align*}
\Phi(u)&:=M_V(-u)q_V(x_3,u),
&J(u)&:=r(u,x_3)q_N(u,x_1),\\
U(u)&:=p_N(x_2,u)q_N(u,x_1),\\
G(u)&:=\frac{q_V(u,x_1)M_V(-x_1)-q_V(u,x_3)M_V(-x_3)}
              {d_{+-}^{13}}.
\end{align*}
In this subsection, set
\[
g(u):=
\frac{(1-2u)\widehat A_N(u)\widehat B_N(u)-1}{u}
=\frac{\widehat B_N(u)-\widehat B_N(-u)}
       {u\widehat B_N(-u)}.
\]
The series \(g\) is regular at zero.

\begin{Lemma}\label{lem:e1-coefficients}
{We have $C_L(e_1)=L_3+B_t=R_3=C_R(e_1)$.
Here $L_3,R_3$ are the contributions without bubbles and $B_t$
is the bubble contribution on the left, given by}
\begin{align*}
L_3={}&
 \ell_{-+}^{23}\ell_{+-}^{13}/\ell_{+-}^{23}-\ell_{-+}^{23}+2\\
&+\ell_{--}^{13}\ell_{+-}^{13}M_N(-x_1)q_N(x_2,x_1)\\
&+\ell_{-+}^{23}W_N(-x_2)(2+\ell_{--}^{13})p_N(x_2,x_1)\\
 &+\ell_{-+}^{23}\ell_{--}^{13}\ell_{+-}^{13} M_N(-x_1)W_N(-x_2)\bigl\{A_V(x_3)q_V(x_3,x_3)U(-x_3)\\
&\qquad +U(-x_3)/(1-2x_3)+U(x_3)-U[x_3,-x_3]\bigr\}\displaybreak[1]\\
 &-\ell_{-+}^{23}\ell_{+-}^{13}W_N(-x_2)p_N(x_2,-x_3)/(1-2x_3)\\
&+\ell_{-+}^{23}\ell_{--}^{13}\ell_{+-}^{13} M_N(-x_1)\bigl(q_N(\cdot,x_1)[-x_2,x_3]-q_N(\cdot,x_1)[x_2,-x_2,x_3]\bigr)\\
 &+\ell_{-+}^{23}\ell_{+-}^{13}W_N(-x_2)\bigl\{p_N(x_2,-x_1)+p_N(x_2,\cdot)[-x_1,x_3]-p_N(x_2,\cdot)[x_1,-x_1]\\
&\qquad -p_N(x_2,\cdot)[x_1,-x_1,x_3]
                    -p_N(x_2,x_3)+p_N(x_2,\cdot)[x_3,-x_3]\bigr\}\displaybreak[1]\\
 &+\ell_{-+}^{23}W_N(-x_2)p_N(x_2,-x_3)\bigl((2\ell_{--}^{13}-\ell_{-+}^{13})A_V(x_3)q_V(x_3,x_3)\\
&\qquad -\ell_{-+}^{13}\Phi[x_1,x_3]\bigr)\\
 &-4M_V(-x_3)q_V(x_3,x_3)+2\ell_{-+}^{13}\Phi[x_1,x_3]-2\ell_{-+}^{23}W_N(-x_2)p_N(x_2,\cdot)[x_3,-x_3]
       \\
&+2\ell_{+-}^{13}\ell_{--}^{13}r(\cdot,x_3)[x_2,-x_1]\\
&+\ell_{-+}^{13}\ell_{-+}^{23}W_N(-x_2)p_N(x_2,\cdot)[-x_3,x_1,x_3]\\
 &+2\ell_{+-}^{13}A_V(x_3)q_V(x_2,x_3)+\ell_{+-}^{13}\ell_{-+}^{23}A_V(x_3) q_V(\cdot,x_3)[x_2,x_3]\\
&-2\ell_{+-}^{13}\ell_{--}^{13}M_N(-x_1)A_V(x_3)q_V(x_3,x_3)q_N(-x_3,x_1)
                   \\
&-\ell_{+-}^{13}\ell_{--}^{13}\ell_{+-}^{23}M_N(-x_1)A_V(x_3)q_N(-x_3,x_1)q_V(\cdot,x_3)[x_2,x_3]\\
 &-2\ell_{+-}^{23}M_V(-x_3)q_V(\cdot,x_3)[x_2,x_3]+\ell_{+-}^{23}\ell_{-+}^{13}G[x_2,x_3]\\
&+2\ell_{+-}^{13}\ell_{--}^{13}M_N(-x_1)J[x_3,-x_3]+\ell_{+-}^{13}\ell_{--}^{13}\ell_{+-}^{23}M_N(-x_1)J[-x_3,x_2,x_3]\\
 &+\ell_{+-}^{13}\ell_{--}^{13}\ell_{+-}^{23}r(\cdot,x_3)[-x_3,x_2,x_3,-x_1]+2\ell_{+-}^{13}\ell_{--}^{13}r(\cdot,x_3)[x_3,-x_3,-x_1],
\end{align*}
and
\begin{align*}
R_3={}&\ell_{++}^{23}
 -2x_2\ell_{--}^{13}M_N(-x_1)q_N(x_2,x_1),
\end{align*}
and %
\[
B_t=\frac12\ell_{-+}^{23}\ell_{--}^{13}\ell_{+-}^{13}
M_N(-x_1)W_N(-x_2)g[-x_1,-x_2,x_3].
\]
Here $g[\cdots]$ denotes the second divided difference of $g$.
\end{Lemma}

\begin{proof}
By Lemmas~\ref{ou}(2) and \ref{gamma-move1},
\[
 \widehat B_N(u)=\frac{1+u-x_3}{1+u+x_3}\widehat B_V(u),\qquad
 \widehat A_N(u)=\frac{1-u+x_3}{1-u-x_3}\widehat A_V(u),
\]
with the  coefficientwise action of $E_0$ understood.
Direct substitution of Definition~\ref{notation-1}, these
transition identities, and
$\widehat A_X(u)\widehat B_X(-u)=(1-2u)^{-1}$ gives
\[
 L_3-R_3=-\frac12\ell_{-+}^{23}\ell_{--}^{13}\ell_{+-}^{13}
 M_N(-x_1)W_N(-x_2)g[-x_1,-x_2,x_3].
\]

{For the bubble, keep the finite representatives $A_N,B_N$.
Their product has the form
\[
 A_N(u)B_N(u)=\frac{(-1)^{c+1}}2m_{N,0}(u)f_N^+(u)f_N^-(u)
 \bigl(1+O(u^k)\bigr),
\]
where $O(u^k)$ is a polynomial divisible by $u^k$, arising from the
truncations of $\Gamma_N$ and $\Gamma_N^{-1}$.}
The factorization in the proof of Lemma~\ref{nbone:lemma:R4} yields
\[
 A_N(u)B_N(u)\frac{O_N(u)}u
 -\frac{(-1)^{c+1}}2\frac{t_{N,0}(u)}{u^{\varepsilon_0(N)+1}}
 \in u^{k-\varepsilon_0(N)-1}Z_{\lambda,N}[u].
\]
Since $\deg t_{N,0}=\varepsilon_0(N)$, the rational term has no
$u^{-1}$-coefficient after division by
$(u+x_1)(u+x_2)(u-x_3)$.
The polynomial remainder contributes its second divided difference,
which vanishes at the nilpotent dots for sufficiently large $k$ by
Lemma~\ref{lem:finite-dd-evaluation}.
As $A_N(0)B_N(0)=1$, the remaining terms give
\[
 \left[
 \frac{A_N(u)B_N(u)\bigl(O_N(u)/u-1+1/(2u)\bigr)}
 {(u+x_1)(u+x_2)(u-x_3)}
 \right]_{{u:-1}}
 =\frac12g[-x_1,-x_2,x_3].
\]
Indeed, $-A_NB_N+A_NB_N/(2u)$ differs from
$((1-2u)A_NB_N-1)/(2u)$ by $1/(2u)$, which contributes nothing.
The same lemma allows replacement by the normalized series in this
divided difference. Multiplying by the exterior factors in the bubble
diagram gives $B_t$ as stated. Hence
$C_L(e_1)-C_R(e_1)=L_3+B_t-R_3=0$.
\end{proof}

\subsubsection{The $e_2$-term}
\label{app:zero-coeff-4}
For $d=e_2$, we choose the positions of dots as     $
\begin{tikzpicture}
[baseline=(current bounding box.center),
  scale=0.8,line cap=round,line join=round,
  font={\fontsize{10pt}{12pt}\selectfont}]
  \begin{scope}[xscale=-1,yscale=1]
  \draw[-,thick] (0,-0.6) to[out=90, in=0] (-.2,-0.3);
  \draw[-,thick] (-.2,-0.3) to[out = 180, in = 90] (-0.4,-.6);
  \node at (-.5,.5){$\scriptstyle 2$};
  \node at (-0.3,.4){$\scriptstyle\bullet$};
  \node at (-.4,-.5){$\scriptstyle 1$};
  \node at (-.3,-.4){$\scriptstyle\bullet$};
  \node at (.6,-.5){$\scriptstyle 3$};
  \node at (.4,-.4){$\scriptstyle\bullet$};
  \draw[-,thick] (.4,.6) to (.4,-.6);
  \draw[-,thick] (0,.6) to[out=-90, in=0] (-.2,0.3);
  \draw[-,thick] (-.2,0.3) to[out = -180, in = -90] (-0.4,.6);
  \end{scope}
\end{tikzpicture}$.
The variables \(x_1,x_2\) on the two arcs are independent,
and \(x_3\) is the variable on the through strand.
The lower and upper boundary variables are
\((x_3,-x_1,x_1)\) and \((x_3,-x_2,x_2)\), respectively.
Define
\begin{align*}
\Phi(u)&:=M_N(-u)q_N(-x_3,u),
&J(u)&:=r(x_3,u)q_V(-u,x_1),\\
U(u)&:=p_V(x_2,u)q_V(u,x_1).
\end{align*}
In this subsection, set
\[
g(u):=
\frac{1-(1+2u)\widehat A_V(-u)\widehat B_V(-u)}{u}
=\frac{1-\widehat B_V(-u)/\widehat B_V(u)}{u}.
\]
The series \(g\) is regular at zero.

\begin{Lemma}\label{lem:e2-coefficients}
{We have $C_L(e_2)=L_4=R_4+B_r=C_R(e_2)$.
Here $L_4,R_4$ are the contributions without bubbles and $B_r$
is the bubble contribution on the right, given by}
\begin{align*}
L_4={}&\ell_{+-}^{23}
 +(1-2x_2\ell_{-+}^{13})M_V(-x_1)q_V(x_2,x_1)\\
&-\ell_{++}^{33}
 \frac{M_V(-x_1)q_V(x_2,x_1)-M_V(x_3)q_V(x_2,-x_3)}
      {d_{++}^{13}}
 +\ell_{--}^{33}M_V(x_3)
 \frac{q_V(x_2,-x_3)-q_V(-x_3,-x_3)}{d_{++}^{23}}\\
&-\ell_{--}^{33}\ell_{-+}^{13}\frac{1}{d_{++}^{13}}
 \left(
 M_V(-x_1)
 \frac{q_V(x_2,x_1)-q_V(-x_3,x_1)}{d_{++}^{23}}
 -M_V(x_3)
 \frac{q_V(x_2,-x_3)-q_V(-x_3,-x_3)}{d_{++}^{23}}
 \right).
\end{align*}
and
\begin{align*}
R_4={}&
 \ell_{--}^{23}\ell_{++}^{13}/\ell_{++}^{23}-\ell_{--}^{23}+2\\
&+\ell_{++}^{13}\ell_{-+}^{13}M_V(-x_1)q_V(x_2,x_1)\\
&+\ell_{--}^{23}W_V(-x_2)(2+\ell_{-+}^{13})p_V(x_2,x_1)\\
 &+\ell_{--}^{23}\ell_{++}^{13}\ell_{-+}^{13} M_V(-x_1)W_V(-x_2)\bigl\{A_N(-x_3)q_N(-x_3,-x_3)U(x_3)\\
&\qquad +U(x_3)/(1+2x_3)+U(-x_3)-U[x_3,-x_3]\bigr\}\displaybreak[1]\\
 &-\ell_{--}^{23}\ell_{++}^{13}W_V(-x_2)A_N(-x_3)p_V(x_2,x_3)q_N(-x_3,-x_3)\\
&+(2\ell_{++}^{13}A_N(-x_3)-3M_N(x_3))q_N(-x_3,-x_3)\\
 &+\ell_{++}^{23}q_N(\cdot,-x_3)[x_2,-x_3]\bigl(\ell_{++}^{13}A_N(-x_3)\\
&\qquad -\ell_{++}^{13}\ell_{-+}^{13}M_V(-x_1)A_N(-x_3)q_V(x_3,x_1)-M_N(x_3)\bigr)
                   \\
&-2\ell_{++}^{13}\ell_{-+}^{13}M_V(-x_1)A_N(-x_3)q_N(-x_3,-x_3)q_V(x_3,x_1)\\
&+2\ell_{--}^{23}W_V(-x_2)M_N(x_3)p_V(x_2,x_3)q_N(-x_3,-x_3)\\
 &+\ell_{--}^{13}\bigl(1-\ell_{--}^{23}W_V(-x_2)p_V(x_2,x_3)\bigr)\Phi[x_1,-x_3]\\
&-2\ell_{--}^{23}W_V(-x_2)p_V(x_2,\cdot)[x_3,-x_3]+\ell_{--}^{23}\ell_{--}^{13}W_V(-x_2)p_V(x_2,\cdot)[x_3,x_1,-x_3]\\
 &+\ell_{--}^{23}\ell_{++}^{13}W_V(-x_2)\bigl\{p_V(x_2,-x_1)+p_V(x_2,\cdot)[-x_3,-x_1]-p_V(x_2,\cdot)[x_1,-x_1]\\
&\qquad -p_V(x_2,\cdot)[-x_3,-x_1,x_1]
                   -p_V(x_2,-x_3)+p_V(x_2,\cdot)[x_3,-x_3]\bigr\}\displaybreak[1]\\
 &+\ell_{--}^{23}\ell_{++}^{13}\ell_{-+}^{13} M_V(-x_1)\bigl(q_V(\cdot,x_1)[-x_2,-x_3]-q_V(\cdot,x_1)[x_2,-x_2,-x_3]\bigr)\\
 &-2\ell_{++}^{13}\ell_{-+}^{13}r(x_3,\cdot)[x_1,-x_2]-2\ell_{++}^{13}\ell_{-+}^{13}M_V(-x_1)J[x_3,-x_3]\\
&+\ell_{++}^{13}\ell_{-+}^{13}\ell_{++}^{23}M_V(-x_1)J[-x_3,x_3,-x_2]\\
 &-\ell_{++}^{13}\ell_{-+}^{13}\ell_{++}^{23}r(x_3,\cdot)[-x_3,x_3,-x_2,x_1]\\
&+\bigl(2\ell_{++}^{13}\ell_{-+}^{13}-\ell_{--}^{23}\ell_{++}^{13}\ell_{-+}^{13} W_V(-x_2)p_V(x_2,x_3)\bigr)r(x_3,\cdot)[x_3,-x_3,x_1],
\end{align*}
and
\[
B_r=\frac12\ell_{--}^{23}\ell_{++}^{13}\ell_{-+}^{13}
M_V(-x_1)W_V(-x_2)g[x_1,x_2,x_3].
\]
Here $g[\cdots]$ denotes the second divided difference of $g$.
\end{Lemma}

\begin{proof}
Direct substitution of Definition~\ref{notation-1} and the
normalized identities in the proof of
Lemma~\ref{lem:e1-coefficients} gives
\[
 L_4-R_4=\frac12\ell_{--}^{23}\ell_{++}^{13}\ell_{-+}^{13}
 M_V(-x_1)W_V(-x_2)g[x_1,x_2,x_3].
\]
For the bubble, retain the finite representatives until coefficient
extraction. Writing $b_V(u)=O_V(u)/u-1+1/(2u)$, the left-leg resolvent
contributes $-b_V(-u)$ with multiplier $-A_V(-u)B_V(-u)$.
The two minus signs cancel.
The argument in the proof of Lemma~\ref{lem:e1-coefficients}, now
with $V$ in place of $N$ and with nodes $x_1,x_2,x_3$, gives
\[
 \left[
 \frac{A_V(-u)B_V(-u)b_V(-u)}
 {(u-x_1)(u-x_2)(u-x_3)}
 \right]_{{u:-1}}
 =\frac12g[x_1,x_2,x_3].
\]
Here the surviving polynomial part is
$\bigl(1-(1+2u)A_V(-u)B_V(-u)\bigr)/(2u)$;
the remaining $-1/(2u)$ contributes nothing.
This yields the stated $B_r$, and hence
$C_L(e_2)-C_R(e_2)=L_4-R_4-B_r=0$.
\end{proof}

\section{A rational-function identity for the mixed braid relation}
\begin{Lemma}\label{compu-T}
{Let $V$ and $E=E_kV$ be nonzero finitely generated objects,
where $k\in\mathbb N_{>0}$, and put $\delta=\delta_{k,1}$.
With the dot positions used in the proof of \eqref{fw:mixed-braid}, we have
$T(x_1,x_2,x_3)=-\delta A_E(x_1)B_V(-x_3)$, where}
\[
\begin{aligned}
T(x_1,x_2,x_3)
={}&-2\frac{1}{\ell_{--}^{11}}
       \frac{d_{-+}^{23}}{\ell_{+-}^{23}}
       \frac{d_{--}^{23}}{\ell_{++}^{23}}
       \frac{1}{\ell_{--}^{33}}
+2\frac{1}{\ell_{--}^{11}}
       \frac{d_{-+}^{23}}{\ell_{+-}^{23}}
       \frac{d_{--}^{23}}{\ell_{++}^{23}}
       \frac{1}{\ell_{--}^{33}}
       \frac{1}{(d_{+-}^{23})^2}
-\frac{1}{d_{+-}^{12}}
       \frac{1}{\ell_{--}^{11}}
       \frac{d_{--}^{12}}{\ell_{++}^{12}}
       \frac{d_{+-}^{12}}{\ell_{-+}^{12}}
       \frac{1}{(d_{+-}^{23})^2}
\\
&-\frac{1}{\ell_{+-}^{23}}
       \frac{1}{\ell_{--}^{11}}
       \frac{d_{--}^{12}}{\ell_{++}^{12}}
       \frac{1}{\ell_{-+}^{12}}
+\frac{1}{\ell_{+-}^{23}}
       \frac{1}{\ell_{--}^{11}}
       \frac{d_{--}^{12}}{\ell_{++}^{12}}
       \frac{1}{\ell_{-+}^{12}}
       \frac{1}{(d_{+-}^{23})^2}
+\frac{d_{-+}^{23}}{\ell_{+-}^{23}}
       \frac{1}{\ell_{++}^{23}}
       \frac{1}{\ell_{--}^{11}}
       \frac{1}{\ell_{++}^{12}}
\\
&-\frac{d_{-+}^{23}}{\ell_{+-}^{23}}
       \frac{1}{\ell_{++}^{23}}
       \frac{1}{\ell_{--}^{11}}
       \frac{1}{\ell_{++}^{12}}
       \frac{1}{(d_{+-}^{23})^2}
+\frac{1}{d_{+-}^{23}}
       \frac{1}{\ell_{--}^{11}}
       \frac{d_{--}^{12}}{\ell_{++}^{12}}
       \frac{d_{+-}^{12}}{\ell_{-+}^{12}}
       \frac{1}{(d_{+-}^{12})^2}+\bigl(\ell_{+-}^{32}\bigr)^\delta
       B_V(-x_3)p_E(-x_3,x_1)\\
       &
-\bigl(\ell_{+-}^{12}\bigr)^\delta
       A_E(x_1)q_V(-x_3,x_1).
\end{aligned}
\]
\end{Lemma}
\begin{proof}
{Use Definitions~\ref{ration111} and \ref{notation-1}.
Identities involving $A_X,B_X$ are understood in the normalized
formal series of Section~4. Lemmas~\ref{ou}(2) and
\ref{fw:Gamma-positive}(1) give}
\begin{equation}\label{appC:shift}
A_E(u)=A_V(u)\frac{u+x_2}{u-x_2}
                   \frac{(u-x_2+1)^{1-\delta}}{u+x_2+1}.
\end{equation}
{By Lemma~\ref{lem:MW-divided-difference},
$A_V(u)B_V(-u)=(1-2u)^{-1}$. Thus}
\[
{A_E(u)B_V(-u)=H_\delta(u):=
\frac{(u+x_2)(u-x_2+1)^{1-\delta}}
     {(u-x_2)(u+x_2+1)(1-2u)}.}
\]
{Since $(1+u-x_2)^\delta H_\delta(u)=H_0(u)$, the definitions
of $p_E,q_V$ give}
\[
{\begin{aligned}
&(\ell_{+-}^{32})^\delta B_V(-x_3)p_E(-x_3,x_1)
-(\ell_{+-}^{12})^\delta A_E(x_1)q_V(-x_3,x_1)\\
&\qquad=-\delta A_E(x_1)B_V(-x_3)
       +\frac{H_0(x_1)-H_0(x_3)}{x_1-x_3}.
\end{aligned}}
\]
{The first eight terms of $T$ sum to
$-(H_0(x_1)-H_0(x_3))/(x_1-x_3)$, as follows by putting them
over a common denominator. This proves the identity.
Here $x_1,x_3$ and $x_2-k$ are nilpotent. The denominators of
$H_0$ are units at $u=0$, $x_2=k$; its divided difference is
formed before evaluating $x_1,x_3$, by
Lemma~\ref{lem:finite-dd-evaluation}.}
\end{proof}


\begin{thebibliography}{BSWW18}
\setlength{\baselineskip}{14.2pt}

\bibitem[Ari96]{Ariki-dec}
{\scshape S. Ariki},
On the decomposition numbers of the Hecke algebra of $G(m,1,n)$,
{\it J. Math. Kyoto Univ.} \textbf{36} (1996), no. 4, 789--808.

\bibitem[AMR06]{AMR-cyc}
{\scshape S. Ariki, A. Mathas and H. Rui},
Cyclotomic Nazarov--Wenzl algebras,
{\it Nagoya Math. J.} \textbf{182} (2006), 47--134.

\bibitem[BSWW18]{BSWW-icate}
{\scshape H. Bao, P. Shan, W. Wang and B. Webster},
Categorification of quantum symmetric pairs I,
{\it Quantum Topol.} \textbf{9} (2018), no. 4, 643--714.

\bibitem[BW18]{BW-ihes}
{\scshape H. Bao and W. Wang},
A new approach to Kazhdan--Lusztig theory of type $B$ via quantum symmetric pairs,
{\it Ast\'erisque} \textbf{402} (2018).

\bibitem[BD17]{BD-cate}
{\scshape J. Brundan and N. Davidson},
Categorical actions and crystals,
{in {\it Categorification and higher representation theory},
{\it Contemp. Math.} \textbf{683} (2017), 105--147.}

\bibitem[BK09a]{BK-gr}
{\scshape J. Brundan and A. Kleshchev},
Graded decomposition numbers for cyclotomic Hecke algebras,
{\it Adv. Math.} \textbf{222} (2009), no. 6, 1883--1942.

\bibitem[BK09b]{BK-klr}
{\scshape J. Brundan and A. Kleshchev},
Blocks of cyclotomic Hecke algebras and Khovanov--Lauda algebras,
{\it Invent. Math.} \textbf{178} (2009), {no. 3,} 451--484.

\bibitem[BSW20]{BSW-heisenberg}
{\scshape J. Brundan, A. Savage and B. Webster},
Heisenberg and Kac--Moody categorification,
{\it Selecta Math.} \textbf{26} (2020), no. 5, Paper No. 74.

\bibitem[BSW26a]{BSW}
{\scshape J. Brundan, A. Savage and B. Webster},
An update on Heisenberg and Kac--Moody categorification,
{\it J. Pure Appl. Algebra} \textbf{230} (2026), no. 8, Paper No. 108315.

\bibitem[BSW26b]{BSW-Brauer}
{\scshape J. Brundan, A. Savage and B. Webster},
Brauer and iquantum categorification,
\href{https://arxiv.org/abs/2610.25922v1}{arXiv:2610.25922v1}, 2026.

\bibitem[BS24]{BS-semi}
{\scshape J. Brundan and C. Stroppel},
Semi-infinite highest weight categories,
{\it Mem. Amer. Math. Soc.} \textbf{293} (2024), no. 1459, vii+152 pp.

\bibitem[BV22]{BV}
{\scshape J. Brundan and M. Vargas},
A new approach to the representation theory of the partition category,
{{\it J. Algebra} \textbf{601} (2022), 198--279.}

\bibitem[BWW24]{BWW24}
{\scshape J. Brundan, W. Wang and B. Webster},
The nil-Brauer category,
{\it Ann. Represent. Theory} \textbf{1} (2024), no. 1, 21--58.

\bibitem[BWW25]{BWW25}
{\scshape J. Brundan, W. Wang and B. Webster},
Categorification of quasi-split $\imath$-quantum groups,
\href{https://arxiv.org/abs/2505.22929v1}{arXiv:2505.22929v1}, 2025.

\bibitem[BWW26]{BWW25v2}
{\scshape J. Brundan, W. Wang and B. Webster},
Categorification of quasi-split $\imath$-quantum groups,
\href{https://arxiv.org/abs/2505.22929v2}{arXiv:2505.22929v2}, 2026.

\bibitem[ES16]{ES-grading}
{\scshape M. Ehrig and C. Stroppel},
Koszul gradings on Brauer algebras,
{\it Int. Math. Res. Not. IMRN} (2016), no. 13, 3970--4011.

\bibitem[ES18]{ES-NW}
{\scshape M. Ehrig and C. Stroppel},
Nazarov--Wenzl algebras, coideal subalgebras and categorified skew Howe duality,
{{\it Adv. Math.} \textbf{331} (2018), 58--142.}

\bibitem[GRS22]{GRS-Kaff}
{\scshape M. Gao, H. Rui and L. Song},
A basis theorem for the affine Kauffman category and its cyclotomic quotients,
{\it J. Algebra} \textbf{608} (2022), 774--846.

\bibitem[GRS23a]{GRS-tri}
{\scshape M. Gao, H. Rui and L. Song},
Representations of weakly triangular categories,
{\it J. Algebra} \textbf{614} (2023), 481--534.

\bibitem[GRS23b]{GRS-block}
{\scshape M. Gao, H. Rui and L. Song},
Blocks of the Brauer category over the complex field,
{{\it Algebr. Represent. Theory}} \textbf{26} (2023), {no. 5,} 1319--1334.

\bibitem[GRS26]{GRS}
{\scshape M. Gao, H. Rui and L. Song},
Affine and cyclotomic Brauer categorification via $\imath$-Kac--Moody
$2$-categories{:} the half-integral type $\operatorname{AIII}$ case,
{\href{https://arxiv.org/abs/2609.07018}{arXiv:2609.07018}}, 2026.

\bibitem[KL09]{KL-rings}
{\scshape M. Khovanov and A. D. Lauda},
A diagrammatic approach to categorification of quantum groups I,
{\it Represent. Theory} \textbf{13} (2009), 309--347.

\bibitem[R08]{R-2cat}
{\scshape R. Rouquier},
$2$-Kac--Moody algebras,
\href{https://arxiv.org/abs/0812.5023}{arXiv:0812.5023}, 2008.

\bibitem[RS19]{RS-cyc}
{\scshape H. Rui and L. Song},
Affine Brauer category and parabolic category $\mathcal O$ in types $B,C,D$,
{\it Math. Z.} \textbf{293} (2019), no. 1--2, 503--550.

\bibitem[RS20]{RS-brauer}
{\scshape H. Rui and L. Song},
Representations of Brauer category and categorification,
{\it J. Algebra} \textbf{557} (2020), 1--36.

\bibitem[RX09]{RX}
{\scshape H. Rui and J. Xu},
The representations of cyclotomic BMW algebras,
{\it J. Pure Appl. Algebra} \textbf{213} (2009), no. 12, 2262--2288.

\bibitem[SW26]{SW}
{\scshape A. Savage and B. Webster},
Bubbles in the affine Brauer and Kauffman categories,
{\it J. Aust. Math. Soc.} \textbf{120} (2026), no. 2, 243--286.

\end{thebibliography}
\end{document}